\documentclass[11pt, reqno]{amsart}
\usepackage{standalone}
\usepackage{subfiles}

\usepackage{bm}
\usepackage[hidelinks]{hyperref}
\usepackage{multirow}
\usepackage{tikz}
\usepackage{caption}
\usepackage{textgreek}
\usepackage{float}
\usetikzlibrary{
 decorations.pathmorphing,%
 decorations.markings,%
 decorations.pathreplacing,%
 decorations.text,%
 calc,%
 patterns,%
 shapes.geometric,%
 arrows,
 decorations.shapes,
 positioning,
 plotmarks,
 shadings,
 math,
 intersections
}
\usetikzlibrary{arrows.meta,decorations.markings,angles,quotes}

\tikzset{>=latex}
\tikzset{font=\small}
\tikzset{mark size=1.5pt, mark options=thin}
\tikzset{pin distance=4pt,
 every pin edge/.style={<-, thin, shorten <= -2pt}}
\usepackage{array,booktabs,tabularx}
\usepackage{graphicx}
\usepackage{amssymb}
\usepackage{arydshln}
\usepackage{enumerate}
\usepackage{color}
\usepackage{esint}
\usepackage{mathtools}
\usepackage{dsfont}
\usepackage{ esint }
\usepackage{placeins}
\usepackage[show]{ed}
\usepackage{mathrsfs}
\usepackage{changepage}   
\usepackage{quiver}
\usepackage{comment}
\usepackage{enumitem}
\usepackage{soul}
\usepackage{needspace}
\usepackage[normalem]{ulem}

\usepackage{multicol}
\usepackage{tabto}
\usepackage{hyperref}

\usepackage{bbm}

\usetikzlibrary{shapes.geometric, arrows}
\tikzstyle{box} = [rectangle, 
minimum width=2.75cm, 
minimum height=1cm, 
text centered, 
text width=2.75cm, 
draw=black]
\tikzstyle{impbox} = [rectangle, double,
minimum width=2.75cm, 
minimum height=1cm, 
text centered, 
text width=2.75cm, 
draw=black]

\definecolor{uipoppy}{RGB}{221,128,71}
\definecolor{uipaleblue2}{RGB}{179,196,215}
\definecolor{uiviolet}{RGB}{86,86,99}
\definecolor{uiblack}{RGB}{0, 0, 0}
\definecolor{azul}{RGB}{0,128,255}
\definecolor{verde}{RGB}{50,180,50}
\definecolor{pale-verde}{RGB}{155,207,145}
\definecolor{uipaleblue}{RGB}{108,199,220}
\definecolor{light-gray}{RGB}{198,198,198}

\newcommand{\Ell}{\operatorname{Ell}}
\newtheorem{lemma}{Lemma}
\newtheorem{theorem}[lemma]{Theorem}

\newtheorem{corollary}[lemma]{Corollary}
\newtheorem{proposition}[lemma]{Proposition}

\theoremstyle{definition}
\newtheorem{definition}[lemma]{Definition}
\newtheorem{remark}[lemma]{Remark}

\usepackage{dsfont}
\newcommand{\1}{\mathds{1}}

\newcommand{\W}{\Omega}
\newcommand{\R}{{\mathbb R}}

\newcommand{\C}{{\mathbb C}}

\renewcommand{\P}{{\mathbb P}}
\renewcommand{\O}{{\Omega}}

\newcommand{\Id}{\operatorname{I}}
\newcommand{\mc}[1]{\mathcal{#1}}

\newcommand{\WFh}{\operatorname{WF}_h}

\newcommand{\intw}{\mathcal{I}_\omega}

\usepackage{multicol}
\usepackage{enumitem}

\newcommand{\abs}[1]{\ensuremath{\left| #1 \right|}}

\newcommand{\comp}{\operatorname{comp}}

\renewcommand{\a}{{\bf a}}

\def\XXint#1#2#3{{\setbox0=\hbox{$#1{#2#3}{\int}$} \vcenter{\hbox{$#2#3$}}\kern-.5\wd0}}

\DeclareMathOperator{\supp}{supp}

\newcommand{\ran}{\operatorname{ran}}
\newcommand{\e}{\epsilon}

\numberwithin{equation}{section}
\numberwithin{lemma}{section}

\usepackage{subcaption}

\newcommand{\spec}{\operatorname{spec}}

\newcommand{\rr}{\mathbf{r}}
\newcommand{\dd}{\mathbf{d}}

\newcommand{\hardy}{\mathfrak{h}}
\newcommand{\hardyr}{\mathfrak{h}^{r_0}}
\newcommand{\rardy}{\mathfrak{r}}
\newcommand{\rardyr}{\mathfrak{r}^{r_0}}
\newcommand{\rardya}{\mathfrak{r}_\alpha}

\newcommand{\loc}{\operatorname{loc}}
\newcommand{\Diff}{\operatorname{Diff}}
\newcommand{\Diffhom}{\prescript{\bullet}{}{\Diff}}

\newcommand{\dist}{\operatorname{dist}}

\newcommand{\Op}{\operatorname{Op}}

\newcommand{\dtn}{DtN }
\newcommand{\acs}{AC smooth }
\newcommand{\acl}{AC Lipschitz }

\renewcommand{\div}{\operatorname{div}}
\newcommand{\WF}{\operatorname{WF}}
\newcommand{\ad}{\operatorname{ad}}
\newcommand{\op}{\operatorname{Op}}
\newcommand{\oph}{\operatorname{Op}_h}

\newcommand{\edra}{\begin{color}{roColor}}
\newcommand{\edrb}{\end{color}}
\usepackage{marginnote}

\newcommand{\xbd}{{x_0}}
\newcommand{\opsi}{O_{\Psi}}
\newcommand{\rre}{\langle \bar x\rangle}
\newcommand{\dde}{\langle x_0\rangle}
\newcommand{\fbdy}{f_\partial}
\newcommand{\oxy}{O_{\div}^{X\rightarrow Y}}
\newcommand{\osing}{\Omega_{\operatorname{sing}}}
\newcommand{\oreg}{\Omega_{\operatorname{reg}}}
\newcommand{\om}{\sqrt{\omega^2+M^2}}
\newcommand{\ohm}{\sqrt{\omega^2+h^2M^2}}

\newcommand{\dom}{\operatorname{Dom}}
\newcommand{\cl}{\operatorname{cl}}

\newcommand{\gplus}{{g_{+}}}

\title[]{Scattering for the Dirichlet-to-Neumann map and the Steklov spectrum}

\author{Jeffrey Galkowski}
\address{Department of Mathematics, University College London, London, UK}
\email{j.galkowski@ucl.ac.uk}

\author{Marcello Malagutti}
\address{Department of Mathematics, University of Warwick, Coventry, UK}
\email{Marcello.Malagutti@warwick.ac.uk}

\author{Ruoyu P. T. Wang}
\address{Department of Mathematics, Yale University, New Haven, CT, USA}
\email{ruoyu.wang@yale.edu}

\begin{document}
\begin{abstract}
We develop scattering theory for Dirichlet-to-Neumann (DtN) maps on asymptotically conic, unbounded domains $\Omega\subset \mathbb{R}^d$, that are only Lipschitz inside a compact set and whose boundary extends to infinity. Specifically, we prove a limiting absorption principle and demonstrate existence of the scattering matrix. 

As an application of the theory, we show that the eigenvalues of Steklov-type problems on curvilinear polygons, $\mc{P}$, are well approximated by the eigenvalues of a Dirac operator on $\partial\mc{P}$ whose transmission coefficients between edges meeting at a corner are given by the scattering matrix for DtN maps on asymptotically conic domains associated to the corners. Up to an arbitrary loss in the power of the eigenvalue, the error in the approximation is optimal among \emph{all} possible asymptotics that use only local information about $\partial\mc{P}$.
\end{abstract}

\maketitle
\setcounter{tocdepth}{1}
\tableofcontents


\section{Introduction}

Let $\mc{P}\subset \mathbb{R}^d$ be a bounded domain with piecewise smooth, Lipschitz boundary. We study the following Steklov-type spectral problem on $\mc{P}$. Fix $M\geq 0$, and consider the non-trivial solutions $u_\sigma\in H^1(\mc{P})$ of
\begin{equation}
\label{e:stekProblem}
\begin{cases}(-\Delta+M^2)u_{\sigma}=0&\text{ in }\mc{P}\\
\partial_\nu u_\sigma=\sigma u_\sigma &\text{ on }\partial\mc{P},
\end{cases}
\end{equation}
where $\nu$ is the outward pointing unit normal to $\partial \mc{P}$. For a Lipschitz domain, $\mc{P}$, there is a sequence $0\leq \sigma_0(\mc{P},M)\leq \sigma_1(\mc{P},M)\leq \dots$, $\sigma_j(\mc{P},M)\rightarrow \infty$ as $j\to \infty$ such that~\eqref{e:stekProblem} has a nontrivial solution if and only if $\sigma=\sigma_j(\mc{P},M)$ for some $j$: see for example ~\cite{ArMa:12}. In fact, $\{\sigma_j(\mc{P},M)\}_{j=0}^\infty$ are the eigenvalues of the associated Dirichlet-to-Neumann (DtN) map, $\Lambda_M^{\mc{P}}:H^1(\partial\mc{P})\to L^2(\partial\mc{P})$ with the corresponding eigenfunctions $u_{\sigma_j(M)}|_{\partial\mc{P}}$.

While we develop tools which will be useful for the study of~\eqref{e:stekProblem} in any dimension, our main application to this problem is in dimension $2$. 
\begin{definition}[Curvilinear polygons]
    We say $\mc{P}\subset \mathbb{R}^2$ is a \emph{curvilinear polygon} if it is a bounded, piecewise smooth, Lipschitz domain. We then decompose the boundary as $\partial\mathcal{P}=\sqcup_{l=1}^L\Gamma^l$ and $\Gamma^l=\cup_{j=0}^{J_l-1} \Gamma_{j}^l$. We then identify $\Gamma_j^{l}$ with $[x^l_j,x^{l}_{j+1}]$ using the arc length parametrized map $\gamma^l_j:[x^{l}_j,x^{l}_{j+1}]\to \Gamma_j^l$ and denote by $\alpha_j^l$ the angle between $-(\gamma^l_j)'(x_{j}^l)$ and $(\gamma^l_{j-1})'(x_{j}^l)$. Let
\begin{equation}
    \alpha_+=\max\{\alpha_j^l: 1\le l\le L,1\le j\le J_l\}
\end{equation}
be the maximal corner angle of $\partial\mc{P}$. 
\end{definition}

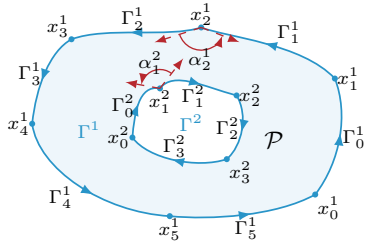
\begin{figure}[ht]
\centering
\begin{tikzpicture}[
    x=1cm,
    y=1cm,
    >=latex,
    font=\small,
    every node/.style={
        font=\fontsize{7}{8}\selectfont,
        inner sep=0.7pt
    },
    mark size=1.5pt,
    mark options={thin},
    boundary/.style={
        draw=plotTertiary,
        line width=0.65pt
    },
    oriented/.style={
        postaction={decorate},
        decoration={
            markings,
            mark=at position 0.55 with {\arrow{latex}}
        }
    },
    vertex/.style={
        circle,
        fill=plotTertiary,
        inner sep=0.8pt
    },
    tangent/.style={
        draw=plotPrimary,
        dashed,
        line width=0.5pt,
        ->
    },
    corner/.style={
        draw=plotPrimary,
        line width=0.5pt,
        ->
    }
]

\definecolor{colorAKANE}{RGB}{183,40,46}
\definecolor{colorUSUAI}{RGB}{40,148,200}
\colorlet{plotPrimary}{colorAKANE}
\colorlet{plotTertiary}{colorUSUAI}

\path[use as bounding box] (0,0) rectangle (7.5,3.2);
\begin{scope}[shift={(3.79,1.58)},scale=0.52]

\coordinate (P0) at ( 3.2,-1.80);
\coordinate (P1) at ( 3.7, 1.15);
\coordinate (P2) at ( 0.3, 2.45);
\coordinate (P3) at (-3.0, 2.15);
\coordinate (P4) at (-4.0, 0.0);
\coordinate (P5) at (-0.5,-2.35);

\coordinate (Q0) at (-1.45,-0.35);
\coordinate (Q1) at (-0.75, 0.90);
\coordinate (Q2) at ( 1.20, 0.72);
\coordinate (Q3) at ( 0.95,-0.90);

\begin{scope}
    \fill[fill=plotTertiary!8,even odd rule]
        (P0)
          .. controls ( 3.85,-0.95) and ( 4.05, 0.25) .. (P1)
          .. controls ( 2.85, 2.05) and ( 0.95, 2.20) .. (P2)
          .. controls (-0.55, 2.18) and (-2.10, 2.48) .. (P3)
          .. controls (-3.55, 1.85) and (-3.90, 0.95) .. (P4)
          .. controls (-3.85,-1.15) and (-2.20,-2.25) .. (P5)
          .. controls ( 0.55,-2.45) and ( 2.20,-2.35) .. cycle
        (Q0)
          .. controls (-1.40, 0.25) and (-1.30, 1.00) .. (Q1)
          .. controls (-0.40, 1.35) and ( 0.35, 0.95) .. (Q2)
          .. controls ( 1.50, 0.22) and ( 1.35,-0.55) .. (Q3)
          .. controls ( 0.20,-1.15) and (-1.05,-0.90) .. cycle;
\end{scope}

\begin{scope}
    \draw[boundary,oriented]
        (P0) .. controls (3.85,-0.95) and (4.05,0.25) ..
        node[pos=.50,right] {$\Gamma^1_0$} (P1);

    \draw[boundary,oriented]
        (P1) .. controls (2.85,2.05) and (0.95,2.20) ..
        node[pos=.44,above right] {$\Gamma^1_1$} (P2);

    \draw[boundary,oriented]
        (P2) .. controls (-0.55,2.18) and (-2.10,2.48) ..
        node[pos=.53,above] {$\Gamma^1_2$} (P3);

    \draw[boundary,oriented]
        (P3) .. controls (-3.55,1.85) and (-3.90,0.95) ..
        node[pos=.52,left] {$\Gamma^1_3$} (P4);

    \draw[boundary,oriented]
        (P4) .. controls (-3.85,-1.15) and (-2.20,-2.25) ..
        node[pos=.48,below left] {$\Gamma^1_4$} (P5);

    \draw[boundary,oriented]
        (P5) .. controls (0.55,-2.45) and (2.20,-2.35) ..
        node[pos=.52,below] {$\Gamma^1_5$} (P0);
\end{scope}

\begin{scope}
    \draw[boundary,oriented]
        (Q0) .. controls (-1.40,0.25) and (-1.30,1.00) ..
        node[pos=.45,left] {$\Gamma^2_0$} (Q1);

    \draw[boundary,oriented]
        (Q1) .. controls (-0.40,1.35) and (0.35,0.95) ..
        node[pos=.53,below] {$\Gamma^2_1$} (Q2);

    \draw[boundary,oriented]
        (Q2) .. controls (1.50,0.22) and (1.35,-0.55) ..
        node[pos=.48,left] {$\Gamma^2_2$} (Q3);

    \draw[boundary,oriented]
        (Q3) .. controls (0.20,-1.15) and (-1.05,-0.90) ..
        node[pos=.50,above] {$\Gamma^2_3$} (Q0);
\end{scope}

\begin{scope}
    \foreach \P in {P0,P1,P2,P3,P4,P5,Q0,Q1,Q2,Q3}
        \node[vertex] at (\P) {};

    \node[below right] at (P0) {$x^1_0$};
    \node[right]       at (P1) {$x^1_1$};
    \node[above]       at (P2) {$x^1_2$};
    \node[above left]  at (P3) {$x^1_3$};
    \node[left]        at (P4) {$x^1_4$};
    \node[below]       at (P5) {$x^1_5$};

    \node[left]        at (Q0) {$x^2_0$};
    \node[below]       at (Q1) {$x^2_1$};
    \node[right]       at (Q2) {$x^2_2$};
    \node[below right] at (Q3) {$x^2_3$};
\end{scope}

\begin{scope}
    \coordinate (Tprevious) at ( 1.275,2.075);
    \coordinate (Tnext)     at (-0.890,2.072);

    \draw[tangent] (P2) -- (Tprevious);
    \draw[tangent] (P2) -- (Tnext);

    \pic[
        corner,
        "$\alpha^1_2$",
        angle radius=3.0mm,
        angle eccentricity=1.40
    ] {angle=Tnext--P2--Tprevious};

    \coordinate (Uprevious) at (-1.630,1.060);
    \coordinate (Unext)     at (-0.155,1.665);

    \draw[tangent] (Q1) -- (Uprevious);
    \draw[tangent] (Q1) -- (Unext);

    \pic[
        corner,
        "$\alpha^2_1$",
        angle radius=2.5mm,
        angle eccentricity=1.45
    ] {angle=Unext--Q1--Uprevious};
\end{scope}

\begin{scope}
    \node[text=plotTertiary]
        at (-2.55,-0.15) {$\Gamma^1$};

    \node[text=plotTertiary]
        at (0.05,0.02) {$\Gamma^2$};

    \node[font=\fontsize{9}{10}\selectfont]
        at (2.15,-0.40) {$\mathcal P$};
\end{scope}

\end{scope}
\end{tikzpicture}
\caption{A curvilinear polygon $\mc{P}$ with $\partial\mathcal P=\Gamma^1\sqcup\Gamma^2$ as parametrized arcs.}
\end{figure}

On the boundary of a curvilinear polygon, we define the Steklov Dirac operator.
\begin{definition}[The Steklov Dirac operator]
Let $\mc{P}$ be a curvilinear polygon and $\{\Sigma^l_j(\lambda):1\le l\le L,1\le j\le J_l\}$ be $2\times 2$ unitary matrices such that for all $\lambda\in\mathbb{R}$,
\begin{equation}
    \sigma_3 \overline{\Sigma_j^l(\lambda)}\sigma_3=-\Sigma_j^l(-\lambda),\qquad \sigma_3=\begin{pmatrix}1&0\\0&-1\end{pmatrix}.
\end{equation}
Define the following $\lambda$-dependent family of Dirac operators on $\partial\mc{P}$: $\mathscr{D}(\mathcal{P},\Sigma_{j}^l):L^2(\partial\mc{P};\mathbb{C}^2)\to L^2(\partial\mc{P};\mathbb{C}^2)$ is the unbounded operator defined by
\begin{gather}
\mathscr{D}_\lambda(\mathcal{P},\Sigma_{j}^l)\begin{pmatrix}u_+\\u_-\end{pmatrix}= \begin{pmatrix}D_x & 0\\ 0 & -D_x\end{pmatrix}\begin{pmatrix}u_+\\u_-\end{pmatrix},\\
 \dom(\mathscr{D}_\lambda(\mathcal{P},\Sigma_{j}^l))=\Big\{ u\in \bigoplus_{j,l} H^1((x^l_j,x^l_{j+1});\mathbb{C}^2): \Sigma_{j}^l(\lambda)
 \begin{pmatrix}u_+(x^l_j-0)\\ u_-(x^l_j+0)\end{pmatrix}=\begin{pmatrix}u_-(x^l_j-0)\\u_+(x^l_j+0)\end{pmatrix}\Big\}.
\end{gather}
Denote
\begin{equation}
    \operatorname{char}(\mathscr{D}(\mathcal{P},\Sigma_j^l))=\{\lambda\in \mathbb{R}\,:\, \ker \big(\mathscr{D}_\lambda(\mathcal{P},\Sigma_j^l)-\lambda\big)\neq \{0\}\}.
\end{equation}
We say $\mathscr{D}(\mathcal{P},\Sigma_j^l)$ is a \emph{Steklov Dirac operator with matching $\Sigma^l_j$}, if we can arrange
\begin{equation}
    \operatorname{char}(\mathscr{D}(\mathcal{P},\Sigma_j^l))=\{\cdots\le \lambda_{-2}(\mathscr{D}(\mathcal{P},\Sigma_j^l))\le\lambda_{-1}(\mathscr{D}(\mathcal{P},\Sigma_j^l))\leq \lambda_0(\mathscr{D}(\mathcal{P},\Sigma_j^l))\leq \lambda_1(\mathscr{D}(\mathcal{P},\Sigma_j^l))\leq \cdots\}
\end{equation}
as a countable bi-sequence where $\lambda_n$ is repeated with multiplicity so that $\lambda_n=-\lambda_{-n-1}$ for $n\ge 0$. 

If 
\begin{equation}
    \Sigma_{j}^l(\lambda)=\Sigma_{1,0}^{\mathcal{W}_{\alpha^l_j}}, \ \Sigma_{1,0}^{\mathcal{W}_\alpha}=\begin{pmatrix}-i\cos \frac{\pi^2}{2\alpha}&\sin\frac{\pi^2}{2\alpha}\\\sin\frac{\pi^2}{2\alpha}&-i\cos\frac{\pi^2}{2\alpha}\end{pmatrix},
\end{equation}
where $\Sigma_{1,0}^{\mathcal{W}_\alpha}$ is the scattering matrix for the infinite wedge of angle $\alpha$ with $M=0$ (see Theorem~\ref{t:simpleScatteringMatrix}), then $\mathscr{D}_\lambda(\mathcal{P},\Sigma_{j}^l)$ is also $\lambda$-independent and we denote $\mathscr{D}^0=\mathscr{D}_\lambda(\mathcal{P},\Sigma_{j}^l)$ and call it the \emph{basic Steklov Dirac operator on $\mathcal{P}$}. In this case, we have $\spec(\mathscr{D}^0)=\operatorname{char}(\mathscr{D}(\mathcal{P},\Sigma_j^l))$. 

\end{definition}

\begin{theorem}
\label{t:boundedSteklov}
Let $\mc{P}\Subset \mathbb{R}^2$ be a curvilinear polygon, $M_0\geq 0$, and $k_0\in\mathbb{N}$. Then for all $\e>0$, there is $C>0$ such that for $n\geq 0$, $M\in[0,M_0]$ we have
\begin{equation}
\label{e:basicCompare}
|\sigma_n(\mc{P})-\lambda_{n}(\mathscr{D}^0)|< C\langle\sigma_n(\mc{P})\rangle^{ -1+\e},
\end{equation}
and there are $\{\Sigma_{j}^{l,k_0}(\lambda)\}$, with $\Sigma_{j}^{l,k_0}(\lambda)$ depending only on the first $(k_0+1)$ derivatives of $\gamma_j^l $ and $\gamma_{j-1}^l$ at $x_j^l$ such that for $n\ge 0$, $M\in[0,M_0]$ we have
\begin{equation}
\label{e:highOrder}
\Big|\sigma_n(\mc{P})-\sqrt{\lambda_n^2(\mathscr{D}(\mathcal{P},\Sigma_{j}^{l,k_0}))+M^2}\Big|\leq C\big(\langle\sigma_n(\mathcal{P})\rangle^{-\frac{2\pi}{\alpha_+}+\e}+\langle\sigma_n(\mathcal{P})\rangle^{-1-k_0+\e}+M^{2-\e}\langle\sigma_{n}(\mathcal{P})\rangle^{-2+\e}\big).
\end{equation}
Furthermore we have: 
\begin{enumerate}[wide]
    \item[(i)] If the derivatives of order $2\leq k\leq k_0+1$ of $\gamma_j^l$ and $\gamma_{j-1}^l$ at $x_j^l$ vanish, then $\Sigma_{j}^{l,k_0}=\Sigma_{1,0}^{\mc{W}_{\alpha_j^l}}$.
    \item[(ii)] The estimate \begin{equation}
     \ \|\Sigma_{j}^{l,k_0}(\lambda)-\Sigma_{j}^{l,k_0'}(\lambda)\|\leq C \langle \lambda\rangle^{-1-k_0+\e},\qquad \big\|\partial_{\lambda}\big(\Sigma_{j}^{l,k_0}-\Sigma_{j}^{l,k_0'}\big)(\lambda)\big\|\leq C \langle \lambda\rangle^{-(1+k_0)/2+\e}
\end{equation}
holds for any $0\leq k_0<k_0'\leq \lceil\frac{2\pi}{\alpha_+}\rceil-1$. 
\end{enumerate}
\end{theorem}
Theorem~\ref{t:boundedSteklov} implies that, when $M=0$, there is a Dirac operator on $\partial\mc{P}$ with matching conditions determined by the curvatures at corners so that the error in~\eqref{e:highOrder} is $O(\sigma^{-\frac{2\pi}{\alpha_+}+0})$. Moreover, if the curvatures vanish to high enough order at each corner, then the error in~\eqref{e:basicCompare} can be improved to  $O(\sigma^{-\frac{2\pi}{\alpha_+}+0})$. 

\begin{figure}
    \centering

    \includegraphics[width=\linewidth]{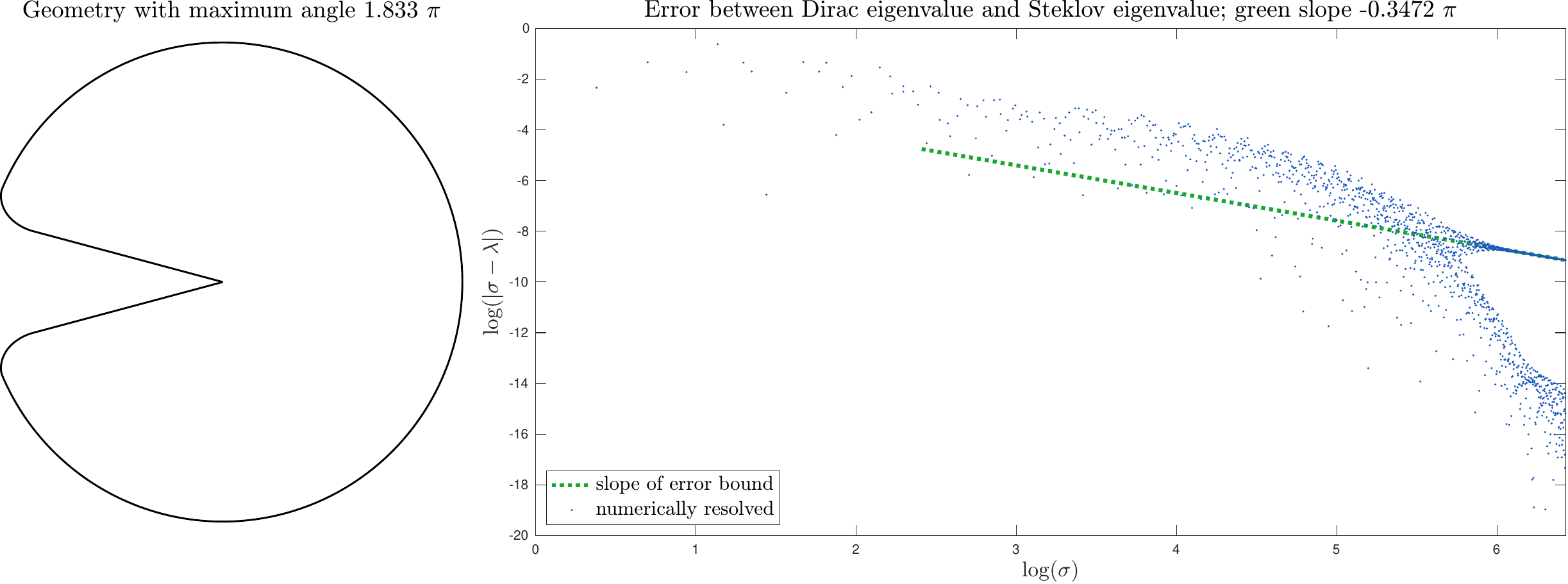}

     \includegraphics[width=\linewidth]{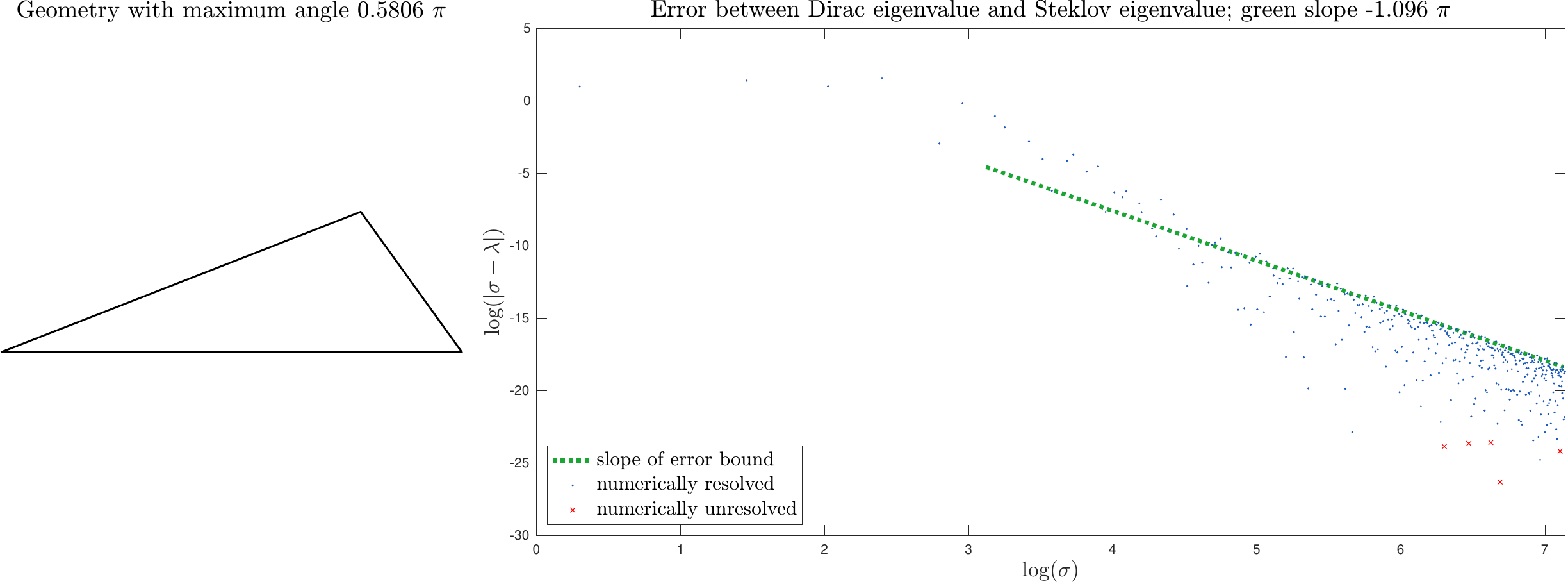}
     
 \includegraphics[width=\linewidth]{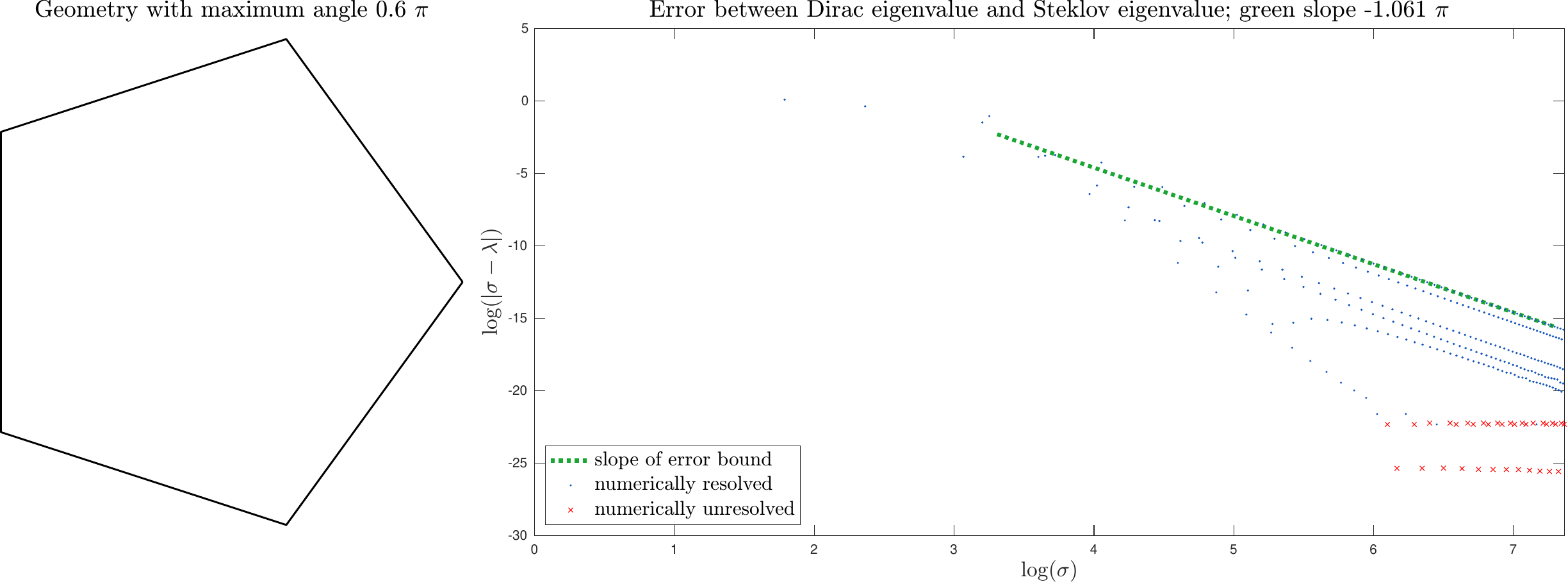}

    \caption{The left column depicts a curvilinear polygon and the right column shows the corresponding plot of $\log |\sigma_n(\mc{P})-\lambda_n(\mathscr{D}^0)|$ v. $\log |\sigma_n(\mc{P})|$ (with $M=0$) in the blue dots. Each curvilinear polygon has one corner at which all curvatures vanish. Any red x's are computed errors where the numerical resolution error was likely higher than the gap between the two spectra. In each figure, the green, dotted line is the upper bound on this error in Theorem~\ref{t:boundedSteklov}.}
    \label{f:numerics1}
\end{figure}

In a recent paper, the first and last authors show that when $M=0$~\eqref{e:highOrder} is optimal (modulo the $\e$ loss in the power) among \emph{all} possible asymptotic expressions for $\sigma_n$ that use only local information about the boundary of $\mc{P}$~\cite[Theorem 1.2]{GaWa:26}. In fact, the authors construct a curvilinear polygon $\mc{P}_{\alpha_+}$ with a single corner of angle $\alpha_+$ such that there are $r>0$, $x_0,y_0\in\mc{P}_{\alpha_{+}}$ satisfying 
\begin{equation}    
\limsup_{n\to \infty}\Big|\sigma_n\big(\mc{P}_{\alpha_+}\setminus B(x_0,r)\big)-\sigma_n\big(\mc{P}_{\alpha_+}\setminus B(y_0,r)\big)\Big|\sigma_n^{\frac{2\pi}{\alpha_+}}>0,
\end{equation}
where $B(x,r)$ denotes the ball of radius $r$ centered at $x$.
Figure~\ref{f:numerics1} shows numerically computed differences between the DtN spectrum and Steklov Dirac spectrum in a variety of curvilinear polygons.

\paragraph{\bf Previous work on the Steklov spectrum for curvilinear polygons} \ 

When $M=0$, and the angles of the curvilinear polygon are less than $\pi$,~\eqref{e:stekProblem} was extensively studied in~\cite{LePaPoSh:22a,LePaPoSh:22b}. Under these additional restrictions the authors proved a substantially weaker version of Theorem~\ref{t:boundedSteklov} where the error was given by $\sigma_n^{-\frac{1}{2}(\min(\frac{\pi}{\alpha_+},\frac{3}{2})-1)+0}$. In order to do this, the authors of~\cite{LePaPoSh:22a,LePaPoSh:22b} used a classical result on the sloping beach problem from~\cite{Pe:50} to construct a large family of quasimodes for~\eqref{e:stekProblem} and, through a series of ingenious comparisons, were able to show that these quasimodes, in fact, form a basis for $L^2(\partial\mc{P})$ modulo a finite dimensional set. Thus, they concluded that all large enough Steklov eigenvalues on $\mc{P}$ lie close to one quasieigenvalue. 

The method of proof used in~\cite{LePaPoSh:22a,LePaPoSh:22b} is highly tuned to exactly the Steklov problem and small perturbations of that problem destroy much of the proof. For instance, their methods do not provide any information when one takes $M>0$. Moreover, there are a number of places in the proof which appear to be `miracles'; for instance some poorly understood cancellations in certain quasimodes. Even for curvilinear polygons with $M=0$, this sensitivity means that the methods of~\cite{LePaPoSh:22a,LePaPoSh:22b}  could not treat angles $\alpha> \pi$.  

To move toward a study of the Steklov problem on piecewise smooth Lipschitz domains in higher dimensions, more robust methods are needed. For instance, if one considers a domain with an $\mathbb{S}^1$ symmetry in $\mathbb{R}^3$, separation of variables with angular frequency, $n$, immediately leads to the Steklov problem with $M=|n|$. The methods developed in this paper are applicable much more broadly than those in~\cite{LePaPoSh:22a,LePaPoSh:22b}. Specifically for $M>0$, they allow us to treat any planar, curvilinear polygon, and, although we do not pursue them here, even some other types of problems such as the interface problem of~\cite[(1.1)]{Bh:Pe:13}, \cite[(5.3)]{AuKl:11} (see also~\cite[(11)]{Bo:92}) and Steklov problems on Riemannian manifolds.  

\begin{figure}
\begin{tikzpicture}
\def \r{.8};
\def \rb{1.7};

\coordinate (A) at (0, 0);
\coordinate (B) at (3, 1);
\coordinate (C) at (5, 3);

\coordinate (D) at (1, 6);

\coordinate (E) at (8,0);
\coordinate (F) at (9,4);
\coordinate (G) at (7,7);
\coordinate (H) at (-4,6);
\coordinate(I) at (1.5,3.5);

%
%

\node[above] at ($(A)+(0,.3)$) {$\alpha_0$};
\node[above left] at (B) {$\alpha_1$};
\node[above left] at (C) {$\alpha_2$};
\node[below] at (D) {$\alpha_3$};
\node at (I) {$\mc{P}$};

\draw [dashed,thick] (A) circle (\r);
\draw [->,thick]  ($(E)+(-.2,0)$) to[out=190, in=-10]($(A)+({.1*\r},{-.1*\r})$) ;

\begin{scope}
\clip (E) circle (\rb-.1);
\draw [thick,fill=light-gray] ($(E)+(120:\rb)$)--(E)--($(E)+(60:\rb)$) arc (60:120:\rb)-- cycle;
\draw node at ($(E)+(90:{\rb/2})$){$\mc{W}_{\alpha_0}$};

\end{scope}

\draw [dashed,thick] (B) circle (\r);
\begin{scope}
\clip (F) circle (\rb-.1);
\draw [thick,fill=light-gray] ($(F)+(90:\rb)$)--(F)--($(F)+(225:\rb)$)arc (225:90:\rb)-- cycle;
\draw node at ($(F)+(165:{\rb/2})$){$\mc{W}_{\alpha_1}$};
\end{scope}
\draw [->,thick]  ($(F)+(.2,0)$) to[out=0, in=-10]($(B)+({.1*\r},{-.1*\r})$) ;

\draw [dashed,thick] (C) circle (\r);
\begin{scope}
\clip (G) circle (\rb-.1);
\draw [thick, fill=light-gray] ($(G)+(220:\rb)$)--(G)--($(G)+(60:\rb)$) arc (60:220:\rb)--cycle;
\draw node at ($(G)+(140:{\rb/2})$){$\mc{W}_{\alpha_2}$};
\end{scope}
\draw [->,thick]  ($(G)+(0,-.2)$) to[out=250, in=0]($(C)+({.1*\r},{-.1*\r})$) ;

\draw [dashed,thick] (D) circle (\r);
\begin{scope}
\clip (H) circle (\rb-.1);
\draw [thick,fill=light-gray] ($(H)+(0:\rb)$)--(H)--($(H)+(110:\rb)$) arc (110:360:\rb)--cycle;
\draw node at ($(H)+(235:{\rb/2})$){$\mc{W}_{\alpha_3}$};
\end{scope}

\draw [->,thick]  ($(H)+(.1,.2)$) to[out=90, in=90]($(D)+(50:{.1*\r})$) ;


\draw[thick, smooth] (A) to[out=60, in=225] (B) to[out=90, in=220] (C) to [out=60, in=0](D) to [out=110, in =120](A);

\end{tikzpicture}

\caption{\label{f:1}A curvilinear polygon $\mc{P}$ having angles $\alpha_0,\dots,\alpha_3$ with the infinite wedge of angle $\alpha_j$, $\mc{W}_{\alpha_j}$, inserted at the $j^{\text{th}}$ corner.}
\end{figure}
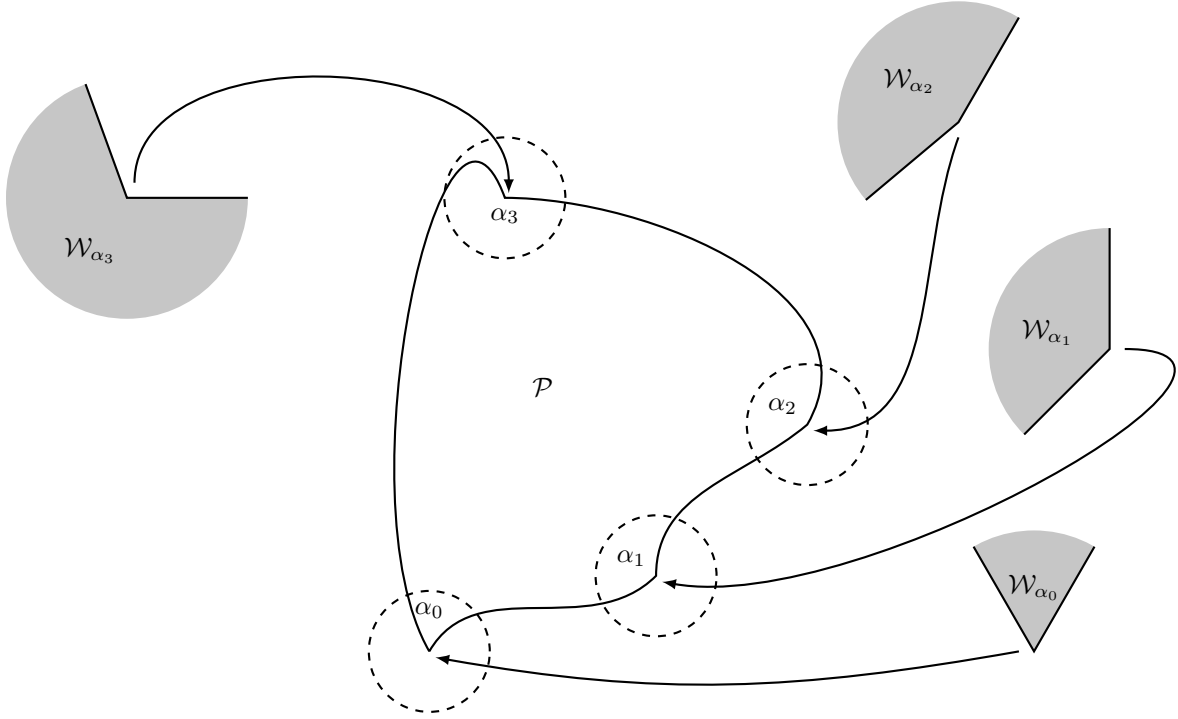

Our methods are based on the scattering theory for the DtN map in non-compact domains, $\Omega$, with boundary extending to infinity. We view~\eqref{e:stekProblem} for large $\sigma$ as a problem on smooth, open curves connected by the scattering problem for the DtN map on an infinite wedge at each corner of the domain (see Figure~\ref{f:1}). While the DtN map is \emph{far} from a local operator, it is nevertheless possible to separate the effects of the corners and the smooth sides at high energy. 

\noindent{ \bf{Scattering theory for the Dirichlet-to-Neumann Map}} 

 To describe the geometry in which we work, we need two definitions.
\begin{definition}
\label{d:scatteringManifold}
Let $(M,g)$ be a complete Riemannian manifold of dimension $d-1$, such that $M=\cup_{l=1}^L M_l\cup K$ where $K$ is compact, $M_l$ are diffeomorphic to $(0,\infty)_r\times Y_l$ for some smooth, compact $(d-2)$-dimensional manifolds, $Y_l$. Then, $M$ can be compactified to $\overline{M}$ using the coordinate $\rho=1/r$ in $M_l$. We allow $K$ to be Lipschitz, that is, locally defined using charts that are only Lipschitz. We say that $M$ is \emph{a scattering manifold} if there is a choice $\rho$ of defining function for $\partial\overline{M}_l$ and a collar neighborhood $[0,\e)_\rho\times Y_l$ of $\partial\overline{M}_l$ such that
\begin{equation}
    g=\frac{d\rho^2}{\rho^4}+\frac{h(\rho)}{\rho^2},
\end{equation}
where $h:[0,\e)\to \operatorname{Sym}^2(T^* Y_l)$ is a smooth family of Riemannian metrics on $Y_l$. 
\end{definition}

\begin{definition}\label{d:ACLip}
We say that a domain $\Omega$ in $\R^d$ with boundary $\partial\Omega$ is an \emph{asymptotically conic Lipschitz domain with cross section $Y$} if $\Omega$ is a Lipschitz domain, $Y\subset \mathbb{S}^{d-1}$ is open with smooth boundary $\partial Y$, and there are $R_0>0$ and $\gamma\in C^\infty((0,1]\times \overline{Y}; \mathbb{S}^{d-1})$, 
with $\gamma(\rho,\cdot)$ as an injective map, satisfying the following conditions:
\begin{enumerate}
\item $\partial\Omega \setminus B(0,R_0) =\{ \rho^{-1}\gamma(\rho,\omega)\,:\, \rho\in(0,1],\, \omega\in \partial Y\}\cap \mathbb{R}^d\setminus B(0,R_0),$
    \item When $d=2$, we ask $\gamma\in S^0$, that is for $\rho \in (0, 1]$, $\omega\in\overline{Y}$, 
    \begin{equation}
        |(\rho\partial_\rho)^\alpha \partial_\omega^\beta \gamma (\rho,\omega)|\leq C_{\alpha\beta}, \ \lim_{\rho\to 0} \|D(\rho^{-1}\gamma)(\rho,\omega) - D(\rho^{-1}\iota)(\rho,\omega)\|_{L^\infty(Y)}=0.
    \end{equation}
    \item When $d>2$,
    \begin{equation}
        \gamma\in C^\infty([0,1]\times \overline{Y}; \mathbb{S}^{d-1}), \ \gamma(0,\omega)=\omega\text{ for all }\omega\in\overline{Y}. 
    \end{equation}
\end{enumerate}
\end{definition}

\begin{remark}
Notice that when $d=2$, the requirement that $(\partial\Omega,g)$ be asymptotically conic is always satisfied since $\partial\Omega$ is 1-dimensional.
\end{remark}

We now give two consequences of our second main theorem whose statement we postpone to Section~\ref{s:microlocalDtN}. The first proves the existence of the limiting absorption resolvent and hence absolute continuity of the spectrum away from eigenvalues.
\begin{theorem}
\label{t:simpleResolve}
Let $\Omega\subset \mathbb{R}^d$ be an asymptotically conic domain and $\delta>0$. Then there is $N\in\{0,1,\dots\}\cup \{\infty\}$ such that the set of eigenvalues, $\{\omega_n\}_{n=1}^N$ of $\Lambda_{M}^\Omega$ in $[M+\delta,\infty)$ is discrete. Moreover, for $\omega\notin\{\omega_n\}_{n=1}^N$, the limiting absorption resolvent, $(\Lambda_{M}^{\Omega}-\omega\pm i0)^{-1}:L^2_{\comp}\to L^2_{\loc}$ exists. Furthermore, if $d=2$ and $\Omega$ is simply connected and smooth away from a single point, or $d>2$ and $\Omega$ is scaling invariant, then $N<\infty$. 
\end{theorem}

Our next theorem demonstrates the existence of the scattering matrix.
\begin{theorem}
\label{t:simpleScatteringMatrix}
Let $\Omega\subset \mathbb{R}^d$ be an asymptotically conic Lipschitz domain with cross-section $Y$ and $\omega>0$. Then for all $g^-\in L^2(\partial Y)$, there is a solution $u$ to 
\begin{equation}
    (\Lambda_M^\Omega-\sqrt{\omega^2+M^2})u=0
\end{equation}
such that, in the coordinates of Definition~\ref{d:scatteringManifold}, 
\begin{equation}
    u- \rho^{\frac{d-2}{2}}(e^{-i\omega/\rho}g^-(y) +e^{i\omega/\rho}g^+(y)) = O_{L^2(\partial \mc{W}_Y)}(\rho^{\frac{d}2}),
\end{equation}
for some $g^+\in L^2(\partial Y)$. 
Moreover, letting $\Pi_\omega$ denote the orthogonal projector onto the eigenspace of $\Lambda_M^\Omega$ at $\sqrt{\omega^2+M^2}$, the function $(I-\Pi_\omega)u$ is unique. In addition, the \emph{scattering matrix}, $\Sigma^{\Omega}_{\omega,M}:L^2(\partial Y)\to L^2(\partial Y)$ given by $\Sigma^{\Omega}_{\omega,M}: g^-\mapsto g^+$, is unitary.
\end{theorem}

\begin{remark}
    Theorems~\ref{t:simpleResolve} and~\ref{t:simpleScatteringMatrix} are simplified versions of the main theorems from Part 1 on scattering theory, Proposition~\ref{p:LAP} and Proposition~\ref{p:existscmatrix} respectively.
\end{remark}

\subsection{Sketch of the strategy and outline of the paper} \ \\

We now give a broad overview of the strategy of proof. We emphasize that, in addition, each of the sections of the main body starts with a sketch of what is done in the corresponding section.

\noindent{{\bf Scattering for the Dirichlet-to-Neumann map on \acl domains in any dimension}}

The first part of the article, Sections~\ref{s:harm} to~\ref{s:scatMat}, is dedicated to the solution theory for the Dirichlet problem and the scattering theory of the DtN map on \acl domains, $\Omega$. When $M>0$, it is clear that the equation 
\begin{equation} \label{e:basicHarmonic}
(-\Delta+M^2)u=f\text{ in }\Omega,\qquad u|_{\partial\Omega}=0
\end{equation}
has a unique solution in $\cup_N\langle x\rangle^NH^1$ for any  $f\in \cup_N\langle x\rangle^NH^{-1}$. However, when $M=0$, this is no longer the case and, indeed, even solutions with $f\in H^{-1}$ are problematic. 

Section~\ref{s:harm} studies solutions to~\eqref{e:basicHarmonic}, for $M\in [0,M_0]$ on uniformly Lipschitz domains (i.e. those for which the Lipschitz norm of the boundary is controlled uniformly to infinity). To do this, we introduce spaces 
\begin{equation}
\begin{gathered}
    \dot X_m^l(\Omega)=\Big\{u\in H^1_{\rm\loc}(\Omega)\,:\, \big\langle d(x,\partial\Omega)\big\rangle ^{m-l}\langle x\rangle^l  \nabla u\in L^2,\, \big\langle d(x,\partial\Omega)\big\rangle ^{m-l-1}\langle x\rangle^l   u\in L^2,\, u|_{\partial\Omega}=0\Big\},\\ Y_{-m}^{-l} =(\dot X_{m}^l)^*.
    \end{gathered}
\end{equation}
and show that for certain ranges of $m$ and $l$ depending on the geometry of $\Omega$ at infinity, 
$$
(-\Delta+M^2):\{ u\in \dot X_m^l\,:\, (-\Delta +M^2)u\in Y_m^l\}\to Y_{m}^l
$$ is invertible for all $M \geq 0$ and has uniform estimates for $M\in[0, M_0]$ for any fixed $M_0>0$. We also show that the solution operator with $M=0$, $\mathcal{R}_0^\Omega: Y_m^l\to\dot X_m^l$, is a weak limit of those $\mathcal{R}_M^\Omega:Y_m^l\to\dot X_m^l$ with $M>0$ and define the corresponding harmonic extension operators solving
\begin{equation}
\label{e:basicHarmonicExtend}
    (-\Delta +M^2)\mathcal{H}_M^\Omega \fbdy=0,\qquad \mathcal{H}_M^\Omega \fbdy|_{\partial\Omega}=\fbdy.
\end{equation}
We then prove boundedness of $\Lambda_M^\Omega: H^{\frac12,l}(\partial\Omega)\rightarrow H^{-\frac12,l}(\partial\Omega)$. The estimates on the solution operators $\mathcal{R}_M^\Omega$ (see Lemma~\ref{l:stableEstimates}) rely on Hardy inequalities controlling $\||x|^{-1}u\|_{L^2}$ and $\|d(x,\partial\Omega)^{-1}u\|_{L^2}$ by $\|\nabla u\|_{L^2}$ and the construction of a delicate Carleman weight. The range of allowed values of $m$ and $l$ is related to the optimal constants in these Hardy inequalities and we show that these ranges are optimal. For future reference, note that on \acl domains, these ranges are
\begin{equation}\label{e:0l1}
    |m|<\rardy_\Omega^\infty,\qquad |m-l|<\tfrac{1}{2},\qquad \rardy_{\Omega}^\infty=\lim_{r_0\to \infty}\inf\Big\{\frac{\|\nabla u\|_{L^2(\Omega)}}{\|\frac{u}{|z|+r_0}\|_{L^2(\Omega)}}: u\in \dot{H}^{1}(\Omega)\setminus\{0\}\Big\}.
\end{equation}
An important example to keep in mind is the half space $\mathbb{R}^d_+=\{ (x_0,x)\in\mathbb{R}^d: x_0>0\}$, where the DtN map $\Lambda_M^{\mathbb{R}^d_+}$ is a Fourier multiplier
\begin{equation}\label{e:0l2}
    \Lambda_M^{\mathbb{R}^d_+}=\sqrt{D^2+M^2}.
\end{equation}
It is then easy to check that
\begin{equation}
    \Lambda_0^{\mathbb{R}^d_+}: H^{s+1, l}(\mathbb{R}^{d-1})\to H^{s, l}(\mathbb{R}^{d-1}),\ s\in\mathbb{R}, \ |l|<
    \tfrac{d+1}{2}
    =\rardy_{\mathbb{R}^{d}_+}^\infty+\tfrac12,
\end{equation}
and cannot be extended to a wider range of $l$. In particular, this shows that the thresholds for $m$ and $l$ in \eqref{e:0l1} cannot be improved.

With the basic solution theory for the Dirichlet problem in hand, Section~\ref{s:3} is devoted to the construction of a second microlocal pseudodifferential approximation to the DtN map on an \acl domain. 

Indeed, \eqref{e:0l2} becomes singular at $|\xi|=0$, $|x|=\infty$ as $M\to 0$ and, as a result the DtN map cannot be pseudodifferential modulo a smoothing and rapidly decaying error. Because of this, we add an artificial mass term to $-\Delta$ and, for any $0\leq \delta <1$, construct a parametrix modulo a smoothing and decaying operator for an artificial DtN map, $P$, corresponding to 
\begin{equation}
    -\Delta+M^2+\langle x\rangle^{-2\delta}. 
\end{equation}
In particular, using a factorization argument, we show that $P$ is a pseudodifferential operator with symbol 
\begin{equation}
    \sigma(P)=\langle x\rangle^{-\delta} \big\langle\langle x\rangle^{\delta}\sqrt{|\xi|_{g}^2+M^2}\big\rangle,
\end{equation}
where $g$ is the metric induced on $\partial\Omega$. (For the construction of the second microlocal calculus in which this lives, see Appendix~\ref{a:1}.) The solution theory from Section~\ref{s:harm} then allows us to control the error between $P$ and the true DtN map. 
One consequence of this approximation is that the harmonic extension and DtN maps on \acl domains are, in fact, pseudodifferential away from the zero frequency and away from the Lipschitz singularities. Thus, they can be extended to data that may have any polynomial growth in $x$ so long as it is away from the zero frequency.

Section~\ref{s:improvedHighFreq} contains a more detailed analysis of solutions to ~\eqref{e:basicHarmonic} and~\eqref{e:basicHarmonicExtend} on \acl domains. More precisely, we study the behavior of $\mathcal{R}_M^\Omega$ acting on functions that have non-trivial oscillations in the directions tangent to the boundary. To do this, we commute vector fields $\mathcal{V}$ that are tangent to $\partial\Omega$ and whose coefficients may grow like $\langle x\rangle $ through $\mathcal{R}_M^{\Omega}$. This analysis is ultimately used in Proposition~\ref{p:dtN} to control the error $\chi (\Lambda_M^\Omega-P)$. Specifically, we show that this error is rapidly decaying on high frequencies, maps high frequencies to high frequencies, low frequencies to low frequencies, and decays like $\langle x\rangle^{-\rardy_{\Omega}^\infty-\frac{1}{2}}$ at low frequencies. 
\begin{remark}
    One can check $\langle x\rangle^{-\rardy_{\Omega}^\infty-\frac{1}{2}}$ is the best possible decay rate for the low frequencies using explicit computations on the infinite planar wedge of any angle $0<\alpha<2\pi$. 
\end{remark}

At this point, we have a very accurate description of $\Lambda_M^\Omega$ microlocally near physical infinity and away from frequency $0$. Except for the non-local part of $\Lambda_M^\Omega$ near 0 frequency, the Fredholm analysis of $\Lambda_M^\Omega-\omega$ is close to that of the Laplacian on an asymptotically flat manifold. In Section~\ref{s:microlocalDtN} we use our microlocal description together with radial point estimates to obtain Fredholm estimates and a limiting absorption resolvent in certain variable order Sobolev spaces (see Proposition~\ref{p:LAP}). Then, in Section~\ref{s:scatMat}, we construct the scattering matrix associated to the DtN map (see Proposition~\ref{p:existscmatrix}). 

\bigskip

\noindent{\bf{Analysis of the Steklov problem on piecewise smooth, planar domains}}

 Sections~\ref{s:scatteringatrixConvergence} to~\ref{s:spectralAsymptotics} apply the scattering theory developed in part 1 to the Steklov problem on piecewise smooth Lipschitz domains -- resulting in Theorem~\ref{t:boundedSteklov}. 

The starting point is Section~\ref{s:scatteringatrixConvergence}, which contains a more detailed study of the scattering matrix on simply connected \acl domains in $\mathbb{R}^2$ whose boundary is smooth away from one single point (i.e. infinite wedges with possibly curved sides). 

We start by showing that, at high energy, $\omega$, the scattering matrix on a curved wedge $\mc{W}^1_\alpha$ with angle $\alpha$ is close to that of any other wedge $\mathcal{W}^2_\alpha$ of the same angle and that, if enough curvatures of $\mc{W}^{1}_{\alpha}$ and $\mc{W}^2_{\alpha}$ agree at the singularity, this difference is controlled by $\omega^{-\frac{2\pi}{\alpha}+0}+M^2\omega^{-2+0}$. To do this, we transplant the distorted plane waves from $\mc{W}^{1}_{\alpha}$ onto $\mc{W}^2_\alpha$. More precisely, let $u_{1}$ satisfy
\begin{equation}
    (\Lambda_{M}^{\mc{W}^1_\alpha}-\om)u_{1}=0,
\end{equation}
and let 
\begin{equation}
u_{1}=u_{1}^{\textnormal{high}}+u_{1}^{\textnormal{low}}
\end{equation}
be the decomposition of $u_{1}$ into high and low frequencies. We then set
\begin{equation}
    u_{2}(x)= \chi(x) \big(F^*u_{1}^{\textnormal{low}}\big)(x) +\chi(\omega^{1-\delta}x)\big(F^*u_{1}^{\textnormal{high}}\big)(x)+
    (1-
    \chi(\omega^{1-\delta}x))u_{\textnormal{match}}^{\textnormal{high}}(x)
    ,
\end{equation}
where $F:\mathcal{W}^2_\alpha\to \mathcal{W}^1_\alpha$, and $u_{\textnormal{match}}^{\textnormal{high}}$ are the distorted plane waves on the smooth edges of $\mathcal{W}_\alpha^2$ with scattering data matching that of $F^*u_{1}^{\textnormal{high}}$.
That is, we truncate the low frequencies of $u_{1}$ at distance from the corner $1$ and match the high frequency parts at distance $\omega^{-1+\delta}$ for some $0<\delta\ll 1$. We then show that 
$
(\Lambda_{M}^{\mc{W}^2_\alpha}-\om)u_{2}
$ is small and apply the limiting absorption resolvent to correct $u_{2}$ to an exact distorted plane wave and hence show that the scattering matrix on $\mc{W}^2_\alpha$ is close to that on $\mc{W}^1_\alpha$.

Next, we compare the scattering matrix on the wedge ${\mc{W}_\alpha}$ with some positive mass $M>0$ to that with $M=0$. In this analysis, we use a change of domain formula
\begin{equation}
\chi\Lambda^{\mc{W}_\alpha}_M= \Lambda_M^{\mathbb{R}^2_+} \chi+\partial_\nu \mc{R}_M^{\mathbb{R}^2_+}[-\Delta, \chi]\mc{H}^{\mc{W}_\alpha}_M,
\end{equation}
proved in Section~\ref{s:improvedHighFreq} together with very precise knowledge of $\Lambda_M^{\Omega}$, when $\Omega$ is a smooth asymptotically conic domain, to study the perturbation $\Lambda_M^{\mc{W}_\alpha}-\Lambda_0^{\mc{W}_\alpha}$.

Section~\ref{s:aprioriSteklov} determines the a priori structure of Steklov eigenfunctions on curvilinear polygons, showing in particular that they are, up to small low frequency errors (of size $\sigma^{-\frac{\pi}{\alpha}-1+0}$), (distorted) plane waves. This analysis proceeds in two steps. First, using the fact that, away from the corners of $\mc{P}$, the DtN map is a pseudodifferential operator with symbol $|\xi|_g$, we are able to show that Steklov eigenfunctions are (distorted) plane waves modulo an error of size $\sigma^{-1}$ at low frequency. With this in hand, similar to how we constructed distorted plane waves on one infinite wedge from another, we are able to construct a distorted plane wave on a curved, infinite wedge that agrees with the domain near the corner from the Steklov eigenfunction near a corner. Our Fredholm estimates then allow us to obtain precise control on this distorted plane wave and hence on $u$. In particular, they show that parts of the eigenfunction that are `incoming' and `outgoing' from the corner are nearly related by the scattering matrix for the infinite wedge.

In Section~\ref{s:grushin}, we construct a Grushin problem for the DtN map on a curvilinear polygon so that the effective Hamiltonian (approximately) characterizes the eigenvalues of the Dirac operator from Theorem~\ref{t:boundedSteklov}. Since we know apriori that Steklov eigenfunctions behave like plane waves on smooth parts of the boundary, the key component in the analysis is using the relationship between the waves incoming to a corner and those outgoing from that corner. This is achieved by the same quasimode construction. 

Section~\ref{s:spectralAsymptotics} then shows that this error does not significantly change the implied spectrum and uses a deformation argument to show that the labeling of the Dirac operator's eigenvalues exactly matches that of the Steklov eigenvalues. 

\noindent{\bf{Appendices}}

Section~\ref{a:1} is devoted to the construction of the second microlocal calculus $\Psi^{s,l}_{k,\delta}$ used in this paper. In Section~\ref{a:2} we prove a few results on the scattering geometry and properties of the Hardy and radial Hardy constants. Section~\ref{a:3} is devoted to the construction of \dtn quasimodes in 2D \acs domains, and Section~\ref{a:4} is devoted to scattering matrices on 2D straight wedges.

\medskip 

\noindent{\textsc{Acknowledgements:}} The authors would like to thank Leonid Parnovski for many useful discussions and comments during the writing of this article, especially about the works~\cite{LePaPoSh:22a,LePaPoSh:22b}. Thanks also to Jean Lagac\'e for the reference~\cite{BuGiTr:20} for continuity of the Steklov eigenvalues. JG is grateful to the EPSRC for support under grants EP/V001760/1 and EP/V051636/1, to the Leverhulme Trust for support under Research Project Grant  RPG-2023-325, and to the ERC for support under the Synergy grant PSINumScat 101167139. MM and RPTW were partially supported by EPSRC grant EP/V001760/1 during the writing of this paper. \\

\noindent{\textsc{AI statement:}} This project began in 2023, before AI tools were widely available. In the later stages of the project ChatGPT 5 and 6 were used for proofreading, generating the numerical codes which produce Figure~\ref{f:numerics1}, and generating the TikZ code for most figures. In addition, in an earlier version of the proof the Carleman estimates in Lemma~\ref{l:stableEstimates} were based only on Hardy inequalities of the form $\|u/\operatorname{dist}(x,\partial\Omega)\|_{L^2}\leq C\|\nabla u\|_{L^2}$. ChatGPT suggested replacing $1/\operatorname{dist}(x,\partial\Omega)$ by $1/|x|$. The authors then realized that this could be further improved to the current, optimal version. In addition, ChatGPT pointed out that one could replace Rouché's theorem, which was originally used in the proof of Lemma~\ref{l:nonUniformEps0} when $k_0=0$, by homotopy invariance of the winding number \cite[p. 134 property (c)]{Hatcher2002}. The manuscript was entirely written by the authors.

\part{The scattering theory for the DtN map on \acl domains}

\section{The harmonic extension in uniformly Lipschitz domains}
\label{s:harm}

The goal of this section is to study solutions of 
\begin{equation}\label{e:toSolveSection2}
(-\Delta+M^2)u=f,\text{ in }\Omega \qquad u|_{\partial\Omega}=g,
\end{equation} 
uniformly as $M\to 0$, either for domains $\Omega$ that may grow or even be unbounded. 

\noindent{\bf{Outline of the strategy}}

We work on domains whose boundaries are uniformly Lipschitz in a sense made precise in Definition~\ref{d:uniformLip}. When $M>0$ or $\Omega$ is bounded, the solution theory for~\eqref{e:toSolveSection2} is simple. Indeed~\eqref{e:toSolveSection2} is uniquely solvable for $f$ and $g$ with any amount of polynomial growth/decay. However, as $\Omega$ grows and $M\to 0$, the implicit constants in this statement degenerate. 

To develop a uniform solution theory, we introduce spaces in which derivatives of $u$ decay faster than $u$ itself (see Definition~\ref{d:Xspaces}). Hardy inequalities of the form
\begin{equation}
\label{e:hardyType}
\Big\|\frac{u}{f}\Big\|_{L^2(\Omega)}\leq \hardy_f\|\nabla u\|_{L^2(\Omega)}
\end{equation}
for different choices of (typically growing) $f$  than allow us to prove apriori estimates in these spaces in Lemma~\ref{l:stableEstimates} which in turn result in a uniform solution theory for~\eqref{e:toSolveSection2} with $g=0$ (see Proposition~\ref{p:extendNoMass}). These estimates are proved using a carefully designed Carleman-type weight that relates the sharp constants in the estimates~\eqref{e:hardyType} to the solution theory of~\eqref{e:toSolveSection2}. Specifically, the maximal growth/decay of $f$ for which one can solve~\eqref{e:toSolveSection2}uniquely are determined by these constants. We are then able to understand the case of $g\neq 0$ using arbitrary extension operators (see Proposition~\ref{p:harmonic}).

\subsection{Geometric Preliminaries and definition of solution operator}

We begin with a few geometric definitions
\begin{definition}[Uniformly Lipschitz domains]
\label{d:uniformLip}
We say a Lipschitz domain $\Omega\subset \mathbb{R}^d$, is \emph{uniformly Lipschitz} if its boundary consists of a countable family of codimension-1 hypersurfaces, $\partial \Omega=\cup_j\Gamma_j$. Here we ask each $\Gamma_j$ to be, up to a rotation and translation, a Lipschitz hypograph with uniform control: that is, we assume there exist $c>0$, $R_j>c$, $f_{j}:\mathbb{R}^{d-1}\to \mathbb{R}$ such that 
\begin{enumerate}
    \item injectivity radius: $\Gamma_j=\{ (f_j(x),x)\,:\, x\in B(0,R_j)\}$, up to a rotation and translation;
    \item uniformly Lipschitz characters:$\|\nabla  f_{j}\|_{L^\infty}\leq C$ independent of $j$;
    \item bounded local degree: $\sup_j \#\{ k\,:\Gamma_j\cap \Gamma_k\neq \emptyset\} <\infty$;
    \item hypograph: $\{(\xbd,x)\in \mathbb{R}^{d}\,:\, x\in B(0,R_{j}),\,f_{j}(x)<\xbd<c\}\subset \Omega$.
\end{enumerate}
\end{definition}
Bounded Lipschitz domains with a uniform Lipschitz constant are uniformly Lipschitz. Note all \acl domains are uniformly Lipschitz. Throughout this section, $\Omega\subset \mathbb{R}^d$ denotes a uniformly Lipschitz domain. 

We use $z\in\mathbb{R}^d$ as a global coordinate for $\Omega$. We assume $\{z=0\}\in \partial\Omega$. We consider set $\mc{V}\subset\partial\Omega$ to be a finite set of vertices that $\{z=0\}\in \mc{V}$. Let the regularized vertex radius function $\rr(z)\in C^{\infty}(\Omega; [\frac18,\infty))$
that
\begin{equation}
    \|\rr(z)-(1+\dist(z,\mc{V})^2)^{\frac12}\|_{L^\infty(\Omega)}\le \frac14, \ \rr\ge \frac14\text{ at }\mc{V}, \ \|\nabla\rr\|_{L^\infty}\le 1, \ \ |\nabla^k\rr(z)|\le C_k\rr^{1-k}(z)
\end{equation}
for $k\in\mathbb{N}$. Let $\dd:C^{\infty}(\Omega; [\frac18,\infty))$ be a regularized distance function to $\partial\Omega$, that
\begin{equation}
    \|\dd(z)-(1+\dist(z,\partial\Omega)^2)^{\frac12}\|_{L^\infty(\Omega)}\le \frac14, \ \dd\ge \frac14\text{ on }\partial\Omega, \ \|\nabla\dd\|_{L^\infty}\le 1, \ \ |\nabla^k\dd(z)|\le C_k\dd^{1-k}(z)
\end{equation}
for $k\in\mathbb{N}$. Without loss of generality $\dist(z,\partial\Omega)\le \dd\le C\rr$ for fixed $C\ge1$. 
For $s,l\in\mathbb{R}$, we denote the Sobolev space by $H^s(\Omega)$, and denote the weighted Sobolev spaces $H^{s,l}(\Omega)=\rr^{-l}H^s(\Omega)$. The trace-free spaces are denoted $\dot H^s(\Omega)$ and $\dot H^{s,l}(\Omega)$. The corresponding spaces over the boundary are $H^s(\partial\Omega)$ and $H^{s,l}(\partial\Omega)=\rr^{-l}H^s(\partial\Omega)$. We define the Dirichlet Laplacian, $-\Delta:\dot{H}^{1}(\Omega )\to H^{-1}(\Omega)$, as the Riesz representation of the quadratic form $Q:\dot{H}^{1}(\Omega)\times \dot{H}^{1}(\Omega)\to \mathbb{C}$, 
\begin{equation}
    Q(u,v):=\langle \nabla u,\nabla v\rangle_{L^2(\Omega)}.
\end{equation}
We understand $\cdot|_{\partial\Omega}$ as the trace if the trace exists, or the non-tangential limit if the trace does not exist, in the latter case of which we abuse notation and call the non-tangential limit the trace. We define the harmonic extension and the Dirichlet resolvent: 
\begin{definition}[Harmonic extension]
    Let $\Omega$ be a uniformly Lipschitz domain. Let $\mc{X}(\partial\Omega)$ and $\mc{Y}(\Omega)$ be Banach spaces of functions on $\partial\Omega$ and $\Omega$. Let $M\ge 0$. We define $\mc{H}^\Omega_M: \mc{X}(\partial\Omega)\rightarrow \mc{Y}(\Omega)$ by sending $\fbdy$ to $\mc{H}_M^\Omega \fbdy$ that
\begin{equation}
\label{e:harmDef}
(-\Delta +M^2)\mc{H}_M^\Omega \fbdy=0\ \ \text{in }\Omega,\ (\mc{H}_M^\Omega \fbdy)|_{\partial \Omega}=\fbdy,
\end{equation}
if such $\mc{H}_M^\Omega \fbdy$ exists and is unique within $\mc{Y}(\Omega)$.
\end{definition}
\begin{definition}[Dirichlet resolvent]
    Let $\Omega$ be a uniformly Lipschitz domain. Let $\mc{X}(\Omega)$ and $\mc{Y}(\Omega)$ be Banach spaces of functions on $\Omega$. Let $M\ge 0$. We define $\mc{R}^\Omega_M: \mc{X}(\Omega)\rightarrow \mc{Y}(\Omega)$ by sending $f$ to $\mc{R}_M^\Omega f$ that
    \begin{equation}
        (-\Delta +M^2)\mc{R}^\Omega_Mf=f\ \ \text{in }\Omega,\ (\mc{R}^\Omega_Mf)|_{\partial \Omega}=0, 
    \end{equation}
    if such $\mc{R}_M^\Omega f$ exists and is unique within $\mc{Y}(\Omega)$. 
\end{definition}
We start with a formula that corrects any extension into a harmonic extension:
\begin{lemma}[Harmonic extension via correction]
Let $\Omega\subset\mathbb{R}^d$ be a uniformly Lipschitz domain. Let $\mc{X}(\partial\Omega)$ and $\mc{Y}(\Omega)$ be Banach spaces of functions on $\partial\Omega$ and $\Omega$. Let $M\ge 0$. Given any extension operator $E: \mc{X}(\partial\Omega)\rightarrow \mc{Y}(\Omega)$ that $E\cdot|_{\partial\Omega}=\Id$, we have
\begin{equation}
\label{e:hForm}
\mc{H}_M^\Omega g=E g-\mc{R}_M^\Omega (-\Delta+M^2)E g,
\end{equation}
if all three terms exist and are unique in $\mc{Y}(\Omega)$. Here we understand the $\cdot|_{\partial\Omega}$ as the non-tangential limit of $\cdot$ approaching $\partial\Omega$, which we require it to agree with the trace when the trace exists.
\end{lemma}
\begin{proof}
    It suffices to note $(-\Delta+M^2)\mc{R}_M^\Omega=\Id$ by definition.  
\end{proof}
Our main goals in this section are to establish uniform boundedness of $\mc{R}_M^\Omega: Y_m^l(\Omega)\rightarrow \dot X_m^l(\Omega)$ in Proposition~\ref{p:extendNoMass}, and $\mc{H}_M^\Omega: H^{\frac12,l}(\partial\Omega)\rightarrow X_m^l(\Omega)$ in Proposition~\ref{p:harmonic}, as $M\rightarrow 0$, where $X_m^l$, $Y_m^l$ (defined in Definition~\ref{d:Xspaces}) are suitably weighted Sobolev-type spaces with weight parameters $m,l$ in intervals controlled by the Hardy constants of $\Omega$.

\subsection{The elliptic theory on bounded uniformly Lipschitz domains}
In this subsection, we include the properties of harmonic extension and Dirichlet/Neumann resolvents on bounded uniformly Lipschitz domains. 
\begin{lemma}[Dirichlet resolvent on bounded domains]
\label{l:reg2}
Let $\Omega\Subset \mathbb{R}^d$ be a bounded uniformly Lipschitz domain. Let $M_0>0$. Let $0<s< 1$. There is $C>0$ such that for all $f\in H^{s-\frac32}(\Omega)$ and $\fbdy\in H^s(\partial \Omega)$, there is a unique $u\in H^{s+\frac{1}{2}}(\Omega)$ satisfying
\begin{equation}
    (-\Delta+M^2) u=f,\ u|_{\partial\Omega}=\fbdy,
\end{equation}
and
\begin{equation}
    \|u\|_{H^{s+\frac{1}{2}}(\Omega)}\leq C\|f\|_{H^{s-\frac32}(\Omega)}+ C\|\fbdy\|_{H^{s}(\partial\Omega)},
\end{equation}
uniformly for all such $f, \fbdy$ and for all $M\in [0,M_0]$. When $f=0$, the above is true for any $0\le s\le 1$. This implies $\mc{H}_M^\Omega: H^{s}(\partial\Omega)\rightarrow H^{s+\frac12}(\Omega)$ for $0\le s\le 1$, and $\mc{R}_M^\Omega: H^{s-\frac32}(\Omega)\rightarrow H^{s+\frac12}(\Omega)$ for $0<s<1$. 
\end{lemma}
\begin{proof}
1. Let $M=0$. The $0<s<1$ case follows from \cite[Theorem~0.5(b)]{JeKe:95} and the extension operators from $H^s(\partial\Omega)$ to $H^{s+\frac12}(\Omega)$ is well-defined. For $f=0$, the endpoint case $s=0$ follows from \cite{JeKe:81} and \cite[Corollary~5.5]{JeKe:95}, where we understand $u|_{\partial\Omega}$ as the non-tangential limit of $u$ towards $\partial\Omega$. The case $s=1$ follows from \cite{JeKe:81b}. 

2. Let $M>0$. We first note that Lax--Milgram gives that $\mc{R}_M^\Omega: H^{-1}(\Omega)\rightarrow \dot H^1(\Omega)$ uniformly in $M\in [0, M_0]$. The correction formula \eqref{e:hForm} gives a uniform  bound in $\mc{H}_M^\Omega=\mc{H}_0^\Omega-M^2\mc{R}_M^\Omega\mc{H}_0^\Omega: H^\frac12(\partial\Omega)\rightarrow  H^1(\Omega)$. Now for $f,\fbdy$ as stated, note $\mc{H}_0^\Omega f_\partial, \mc{R}_0^\Omega f\in H^{s+\frac12}(\Omega)\subset H^{\frac12}(\Omega)$. The basic bounds above then give
\begin{equation}
    u=\mc{R}_M^\Omega f+\mc{H}_M^\Omega f_\partial=\mc{R}_0^\Omega f-M^2\mc{R}_M^\Omega\mc{R}_0^\Omega f+\mc{H}_0^\Omega f_\partial-M^2\mc{R}_M^\Omega\mc{H}_0^\Omega f_\partial\in H^{\frac12}(\Omega)
\end{equation}
to exist. Since \cite[Theorem~B]{JeKe:95} implies $\mc{R}_0^\Omega: L^2(\Omega)\rightarrow \dot H^{\frac32}(\Omega)$ bounded, we have $\mc{R}_M^\Omega=\mc{R}_0^\Omega-M^2 \mc{R}_0^\Omega \mc{R}_M^\Omega$ also bounded on $L^2(\Omega)\rightarrow \dot H^{\frac32}(\Omega)$. Then
\begin{equation}
    \mc{H}_M^\Omega f_\partial=\mc{H}_0^\Omega f_\partial-M^2\mc{R}_0^\Omega\mc{H}_M^\Omega f_\partial\in H^{s+\frac12}(\Omega)
\end{equation}
for any $0\le s\le 1$. Similarly we have
\begin{equation}
    \mc{R}_M^\Omega f=\mc{R}_0^\Omega f-M^2 \mc{R}_0^\Omega \mc{R}_M^\Omega f\in H^{s+\frac12}(\Omega)
\end{equation}
for any $0<s<1$. They conclude the proof.  
\end{proof}

\begin{lemma}[Neumann resolvent on bounded domains]
\label{l:reg3}
Let $\Omega\Subset \mathbb{R}^d$ be a bounded uniformly Lipschitz domain. Let $M_0>0$. Let $0<s< 1$. Let $u\in H^{s+\frac{1}{2}}(\Omega)$, $f\in \dot H^{s-\frac32}(\Omega)$, and $\fbdy\in H^{s-1}(\partial \Omega)$, satisfying
\begin{equation}
    (-\Delta+M^2) u=f,\ \partial_\nu u|_{\partial\Omega}=\fbdy. 
\end{equation}
Then 
\begin{equation}\label{e:2l18}
    \|u\|_{H^{s+\frac{1}{2}}(\Omega)}\leq C\|f\|_{H^{s-\frac32}(\Omega)}+C\|\fbdy\|_{H^{s-1}(\partial\Omega)}+C\|u\|_{H^{-\frac12+0}(\Omega)},
\end{equation}
uniformly for all such $f, \fbdy$ and for all $M\in [0,M_0]$. When $f=0$, the above is true for any $0\le s\le 1$. 
\end{lemma}
\begin{remark}
    One here can replace the last term in \eqref{e:2l18} by $C|\langle u,1\rangle|$ and show 
    \begin{equation}
        \|u\|_{H^{s}(\partial\Omega)}\leq C\|f\|_{H^{s-\frac32}(\Omega)}+C\|f_\partial\|_{H^{s-1}(\partial\Omega)}+C\|u\|_{H^{-1}(\partial\Omega)},
    \end{equation}
    in the above-mentioned scenarios, by using \cite[Theorem~3.6, Corollary~5.7]{Behrndt2025SharpDomains}. 
\end{remark}
\begin{proof}
1. We first address $M=0$. The cases $f=0$, and $f_\partial=0$ are known to hold. The homogeneous problem $s=0$ was treated in \cite[Proposition~4.2]{Verchota2005TheDomains}, and $s=1$ was treated in \cite{JeKe:81b} under the condition $\langle f_\partial, 1\rangle_{\partial\Omega}=0$. The inhomogeneous problem was treated in \cite[Theorem~9.2]{Fabes1998BoundaryDomains} for $\langle f, 1\rangle=0$. We will use those results to allow general $f,\fbdy$ with $\langle f, 1\rangle+\langle\fbdy, 1\rangle_{\partial\Omega}=0$ when $0<s<1$. Denoting the volume of $\Omega$ by $|\Omega|$, we build the average corrector $h(z)\in C_c^\infty(\overline\Omega)$ by
\begin{equation}
    h(z)=-\frac1{2d}|z|^2, \ -\Delta h=1, \ h_\partial=\partial_\nu h\in L^2(\partial\Omega), \ \langle h_\partial,1\rangle_{\partial\Omega}=\langle \Delta h,1\rangle=-|\Omega|. 
\end{equation}
The cases discussed above imply there exist $u_1,u_2\in H^{s+\frac12}(\Omega)$ that
\begin{gather}
    -\Delta u_1=f-\frac1{|\Omega|}\langle f, 1\rangle, \ \partial_\nu u_1=0,\\
    -\Delta u_2=0, \ \partial_\nu u_2=f_\partial-\frac1{|\Omega|}\langle f,1\rangle h_\partial.\label{e:2l12}
\end{gather}
Now setting $\tilde u=u_1+u_2+\frac{1}{|\Omega|}\langle f, 1\rangle h-\frac1{|\Omega|}\langle u_1+u_2+\frac{1}{|\Omega|}\langle f, 1\rangle h, 1\rangle$ gives $u=\tilde u+\frac1{|\Omega|}\langle u,1\rangle$ and thus
\begin{equation}
    \|u\|_{H^{s+\frac{1}{2}}(\Omega)}\leq C\|f\|_{H^{s-\frac32}(\Omega)}+C\|\fbdy\|_{H^{s-1}(\partial\Omega)}+C\|u\|_{H^{-\frac12+0}(\Omega)}.
\end{equation}

2a. For $M>0$, note $-\Delta u=f-M^2 u$ and applying the case $M=0$ gives the proof for $0<s<1$. It remains to address the endpoint cases. 

2b. Let $f=0$ and $s=1$. We have $-\Delta u=-M^2 u$ for $u\in H^{\frac32}$. First let $u_3=-M^2\mc{R}_0^\Omega u\in \dot H^{\frac32}(\Omega)$. Then Step 1 with $M=0, s=\frac12$ together with \cite[Theorem~5.4]{Behrndt2025SharpDomains} implies
\begin{gather}
    \|u_3\|_{H^{\frac32}}\le CM^2\|u\|\le CM^2\|f_\partial\|_{H^{-\frac12}(\partial\Omega)}+CM^2\|u\|_{H^{-\frac12+0}}, \label{e:2l10}\\
    \|\partial_\nu u_3\|_{L^2(\partial\Omega)}\le C\|u_3\|_{H^{\frac32}}+CM^2\|u\|\le CM^2\|f_\partial\|_{H^{-\frac12}(\partial\Omega)}+CM^2\|u\|_{H^{-\frac12+0}}.\label{e:2l11}
\end{gather}
Set $u_4=u-u_3$, we have $-\Delta u_4=0$, $\partial_\nu u_4=f_\partial-\partial_\nu u_3$. Apply Step 1 with $M=0, s=1$ to see
\begin{equation}
    \|u_4\|_{H^{\frac32}}\le C\|f_\partial\|_{L^2(\partial\Omega)}+C\|u\|_{H^{-\frac12+0}}.
\end{equation}

2c. Let $f=0$ and $s=0$, where $u\in H^{\frac12}$. Similar to \eqref{e:2l12}, there exists a unique $u_5\in H^{\frac12}(\Omega)$ with $\langle u_5, 1\rangle=0$ that 
\begin{equation}
    -\Delta u_5=0, \ \partial_\nu u_5=f_\partial-\frac{1}{|\partial\Omega|}\langle f_\partial, 1\rangle_{\partial\Omega},
\end{equation}
by \cite[Proposition~4.2]{Verchota2005TheDomains}, and we have
\begin{equation}
    \|u_5\|_{H^{\frac12}}\le C\|f_\partial-\frac{1}{|\partial\Omega|}\langle f_\partial, 1\rangle_{\partial\Omega}\|_{H^{-1}(\partial\Omega)}\le C\|f_\partial\|_{H^{-1}(\partial\Omega)}. 
\end{equation}
On another hand, $u_6=u-u_5\in H^{\frac12}$ satisfies
\begin{equation}
    (-\Delta+M^2)u_6=-M^2 u_5\in H^{\frac12}, \ \partial_\nu u_6=\frac{1}{|\partial\Omega|}\langle f_\partial, 1\rangle_{\partial\Omega}\in H^1(\partial\Omega).
\end{equation}
We can use Step 2a with $M>0, s=0+0$ to estimate
\begin{equation}
    \|u_6\|_{H^{\frac12+0}}\le CM^2\| u_5\|_{H^{-\frac32+0}}+\frac{C}{|\partial\Omega|}\|\langle f_\partial, 1\rangle_{\partial\Omega}\|_{H^{-1+0}}+C\|u_6\|_{H^{-\frac12+0}}\le C\|f_\partial\|_{H^{-1}(\partial\Omega)}+C\|u\|_{H^{-\frac12+0}}
\end{equation}
as desired. 
\end{proof}

\subsection{Spaces with Hardy-type weights}
In order to specify the boundedness of $\mc{H}_M^\Omega$ and $\mc{R}_M^\Omega$ on unbounded uniformly Lipschitz domains, we introduce the following space that contains a Hardy-type weight that relates $\nabla$ to the regularized distance $\dd$ to the boundary:
\begin{definition}[Spaces with Hardy-type weights]\label{d:Xspaces}
    Let $\Omega\subset \mathbb{R}^d$ be a uniformly Lipschitz domain. Let $m,l\in \mathbb{R}$. We define
    \begin{gather}
        X_m^l(\Omega;\mc{V})=\{u\in  \rr^{-l}\dd^{l-m+1}H^{1}(\Omega): \rr^{l}\dd^{m-l}\nabla u\in L^2(\Omega),\ \rr^{l}\dd^{m-1-l}u\in L^2(\Omega)\},\\
        \dot X_m^l(\Omega;\mc{V})=\{u\in  \rr^{-l}\dd^{l-m+1}\dot H^{1}(\Omega): \rr ^{l}\dd^{m-l} \nabla u\in L^2(\Omega),\ \rr^{l}\dd^{m-1-l}u\in L^2(\Omega)\},
    \end{gather}
    where $\dot H^{1}(\Omega)$ denote $H^1(\Omega)$-elements with vanishing traces over $\partial\Omega$. We drop $\mc{V}$ and write $X_m^l(\Omega)$ whenever there is no confusion about $\mc{V}$. The spaces are equipped with norm
    \begin{equation}
        \|u\|_{X_m^l}^2=\|\rr^{l}\dd^{m-l-1}u\|_{L^2}^{2}+\|\rr^{l}\dd^{m-l}\nabla u\|_{L^2}^{2}.
    \end{equation}
    We define their dual spaces with respect to $L^2(\Omega)$ by
    \begin{equation}
        Y_m^l(\Omega)=(\dot X_{-m}^{-l}(\Omega))^*, \ \dot Y_m^l(\Omega)=(X_{-m}^{-l}(\Omega))^*.
    \end{equation}
\end{definition}
We also introduce the following Hardy-type constants for $\Omega$:
\begin{definition}[Radial Hardy constants]\label{d:rardy} Let $\Omega\subset \mathbb{R}^d$ be a uniformly Lipschitz domain. For each $r_0\ge0$, we define the \emph{$r_0$-regularized radial Hardy constant} by 
\begin{equation}
\label{e:rardyConstant}
\rardy^{r_0}_\Omega=\inf\Big\{\frac{\|\nabla u\|_{L^2(\Omega)}}{\|\frac{u}{\dist(z,\mc{V})+r_0}\|_{L^2(\Omega)}}: u\in \dot{H}^{1}(\Omega)\setminus\{0\}\Big\}.
\end{equation}
We call $\rardy^{0}_{\Omega}$ \emph{the (strong) radial Hardy constant}. We then define \emph{the regularized radial Hardy constant} by 
\begin{equation}
    \rardy^{\infty}_{\Omega}=\lim_{r_0\to \infty}\rardy^{r_0}_{\Omega}. 
\end{equation}
Note that $\rardy^{r_0}_{\Omega}>0$ is non-decreasing in $r_0$ and hence the limit either exists or is $\infty$. Given $\Omega\subset\tilde\Omega$, we have $\rardy^\infty_{\Omega}\ge \rardy^\infty_{\tilde\Omega}$; we further have $\rardy^{r_0}_{\Omega}\ge \rardy^{r_0}_{\tilde\Omega}$ for each $r_0\ge 0$ if they have same choice of $\mc{V}$. We require any uniformly Lipschitz domain $\Omega$ considered in this paper to satisfy $\rardy^{\infty}_\Omega>0$. 
\end{definition}
\begin{definition}[Hardy constants]\label{d:hardy} Let $\Omega\subset \mathbb{R}^d$ be a uniformly Lipschitz domain, and $U$ an open subset of $\Omega$. For each $r_0\ge0$, we define the \emph{$r_0$-regularized Hardy constant of $U$ with respect to $\partial\Omega$} by 
\begin{equation}
\label{e:hardyConstant}
\hardy^{r_0}_{U, \partial\Omega}=\inf\Big\{\frac{\|\nabla u\|_{L^2(\Omega)}}{\|\frac{u}{\dist(z, \partial\Omega)+r_0}\|_{L^2(\Omega)}}: u\in \dot{H}^{1}(U)\setminus\{0\}, u|_{\partial\Omega}=0\Big\}.
\end{equation}
We define the \emph{regularized Hardy constant of $U$ with respect to $\partial\Omega$} by 
\begin{equation}
    \hardy^{\infty}_{U, \partial\Omega}=\lim_{r_0\to \infty}\hardy^{r_0}_{U, \partial\Omega}.
\end{equation}
Note that $\hardy^{r_0}_{U, \partial\Omega}>0$ is non-decreasing in $r_0$ and hence the limit either exists or is $\infty$. Also for $U\subset\tilde U$, we have $\hardy^{r_0}_{U, \partial\Omega}\ge \hardy^{r_0}_{\tilde U, \partial\Omega}$. We define the \emph{$r_0$-regularized Hardy constant of $\Omega$} by 
\begin{equation}
\label{d:hardyr}
    \hardy^{r_0}_{\Omega}:=\hardy^{r_0}_{\Omega,\partial\Omega}=\inf\{\hardy^{r_0}_{U, \partial\Omega}: U\subset\Omega\}>0. 
\end{equation}
We call $\hardy^{0}_{\Omega}$ the \emph{(strong) Hardy constant} and $\hardy^{\infty}_{\Omega}$ the \emph{regularized Hardy constant}. For each $r_0\ge 0$, we have $\rardy_\Omega^{r_0}\ge \hardy_\Omega^{r_0}$. For each $\gamma>0$ we define the conic collar
\begin{equation}
    U^{r_0}_\gamma=\{z\in\Omega, \dist(z,\partial\Omega)+r_0< \gamma(\dist(z,\mc{V})+r_0)\}.
\end{equation}
The constant $\hardy^{r_0}_{U^{r_0}_\gamma, \partial\Omega}$ will appear as the threshold for $|m-l|$ in this paper. 
\end{definition}

\begin{remark}\label{r:fixr0}
    Throughout the paper, we fix some $r_0\ge 0$ if not specified otherwise. 
\end{remark}

\begin{remark}
    Here $\rardy^\infty_\Omega+\frac12-0$ is the maximal decay one can track at the zero frequency for one specific domain $\Omega$ individually, and $\rardyr_\Omega+\hardy^{r_0}_{U^{r_0}_\gamma, \partial\Omega}-0$ is the maximal decay one can trace at the zero frequency uniformly for all domains with controlled $\rardyr_{\Omega}$, $\gamma$. For $\mc{W}_\alpha$, an infinite planar wedge of corner angle $\alpha\in (0,2\pi)$, we have $\rardy^{0}_{\mc{W}_\alpha}=\frac{\pi}{\alpha}$ (see Lemma~\ref{l:conerardy}). 
\end{remark}

\begin{remark}\label{r:ydiv}
    The Hardy inequality implies
    \begin{equation}
        \dot{X}_m^l=\{u\in  \rr^{-l}\dd^{l-m+1}\dot H^{1}(\Omega): \dd^{m-l}\rr ^{l}\nabla u\in L^2(\Omega)\}, \ Y_m^l=\div(\dd^{l-m}\rr^{-l}L^2(\Omega;T\Omega)).
    \end{equation}
\end{remark}
\begin{lemma}[Properties for $X_m^l$-spaces]\label{l:X-SpacesProp}
Let $\Omega$ be a uniformly Lipschitz domain. Let $m,l\in \mathbb{R}$. The following statements are true:
\begin{gather}
    X_m^l(\Omega)=\rr^{-l} \dd^{l-m}X_{0}^{0}(\Omega), \ Y_m^l(\Omega)=\rr^{-l} \dd^{l-m} Y_{0}^{0}(\Omega),\label{e:2l5}\\
    \rr^{-l}\dd^{l-m} H^1(\Omega)\subset X_m^l(\Omega)\subset \rr^{-l}\dd^{l-m+1} H^1(\Omega), \\
    \rr^{-l}\dd^{l-m-1} H^{-1}(\Omega)\subset Y_m^l(\Omega)\subset \rr^{-l}\dd^{l-m} H^{-1}(\Omega), \\
    X_m^l(\Omega)\subset Y_{m-2}^l(\Omega). 
\end{gather}
The statements are also true for $\dot X_m^l$ and $\dot Y_m^l$ when $H^1, H^{-1}$ are replaced with $\dot H^{1}, \dot H^{-1}$.
\end{lemma}
\begin{proof}
    Note $\rr^l\dd^{m-l}\nabla \rr^{-l}\dd^{-m+l}=\nabla+O(\dd^{-1})$ to see
    \eqref{e:2l5}. The rest follows similarly. 
\end{proof}

\subsection{The Poisson problem}
We study the solution $u\in \dot{H}^{1,l}(\Omega)$ to the Poisson problem
\begin{equation}
\label{e:posM}
(-\Delta+M^2)u=f,
\end{equation}
for $f\in H^{-1,l}(\Omega)$. For each $M>0$, the solution $u$ to~\eqref{e:posM} exists by the Lax--Milgram Theorem. However, as $M\to 0$, estimates for $u\in \dot{H}^{1,l}$ degenerate. In this subsection, we study uniform estimates for~\eqref{e:posM} in $X_m^l$ instead of $\dot H^{1,l}$ as $M\to 0$. 
\begin{lemma}[Stability for the Poisson problem]
\label{l:stableEstimates}
Let $r_0\geq 0$, $\gamma\in(0,\frac12)$, $\delta_0\in(0,1)$. Then for all $\Omega\subset\mathbb{R}^d$ a uniformly Lipschitz domain with $|m|\le\delta_0\rardy_{\Omega}^{r_0}$, $|m-l|\le \delta_0\hardy_{U^{r_0}_{\gamma},\partial\Omega}^{r_0}$, $m_1, m_2\in \mathbb{R}$, we have
\begin{equation}
\|\langle M\rr\rangle^{m_1}\langle M\dd\rangle^{m_2}u\|_{X_{m}^{l}}+M\|\langle M\rr\rangle ^{m_1}\langle M\dd\rangle ^{m_2} \rr^l \dd^{m-l} u\|_{L^2}
\leq C\|\langle M\rr\rangle ^{m_1}\langle M\dd\rangle ^{m_2}(-\Delta+M^2)u\|_{Y_m^l}.
\end{equation}
uniformly for $u\in \langle M\rr\rangle^{-m_1}\langle M\dd\rangle^{-m_2}\rr^{-l}\dd^{l-m} \dot{H}^{1}(\Omega)$ and all $M\ge0$.
\end{lemma}
\begin{proof}
1. Fix $r_0'>0$ that $r_0'\ge r_0$. Only in this proof, we use the following notation for convenience:
\begin{equation}
    r=\dist(z,\mc{V})+r_0',\ d=\dist(z,\partial\Omega)+r_0'.
\end{equation}
Note $r,d$ are only $C^{0,1}(\Omega)$, and $r\sim\rr$, $d\sim\dd$. Let $\delta_1=1-\delta_0\in (0,1)$ for convenience. By Definition~\ref{d:hardy}, Definition~\ref{d:rardy}, and $\rardy_{\Omega}^{r_0'}\ge \rardy_{\Omega}^{r_0}$, $\hardy_{U^{r_0'}_\gamma,\partial\Omega}^{r_0'}\ge \hardy_{U^{r_0}_\gamma,\partial\Omega}^{r_0}$, we have
\begin{gather}\label{e:2l7}
    (1-\frac14\delta_1)^2 (\rardy_{\Omega}^{r_0})^2\|\frac{u}{r}\|_{L^2(\Omega)}^2\le (1-\frac18\delta_1)\|\nabla u\|_{L^2(\Omega)}^2,\ u\in \dot H^1(\Omega)\\
    (1-\frac14\delta_1)^2 (\hardy_{U^{r_0}_\gamma,\partial\Omega}^{r_0})^2\|\frac{u}{d}\|_{L^2(\Omega)}^2\le (1-\frac18\delta_1)\|\nabla u\|_{L^2(\Omega)}^2,\ u\in \dot H^1(U^{r_0}_\gamma), \ u|_{\partial\Omega}=0.
\end{gather}
All constants in the lemma only depend on $\delta_0, r_0,\gamma$, the choices of $m, l,m_1, m_2$, and lower bounds on $\hardy_{U^{r_0}_\gamma,\partial\Omega}^{r_0}$, $\rardy_\Omega^{r_0}$.

2. We construct the weight functions $\varphi_M\in C^{0,1}(\Omega)$. Set $K>1$, $K_1>\gamma^{-1}$ large to be determined later. Let 
\begin{equation}
\label{e:weightDef}
\varphi_M = m\log r+(m-l)\log (\frac{d/r}{1+Kd/r})
+ m_1\log (Mr+K(|m_1|+1))+m_2\log (Md+K(|m_2|+1)).
\end{equation}
Noting $1+Kd/r$ is bounded from above and below by positive constants we have
\begin{equation}\label{e:2l8}
    e^{\varphi_M}\underset{r_0,\delta_0}{\sim} \langle M\rr\rangle ^{m_1}\langle M\dd\rangle ^{m_2} \rr^l\dd^{m-l} ,\ |\nabla e^{\varphi_M}|\underset{r_0,\delta_0}{\le} \langle M\rr\rangle ^{m_1}\langle M\dd\rangle ^{m_2}\rr^l\dd^{m-l-1}.
\end{equation}
Let $\chi_0, \chi_1\in C^{0,1}(\Omega)$ that $\chi_0^2+\chi_1^2\equiv 1$ and $\chi_1\equiv 0$ on $\{\log(\frac{r}d)\ge K_1+\log(K_1)\}=\{\frac{d}r\le K_1^{-1}e^{-K_1}\}$ and $\chi_1\equiv 1$ on $\{\log(\frac{r}d)\le \log(K_1)\}=\{\frac{d}r\ge K_1^{-1}\}$, and 
\begin{equation}\label{e:2l19}
    |\nabla \chi_0|\le \frac{C}{K_1 d}, \ |\nabla \chi_1|\le \frac{C}{K_1 d}
\end{equation}
for some $C>0$ uniformly. This makes $\supp\chi_0\subset\{d\le \gamma r\}$. We claim choosing large $K, K_1$ gives
\begin{equation}
\label{e:weight}
|\nabla \varphi_M|^2\leq(1-\frac14\delta_1)^2\Big(\frac{\chi_0^2(\hardy_{U^{r_0}_\gamma,\partial\Omega}^{r_0})^2}{d^{2}}+\frac{\chi_1^2(\rardy_\Omega^{r_0})^2}{r^{2}}\Big)+\frac14 M^2
\end{equation}
in this step. To see it, we compute
\begin{equation}
    \nabla \varphi_M=\frac{m\nabla r}{r}+\frac{m-l}{1+Kd/r}\nabla \log(\frac{d}r)+\frac{m_1M\nabla r}{Mr+K(|m_1|+1)}+\frac{m_2M\nabla d}{Md+K(|m_2|+1)},
\end{equation}
the last two terms of which can be made small by a choosing a large $K$ since 
\begin{equation}
    \Big|\frac{m_1M}{Mr+K(|m_1|+1)}\nabla r\Big|+\Big|\frac{m_2M}{Md+K(|m_2|+1)}\nabla d\Big|\leq \frac{2M}{K}.
\end{equation}
We now analyze the first two terms. On $\supp\chi_0$ we have $\frac {d}{r}\le K_1^{-1}$ and thus
\begin{multline}
    \chi_0^2|\nabla \varphi_M|^2\le \chi_0^2\Big((1+\e)\Big|\frac{m\nabla r}{r}+\frac{m-l}{1+Kd/r}\nabla \log(\frac{d}r)\Big|^2+\frac{C_\e M^2}{K^2}\Big)\\
    \le \chi_0^2\Big((1+2\e)(1+K_1^{-1})^2\frac{(m-l)^2}{d^2}+\frac{C_\e m^2}{K_1^2d^2}+\frac{C_\e M^2}{K^2}\Big),
\end{multline}
where we used $|\nabla d|\le1$, $|\nabla r|\le1$ and on $\supp\chi_0$ we have
\begin{equation}
    \Big|\frac{\nabla r}{r}\Big|\le \frac{1}{K_1d},\ \Big|\frac{1}{1+Kd/r}\nabla \log(\frac{d}r)\Big|\le \Big|\frac{\nabla d}{d}-\frac{\nabla r}{r}\Big|\le \frac1d+\frac1r\le \frac{1+K_1^{-1}}{d}. 
\end{equation}
On $\supp\chi_1$ we have $\frac {d}{r}\ge K_1^{-1}e^{-K_1}$ and
\begin{multline}
    \chi_1^2|\nabla \varphi_M|^2\le \chi_1^2\Big((1+\e)\Big|\frac{m\nabla r}{r}+\frac{m-l}{1+Kd/r}\nabla \log(\frac{d}r)\Big|^2+\frac{C_\e M^2}{K^2}\Big)\\
    \le \chi_1^2\Big((1+2\e)\frac{m^2}{r^2}+\frac{C_\e(1+K_1e^{K_1})^2 (m-l)^2}{(1+K K_1^{-1}e^{-K_1})^2r^2}+\frac{C_\e M^2}{K^2}\Big),
\end{multline}
where we used
\begin{multline}
    \Big|\frac{1}{1+Kd/r}\nabla \log(\frac{d}r)\Big|\le \frac{1}{1+K K_1^{-1}e^{-K_1}}\Big|\frac{\nabla d}{d}-\frac{\nabla r}{r}\Big|\le \frac{1}{1+K K_1^{-1}e^{-K_1}}(\frac{1}{d}+\frac{1}{r})\\
    \le \frac{1+K_1e^{K_1}}{(1+K K_1^{-1}e^{-K_1})r}.
\end{multline}
Now sequentially choosing small $\e$, large $K_1$, large $K$ gives \eqref{e:weight}. 

3. Consider the conjugated operator
\begin{equation}
    P_{\varphi_M}:=e^{\varphi_M}(-\Delta+M^2)e^{-\varphi_M}.
\end{equation}
We claim $P_{\varphi_M}$ is elliptic in the sense that there is $c>0$ small that
\begin{equation}\label{e:2l6}
    \Re \langle P_{\varphi_M} v,v\rangle\geq c(\|v\|^2_{X^0_0}+M^2\|v\|_{L^2}^2)
\end{equation}
uniformly for all $v\in \dot{H}^{1}(\Omega)$ and $M\ge 0$. Indeed, 
\begin{multline}
2\Re \langle P_{\varphi_M} v,v\rangle =2 \|\nabla v\|^2 +2\langle \langle \nabla \varphi_M,\nabla v\rangle, v\rangle +2\langle v,\langle \nabla \varphi_M,\nabla v\rangle \rangle +2\langle (M^2-\nabla^*\nabla \varphi_M -|\nabla \varphi_M|^2 )v,v\rangle \\
=2 \|\nabla v\|^2 +2\langle (M^2-|\nabla \varphi_M|^2 )v,v\rangle=2\|\nabla v\|^2+2M^2\|v\|^2-2\||\nabla\varphi_M|v\|^2.
\end{multline}
Making $K_1>\gamma^{-1}$ larger if necessary, we note Step 1 and 2 give
\begin{multline}
    \||\nabla\varphi_M|v\|^2\le (1-\frac14\delta_1)^2\Big((\hardy^{r_0}_{U^{r_0}_\gamma,\partial\Omega})^2\Big\|\frac{\chi_0 v}{d}\Big\|^2+(\rardy^{r_0}_{\Omega})^2\Big\|\frac{\chi_1 v}{r}\Big\|^2)+\frac14 M^2\|v\|^2\\\le (1-\frac18\delta_1)\Big(\Big\|\nabla (\chi_0v)\Big\|^2+\Big\|\nabla (\chi_1v)\Big\|^2)+\frac14 M^2\|v\|^2\le (1-\frac18\delta_1)\|\nabla v\|^2+\frac{C}{K_1}\Big\|\frac{v}{d}\Big\|^2+\frac14 M^2\|v\|^2\\
    \le (1-\frac1{16}\delta_1)\|\nabla v\|^2+\frac14 M^2\|v\|^2,
\end{multline}
where we used \eqref{e:2l19}, \eqref{e:weight}, and $\chi_0 v, (\nabla \chi_0)v, (\nabla \chi_1)v$ are supported inside $U^{r_0}_\gamma\cup (\overline{U^{r_0}_\gamma}\cap\partial\Omega)$.
We now have for some $C>0$ large independent of $M$ that
\begin{equation}
2\Re \langle P_{\varphi_M} v,v\rangle \geq \frac{\delta_1}8\|\nabla v\|^2+\frac{3M^2}{2}\|v\|^2\ge C^{-1}\Big(\|\nabla v\|^2_{L^2}+\Big\|\frac{v}{\dd}\Big\|_{L^2}^2+M^2\|v\|_{L^2}^2\Big).
\end{equation}
This gives \eqref{e:2l6} as desired by noting $\|v\|_{X_{0}^0}^2=\|\nabla v\|_{L^2}^2+\|\dd^{-1}v\|_{L^2}^2$.

4. Now we take $v=e^{\varphi_M}u\in \dot{H}^{1}(\Omega)$. Step 3 gives
\begin{equation}
\|v\|^2_{X_{0}^0}+M^2\|v\|^2_{L^2}\leq C|\langle e^{\varphi_{M}}(-\Delta+M^2)u,v\rangle|\leq C\e^{-1}\| e^{\varphi_{M}}(-\Delta+M^2)u\|_{Y_0^0}^2+\e\|v\|^2_{X_{0}^0}.
\end{equation}
The absorption of the last term on the right concludes the proof with \eqref{e:2l8}. 
\end{proof}

We now show $\mc{R}_M^\Omega$, the Dirichlet resolvent, is uniformly bounded in $M\ge 0$ as an operator from $Y_m^l$ to $\dot X^l_m$, and show $\mc{R}_M^\Omega\rightarrow \mc{R}_0^\Omega$ as a weak limit.

\begin{proposition}[Zero mass limit of Dirichlet resolvent]
\label{p:extendNoMass}
Let $r_0\geq 0$, $\gamma\in(0,\frac12)$. Let $m,l\in\mathbb{R}$. Let $\Omega\subset\mathbb{R}^d$ be any uniformly Lipschitz domain with $|m|<\rardyr_\Omega$, $|m-l|<\hardyr_{U_{\gamma}^{r_0},\partial\Omega}$. Let $m_1,m_2\in\mathbb{R}$. Then 
\begin{equation}\label{e:2l21}
    \mc{R}^\Omega_M:\langle M\rr\rangle^{-m_1}\langle M\dd\rangle^{-m_2} Y_{m}^{l}(\Omega)\to \langle M\rr\rangle^{-m_1}\langle M\dd\rangle^{-m_2} \dot X_m^l(\Omega)
\end{equation}
exists and is uniformly bounded for all $M\ge 0$, where $\mc{R}^\Omega_M f$ is the unique element in that class such that $(-\Delta+M^2)\mc{R}^\Omega_M f=f$. Indeed, there is $C>0$ such that for all $M\ge 0$, all $f\in \langle M\rr\rangle^{-m_1}\langle M\dd\rangle^{-m_2} Y_{m}^{l}(\Omega)$, we have
\begin{equation}\label{e:2l20}
 \|\langle M\rr\rangle^{m_1}\langle M\dd\rangle^{m_2}\mc{R}^\Omega_Mf \|_{X_{m}^{l}}+M\|\langle M\rr\rangle^{m_1}\langle M\dd\rangle^{m_2}\rr^{l}\dd^{m-l}\mc{R}^\Omega_Mf \|_{L^2}
 \leq C\|\langle M\rr\rangle^{m_1}\langle M\dd\rangle^{m_2}f\|_{Y_m^l}.
\end{equation}
Moreover $\mc{R}_0^\Omega: Y_m^l(\Omega)\rightarrow \dot X_m^l(\Omega)$ agrees with the zero mass limit of $\mc{R}_M^\Omega: Y_m^l(\Omega)\rightarrow \dot X_m^l(\Omega)$ as $M\rightarrow 0$, 
\begin{equation}
    Y_{m}^{l}(\Omega)\ni f\mapsto \mc{R}^\Omega_0f=\lim_{M\to 0}\mc{R}^\Omega_Mf\in \dot X_m^l(\Omega)
\end{equation}
as weak convergence in $\dot X_{m}^{l}$, and in the norm sense in $\dot X_{m-0}^{l-0}$. 
The only dependence of bounds on $\Omega$ are through the lower bounds on $\hardyr_{U_{\gamma}^{r_0},\partial\Omega},\rardyr_{\Omega}$. 
\end{proposition}
\begin{remark}
    For $u\in \dot X_m^l(\Omega)$ within the specified limits for $l,m, \rardyr_{\Omega}, \hardyr_{U_{\gamma}^{r_0},\partial\Omega}$, we have $\mc{R}^\Omega_M(-\Delta+M^2)u=u$ by definition, for all $M\ge 0$.
\end{remark}
\begin{proof}
1. Set $\varphi_M$ as in~\eqref{e:weightDef}. For each $M>0$, $\rr\sim \langle M\rr\rangle$, $\dd\sim \langle M\dd\rangle$. Let 
\begin{equation}
    f\in \langle M\rr\rangle^{-m_1}\langle M\dd\rangle^{-m_2} Y_m^l\subset \langle M\rr\rangle^{-m_1}\langle M\dd\rangle^{-m_2} \rr^{-l}\dd^{-m+l} H^{-1}.
\end{equation}
Then, we have $e^{\varphi_M}f \in H^{-1}(\Omega)$ per \eqref{e:2l8}. By the Lax--Milgram theorem and \eqref{e:2l6}, there is a unique solution $v\in {\dot{H}^{1}(\Omega)}$ to $P_{\varphi_M}v= e^{\varphi_M}f$. Now, setting $u=e^{-\varphi_M}v$, we have 
\begin{equation}
    (-\Delta+M^2)u=f,\ u\in \langle M\rr\rangle^{-m_1}\langle M\dd\rangle^{-m_2}\dd^{l-m}\rr^{-l} \dot H^{1}.
\end{equation}
Thus, $u=\mc{R}^\Omega_Mf$, and applying Lemma~\ref{l:stableEstimates}, we obtain \eqref{e:2l20} and \eqref{e:2l21}.

2. Let $\delta>0$, $f\in Y_m^l$ and let $u_M=\mc{R}^\Omega_M f$ for $M\ge 0$. By \eqref{e:2l21} we have $u_M\in \dot X_m^l$ and
\begin{equation}
\sup_{M}\|u_M\|_{X_{m}^{l}}\le C\|f\|_{Y_m^l}
\end{equation}
uniformly. There is a sequence of $M\to 0$ and there is $u_0\in \dot X_{m}^{l}$ such that 
\begin{equation}\label{eq:2l1}
u_{M}\rightharpoonup u_0 \text{ weakly in }X_{m}^l,\ \dd^{m-l-1} \rr^{l-\delta}u_{M}\rightarrow\dd^{m-l-1}\rr^{l-\delta}u_0 \text{ in }L^2. 
\end{equation}
It follows that $-\Delta u_0=f$ and $u_0=\mc{R}^\Omega_0 f$ via uniqueness. Here we used the fact that the embedding $\dot X_m^l\hookrightarrow \dd^{l-m+1}\rr^{-l+\delta} L^2$ is compact. Indeed, note that $\dd\le C\rr$, we have that $\dd^{l-m}\rr^{-l}\le C\dd^{l-m-\frac12\delta}\rr^{-l+\frac12\delta}$. We then have
\begin{equation}
    \dot X_m^l=\dd^{l-m}\rr^{-l}\dot X_0^0\subset \dd^{l-m+1}\rr^{-l}\dot H^1\hookrightarrow \dd^{l-m+1}\rr^{-l+\delta}L^2
\end{equation}
compactly by the Rellich theorem. It remains to show the norm convergence of $u_{M}\rightarrow u_0$ in $\dot{X}_{m-\delta}^{l-\delta}$ in each $\delta>0$, by showing $\dd^{m-l}\rr^{l-\delta} \nabla u_{M}\rightarrow\dd^{m-l}\rr^{l-\delta}\nabla u_0$ in $L^2$. This is equivalent to showing $\dd^{m-l}\rr^{l-\delta} \nabla w_M\rightarrow 0$ for $w_M=u_M-u_0$. 

2a. For each $R>1$, let $\chi_R(z)\in C^\infty(\mathbb{R}^d)$ that $\chi_{R}\equiv 1$ in $\{\rr\le R-1\}$ and $\chi_R\equiv 0$ in $\{\rr\ge R+1\}$ such that $|\nabla \chi_R|\le C$ independently of $R$. In this step, we claim for each $\epsilon>0$, there is $R>1$ that $\|\dd^{m-l}\rr^{l-\delta} (1-\chi_R)\nabla w_M\|_{L^2}<\epsilon$ for all $M$ along the sequence. Indeed, 
\begin{equation}
    \|\dd^{m-l}\rr^{l-\delta} (1-\chi_R) \nabla  w_M\|\le CR^{-\delta}\|\dd^{m-l}\rr^{l} \nabla w_M\|\le CR^{-\delta}\|w_M\|_{X_m^l}\le CR^{-\delta}\|f\|_{Y_m^l}
\end{equation}
yielding the claim. 

2b. We now claim for such $\epsilon$ and its corresponding $R$ given in Step 2a, 
\begin{equation}
    \|\dd^{m-l}\rr^{l-\delta} \chi_R \nabla  w_M\|\rightarrow 0, \ M\rightarrow 0. 
\end{equation}
Showing this would imply the norm convergence in $\dot X_{m-\delta}^{l-\delta}$ as desired. For so, compute
\begin{equation}\label{e:2l9}
\|\rr^{l-\delta}\dd^{m-l}\chi_R\nabla w_M\|_{L^2}^2
=|\langle \nabla w_M,\nabla  \rr^{2l-2\delta} \dd^{2m-2l}\chi_R^2w_M\rangle
+\langle \nabla w_M,[ \rr^{2l-2\delta}\dd^{2m-2l}\chi_R^2,\nabla]w_M\rangle |.
\end{equation}
Noting that
\begin{equation}
[ \rr^{2l-2\delta}\dd^{2m-2l}\chi_R^2,\nabla]=O(\rr^{2l-2\delta}\dd^{2m-2l-1}),
\end{equation}
we can estimate the second term of \eqref{e:2l9} by 
\begin{multline}
|\langle \nabla w_M,[ \rr^{2l-2\delta}\dd^{2m-2l}\chi_R^2,\nabla]w_M\rangle |\le C\|\rr^{l-\delta}\dd^{m-l}\nabla w_M\|_{L^2}\|\rr^{l-\delta} \dd^{m-l-1}w_M\|_{L^2}\\
\le C\|f\|_{Y_m^l}\|\rr^{l-\delta} \dd^{m-l-1}w_M\|_{L^2}\rightarrow 0
\end{multline}
as $M\rightarrow 0$. For the first term, we note $-\Delta w_M=-M^2u_M$ gives
\begin{multline}
    |\langle \nabla w_M,\nabla  \rr^{2l-2\delta} \dd^{2m-2l}\chi_R^2w_M\rangle|\le CR^2M^2\|\rr^{l-\delta} \dd^{m-l-1}u_M\|\|\rr^{l-\delta} \dd^{m-l-1}w_M\|\\
    \le CR^2M^2\|f\|_{Y_m^l}\|\rr^{l-\delta} \dd^{m-l-1}w_M\|\rightarrow 0. 
\end{multline}
Now we have $u_{M}\rightarrow u_0$ in $\dot X_{m-\delta}^{l-\delta}$ for every $\delta>0$.
\end{proof}

\subsection{The harmonic extension}
We are now in a position to study the harmonic extension operator for $-\Delta +M^2$. 

\begin{lemma}[Extension operator]
\label{l:babyExtend}
Let $\Omega$ be a uniformly Lipschitz domain.  Then there is an extension operator $E_0$ such that for all $0\le s\le 1$, $l\in\mathbb{R}$, and $N\in\mathbb{R}$, $E_0: H^{s,l}(\partial \Omega)\to  \dd^{-N}H^{s+\frac{1}{2},l}(\Omega)$ is continuous, and for all $\fbdy\in H^{s,l}(\partial\Omega)$, $(E_0\fbdy)|_{\partial \Omega}=\fbdy$ and $\supp E_0\fbdy\subset \{\dd\le 2\}$. 
\end{lemma}

\begin{proof}
Consider a partition of unity $\chi_k\in C^\infty_c(\partial\Omega)$ of bounded local degree, with $\tilde\chi_k\in C^\infty_c(\partial\Omega)$ that $\supp\chi_k\cap\supp(1-\tilde\chi_k)=\emptyset$, $\supp\tilde\chi_k\subset\{\dd\le 2\}$. We do not distinguish functions on $\partial\Omega$ with their extensions in $\Omega$, and without loss of generality assume $\supp\tilde\chi_k\subset \tilde\Omega_k$ for some bounded uniformly Lipschitz domain $\tilde\Omega_k$ that $\partial\Omega=\partial\tilde\Omega_k$ over $\supp\tilde\chi_k$, where the Lipschitz constants and diameters of such $\tilde\Omega_k$ does not depend on $k$. Note the above makes $\mc{H}_M^{\tilde\Omega_k}: H^s(\partial\tilde\Omega_k)\rightarrow H^{s+\frac12}(\tilde\Omega_k)$ is bounded independently of $k$. Now set up $E_M=\sum_{k}\tilde\chi_k \mc{H}_M^{\tilde\Omega_k} \chi_k$ to conclude with $E_0$. 
\end{proof}

We now show $\mc{H}_M^\Omega$ exists for each $M\ge 0$, and $\mc{H}_0^\Omega$ agrees with $\lim_{M\rightarrow 0}\mc{H}_M^\Omega$.

\begin{proposition}[Zero mass limit of harmonic extension]
\label{p:harmonic}
Let $r_0\geq 0$, $\gamma\in(0,\frac12)$. Let $m,l\in\mathbb{R}$. Let $\Omega \subset \mathbb{R}^d$ be a uniformly Lipschitz domain with $|m|<\rardyr_{\Omega}$, $|m-l|<\hardyr_{U_{\gamma}^{r_0},\partial\Omega}$. Let $m_1,N\in\mathbb{R}$. Let $M_0>0$. Then
\begin{equation}\label{e:2l23}
    \mc{H}^\Omega_M=E_0-\mc{R}_M^\Omega (-\Delta+M^2)E_0:H^{\frac12, l}(\partial\Omega)\to X_m^l(\Omega)
\end{equation}
is bounded uniformly in $M\in[0,M_0]$. Indeed, there is $C>0$ such that
\begin{equation}\label{e:2l22}
 \|\langle M\rr\rangle^{m_1}\langle M\dd\rangle^{N}\mc{H}_M^\Omega \fbdy \|_{X_{m}^{l}}+M\|\langle M\rr\rangle^{m_1}\langle M\dd\rangle^{N}\rr^{l}\dd^{m-l}\mc{H}_M^\Omega \fbdy \|_{L^2}
 \leq C\|\langle M\rr\rangle^{m_1}\fbdy\|_{H^{\frac12,l}(\partial\Omega)},
\end{equation}
uniformly in $\fbdy\in \langle M\rr\rangle^{-m_1}H^{\frac12, l}(\partial\Omega)$ and all $M\in[0,M_0]$. Moreover, $\mc{H}_0^\Omega=\lim_{M\to 0}\mc{H}_M^\Omega$ in the operator norm $H^{\frac12, l}(\partial\Omega)\to X_{m}^{l}(\Omega)$, with
\begin{equation}\label{e:2l4}
     \|(\mc{H}^\Omega_M-\mc{H}^\Omega_0)\fbdy\|_{X_{m}^{l}}\leq CM^{\min(\rardyr_\Omega-m-0, \hardyr_{U_{\gamma}^{r_0},\partial\Omega}-(m-l')-0,2)}\|\fbdy\|_{H^{\frac12, l'}(\partial\Omega)}
\end{equation}
uniformly for all $M\in[0, M_0]$, where $l'\ge l$ satisfies $|l'|<\rardyr_{\Omega}+\hardyr_{U_{\gamma}^{r_0},\partial\Omega}$. 
\end{proposition}
\begin{remark}\label{r:pickm}
    Given fixed $|l|<\rardyr_{\Omega}+\hardyr_{U_{\gamma}^{r_0},\partial\Omega}$, one can choose $m,l'$ stated above that
    \begin{equation}
        \min(\rardyr_\Omega-m-0, \hardyr_{U_{\gamma}^{r_0},\partial\Omega}-(m-l')-0,2)=\min(2\rardyr_\Omega-0,  \rardyr_\Omega+\hardyr_{U_{\gamma}^{r_0},\partial\Omega}-l-0, 2).
    \end{equation}
\end{remark}
\begin{proof}
1. Note that Lemma~\ref{l:babyExtend} gives
\begin{equation}
    (-\Delta+M^2) E_0\fbdy\in \langle M \rr\rangle^{-m_1}\dd^{-N-1} H^{-1,l}\subset \langle M \rr\rangle^{-m_1}Y_{l+N}^{l},
\end{equation}
uniformly in $M\in[0, M_0]$. Applying Proposition~\ref{p:extendNoMass} to \eqref{e:2l23} gives \eqref{e:2l22} and \eqref{e:2l23} bounded. 

2. It remains to show the convergence of $\mc{H}_M^\Omega\rightarrow \mc{H}_0^\Omega$ in the operator norm as $M\rightarrow 0$. Consider
\begin{equation}
    -\Delta(\mc{H}_{M}^\Omega-\mc{H}_{0}^\Omega)\fbdy=-M^2\mc{H}_{M}^\Omega \fbdy,\  (\mc{H}_{M}^\Omega-\mc{H}_{0}^\Omega)\fbdy|_{\partial\Omega}=0.
\end{equation}
By the uniqueness of $\mc{R}^\Omega_0$, we have
\begin{equation}
   \label{e:2l2}
(\mc{H}_{M}^\Omega-\mc{H}_{0}^\Omega)\fbdy=-\mc{R}_{0}^\Omega M^2\mc{H}_{M}^\Omega \fbdy
\end{equation}
Let $\e=\min(\rardyr_\Omega-m-0, \hardyr_{U_{\gamma}^{r_0},\partial\Omega}-(m-l')-0,2)$. Now, for $|m|<\rardyr_{\Omega}$, $|m-l|<\hardyr_{U_{\gamma}^{r_0},\partial\Omega}$, we have some $m'>m+\epsilon$ such that $|m'|<\rardyr_{\Omega}$, $|m'-l'|<\hardyr_{U_{\gamma}^{r_0},\partial\Omega}$ as well. Note that $X_{m'+2-\epsilon}^{l'}\subset X_{m+2}^{l}\subset Y_m^l$ implies
\begin{equation}
    \|M^{2-\epsilon}\mc{H}_{M}^\Omega \fbdy\|_{Y_m^l}\le C\|M^{2-\epsilon}\mc{H}_{M}^\Omega \fbdy\|_{X_{m'+2-\epsilon}^{l'}}.
\end{equation}
We can use \eqref{e:2l22} to further estimate
\begin{equation}
    \|M^{2-\epsilon}\mc{H}_{M}^\Omega \fbdy\|_{X_{m'+2-\epsilon}^{l'}}\le C\|\langle M\dd\rangle^{2-\epsilon}\mc{H}_{M}^\Omega \fbdy\|_{X_{m'}^{l'}}\le C\|\fbdy\|_{H^{\frac12,l'}_\partial}.
\end{equation}
Then using Proposition~\ref{p:extendNoMass}, we can estimate \eqref{e:2l2} by
\begin{equation}
    \|(\mc{H}_{M}^\Omega -\mc{H}_{0}^\Omega)\fbdy\|_{X_{m}^l}\leq \frac{CM^2}{M^{2-\epsilon}}\|M^{2-\epsilon}\mc{H}_{M}^\Omega \fbdy\|_{Y_m^l}\leq CM^\epsilon\|\fbdy\|_{H^{\frac12,l'}_\partial}.
\end{equation}
This is \eqref{e:2l4}. We have $\|\mc{H}_M^\Omega- \mc{H}_0^\Omega\|_{H^{\frac12,l}_\partial\to X_{m}^{l}}$ goes to 0 as $M\rightarrow 0$ as desired, by taking $l'=l$. 
\end{proof}

\subsection{The \dtn map and elliptic regularity}
As an immediate corollary of Proposition~\ref{p:harmonic}, we give a definition of the \dtn map.   
\begin{definition}[\dtn map]\label{d:dtn}
    Let $\Omega\subset\mathbb{R}^d$ be a uniformly Lipschitz domain, $|l|<\rardyr_\Omega+\hardyr_{U_{\gamma}^{r_0},\partial\Omega}$. For each $M\ge 0$, we define the Dirichlet-to-Neumann map, $\Lambda_M^\Omega:H^{\frac12,l}(\partial\Omega)\to H^{-\frac12,l}(\partial\Omega)$ as a bounded map given by the pairing
\begin{equation}
\label{e:quadraticDtN}
\langle \Lambda_M^\Omega u,\phi\rangle_{\partial\Omega}= \langle \nabla\mc{H}^\Omega_Mu,\nabla\mc{H}^\Omega_M\phi\rangle_{\Omega} + M^2\langle \mc{H}_M^\Omega u,\mc{H}_M^\Omega \phi\rangle_{\Omega}, \ u\in H^{\frac12, l}, \ \phi\in H^{\frac12,-l}, 
\end{equation}
or equivalently $\Lambda_M^\Omega=\partial_\nu \mc{H}_M^\Omega$. Here $\Lambda_M^\Omega$ is self-adjoint on $L^2(\partial\Omega)$. Moreover $\Lambda_0^\Omega=\lim_{M\rightarrow 0}\Lambda_M^\Omega$ in norm. 
\end{definition}

We prove some useful elliptic estimates to be used later. 
\begin{lemma}
\label{l:nearBoundary}
Let $\Omega\subset \mathbb{R}^d$ be a uniformly Lipschitz domain. Let $M_0>0$. Let $\psi,\tilde{\psi} \in C_c^\infty(\mathbb{R})$ with $\supp \psi\cap \supp (1-\tilde{\psi})=\emptyset$ and $\psi(\dd)\equiv 1$ near $\{\dd\le 2\}$. 
Then, for $0<s<1$, and $l\in\mathbb{R}$,
\begin{equation}\label{e:4l29}
    \|\rr^l\psi(\dd) u\|_{H^{s+\frac12}(\Omega)}\leq C\|\rr^l \tilde\psi(\dd)(-\Delta+M^2)u\|_{H^{s-\frac32}(\Omega)}+C\|\rr^l u\|_{H^s(\partial\Omega)}+C\|\rr^l\tilde{\psi}(\dd) u\|_{H^{-\frac12+0}(\Omega)},
\end{equation}
uniformly for all $u\in \rr^{\infty}H^{\frac12+0}(\Omega)$. When $s=1$, we have
\begin{equation}
    \|\rr^l\psi(\dd) u\|_{H^{\frac32}(\Omega)}\leq C\|\rr^l \tilde\psi(\dd)(-\Delta+M^2)u\|_{L^2(\Omega)}+C\|\rr^l u\|_{H^1(\partial\Omega)}+C\|\rr^l\tilde{\psi}(\dd) u\|_{H^{-\frac12+0}(\Omega)},
\end{equation}
for all $u\in \rr^{\infty}H^{\frac12+0}(\Omega)$. When $s=0$, we have
\begin{equation}
    \|\rr^l\psi(\dd) u\|_{H^{\frac12}(\Omega)}\leq C\|\rr^l \tilde\psi(\dd)(-\Delta+M^2)u\|_{H^{-\frac32+0}(\Omega)}+C\|\rr^l u\|_{L^2(\partial\Omega)}+C\|\rr^l\tilde{\psi}(\dd) u\|_{H^{-\frac12+0}(\Omega)},
\end{equation}
for all $u\in \rr^{\infty}H^{\frac12}(\Omega)$ with non-tangential limits almost everywhere on $\partial\Omega$. The above estimates are uniform in $M\in [0,M_0]$. Moreover, the above are true for all $0\le s\le 1$ by replacing $\|\rr^l u\|_{H^s(\partial\Omega)}$ on the right by $\|\rr^l \partial_\nu u\|_{H^{s-1}(\partial\Omega)}$.   
\end{lemma}
\begin{proof}
1. Since $\Omega$ is uniformly Lipschitz, we may find $\chi_k$ uniformly bounded in $C_c^\infty(\overline{\Omega})$ a locally finite partition of unity on $\supp \psi(\dd)$ with uniformly bounded multiplicity and diameters of $\supp\chi_k$. There are also $\tilde{\chi}_k\in C_c^\infty(\overline{\Omega})$ that has uniformly bounded multiplicities and diameters of $\supp\tilde\chi_k$, that $\supp\chi_k\cap\supp(1-\tilde\chi_k)=\emptyset$. Fix $k$ and let $\tilde\Omega$ be a bounded uniformly Lipschitz domain such that $\partial\tilde\Omega$ agrees with $\partial\Omega$ on $\supp\tilde\chi_k$. Let $f_\partial=\rr^l\psi(\dd)\chi_ku|_{\partial\Omega}\in \rr^{l}H^s(\partial\tilde\Omega)$ in the sense of non-tangential limits: see further in \cite[Corollary~5.5]{JeKe:95}. Then we aim to estimate
\begin{equation}\label{e:4l30}
    \|\rr^l\psi(\dd)\chi_ku\|_{H^{s+\frac12}}\le C\|\mc{H}_M^{\tilde\Omega}f_\partial\|_{H^{s+\frac12}}+\|\mc{R}_M^{\tilde\Omega}\rr^l\psi(\dd)\chi_k(-\Delta+M^2)u+\mc{R}_M^{\tilde\Omega}[-\Delta, \rr^l\psi(\dd)\chi_k]u\|_{H^{s+\frac12}},
\end{equation}
the first term of which can be bounded by Lemma~\ref{l:reg2}, 
\begin{equation}
    \|\mc{H}_M^{\tilde\Omega}f_\partial\|_{H^{s+\frac12}}\le C\|f_\partial\|_{H^{s}_\partial}=C\|\rr^l\chi_ku \|_{H^{s}_\partial}.
\end{equation}
For $0<s<1$, we utilize $\mc{R}_M^{\tilde\Omega}: H^{s-\frac32}\rightarrow H^{s+\frac12}$ bounded from Lemma~\ref{l:reg2} to see
\begin{equation}\label{e:4l33}
    \|\rr^l\psi(\dd)\chi_ku\|_{H^{s+\frac12}}\le C\| \rr^l\psi(\dd) \chi_k(-\Delta+M^2) u\|_{H^{s-\frac32}}+C\|\rr^l \chi_ku\|_{H^s_\partial}+C\|\rr^l\tilde\psi(\dd)\tilde\chi_k u\|_{H^{s-\frac12}}.
\end{equation}
One can choose $\chi_k$ that the diameters and Lipschitz constants of $\tilde\Omega$ are uniformly bounded for each $k$: this makes all constants $C$ uniformly in $k$. Now since $\chi_k$ have uniformly bounded multiplicity, we have
\begin{equation}
    \|\rr^l\psi(\dd) u\|_{H^{s+\frac12}}\leq C\|\rr^l\psi(\dd) (-\Delta+M^2)u\|_{H^{s-\frac32}}+C\|\rr^l u\|_{H^s_\partial}+C\|\rr^l\tilde{\psi}(\dd) u\|_{H^{s-\frac12}}.
\end{equation}
Applying it first to $s=0+0$ and then to any $s<1$ justifies \eqref{e:4l29}. 

2. When $s=0$, 
by Lemma~\ref{l:reg2} the second term of \eqref{e:4l30} can be bounded by
\begin{gather}
    \|\mc{R}_M^{\tilde\Omega}\rr^l\psi(\dd)\chi_k(-\Delta+M^2)u\|_{H^{\frac12}}\le C\|\rr^l\psi(\dd)\chi_k(-\Delta+M^2)u\|_{H^{-\frac32+0}},\\
    \|\mc{R}_M^{\tilde\Omega}[-\Delta, \rr^l\psi(\dd)\chi_k]u\|_{H^{\frac12}}\le C\|[-\Delta, \rr^l\psi(\dd)\chi_k]u\|_{H^{-\frac32+0}}\le C\|\rr^l\tilde\psi(\dd)\tilde\chi_k u\|_{H^{-\frac12+0}}.
\end{gather}
Similar arguments gluing up estimates in $k$ yield the case. When $s=1$, we have $\mc{R}_M^{\tilde\Omega}: L^2\rightarrow H^{\frac32}$ from Lemma~\ref{l:reg2} instead
\begin{gather}
    \|\mc{R}_M^{\tilde\Omega}\rr^l\psi(\dd)\chi_k(-\Delta+M^2)u\|_{H^{\frac32}}\le C\|\rr^l\psi(\dd)\chi_k(-\Delta+M^2)u\|_{L^2},\\
    \|\mc{R}_M^{\tilde\Omega}[-\Delta, \rr^l\psi(\dd)\chi_k]u\|_{H^{\frac32}}\le C\|[-\Delta, \rr^l\psi(\dd)\chi_k]u\|_{L^2}\le C\|\rr^l\tilde\psi(\dd)\tilde\chi_k u\|_{H^{1}}.
\end{gather}
One can apply \eqref{e:4l29} with $s=\frac12$ to the second line to conclude the proof for the Dirichlet case.

3. We address the Neumann cases where $0<s<1$. We further ask that $\partial_\nu \chi_k\equiv 0$. For each $k$, the uniformly bounded diameter of $\supp\tilde\chi_k$ implies that $C^{-1}K(k)\le\rr\le CK(k)$ on $\supp\tilde\chi_k$ for some $C>0$ independent of $k$. Step 1 in Lemma~\ref{l:reg3} then implies
\begin{multline}
    \|\psi(\dd)\chi_k u\|_{H^{s+\frac12}}\le C\|(-\Delta+M^2)\psi(\dd)\chi_k u\|_{H^{s-\frac32}}+C\|\chi_k \partial_\nu u\|_{H^{s-1}_{\partial}}+C\|\psi(\dd)\chi_k u\|_{H^{-\frac12+0}}\\
    \le C\|\psi(\dd)\chi_k (-\Delta+M^2) u\|_{H^{s-\frac32}}+C\|\chi_k \partial_\nu u\|_{H^{s-1}_{\partial}}+C\|\tilde\psi(\dd)\tilde\chi_k u\|_{H^{s-\frac12}},
\end{multline}
by looking at $[-\Delta, \psi(\dd)\chi_k]$. Now noting that $\rr$ is comparable to a constant on $\supp\tilde\chi_k$, we have
\begin{multline}
    \|\rr^l\psi(\dd)\chi_k u\|_{H^{s+\frac12}}\le CK^l\|\psi(\dd)\chi_k u\|_{H^{s+\frac12}}\le CK^l\|\psi(\dd)\chi_k (-\Delta+M^2) u\|_{H^{s-\frac32}}+CK^l\|\chi_k \partial_\nu u\|_{H^{s-1}_{\partial}}\\
    +CK^l\|\tilde\psi(\dd)\tilde\chi_k u\|_{H^{s-\frac12}}\le C\|\rr^l\psi(\dd)\chi_k (-\Delta+M^2) u\|_{H^{s-\frac32}}+C\|\rr^l\chi_k \partial_\nu u\|_{H^{s-1}_{\partial}}    +C\|\rr^l\tilde\psi(\dd)\tilde\chi_k u\|_{H^{s-\frac12}},
\end{multline}
for $C>0$ independent of $k$. This is the counterpart to \eqref{e:4l33} and the rest follows similarly. The cases $s=0$ and $s=1$ are also similar. 
\end{proof}

As a consequence of Lemma~\ref{l:nearBoundary}, we can control the finer regularity gain of $\mc{R}_M^\Omega$. 
\begin{lemma}
\label{l:basicReg}
Let $\Omega\subset\mathbb{R}^d$ be a uniformly Lipschitz domain with $|m|<\rardyr_{\Omega}$, $|m-l|<\hardyr_{U_{\gamma}^{r_0},\partial\Omega}$. Let $0<s<1$. Let $M_0>0$. There is $C>0$ such that
\begin{equation}\label{e:4l8}
    \|\rr^{l}\dd^{m-1-l}\mc{R}_M^\Omega f\|_{H^{s+\frac12}}\le C\|\rr^l \dd^{m+1-l}f\|_{H^{s-\frac32}},
\end{equation}
uniformly for all $f\in \rr^{-l}\dd^{-m-1+l}H^{s-\frac32}(\Omega)$ and $M\in[0,M_0]$. 
\end{lemma}
\begin{proof}
1. Consider first the case $s\ge\frac12$. We begin with $f\in \mc{S}(\Omega)$ consisting of the restriction to $\Omega$ of elements in $\mc{S}(\mathbb{R}^d)$. Since $\mc{R}_M^\Omega: Y_m^l\rightarrow \dot X_m^l$ bounded,
Lemma~\ref{l:nearBoundary} implies
\begin{multline}\label{e:4l7}
    \|\rr^l\psi(\dd) \mc{R}_M^\Omega f\|_{H^{s+\frac12}}\leq C\|\rr^l\tilde\psi(\dd) f\|_{H^{s-\frac32}}+C\|\rr^l\tilde{\psi}(\dd) \mc{R}_M^\Omega f\|_{H^{-\frac12+0}}\\
    \le C\|\rr^l\tilde\psi(\dd) f\|_{H^{s-\frac32}}+C\|f\|_{Y_m^l}\le C\|\rr^l \dd^{m+1-l}f\|_{H^{s-\frac32}},
\end{multline}
where we used $\mc{R}_M^\Omega f\in \rr^\infty\dot H^{\frac12+0}$, $\rr^{-l}\dd^{-m-1+l} H^{-1}\subset Y_m^l$ from Lemma \ref{l:X-SpacesProp}, $s\ge \frac12$. We also have
\begin{equation}   
\|\rr^l(1-\psi(\dd))\dd^{m-1-l}\mc{R}_M^\Omega f\|_{H^{1}}\le \|(1-\psi(\dd))\mc{R}_M^\Omega f\|_{X_m^l}\le C\|f\|_{Y_m^l}\le C\|\rr^l\dd^{m+1-l} f\|_{H^{-1}}.
\end{equation}
Noting
\begin{equation}
    (1-\psi(\dd))\nabla\mc{R}_M^\Omega f=\mc{R}_M^\Omega (1-\psi(\dd))\nabla f+\mc{R}_M^\Omega[\Delta, \psi(\dd)]\nabla \mc{R}_M^\Omega f
\end{equation}
and $[\Delta, \psi(\dd)]\nabla: X_m^l\rightarrow Y_m^l$ bounded, we have
\begin{multline}   
\|(1-\psi(\dd))\nabla\mc{R}_M^\Omega f\|_{X_m^l}\leq \|(1-\psi(\dd)) \nabla f\|_{Y_m^l}+\|f\|_{Y_m^l}\\
\leq \|\rr^l(1-\psi(\dd))\dd^{m+1-l} \nabla f\|_{H^{-1}}+\|\rr^{l}\dd^{m+1-l}f\|_{H^{-1}}\le C\|\rr^l\dd^{m+1-l} f\|_{H^{0}}.
\end{multline}
We end up with
\begin{equation}
    \|\rr^l(1-\psi(\dd))\dd^{m-1-l}\mc{R}_M^\Omega f\|_{H^2}\le C\|\rr^l\dd^{m+1-l} f\|_{H^{0}}.
\end{equation}
Note $H^\cdot(\Omega)$ is an interpolation scale: see for example \cite[Lemma~4.2]{Chandler-Wilde2015InterpolationCounterexamples}. Interpolation now gives
\begin{equation}
    \|(1-\psi(\dd))\rr^l\dd^{m-1-l}\mc{R}_M^\Omega f\|_{H^{s+\frac12}}\le C\|\rr^l\dd^{m+1-l} f\|_{H^{s-\frac32}},
\end{equation}
and with \eqref{e:4l7},
\begin{equation}
    \|\rr^{l}\dd^{m-1-l}\mc{R}_M^\Omega f\|_{H^{s+\frac12}}\le C\|\rr^l \dd^{m+1-l}f\|_{H^{s-\frac32}},
\end{equation}
uniformly for each $f\in\mc{S}'(\Omega)$. Now by the density of $\mc{S}(\Omega)$ in $\mc{S}'(\Omega)$, we have \eqref{e:4l8} as desired for $\frac12\le s<1$.

2. The case $0<s<\frac12$ remains. For any $\phi\in \mc{S}$, since $\mc{R}_M^\Omega$ is self-adjoint on $L^2$ we have
\begin{equation}
    \langle \rr^{l}\dd^{m-1-l}\mc{R}_M^\Omega f, \phi\rangle=\langle f, \mc{R}_M^\Omega\rr^{l}\dd^{m-1-l}\phi\rangle=\langle \rr^l\dd^{m+1-l} f, \rr^{-l}\dd^{l-m-1}\mc{R}_M^\Omega\rr^{l}\dd^{m-1-l}\phi\rangle,
\end{equation}
where we can use Step 1 to estimate
\begin{equation}
    \|\rr^{-l}\dd^{l-m-1}\mc{R}_M^\Omega\rr^{l}\dd^{m-1-l}\phi\|_{H^{\frac32-s}}\le C\|\phi\|_{H^{-s-\frac12}}
\end{equation}
by noting that $\frac32-s>1$. This implies 
\begin{equation}
    \|\rr^{l}\dd^{m-1-l}\mc{R}_M^\Omega f\|_{H^{s+\frac12}}\le C\|\rr^l\dd^{m+1-l} f\|_{H^{s-\frac32}}
\end{equation}
as desired. 
\end{proof}

We can now improve the mapping properties of $\mc{H}^\Omega_M$.
\begin{lemma}[Elliptic regularity near the Lipschitz boundary]\label{l:ellregunbdd}
Let $\Omega\subset \mathbb{R}^d$ be a uniformly Lipschitz domain. Let $0\le s\le 1$. Let $M_0>0$. Let $|l|<\rardyr_{\Omega}+\hardyr_{U_{\gamma}^{r_0},\partial\Omega}$. Let $\chi \in C_c^\infty(\overline{\Omega})$. Let $\tilde\chi\in C^\infty_c(\overline{\Omega})$ that $\supp\chi\cap\supp(1-\tilde\chi)=\emptyset$. Then there is $C>0$ such that 
\begin{equation}\label{e:4l25}
    \|\chi \mc{H}_M^\Omega f\|_{H^{s+\frac12}(\Omega)}\leq C\|\tilde\chi f\|_{H^{s}(\partial\Omega)}+C\|f\|_{H^{0,l}(\partial\Omega)},
\end{equation}
uniformly for all $f\in H^{s,l}(\partial\Omega)$ and $M\in[0,M_0]$.  If $\partial_\nu \chi\equiv 0$, we further have for all $f\in H^{s,l}(\partial\Omega)$, 
\begin{equation}\label{e:4l26}
    \|\chi \mc{H}_M^\Omega f\|_{H^{s+\frac12}(\Omega)}\leq C\|\tilde\chi 
    \Lambda_M^\Omega f\|_{H^{s-1}(\partial\Omega)}+C\|f\|_{H^{0,l}(\partial\Omega)}.
\end{equation}
\end{lemma}
\begin{proof}
1. We first work with \eqref{e:4l25}. Let $\psi\in C_c^\infty(\mathbb{R})$ that $\supp\chi\cap\supp(1-\psi(\dd))=\emptyset$. For $s\in [0,1)$, Lemma~\ref{l:nearBoundary} 
implies
\begin{multline}\label{e:4l27}
\|\psi(\dd)\chi \mc{H}^\Omega_M f\|_{H^{s+\frac12}}\le C\|\chi f\|_{H^{s}_\partial}+C\|[-\Delta ,\chi]\mathcal{H}_M^\Omega f\|_{H^{s-\frac32+0}}+C\|\tilde\psi(\dd) \chi \mathcal{H}_M^\Omega f\|_{H^{-\frac12+0}}
\\
\leq C\|\chi f\|_{H^{s}_\partial}+C\|\tilde\psi(\dd)\tilde\chi\mathcal{H}_M^\Omega f\|_{H^{s-\frac12+0}}+C\|\tilde\psi(\dd) \chi \mathcal{H}_M^\Omega f\|_{H^{-\frac12+0}},
\end{multline}
repeated use of which gives
\begin{equation}\label{e:4l31}
    \|\psi(\dd)\chi \mc{H}^\Omega_M f\|_{H^{s+\frac12}}\leq C\|\tilde\chi f\|_{H^{s}_\partial}+C\|\tilde\psi(\dd) \tilde\chi \mathcal{H}_M^\Omega f\|_{H^{-\frac12+0}}.
\end{equation}
We aim to bound the last term here. We remind the reader of the extension $E_0:H^{0,l}_\partial\to H^{\frac12,l}$ with $\supp E_0 f\subset \{\dd\le 2\}$ constructed in Lemma~\ref{l:babyExtend}, by $E_0=\sum_k \tilde\chi_k\mc{H}_0^{\tilde\Omega_k}\chi_k$. Note we have
\begin{equation}
    \mc{H}^\Omega_M f= E_0 f-\mc{R}^\Omega_M(-\Delta+M^2)E_0 f=E_0 f-M^2\mc{R}^\Omega_ME_0 f-\sum_k \mc{R}^\Omega_M[-\Delta, \tilde\chi_k]\mc{H}^{\tilde\Omega_k}_0\chi_k f
\end{equation}
from the correction formula \eqref{e:hForm}. Then $\mc{R}_M^\Omega: Y_m^l\rightarrow \dot X_m^l$ bounded and compact $\supp\tilde\chi$ implies 
\begin{equation}
    \|\tilde\psi(\dd)\tilde\chi\sum_k\mc{R}^\Omega_M[-\Delta, \tilde\chi_k]\mc{H}^{\tilde\Omega_k}_0\chi_k f\|_{H^{-\frac12+0}}\le C \|\rr^l\dd^{-N} \sum_k[-\Delta, \tilde\chi_k]\mc{H}^{\tilde\Omega_k}_0\chi_k f\|_{H^{-1}}\le C\|f\|_{H^{0,l}_\partial}.
\end{equation}
Bounding other two terms similarly, we have
\begin{equation}
    \|\tilde\psi(\dd)\tilde\chi\mc{H}^\Omega_M f\|_{H^{-\frac12+0}}\le C\|f\|_{H_{\partial}^{0,l}},
\end{equation}
bringing which into \eqref{e:4l31} concludes the case $s\in [0,1)$. For the case $s=1$, \eqref{e:4l27} now reads
\begin{equation}
\|\psi(\dd)\chi \mc{H}^\Omega_M f\|_{H^{\frac32}}\leq C\|\chi f\|_{H^{1}_\partial}+C\|\tilde\psi(\dd)\tilde\chi\mathcal{H}_M^\Omega f\|_{H^1}+C\|\tilde\psi(\dd) \chi \mathcal{H}_M^\Omega f\|_{H^{-\frac12+0}}\le C\|\tilde\chi f\|_{H^{1}_\partial}+C\|f\|_{H_{\partial}^{0,l}},
\end{equation}
by applying the case $s=\frac12$ to get rid of the $H^1$-term and renaming the cutoffs if necessary. This shows \eqref{e:4l25}. 

2. We show \eqref{e:4l26}. Lemma~\ref{l:nearBoundary} gives the alternative to \eqref{e:4l27}, in the case of $0\le s<1$, 
\begin{equation}
    \|\chi \mc{H}^\Omega_M f\|_{H^{s+\frac12}}\le C\|\chi \Lambda_M^\Omega f\|_{H^{s-1}_\partial}+C\|[-\Delta ,\chi]\mathcal{H}_M^\Omega f\|_{H^{s-\frac32+0}}+C\|\tilde\psi(\dd) \chi \mathcal{H}_M^\Omega f\|_{H^{-\frac12+0}},
\end{equation}
and the one with $C\|[-\Delta ,\chi]\mathcal{H}_M^\Omega f\|_{L^2}$ in place of $C\|[-\Delta ,\chi]\mathcal{H}_M^\Omega f\|_{H^{s-\frac32+0}}$ for $s=1$. Repeat Step 1 and 2 to see
\begin{equation}\label{e:4l28}
    \|\chi \mc{H}^\Omega_M f\|_{H^{s+\frac12}}\le C\|\tilde\chi\Lambda_M^\Omega f\|_{H^{s-1}_\partial}+C\|\tilde\psi(\dd) \tilde\chi \mathcal{H}_M^\Omega f\|_{H^{-\frac12+0}}, 
\end{equation}
which gives \eqref{e:4l26}. 
\end{proof}

\begin{remark}
All constants in this section can depend on $\Omega$ in the following way: the upper bound on uniform Lipschitz characters of $\partial\Omega$, the upper bound on the choice of $r_0\ge 0$, the lower bound on the choice of $\gamma>0$, and the lower bound on Hardy constants $\rardyr_{\Omega}, \hardyr_{U_{\gamma}^{r_0},\partial\Omega}$. We show in Lemma~\ref{l:fixhardy} that for all \acl domains and all $\e>0$ there are $\gamma,r_0$ such that $\hardyr_{U_{\gamma}^{r_0},\partial\Omega}>\frac{1}{2}-\e$. Moreover, if the only Lipschitz singularity of $\partial\Omega$ is a conic singularity with smooth cross-sections that occurs at $0$, then for all $\e>0$ there is $\gamma>0$ such that for \emph{all} $r_0\geq 0$ we have $\hardyr_{U_{\gamma}^{r_0},\partial\Omega}>\frac{1}{2}-\e$. In particular, this means that in the rest of the paper, we may replace $\hardyr_{U_{\gamma}^{r_0},\partial\Omega}$ by $\frac{1}{2}$ and disregard the dependency on $\gamma$, without any loss of generality or control on constants, without further comment.
\end{remark}

\section{Structure of the \dtn map on \acl domains}\label{s:3}
As a preparation for our study of the \dtn map on \acl domains, we first find a microlocal approximation of it on \acl domains. 

\noindent{\bf{Outline of the strategy}}

It is easy to see that on $\mathbb{R}_+^d$, the DtN map at $M=0$ is given by the Fourier multiplier $\Lambda_0^{\mathbb{R}^d_+}=|D|$ and therefore is a bounded operator from $H^{s,l}(\mathbb{R}^{d-1})$ to $H^{s-1,l}(\mathbb{R}^{d-1})$ only for $|l|<\frac{d-1}{2}+1$. Moreover, it is easy to check that $|D|$ is far from local. For example, when $d=2$ and $\psi_\pm \in C^\infty( \mathbb{R})$ with $\supp \psi_{\pm} \subset \{ \pm x>0\}$ and $\supp (1-\psi_{\pm})\subset \{\pm x< 1\}$, we have
\begin{equation}
    \psi_-|D|\psi_+:H^{-N,l}\to H^{N,l+1}
\end{equation}
for any $N$ and 
 $-\frac{3}{2}<l<l+1<\frac{3}{2}$, but not better. Both of these remarks demonstrate that $|D|$ is far from pseudodifferential. This failure to be pseudodifferential is caused by the singularity of $|\xi|$ at $\xi=0$ and, in fact, one can show that away from $\xi=0$, $\Lambda_0^{\mathbb{R}_+}$ is pseudodifferential.

 Motivated by the above discussion, the goal of this section is to find a (second microlocal) pseudodifferential approximation to the DtN map on an \acl domain that is accurate as close to $\xi=0$ as possible.

In Fermi normal coordinates $\bar{x}=(x_0,x)$, the Laplacian takes the form 
\begin{equation}
    -\Delta=D_{x_0}^2+b(\bar{x})D_{x_0}+R(\bar{x},D_x).
\end{equation}

One well-known method for finding an approximation to the DtN map with mass M is to factor
\begin{equation}
    D_{x_0}^2+b(\bar{x})D_{x_0}+R(\bar{x},D_x)+M^2= (D_{x_0}+iP+b(\bar{x}))(D_{x_0}-iP)+\text{small},
\end{equation}
where $P$ is pseudodifferential with symbol $\sqrt{R(\bar{x}, \xi)+M^2}$ (see e.g.~\cite{lee1989determining}). This is straightforward when $M>0$. However, when $M=0$, one immediately arrives at the problem that $R(\bar{x},D_x)$ is \emph{not} elliptic in the scattering calculus. In particular, its symbol vanishes at $\{\xi=0\}\cap \{ |x|=\infty\}$. As a result, $\sqrt{R(\bar{x}, \xi)}$ is \emph{not} a scattering symbol and the factorization fails. 

To remedy this, we allow an error in the factorization that has a fixed, polynomial decay in $x$  at $\xi=0$ but rapidly becomes residual in both $x$ and $\xi$ away from $\xi=0$ (see Lemma~\ref{l:factor}). One way to think of this is that we add an artificial damping to the Laplacian. Thus, instead of factoring the Laplacian, we factor
\begin{equation}
    D_{x_0}^2+b(\bar{x})D_{x_0}+R(\bar{x},D_{x}) + \Op( \langle x\rangle^{-2\delta}\langle \langle x\rangle^{\delta}|\xi|\rangle^{-N}),
\end{equation}
for an arbitrarily large $N$. 

This factorization involves a pseudodifferential calculus that is second microlocalized at $\{\xi=0\}\cap \{|x|=\infty\}$ (see~\eqref{e:SymbolDef} for the symbols in this class and Appendix~\ref{a:1} for the construction of the calculus). 

Once we have factored $-\Delta$ appropriately, in Lemma~\ref{l:dtnParametrixExtension}, we use energy estimates to construct the corresponding approximate Poisson extension operator, $E_P$, satisfying $(D_{x_0}-iP)E_P=0$. 

In section~\ref{s:dtnStructure}, we finally show that $E_P$ is a good approximation of $\mc{H}_M^{\Omega}$ (and hence also $P$ is a good approximation of $\Lambda_M^{\Omega}$) using microlocally weighted estimates on $\mc{R}_M^{\Omega}$ in $X_{m}^l$ spaces. Building on the theory developed in Section~\ref{s:harm}, these estimates are proved in Lemma~\ref{l:resgammaBounds}.

\subsection{Pseudodifferential calculus and symbol class}
Given an \acs domain $\Omega\subset\mathbb{R}^d$, its boundary $\partial \Omega$ is a scattering manifold of dimension $(d-1)$, embedded in $\mathbb{R}^d$. Denote the half-space collar of $\partial\Omega$ by $\Omega_+=\{(x_0,x): x_0\in [0,\infty), x\in\partial\Omega\}$, where $x$ is the local coordinates on $\partial\Omega$. We will use the second-microlocal pseudodifferential calculus $\Psi^{s,l}_{k,\delta}(\Omega_+)$ with $M_0>0$ fixed, developed in the Appendix~\ref{a:1}. To make it easier to read, we briefly describe the calculus in this subsection. 

Consider the phase space $(\bar x,\xi)\in \Omega_+\times \mathbb{R}^{d-1}$, with position variable $\bar x=(x_0,x)\in\Omega_+$ and momentum variable $\xi\in \mathbb{R}^{d-1}$ dual to $x$. One heuristically sees $x_0$ as the normal variable to $\partial\Omega_x$ and the phase space is $\mathbb{R}_{x_0}\times T^*\partial\Omega_{x}$ equivalently. 

Let $s,l,k\in\mathbb{R}$ and $\delta\in[0,1)$. Fix $M_0>0$. We consider the symbol class $S^{s,l}_{k,\delta}(\Omega_{+}; \mathbb{R}^{d-1}_\xi)$ that contains all $a(x_0,x,\xi;h)\in C^{\infty}(\Omega_+\times \mathbb{R}^{d-1}_\xi\times(0,h_0))$ that for some $C_{\alpha\beta}>0$,
\begin{equation} \label{e:SymbolDef}
\abs{(\langle x\rangle\partial_{\bar x})^\alpha (\langle x\rangle^{-\delta} \langle \langle x\rangle^{\delta}\sqrt{|\xi|^2+h^2M^2}\rangle\partial_\xi)^\beta a(\bar x,\xi;h)}\le C_{\alpha\beta}\langle \xi\rangle^{s}\langle x\rangle^{l}\langle \langle x\rangle^{\delta}\sqrt{|\xi|^2+h^2M^2}\rangle^{k},
\end{equation}
uniform in all $(\bar x,\xi;h)\in \Omega_+\times \mathbb{R}^{d-1}_\xi\times(0, h_0)$ and $M\in [0, M_0]$. This is essentially a tangential calculus in $x$ with global control in $x_0$, and a second microlocalization in $\langle x\rangle^\delta\xi$. Here we interpret $s$ as the regularity in $x$, $l$ as the decay in $x$, and $k$ as the weighted regularity in $x$. The residual symbol class is denoted $S_{-\infty,\delta}^{-\infty,-\infty}(\Omega_+; \mathbb{R}^{d-1}_\xi)$. 

We have an associated quantization map $\Op_h: S^{s,l}_{k,\delta}(\Omega_+; \mathbb{R}^{d-1}_\xi)\rightarrow \Psi^{s,l}_{k,\delta}(\Omega_+)$ locally given by
\begin{equation}
    \Op_h(a)(x_0,x,y_0,y)=\frac{\delta(y_0-x_0)}{(2\pi h)^{d-1}}\int_{\mathbb{R}^{d-1}_{\xi}}e^{\frac{i}{h}\langle x-y,\xi\rangle}a(x_0,x,\xi)\,d\xi, 
\end{equation}
and a principal symbol map $\sigma_h: \Psi^{s,l}_{k,\delta}\rightarrow S^{s,l}_{k,\delta}/hS^{s,l-(1-\delta)}_{k-1,\delta}$. The residual operators are $h^\infty\Psi_{-\infty,\delta}^{-\infty,-\infty}(\Omega_+):\langle x\rangle^N L^2_{x_0}H^{-N}_{x}\rightarrow h^N\langle x\rangle^{-N}L^2_{x_0}H^{N}_{x}$ for any $N\ge 0$. For every $A\in\Psi^{s,l}_{k,\delta}(\Omega_+)$, there exists $a\in S^{s,l}_{k,\delta}(\Omega_+; \mathbb{R}^{d-1}_\xi)$ that $A=\Op_h(a)+O(h^\infty\Psi^{-\infty,-\infty}_{-\infty,\delta})$, unique modulo $S^{-\infty,-\infty}$. 

We say $\chi(x_0, x)\in S^{l}(\Omega_+)$ if $\langle x\rangle^{-l}\chi(x_0,x)$ extends to a function in $ C^\infty([0,\infty)_{x_0}\times\overline{\partial\Omega_x})$, where $\overline{\partial\Omega_x}$ is the radial compactification of $\partial\Omega_x$ using the boundary defining function $\langle x\rangle^{-1}$. This is equivalent to say there are $C_\alpha>0$ that
\begin{equation}
    |\partial_{\bar x}^{\alpha}\chi(\bar x)|\le C_\alpha\langle x\rangle^{l-|\alpha|}
\end{equation}
uniformly in $\bar x\in \Omega_+$ and $\chi(x_0,x)$ is polyhomogeneous in $x$. 

We introduce a convenient notation
\begin{equation}
    \opsi=O(\rre^{-N}\Psi_{-\infty,\delta}^{-\infty,-\infty}(\Omega_+)): \langle \bar x\rangle^N L^2_{x_0}H^{-N}_{x}\rightarrow \langle \bar x\rangle^{-N}L^2_{x_0}H^{N}_{x}. 
\end{equation}
For $\chi(\bar x)\in S^0(\Omega_+)$ such that $\supp\chi\subset\{x_0\le C\langle x\rangle\}$ for some $C>0$, we have an improved remainder estimate:
\begin{equation}
    A\in \Psi^{-\infty,-\infty}_{-\infty,\delta}\Rightarrow A\chi=\opsi\chi,
\end{equation}
by that $\rre\sim \langle x\rangle$ on the support of $\chi$.

The operators in $\Psi^{s,l}_{k,\delta}(\Omega_+)$ are by nature tangential to $\partial\Omega_x=\{x_0=0\}$, in the sense that they are bounded on $\mc{S}'(\partial\Omega_x)$. We denote their restriction on $\partial\Omega_x$ by $\Psi^{s,l}_{k,\delta}(\partial\Omega)$. 

\begin{remark}
    A heuristic way to understand why we use such a calculus, is to see $\Lambda^{\mathbb{R}^2_+}_M=\sqrt{D^2+M^2}$ converging to $\Lambda^{\mathbb{R}^2_+}_0=|D|$ as Fourier multipliers. The limit $|D|$ has a $H^{\frac32-0}$-singularity at $\xi=0$ and thus maps $H^{s,l}$ to $H^{s-1,l}$ for $|l|\le \frac32-0$, $s\in\mathbb{R}$. This decay threshold only happens at $\xi=0$, and our calculus captures that key feature. 
\end{remark}

We denote by $\Gamma^{s,l}_{k,\delta}\in \Psi^{s,l}_{k,\delta}(\Omega_+)$ the scaling operators built in Lemma~\ref{l:GammaInverse} in Appendix~\ref{a:1} and corresponding second-microlocal weighted Sobolev spaces $H^{s,l}_{k,\delta}(\Omega_+)=\Gamma^{-s,-l}_{-k,\delta}L^2(\Omega_+)$, and $H^{s,l}_{k,\delta}(\partial\Omega)=\Gamma^{-s,-l}_{-k,\delta}L^2(\partial\Omega)$. We brief the readers that $\Gamma^{s,l}_{k,\delta}$ has classical principal symbol $\langle h\xi\rangle^{s}\langle x\rangle^{l}\langle h\langle x\rangle^{\delta}(|\xi|^2+M^2)^{\frac12}\rangle^{k}$. Note that for $k\ge 0$, we have
\begin{equation}
H^{s+k,l+\delta k}_{0, \delta}\subset H^{s,l}_{k, \delta}. 
\end{equation}

\subsection{The boundary factorization}\label{s:3-2}
In this subsection, we factor the Laplacian near $\partial\Omega$ and derive its tangential part, $P\in\Psi^{0,-\delta}_{1,\delta}$, as a microlocal representation of $\Lambda_M^\Omega$. We also identify an extension operator $E_{P}$ such that $(D_{x_0}-iP)E_{P}=0$. 

\begin{figure}[t]
  \centering
  \begin{subfigure}[t]{0.5\textwidth}
        \centering
        \begin{tikzpicture}[x=1cm,y=1cm]
  \definecolor{colorAKANE}{RGB}{183,40,46}
  \definecolor{colorUSUAI}{RGB}{40,148,200}
  \colorlet{plotPrimary}{colorAKANE}
  \colorlet{plotTertiary}{colorUSUAI}
  \colorlet{plotAxis}{black!72}
  \tikzset{
    >=latex, font=\small, mark size=1.5pt, mark options={thin},
    line cap=round, line join=round,
    collar fill/.style={fill=plotTertiary!16,draw=none},
    collar edge/.style={draw=plotAxis,line width=0.75pt},
    boundary curve/.style={draw=plotPrimary,line width=0.6pt},
    coordinate axis/.style={draw=plotAxis,line width=0.6pt,->},
    figure label/.style={text=black!84,inner sep=1.3mm},
    region label/.style={figure label,draw=none,fill=none}
  }

  \path[use as bounding box] (0,0) rectangle (7.5,3.2);
  \fill[fill=white] (0,0) rectangle (7.5,3.2);
  \begin{scope}[shift={(3.5,0.8)},scale=1.6]
    \def\collarOpening{0.35}
    \coordinate (C) at (2,{1+\collarOpening});
    \coordinate (D) at (-2,{1+\collarOpening});
    \coordinate (Center) at (0,0.5);
    \coordinate (lowerRight) at (2,0.5);
    \coordinate (lowerCenter) at (0,0);
    \coordinate (lowerLeft) at (-2,0.5);
    \coordinate (lowerControlLeft) at ({-2/3},{-1/6});
    \coordinate (lowerControlRight) at ({2/3},{-1/6});
    \path (lowerControlLeft) ++(0,{0.5-\collarOpening/3})
      coordinate (upperControlLeft);
    \path (lowerControlRight) ++(0,{0.5-\collarOpening/3})
      coordinate (upperControlRight);
    \coordinate (axisTop) at (0,1);
    \coordinate (regionLabel) at (-0.5,0.28);
    \coordinate (boundaryLabel) at (-1.7,0.3);
    \coordinate (normalLabel) at (0.3,1);
    \coordinate (tangentialLabel) at (2.1,0.65);

    \begin{scope}
      \path[collar fill]
        (lowerRight) -- (C)
        .. controls (upperControlRight) and (upperControlLeft) .. (D)
        -- (lowerLeft)
        .. controls (lowerControlLeft) and (lowerControlRight) .. (lowerRight)
        -- cycle;
      \draw[collar edge]
        (C) .. controls (upperControlRight) and (upperControlLeft) .. (D);
    \end{scope}

    \begin{scope}
      \draw[coordinate axis] (lowerCenter) -- (axisTop);
      \draw[boundary curve,->]
        (lowerLeft)
        .. controls (lowerControlLeft) and (lowerControlRight) .. (lowerRight);
    \end{scope}

    \begin{scope}
      \node[figure label,text=plotPrimary,anchor=north]
        at (boundaryLabel) {$\partial \Omega$};
      \node[region label] at (regionLabel) {$U$};
      \node[figure label] at (normalLabel) {$x_0$};
      \node[figure label] at (tangentialLabel) {$x$};
    \end{scope}
  \end{scope}
\end{tikzpicture}%
        \caption{$\Omega$}
  \end{subfigure}\hfill
  \begin{subfigure}[t]{0.5\textwidth}
        \centering
        \begin{tikzpicture}[x=1cm,y=1cm]
  \definecolor{colorAKANE}{RGB}{183,40,46}
  \definecolor{colorUSUAI}{RGB}{40,148,200}
  \colorlet{plotPrimary}{colorAKANE}
  \colorlet{plotTertiary}{colorUSUAI}
  \colorlet{plotAxis}{black!72}
  \tikzset{
    >=latex, font=\small, mark size=1.5pt, mark options={thin},
    line cap=round, line join=round,
    collar fill/.style={fill=plotTertiary!16,draw=none},
    collar edge/.style={draw=plotAxis,line width=0.75pt},
    coordinate axis/.style={draw=plotAxis,line width=0.6pt,->},
    figure label/.style={text=black!84,inner sep=1.3mm},
    region label/.style={figure label,draw=none,fill=none}
  }

  \path[use as bounding box] (0,0) rectangle (7.5,3.2);
  \fill[fill=white] (0,0) rectangle (7.5,3.2);
  \begin{scope}[shift={(3.5,0.8)},scale=1.6]
    \coordinate (A) at (-2,0);
    \coordinate (B) at (2,0);
    \coordinate (C) at (2,1);
    \coordinate (D) at (-2,1);
    \coordinate (Center) at (0,0.5);
    \coordinate (axisOrigin) at (0,0);
    \coordinate (axisTop) at (0,1);
    \coordinate (regionLabel) at (-0.5,0.28);
    \coordinate (boundaryLabel) at (-1.7,-0.2);
    \coordinate (normalLabel) at (0.3,1);
    \coordinate (tangentialLabel) at (2.1,0.15);

    \begin{scope}
      \path[collar fill]
        (A) -- (B) -- (C)
        to[out=-150,in=0] (Center)
        to[out=180,in=-30] (D) -- cycle;
      \draw[collar edge]
        (C) to[out=-150,in=0] (Center)
        to[out=180,in=-30] (D);
    \end{scope}

    \begin{scope}
      \draw[coordinate axis,draw=plotPrimary] (A) -- (B);
      \draw[coordinate axis] (axisOrigin) -- (axisTop);
      \node[figure label,text=plotPrimary] at (boundaryLabel) {$\partial \Omega$};
      \node[region label] at (regionLabel) {$U$};
      \node[figure label] at (normalLabel) {$x_0$};
      \node[figure label] at (tangentialLabel) {$x$};
    \end{scope}
  \end{scope}
\end{tikzpicture}%
        \caption{$\Omega_+$}
  \end{subfigure}
  \caption{Illustration for $\Omega$ and $\Omega_+$.}
  \label{f:3f1}
\end{figure}
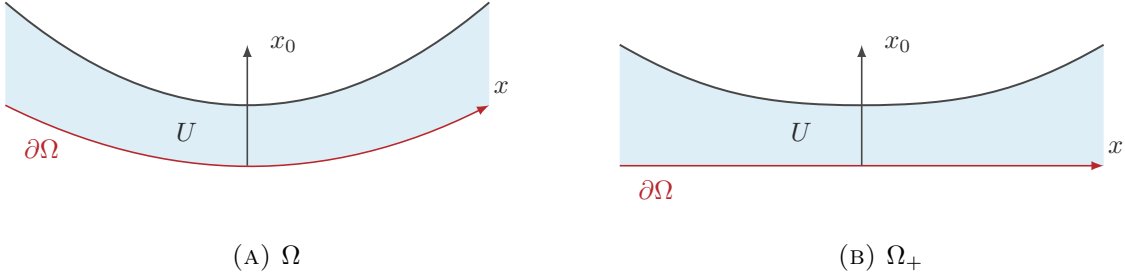

We first identify a conic collar of $\partial\Omega$ in $\Omega$ by that in $\Omega_+$. Fix some $c>0$ small, we choose a finite set of local Fermi normal coordinates $\bar x=(\xbd,x)$ in an asymptotically conic set $U$ identified as $\{(\xbd,x)\in \Omega_+: \xbd <c\langle x\rangle\}$ near $\partial\Omega$, such that $\partial\Omega$ is identified as $\{\xbd=0\}$, and the metric on $\Omega$ in those coordinates, $\gplus$, satisfies $(\gplus)_{00}(x_0, x)\equiv 1$, $(\gplus)_{0j}(x_0, x)\equiv0$ for each $1\le j\le d-1$. Such $U$ is an open subset of $\overline{\Omega}$. See Figure~\ref{f:3f1}.

Here we interpret $\xbd$ as the normal coordinates to the $\partial\Omega$ and $x$ as the local coordinates that parametrize $\partial\Omega$. Our strategy it to perform analysis in $U\subset \Omega_+$ using the pseudodifferential tools, and then transplant the analysis to $U\subset \Omega$. We say $\chi(\bar x)\in S^{l}(U)$ if $\chi(\bar x)\in S^l(\Omega_+)$ and $\supp \chi\subset U$. 

Furthermore $\gplus$ satisfies the symbolic estimate in $\bar x$: 
\begin{equation}\label{e:3l26}
    |\partial_{\bar x}^{\alpha}(\gplus)_{ij}(\bar x)|\le C_\alpha\langle x\rangle^{-|\alpha|}, \ 0\le i,j\le d-1, 
\end{equation}
for all $\bar x\in U$. We remind the reader that $\langle \bar x\rangle\sim \langle x\rangle$ on $U$, and we can extend the metric into $(\gplus)_{ij}\in S^0(\Omega_+)$ that retains $(\gplus)_{00}\equiv 1$ and $(\gplus)_{0j}\equiv 0$, and $\det((\gplus)_{ij})\ge C^{-1}$ uniformly for some $C>0$. In particular, one can extend in the way that there is $b(\bar x)\in S^{-1}(\Omega_+)$ that
\begin{equation}\label{e:boundarydecomp}
    -\Delta_{\gplus}= D_{x_0}^2+b(\bar x)D_{x_0}+R_2(\bar x,D_{x})+R_1(\bar x,D_{x}),
\end{equation}
where $R_k(\bar x,\xi)$ are homogeneous polynomials of degree $k$ in $\xi$ and 
\begin{equation}
    |\partial_{\bar x}^{\alpha}R_k(\bar x,\xi)|\leq C_{\alpha}\langle x\rangle^{-2+k-|\alpha|}|\xi|^k,\ R_2(\bar x,\xi)=\sum_{i,j=1}^{d-1}(\gplus)^{ij}(\bar x)\xi_i\xi_j=|\xi|_{\gplus(\bar x)}^2,
\end{equation}
where $(g_+^{ij})_{i,j=0}^{d-1}=((\gplus)_{ij})^{-1}$. Here we abbreviated $|(0, \xi_1, \cdots, \xi_{d-1})|_{\gplus(x_0,x)}^2=|\xi|_{\gplus(x_0,x)}^2$. At $x_0=0$, we have $|\xi|_{\gplus(0,x)}^2=|\xi|_g^2$. Note that $R_1(\bar x, D_x)\in h^{-1}\Psi^{0, -1-\delta}_{1,\delta}$, $R_2(\bar x, D_x)\in h^{-2}\Psi^{0, -2\delta}_{2, \delta}$. 

Our first goal is to factor $-\Delta_\gplus+M^2$ uniformly in $M\geq 0$ inside $U$, into tangential parts and conormal parts with respect to $\partial\Omega$. However, at $M=0$, $-\Delta_\gplus$ is not uniformly elliptic in the scattering calculus or its second-microlocal variant. Therefore one expects to accept an error in the factorization in order to work in a well-behaved calculus of pseudodifferential operators. 

\begin{lemma}[Factorization of $-\Delta_\gplus+M^2$ on $\Omega_+$]
\label{l:factor}
Let $M_0>0$. Let $\delta\in[0,1)$. Let $h_0>0$ be small enough. There then exists a family of invertible operators $P\in h^{-1}\Psi^{0,-\delta}_{1,\delta}(\Omega_+)$ smoothly dependent on the parameter $h\in(0, h_0]$ and $M\in[0,M_0]$, such that the remainder 
\begin{equation}
    \mc{E}_0^\circ=(-\Delta_\gplus+M^2)-(D_{x_0}+iP+b(x_0,x))(D_{x_0}-iP)\in h^{-2}\Psi^{0,-2\delta}_{-\infty,\delta}(\Omega_+)\subset h^{-N-2}\Psi^{-N,-2\delta}_{-\infty,\delta}(\Omega_+),
\end{equation}
for each $N\ge0$. Moreover we have
\begin{equation}
    \sigma_h(P)= h^{-1}\langle x\rangle^{-\delta}\big\langle \langle  x\rangle^\delta \sqrt{|\xi|^2_{\gplus(\bar x)}+h^2M^2}\big\rangle,\ \sigma(P)=h^{-1}\langle x\rangle^{-\delta}\big\langle h\langle  x\rangle^\delta \sqrt{|\xi|^2_{\gplus(\bar x)}+M^2}\big\rangle.
\end{equation}
The construction of $P$ only depends locally on $g_+$. 
\end{lemma}
\begin{proof}
1. Set
\begin{equation}
    P_0=h^{-1}\Op_h\big(\langle x\rangle^{-\delta}\langle\langle x\rangle^\delta \sqrt{|\xi|^2_{\gplus}+h^2M^2}\rangle\big)\in h^{-1}\Psi^{0,-\delta}_{1,\delta}, \ a_0=\sigma_h(P_0)\in h^{-1}S^{0,-\delta}_{1,\delta},
\end{equation}
where $P_0$ is elliptic. We have from \eqref{e:boundarydecomp}
\begin{equation*}
(D_{x_0}+iP_0+b)(D_{x_0}-iP_0)=D_{x_0}^2+bD_{x_0}-ibP_0+i[P_0,D_{x_0}]+P_0^2=-\Delta_\gplus+M^2+\mc{E}_1^\circ,
\end{equation*}
where
\begin{equation}
    \mc{E}_1^\circ=(-R_2(x_0,x, D_{x})-M^2+P_0^2)-ibP_0+i[P_0, D_{x_0}]-R_1(x_0, x, D_{x})\in h^{-2}\Psi^{0,-2\delta}_{1,\delta},
\end{equation}
since 
\begin{equation}
    [P_0, D_{x_0}]\in h^{-1}\Psi^{0,-1-\delta}_{1,\delta}, \ bP_0\in h^{-1}\Psi^{0,-1-\delta}_{1,\delta}, \ -R_2(x_0,x, D_{x})-M^2+P_0^2\in h^{-2}\Psi^{0,-2\delta}_{1,\delta}.
\end{equation}

2. Iteratively, assume that we have $a_j$ for $j=0,\dots N-1$, $a_j\in h^{-1}S^{0,-\delta}_{1-j,\delta}$, such that with $P_{N-1}=\sum_{j=0}^{N-1}\Op_h(a_j)$, we have
\begin{equation}
    (D_{x_0}+iP_{N-1}+b)(D_{x_0}-iP_{N-1})=-\Delta_\gplus+M^2+\mc{E}_N^\circ,
\end{equation}
with $\mc{E}_N^\circ\in h^{-2}\Psi^{0,-2\delta}_{2-N, \delta}$. Then, 
\begin{multline}
(D_{x_0}+iP_{N-1}+i\Op_h(a_N)+b)(D_{x_0}-iP_{N-1}-i\Op_h(a_N))
=-\Delta_\gplus+M^2+\mc{E}_N^\circ\\-i(D_{x_0}+iP_{N-1}+b)\Op_h(a_N)+i\Op_h(a_N)(D_{x_0}-iP_{N-1}-i\Op_h(a_N))
=-\Delta_\gplus+M^2+\mc{E}_N^\circ\\+2P_{N-1}\Op_h(a_N)+[\Op_h(a_N),P_{N-1}]-i[D_{x_0},\Op_h(a_N)]-ib\Op_h(a_N)+\Op_h(a_N)^2.
\end{multline}
Setting 
\begin{equation}
    a_N=-\frac{\sigma_h(\mc{E}_N^\circ)}{2a_0}\in h^{-1}S^{0,-\delta}_{1-N,\delta},
\end{equation}
we have 
\begin{equation}
    \mc{E}_N^\circ+2P_{N-1}\Op_h(a_N), [\Op_h(a_N),P_{N-1}], [D_{x_0},\Op_h(a_N)]\in h^{-2}\Psi^{0,-1-\delta}_{1-N,\delta}, \ \Op_h(a_N)^2\in h^{-2}\Psi^{0,-2\delta}_{2-2N,\delta}.
\end{equation}
Then, with $P_N=P_{N-1}+\Op_h(a_N)$, we have
\begin{equation}
    (D_{x_0}+iP_{N}+b)(D_{x_0}-iP_{N})=-\Delta_\gplus+M^2+\mc{E}_{N+1}^\circ,
\end{equation} 
with $\mc{E}_{N+1}^\circ\in h^{-2}\Psi^{0,-2\delta}_{1-N,\delta}$, and the proof is complete by induction. We end up with the asymptotic sum
\begin{equation}
    P=\sum_{j=0}^\infty \Op_h(a_j), \ \mc{E}_0^\circ=((-\Delta_\gplus+M^2) -(D_{x_0}+iP+b)(D_{x_0}-iP)), 
\end{equation}
as desired.  
\end{proof}

Now that we have factored $-\Delta_\gplus+M^2$, we obtain energy estimates for $(D_{x_0}+iP+b)$. Here we denote $L^2_{\bar x}=L^2(\Omega_+)$, and $L^\infty_{x_0}L^2_x=L^\infty([0,\infty)_{x_0};L^2(\partial\Omega_x))$. 
\begin{lemma}[Energy estimate for $D_{x_0}-iP$]
\label{l:energy}
Let $\delta\in [0,1)$. Let $N>0$. Let $\e>0$ small. Let $J\ge 0$. There exists $h_0>0$ and $C>0$ such that for all $h\in(0,h_0]$, for all $|s|,|l|,|k|\le N$, we have
\begin{multline}\label{e:3l2}
\|\langle x_0\rangle^{-\frac{1}{2}-\e}\Gamma^{s+\frac{1}{2},l}_{k,\delta} (x_0 P)^J v\|_{L^2_{\bar x}}^2+\|\Gamma^{s,l}_{k,\delta}(x_0 P)^J v\|^2_{L^\infty_{x_0}L^2_{x}}+ h^{-1}\|\Gamma^{s,l-\frac\delta2}_{k+\frac12,\delta}(x_0 P)^J v\|_{L^2_{\bar x}}^2\\
\leq C\sum_{j=0}^J( \| \langle x_0\rangle^{1/2+\e}\Gamma^{s-\frac12,l}_{k,\delta}(x_0 P)^j (D_{x_0}-iP)v\|^2_{L^2_{\bar x}}+\|\Gamma^{s,l}_{k,\delta}v(0)\|_{L^2_{x}}^2),
\end{multline}
and
\begin{multline}\label{e:3l3}
\|\langle x_0\rangle^{-\frac{1}{2}-\e}\Gamma^{s+\frac{1}{2},l}_{k,\delta}v\|_{L^2_{\bar x}}^2+\|\Gamma^{s,l}_{k,\delta}v\|^2_{L^\infty_{x_0}L^2_{x}}+h^{-1}\|\Gamma^{s,l-\frac\delta2}_{k+\frac12,\delta}v\|_{L^2_{\bar x}}^2\\
 \leq C( \| \langle x_0\rangle^{\frac12+\e}\Gamma^{s-\frac12,l}_{k,\delta}(D_{x_0}+iP)v\|^2_{L^2}+\liminf_{x_0\rightarrow \infty}\|\Gamma^{s,l}_{k,\delta}v(x_0)\|_{L^2_{x}}^2),
\end{multline}
uniformly for all $v\in\mc{S}'(\Omega_{+})$. Here we denoted $x\mapsto v(x_0,x)$ by $v(x_0)$.
\end{lemma}
\begin{proof}
0. Through the proof, all we assumed is $P\in h^{-1}\Psi^{0,-\delta}_{1,\delta}$ being elliptic, that is for some $c>0$, 
\begin{equation}
    \Re \sigma_h(P)(\bar x, \xi)>ch^{-1}\langle x\rangle^{-\delta}\langle \langle x\rangle^\delta\sqrt{|\xi|^2+h^2M^2}\rangle>0, \ (\bar x,\xi)\in \Omega_+\times \mathbb{R}^{d-1}.
\end{equation}
Under this assumption, for $h>0$ small, $P$ is invertible with its inverse $P^{-1}\in h\Psi^{0,\delta}_{-1,\delta}$. All results here still hold for any elliptic operators other than $P$ in $h^{-1}\Psi^{0,-\delta}_{1,\delta}$. 

1. We first show \eqref{e:3l2} in the case of $(s,l,k)=(0,0,0)$, $J=0$. Consider
\begin{multline}\label{e:3l1}
\partial_{x_0}\|v(x_0)\|_{L^2_{x}}^2= -2\Im\langle (D_{x_0}-iP)v(x_0),v(x_0)\rangle_{L^2_{x}}  -2\Re \langle P v(x_0),v(x_0)\rangle_{L^2_{x}}\\
\leq 2|\Im\langle (D_{x_0}-iP)v(x_0),v(x_0)\rangle_{L^2_{x}}| -ch^{-1}\|\Gamma^{0, -\frac\delta2}_{\frac12, \delta}v(x_0)\|_{L^2_{x}}^2
\end{multline}
Now integrate \eqref{e:3l1} in $x_0$ from $0$ to $T$ to see
\begin{multline}\label{e:3l35}
    \|v(T)\|_{L^2_{x}}^2+ch^{-1}\int_0^T\|\Gamma^{0, -\frac\delta2}_{\frac12,\delta} v(x_0)\|_{L^2_{x}}^2\,dx_0
    \le K\int_0^T\|\langle x_0\rangle^{\frac12+\e}\langle D_{x}\rangle^{-\frac12}(D_{x_0}-iP)v(x_0)\|_{L^2_{x}}^2\,dx_0\\
    + \frac{1}{4K}\int_0^T \|\langle x_0\rangle^{-\frac12-\e}\langle D_{x}\rangle^{\frac12}v(x_0)\|_{L^2_{x}}^2\,dx_0+\|v(0)\|_{L^2_{x}}^2.
\end{multline}
By noting $\Id+\Gamma^{0,-\frac\delta2}_{\frac12,\delta}$ is elliptic over the microsupport of $\langle D_x\rangle^\frac12\in \Psi^{0,0}_{\frac12, \delta}$, We can estimate 
\begin{multline}\label{e:3l36}
     \frac{1}{4K}\int_0^T  
    \|\langle x_0\rangle^{-\frac12-\e}\langle D_{x}\rangle^{\frac12}v\|_{L^2_{x}}^2\,dx_0\le \frac{C}{4K}\|\langle x_0\rangle^{-\frac12-\e}v\|_{L^2_{\bar x}}^2+\frac{C}{4K}\|\langle x_0\rangle^{-\frac12-\e}\Gamma^{0, -\frac\delta2}_{\frac12,\delta} v\|_{L^2_{\bar x}}^2\\
    \le \frac12 \|v\|_{L^\infty_{x_0}L^2_{x}}^2+\frac{c}{2} h^{-1}\|\Gamma^{0, -\frac\delta2}_{\frac12,\delta} v\|_{L^2_{\bar x}}^2,
\end{multline}
by choosing some $K$ large. Then taking the supremum in $T\in [0,\infty)$ in \eqref{e:3l35} we have
\begin{equation}\label{e:3l4}
    \frac{1}{2}\|v\|_{L^\infty_{x_0}L^2_{x}}^2+\frac12 ch^{-1}\|\Gamma^{0,-\frac\delta2}_{\frac12,\delta} v\|^2_{L^2_{\bar x}}\le C\|\langle x_0\rangle^{\frac12+\e}\langle D_{x}\rangle^{-\frac12}(D_{x_0}-iP)v\|_{L^2_{\bar x}}^2+\|v(0)\|_{L^2_{x}}^2. 
\end{equation}
Note \eqref{e:3l36} allows us to bound $\|\dde^{-\frac12-\e}\Gamma^{\frac12,0}_{0,\delta} v\|_{L^2_{\bar x}}^2$. This is \eqref{e:3l2} when $(s,l,k)=(0,0,0)$, $J=0$. 

2. We now consider $J\ge 1$. Consider
    \begin{multline}
        \partial_{x_0} \|(x_0 P) v(x_0)\|_{L^2_x}^2\le -ch^{-1}\|\Gamma^{0,-\frac\delta2}_{\frac12,\delta} (x_0 P) v(x_0)\|_{L^2_x}^2+2\langle P v(x_0), (x_0 P)v(x_0)\rangle_{L^2_{x}}\\
        +2|\Im\langle (x_0 P)(D_{x_0}-i P)v(x_0), (x_0 P) v(x_0)\rangle_{L^2_{x}}|
        +2|\Re\langle x_0[\partial_{x_0},P]v(x_0), (x_0 P) v(x_0)\rangle_{L^2_{x}}|.
    \end{multline}
    We estimate term by term: 
    \begin{equation}
        2\langle P v(x_0), (x_0 P)v(x_0)\rangle_{L^2_{x}}\le Ch^{-1}\|\Gamma^{0, -\frac\delta2}_{\frac12,\delta}v(x_0)\|_{L^2_x}^2+\frac{c}{4h}\|\Gamma^{0, -\frac\delta2}_{\frac12,\delta} (x_0 P)v(x_0)\|_{L^2_x}^2,
    \end{equation}
    and
    \begin{equation}
        2|\Re\langle x_0[\partial_{x_0},P]v(x_0), (x_0 P) v(x_0)\rangle_{L^2_x}|\le C\|\Gamma^{0,-\frac12}_{0,\delta}(x_0 P) v(x_0)\|_{L^2_x}^2
    \end{equation}
    by that $[\partial_{x_0}, P]\in h^{-1}\Psi^{0, -\delta-1}_{1, \delta}$ 
    and $P\in h^{-1}\Psi^{0,-\delta}_{1,\delta}$ is elliptic, and $\Psi^{0,0}_{0,\delta}$ is bounded on $L^2_x$.
    So
    \begin{multline}
        \partial_{x_0} \|(x_0 P) v(x_0)\|_{L^2_x}^2\le -\frac{c}2h^{-1}\|\Gamma^{0,-\frac\delta2}_{\frac12,\delta} (x_0 P) v(x_0)\|_{L^2_x}^2+Ch^{-1}\|\Gamma^{0, -\frac\delta2}_{\frac12,\delta}v(x_0)\|_{L^2_x}^2\\+K\|\langle x_0\rangle^{\frac12+\e}\langle D_{x}\rangle^{-\frac12}(x_0P)(D_{x_0}-iP)v(x_0)\|_{L^2_{x}}^2+\frac{1}{4K}\|\langle x_0\rangle^{-\frac12-\e}\langle D_{x}\rangle^{\frac12}(x_0P)v(x_0)\|_{L^2_{x}}^2
    \end{multline}
    After integrating and taking supremum in $x_0\in[0,\infty)$, noting $(x_0 P)v(x_0)=0$ at $x_0=0$, we have
    \begin{multline}
        \|(x_0 P) v\|_{L^\infty_{x_0}L^2_x}^2+\frac{c}2h^{-1}\|\Gamma^{0,-\frac\delta2}_{\frac12,\delta} (x_0 P) v\|_{L^2_{\bar x}}^2
        \le Ch^{-1}\|\Gamma^{0, -\frac\delta2}_{\frac12,\delta}v\|_{L^2_{\bar x}}^2\\
        + K\|\langle x_0\rangle^{\frac12+\e}\langle D_{x}\rangle^{-\frac12}(x_0P)(D_{x_0}-iP)v\|_{L^2_{\bar x}}^2+\frac{1}{4K}\|\langle x_0\rangle^{-\frac12-\e}\langle D_{x}\rangle^{\frac12}(x_0P)v\|_{L^2_{\bar x}}^2.
    \end{multline}
    Now similar to \eqref{e:3l36} we have for large $K>0$, 
    \begin{equation}
        \frac{1}{4K}\int_0^T  
    \|\langle x_0\rangle^{-\frac12-\e}\langle D_{x}\rangle^{\frac12}(x_0P)v\|_{L^2_{x}}^2\,dx_0\le \frac12 \|(x_0P)v\|_{L^\infty_{x_0}L^2_{x}}^2+\frac{c}{4} h^{-1}\|\Gamma^{0, -\frac\delta2}_{\frac12,\delta} (x_0P)v\|_{L^2_{\bar x}}^2,
    \end{equation}
    and using \eqref{e:3l4} to estimate $Ch^{-1}\|\Gamma^{0, -\frac\delta2}_{\frac12,\delta}v\|_{L^2_{\bar x}}^2$, we arrive at
    \begin{multline}
        \|(x_0 P) v\|_{L^\infty_{x_0}L^2_x}^2+\frac{c}4h^{-1}\|\Gamma^{0,-\frac\delta2}_{\frac12,\delta} (x_0 P) v\|_{L^2_{\bar x}}^2\le C\sum_{j=0}^1\|\langle x_0\rangle^{\frac12+\e}\langle D_{x}\rangle^{-\frac12}(x_0P)^j(D_{x_0}-iP)v\|_{L^2_{\bar x}}^2\\
        +C\|v(0)\|_{L^2_x}^2.
    \end{multline}
    Iterate this procedure to observe
    \begin{multline}
        \frac12 \|(x_0 P)^J v\|_{L^\infty_{x_0}L^2_x}^2+\frac{c}{2^{J+1}}h^{-1}\|\Gamma^{0,-\frac\delta2}_{\frac12,\delta} (x_0 P)^J v\|_{L^2_{\bar x}}^2\le C\sum_{j=0}^J\|\langle x_0\rangle^{\frac12+\e}\langle D_{x}\rangle^{-\frac12}(x_0P)^j(D_{x_0}-iP)v\|_{L^2_{\bar x}}^2\\
        +C\|v(0)\|_{L^2_x}^2,
    \end{multline}
    this is \eqref{e:3l2} when $(s,l,k)=(0,0,0)$.

3. We now show \eqref{e:3l3} in the case $(s,l,k)=(0,0,0)$. Replacing $P$ by $-P$ in \eqref{e:3l1} gives
\begin{equation}
    -\partial_{x_0}\|v(x_0)\|_{L^2_{x}}^2
\leq 2|\Im\langle (D_{x_0}+iP)v(x_0),v(x_0)\rangle_{L^2_{x}}| -ch^{-1}\|\Gamma^{0, -\frac\delta2}_{\frac12, \delta}v(x_0)\|_{L^2_{x}}^2.
\end{equation}
Let $R$ be a large number to be set. Similarly the integration in $x_0$ from $R$ back to some $T\in [0, R]$ and taking the supremum of $T\in [0, R]$ gives
\begin{multline}
\sup_{T\in[0,R]}\|v(T)\|_{L^2_{x}}^2+ch^{-1}\int_0^R\|\Gamma^{0, -\frac\delta2}_{\frac12,\delta} v(x_0)\|_{L^2_{x}}^2\,dx_0
    \le C\int_0^R\|\langle x_0\rangle^{1/2+\e}\langle D_{x}\rangle^{-\frac12}(D_{x_0}+ iP)v(x_0)\|_{L^2_{x}}^2\,dx_0\\
    +\frac12 \sup_{T\in [0, R]}\|v(T)\|_{L^2_{x}}^2+\frac{c}{2} h^{-1}\int_{0}^R\|\Gamma^{0, -\frac\delta2}_{\frac12,\delta} v(x_0)\|_{L^2_x}^2\,dx_0+\|v(R)\|_{L^2_{x}}^2,
\end{multline}
where $C$ is chosen independent of $R$. Now taking a sequence of $R\rightarrow \infty$ where $\|v(R)\|_{L^2_{x}}\rightarrow \liminf_{x_0\rightarrow \infty}\|v(x_0)\|_{L^2_{x}}$ gives \eqref{e:3l3} similarly when $(s,l,k)=(0,0,0)$, $J=0$. The procedure for the case $J\ge 1$ is similar to Step 2. 

4. For the general case $(s,l,k)$, one notes 
\begin{equation}
   [\Gamma^{s,l}_{k,\delta}, D_{x_0}]\Gamma^{-s,-l}_{-k,\delta}\in \Psi^{0,-1}_{0,\delta},\ [\Gamma^{s,l}_{k,\delta}, P]\Gamma^{-s,-l}_{-k,\delta}\in \Psi^{0,-1}_{0,\delta}
\end{equation}
are lower order terms compared to $P\in h^{-1}\Psi^{0,-\delta}_{1,\delta}$ elliptic, and thus
\begin{equation}
    \Re\sigma_h(P\pm i[\Gamma^{s,l}_{k,\delta}, D_{x_0}]\Gamma^{-s,-l}_{-k,\delta}+[\Gamma^{s,l}_{k,\delta}, P]\Gamma^{-s,-l}_{-k,\delta})=\Re\sigma_h(P)
\end{equation}
also satisfies the ellipticity assumption specified in Step 0. Noting that
\begin{equation}
     \Gamma^{s,l}_{k,\delta}(D_{x_0}-iP)\Gamma^{-s,-l}_{-k,\delta}=D_{x_0}-i(P+ i[\Gamma^{s,l}_{k,\delta}, D_{x_0}]\Gamma^{-s,-l}_{-k,\delta}+[\Gamma^{s,l}_{k,\delta}, P]\Gamma^{-s,-l}_{-k,\delta}),
\end{equation}
we arrive at \eqref{e:3l4} with $P$ replaced with $P+ i[\Gamma^{s,l}_{k,\delta}, D_{x_0}]\Gamma^{-s,-l}_{-k,\delta}+[\Gamma^{s,l}_{k,\delta}, P]\Gamma^{-s,-l}_{-k,\delta}$:
\begin{equation}
    \|v\|_{L^\infty_{x_0}L^2_{x}}^2+\frac12 ch^{-1}\|\Gamma^{0,-\frac\delta2}_{\frac12,\delta} v\|^2_{L^2_{\bar x}}\le C\|\langle x_0\rangle^{\frac12+\e}\langle D_{x}\rangle^{-\frac12}\Gamma^{s,l}_{k,\delta}(D_{x_0}-iP)\Gamma^{-s,-l}_{-k,\delta}v\|_{L^2_{\bar x}}^2+\|v(0)\|_{L^2_{x}}^2. 
\end{equation}
Replace $v$ by $\Gamma^{s,l}_{k,\delta} v$ to observe \eqref{e:3l2}. The proof of cases $J\ge 1$ and \eqref{e:3l3} are similar. 
\end{proof}

Next, we construct an extension operator $E_{P}:\mc{S}'(\partial\Omega)\rightarrow \mc{S}'(\Omega_+)$ that solves $(D_{x_0}-iP)E_{P}=0$. 
\begin{lemma}[Poisson extension of $D_{x_0}-iP$]
\label{l:dtnParametrixExtension}
Let $\delta\in [0,1)$. Let $N>0$. Let $\e>0$ small. 
There exists an extension operator $E_{P}:\mc{S}'(\partial\Omega)\rightarrow \mc{S}'(\Omega_+)$ such that for any $u_0\in\mc{S}'(\partial\Omega)$, 
\begin{equation}
    (D_{x_0}-iP)E_{P}u_0=0,\ E_{P}u_0|_{x_0=0}=u_0,\ E_{P}u_0\in C^\infty([0,\infty)_{x_0}; \langle x_0\rangle^{\frac{1}{2}+\e}\mc{S}'(\partial\Omega)).
\end{equation}
For such extension operator $E_{P}$, there exist $h_0>0$ and $C>0$ such that
\begin{equation}\label{e:3l28}
    \|\langle x_0\rangle^{-\frac{1}{2}-\e}\Gamma^{-s+\frac12,-l}_{-k,\delta}E_{P}\Gamma^{s, l}_{k,\delta}\|_{L^2_{x}\to L^2_{\bar x}}+\|\Gamma^{-s,-l}_{-k,\delta}E_{P}\Gamma^{s,l}_{k,\delta}\|_{L^2_x\to L^\infty_{x_0} L^2_x}+\|\Gamma^{-s,-l-\frac{\delta}{2}}_{-k+\frac{1}{2},\delta}E_{P}\Gamma^{s,l}_{k,\delta}\|_{L^2_{x}\to L^2_{\bar x}}<C,
\end{equation}
uniformly for all $h\in(0,h_0]$, for all $|s|,|l|,|k|\le N$.
\end{lemma}
\begin{proof}
1. We follow the approach of \cite[Theorem 23.1.2]{HorV3}. Let $u_0\in \Gamma^{-s,-l}_{-k,\delta}L^2(\partial\Omega)$. We want to show that there exists a 
unique solution $u$ such that $\langle x_0\rangle^{-\frac{1}{2}-\e}\Gamma^{s+\frac{1}{2},l}_{k,\delta} u\in L^2(\Omega_+)$ to 
\begin{equation} \label{e:hyperbPb_-}
(D_{x_0}-iP)u=0,\ u|_{x_0=0}=u_0, \ u=E_{P}u_0.
\end{equation}
The uniqueness of $u$ follows from the energy estimate \eqref{e:3l2}. We show the existence of $u$. Observe that $u$ satisfies \eqref{e:hyperbPb_-}, if and only if for all $T>0$ and all $v\in C_c^\infty(\{0\le x_0\le T\})$, we have
\begin{equation} \label{e:weakeq}
\int_0^\infty\langle u,(D_{x_0}+iP^*)v\rangle_{L^2_{x}}\,dx_0=-i \langle u_0,v(0)\rangle_{L^2_{x}}.
\end{equation}

1a. We first find $u\in \langle x_0\rangle^{\frac{1}{2}+\e}\Gamma_{-k,\delta}^{-s-\frac{1}{2},-l}L^2_{\bar x}$ such that \eqref{e:weakeq} holds. Applying the energy estimate \eqref{e:3l3} to $v$ with $P$ replaced by $P^*$ and noting that $\Re\sigma(P^*)=\Re\sigma(P)$ gives
\begin{equation}
\|\Gamma^{-s,-l}_{-k,\delta}v\|^2_{L^\infty_{x_0}L^2_{x}}\leq C \|\langle x_0\rangle^{\frac12+\e}\Gamma^{-s-\frac12,-l}_{-k,\delta}(D_{x_0}+iP^*)v\|^2_{L^2_{\bar x}} 
\end{equation}
which implies
\begin{equation}
|\langle u_0,v(0)\rangle_{L^2_{x}}| \leq C \| \langle x_0\rangle^{\frac12+\e}\Gamma^{-s-\frac12,-l}_{-k,\delta}(D_{x_0}+iP^*)v\|_{L^2_{\bar x}}\|\Gamma^{s,l}_{k,\delta}u_0\|_{L^2_{x}}. 
\end{equation}
Hence, by the Hahn-Banach theorem, we can extend the functional
\begin{equation}
    (D_{x_0}+iP^*)v\mapsto -i\langle u_0,v(0)\rangle_{L^2_{x}}
\end{equation}
to any $(D_{x_0}+iP^*)v\in \Gamma^{s+\frac12,l}_{k,\delta}\langle x_0\rangle^{-\frac12-\e}L^2_{\bar x}$. By the Riesz Theorem there is a $u\in \langle x_0\rangle^{\frac12+\e} \Gamma^{-s-\frac12,-l}_{-k,\delta}L^2_{\bar x}$ satisfying \eqref{e:weakeq} and thus \eqref{e:hyperbPb_-}, as desired. 

1b. Since $\langle x_0\rangle^{-\frac{1}{2}-\e}\Gamma_{k,\delta}^{s+\frac{1}{2},l}u\in L^2_{\bar x}$, we have by the energy estimate \eqref{e:3l2} that
\begin{equation}
\|\langle x_0\rangle^{-\frac{1}{2}-\e}\Gamma^{s+\frac{1}{2},l}_{k,\delta}u\|_{L^2_{\bar x}}^2+\|\Gamma^{s,l}_{k,\delta}u\|^2_{L^\infty_{x_0}L^2_{x}}+ h^{-1}\|\Gamma^{s,l-\frac{\delta}{2}}_{k+\frac{1}{2},\delta}u\|_{L^2_{\bar x}}^2\leq C\|\Gamma^{s,l}_{k,\delta}u_0\|_{L^2_{x}}^2,
\end{equation}
yielding \eqref{e:3l28} and $u\in L^\infty([0,\infty)_{x_0}; \Gamma^{-s,-l}_{-k,\delta}L^2_x)$. 

2. We show $u$ is $C^\infty$ in $x_0$. Note $P\in h^{-1}\Psi^{0,-\delta}_{1,\delta}$ and thus \eqref{e:hyperbPb_-} implies for almost every $x_0\in[0,\infty)$. 
\begin{equation}\label{e:3l27}
    D_{x_0}\Gamma_{k,\delta}^{s,l}u=i\Gamma_{k,\delta}^{s,l}Pu+[D_{x_0},\Gamma_{k,\delta}^{s,l}]u\in \Gamma^{0, -\delta}_{1,\delta}L^2_{x},
\end{equation}
where the commutator terms are of lower order. Then $u\in H^1_{\operatorname{loc}}([0,\infty)_{x_0};\Gamma_{-k+1,\delta}^{-s,-l-\delta}L^2_x)$ and is $C^\frac12$ in $x_0$. Now note
\begin{equation}
    D_{x_0}^2 u=-P^2 u+i[D_{x_0},P]u, \ [D_{x_0},P]\in h^{-1}\Psi^{0,-1-\delta}_{1,\delta}.
\end{equation}
One can then repeat \eqref{e:3l27} for $D_{x_0}^2$ to conclude $u\in H^2_{\operatorname{loc}}([0,\infty)_{x_0};\Gamma_{-k+2,\delta}^{-s,-l-2\delta}L^2_x)$ embeds into $C^{\frac32}([0,\infty)_{x_0};\Gamma_{-k+2,\delta}^{-s,-l-2\delta}L^2_x)$ continuously. Iteration yields the desired regularity in $x_0$. 
\end{proof}

\begin{remark}\label{r:fixe}
    From this point onward, we fix some $h>0$ such that Lemma~\ref{l:energy} and Lemma~\ref{l:dtnParametrixExtension} hold, and disregard all further dependency on $h$. This also makes $P\in \Psi^{0,-\delta}_{1,\delta}$ invertible with $P^{-1}\in \Psi^{0,\delta}_{-1,\delta}$. The choice of such $h$ only depends on upper bounds $C_\alpha$ of all derivatives of the metric $g_+$ in \eqref{e:3l26}, choice of $\delta$, $M_0$, and the upper bounds $N$ of indices in Lemma~\ref{l:energy} and Lemma~\ref{l:dtnParametrixExtension}. 
\end{remark}

\begin{figure}[t]
  \centering
  \begin{subfigure}[t]{0.5\textwidth}
        \centering
        \begin{tikzpicture}[x=1cm,y=1cm]
  \definecolor{colorAKANE}{RGB}{183,40,46}
  \definecolor{colorUSUAI}{RGB}{40,148,200}
  \colorlet{plotPrimary}{colorAKANE}
  \colorlet{plotTertiary}{colorUSUAI}
  \colorlet{plotAxis}{black!72}
  \tikzset{
    >=latex, font=\small, mark size=1.5pt, mark options={thin},
    line cap=round, line join=round,
    support fill/.style={fill=plotTertiary!16,draw=none},
    collar edge/.style={draw=plotAxis,line width=0.75pt},
    support edge/.style={draw=plotTertiary,line width=0.6pt},
    coordinate axis/.style={draw=plotAxis,line width=0.6pt,->},
    figure label/.style={text=black!84,inner sep=1.3mm},
    support label/.style={figure label,text=plotTertiary,draw=none,fill=none}
  }

  \path[use as bounding box] (0,0) rectangle (7.5,3.2);
  \fill[fill=white] (0,0) rectangle (7.5,3.2);
  \begin{scope}[shift={(3.5,0.8)},scale=1.6]
    \def\supportOpening{0.35}
    \coordinate (A) at (-2,0);
    \coordinate (B) at (2,0);
    \coordinate (C) at (2,{1+\supportOpening});
    \coordinate (D) at (-2,{1+\supportOpening});
    \coordinate (Center) at (0,0.5);
    \coordinate (lowerRight) at (2,0.6);
    \coordinate (lowerCenter) at (0,0.1);
    \coordinate (lowerLeft) at (-2,0.6);
    \coordinate (lowerControlLeft) at ({-2/3},{-1/15});
    \coordinate (lowerControlRight) at ({2/3},{-1/15});
    \path (lowerControlLeft) ++(0,{0.4-\supportOpening/3})
      coordinate (upperControlLeft);
    \path (lowerControlRight) ++(0,{0.4-\supportOpening/3})
      coordinate (upperControlRight);
    \coordinate (axisOrigin) at (0,0);
    \coordinate (axisTop) at (0,1);
    \coordinate (supportLabel) at (-0.5,0.32);
    \coordinate (boundaryLabel) at (-1.7,-0.2);
    \coordinate (normalLabel) at (0.3,1);
    \coordinate (tangentialLabel) at (2.1,0.15);

    \begin{scope}
      \path[support fill]
        (lowerRight) -- (C)
        .. controls (upperControlRight) and (upperControlLeft) .. (D)
        -- (lowerLeft)
        .. controls (lowerControlLeft) and (lowerControlRight) .. (lowerRight)
        -- cycle;
      \draw[collar edge]
        (C) .. controls (upperControlRight) and (upperControlLeft) .. (D);
      \draw[support edge]
        (lowerLeft)
        .. controls (lowerControlLeft) and (lowerControlRight) .. (lowerRight);
    \end{scope}

    \begin{scope}
      \draw[coordinate axis,draw=plotPrimary] (A) -- (B);
      \draw[coordinate axis] (axisOrigin) -- (axisTop);
      \node[support label] at (supportLabel) {$\mathrm{supp}\,\chi$};
      \node[figure label,text=plotPrimary] at (boundaryLabel) {$\partial \Omega$};
      \node[figure label] at (normalLabel) {$x_0$};
      \node[figure label] at (tangentialLabel) {$x$};
    \end{scope}
  \end{scope}
\end{tikzpicture}%
        \caption{Lemma~\ref{l:extensiongainsdecay}}
  \end{subfigure}\hfill
  \begin{subfigure}[t]{0.5\textwidth}
        \centering
        \begin{tikzpicture}[x=1cm,y=1cm]
  \definecolor{colorAKANE}{RGB}{183,40,46}
  \definecolor{colorTOKIWA}{RGB}{0,123,67}
  \definecolor{colorUSUAI}{RGB}{40,148,200}
  \colorlet{plotPrimary}{colorAKANE}
  \colorlet{plotSecondary}{colorTOKIWA}
  \colorlet{plotTertiary}{colorUSUAI}
  \colorlet{plotAxis}{black!72}

  \tikzset{
    >=latex,
    font=\small,
    line cap=butt,
    line join=round,
    wedge edge/.style={draw=plotAxis,line width=0.95pt},
    old arc/.style={draw=plotPrimary,line width=0.95pt},
    marked edge/.style={draw=plotSecondary,line width=1.25pt,line cap=butt},
    marked point/.style={fill=plotSecondary,draw=none},
    collar region/.style={fill=plotTertiary!14,draw=none},
    collar curve/.style={
      draw=plotTertiary,
      line width=1pt,
      line cap=round,
      line join=round
    },
    figure label/.style={text=black!84,inner sep=0.8pt}
  }

  \pgfmathsetmacro{\figureWidth}{7.5}
  \pgfmathsetmacro{\figureHeight}{3.2}
  \pgfmathsetmacro{\frameBottom}{-0.28}
  \pgfmathsetmacro{\frameTop}{\figureHeight+\frameBottom}
  \pgfmathsetmacro{\edgeDx}{4.26}
  \pgfmathsetmacro{\edgeDy}{2.48}
  \pgfmathsetmacro{\tipY}{0.34}
  \pgfmathsetmacro{\edgeLength}{veclen(\edgeDx,\edgeDy)}
  \pgfmathsetmacro{\tipX}{(\figureWidth-\edgeLength)/2}
  \pgfmathsetmacro{\wedgeAngle}{atan2(\edgeDy,\edgeDx)}
  \pgfmathsetmacro{\arcDist}{3.15}
  \pgfmathsetmacro{\arcRadius}{\arcDist*tan(\wedgeAngle/2)}
  \pgfmathsetmacro{\braceInnerClearance}{0.36}
  \pgfmathsetmacro{\braceOuterClearance}{0.60}
  \pgfmathsetmacro{\braceRadius}{\arcRadius-\braceInnerClearance}
  \pgfmathsetmacro{\domainLabelBaseline}{-0.08}

  \path[use as bounding box] (0,{\frameBottom}) rectangle ({\figureWidth},{\frameTop});
  \fill[fill=white] (0,{\frameBottom}) rectangle ({\figureWidth},{\frameTop});

  \coordinate (tip) at ({\tipX},{\tipY});
  \coordinate (top) at ({\tipX+\edgeDx},{\tipY+\edgeDy});
  \coordinate (bottom) at ({\tipX+\edgeLength},{\tipY});
  \path (tip) ++(\wedgeAngle:\arcDist) coordinate (arcTop);
  \path (tip) ++({\arcDist},{\arcDist*tan(\wedgeAngle/2)}) coordinate (arcCenter);
  \path (arcCenter) ++(270:\arcRadius) coordinate (arcBottom);
  \path (arcCenter) ++({90+\wedgeAngle}:\braceRadius) coordinate (braceTop);
  \path (arcCenter) ++(270:\braceRadius) coordinate (braceBottom);
  \path (top) ++({\wedgeAngle-90}:\braceOuterClearance) coordinate (braceTopOuter);
  \path (tip) ++({\edgeLength},{\braceOuterClearance}) coordinate (braceBottomOuter);
  \path (arcBottom) -- coordinate[pos=0.07] (tildeOmegaLabel) (bottom);
  \path (arcBottom) --
    coordinate[pos=0.34] (bottomMarkStart)
    (bottom);
  \path (arcTop) --
    coordinate[pos=0.34] (topMarkStart)
    (top);

  \fill[collar region]
    (arcTop) -- (top) -- (braceTopOuter) -- (braceTop)
    arc[
      start angle={90+\wedgeAngle},
      end angle=270,
      radius=\braceRadius cm
    ]
    -- (braceBottomOuter) -- (bottom) -- (arcBottom)
    arc[
      start angle=270,
      end angle={90+\wedgeAngle},
      radius=\arcRadius cm
    ]
    -- cycle;

  \draw[wedge edge] (tip) -- (top);
  \draw[wedge edge] (tip) -- (bottom);
  \draw[old arc] (arcTop) -- (top);
  \draw[old arc] (arcBottom) -- (bottom);
  \draw[old arc]
    (arcTop) arc[
      start angle={90+\wedgeAngle},
      end angle=270,
      radius=\arcRadius cm
    ];
  \draw[marked edge] (topMarkStart) -- (top);
  \path[marked point] (topMarkStart) circle[radius=1.45pt];
  \draw[marked edge] (bottomMarkStart) -- (bottom);
  \path[marked point] (bottomMarkStart) circle[radius=1.45pt];
  \draw[collar curve]
    (braceTopOuter) -- (braceTop)
    arc[
      start angle={90+\wedgeAngle},
      end angle=270,
      radius=\braceRadius cm
    ]
    -- (braceBottomOuter);

  \node[figure label,text=plotSecondary,anchor=south east]
    at ([xshift=-0.10cm,yshift=0.08cm]topMarkStart)
    {$\operatorname{supp}\chi$};
  \node[figure label,anchor=base west]
    at ({\tipX+0.16},{\domainLabelBaseline}) {$\Omega$};
  \node[figure label,text=plotPrimary,anchor=base]
    at (tildeOmegaLabel |- 0,{\domainLabelBaseline})
    {$\tilde\Omega$};
  \node[figure label,text=plotTertiary,anchor=west]
    at ([xshift=0.22cm,yshift=-0.30cm]braceBottomOuter) {$U$};
\end{tikzpicture}%
        \caption{Remark~\ref{r:aclcollar}}
  \end{subfigure}
  \caption{Illustration for Lemma~\ref{l:extensiongainsdecay} and Remark~\ref{r:aclcollar}.}
  \label{f:3f2}
\end{figure}
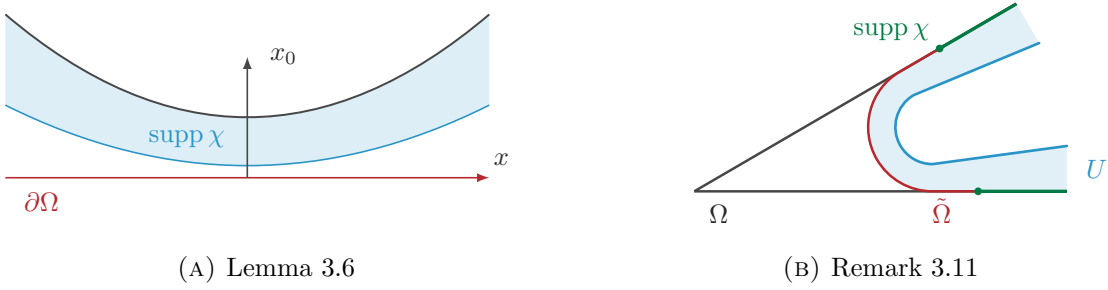

\begin{lemma}[Nonlocality of $E_P$ gives decay and regularity]\label{l:extensiongainsdecay}
Consider $\chi(\bar x), \varphi(\bar x)\in S^0(\Omega_+)$ that $\supp\chi\cap\supp\varphi=\emptyset$. Let $N>0$. Then, there is $C_N>0$ that 
    \begin{equation}
        \|\langle x_0\rangle^{N}\Gamma^{N,N}_{N,\delta} \chi E_{P} \varphi f\|_{L^2(\Omega_+)}\le C_N\|\langle x_0\rangle^{-N}\Gamma^{-N,-N}_{-N,\delta} \varphi f\|_{L^2(\partial\Omega)},
    \end{equation}
    uniformly for $f\in \mc{S}'(\partial\Omega)$. See Figure~\ref{f:3f2}(A).
\end{lemma}
\begin{proof}
It suffices to split the cases into $\supp\chi\subset \{x_0<2c\langle x\rangle\}$ and $\supp\chi\cap \{x_0>c\langle x\rangle\}$ for some $c>0$. Let $\tilde\chi\in S^0(\Omega_+)$ that satisfy all conditions for $\chi$ and $\supp\chi\cap \supp(1-\tilde\chi)=\emptyset$. 

1. Assume $\supp\chi\subset \{x_0<2c\langle x\rangle\}$. Note $\chi E_P \varphi|_{x_0=0}\equiv 0$ and
\begin{equation}
    (D_{x_0}-iP)\chi E_P \varphi f=[D_{x_0}-iP,\chi] \tilde\chi E_P \varphi f+\opsi E_P\varphi f, \ [D_{x_0}-iP,\chi]\in \Psi^{0,-1}_{0,\delta}.
\end{equation}
We remind the reader of the notation $\opsi=O(\rre^{-N}\Psi_{-\infty,\delta}^{-\infty,-\infty}(\Omega_+))$. Lemma~\ref{l:energy} implies
\begin{multline}\label{e:3l43}
\|\Gamma^{s,l}_{k,\delta}\chi E_P \varphi f\|_{L^2_{\bar x}}
\leq C \| \langle x_0\rangle^{1/2+\e}\Gamma^{s-\frac12,l+\frac\delta2}_{k-\frac12,\delta}(D_{x_0}-iP)\chi E_P \varphi f\|_{L^2_{\bar x}}\\\le C \| \langle x_0\rangle^{1/2+\e}\tilde\chi\Gamma^{s-\frac12,l+\frac\delta2}_{k-\frac12,\delta}[D_{x_0}-iP,\chi] \tilde\chi E_P \varphi f\|_{L^2_{\bar x}}+C \| \langle x_0\rangle^{1/2+\e}O(\Psi^{-\infty,-\infty}_{-\infty,\delta})\tilde\chi E_P \varphi f\|_{L^2_{\bar x}}\\
\le C \| \Gamma^{s-\frac12,l+\frac\delta2+\e-\frac12}_{k-\frac12,\delta} \tilde\chi E_P \varphi f\|_{L^2_{\bar x}}+C \| \opsi\tilde\chi E_P \varphi f\|_{L^2_{\bar x}}\le C \| \Gamma^{s-\frac12,l+\frac\delta2+\e-\frac12}_{k-\frac12,\delta} \tilde\chi E_P \varphi f\|_{L^2_{\bar x}},
\end{multline}
where we used $x_0<2c\langle x\rangle$. Now we choose $\e>0$ small that $\frac\delta2+\e-\frac12<0$. Repeating \eqref{e:3l43} by increasing the support of $\tilde\chi$ we have
\begin{equation}
    \|\Gamma^{s,l}_{k,\delta}\chi E_P \varphi f\|_{L^2_{\bar x}}\le C \| \Gamma^{s-\frac12N,l+(\frac\delta2+\e-\frac12)N}_{k-\frac12N,\delta} \tilde\chi E_P \varphi f\|_{L^2_{\bar x}}\le C \| \Gamma^{s-\frac12N,l+\frac\delta2+(\frac\delta2+\e-\frac12)N}_{k-\frac12-\frac12N,\delta} \varphi f\|_{L^2_{x}}
\end{equation}
through Lemma~\ref{l:dtnParametrixExtension}. This concludes the proof for the case $\supp\chi\subset \{x_0<2c\langle x\rangle\}$. 

2. Assume $\supp\chi\subset \{x_0>c\langle x\rangle\}$ and $\supp(1-\tilde\chi)\subset\{x_0<2c\langle x\rangle\}$ without loss of generality. Note $P^{-1}\in \Psi^{0,\delta}_{-1,\delta}$ per Remark~\ref{r:fixe} and
    \begin{multline}
        \Gamma^{s,l}_{k,\delta}\chi E_{P} \varphi=\tilde\chi \Gamma^{s,l}_{k,\delta}\chi P^{-N}P^N E_{P} \varphi+(1-\tilde\chi)O(\Psi^{-\infty,-\infty}_{-\infty,\delta})E_{P} \varphi\\
        =\tilde\chi \langle x\rangle^{N(\frac12+\frac\delta{2})} O(\Psi^{s,l-\frac{N(1-\delta)}2}_{k-N,\delta})P^NE_{P} \varphi
        +\opsi E_{P} \varphi.
    \end{multline}
    With that $x_0>c \langle x\rangle$ on $\supp\tilde\chi$, we have
    \begin{multline}
    \|x_0^{\frac{N(1-\delta)}{2}}\Gamma^{s,l}_{k,\delta} \chi E_{P}\varphi f\|
    \le C\|\tilde\chi\Gamma^{s,l-\frac{N(1-\delta)}2}_{k-N,\delta}(x_0 P)^N E_{P}\varphi f\|+C\|\opsi E_{P}\varphi f\|_{L^2_{x_0}H^{-N}_x}\\
    \le C\|\Gamma^{s, l+\frac\delta2-\frac{N(1-\delta)}2}_{k-\frac12-N,\delta}\varphi f\|_{L^2_x}, 
\end{multline}
where we used Lemma \ref{l:energy} and Lemma \ref{l:dtnParametrixExtension}. Renumber $N$ and note the $L^2$-boundedness of $\Gamma_{-\frac N2,\delta}^{\frac N2, 0}\subset \Psi^{0,0}_{0,\delta}$ to conclude the proof. 
\end{proof}

\subsection{Boundedness properties on spaces with Hardy-type weights}
\label{s:boundedHardyWeight}
Let $\Omega\subset\mathbb{R}^d$ be an \acs domain. We now relate the action of the Poisson operators $\mc{R}_M^\Omega$ defined on $\Omega$ conjugated with operators in $\Gamma^{s,l}_{k,\delta}$ defined on the model space $\Omega_+$, with respect to the spaces $X_{m}^{l}(\Omega_+)$ defined in Definition~\ref{d:Xspaces}. Similarly we consider
$\dot Y_m^l(\Omega_+)=(X_{-m}^{-l}(\Omega_+))^*$ and $Y_m^l(\Omega_+)=(\dot X_{-m}^{-l}(\Omega_+))^*$.

Within this whole section, we abuse notation by using $\rr$, $\dd$ for both weights on $X_m^l(\Omega)$ and the preferred weights $\rr=\langle\bar x\rangle$, $\dd=\langle x_0\rangle$ on $X_m^l(\Omega_+)$. One handy feature about connecting $X_m^l(\Omega_+)$ to $X_m^l(\Omega)$ is that given $\chi(x_0,x)\in S^0(U)$, we have $\rr\sim\langle\bar x\rangle$ and $\dd\sim\langle x_0\rangle$ on $U$, and thus $\chi X_m^l(\Omega_+)\subset X_m^l(\Omega)$ and $\chi X_m^l(\Omega)\subset X_m^l(\Omega_+)$. 

Before stating that result we need to state couple of lemmas that will be relevant in its proof. We first give a decomposition formula for the commutator of the Laplacian with an operator in the $\Psi^{s,l}_{k,\delta}$-calculus. 
\begin{lemma}[Commutator decomposition]\label{l:CommDec}
Let $s,l,k\in\R$. Let $\delta\in[0,1)$. Let $A\in\Psi^{s,l}_{k,\delta}(\Omega_+)$. Then there exist
\begin{equation}
    (A_0)_{ij}\in \Psi^{s,l-1+\delta}_{k-1,\delta}(\Omega_+), \ (A_1)_j\in \Psi^{s,l-1}_{k,\delta}(\Omega_+), \ A_2\in \Psi^{s,l-2}_{k,\delta}(\Omega_+)
\end{equation}
for $0\le i,j\le d-1$ such that
\begin{equation} \label{e:LaplacianCommutator}
[\Delta_\gplus,A]=\operatorname{div}A_0\nabla+\operatorname{div}A_1+A_2.
\end{equation}
\end{lemma}

\begin{proof}
Locally, the commutator reads
\begin{multline}
    [\sqrt{\gplus}^{-1}\partial_{x_i}\sqrt{\gplus}g_+^{ij}\partial_{x_j}, A]=[\sqrt{\gplus}^{-1},A]\partial_{x_i}\sqrt{\gplus}g_+^{ij}\partial_{x_j}+\sqrt{\gplus}^{-1}[\partial_{x_i},A]\sqrt{\gplus}g_+^{ij}\partial_{x_j}\\+\sqrt{\gplus}^{-1}\partial_{x_i}[\sqrt{\gplus}g_+^{ij},A]\partial_{x_j}
    +\sqrt{\gplus}^{-1}\partial_{x_i}\sqrt{\gplus}g_+^{ij}[\partial_{x_j},A]
\end{multline}
We identify the above four terms individually. Note 
\begin{multline}
[\sqrt{\gplus}^{-1},A]\partial_{x_i}\sqrt{\gplus}g_+^{ij}\partial_{x_j}=\partial_{x_i}[\sqrt{\gplus}^{-1},A]\sqrt{\gplus}g_+^{ij}\partial_{x_j}+\partial_{x_j}[[\sqrt{\gplus}^{-1},A],\partial_{x_i}]\sqrt{\gplus}g_+^{ij}\\
+[[[\sqrt{\gplus}^{-1},A],\partial_{x_i}]\sqrt{\gplus}g_+^{ij},\partial_{x_j}]
\in \partial_{x_i}\Psi^{s,l-1+\delta}_{k-1,\delta}\partial_{x_j}+\partial_{x_j}\Psi^{s,l-2+\delta}_{k-1,\delta}+\Psi^{s,l-3+\delta}_{k-1,\delta},
\end{multline}
since $[\partial_{x_i}, \Psi^{s,l}_{k,\delta}]\subset \Psi^{s,l-1}_{k,\delta}$. Similarly, we have
\begin{multline}
    \sqrt{\gplus}^{-1}[\partial_{x_i},A]\sqrt{\gplus}g_+^{ij}\partial_{x_j}=\partial_{x_j}\sqrt{\gplus}^{-1}[\partial_{x_i},A]\sqrt{\gplus}g_+^{ij}+[\sqrt{\gplus}^{-1}[\partial_{x_i},A]\sqrt{\gplus}g_+^{ij},\partial_{x_j}]\\
    \in \partial_{x_j}\Psi^{s, l-1}_{k, \delta}+\Psi^{s, l-2}_{k, \delta},
\end{multline}
\begin{multline}
    \sqrt{\gplus}^{-1}\partial_{x_i}\sqrt{\gplus}g_+^{ij}[\partial_{x_j},A]=\partial_{x_i}\sqrt{\gplus}^{-1}\sqrt{\gplus}g_+^{ij}[\partial_{x_j},A]+[\sqrt{\gplus}^{-1}, \partial_{x_i}]\sqrt{\gplus}g_+^{ij}[\partial_{x_j},A]\\
    \in \partial_{x_i} \Psi^{s,l-1}_{k,\delta}+\Psi^{s, l-2}_{k, \delta},
\end{multline}
and
\begin{multline}
    \sqrt{\gplus}^{-1}\partial_{x_i}[\sqrt{\gplus}g_+^{ij},A]\partial_{x_j}=\partial_{x_i}\sqrt{\gplus}^{-1}[\sqrt{\gplus}g_+^{ij},A]\partial_{x_j}+\partial_{x_j}[\sqrt{\gplus}^{-1},\partial_{x_i}][\sqrt{\gplus}g_+^{ij},A]\\
    +[[\sqrt{\gplus}^{-1},\partial_{x_i}][\sqrt{\gplus}g_+^{ij},A],\partial_{x_j}]
    \in \partial_{x_i}\Psi^{s,l-1+\delta}_{k-1,\delta}\partial_{x_j}+\partial_{x_j}\Psi^{s, l-2+\delta}_{k-1,\delta}+\Psi^{s, l-3+\delta}_{k-1,\delta}. 
\end{multline}
Thus for each pair $i,j$, we have
\begin{equation}
    [\sqrt{\gplus}^{-1}\partial_{x_i}\sqrt{\gplus}g_+^{ij}\partial_{x_j}, A]\in \partial_{x_i} \Psi^{s, l-1+\delta}_{k-1,\delta}\partial_{x_j}+\partial_{x_i}\Psi^{s,l-1}_{k,\delta}+\partial_{x_j}\Psi^{s,l-1}_{k,\delta}+\Psi^{s,l-2}_{k,\delta}.
\end{equation}
Summing up indices $i,j$ concludes the proof. 
\end{proof}

Another lemma is needed to establish the boundedness property of $\Gamma^{s,l}_{k,\delta}$ with respect to the spaces $X_{m}^{l}(\Omega_+)$.
\begin{lemma}[Boundedness of $\Psi^{s,l}_{k,\delta}$ on $X_{m}^{l}(\Omega_+)$]\label{l:xmlbdd}
Let $\delta\in[0,1)$. Let $s,l',k,l,m\in\R$. Let $A\in\Psi^{s,l'}_{k,\delta}$. Let $\chi(x_0,x)\in S^0(U)$. For any $q\ge 0$, there then exists $C>0$ such that
\begin{gather}\label{e:GammaBoundedness_4}
\|\rre^{l}\dde^{m-l+q-1}A\chi u\|_{L^2(\Omega_+)}\le C\|\rre^l\dde^{m-l-1}\Gamma^{s,l'+q}_{k, \delta}\chi u\|_{L^2(\Omega_+)},\\
\label{e:GammaBoundedness_1}
\|\rre^{l}\dde^{m-l+q}A\chi u\|_{L^2(\Omega_+)}\leq C\|\Gamma^{s,l'+1+q}_{k-1,\delta}\chi u\|_{X_m^l(\Omega_+)},\\
\label{e:GammaBoundedness_2}
\|\rre^{l}\dde^{m-l+q}\nabla A\chi u\|_{L^2(\Omega_+)}\leq C\|\Gamma^{s,l'+q}_{k,\delta}\chi u\|_{X_m^l(\Omega_+)},\\\label{e:GammaBoundedness_3}
\|A\chi u\|_{X_{m+q}^l(\Omega_+)}\le C\|\Gamma_{k,\delta}^{s,l'+q}\chi u\|_{X_m^l(\Omega_+)},
\end{gather}
for all $u\in \mc{S}'(\Omega_+)$. Moreover, for $B\in \Psi^{0,0}_{0,\delta}$, we have for any $N>0$ large, 
\begin{gather}\label{e:GammaBoundedness_5}
    \|BA\chi u\|_{X_{m}^l(\Omega_+)}\le C\|A\chi u\|_{X_m^l(\Omega_+)}+C\|\opsi\chi u\|_{X_{-N}^{-N}(\Omega_+)},\\
    \label{e:GammaBoundedness_7}
    \|Bu\|_{X_{N}^N(\Omega_+)}\le C\|u\|_{X_{2N+1}^{2N}(\Omega_+)}, \ 
    \|Bu\|_{Y_{-2N-1}^{-2N}(\Omega_+)}\le C\|u\|_{Y_{-N}^{-N}(\Omega_+)}.
\end{gather}
\end{lemma}
\begin{proof}
0. We begin with some error estimates. We have for any $N>0$, any $s,l',k\in\R$,
\begin{multline}
    \|\opsi u\|\le C\| \rre^{-N}O(\Psi^{-\infty,-\infty}_{-\infty,\delta}) u\|\le C\| \dde^{-N}O(\Psi^{-\infty,-\infty}_{-\infty,\delta}) u\|\\
    \le C\| O(\Psi^{-\infty,-\infty}_{-\infty,\delta})\dde^{-N}\langle x\rangle^{-N}\Gamma^{s,l'}_{k,\delta} u\|\le C\| \rre^{-N}\Gamma^{s,l'}_{k,\delta} u\|\le C\|\Gamma^{s,l'}_{k,\delta} u\|_{X_{-N}^{-N}}. 
\end{multline}
by noting $\rre^2\le \dde^2\langle x\rangle^2$. Similarly we have
\begin{equation}
    \|\nabla\opsi u\|\le \|\opsi\nabla \Gamma^{s,l'}_{k,\delta}u\|+\|\opsi\Gamma^{s,l'}_{k,\delta}u\|\le C\|\Gamma^{s,l'}_{k,\delta} u\|_{X_{-N}^{-N}}. 
\end{equation}
since $[\nabla, \opsi]=\opsi$.

1. Let $\chi_1(\bar x), \chi_2(\bar x)\in S^0(U)$ that $\chi_1\chi=\chi$, $\chi_2\chi_1=\chi_1$. For \eqref{e:GammaBoundedness_4}, note first 
\begin{equation}
    \|\rre^{l}\dde^{m-l+q-1}A\chi u\|\le C\|\chi_1\rre^{l}\dde^{m-l+q-1}A\chi u\|+C\|\rre^{l}\dde^{m-l+q-1}\opsi\chi u\|.
\end{equation}
Noting $\rre\sim \langle x\rangle$ on $U$ and $\chi_2A\Gamma^{-s,-l'}_{-k, \delta}\chi_2\rre^l\dde^{m-l-1}\Gamma^{s,l'+q}_{k, \delta}$ elliptic over the microsupport of $\chi_1\rre^{l+q}\dde^{m-l-1}A$, we have
\begin{multline}
    \|\chi_1\rre^{l}\dde^{m-l+q-1}A\chi u\|\le C\|\chi_2A\Gamma^{-s,-l'}_{-k, \delta}\chi_2\rre^l\dde^{m-l-1}\Gamma^{s,l'+q}_{k, \delta}\chi u\|+C\|\opsi\chi u\|\\
    \le C\|\rre^l\dde^{m-l-1}\Gamma^{s,l'+q}_{k, \delta}\chi u\|+C\|\opsi\chi u\|\le C\|\rre^l\dde^{m-l-1}\Gamma^{s,l'+q}_{k, \delta}\chi u\|,
\end{multline}
yielding \eqref{e:GammaBoundedness_4}. 

2. We claim we have the symbolic decomposition
\begin{equation}\label{e:3l29}
1=B_0+\pmb{B}\cdot \nabla_{x}, \ \pmb{B}=(B_j)_{j=1}^{d-1}, \ B_0\in \Psi^{0,0}_{-1,\delta}, \ B_j\in\Psi^{0,\delta}_{-1,\delta} \text{ for }j\ge 1.
\end{equation}
Indeed, consider for $\psi(\xi)\equiv 1$ near 0 and compactly supported, 
\begin{equation}
\pmb{B}=(\op_h(b_j))_{j=1}^{d-1}, \ b_j=h\abs{\xi}^{-2}\xi_j(1-\psi(\langle x \rangle^\delta\abs{\xi}))\in S^{0,\delta}_{-1,\delta}
\end{equation}
and 
\begin{equation}
\tilde b_0=1-h^{-1}\sum_{j=1}^{d-1}\xi_jb_j\in S^{0,0}_{-\infty,\delta}
\end{equation}
One can check the vanishing of the $S^{0,0}_{0,\delta}$-principal symbol of 
\begin{equation}
F=1-\Op_h(\tilde b_0)+i\operatorname{div}_x\pmb{B}\in \Psi^{0, -1+\delta}_{-1, \delta}.
\end{equation}
Set $B_0=\Op_h(\tilde b_0)+F\in \Psi^{0, 0}_{-1,\delta}$.
We now have the symbolic decomposition 
\begin{equation}
1= B_0-i\operatorname{div}_x\pmb{B},
\end{equation}
whose adjoint gives \eqref{e:3l29} after renaming the operators. 

3. We show \eqref{e:GammaBoundedness_1}. Decompose, with $A_0\in \Psi^{0,-1}_{1,\delta}$ that
\begin{equation}\label{e:3l37}
A=\langle x\rangle^{-q}A_0\Gamma^{s,l'+1+q}_{k-1,\delta}=\langle x\rangle^{-q}A_0(B_0+\pmb{B}\cdot \nabla_{x})\Gamma^{s,l'+1+q}_{k-1,\delta}.
\end{equation}
Let $v=\Gamma^{s,l'+1+q}_{k-1,\delta}\chi u$ and we have 
\begin{equation}
\rre^{l}\dde^{m-l+q}A\chi u=\rre^{l}\dde^{m-l+q}\langle x\rangle^{-q}A_0B_0v+\rre^{l}\dde^{m-l+q}\langle x\rangle^{-q}A_0\pmb{B}\cdot \nabla_{x}v,
\end{equation}
and
\begin{multline}
\rre^{l}\dde^{m-l+q}\langle x\rangle^{-q}A_0B_0\Gamma^{s,l'+1+q}_{k-1,\delta}\chi=\rre^{l}\dde^{m-l+q}\langle x\rangle^{-q}\chi_1 A_0B_0\Gamma^{s,l'+1+q}_{k-1,\delta}\chi+\opsi\chi\\
=A_0B_0\rre^{l}\dde^{m-l+q}\langle x\rangle^{-q}\chi_1\Gamma^{s,l'+1+q}_{k-1,\delta}\chi+[\rre^{l}\langle x\rangle^{-q}\chi_1,A_0B_0]\dde^{m-l+q}\Gamma^{s,l'+1+q}_{k-1,\delta}\chi+\opsi\chi.
\end{multline}
Since $[\rre^{l}\langle x\rangle^{-q}\chi_1, A_0B_0]\in \Psi^{0, l-q-2+\delta}_{-1,\delta}$, we can estimate using $\rre^{-q}\le \dde^{-q}$ that
\begin{multline}
\|\rre^{l}\dde^{m-l+q}\langle x\rangle^{-q}A_0B_0v\|\le\|A_0B_0\rre^{l}\dde^{m-l+q}\langle x\rangle^{-q}\chi_1 v\|+\|[\rre^{l}\langle x\rangle^{-q}\chi_1,A_0B_0]\dde^{m-l+q}v\|\\
+C\|\opsi\chi u\|
\le C\|\rre^{l}\dde^{m-l-1}v\|+C\|\opsi\chi u\|.
\end{multline}
Similarly we have
\begin{multline}
\rre^l\dde^{m-l+q}\langle x\rangle^{-q}A_0\pmb{B}\cdot \nabla_{x}\Gamma^{s,l'+1+q}_{k-1,\delta}\chi =\rre^l\dde^{m-l+q}\langle x\rangle^{-q}\chi_1 A_0\pmb{B}\cdot \nabla_{x}\Gamma^{s,l'+1+q}_{k-1,\delta}\chi+\opsi\chi\\
=A_0\pmb{B}\rre^l\dde^{m-l+q}\langle x\rangle^{-q}\chi_1 \cdot \nabla_{x}\Gamma^{s,l'+1+q}_{k-1,\delta}\chi+[\rre^l\langle x\rangle^{-q}\chi_1, A_0\pmb{B}]\cdot \dde^{m-l+q}\nabla_{x}\Gamma^{s,l'+1+q}_{k-1,\delta}\chi+\opsi\chi,
\end{multline}
and
\begin{multline}
\|\rre^{l}\dde^{m-l+q}\langle x\rangle^{-q}A_0\pmb{B}\cdot \nabla_{x}v\|\le\|A_0\pmb{B}\rre^l\dde^{m-l+q}\langle x\rangle^{-q}\chi_1 \cdot \nabla_{x}v\|\\
+\|[\rre^l\langle x\rangle^{-q}\chi_1 , A_0\pmb{B}]\cdot \dde^{m-l+q}\nabla_{x}v\|+C\|\opsi\chi v\|\le C\|\rre^{l-1+\delta}\dde^{m-l}\nabla_{x}v\|+C\|\opsi\chi u\|,
\end{multline}
since $[\rre^l\langle x\rangle^{-q}\chi_1, A_0\pmb{B}]\in \Psi^{0,l-q-2+2\delta}_{-1,\delta}$. Thus we have 
\begin{multline}
\|\rre^{l}\dde^{m-l+q}A\chi u\|\le C\|\rre^{l}\dde^{m-l-1} v\|+C\|\rre^{l-1+\delta}\dde^{m-l}\nabla_{x}v\|+C\|\opsi\chi u\|\\
\le C\|\Gamma^{s,l'+1+q}_{k-1,\delta}\chi u\|_{X_m^l}+C\|\opsi\chi u\|\le C\|\Gamma^{s,l'+1+q}_{k-1,\delta}\chi u\|_{X_m^l}.
\end{multline}
This is \eqref{e:GammaBoundedness_1}. 

4. For \eqref{e:GammaBoundedness_2}, we employ a decomposition similar to \eqref{e:3l37} to see for some $A_0\in \Psi^{0,0}_{0,\delta}$, 
\begin{equation}
\nabla A \chi=\chi_1 \langle x\rangle^{-q} A_0\nabla\Gamma^{s,l'+q}_{k,\delta}\chi+[\nabla, \chi_1\langle x\rangle^{-q} A_0]\Gamma^{s,l'+q}_{k,\delta}\chi+\nabla\opsi\chi.
\end{equation}
We estimate for $v=\Gamma^{s,l'+q}_{k,\delta}\chi u$, since $[\rre^{l} \langle x\rangle^{-q} \chi_1 ,A_0]\in \Psi^{0, l-1+\delta-q}_{-1, \delta}$,
\begin{multline}
\|\rre^{l}\dde^{m-l+q}\chi_1 \langle x\rangle^{-q} A_0\nabla\Gamma^{s,l'+q}_{k,\delta}\chi u\|\le C\|A_0\rre^{l-q}\dde^{m-l+q}\chi_1 \nabla v\|\\
+C\|[\rre^{l} \langle x\rangle^{-q} \chi_1 ,A_0]\dde^{m-l+q}\nabla v\|
\le C\|\rre^{l}\dde^{m-l} \nabla v\|+C\|\rre^{l-1+\delta}\dde^{m-l}\nabla v\|\\
\le C\|\Gamma^{s,l'+q}_{k,\delta}\chi u\|_{X_m^l}.
\end{multline}
Note that $[\nabla, \chi_1\langle x\rangle^{-q} A_0]\Gamma^{s,l'+q}_{k,\delta}\in \Psi^{s,l'-1}_{k, \delta}$ and applying \eqref{e:GammaBoundedness_4} to it gives
\begin{equation} 
\|\rre^{l}\dde^{m-l+q}  [\nabla, \chi_1\langle x\rangle^{-q} A_0]\Gamma^{s,l'+q}_{k,\delta}\chi u\|\le C\|\Gamma^{s, l'+q}_{k, \delta}\chi u\|_{X_m^l}+C\|\opsi\chi u\|.
\end{equation}
Summing up above two estimates and add $\nabla \opsi\chi u$ to obtain \eqref{e:GammaBoundedness_2} as desired. 

5. So far we have obtained \eqref{e:GammaBoundedness_3}, 
\begin{equation}
\|A\chi u\|_{X_{m+q}^l(\Omega_+)}\le C\|\Gamma_{k,\delta}^{s,l'+q}\chi u\|_{X_m^l(\Omega_+)}+C\|\opsi\chi u\|_{L^2(\Omega_+)}+C\|\opsi\nabla\chi u\|_{L^2(\Omega_+)},
\end{equation}
where the two remainder terms are controlled by $C\|\Gamma_{k,\delta}^{s,l'+q}\chi u\|_{X_m^l(\Omega_+)}$ by Step 1. Note
\begin{equation}
    BA\chi u=B\chi_1A\chi u+B(1-\chi_1)A\chi u, \ B(1-\chi_1)A\chi=O(\Psi^{-\infty,-\infty}_{-\infty,\delta})\chi=\opsi\chi
\end{equation}
implies \eqref{e:GammaBoundedness_5}, 
\begin{equation}
    \|BA\chi u\|_{X_m^l}\le  C\|A\chi u\|_{X_m^l}+ C\|\opsi\chi u\|_{X_{-N}^{-N}}.
\end{equation}
For \eqref{e:GammaBoundedness_7}, note
\begin{equation}
    \|\rre^N\dde^{-1}B u\|\le \|\dde^{N-1}(\langle x\rangle^NB\langle x\rangle^{-N})\langle x\rangle^N u\|\le C\|\dde^{N-1}\langle x\rangle^N u\|\le C\|\rre^{2N}\dde^{-1} u\|
\end{equation}
and similarly
\begin{equation}
    \|\rre^N\nabla B u\|\le C\|\rre^N B \nabla u\|+C\|\rre^N O(\Psi^{0,-1}_{0,\delta}) u\|\le C\|\rre^{2N}\nabla u\|+C\|\rre^{2N}u\|
\end{equation}
to observe $B$ is bounded from $X_{2N+1}^{2N}$ to $X_{N}^{N}$ and also bounded from $\dot X_{2N+1}^{2N}$ to $\dot X_{N}^{N}$. Duality arguments then give the second half of \eqref{e:GammaBoundedness_7}. 
\end{proof}

\begin{lemma}
\label{l:resgammaBounds}
Let $\delta\in [0,1)$. Let $M_0>0$. Let $|m|<\rardyr_{\Omega}$ and $|m-l|<\frac12$. Fix $\chi, \varphi\in S^0(U)$. Then for all $s,k\in\mathbb{R}$,
\begin{equation}
    \Gamma^{s,0}_{k,\delta}\chi \mathcal{R}^\Omega_M\varphi \Gamma^{-s,0}_{-k,\delta}:Y_m^l(\Omega_+)\to \dot X_m^l(\Omega_+)
\end{equation}
is bounded uniformly in $M\in[0,M_0]$. 
\end{lemma}
\begin{proof}
1. Let $\chi_1(\bar x), \chi_2(\bar x), \chi_3(\bar x)\in S^0(U)$ that $\supp\chi\cap\supp(1-\chi_1)=\supp\chi_1\cap\supp(1-\chi_2)=\supp\chi_2\cap\supp(1-\chi_3)=\emptyset$. Let $\varphi_1(\bar x), \varphi_2(\bar x)\in S^0(U)$ that $\supp\varphi\cap\supp(1-\varphi_1)=\supp\varphi_1\cap\supp(1-\varphi_2)=\emptyset$. We claim that for any $N>0$ there exists a constant $C>0$ such that for all $v\in \mathcal{S}'(\Omega_+)$,
\begin{equation}\label{e:3l13}
    \|\Gamma^{s,0}_{k,\delta}\chi\mathcal{R}^\Omega_M\varphi v\|_{X_{m}^l}\le C\|\chi_2\Gamma^{s,0}_{k,\delta} v\|_{Y_m^l}+C\|\Gamma^{s,0}_{k,\delta} v\|_{Y_{-N}^{-N}}+C\|\Gamma^{-N,0}_{-N,\delta}\chi_1\mathcal{R}^\Omega_M\varphi v\|_{X_m^l}
\end{equation}
We will prove the claim in this step. Begin with noting that $\chi_1\Gamma^{s,0}_{k,\delta}\chi \mathcal{R}^\Omega_M\varphi v|_{x_0=0}\equiv 0$ and
\begin{equation}
 (-\Delta+M^2)\chi_1\Gamma^{s,0}_{k,\delta}\chi\mathcal{R}^\Omega_M\varphi v= \chi_1\Gamma^{s,0}_{k,\delta}\chi\varphi v+ [-\Delta, \chi_1\Gamma^{s,0}_{k,\delta} \chi] \mathcal{R}^\Omega_M\varphi v.
\end{equation}
The definition of $\mc{R}_M^\Omega$ implies 
\begin{equation}
    \chi_1\Gamma^{s,0}_{k,\delta}\chi\mathcal{R}^\Omega_M\varphi v=\chi_2\mathcal{R}^\Omega_M\chi_1\Gamma^{s,0}_{k,\delta}\chi\varphi v+ \chi_2\mathcal{R}^\Omega_M[-\Delta, \chi_1\Gamma^{s,0}_{k,\delta} \chi] \mathcal{R}^\Omega_M\varphi v,
\end{equation}
the first term of which can be estimated by
\begin{equation}
    \|\chi_2\mathcal{R}^\Omega_M\chi_1\Gamma^{s,0}_{k,\delta}\chi\varphi v\|_{X_m^l(\Omega_+)}\le C\|\chi_1\Gamma^{s,0}_{k,\delta}\chi\varphi v\|_{Y_m^l(\Omega)}\le C\|\chi_2 \Gamma^{s,0}_{k,\delta} v\|_{Y_m^l(\Omega_+)}+C\|\Gamma^{s,0}_{k,\delta} v\|_{Y_{-N}^{-N}(\Omega_+)}.
\end{equation}
Here we used the boundedness of $\mathcal{R}_{M}^\Omega: Y_m^l(\Omega)\rightarrow \dot X_m^l(\Omega)$ in Proposition \ref{p:extendNoMass}, and for the last inequality, consider 
\begin{equation}
    \chi_1\Gamma^{s,0}_{k,\delta}\chi\varphi v=\chi_1\Gamma^{s,0}_{k,\delta}\chi\varphi \Gamma^{-s,0}_{-k,\delta}\chi_2 \Gamma^{s,0}_{k,\delta}v+\chi_1\Gamma^{s,0}_{k,\delta}\chi\varphi \Gamma^{-s,0}_{-k,\delta}(1-\chi_2) \Gamma^{s,0}_{k,\delta}v,
\end{equation}
which implies for each $\phi\in \dot X_{-m}^{-l}(\Omega)$,
\begin{multline}
    |\langle \chi_1\Gamma^{s,0}_{k,\delta}\chi\varphi v, \phi\rangle|\le \|\chi_2\Gamma^{s,0}_{k,\delta} v\|_{Y_m^l(\Omega_+)}\|\Gamma^{-s,0}_{-k,\delta}\varphi\chi\Gamma^{s,0}_{k,\delta}\chi_1\phi\|_{X_{-m}^{-l}(\Omega_+)}\\
    +\|\Gamma^{s,0}_{k,\delta} v\|_{Y_{-N}^{-N}(\Omega_+)}\|(1-\chi_2)\Gamma^{-s,0}_{-k,\delta}\varphi\chi\Gamma^{s,0}_{k,\delta}\chi_1\phi\|_{X_{N}^{N}(\Omega_+)}\\
    \le C\big(\|\chi_2\Gamma^{s,0}_{k,\delta} v\|_{Y_m^l(\Omega_+)}+\|\Gamma^{s,0}_{k,\delta} v\|_{Y_{-N}^{-N}(\Omega_+)}\big)\|\phi\|_{X_{-m}^{-l}(\Omega)},
\end{multline}
by noting $\Gamma^{-s,0}_{-k,\delta}\varphi\chi\Gamma^{s,0}_{k,\delta}\in\Psi^{0,0}_{0,\delta}$ and $(1-\chi_2)\Gamma^{-s,0}_{-k,\delta}\varphi\chi\Gamma^{s,0}_{k,\delta}\chi_1=\opsi\chi_1$ to apply Lemma~\ref{l:xmlbdd}. To estimate the second term, Lemma \ref{l:CommDec} gives a decomposition of the commutator:
\begin{equation}
    [-\Delta, \chi_1\Gamma^{s,0}_{k,\delta} \chi]=\chi_2\big(\operatorname{div}A_0\nabla+\operatorname{div}A_1+A_2\big) \chi_1
\end{equation}
with $A_0\in \Psi^{s,-1+\delta}_{k-1,\delta}$, $A_1\in \Psi^{s,-1}_{k,\delta}$, $A_2\in \Psi^{s,-2}_{k,\delta}$. We can then bound
\begin{multline}
    \|\chi_2\mathcal{R}^\Omega_M[-\Delta, \chi_1\Gamma^{s,0}_{k,\delta} \chi] \mathcal{R}^\Omega_M\varphi v\|_{X_m^l(\Omega_+)}\le C\|\chi_2\operatorname{div}(A_0\nabla+A_1)\chi_1 \mathcal{R}^\Omega_M\varphi v\|_{Y_m^l(\Omega_+)}\\
    +C\|\chi_2 A_2\chi_1\mathcal{R}^\Omega_M\varphi v\|_{Y_m^l(\Omega_+)}.
\end{multline}
We first estimate the $A_2$ term. For each $\phi\in \dot X_{-m}^{-l}(\Omega_+)$, we have $\chi_2\phi\in \dot X^{-l}_{-m}(\Omega_+)$. By Lemma~\ref{l:xmlbdd}, 
\begin{equation}
    \|\rre^{-l}\dde^{-m+l}\Gamma^{0,-1}_{1,\delta}\chi_2\phi\|\le C\|\chi_2\phi\|_{X_{-m}^{-l}(\Omega_+)}\le C\|\phi\|_{X_{-m}^{-l}(\Omega_+)}.
\end{equation}
We then can bound
\begin{multline}
    |\langle \chi_2 A_2\chi_1\mathcal{R}^\Omega_M\varphi v, \phi\rangle|\le C\|\rre^{l}\dde^{m-l}\Gamma^{0,1}_{-1,\delta}A_2\chi_1\mathcal{R}^\Omega_M\varphi v\|\|\rre^{-l}\dde^{-m+l}\Gamma^{0,-1}_{1,\delta}\chi_2\phi\|\\
    \le C\|\rre^{l}\dde^{m-l}\Gamma^{0,1}_{-1,\delta}A_2\chi_1\mathcal{R}^\Omega_M\varphi v\|\|\phi\|_{X_{-m}^{-l}(\Omega_+)}.
\end{multline}
This implies 
\begin{multline}
    \|\chi_2 A_2\chi_1\mathcal{R}^\Omega_M\varphi v\|_{Y_m^l(\Omega_+)}\le C\|\rre^{l}\dde^{m-l}\Gamma^{0,1}_{-1,\delta}A_2\chi_1\mathcal{R}^\Omega_M\varphi v\|\\
    \le C\|\rre^{l}\dde^{m-l}\chi_2\Gamma^{0,1}_{-1,\delta}A_2\chi_1\mathcal{R}^\Omega_M\varphi v\|+C\|\opsi\chi_1\mathcal{R}^\Omega_M\varphi v\|
\end{multline}
By the ellipticity of $\rre^{l}\dde^{m-l}\chi_3\Gamma^{s,-1}_{k-1,\delta}$ over the microsupport of $\rre^{l}\dde^{m-l}\chi_2\Gamma^{0,1}_{-1,\delta}A_2$, we have 
\begin{multline}
    \|\chi_2 A_2\chi_1\mathcal{R}^\Omega_M\varphi v\|_{Y_m^l(\Omega_+)}\le C\|\rre^{l}\dde^{m-l}\chi_3\Gamma^{s,-1}_{k-1,\delta}\chi_1\mathcal{R}^\Omega_M\varphi v\|+C\|\opsi\chi_1\mathcal{R}^\Omega_M\varphi v\|\\
    \le C\|\Gamma^{s,0}_{k-1,\delta}\chi_1\mathcal{R}^\Omega_M\varphi v\|_{X_m^l(\Omega_+)}+C\|\opsi\chi_1\mathcal{R}^\Omega_M\varphi v\|_{X_{-N}^{-N}(\Omega_+)}\le C\|\Gamma^{s,0}_{k-1,\delta}\chi_1\mathcal{R}^\Omega_M\varphi v\|_{X_m^l(\Omega_+)}.
\end{multline}
We now deal with the $(A_0\nabla+A_1)$-term. Then
\begin{multline}
    |\langle \chi_2\operatorname{div}(A_0\nabla+A_1)\chi_1 \mathcal{R}^\Omega_M\varphi v, \phi\rangle|\le C\|\rre^l\dde^{m-l}(A_0\nabla+A_1)\chi_1 \mathcal{R}^\Omega_M\varphi v\|\|\rre^{-l}\dde^{-m+l}\nabla\phi\|\\
    \le C\big(\|\rre^l\dde^{m-l}A_0\nabla \chi_1 \mathcal{R}^\Omega_M\varphi v\|+\|\rre^l\dde^{m-l}A_1 \chi_1 \mathcal{R}^\Omega_M\varphi v\|\big)\|\phi\|_{X_{-m}^{-l}(\Omega_+)},
\end{multline}
With $A_0\nabla=\nabla A_0+O(\Psi^{s,-2+\delta}_{k-1,\delta})$, Lemma~\ref{l:xmlbdd} then implies
\begin{equation}
    \|\chi_2\operatorname{div}(A_0\nabla+A_1)\chi_1 \mathcal{R}^\Omega_M\varphi v\|_{Y_m^l(\Omega_+)}\le C\|\Gamma^{s,0}_{k-1,\delta}\chi_1 \mathcal{R}^\Omega_M\varphi v\|_{X_{m}^{l}(\Omega_+)}.
\end{equation}
Summing up the estimates above and noting $\Gamma^{\frac12,0}_{-\frac12,\delta}\in\Psi^{0,0}_{0,\delta}$ is bounded on $X_{m}^{l}(\Omega_+)$, we have
\begin{equation}\label{e:3l38}
    \|\Gamma^{s,0}_{k,\delta}\chi\mathcal{R}^\Omega_M\varphi v\|_{X_{m}^l}\le C\|\chi_2\Gamma^{s,0}_{k,\delta} v\|_{Y_m^l}+C\|\Gamma^{s,0}_{k,\delta} v\|_{Y_{-N}^{-N}}+C\|\Gamma^{s-\frac12,0}_{k-\frac12,\delta}\chi_1\mathcal{R}^\Omega_M\varphi v\|_{X_m^l}.
\end{equation}
Repeated use of \eqref{e:3l38} gives \eqref{e:3l13} with cutoff functions renamed.

2. Apply \eqref{e:3l13} with $v=\Gamma^{-s,0}_{-k,\delta}u$, we have
\begin{equation}\label{e:3l14}
    \|\Gamma^{s,0}_{k,\delta}\chi\mathcal{R}^\Omega_M\varphi \Gamma^{-s,0}_{-k,\delta}u\|_{X_{m}^l}\le C\|\chi_2u\|_{Y_m^l}+C\|u\|_{Y_{-N}^{-N}}\\
    +C\|\Gamma^{-N,0}_{-N,\delta}\chi_1\mathcal{R}^\Omega_M\varphi \Gamma^{-s,0}_{-k,\delta}u\|_{X_m^l}.
\end{equation}
The first two terms are bounded by $C\|u\|_{Y_m^l}$. We will control the last term. We claim
\begin{equation}
    T=\Gamma^{-N,0}_{-N,\delta}\chi_1\mathcal{R}^\Omega_M\varphi \Gamma^{-s,0}_{-k,\delta}: Y_m^l(\Omega_+)\rightarrow \dot X_m^l(\Omega_+)
\end{equation}
is bounded. This is equivalent to show its adjoint
\begin{equation}
    T^*=\Gamma^{-s,0}_{-k,\delta}\varphi\mathcal{R}^\Omega_M\chi_1 \Gamma^{-N,0}_{-N,\delta}: Y_{-m}^{-l}(\Omega_+)\rightarrow \dot X_{-m}^{-l}(\Omega_+)
\end{equation}
is bounded. Here we used the fact that the adjoint of $\mathcal{R}^\Omega_M$ coincides with itself on the spaces of inverted parameters. Now for $w\in Y_{-m}^{-l}(\Omega_+)$, apply \eqref{e:3l13} with inverted parameters and $v=\Gamma^{-N,0}_{-N,\delta}w$ to see for $N>0$ large, we have
\begin{multline}
    \|T^*w\|_{X_{-m}^{-l}}=\|\Gamma^{-s,0}_{-k,\delta}\varphi\mathcal{R}^\Omega_M\chi_1 \Gamma^{-N,0}_{-N,\delta} w\|_{X_{-m}^{-l}(\Omega_+)}\le C\|\varphi_2 \Gamma^{-s-N,0}_{-k-N,\delta}w\|_{Y_{-m}^{-l}(\Omega_+)}+C\|\Gamma^{-s-N,0}_{-k-N,\delta}w\|_{Y_{-N}^{-N}(\Omega_+)}\\+C\|\Gamma^{-N,0}_{-N,\delta}\varphi_1\mathcal{R}^\Omega_M\chi_1\Gamma^{-N,0}_{-N,\delta} w\|_{X_{-m}^{-l}(\Omega_+)}\le C\|w\|_{Y_{-m}^{-l}(\Omega_+)}.
\end{multline} 
Note here we used by Lemma~\ref{l:xmlbdd}, for $N>0$ large we have
\begin{gather}
    \|\varphi_2 \Gamma^{-s-N,0}_{-k-N,\delta}w\|_{Y_{-m}^{-l}(\Omega_+)}\le C\|w\|_{Y_{-m}^{-l}(\Omega_+)}, \ \|\Gamma^{-s-N,0}_{-k-N,\delta}w\|_{Y_{-N}^{-N}(\Omega_+)}\le C\|w\|_{Y_{-N/2}^{-N/2}(\Omega_+)},
\end{gather}
and since $\mc{R}^\Omega_M: Y_m^l(\Omega)\rightarrow \dot X_m^l(\Omega)$ is bounded, we also have
\begin{multline}
    \|\Gamma^{-N,0}_{-N,\delta}\varphi_1\mathcal{R}^\Omega_M\chi_1\Gamma^{-N,0}_{-N,\delta} w\|_{X_{-m}^{-l}(\Omega_+)}\le C\|\varphi_1\mathcal{R}^\Omega_M\chi_1\Gamma^{-N,0}_{-N,\delta} w\|_{X_{-m}^{-l}(\Omega)}\\
    \le C\|\chi_1\Gamma^{-N,0}_{-N,\delta} w\|_{Y_{-m}^{-l}(\Omega)}\le C\|w\|_{Y_{-m}^{-l}(\Omega_+)}.
\end{multline}
Now $T$ and $T^*$ are bounded. Revisit \eqref{e:3l14} to conclude the proof.
\end{proof}

\subsection{Structure of the \dtn map on \acl domains}
\label{s:dtnStructure}
\begin{definition}
    We say two \acl domains $\Omega_1,\Omega_2$ share the common boundary near a closed set $K\subset\partial\Omega_1\cap\partial\Omega_2$, if there is $U\subset \Omega_1\cap\Omega_2$ as an open subset of the manifolds with boundary $\Omega_1,\Omega_2$, that $U$ contains $K$. Furthermore, we say they share the common smooth boundary, if $\partial\Omega_1\cap\partial\Omega_2$ is further smooth.
\end{definition}
\begin{remark}[Geometric setup of \acl domains]\label{r:aclcollar}
    Let $\Omega$ be an \acl domain. We define the regular set $\partial\oreg$ consists of $x\in\partial\Omega$ that $\partial\Omega$ is smooth near $x$, and the singular set $\partial\osing=\partial\Omega\setminus \partial\oreg$. Near any closed $K\subset\partial\oreg$, there exists an \acs domain $\tilde\Omega\subset\Omega$ that $\Omega, \tilde\Omega$ share the common smooth boundary near $K$, associated with its conic collar $U$, half-space collar $\tilde\Omega_+$, $P$, $E_P$ given in Lemma~\ref{l:factor} and Lemma~\ref{l:dtnParametrixExtension}. We call such $U$ a regular conic collar of $\Omega$. See Figure~\ref{f:3f2}(B).
\end{remark}

In this setup, we now show the main result in this section: $P$ describes the structure of $\Lambda^\Omega_M$, the Dirichlet to Neumann map on $\Omega$, with an error that is arbitrarily smooth and decays fast up to a threshold at the zero frequency. 

    Let $\chi(x),\tilde\chi(x)\in S^0(\partial\oreg)$ with $\supp\chi\cap\supp(1-\tilde\chi)=\emptyset$. Here we confuse $\chi, \tilde\chi$ with any of their extensions in $S^0(U)$ that $\supp\tilde\chi\subset\overline{\Omega}$, $\supp\chi\cap\supp(1-\tilde\chi)=\emptyset$ and $\partial_\nu \tilde\chi=\partial_\nu\chi=0$. 
\begin{proposition}[Structure of $\Lambda^\Omega_M$ and $\mathcal{H}^\Omega_M$]
\label{p:dtNParametrix}
    Let $\Omega\subset\mathbb{R}^d$ be an \acl domain with a regular collar $U$. Let $M_0>0$. Let $\delta\in [0,1)$. Let $\chi,\tilde\chi, \tilde{\tilde{\chi}}\in S^0(U)$ with $\supp\chi\cap\supp(1-\tilde\chi)=\supp\tilde\chi\cap\supp(1-\tilde{\tilde{\chi}})=\emptyset$. Let $l,m, L\in\mathbb{R}$ such that $|m|< \rardyr_{\Omega}$, $|l-m|<\frac12$ and
    \begin{equation}
        L\ge \max(l, m+1)-\frac32\delta.
    \end{equation}
    For any $N>0$, we have
    \begin{equation}\label{e:DtN_Approx}
        \|\tilde\chi \Gamma^{N,l}_{N,\delta}\chi (\Lambda^\Omega_M-P)\tilde\chi f\|_{L^2(\partial\Omega)}\le C_N\|\tilde{\tilde{\chi}}\Gamma^{-N,L}_{-N,\delta}\tilde\chi f\|_{L^2(\partial\Omega)}+C_N\|\tilde\chi f\|_{H^{-N,-N}(\partial\Omega)},
    \end{equation}
    and
    \begin{equation}\label{e:3l12}
    \|\tilde\chi\Gamma^{N,0}_{N,\delta}\chi(\mathcal{H}^\Omega_M-E_P)\tilde\chi f\|_{X_{m}^{l}(\Omega)}\le C_N\|\tilde{\tilde{\chi}}\Gamma^{-N, L}_{-N, \delta}\tilde\chi f\|_{L^2(\partial\Omega)}+C_N\|\tilde\chi f\|_{H^{-N,-N}(\partial\Omega)},
    \end{equation}
    uniformly for all $M\in[0,M_0]$ and $f\in \mc{S}'(\partial\Omega)$. 
\end{proposition}
\begin{remark}
    Here the error $\mc{E}_0$, between $\Lambda_M^\Omega$ and its approximation $P$, lives on the boundary and consists of the interior error $\mc{E}_0^\circ$ induced in the boundary decomposition, and other boundary correction terms. 
\end{remark}

\begin{proof}
    0. It suffices to show the claim $f\in \Gamma^{-N,L}_{-N,\delta}L^2(\partial\Omega)\cap \mc{S}(\partial\Omega)$ and extend it to $\mc{S}(\partial\Omega)$ using density. 

    1. Let $\chi_1(x_0,x),\chi_2(x_0,x)\in S^0(U)$ that $\supp\tilde\chi\cap\supp(1-\chi_1)=\supp\chi_1\cap\supp(1-\chi_2)=\emptyset$, and $\supp\chi_2\subset\overline{\Omega}$. Note $\chi_1 E_{P}\tilde\chi f$ extends $\tilde\chi f$ into $\Omega$, and the extension correction formula \eqref{e:hForm} implies
    \begin{multline}
        \chi\mathcal{H}_M^\Omega \tilde\chi f=\chi E_{P}\tilde\chi f-\chi\mathcal{R}_M^\Omega (-\Delta+M^2)\chi_1 E_{P}\tilde\chi f=\chi E_{P} \tilde\chi f-\chi \mathcal{R}_M^\Omega\chi_1 (-\Delta+M^2)\chi_2 E_{P}\tilde\chi f-\chi\mathcal{R}_M^\Omega [-\Delta,\chi_1]E_{P}\tilde\chi f  \\
        =\chi E_{P}\tilde\chi f-\chi\mathcal{R}_M^\Omega \chi_1(D_{x_0}+iP+b)(D_{x_0}-iP)\chi_2 E_{P}\tilde\chi f-\chi\mathcal{R}_M^\Omega[-\Delta,\chi_1]E_{P}\tilde\chi f-\chi\mathcal{R}_M^\Omega\chi_1 \mathcal{E}_0^\circ\chi_2 E_{P}\tilde\chi f\\
        =\chi E_{P}\tilde\chi f -
        \chi\mathcal{R}_M^\Omega[-\Delta,\chi_1]E_{P}\tilde\chi f-\chi\mathcal{R}_M^\Omega\chi_1 \mathcal{E}_0^\circ\chi_2 E_{P}\tilde\chi f 
         +
        \chi\mathcal{R}_M^\Omega \chi_1(D_{x_0}+iP+b)(D_{x_0}-iP)(1-\chi_2) E_{P}\tilde\chi f,
    \end{multline}
    where we used that $(D_{x_0}-iP)E_{P}=0$. Now for the last term,
    \begin{multline}
        \chi_1(D_{x_0}+iP+b)(D_{x_0}-iP)(1-\chi_2)=\chi_1(D_{x_0}^2+b(x_0, x)D_{x_0}+i[P,D_{x_0}]+P^2-ib(x_0,x)P)(1-\chi_2)\\
        =\chi_1O(\Psi^{-\infty,-\infty}_{-\infty,\delta})
    \end{multline}
    since $\chi_1$ and $1-\chi_2$ have disjoint support. Note $\mathcal{E}_0^\circ\in \Psi^{-\infty,-2\delta}_{-\infty,\delta}$ and we arrive at
    \begin{equation}
        \chi(\mathcal{H}_M^\Omega-E_{P})\tilde\chi f=\chi\mathcal{R}_M^\Omega \chi_2 v, \ v=(-[-\Delta, \chi_1]+\chi_1O(\Psi^{-\infty, -2\delta}_{-\infty,\delta}))E_{P}\tilde\chi f.
    \end{equation}
    Here and later we repeatedly used for any $s,l',k\in\mathbb{R}$ that
    \begin{equation}\label{e:3l34}
        \|\Gamma^{s,l'}_{k,\delta} E_{P}\tilde\chi f\|\le C\|\Gamma^{s,l'+\frac\delta2}_{k-\frac12,\delta} \tilde\chi f\|_{L^2(\partial\tilde\Omega)}, \ \|O(\Psi^{-\infty,-\infty}_{-\infty,\delta}) E_{P}\tilde\chi f\|
    \le C\|\Gamma^{-N, -N}_{-N,\delta}\tilde\chi f\|_{L^2(\partial\tilde\Omega)}
    \end{equation}
    from Lemma~\ref{l:dtnParametrixExtension}. We will focus on estimating the size of $v$. 

2. In this step, we estimate the size of $\Gamma^{s,0}_{k,\delta} v$ in $L^2(\tilde\Omega_+)$ and in $Y_m^l(\tilde\Omega_+)$.
We claim for any $N>0$, $s, l', k\in\mathbb{R}$, we have
\begin{equation}\label{e:3l10}
    \|\Gamma^{s,l'}_{k,\delta} v\|_{Y_m^l(\tilde\Omega_+)}\le C\|\Gamma^{-N, l'+\max(l,m+1)-\frac{3}2\delta}_{-N,\delta} \tilde\chi f\|_{L^2(\partial\tilde\Omega)}\le C\|\tilde{\tilde{\chi}}\Gamma^{-N, l'+L}_{-N,\delta} \tilde\chi f\|_{L^2(\partial\Omega)}+C\|\tilde\chi f\|_{H^{-N, -N}(\partial\Omega)},
\end{equation}
and
\begin{equation}\label{e:3l11}
    \|\Gamma^{s,l'}_{k,\delta}v\|_{L^2(\tilde\Omega_+)}\le C\|\Gamma^{-N, l'-\frac{3}2\delta}_{-N,\delta}\tilde\chi f\|_{L^2(\partial\tilde\Omega)}\le C\|\tilde{\tilde{\chi}}\Gamma^{-N, l'-\frac{3}2\delta}_{-N,\delta}\tilde\chi f\|_{L^2(\partial\Omega)}+C\|\tilde\chi f\|_{H^{-N, -N}(\partial\Omega)}.
\end{equation}
Note that \eqref{e:3l10} implies with $l'=0$, by Lemma \ref{l:resgammaBounds},
\begin{equation}
    \|\tilde\chi\Gamma^{s,0}_{k,\delta}\chi(\mathcal{H}^\Omega_M-E_{P})\tilde\chi f\|_{X_m^l(\Omega)}\le C\|\Gamma^{s,0}_{k,\delta}\chi\mathcal{R}^\Omega_M \chi_2 v\|_{X_m^l(\tilde\Omega_+)}\le C\|\Gamma^{s,0}_{k,\delta} v\|_{Y_m^l(\tilde\Omega_+)}\le C\|\Gamma^{-N, L}_{-N, \delta} \tilde\chi f\|_{L^2(\partial\tilde\Omega)},
\end{equation}
that is \eqref{e:3l12}. The rest of Step 2 is devoted to prove those two estimates. 

2a. We first estimate the size of $\Gamma^{s,l'}_{k,\delta}[-\Delta, \chi_1]E_{P}\tilde\chi f$. Using that $(D_{x_0}-iP)E_{P}=0$, we have 
\begin{multline}
    [-\Delta, \chi_1]E_{P}\tilde\chi f=-(\Delta\chi_1) E_{P}\tilde\chi f-2i(\partial_{x_0}\chi_1)D_{x_0} E_{P}\tilde\chi f-2\sum_{j,k=1}^{d-1}(\partial_{x_j} \chi_1)g_+^{jk}\partial_{x_k} E_{P}\tilde\chi f\\
    =-(\Delta\chi_1) E_{P}\tilde\chi f+2\big((\partial_{x_0}\chi_1)P-\sum_{j,k=1}^{d-1}(\partial_{x_j} \chi_1)g_+^{jk}\partial_{x_k}\big)E_{P}\tilde\chi f=BE_{P}\tilde\chi f,
\end{multline}
where $B\in\Psi^{0,-1-\delta}_{1,\delta}$ and $B=\varphi_0 B\varphi_0+\varphi_0O(\Psi^{-\infty,-\infty}_{-\infty,\delta})$ for some $\varphi_0(x_0,x)\in S^0(U)$ supported in some $\{x_0< c\langle x\rangle\}$ for some $c>0$, and $\varphi_0\nabla \chi_1=\nabla \chi_1$. Let $\varphi\in S^0(U)$ satisfy all conditions for $\varphi_0$ with $\supp\varphi_0\cap\supp(1-\varphi)=\emptyset$. Then $(1-\varphi)\Gamma^{s,l'}_{k,\delta}\varphi_0\in O(\Psi^{-\infty,-\infty}_{-\infty,\delta})$ implies
\begin{equation}\label{e:3l32}
    \|\Gamma^{s,l'}_{k,\delta} [-\Delta,\chi_1]E_{P}\tilde\chi f\|\le C\|\varphi\Gamma^{s,l'}_{k,\delta}B\varphi_0E_{P}\tilde\chi f\|+\|O(\Psi^{-\infty,-\infty}_{-\infty,\delta}) \varphi E_{P}\tilde\chi f\|\le C\|\Gamma^{-N, -N}_{-N,\delta}\tilde\chi f\|_{L^2(\partial\tilde\Omega)},
\end{equation}
since $\varphi_0E_{P}\tilde\chi\in \opsi$ from Lemma~\ref{l:extensiongainsdecay} and $O(\Psi^{-\infty,-\infty}_{-\infty,\delta}) \varphi=\opsi$. We have
\begin{multline}
    \|\rre^{l}\dde^{m+1-l} \Gamma^{s,l'}_{k,\delta} [-\Delta,\chi_1]E_{P}\tilde\chi f\|\le C\|\varphi\Gamma^{s,l'+m+1}_{k,\delta}BE_{P}\tilde\chi f\|+\|O(\Psi^{-\infty,-\infty}_{-\infty,\delta}) \varphi E_{P}\tilde\chi f\|\le\\
    C\|\Gamma^{-N, -N}_{-N,\delta}\tilde\chi f\|_{L^2(\partial\tilde\Omega)}.
\end{multline}
and similarly 
\begin{multline}
    \|\rre^{l}\dde^{m+2-l}\nabla_{x}\Gamma^{s,l'}_{k,\delta}[-\Delta, \chi_1]E_{P}\tilde\chi f\|\le C\|\varphi\Gamma^{s,l'+m+2}_{k,\delta}\nabla_{x}BE_{P}\tilde\chi f\|+\|O(\Psi^{-\infty,-\infty}_{-\infty,\delta}) \varphi E_{P}\tilde\chi f\|\\
    \le C\|\Gamma^{-N, -N}_{-N,\delta}\tilde\chi f\|_{L^2(\partial\tilde\Omega)}.
\end{multline}
Now noting $(D_{x_0}-iP)E_{P}=0$ and $[D_{x_0}, \Gamma^{s,l'}_{k,\delta}]=0$, we have
\begin{equation}
    D_{x_0}\Gamma^{s,l'}_{k,\delta}[-\Delta, \chi_1]E_{P}\tilde\chi f=D_{x_0}\Gamma^{s,l'}_{k,\delta}BE_{P}\tilde\chi f=i\Gamma^{s,l'}_{k,\delta}BPE_{P}\tilde\chi f+\Gamma^{s,l'}_{k,\delta}[D_{x_0}, B]E_{P}\tilde\chi f,
\end{equation}
where $\Gamma^{s,l'}_{k,\delta}[D_{x_0}, B]\in \Psi^{s,l'-2-\delta}_{k+1,\delta}$ 
has microsupport inside $\{\varphi\equiv 1\}$. Thus we can similarly estimate
\begin{gather}
    \|\rre^{l}\dde^{m+2-l} \Gamma^{s,l'}_{k,\delta}BPE_{P}\tilde\chi f\|\le C\|\Gamma^{-N, -N}_{-N,\delta}\tilde\chi f\|_{L^2(\partial\tilde\Omega)}, \\
    \|\rre^{l}\dde^{m+2-l} \Gamma^{s,l'}_{k,\delta}[D_{x_0}, B]E_{P}\tilde\chi f\|\le C\|\Gamma^{-N, -N}_{-N,\delta}\tilde\chi f\|_{L^2(\partial\tilde\Omega)},
\end{gather}
and ultimately we have
\begin{equation}
    \|\rre^{l}\dde^{m+2-l} D_{x_0}\Gamma^{s,l'}_{k,\delta}[-\Delta, \chi_1]E_{P}\tilde\chi f\|\le C\|\Gamma^{-N, -N}_{-N,\delta}\tilde\chi f\|_{L^2(\partial\tilde\Omega)}.
\end{equation}
Now we have
\begin{multline}\label{e:3l30}
    \|\Gamma^{s,l'}_{k,\delta}[-\Delta, \chi_1]E_{P}\tilde\chi f\|_{Y_m^l(\tilde\Omega_+)}\le C\|\Gamma^{s,l'}_{k,\delta}[-\Delta, \chi_1]E_{P}\tilde\chi f\|_{X_{m+2}^l(\tilde\Omega_+)}\\\le C\|\rre^{l}\dde^{m+2-l}\nabla \Gamma^{s,l'}_{k,\delta} [-\Delta, \chi_1]E_{P}\tilde\chi f\|
    +C\|\rre^{l}\dde^{m+1-l}\Gamma^{s,l'}_{k,\delta}[-\Delta, \chi_1]E_{P}\tilde\chi f\|\le C\|\Gamma^{-N, -N}_{-N,\delta}\tilde\chi f\|_{L^2(\partial\tilde\Omega)},
\end{multline}
for any $N>0$. 

2b. We estimate the size of $\Gamma^{s,l'}_{k,\delta}\chi_1 O(\Psi^{-\infty, -2\delta}_{-\infty,\delta}) E_{P}\tilde\chi f$. Note 
\begin{equation}\label{e:3l33}
    \|\Gamma^{s,l'}_{k,\delta}\chi_1 O(\Psi^{-\infty,-2\delta}_{-\infty,\delta}) E_{P}\tilde\chi f\|\le C\|\Gamma^{-N,l'-2\delta}_{-N,\delta} E_{P}\tilde\chi f\|\le C\|\Gamma^{-N,l'-\frac32\delta}_{-N,\delta} \tilde\chi f\|_{L^2(\partial\tilde\Omega)},
\end{equation}
in the last step of which we used \eqref{e:3l34}. Consider that for each $\phi\in \dot X_{-m}^{-l}(\tilde\Omega_+)$, 
\begin{multline}
    |\langle \Gamma^{s,l'}_{k,\delta}\chi_1O(\Psi^{-\infty, -2\delta}_{-\infty,\delta}) E_{P}\tilde\chi f, \phi\rangle|=|\langle \rre^{l}\dde^{m-l+1}\Gamma^{s,l'}_{k,\delta}\chi_1O(\Psi^{-\infty, -2\delta}_{-\infty,\delta}) E_{P}\tilde\chi f, \rre^{-l}\dde^{-m+l-1}\phi\rangle|\\
    \le C\|\rre^{l}\dde^{m-l+1}\Gamma^{s,l'}_{k,\delta}\chi_1O(\Psi^{-\infty, -2\delta}_{-\infty,\delta}) E_{P}\tilde\chi f\|\|\phi\|_{X_{-m}^{-l}(\tilde\Omega_+)}.
\end{multline}
Noting $\rre^l\dde^{m+1-l}\le \langle x\rangle^{\max(l,m+1)}$ over $\supp\chi_2$, we have by \eqref{e:3l34},
\begin{multline}\label{e:3l31}
    \|\Gamma^{s,l'}_{k,\delta}\chi_1O(\Psi^{-\infty, -2\delta}_{-\infty,\delta}) E_{P}\tilde\chi f\|_{Y_m^l(\tilde\Omega_+)}\le C\|\rre^{l}\dde^{m+1-l}\Gamma^{s,l'}_{k,\delta}\chi_1O(\Psi^{-\infty, -2\delta}_{-\infty,\delta}) E_{P}\tilde\chi f\|\\ 
    \le C\|O(\Psi^{-\infty, l'+\max(l,m+1)-2\delta}_{-\infty,\delta}) E_{P}\tilde\chi f\|+C\|\opsi E_{P}\tilde\chi f\|
    \le C\|\Gamma^{-N, l'+\max(l,m+1)-\frac{3}2\delta}_{-N,\delta} \tilde\chi f\|_{L^2(\partial\tilde\Omega)},
\end{multline}
where we used $(1-\chi_2)\Gamma^{s,l'}_{k,\delta}\chi_1O(\Psi^{-\infty, -2\delta}_{-\infty,\delta})=\opsi$. Sum up \eqref{e:3l30}, \eqref{e:3l31} to observe \eqref{e:3l10}, and sum up \eqref{e:3l32}, \eqref{e:3l33} to observe \eqref{e:3l11}. This proves \eqref{e:3l12}. 

3. It remains to show \eqref{e:DtN_Approx}. Noting $(\mathcal{H}_M^\Omega-E_{P})\tilde\chi f$ has zero trace, one can estimate 
\begin{multline}\label{e:3l8}
    \|\tilde\chi \Gamma^{s,l}_{k,\delta}D_{x_0}\dde^{-N}\chi (\mathcal{H}_M^\Omega-E_{P})\tilde\chi f\|_{H^{-\frac12}(\partial\tilde\Omega)}=\|D_{x_0} \dde^{-N} \tilde\chi  \Gamma^{s,l}_{k,\delta}\chi (\mathcal{H}_M^\Omega-E_{P})\tilde\chi f\|_{H^{-\frac12}(\partial\tilde\Omega)}\\
    \le \|\Delta \dde^{-N} \tilde\chi  \Gamma^{s,l}_{k,\delta}\chi (\mathcal{H}_M^\Omega-E_{P})\tilde\chi f\|_{L^2}
    +\|\nabla \dde^{-N} \tilde\chi  \Gamma^{s,l}_{k,\delta}\chi (\mathcal{H}_M^\Omega-E_{P})\tilde\chi f\|_{L^2}\\
    =\|\Delta \dde^{-N} \tilde\chi \Gamma^{s,l}_{k,\delta}\chi\mathcal{R}_M^\Omega \chi_2 v\|+\|\nabla \dde^{-N} \tilde\chi  \Gamma^{s,l}_{k,\delta}\chi\mathcal{R}_M^\Omega\chi_2 v\|,
\end{multline}
where $N>0$ will be chosen such that $-N+l+2\le m$ to embed spaces into $X_m^l(\Omega_+)$ in this step. 

3a. We take care of the term $\nabla \dde^{-N}\tilde\chi  \Gamma^{s,l}_{k,\delta}\chi\mathcal{R}_M^\Omega \chi_2 v$ in this step. Consider
\begin{multline}
    \|\nabla \dde^{-N} \tilde\chi \Gamma^{s,l}_{k,\delta}\chi\mathcal{R}_M^\Omega \chi_2 v\|\le C\|\nabla \rre^l\dde^{-N}\tilde\chi O(\Psi^{s,0}_{k,\delta})\chi\mathcal{R}_M^\Omega \chi_2 v\|\le C\|\Gamma^{s,0}_{k,\delta}\chi\mathcal{R}_M^\Omega \chi_2 v\|_{X_{{-N}+l}^l(\tilde\Omega_+)}\\
    \le C\|\Gamma^{s,0}_{k,\delta}\chi\mathcal{R}_M^\Omega \chi_2 v\|_{X_m^l(\tilde\Omega_+)}
\end{multline}
for large $N>0$, where we used $\nabla \rre^l\dde^{-N}=\rre^l\dde^{-N}\nabla + O(\rre^l\dde^{-N-1})$ and the boundedness of $\nabla O(\Psi^{s,0}_{k,\delta})$ on $X_m^l$ spaces as in Lemma \ref{l:xmlbdd}. 

3b. We take care of the term $\Delta \dde^{-N} \tilde\chi \Gamma^{s,l}_{k,\delta}\chi\mathcal{R}_M^\Omega \chi_2 v$ in this step. Consider 
\begin{multline}\label{e:3l39}
    -\Delta \dde^{-N} \tilde\chi \Gamma^{s,l}_{k,\delta}\chi\mathcal{R}_M^\Omega \chi_2 v=\dde^{-N} \tilde\chi \Gamma^{s,l}_{k,\delta}\chi(-\Delta)\mathcal{R}_M^\Omega \chi_2 v+[-\Delta,  \dde^{-N} \tilde\chi \Gamma^{s,l}_{k,\delta}\chi ]\chi_1\mathcal{R}_M^\Omega \chi_2 v\\
    =\dde^{-N} \tilde\chi \Gamma^{s,l}_{k,\delta}\chi v-M^2 \dde^{-N} \tilde\chi \Gamma^{s,l}_{k,\delta}\chi\mathcal{R}_M^\Omega \chi_2 v+[-\Delta,  \dde^{-N} \tilde\chi \Gamma^{s,l}_{k,\delta}\chi ]\chi_1 \mathcal{R}_M^\Omega\chi_2 v.
\end{multline}
For the second term of \eqref{e:3l39}, we have from elliptic estimates
\begin{multline}
    \|\dde^{-N} \tilde\chi \Gamma^{s,l}_{k,\delta}\chi\mathcal{R}_M^\Omega \chi_2 v\|\le C\|\rre^l\dde^{-N}\chi_2\Gamma^{s,0}_{k,\delta} \chi\mathcal{R}^\Omega_M \chi_2 v\|+C\|\opsi \chi\mathcal{R}^\Omega_M \chi_2 v\|\\
    \le C\|\Gamma^{s,0}_{k, \delta}\chi \mathcal{R}_M^\Omega \chi_2 v\|_{X_{-N+l+1}^l} \le C\|\Gamma^{s,0}_{k, \delta}\chi \mathcal{R}_M^\Omega \chi_2 v\|_{X_{m}^l}.
\end{multline}
For the last term of \eqref{e:3l39}, note
\begin{equation}\label{e:3l40}
    [-\Delta,  \dde^{-N} \tilde\chi \Gamma^{s,l}_{k,\delta}\chi ]\chi_1 \mathcal{R}_M^\Omega\chi_2 v=[-\Delta,  \dde^{-N} \tilde\chi ]\Gamma^{s,l}_{k,\delta}\chi\mathcal{R}_M^\Omega\chi_2 v+ \dde^{-N} \tilde\chi [-\Delta, \Gamma^{s,l}_{k,\delta}\chi ]\chi_1 \mathcal{R}_M^\Omega\chi_2 v
\end{equation}
the first term of which can be bounded by Lemma~\ref{l:xmlbdd}: for large $N>0$, 
\begin{multline}
    \|[-\Delta,  \dde^{-N} \tilde\chi ]\Gamma^{s,l}_{k,\delta}\chi \mathcal{R}_M^\Omega \chi_2 v\|\le C\|\dde^{{-N}-1}\chi_1\nabla \Gamma^{s,l}_{k,\delta}\chi \mathcal{R}_M^\Omega \chi_2 v\|+C\|\dde^{{-N}-2}\chi_1 \Gamma^{s,l}_{k,\delta}\chi \mathcal{R}_M^\Omega \chi_2 v\|\\
    \le C\|\rre^{l}\dde^{{-N}-1}\chi_1 \nabla O(\Psi^{s,0}_{k,\delta})\chi \mathcal{R}_M^\Omega \chi_2 v\|+ C\|\rre^{l-1}\dde^{{-N}-1}\chi_1 O(\Psi^{s,0}_{k,\delta})\chi \mathcal{R}_M^\Omega \chi_2 v\|\\
    +C\|\rre^{l}\dde^{{-N}-2}\chi_1 O(\Psi^{s,0}_{k,\delta})\chi \mathcal{R}_M^\Omega \chi_2 v\|
    \le C\|\Gamma^{s,0}_{k,\delta} \chi \mathcal{R}_M^\Omega \chi_2 v\|_{X_{m}^{l}}.
\end{multline}
Noting that $-\Delta_\gplus=D_{x_0}^2+b(\bar x)D_{x_0}+O(\Psi^{0,-1-\delta}_{1,\delta})+O(\Psi^{0,-2\delta}_{2,\delta})$ and $[D_{x_0}, \Gamma^{s,l}_{k,\delta}]=0$, we have 
\begin{multline}
    \dde^{-N}\tilde\chi [-\Delta, \Gamma^{s,l}_{k,\delta}\chi ]\chi_1=\dde^{-N}\tilde\chi (O(\Psi^{s,l-1}_{k,\delta})D_{x_0}+O(\Psi^{s,l-1-\delta}_{k+1,\delta}))\chi_1\\
    =\dde^{-N}\big(\rre^{l-1}\tilde\chi  D_{x_0}O(\Psi^{s,0}_{k,\delta})
    +
    \rre^{l-2}\tilde\chi  O(\Psi^{s,0}_{k,\delta}) 
    +\rre^{l-1-\delta}\tilde\chi O(\Psi^{s,0}_{k+1,\delta})\big)\chi_1
\end{multline}
and the second term of \eqref{e:3l40} can be estimated by Lemma~\ref{l:xmlbdd}: for large $N>0$, 
\begin{multline}
    \|\dde^{-N} \tilde\chi [-\Delta, \Gamma^{s,l}_{k,\delta}\chi ]\chi_1\mathcal{R}_M^\Omega\chi_2 v\|\le C\|\Gamma^{s,0}_{k,\delta}\chi_1 \mathcal{R}_M^\Omega\chi_2 v\|_{X_{-N+l-1}^{l-1}}+C\|\Gamma^{s,0}_{k+1,\delta}\chi_1 \mathcal{R}_M^\Omega\chi_2 v\|_{X_{-N+l-\delta}^{l-1-\delta}}\\
    \le
    C\|\Gamma^{s,0}_{k,\delta}\chi_1 \mathcal{R}_M^\Omega \chi_2 v\|_{X_m^l}+C\|\Gamma^{s,0}_{k+1,\delta}\chi_1 \mathcal{R}_M^\Omega \chi_2 v\|_{X_m^l}.
\end{multline}
Eventually, we have the estimate of \eqref{e:3l39} by summing up those estimates above, and hence back in \eqref{e:3l8} we have
\begin{multline}\label{e:3l9}
    \|\tilde\chi \Gamma^{s,l}_{k,\delta}D_{x_0}\dde^{-N}\chi (\mathcal{H}_M^\Omega-E_{P})\tilde\chi f\|_{H^{-\frac12}(\partial\tilde\Omega)}\le C\|\dde^{-N} \tilde\chi \Gamma^{s,l}_{k,\delta}\chi v\|
    +C(1+M^2)\|\Gamma^{s,0}_{k, \delta}\chi \mathcal{R}_M^\Omega \chi_2 v\|_{X_{m}^l}\\
    +C\|\Gamma^{s,0}_{k,\delta}\chi_1\mathcal{R}_M^\Omega \chi_2 v\|_{X_m^l}
+C\|\Gamma^{s,0}_{k+1,\delta}\chi_1\mathcal{R}_M^\Omega \chi_2 v\|_{X_m^l}
\end{multline}

4. We now conclude the proof by estimating each term in \eqref{e:3l9}. Note per \eqref{e:3l11}, we have
\begin{equation}
    \|\dde^{-N} \tilde\chi \Gamma^{s,l}_{k,\delta}\chi v\|\le C\|\Gamma^{s,l}_{k,\delta}\chi v\|\le C\|\chi \Gamma^{s,l}_{k,\delta}v\|+C\|O(\Psi^{s,l-1+\delta}_{k-1,\delta})v\|\le C\|\Gamma^{-N, l-\frac{3}2\delta}_{-N,\delta}\tilde\chi f\|_{L^2(\partial\tilde\Omega)}.
\end{equation}
For the other terms, we have
\begin{equation}
    \|\Gamma^{s,0}_{k, \delta}\chi \mathcal{R}_M^\Omega \chi_2 v\|_{X_{m}^l}\le C\|\Gamma^{s,0}_{k, \delta}v\|_{Y_m^l}\le C\|\Gamma^{-N, \max(l,m+1)-\frac32\delta}_{-N,\delta} \tilde\chi f\|_{L^2(\partial\tilde\Omega)}\le C\|\Gamma^{-N, L}_{-N,\delta} \tilde\chi f\|_{L^2(\partial\tilde\Omega)}
\end{equation}
by \eqref{e:3l10} and that $\Gamma^{s,0}_{k,\delta}\chi\mathcal{R}^\Omega_M \chi_2\Gamma^{-s,0}_{-k,\delta}: Y_m^l(\tilde\Omega_+)\rightarrow \dot X_m^l(\tilde\Omega_+)$ is bounded as in Lemma \ref{l:resgammaBounds}. This is the place we used the assumption that $|m|<\rardyr_{\Omega}$ and $|m-l|<\frac12$. All terms of such shape can be estimated similarly. For any $s,k\in\mathbb{R}$, any $N>0$ large, we end up with 
\begin{equation}
    \|\tilde\chi \Gamma^{s,l}_{k,\delta}D_{x_0}\dde^{-N}\chi (\mathcal{H}_M^\Omega-E_{P})\tilde\chi f\|_{H^{-\frac12}(\partial\tilde\Omega)}\le C\|\Gamma^{-N, l-\frac{3}2\delta}_{-N,\delta}\tilde\chi f\|_{L^2(\partial\tilde\Omega)}+C\|\Gamma^{-N, L}_{-N,\delta}\tilde\chi f\|_{L^2(\partial\tilde\Omega)}\le C\|\Gamma^{-N, L}_{-N,\delta}\tilde\chi f\|_{L^2(\partial\tilde\Omega)},
\end{equation}
Eventually noting that $D_{x_0}\dde=0$ on $\{x_0=0\}$, we have
\begin{multline}
    \|\tilde\chi\Gamma^{s,l}_{k,\delta}\chi (\Lambda^\Omega_M-P)\tilde\chi f\|_{L^2(\partial\tilde\Omega)}\le C\|\tilde\chi\Gamma^{s+\frac12,l}_{k,\delta}D_{x_0}\dde^{-N}\chi(\mathcal{H}^\Omega_M-E_{P})\tilde\chi f\|_{H^{-\frac12}(\partial\tilde\Omega)}\\+C\|\Gamma^{-N, -N}_{-N,\delta}\tilde\chi f\|_{L^2(\partial\tilde\Omega)}\le C\|\Gamma^{-N, L}_{-N,\delta}\tilde\chi f\|_{L^2(\partial\tilde\Omega)}
\end{multline}
as desired.
\end{proof}
\begin{corollary}[Bounds for $\mathcal{H}_M^\Omega$]\label{c:harmonicBounds}
Let $\Omega\subset\mathbb{R}^d$ be an \acl domain with a regular collar $U$. Let $M_0>0$. Let $\delta\in [0,1)$. Let $\chi,\tilde\chi, \tilde{\tilde{\chi}}\in S^0(U)$ with $\supp\chi\cap\supp(1-\tilde\chi)=\supp\tilde\chi\cap\supp(1-\tilde{\tilde{\chi}})=\emptyset$. Let $l,m\in\mathbb{R}$ such that $|m|< \rardyr_{\Omega}$, $|l-m|<\frac12$. Let $s,k\in \mathbb{R}$. For each $N>0$, we have
\begin{multline}
    \|\tilde\chi \Gamma^{s,0}_{k,\delta}\chi\mathcal{H}^\Omega_M\tilde\chi f\|_{X_m^l(\Omega)}\le C\|\tilde{\tilde{\chi}}\Gamma^{s-\frac12,\max(l, m-\frac{1}2+0)}_{k,\delta}\tilde\chi f\|_{L^2(\partial\Omega)}+C\|\tilde{\tilde{\chi}}\Gamma^{s,\max(l,m)-\frac{\delta}{2}}_{k+\frac12,\delta} \tilde\chi f\|_{L^2(\partial\Omega)}\\
    +C\|\tilde{\tilde{\chi}}\Gamma^{-N, \max(l,m+1)-\frac32\delta}_{-N, \delta} \tilde\chi f\|_{L^2(\partial\Omega)}+C_N\|\tilde\chi f\|_{H^{-N,-N}(\partial\Omega)},
\end{multline}
uniformly for all $M\in[0,M_0]$ and $f\in \mc{S}'(\partial\Omega)$. 
\end{corollary}

\begin{proof}
    Note that Proposition~\ref{p:dtNParametrix} implies
    \begin{equation}\label{e:3l15} \|\Gamma^{s,0}_{k,\delta}\chi\mathcal{H}^\Omega_M\tilde\chi f\|_{X_m^l(\tilde\Omega_+)}\le \|\Gamma^{s,0}_{k,\delta}\chi E_{P}\tilde\chi f \|_{X_m^l(\tilde\Omega_+)}+C\|\Gamma^{-N, \max(l,m+1)-\frac32\delta}_{-N, \delta} \tilde\chi f\|_{L^2(\partial\tilde\Omega)}.
    \end{equation}
    We will estimate the $E_{P}\tilde\chi f$ term. Consider $D_{x_0}E_{P}=iPE_{P}$ with $P\in \Psi^{0,-\delta}_{1,\delta}$ and 
    \begin{multline}\label{e:3l41}
        \|\Gamma^{s,0}_{k,\delta}\chi E_{P}\tilde\chi f \|_{X_m^l(\tilde\Omega_+)}\le C\|\rre^{l}\dde^{m-l-1}\Gamma^{s,0}_{k,\delta}\chi E_{P}\tilde\chi f\|+C\|\rre^{l}\dde^{m-l}\nabla_{x}\Gamma^{s,0}_{k,\delta}\chi E_{P}\tilde\chi f\|
        \\
        +C\|\rre^{l}\dde^{m-l}\Gamma^{s,0}_{k,\delta}\chi PE_{P}\tilde\chi f\|+C\|\rre^{l}\dde^{m-l}\Gamma^{s,0}_{k,\delta}[D_{x_0},\chi] E_{P}\tilde\chi f\|.
    \end{multline}
    Since $\rre^{l}\dde^{m-l-\frac12+0}\le C\langle x\rangle^{\max(l, m-\frac{1}2+0)}$ on $U$, the first term of \eqref{e:3l41} can be estimated by
    \begin{multline}
        \|\rre^{l}\dde^{m-l-1}\Gamma^{s,0}_{k,\delta}\chi E_{P}\tilde\chi f\|\le C\|\tilde \chi\dde^{-\frac{1}2-0}O(\Psi^{s,\max(l, m-\frac{1}2+0)}_{k,\delta})\chi E_{P}\tilde\chi f\|+C\|\opsi\chi E_{P}\tilde\chi f\|\\
        \le C\|\Gamma^{s-\frac12,\max(l, m-\frac{1}2+0)}_{k,\delta}\tilde\chi f\|_{L^2(\partial\tilde\Omega)}
    \end{multline}
    where we used the boundedness of $E_{P}$ as in Lemma~\ref{l:dtnParametrixExtension}. Similarly we have
    \begin{multline}
        \|\rre^{l}\dde^{m-l}\nabla_{x}\Gamma^{s,0}_{k,\delta}\chi E_{P}\tilde\chi f\|
        +\|\rre^{l}\dde^{m-l}\Gamma^{s,0}_{k,\delta}\chi PE_{P}\tilde\chi f\|\le C\|\tilde\chi \Gamma^{s,\max(l,m)-\delta}_{k+1,\delta} E_P\tilde\chi f\|\\
        +C\|\opsi E_P \tilde\chi f\|\le C\|\Gamma^{s,\max(l,m)-\frac{\delta}{2}}_{k+\frac12,\delta} \tilde\chi f\|_{L^2(\partial\tilde\Omega)}, 
    \end{multline}
    and we have
    \begin{equation}
        \|\rre^{l}\dde^{m-l}\Gamma^{s,0}_{k,\delta}[D_{x_0},\chi]E_P\tilde\chi f\|
        \le C\|\Gamma^{s,\max(l,m)
        +\frac{\delta}{2}-1}_{k
        -\frac12,\delta} \tilde\chi f\|_{L^2(\partial\tilde\Omega)}.
    \end{equation}
    Thus we obtain
    \begin{equation}
        \|\Gamma^{s,0}_{k,\delta}\chi E_P\tilde\chi f \|_{X_{m}^{l}(\tilde\Omega_+)}\le C\|\Gamma^{s-\frac12,\max(l, m-\frac{1}2+0)}_{k,\delta}\tilde\chi f\|_{L^2(\partial\tilde\Omega)}+C\|\Gamma^{s,\max(l,m)-\frac{\delta}{2}}_{k+\frac12,\delta}\tilde\chi  f\|_{L^2(\partial\tilde\Omega)}.
    \end{equation}
    Revisit \eqref{e:3l15} to see
    \begin{multline}\label{e:3l42}
    \|\Gamma^{s,0}_{k,\delta}\chi\mathcal{H}^\Omega_M\tilde\chi f\|_{X_m^l(\tilde\Omega_+)}\le C\|\Gamma^{s-\frac12,\max(l, m-\frac{1}2+0)}_{k,\delta}\tilde\chi f\|_{L^2(\partial\tilde\Omega)}+C\|\Gamma^{s,\max(l,m)-\frac{\delta}{2}}_{k+\frac12,\delta} \tilde\chi f\|_{L^2(\partial\tilde\Omega)}\\
    +C\|\Gamma^{-N, \max(l,m+1)-\frac32\delta}_{-N, \delta} \tilde\chi f\|_{L^2(\partial\tilde\Omega)}.
    \end{multline}
    Localize with $\tilde\chi$ and $\tilde{\tilde{\chi}}$ to conclude the proof.
\end{proof}

\begin{remark}
All constants in this sections can depend on $\Omega$ in the following way: the upper bound on the choice of $r_0\ge 0$, and the lower bound on the radial Hardy constant $\rardyr_{\Omega}$, the upper bound on the diameter of $\partial\Omega\setminus\partial\tilde\Omega$ in Remark~\ref{r:aclcollar}, for each multi-index $\alpha$, the upper bounds on $C_\alpha$ of derivatives of the metric $\gplus$ in \eqref{e:3l26}, and lower bounds on $|\det \gplus|$.
\end{remark}

\section{Refined estimates near infinity on \acl domains}
\label{s:improvedHighFreq}

The goal of this section is to study the solution operators $\mathcal{R}_M^\Omega$ and $\mc{H}_M^{\Omega}$ in more detail. Specifically, we show that they are pseudolocal when acting on functions with non-trivial oscillation near infinity. We then use this to understand the structure of the DtN map on \acl domains-- showing that it is almost local at high frequencies and giving precise estimates on how varying the domain changes the DtN map. 

\noindent{\bf{Outline of the strategy}}

As we saw in Section~\ref{s:harm}, $\mc{R}_M^\Omega:Y_{m}^l \to X_m^l$ is well defined only for 
\begin{equation}
    \label{e:thresholdAgain}
    |m-l|<\frac{1}{2},\qquad |m|<\rardyr_{\Omega}.
\end{equation}  
The main idea in this section is to commute certain ($b$)- vector fields that are tangent to $\partial\Omega$ and have linearly growing coefficients through $\mc{R}_M^{\Omega}$ to improve these thresholds when acting on functions that oscillate. 

Roughly speaking, we will use the fact that 
\begin{equation}
    [-\Delta,\langle x\rangle]=\langle x\langle x\rangle^{-1},\nabla\rangle +O(\langle x\rangle^{-1})
\end{equation}
to see that 
\begin{equation}
    xD_x\mc{R}_M^\Omega= \mc{R}_M^\Omega xD_x + \text{bounded} 
\end{equation}
and hence that, if $f=D_x\tilde{f}+decaying$, then we can apply $\mc{R}_M^\Omega$ to $f$ even when $f$ lies in $Y_{m-1}^{l-1}$ with $(m,l)$ satisfying~\eqref{e:thresholdAgain}.  Similarly, if $f$ oscillates in the sense above and $f\in Y_{m+1}^{l+1}$, with $(m,l)$ satisfying~\eqref{e:thresholdAgain}, then $\mc{R}_M^{\Omega}f\in X_{m+1}^{l+1}$; i.e. the requirements on $(m,l)$ relax wehn acting on oscillating functions.

Functions that can be written as $f=D_x^N\tilde{f}_N+\text{decaying}$ have non-trivial oscillations at $|x|=\infty$ and hence by repeating the commutator argument, we can show that $\mc{R}_M^\Omega$ acts on functions with any amount of polynomial growth/decay, so long as these functions are oscillating. 

By further localizing this argument we, in fact, show that $\mc{R}_M^\Omega$ and $\mc{H}_M^\Omega$ are pseudolocal when acting on oscillating functions. We finally use this in Section~\ref{s:parametrix} to show that, on oscillatory functions and near infinity, $\Lambda_M^\Omega$ depends only on the local geometry of $\Omega$.

\begin{remark}
    Throughout this section, if not stated otherwise, we assume $r_0\ge 0$, $m,l\in\mathbb{R}$, and $\Omega\subset \mathbb{R}^d$ be any \acl domain with $M_0>0$.
\end{remark}

\subsection{Spaces with $b$-regularity}

In this section, we will study spaces on an \acl domain $\Omega\subset\mathbb{R}^d$. We first discuss the idea of $b$-regularity. Let $z$ be the Euclidean coordinates in $\Omega\subset\mathbb{R}^d$. We say $\chi\in S^l(\Omega)$ if $\chi\in C^\infty(\Omega)$ and $|\partial_{z}^\alpha \chi|\le C_\alpha\langle z\rangle^{l-|\alpha|}$. We say $\gamma\in S^l(\partial\Omega)$ if it extends to a function in $S^l(\Omega)$. Denote the order $(1,l)$ vector fields by
\begin{equation}
    \Diffhom^{1,l}(\Omega):=\Big\{ \sum_{|\alpha|= 1}a_\alpha(z)D_{z}^\alpha, \ a_\alpha\in S^{l}(\Omega)\Big\}, \ \Diff^{1,l}(\Omega):=\Big\{ \sum_{|\alpha|\le 1}a_\alpha(z)D_{z}^\alpha, \ a_\alpha\in S^{l-1+|\alpha|}(\Omega)\Big\}.
\end{equation}
In order to understand the boundary regularities, we define the order $(1,l)$ $b$-vector fields:
\begin{gather}
    \Diffhom_b^{1,l}(\Omega):=\Big\{\mc{V}\in \Diffhom^{1,l}(\Omega): \mc{V}|_{z}\in T_{z}\partial\Omega\text{ for all }z\in \partial\Omega\Big\},\\
    \Diff_b^{1,l}(\Omega):=\Big\{\mc{V}+f: f\in S^{l-1}(\Omega), \mc{V}\in \Diffhom_b^{1,l}(\Omega), \mc{V}|_{z}\in T_{z}\partial\Omega\text{ for all }z\in \partial\Omega \Big\}.
\end{gather}
Note by definition, $\mc{V}\in \Diffhom_b^{1,l}(\Omega)$ must be 0 on $\partial\osing$. Those $b$-vector fields are tangential to $\partial\Omega$. We define $\Diff^{1,l}_b(\partial\Omega)$ by the restriction of $\Diff^{1,l}_b(\Omega)$ onto $\partial\Omega$, and $\Diffhom^{1,l}_b(\partial\Omega)$ similarly. We then move on to define for $s\in\mathbb{N}$ that
\begin{gather}
    \Diff^{s,sl}(\Omega)=\operatorname{span}\Big\{ \mc{V}_1\cdots \mc{V}_k: k\le s, \mc{V}_1,\dots, \mc{V}_k\in \Diff^{1,l}(\Omega)\Big\}, \\
    \Diff^{s,sl}_b(\Omega)=\operatorname{span}\Big\{ \mc{V}_1\cdots \mc{V}_k: k\le s, \mc{V}_1,\dots, \mc{V}_k\in \Diff_b^{1,l}(\Omega)\Big\}. 
\end{gather}
We now define $X_m^l$-spaces with order $k$ $\Diff_b^{1,1}$-regularity and order $s$ $\Diff_b^{1,0}$-regularity. 
\begin{definition}[$X^{l,k}_{m,s}$]
For $k\in \mathbb{N}_0$, we define the $(m,l,k)$-spaces by
\begin{equation}
    X_{m}^{l,k}(\Omega)=\big\{u\in X_{m}^{l}(\Omega): \mc{V}_1\cdots\mc{V}_{k'} u\in X_{m}^{l}(\Omega), \mc{V}_j\in \Diffhom_b^{1,1}(\Omega), 1\le k'\le k\big\},
\end{equation}
\begin{equation}
    Y_{m}^{l,k}(\Omega)=\big\{u\in Y_{m}^{l}(\Omega): \mc{V}_1\cdots\mc{V}_{k'} u\in Y_{m}^{l}(\Omega), \mc{V}_j\in \Diffhom_b^{1,1}(\Omega), 1\le k'\le k\big\},
\end{equation}
\begin{equation}
    \dot X_{m}^{l,k}(\Omega)=\big\{u\in \dot X_{m}^{l}(\Omega): \mc{V}_1\cdots\mc{V}_{k'} u\in \dot X_{m}^{l}(\Omega), \mc{V}_j\in \Diffhom_b^{1,1}(\Omega), 1\le k'\le k\big\},
\end{equation}
\begin{equation}
    \dot{Y}_{m}^{l,k}(\Omega)=\big\{u\in \dot Y_{m}^{l}(\Omega): \mc{V}_1\cdots\mc{V}_{k'} u\in \dot Y_{m}^{l}(\Omega), \mc{V}_j\in \Diffhom_b^{1,1}(\Omega), 1\le k'\le k\big\}.
\end{equation}
To allow the negative indices, we let
\begin{equation}
    X_{m}^{l,-k}(\Omega)=(\dot Y_{-m}^{-l,k}(\Omega))^*, \ Y_{m}^{l,-k}(\Omega)=(\dot X_{-m}^{-l,k}(\Omega))^*,
\end{equation}
and $\dot X_{m}^{l,-k}$, $\dot Y_{m}^{l,-k}$ as the dual to $Y_{-m}^{-l,k}$, $X_{-m}^{-l,k}$ correspondingly. The above allows the definition of $X_{m}^{l,k}(\Omega)$, $Y_{m}^{l,k}(\Omega)$ for any $k\in\mathbb{Z}$. We use $X_{m}^{l,k}(\Omega)$ as the pivot space to define $(m,l,k,s)$-spaces. Define for $s\in\mathbb{N}_0$, 
\begin{equation}
    X_{m,s}^{l,k}(\Omega)=\big\{u\in X_{m}^{l,k}(\Omega): \mc{U}_1\cdots\mc{U}_{s'} u\in X_{m}^{l,k}(\Omega), \mc{U}_j\in \Diffhom_b^{1,0}(\Omega), 1\le s'\le s\big\},
\end{equation}
\begin{equation}
    Y_{m,s}^{l,k}(\Omega)=\big\{u\in Y_{m}^{l,k}(\Omega): \mc{U}_1\cdots\mc{U}_{s'} u\in Y_{m}^{l,k}(\Omega), \mc{U}_j\in \Diffhom_b^{1,0}(\Omega), 1\le s'\le s\big\},
\end{equation}
\begin{equation}
    \dot X_{m,s}^{l,k}(\Omega)=\big\{u\in \dot X_{m}^{l,k}(\Omega): \mc{U}_1\cdots\mc{U}_{s'} u\in \dot X_{m}^{l,k}(\Omega), \mc{U}_j\in \Diffhom_b^{1,0}(\Omega), 1\le s'\le s\big\},
\end{equation}
\begin{equation}
    \dot Y_{m,s}^{l,k}(\Omega)=\big\{u\in \dot Y_{m}^{l,k}(\Omega): \mc{U}_1\cdots\mc{U}_{s'} u\in \dot Y_{m}^{l,k}(\Omega), \mc{U}_j\in \Diffhom_b^{1,0}(\Omega), 1\le s'\le s\big\}.
\end{equation}
\begin{equation}
    X_{m,-s}^{l,k}(\Omega)=(\dot Y_{-m,s}^{-l,-k}(\Omega))^*, \ Y_{m,-s}^{l,k}(\Omega)=(\dot X_{-m,s}^{-l,-k}(\Omega))^*.
\end{equation}
The above defines $X_{m,s}^{l,k}(\Omega)$ for all $k,s\in\mathbb{Z}$.
\end{definition}
We identify a basis $\mc{W}_0, \mc{W}_1,\dots, \mc{W}_{d-1}\in \Diffhom_b^{1,0}(\Omega)$, that $\{\mc{W}_j|_{z}\}_{j=0}^{d-1}$ spans $T_z\Omega$ for each $z$ in the interior of $\Omega$. We further ask, for each $z\in \partial\oreg$, locally near $z$ inside $\overline{\Omega}$ in the Fermi normal coordinates $(x_0,x)$ where $x_0$ is the variable normal to $\partial\Omega$, we have
\begin{equation}
    \mc{W}_0=a_0(x_0,x) x_0 D_{x_0}, \ \mc{W}_j=a_j(x_0,x) D_{x_j}\text{ for }1\le j\le d-1,
\end{equation}
for some $a_0, a_1, \dots, a_{d-1}$ locally $C^\infty$ near $z$ inside $\overline{\Omega}$. 
Given multi-index $\alpha=(\alpha_0, \dots, \alpha_{d-1})$ with $\alpha_j\in\mathbb{N}_0$, we define the multi-index operator $\mc{W}^\alpha$ by
\begin{equation}
    \mc{W}^\alpha=\mc{W}_0^{\alpha_0}\mc{W}_1^{\alpha_1}\cdots \mc{W}_{d-1}^{\alpha_{d-1}}.
\end{equation}
Given a word $\alpha=(\alpha_1, \dots, \alpha_{k})$ of length $k$, with $\alpha_j\in \{0,\dots, d-1\}$ for each $1\le j\le k$, we define the operator word $\mc{W}_\alpha$ by
\begin{equation}
    \mc{W}_\alpha=\mc{W}_{\alpha_1}\cdots \mc{W}_{\alpha_{k}}.
\end{equation}
In this section, $\gamma$ cutoffs live on $\partial\Omega$ and $\chi$ cutoffs live in $\Omega$.

\subsection{Regularity of the resolvent near infinity}

Lemma~\ref{l:reg2} allows us to obtain higher regularity inside a compact set. Our next estimates provide increased regularity near infinity. Before, proceeding, we record two lemmas on the interaction between commutators with the Laplacian and various vector fields.

We define the remainder class of divergence form operators, 
\begin{equation}
    \oxy=O(\div_{z}B_1\nabla_{z}+B_2\nabla_{z}+B_3), \ (B_1)_{ij}\in S^0(\Omega), (B_2)_j\in S^{-1}(\Omega), B_3\in S^{-2}(\Omega).
\end{equation}
We show any such remainder in $\oxy$ is bounded from $X_m^l(\Omega)$ to $Y_{m}^l(\Omega)$, where we remind the reader that $Y_m^l=(\dot X_{-m}^{-l})^*$.  
\begin{lemma}
\label{l:starEstimate}
Let $\Omega$ be a uniformly Lipschitz domain. Let $B_1=\sum_{i,j=1}^d \partial_{z_i}b_{1,i}^j(z)\partial_{z_j}$ and $B_2=\sum_{j=1}^d b_{2}^j(z)\partial_{z_j}$, $B_3=b_{3}(z)$,
where the $L^\infty(\Omega)$-coefficients satisfy
\begin{equation}
    |b_{1,i}^j|\leq  C,\ | b_{2}^j|\leq C\dd^{-1},\ | b_{3}|\leq C\dd^{-2}. 
\end{equation}
Then, there is $C>0$ that uniformly for all $u\in X_m^l(\Omega)$ we have
\begin{equation}
    \|B_1 u\|_{Y_m^l(\Omega)}+\|B_2u\|_{Y_m^l(\Omega)}+\|B_3u\|_{Y_m^l(\Omega)}\leq C\| u\|_{X_m^l(\Omega)}.
\end{equation}
\end{lemma}
\begin{proof}
Note for $\phi\in \dot X_{-m}^{-l}$ we have
\begin{equation}
    |\langle b_3(z)u,\phi\rangle|\le C\|\rr^l\dd^{m-l-1}u\|\|\rr^{-l}\dd^{-m+l-1}\phi\|\le C\|u\|_{X_m^l}\|\phi\|_{X_{-m}^{-l}}.
\end{equation}
Observe that if $f=v+\mathrm{div}\,w$, we have
\begin{equation}
    \left\Vert f\right\Vert _{Y_m^l}\leq\Vert  \rr^{l}\dd^{m-l+1}v\Vert _{L^{2}}+\Vert \rr^l\dd^{m-l} w\Vert _{L^{2}}.
\end{equation}
Therefore, 
\begin{equation}
    \|B_2 u\|_{Y_m^l}\leq \Vert  \rr^{l}\dd^{m-l+1}B_2u\Vert _{L^{2}}\leq C\Vert \rr^l\dd^{m-l}\nabla u\Vert_{L^2}\leq C\|u\|_{X_m^l},
\end{equation}
and
\begin{equation}
    \|B_1u\|_{Y_m^l}\leq \sum_{i,j=1}^d\Vert \rr^l\dd^{m-l}b_{1,i}^j\partial_{z_j}u\Vert_{L^{2}}\leq C\Vert \rr^l\dd^{m-l}\nabla u\Vert_{L^2}\leq C\|u\|_{X_m^l},
\end{equation}
 as desired. 
\end{proof}

We show commutators between $-\Delta$, $\Diff^{1,1}(\Omega)$ are remainders in $\oxy$. 
\begin{lemma}
\label{l:algCommute}
Let $V\in \Diffhom^{1,1}(\Omega)$, $\varphi\in S^0(\Omega)$. Separately we have
\begin{equation}
\label{e:1commute}
[-\Delta, V]=\oxy, \ [-\Delta, \varphi]=\oxy, \ [-\Delta, V\varphi]=\oxy, \ [-\Delta, \varphi]V=\oxy,
\end{equation}
and 
\begin{equation}
\label{e:2commute}
[V,\oxy]=\oxy, \ [\varphi, \oxy]=\oxy.
\end{equation}
\end{lemma}
\begin{proof}
Given $V\in \Diffhom^{1,1}(\Omega)$, we have
\begin{equation}
    V=\sum_{k=1}^{d}v^k(z)\partial_{z_k}, \ v^k\in S^1(\Omega)\text{ for }1\le k \le d.
\end{equation}
To obtain \eqref{e:1commute}, we compute
\begin{equation}
[\partial_{z_i}\partial_{z_j}, v^k\partial_{z_k}]=\partial_{z_i}[\partial_{z_j}, v^k]\partial_{z_k}+\partial_{z_j}[\partial_{z_i}, v^k]\partial_{z_k}+[[\partial_{z_i}, v^k]\partial_{z_k},\partial_{z_j}],
\end{equation}
which is in $\div_{z}(S^0(\Omega))_{d\times d}\nabla_{z}+(S^{-1}(\Omega))_{d}\nabla_{z}$. Summing up in $i,j$ to observe $[-\Delta, V]=\oxy$. The rest is similar. 
\end{proof}

We show any such remainder in $\oxy$ is bounded from $X_{m,s}^{l,k}(\Omega)$ to $Y_{m,s}^{l,k}(\Omega)$ for all $k,s\in\mathbb{Z}$. 
\begin{lemma}\label{l:oxybddxmlk}
Let $l,m\in\mathbb{R}$, $k,s\in \mathbb{Z}$. Let $\varphi\in S^0(\Omega)$. Then we have uniformly for all $u\in X_{m,s}^{l,k}(\Omega)$, 
\begin{gather}
    \|\oxy u\|_{Y_{m,s}^{l,k}(\Omega)}\le C\| u\|_{X_{m,s}^{l,k}(\Omega)},\\
    \|[-\Delta, \varphi] u\|_{Y_{m,s+1}^{l+1,k}(\Omega)}\le C\| u\|_{X_{m,s}^{l,k}(\Omega)}.
\end{gather}
\end{lemma}

\begin{proof}
    1. Let $s=0$. We claim for $k\ge 0$ that
    \begin{equation}\label{e:4l18}
        \|\mc{V}_1\cdots\mc{V}_k\oxy u\|_{Y_m^l}\le C\sum_{|\alpha|\le k}\|\mc{V}_\alpha u\|_{X_m^l}.
    \end{equation}
    The case $k=0$ follows from Lemma~\ref{l:starEstimate}. Assume \eqref{e:4l18} holds for cases up to $k-1$, we have
    \begin{equation}
        \|\mc{V}_1\cdots\mc{V}_k\oxy u\|_{Y_m^l}\le\|\oxy \mc{V}_1\cdots\mc{V}_ku\|_{Y_m^l}+\|[\mc{V}_1\cdots\mc{V}_k,\oxy]u\|_{Y_m^l}\le C\sum_{|\alpha|\le k}\|\mc{V}_\alpha u\|_{X_m^l},
    \end{equation}
    as desired in \eqref{e:4l18}. Note here we used the iterative assumption and Lemma~\ref{l:algCommute} for
    \begin{equation}
        [\mc{V}_1\cdots\mc{V}_k,\oxy]=\sum_{|\gamma|< k}(\ad^{\gamma^c}_{\mc{V}}(\oxy))\mc{V}_{\gamma}= \sum_{|\gamma|< k}\oxy\mc{V}_{\gamma},
    \end{equation}
    as a summation over all proper subwords $\gamma$ of $(1,\dots,k)$ where $\gamma^c=(1,\dots,k)\setminus \gamma$. Here we used the notation
    \begin{equation}\label{e:adword}
        \ad_{\mc{V}}^\gamma (\cdot)=[\mc{V}_{\gamma_1},[\mc{V}_{\gamma_2}, \dots, [\mc{V}_{\gamma_j},\cdot]]], \ \gamma=(\gamma_1, \dots, \gamma_j). 
    \end{equation}
    The $k<0$ cases follow from the duality and $(\oxy)^*=\oxy$. 

    2. Now let $s>0$. Assume the claim holds for cases up to $s-1$ and all $k$. We have
    \begin{equation}
        \|\mc{U}_1\cdots\mc{U}_s\oxy u\|_{Y_m^{l,k}}\le\|\oxy \mc{U}_1\cdots\mc{U}_su\|_{Y_m^{l,k}}+\|[\mc{U}_1\cdots\mc{U}_s,\oxy]u\|_{Y_m^{l,k}}\le C\sum_{|\alpha|\le s}\|\mc{U}_\alpha u\|_{X_m^{l,k}},
    \end{equation}
    similarly as in Step 1. We here used Lemma~\ref{l:algCommute} and $\Diffhom^{1,0}_b\subset\Diffhom^{1,1}_b$. The duality gives the cases for $s<0$ and all $k$. 

    3. Note for any $m,l\in\mathbb{R}$, we have $[-\Delta,\varphi]: X_m^l(\Omega)\rightarrow \dot Y_{m,1}^{l+1,0}(\Omega)$. Indeed, note
\begin{equation}
    [-\Delta,\varphi]=\rr^{-1}O(S^0)\nabla+\rr^{-2}O(S^0). 
\end{equation}
For each $u\in X_m^l(\Omega)$ and $\phi\in X^{-l-1}_{-m}(\Omega)$, we have
\begin{gather}
    |\langle \rr^{-1}O(S^0)\nabla u, \phi\rangle|\le C\|\rr^{l}\dd^{-l+m}\nabla u\|\|\rr^{-l-1}\dd^{l-m}\phi\|\le C\|u\|_{X_m^l}\|\phi\|_{X_{-m}^{-l-1}},\\
    |\langle \mc{U}_1\rr^{-1}O(S^0)\nabla u, \phi\rangle|\le C\|\rr^{l}\dd^{-l+m}\nabla u\|\|\rr^{-l-1}\dd^{l-m}\mc{U}_1^*\phi\|\le C\|u\|_{X_m^l}\|\phi\|_{X_{-m-1}^{-l-1}}
\end{gather}
showing $\rr^{-1}O(S^0)\nabla: X_m^l(\Omega)\rightarrow \dot Y_{m,1}^{l+1,0}(\Omega)$. Similarly one can show $\rr^{-2}O(S^0): X_m^l(\Omega)\rightarrow \dot Y_{m,1}^{l+2,0}(\Omega)$. This concludes the claim. Now address the $k$, $s$ indices similarly like in Step 1 and 2. 
\end{proof}

We can now improve the mapping properties of $\mc{R}_M^\Omega$ on \acl domains.
\begin{lemma}[$\mc{R}_M^\Omega$ conserves positive $b$-regularity]
\label{l:longDistanceResolve}
Let $k,s\in\mathbb{N}_0$. Let $|m|<\rardyr_\Omega$, $|m-l|<\frac12$. Let vector fields $\mc{V}_1,\dots \mc{V}_k\in \Diffhom_b^{1,1}(\Omega)$, $\mc{U}_1,\dots \mc{U}_s\in \Diffhom_b^{1,0}(\Omega)$. Consider $\tilde\chi \in S^0(\Omega)$ with $\supp\tilde\chi\cap \partial\Omega$ is smooth, and $\tilde\chi\equiv 1$ near  $(\cap_{k'=1}^k\supp\mc{V}_{k'})\cap(\cap_{s'=1}^s\supp\mc{U}_{s'})$. Then we have
\begin{equation}
    \|\mc{U}_1\dots \mc{U}_s\mc{V}_1\dots \mc{V}_k\mathcal{R}_{M}^\Omega f\|_{X_m^l(\Omega)}\leq C\|f\|_{Y_m^l(\Omega)}+C\sum_{|\alpha|\leq k,|\beta|\le s}\| \tilde{\chi}\mc{W}^{\beta}(\rr\mc{W})^\alpha f\|_{Y_m^l(\Omega)},
\end{equation}
and
\begin{equation}
    \| (-\Delta+M^2)\mc{U}_1\dots \mc{U}_s\mc{V}_1\dots \mc{V}_k\mathcal{R}_{M}^\Omega f\|_{Y_m^l(\Omega)}\leq C\|f\|_{Y_m^l(\Omega)}+C\sum_{|\alpha|\leq k,|\beta|\le s}\| \tilde{\chi}\mc{W}^{\beta}(\rr\mc{W})^\alpha f\|_{Y_m^l(\Omega)},
\end{equation}
uniformly for all $f\in Y_{m,s}^{l,k}(\Omega)$ and $M\in[0,M_0]$. This specifically implies
\begin{equation}
    \mc{R}^\Omega_M: Y_{m,s}^{l,k}(\Omega)\rightarrow \dot X_{m,s}^{l,k}(\Omega)
\end{equation}
is uniformly bounded in $M\in[0,M_0]$. 
\end{lemma}
\begin{proof}
1. Let $\chi\in S^0(\Omega)$ satisfying all requirements for $\tilde\chi$ and $\supp\chi\cap\supp(1-\tilde\chi)=\emptyset$. We start by showing for each $k,s\ge 0$, we have the following decomposition:
\begin{equation}
\label{e:alg}
(-\Delta+M^{2}) \chi \mc{U}_1\cdots \mc{U}_s\mc{V}_1\cdots \mc{V}_k\mathcal{R}_{M}^\Omega=\chi\mc{U}_1\dots \mc{U}_s\mc{V}_1\cdots \mc{V}_k+\sum_{|\alpha|+|\beta|< k+s}\oxy\tilde\chi \mc{U}^\beta\mc{V}^\alpha\mc{R}^\Omega_M,
\end{equation}
as a summation over all possibly improper subwords $\alpha$ of $(1,\dots, k)$ and $\beta$ of $(1,\dots, s)$. This is done by noting
\begin{equation}
    [-\Delta,\chi\mc{U}_1\cdots \mc{U}_s\mc{V}_1\cdots \mc{V}_k]=-\sum_{|\alpha|+|\beta|<k+s} \chi(\ad_\mc{U}^{\beta^c}\ad_\mc{V}^{\alpha^c}(-\Delta))\tilde\chi\mc{U}^\beta\mc{V}^\alpha,
\end{equation}
as a summation in subwords described above. Note $\ad_\mc{U}^{\beta^c}\ad_\mc{V}^{\alpha^c}(-\Delta)=\oxy$ per Lemma~\ref{l:algCommute} to conclude \eqref{e:alg}.  

2. We claim for each $(k,s)\in \mathbb{N}_0\times \mathbb{N}_0$, we have for all $f\in Y_{m,s}^{l,k}$, 
\begin{equation}\label{e:4l10}
    \| (-\Delta+M^2)\chi \mc{U}_1\dots \mc{U}_s\mc{V}_1\dots \mc{V}_k\mathcal{R}_{M}^\Omega f\|_{Y_m^l(\Omega)}\leq  C\|f\|_{Y_m^l(\Omega)}+C\sum_{|\alpha|\leq k,|\beta|\le s}\| \tilde{\chi}\mc{W}^{\beta}(\rr\mc{W})^\alpha f\|_{Y_m^l(\Omega)},
\end{equation}
and 
\begin{equation}\label{e:4l21}
    \| \chi \mc{U}_1\dots \mc{U}_s\mc{V}_1\dots \mc{V}_k\mathcal{R}_{M}^\Omega f\|_{X_m^l(\Omega)}\leq C\|f\|_{Y_m^l(\Omega)}+C\sum_{|\alpha|\leq k,|\beta|\le s}\| \tilde{\chi}\mc{W}^{\beta}(\rr\mc{W})^\alpha f\|_{Y_m^l(\Omega)},
\end{equation}
and we prove it iteratively. Note $-[\Delta, \chi]=\oxy$ and we have
\begin{equation}
    (-\Delta+M^{2}) \chi\mathcal{R}_{M}^\Omega f=\chi f-[\Delta, \chi]\mathcal{R}_{M}^\Omega f\in Y_m^l. 
\end{equation}
The case $(k,s)=(0,0)$ here follows from Lemma~\ref{l:starEstimate}, the boundedness of $\mc{R}_M^\Omega: Y_m^l\rightarrow \dot X_m^l$ per Proposition~\ref{p:extendNoMass}, and that $\chi\mc{R}_M^\Omega f$ has zero trace over the boundary. Now assume \eqref{e:4l10} and \eqref{e:4l21} hold for all cases $(k,s)$ where $k+s<N$. Consider some case where $k+s=N$. Note $f\in Y_{m,s}^{l,k}$, and \eqref{e:alg} implies
\begin{equation}\label{e:4l19}
    (-\Delta+M^2)\chi\mc{U}_1\cdots\mc{U}_s\mc{V}_1\cdots \mc{V}_k\mathcal{R}_{M}^\Omega f
=\chi\mc{U}_1\cdots\mc{U}_s\mc{V}_1\cdots \mc{V}_kf+\sum_{|\alpha|+|\beta|< N}\oxy\tilde\chi \mc{U}^\beta\mc{V}^\alpha\mc{R}^\Omega_M f,
\end{equation}
as a summation over all possibly improper subwords $\alpha$ of $(1,\dots, k)$ and $\beta$ of $(1,\dots, s)$. The iterative assumption \eqref{e:4l21} implies 
\begin{equation}
    \|\sum_{|\alpha|+|\beta|<N}\oxy\tilde\chi \mc{U}^\beta\mc{V}^\alpha\mc{R}^\Omega_M f\|_{Y_m^l}\le C\|f\|_{Y_m^l}+C\sum_{|\alpha|+|\beta|< N}\| \tilde{\chi}\mc{W}^\beta(\rr\mc{W})^\alpha f\|_{Y_m^l}
\end{equation}
over subwords $\alpha,\beta$ specified above. One can estimate the principal term by
\begin{equation}
    \|\chi\mc{U}_1\cdots\mc{U}_s\mc{V}_1\cdots \mc{V}_kf\|_{Y_m^l}\le C\sum_{|\alpha|\leq k,|\beta|\le s}\| \tilde{\chi}\mc{W}^{\beta}(\rr\mc{W})^\alpha f\|_{Y_m^l}. 
\end{equation}
Those two estimates justify \eqref{e:4l10} for $k+s=N$ and the zero trace of $\chi \mc{U}_1\dots \mc{U}_s\mc{V}_1\dots \mc{V}_k\mathcal{R}_{M}^\Omega f$ then gives \eqref{e:4l21}. Iteration concludes the proof.
\end{proof}

Using self-adjointness of $\mc{R}_M^\Omega$, we are able to work in spaces with very low $b$-regularity.
\begin{lemma}[$\mc{R}_M^\Omega$ conserves $b$-regularity]
\label{l:lowRegularityAction}
Let $|m|<\rardyr_\Omega$, $|m-l|<\frac12$. For all $k,s\in\mathbb{Z}$, we have
\begin{equation}
    \mathcal{R}_M^\Omega : Y_{m,s}^{l,k}(\Omega) \to \dot X_{m,s}^{l,k}(\Omega)
\end{equation}
is uniformly bounded for $M\in[0,M_0]$. 
\end{lemma}
\begin{proof}
Let $k\le 0,s=0$. By Lemma~\ref{l:longDistanceResolve}, $\mc{R}_M^\Omega:Y_{-m}^{-l,-k}\to \dot X_{-m}^{-l,-k}$ is bounded. Note for $f\in Y^{l,k}_{m}=(\dot X_{-m}^{-l,-k})^*$, we have
\begin{equation}
    \langle (-\Delta+M^2)(\mathcal{R}_M^\Omega)^*f, \phi\rangle= \langle f, \mathcal{R}_M^\Omega(-\Delta+M^2)\phi\rangle=\langle f, \phi\rangle
\end{equation}
for all $\phi\in \dot X_{-m}^{-l,-k}$, where we used that $\mathcal{R}_M^\Omega(-\Delta+M^2)$ is identity on $\dot X^{-l}_{-m}$. This implies $(\mathcal{R}_M^\Omega)^*f=\mathcal{R}_M^\Omega f$ for $f\in Y_{m}^{l,k}$ and thus $\mathcal{R}_M^\Omega =(\mathcal{R}_M^\Omega)^*: Y_{m}^{l,k}\to \dot X_{m}^{l,k}$ is bounded. The inclusion of negative $s$ is similar. 
\end{proof}

Finally, when the data is away from the measurement, we obtain a smoothing operator.
\begin{lemma}[Nonlocal effect of $\mc{R}_M^\Omega$ on smooth part is $b$-smoothing]
\label{l:smoothing}
Let $|m|<\rardyr_\Omega$, $|m-l|<\frac12$. Let $\chi_1,\chi_2 \in S^0(\Omega)$ with disjoint support. Then for all $N\in\mathbb{N}_0$,
\begin{equation}
    \chi_1\mc{R}_M^\Omega\chi_2: Y_{m,-N}^{l,-N}(\Omega)\to \dot X_{m,N}^{l,N}(\Omega)
\end{equation}
is uniformly bounded for $M\in[0,M_0]$. 
\end{lemma}
\begin{proof}
By Lemma~\ref{l:longDistanceResolve}, for any such $\chi_1$, $\chi_2$ with disjoint support, we have $\chi_1\mc{R}_M^\Omega\chi_2:Y_m^l\to \dot X_{m,N}^{l,N}$. Since $(\mathcal{R}_M^\Omega)^*=\mathcal{R}_M^\Omega$ on $Y^{-l,-N}_{-m,-N}$, we also have $\chi_1\mc{R}_M^\Omega\chi_2: Y_{-m,-N}^{-l,-N}\to \dot X_{-m}^{-l}$. Now, let $\tilde{\chi}_1\in S^0(\Omega)$, $\supp \tilde{\chi}_1\cap \supp \chi_2=\emptyset$, and $\supp\chi_1\cap\supp(1-\tilde\chi_1)=\emptyset$. Then we have
\begin{equation}
\tilde{\chi}_1\mathcal{R}_M^\Omega=\mathcal{R}_M^\Omega\tilde{\chi}_1+\mathcal{R}_M^\Omega[-\Delta,\tilde{\chi}_1]\mathcal{R}_M^\Omega.
\end{equation}
Letting $\varphi\in S^0(\Omega)$, $\supp \varphi$ disjoint from the support of $\chi_1$ and $\chi_2$, and $\varphi\nabla\tilde{\chi}_1= \nabla\tilde{\chi}_1$, we have
\begin{equation}
\chi_1\mc{R}_M^\Omega\chi_2=\chi_1\mc{R}_M^\Omega\varphi[-\Delta,\tilde{\chi}_1]\varphi\mc{R}_M^\Omega\chi_2.
\end{equation}
It suffices to note
\begin{equation}
    \varphi\mc{R}_M^\Omega \chi_2:Y_{m,-N}^{l,-N}\to \dot X_{m}^{l}, \ [-\Delta,\tilde{\chi}_1]:\dot X_{m}^{l}\to Y_m^l, \ \chi_1\mc{R}_M^\Omega\varphi:Y_m^l\to \dot X_{m,N}^{l,N}
\end{equation}
are bounded to conclude the proof. 
\end{proof}

\begin{corollary}\label{c:RimproveM}
Let $|m|<\rardyr_\Omega$, $|m-l|<\frac12$. For all $m_2\in\mathbb{R}$, $s\in\mathbb{Z}$, we have
\begin{gather}\label{e:4l40}
    \langle M\dd\rangle^{m_2}\mathcal{R}_M^\Omega \langle M\dd\rangle^{-m_2}: Y_{m,s}^{l,0}(\Omega) \to \dot X_{m,s}^{l,0}(\Omega)
\end{gather}
is uniformly bounded for $M\in[0,M_0]$. Let $\chi_1,\chi_2 \in S^0(\Omega)$ with disjoint support. Then for all $N\in\mathbb{N}_0$,
\begin{equation}\label{e:4l44}
    \chi_1\langle M\dd\rangle^{m_2}\mc{R}_M^\Omega\langle M\dd\rangle^{-m_2}\chi_2: Y_{m,-N}^{l,0}(\Omega)\to \dot X_{m,N}^{l,0}(\Omega)
\end{equation}
is uniformly bounded for $M\in[0,M_0]$. 
\end{corollary}
\begin{proof}
    1. We briefly show \eqref{e:4l40}. We follow the strategy of Lemma~\ref{l:longDistanceResolve} and show the case when $s=1$. Given $f\in Y_{m,1}^{l,0}(\Omega)$ and $\mc{U}_1\in \Diffhom_b^{1,0}(\Omega)$, consider
\begin{equation}
    (-\Delta+M^2)\mc{U}_1\langle M\dd\rangle^{m_2}\mc{R}_M^\Omega \langle M\dd\rangle^{-m_2}f=\mc{U}_1 f+[-\Delta, \mc{U}_1\langle M\dd\rangle^{m_2}]\langle M\dd\rangle^{-m_2}\langle M\dd\rangle^{m_2}\mc{R}_M^\Omega \langle M\dd\rangle^{-m_2}f.
\end{equation}
Note $\langle M\dd\rangle^{m_2}\mc{R}_M^\Omega \langle M\dd\rangle^{-m_2}: Y_m^l\rightarrow \dot X_m^l$ is bounded per Proposition~\ref{p:extendNoMass} and
\begin{equation}
    [-\Delta, \mc{U}_1\langle M\dd\rangle^{m_2}]\langle M\dd\rangle^{-m_2}=O(\dd^{-1})\nabla^2+O(\dd^{-2})\nabla+O(\dd^{-3}): X_m^l\rightarrow Y_m^l
\end{equation}
is also bounded. This implies
\begin{equation}
    (-\Delta+M^2)\mc{U}_1\langle M\dd\rangle^{m_2}\mc{R}_M^\Omega \langle M\dd\rangle^{-m_2}f\in Y_m^l(\Omega), \ \mc{U}_1\langle M\dd\rangle^{m_2}\mc{R}_M^\Omega \langle M\dd\rangle^{-m_2}f\in X_m^l(\Omega)
\end{equation}
as desired. The cases $s\ge0$ follow from induction, and the cases $s\le 0$ follow from taking the adjoint.

2. The proof of \eqref{e:4l44} is similar to Lemma~\ref{l:smoothing}. 
\end{proof}

\subsection{Regularity of the harmonic extension near infinity}
In this subsection, we use the results obtained early to obtain improved the regularity of $\mc{H}_M^\Omega$ near infinity. 

We show that $\mc{H}_M^\Omega$ maps functions on $\partial\Omega$ with $k$ $\Diff^{1,1}_b$-regularity to $X^{l,k}_m(\Omega)$ for $k\in \mathbb{N}_0$. 

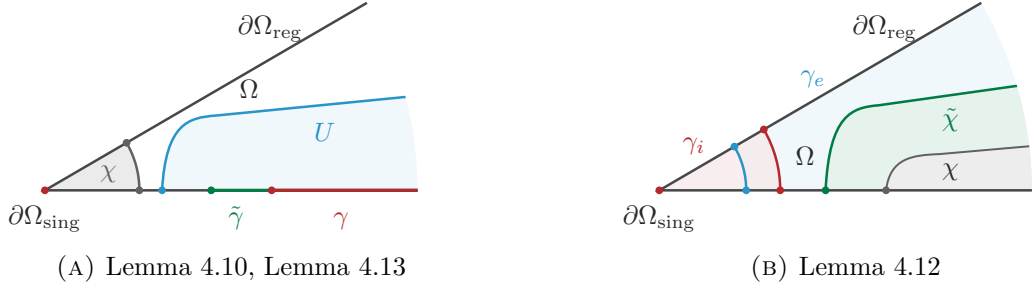
\begin{figure}[t]
  \centering
  \begin{subfigure}[t]{0.5\textwidth}
        \centering
        \begin{tikzpicture}[x=1cm,y=1cm]
  \definecolor{colorAKANE}{RGB}{183,40,46}
  \definecolor{colorTOKIWA}{RGB}{0,123,67}
  \definecolor{colorUSUAI}{RGB}{40,148,200}
  \colorlet{plotPrimary}{colorAKANE}
  \colorlet{plotSecondary}{colorTOKIWA}
  \colorlet{plotTertiary}{colorUSUAI}
  \colorlet{plotAxis}{black!72}

  \tikzset{
    >=latex, font=\small, mark size=1.5pt, mark options={thin},
    line cap=butt, line join=round,
    wedge edge/.style={draw=plotAxis,line width=0.95pt},
    collar edge/.style={draw=plotTertiary,line width=0.95pt},
    outer support/.style={draw=plotSecondary,line width=1.2pt},
    inner support/.style={draw=plotPrimary,line width=1.2pt,solid},
    chi edge/.style={draw=black!62,line width=0.95pt},
    figure label/.style={text=black!84,inner sep=0.8pt}
  }

  \pgfmathsetmacro{\figureWidth}{7.5}
  \pgfmathsetmacro{\figureHeight}{3.2}
  \pgfmathsetmacro{\frameBottom}{-0.28}
  \pgfmathsetmacro{\frameTop}{\frameBottom+\figureHeight}
  \pgfmathsetmacro{\edgeDx}{4.26}
  \pgfmathsetmacro{\edgeDy}{2.48}
  \pgfmathsetmacro{\edgeLength}{veclen(\edgeDx,\edgeDy)}
  \pgfmathsetmacro{\tipX}{(\figureWidth-\edgeLength)/2}
  \pgfmathsetmacro{\tipY}{0.34}
  \pgfmathsetmacro{\wedgeAngle}{atan2(\edgeDy,\edgeDx)}

  \path[use as bounding box]
    (0,\frameBottom) rectangle (\figureWidth,\frameTop);
  \fill[fill=white]
    (0,\frameBottom) rectangle (\figureWidth,\frameTop);
  \coordinate (vertex) at (\tipX,\tipY);
  \coordinate (upperEnd) at ({\tipX+\edgeDx},{\tipY+\edgeDy});
  \coordinate (lowerEnd) at ({\tipX+\edgeLength},\tipY);
  \coordinate (cropBeyond) at ({\tipX+\edgeLength+0.5},\tipY);

  \coordinate (collarStart) at ({\tipX+1.55},\tipY);
  \coordinate (collarControlOne) at ({\tipX+1.60},{\tipY+0.884});
  \coordinate (collarControlTwo) at ({\tipX+1.90},{\tipY+0.962});
  \coordinate (collarShoulder) at ({\tipX+2.25},{\tipY+0.988});
  \coordinate (collarEnd) at
    ({\tipX+\edgeLength+0.5},{\tipY+1.30});
  \coordinate (tildeStart) at ({\tipX+2.20},\tipY);
  \coordinate (gammaStart) at ({\tipX+3.00},\tipY);

  \def\chiRadius{1.25}
  \path (vertex) ++(0:\chiRadius) coordinate (chiLower);
  \path (vertex) ++(\wedgeAngle:\chiRadius) coordinate (chiUpper);

  \begin{scope}
    \clip (vertex) -- (lowerEnd)
      arc[start angle=0,end angle=\wedgeAngle,radius=\edgeLength cm]
      -- cycle;
    \fill[fill=plotTertiary!7]
      (collarStart) .. controls (collarControlOne) and (collarControlTwo)
      .. (collarShoulder) -- (collarEnd) -- (cropBeyond) -- cycle;
    \draw[collar edge]
      (collarStart) .. controls (collarControlOne) and (collarControlTwo)
      .. (collarShoulder) -- (collarEnd);

    \fill[fill=black!8]
      (vertex) -- (chiLower)
      arc[start angle=0,end angle=\wedgeAngle,radius=\chiRadius cm]
      -- cycle;
    \draw[chi edge] (chiLower)
      arc[start angle=0,end angle=\wedgeAngle,radius=\chiRadius cm];
  \end{scope}

  \draw[wedge edge] (upperEnd) -- (vertex) -- (lowerEnd);
  \draw[outer support] (tildeStart) -- (lowerEnd);
  \draw[inner support] (gammaStart) -- (lowerEnd);
  \fill[fill=plotPrimary] (vertex) circle[radius=1.45pt];
  \fill[fill=plotTertiary] (collarStart) circle[radius=1.45pt];
  \fill[fill=plotSecondary] (tildeStart) circle[radius=1.45pt];
  \fill[fill=plotPrimary] (gammaStart) circle[radius=1.45pt];
  \foreach \point in {chiLower,chiUpper}
    \fill[fill=black!62] (\point) circle[radius=1.45pt];

  \path (vertex) ++(\wedgeAngle:3.55)
    ++({\wedgeAngle+90}:0.17) coordinate (regularLabel);
  \node[figure label,anchor=south]
    at (regularLabel) {$\partial\Omega_{\rm reg}$};
  \node[figure label,anchor=base]
    at (\tipX,-0.10) {$\partial\Omega_{\rm sing}$};
  \node[figure label] at ({\tipX+2.70},{\tipY+1.34}) {$\Omega$};
  \node[figure label,text=plotTertiary]
    at ({\tipX+3.70},{\tipY+0.77}) {$U$};
  \node[figure label,text=black!62]
    at ({\tipX+0.84},{\tipY+0.22}) {$\chi$};
  \node[figure label,text=plotSecondary,anchor=base]
    at ({\tipX+2.53},-0.10) {$\tilde\gamma$};
  \node[figure label,text=plotPrimary,anchor=base]
    at ({\tipX+3.92},-0.10) {$\gamma$};
\end{tikzpicture}%
        \caption{Lemma~\ref{l:hnegbreg}, Lemma~\ref{l:4t2}}
  \end{subfigure}\hfill
  \begin{subfigure}[t]{0.5\textwidth}
        \centering
        \begin{tikzpicture}[x=1cm,y=1cm]
  \definecolor{colorAKANE}{RGB}{183,40,46}
  \definecolor{colorTOKIWA}{RGB}{0,123,67}
  \definecolor{colorUSUAI}{RGB}{40,148,200}
  \colorlet{plotPrimary}{colorAKANE}
  \colorlet{plotSecondary}{colorTOKIWA}
  \colorlet{plotTertiary}{colorUSUAI}
  \colorlet{plotAxis}{black!72}

  \tikzset{
    >=latex, font=\small, mark size=1.5pt, mark options={thin},
    line cap=butt, line join=round,
    wedge edge/.style={draw=plotAxis,line width=0.95pt},
    inner edge/.style={draw=plotPrimary,line width=0.95pt},
    exterior edge/.style={draw=plotTertiary,line width=0.95pt},
    measurement edge/.style={draw=plotSecondary,line width=0.95pt,
      line cap=round},
    chi edge/.style={draw=black!62,line width=0.75pt,line cap=round},
    figure label/.style={text=black!84,inner sep=0.8pt}
  }

  \pgfmathsetmacro{\figureWidth}{7.5}
  \pgfmathsetmacro{\figureHeight}{3.2}
  \pgfmathsetmacro{\frameBottom}{-0.28}
  \pgfmathsetmacro{\frameTop}{\frameBottom+\figureHeight}
  \pgfmathsetmacro{\edgeDx}{4.26}
  \pgfmathsetmacro{\edgeDy}{2.48}
  \pgfmathsetmacro{\edgeLength}{veclen(\edgeDx,\edgeDy)}
  \pgfmathsetmacro{\tipX}{(\figureWidth-\edgeLength)/2}
  \pgfmathsetmacro{\tipY}{0.34}
  \pgfmathsetmacro{\wedgeAngle}{atan2(\edgeDy,\edgeDx)}
  \def\exteriorStartRadius{1.15}
  \def\innerEndRadius{1.60}

  \path[use as bounding box]
    (0,\frameBottom) rectangle (\figureWidth,\frameTop);
  \fill[fill=white]
    (0,\frameBottom) rectangle (\figureWidth,\frameTop);
  \coordinate (vertex) at (\tipX,\tipY);
  \coordinate (upperEnd) at ({\tipX+\edgeDx},{\tipY+\edgeDy});
  \coordinate (lowerEnd) at ({\tipX+\edgeLength},\tipY);
  \path (vertex) ++(0:\exteriorStartRadius) coordinate (exteriorLower);
  \path (vertex) ++(\wedgeAngle:\exteriorStartRadius) coordinate (exteriorUpper);
  \path (vertex) ++(0:\innerEndRadius) coordinate (innerLower);
  \path (vertex) ++(\wedgeAngle:\innerEndRadius) coordinate (innerUpper);

  \foreach \prefix/\start/\height in {tilde/2.20/1.48,chi/3.00/0.63} {
    \coordinate (\prefix Start) at ({\tipX+\start},\tipY);
    \coordinate (\prefix ControlOne) at
      ({\tipX+\start+0.05},{\tipY+0.68*\height});
    \coordinate (\prefix ControlTwo) at
      ({\tipX+\start+0.35},{\tipY+0.74*\height});
    \coordinate (\prefix Shoulder) at
      ({\tipX+\start+0.70},{\tipY+0.76*\height});
    \coordinate (\prefix End) at
      ({\tipX+\edgeLength+0.5},{\tipY+\height});
  }
  \coordinate (cropBeyond) at ({\tipX+\edgeLength+0.5},\tipY);

  \begin{scope}
    \clip (vertex) -- (lowerEnd)
      arc[start angle=0,end angle=\wedgeAngle,radius=\edgeLength cm]
      -- cycle;

    \fill[fill=plotTertiary!7]
      (exteriorLower) -- (lowerEnd)
      arc[start angle=0,end angle=\wedgeAngle,radius=\edgeLength cm]
      -- (exteriorUpper)
      arc[start angle=\wedgeAngle,end angle=0,radius=\exteriorStartRadius cm]
      -- cycle;
    \fill[fill=plotPrimary!10,fill opacity=0.75]
      (vertex) -- (innerLower)
      arc[start angle=0,end angle=\wedgeAngle,radius=\innerEndRadius cm]
      -- cycle;

    \fill[fill=plotSecondary!10]
      (tildeStart) .. controls (tildeControlOne) and (tildeControlTwo)
      .. (tildeShoulder) -- (tildeEnd) -- (cropBeyond) -- cycle;
    \fill[fill=black!7]
      (chiStart) .. controls (chiControlOne) and (chiControlTwo)
      .. (chiShoulder) -- (chiEnd) -- (cropBeyond) -- cycle;
    \draw[measurement edge]
      (tildeStart) .. controls (tildeControlOne) and (tildeControlTwo)
      .. (tildeShoulder) -- (tildeEnd);
    \draw[chi edge]
      (chiStart) .. controls (chiControlOne) and (chiControlTwo)
      .. (chiShoulder) -- (chiEnd);
  \end{scope}

  \draw[wedge edge] (upperEnd) -- (vertex) -- (lowerEnd);
  \draw[exterior edge] (exteriorLower)
    arc[start angle=0,end angle=\wedgeAngle,radius=\exteriorStartRadius cm];
  \draw[inner edge] (innerLower)
    arc[start angle=0,end angle=\wedgeAngle,radius=\innerEndRadius cm];
  \foreach \point in {exteriorLower,exteriorUpper}
    \fill[fill=plotTertiary] (\point) circle[radius=1.45pt];
  \foreach \point in {innerLower,innerUpper,vertex}
    \fill[fill=plotPrimary] (\point) circle[radius=1.45pt];
  \fill[fill=plotSecondary] (tildeStart) circle[radius=1.45pt];
  \fill[fill=black!62] (chiStart) circle[radius=1.45pt];

  \path (vertex) ++(\wedgeAngle:0.62)
    ++({\wedgeAngle+90}:0.15) coordinate (innerLabel);
  \path (vertex) ++(\wedgeAngle:2.43)
    ++({\wedgeAngle+90}:0.13) coordinate (exteriorLabel);
  \path (vertex) ++(\wedgeAngle:3.55)
    ++({\wedgeAngle+90}:0.17) coordinate (regularLabel);
  \node[figure label,text=plotPrimary,anchor=south]
    at (innerLabel) {$\gamma_i$};
  \node[figure label,text=plotTertiary,anchor=south]
    at (exteriorLabel) {$\gamma_e$};
  \node[figure label,anchor=south]
    at (regularLabel) {$\partial\Omega_{\rm reg}$};
  \node[figure label,anchor=base]
    at (\tipX,-0.10) {$\partial\Omega_{\rm sing}$};
  \node[figure label] at ({\tipX+1.93},{\tipY+0.45}) {$\Omega$};
  \node[figure label,text=plotSecondary]
    at ({\tipX+3.85},{\tipY+0.94}) {$\tilde\chi$};
  \node[figure label] at ({\tipX+3.85},{\tipY+0.25}) {$\chi$};
\end{tikzpicture}%
        \caption{Lemma~\ref{l:longDistanceHarmonic}}
  \end{subfigure}
  \caption{Illustration for Lemma~\ref{l:hnegbreg}, Lemma~\ref{l:4t2} and Lemma~\ref{l:longDistanceHarmonic}.}
  \label{f:4f1}
\end{figure}

\begin{lemma}[$\mc{H}_M^\Omega$ conserves negative $b$-regularity]\label{l:hnegbreg}
Let $k,s\in\mathbb{N}_0$. Let $m,l\in\mathbb{R}$ that $|m|<\rardyr_\Omega$, $|m-l|<\frac12$. Let $\gamma,\tilde\gamma\in S^0(\partial\Omega)$ with $\supp\tilde\gamma\subset\partial\oreg$ and $\supp\gamma\cap\supp(1-\tilde\gamma)=\emptyset$. Then for all vector fields $\mc{V}_1,\dots, \mc{V}_k\in \Diff^{1,1}_b(\partial\Omega)$, we have
\begin{equation}\label{e:4l13}
    \|\mc{H}^{\Omega}_M\gamma \mc{V}_1\cdots\mc{V}_kf\|_{X_{m,-s}^{l,-k}(\Omega)} \le C\|\tilde\gamma f\|_{H^{\frac{1}{2}-s,l}(\partial\Omega)},
\end{equation}
uniformly for all $f\in \mc{S}'(\partial\Omega)$ and $M\in [0,M_0]$. See Figure~\ref{f:4f1}(A). 
\end{lemma}
\begin{proof}
    0. Let $U\subset\Omega$ be a conic neighbourhood whose closure contains $\supp\tilde\gamma\subset\partial\oreg$, on which we have the Fermi normal coordinates. We confuse $\gamma\in S^0(\partial\Omega)$ with one of its extensions $\gamma\in S^0(U)$ that $\partial_\nu \gamma=0$. This similarly applies to $\tilde\gamma$. We confuse $\mc{V}_j\in \Diff^{1,1}_b(\partial\Omega)$ with their extension $\mc{V}_j\in \Diff^{1,1}_b(\Omega)$ with $\supp \mc{V}_j\subset U$. Assume without loss of generality none of $\mc{V}_j$'s contains a component in $x_0D_{x_0}$: if it does, then $x_0D_{x_0}$ is effectively 0 on $\partial\Omega$. 
    
    1. We claim
\begin{equation}
\mc{V}_1\cdots\mc{V}_k\partial_\nu=\partial_\nu\mc{V}_1\cdots\mc{V}_k+\sum_{|\alpha|\le k,\alpha_0=0}O(S^{-1}(U))(\rre D)^\alpha. 
\end{equation}
The case $k=0$ is trivial. For the case $k=1$, one can check in local coordinates in $U$ that
\begin{equation}\label{e:4l22}
    \mc{V}_1\partial_\nu=\partial_\nu \mc{V}_1+O(S^{-1}(U))+\sum_{j=1}^{d-1}O(S^{-1}(U))\rre D_{x_j},
\end{equation}
by identifying $U$ as $\tilde\Omega_+=\{(x_0,x): x_0>0\}$ as in Remark~\ref{r:aclcollar} and $\partial_\nu$ by $-\partial_{x_0}$. Assume the case $k-1$ holds, and 
\begin{multline}
\mc{V}_1\cdots\mc{V}_k\partial_\nu=\mc{V}_1\big(\partial_\nu\mc{V}_2\cdots\mc{V}_k+\sum_{|\alpha|\le k-1,\alpha_0=0}O(S^{-1}(U))(\rre D)^\alpha\big)=\partial_\nu\mc{V}_1\cdots\mc{V}_k\\
+\sum_{|\alpha|\le 1, \alpha_0=0}O(S^{-1}(U))(\rre D)^\alpha\mc{V}_2\cdots\mc{V}_k+\sum_{|\alpha|\le k-1,\alpha_0=0}\mc{V}_1O(S^{-1}(U))(\rre D)^\alpha
\end{multline}
justifies the iterative claim. 

2. Let multi-index $|\beta|\le s$ with $\beta_0=0$. We claim uniformly for all $\phi\in \dot Y_{-m,s}^{-l,k}(\Omega)$, 
\begin{equation}
    \|D^\beta\mc{V}_1\cdots\mc{V}_k\gamma \partial_\nu \mc{R}^{\Omega}_M \phi\|_{H^{-\frac12,-l}(\partial\Omega)}\le C\|\phi\|_{Y_{-m,s}^{-l,k}}
\end{equation}
in this step. First note $D^\beta$ commutes with $\partial_\nu$ in $U$, and hence Step 1 implies
\begin{equation}\label{e:4l20}
    D^\beta\mc{V}_1\cdots\mc{V}_k\gamma \partial_\nu \mc{R}^{\Omega}_M \phi=\partial_\nu D^\beta\mc{V}_1\cdots\mc{V}_k\gamma\mc{R}^{\Omega}_M\phi+\sum_{|\alpha|\le k,\alpha_0=0}O(S^{-1})D^\beta(\rre D)^\alpha\gamma\mc{R}^{\Omega}_M\phi,
\end{equation}
the second term of which can be estimated via Lemma~\ref{l:longDistanceResolve} by
\begin{multline}
    \|O(S^{-1})D^\beta(\rre D)^\alpha\gamma\mc{R}^{\Omega}_M\phi\|_{H^{-\frac12,-l}_\partial}\le C\sum_{|\alpha|\le k, \alpha_0=0}\|D^\beta(\rre D)^\alpha\gamma\mc{R}^{\Omega}_M\phi\|_{X_{-m}^{-l}}\\
    \le C\|\phi\|_{Y_{-m}^{-l}}+C\sum_{|\alpha|\le k, \alpha_0=0}\|\tilde\gamma D^\beta(\rre D)^\alpha\phi\|_{Y_{-m}^{-l}}\le C\|\phi\|_{Y^{-l,k}_{-m,s}},
\end{multline}
For the first term of \eqref{e:4l20}, we have
\begin{multline}
\|\partial_\nu D^\beta\mc{V}_1\cdots\mc{V}_k\gamma\mc{R}^{\Omega}_M\phi\|_{H^{-\frac12,-l}_\partial}\le C\|\partial_\nu\langle x_0\rangle^{-N}\psi(x_0) D^\beta\mc{V}_1\cdots\mc{V}_k\gamma\mc{R}^{\Omega}_M\phi\|_{H^{-\frac12,-l}_\partial}\\
\le C\|\langle x_0\rangle^{-N}(-\Delta+M^2)D^\beta\mc{V}_1\cdots\mc{V}_k\gamma\mc{R}^{\Omega}_M\phi\|_{H^{-1,-l}}+C\|\psi(x_0)\nabla D^\beta\mc{V}_1\cdots\mc{V}_k\gamma\mc{R}^{\Omega}_M\phi\|_{H^{0,-l}}\\
+C\|\langle x_0\rangle^{-N}D^\beta\mc{V}_1\cdots\mc{V}_k\gamma\mc{R}^{\Omega}_M\phi\|_{H^{-1,-l}},
\end{multline}
for some $\psi\in C^\infty_c(\mathbb{R})$ that $\psi(x_0)\equiv 1$ near $x_0=0$. By Lemma~\ref{l:longDistanceResolve} we have
\begin{multline}
    \|\langle x_0\rangle^{-N}(-\Delta+M^2)D^\beta\mc{V}_1\cdots\mc{V}_k\gamma\mc{R}^{\Omega}_M\phi\|_{H^{-1,-l}}\le C\|(-\Delta+M^2)D^\beta\mc{V}_1\cdots\mc{V}_k\gamma\mc{R}^{\Omega}_M\phi\|_{Y_{-m}^{-l}}\\
    \le C\|\phi\|_{Y_{-m}^{-l}}+ C\sum_{|\alpha|\le k}\|\tilde\gamma D^\beta(\rre D)^\alpha \phi\|_{Y_{-m}^{-l}}\le C\|\phi\|_{Y_{-m,s}^{-l,k}}
\end{multline}
and
\begin{equation}
    \|\psi(x_0)\nabla D^\beta\mc{V}_1\cdots\mc{V}_k\gamma\mc{R}^{\Omega}_M\phi\|_{H^{0,-l}}\le C\|D^\beta\mc{V}_1\cdots\mc{V}_k\gamma\mc{R}^{\Omega}_M\phi\|_{X_{-m}^{-l}}\le C\|\phi\|_{Y_{-m,s}^{-l,k}}.
\end{equation}
This shows the claim. 

3. Noting $(\mc{H}^{\Omega}_M)^*=-\partial_\nu \mc{R}^{\Omega}_M$ we have for any $\phi\in \dot Y_{-m,s}^{-l,k}(\Omega)$,
\begin{multline}
    |\langle \mc{H}_M^\Omega \gamma \mc{V}_1\cdots\mc{V}_kf,\phi\rangle|=|\langle\tilde\gamma f,\mc{V}_k^*\cdots\mc{V}_1^*\gamma \partial_\nu \mc{R}^{\Omega}_M \phi\rangle_{\partial\Omega}|\\
    \le C\|\tilde\gamma f\|_{H^{\frac12-s,l}_\partial}\sum_{|\beta|\le s, \beta_0=0}\| D^\beta\mc{V}_k^*\cdots\mc{V}_1^*\gamma \partial_\nu \mc{R}^{\Omega}_M \phi\|_{H^{-\frac12,-l}}\le C\|\tilde\gamma f\|_{H^{\frac12-s,l}_\partial}\|\phi\|_{Y_{-m,s}^{-l,k}}.
\end{multline}
This demonstrates the desired bounds over $\mc{H}_M^\Omega \gamma \mc{V}_1\cdots\mc{V}_kf$.
\end{proof}

\begin{corollary}\label{c:HimproveM}
Let $m,l\in\mathbb{R}$ that $|m|<\rardyr_\Omega$, $|m-l|<\frac12$. Let $m_2\in\mathbb{R}$. We have
\begin{equation}
    \|\langle M\dd\rangle^{m_2}\mc{H}^{\Omega}_M f\|_{Y_{m-2,-1}^{l,0}(\Omega)} \le C\|f\|_{H^{0,l}(\partial\Omega)},
\end{equation}
uniformly for all $f\in \mc{S}'(\partial\Omega)$ and $M\in [0,M_0]$.
\end{corollary}
\begin{proof}
    1. Let $\gamma,\tilde\gamma\in S^0(\partial\Omega)$ with $\supp\tilde\gamma\subset\partial\oreg$ and $\supp\gamma\cap\supp(1-\tilde\gamma)=\emptyset$. We show
\begin{equation}\label{e:4l41}
    \|\langle M\dd\rangle^{m_2}\mc{H}^{\Omega}_M\gamma f\|_{X_{m,-1}^{l,0}(\Omega)} \le C\|\tilde\gamma f\|_{H^{-\frac12,l}(\partial\Omega)}.
\end{equation}
The proof is similar to that of Lemma~\ref{l:hnegbreg}. Note $[D_{x_j}, \partial_\nu]=0$ for the boundedness of
    \begin{equation}
        D_{x_j}\gamma\partial_\nu\langle M\dd\rangle^{-m_2}\mc{R}^{\Omega}_M \langle M\dd\rangle^{m_2}=\partial_\nu D_{x_j}\gamma\langle M\dd\rangle^{-m_2}\mc{R}^{\Omega}_M \langle M\dd\rangle^{m_2}: Y_{-m, 1}^{-l, 0}\rightarrow H^{-\frac12, -l}_\partial,
    \end{equation}
    where we used Corollary~\ref{c:RimproveM}. Like in Step 3 of Lemma~\ref{l:hnegbreg}, for each $\phi\in \dot Y_{-m, 1}^{-l, 0}(\Omega)$ we have
    \begin{multline}
    |\langle \langle M\dd\rangle^{m_2}\mc{H}_M^\Omega \gamma f,\phi\rangle|=|\langle\langle M\dd\rangle^{m_2}\tilde\gamma f,\gamma \partial_\nu \langle M\dd\rangle^{-m_2}\mc{R}^{\Omega}_M \langle M\dd\rangle^{m_2}\phi\rangle_{\partial\Omega}|\\
    \le C\|\tilde\gamma f\|_{H^{-\frac12,l}_\partial}\big(\|\gamma \partial_\nu \langle M\dd\rangle^{-m_2}\mc{R}^{\Omega}_M \langle M\dd\rangle^{m_2}\phi\|_{H^{-\frac12,-l}_\partial}+\sum_{j=1}^{d-1}\| D_{x_j}\gamma \partial_\nu \langle M\dd\rangle^{-m_2}\mc{R}^{\Omega}_M \langle M\dd\rangle^{m_2}\phi\|_{H^{-\frac12,-l}_\partial}\big)\\\le C\|\tilde\gamma f\|_{H^{-\frac12,l}_\partial}\|\phi\|_{Y_{-m,1}^{-l,0}}
\end{multline}
yielding the desired bounds. 

2. Let $\kappa,\kappa_1, \tilde\kappa\in C^\infty_c(\overline{\Omega})$ be that $\supp(1-\kappa)\cap \partial\Omega \subset\partial\oreg$, and $\supp\kappa\cap\supp(1-\kappa_1)=\supp\kappa_1\cap\supp(1-\tilde\kappa)=\emptyset$. Let $\varphi\in C^\infty_c(\overline{\Omega})$ that $\supp\nabla \kappa_1\cap\supp(1-\varphi)=\emptyset$. We decompose
\begin{equation}
    \mc{H}_M^{\Omega}f=\tilde\kappa\mc{H}_M^{\Omega}f+(1-\tilde\kappa)\mc{H}_M^{\Omega}(1-\kappa)f+(1-\tilde\kappa)\mc{H}_M^{\Omega}\kappa f
\end{equation}
and estimate it term by term. Lemma~\ref{l:ellregunbdd} implies 
\begin{equation}
    \|\langle M\dd\rangle^{m_2}\tilde\kappa\mc{H}_M^{\Omega}f\|_{Y_{N,0}^{N,0}}\le C\|\langle M\dd\rangle^{m_2}\tilde\kappa\mc{H}_M^{\Omega}f\|_{H^{-1}}\le C\|f\|_{H^{0, l}_\partial}.
\end{equation}
For the second term, Step 1 implies
\begin{equation}
    \|\langle M\dd\rangle^{m_2}(1-\tilde\kappa)\mc{H}_M^{\Omega}(1-\kappa)f\|_{X_{m,-1}^{l,0}}\le C\|(1-\kappa)f\|_{H^{-\frac12, l}_\partial}.
\end{equation}
For the last term, utilizing $(1-\tilde\kappa)\mc{H}_M^{\Omega}\kappa=(1-\tilde\kappa)\mc{R}_M^{\Omega}[\Delta, \kappa_1]\mc{H}_M^{\Omega}\kappa$ we have
\begin{multline}
    \|\langle M\dd\rangle^{m_2}(1-\tilde\kappa)\mc{R}_M^{\Omega}[\Delta, \kappa_1]\mc{H}_M^{\Omega}\kappa f\|_{X_{m,-1}^{l,0}}\le C\|\langle M\dd\rangle^{m_2}[\Delta, \kappa_1]\mc{H}_M^{\Omega}\kappa f\|_{Y_{m,-1}^{l,0}}\\
    \le C\|\varphi\mc{H}_M^{\Omega}\kappa f\|_{H^{-1}} \le C\|\kappa f\|_{H_\partial^{0, l}}
\end{multline}
from Corollary~\ref{c:RimproveM} and Lemma~\ref{l:ellregunbdd}. Summing up above gives
\begin{multline}
    \|\langle M\dd\rangle^{m_2}(1-\tilde\kappa)\mc{H}_M^{\Omega}(1-\kappa)f+\langle M\dd\rangle^{m_2}(1-\tilde\kappa)\mc{H}_M^{\Omega}\kappa f\|_{Y^{l,0}_{m-2, -1}}\\
    \le C\|\langle M\dd\rangle^{m_2}(1-\tilde\kappa)\mc{H}_M^{\Omega}(1-\kappa)f+\langle M\dd\rangle^{m_2}(1-\tilde\kappa)\mc{H}_M^{\Omega}\kappa f\|_{X^{l,0}_{m, -1}}\le C\|f\|_{H^{0,l}_\partial},
\end{multline}
which concludes the proof. 
\end{proof}

\begin{lemma}[$\mc{H}_M^\Omega$ conserves positive $b$-regularity]
\label{l:longDistanceHarmonic}
Let $\Omega$ be an \acl domain. Let $k,s\in\mathbb{N}_0$. Let $m,l\in\mathbb{R}$ that $|m|<\rardyr_\Omega$, $|m-l|<\frac12$. Let vector fields $\mc{V}_1,\dots \mc{V}_k\in \Diffhom_b^{1,1}(\Omega)$, $\mc{U}_1,\dots \mc{U}_s\in \Diffhom_b^{1,0}(\Omega)$. Let $\chi,\tilde\chi\in S^0(\Omega)$ with $\supp\chi\cap\supp(1-\tilde\chi)=\emptyset$. Let $\gamma_i\in C_c^\infty(\overline{\Omega})$, $\gamma_e\in S^0(\Omega)$ that $\gamma_i+\gamma_e=1$ on $\Omega$, $\supp\tilde\chi\cap\supp\gamma_i=\emptyset$, $(\supp\gamma_e\cap \partial\Omega)\cap \partial\osing=\emptyset$. For each $N\ge 1$ large, we then have
\begin{multline}
    \|\chi\mc{U}_1\cdots \mc{U}_s \mc{V}_1\cdots \mc{V}_k\mathcal{H}^\Omega_{M}f\|_{X_{m}^l(\Omega)}\le C\sum_{|\alpha|\le k,|\beta|\le s}\|\tilde\chi D^\beta(\langle x\rangle D)^\alpha f\|_{H^{\frac{1}{2},l}(\partial\Omega)}+C\|\gamma_i f\|_{H^{0, -N}(\partial\Omega)}\\
    +C\|\gamma_e f\|_{H^{-N, m-\frac12+0}(\partial\Omega)},
\end{multline}
uniformly for all $f\in \mc{S}'(\partial\Omega)$ and $M\in [0,M_0]$. See Figure~\ref{f:4f1}(B). 
\end{lemma}

\begin{proof}
0. Let $\tilde\gamma_i\in C_c^\infty(\overline{\Omega})$ that $\supp\gamma_i\cap\supp(1-\tilde\gamma_i)=\emptyset$, $\supp\tilde\gamma_i\cap \supp\chi=\emptyset$ and $\tilde\gamma_e\in S^0$ that $\supp\gamma_e\cap\supp(1-\tilde\gamma_e)=\emptyset$, $(\supp\tilde\gamma_e\cap \partial\Omega)\cap \partial\osing=\emptyset$. Let $\kappa\in C^\infty_c(\bar\Omega)$ that $\supp\nabla\tilde\gamma_i\cap\supp(1-\kappa)=\emptyset$ and $\supp\kappa\cap\supp\tilde\chi=\emptyset$. Let $\chi_1,\chi_2\in S^0(\Omega)$ that satisfy all requirements for $\chi$ with $\supp\chi\cap\supp(1-\chi_1)=\supp\chi_1\cap\supp(1-\chi_2)=\supp\chi_2\cap\supp(1-\tilde\chi)=\emptyset$. We decompose
\begin{equation}\label{e:4l32}
    \chi\mc{V}_1\cdots \mc{V}_k\mc{H}^\Omega_Mf=\chi\mc{V}_1\cdots \mc{V}_k\mc{H}^\Omega_M\gamma_e \tilde\chi f+\chi\mc{V}_1\cdots \mc{V}_k\mc{H}^\Omega_M\gamma_e (1-\tilde\chi) f+\chi\mc{V}_1\cdots \mc{V}_k\mc{H}^\Omega_M\gamma_i f
\end{equation}
We will estimate \eqref{e:4l32} term by term in the rest of the proof to conclude
\begin{multline}
    \|\chi\mc{V}_1\cdots \mc{V}_k\mathcal{H}^\Omega_{M}f\|_{X_{m}^l(\Omega)}\le C\sum_{|\alpha|\le k}\|\tilde\chi (\langle x\rangle D)^\alpha f\|_{H^{\frac{1}{2},l}(\partial\Omega)}+C\|\gamma_i f\|_{H^{0, -N}(\partial\Omega)}\\
    +C\|\gamma_e f\|_{H^{-N, m-\frac12+0}(\partial\Omega)},
\end{multline}
for all $k\in\mathbb{N}_0$. Adding $\mc{U}_{1}\cdots\mc{U}_s$ is similar: see Lemma~\ref{l:longDistanceResolve}. This would close the proof. 

1. For the first term in \eqref{e:4l32}, we iteratively show for $k\in\mathbb{N}_0$ that 
\begin{equation}\label{e:4l14}
    \|\chi\mc{V}_1\cdots \mc{V}_k \tilde\gamma_e\mc{H}_M^{\Omega}\gamma_e \tilde\chi f\|_{X_{m}^{l}(\Omega)}\leq C\sum_{|\alpha|\leq k}\|\tilde\chi (\langle x\rangle D)^\alpha f\|_{H^{\frac{1}{2},l}(\partial\Omega)}.
\end{equation}
The case $k=0$ follows from Proposition~\ref{p:harmonic}. 
Assume \eqref{e:4l14} holds for $k-1$. Noting $\chi_1\mc{V}_k\mc{H}_M^{\Omega}\gamma_e \tilde\chi f$ extends $\chi_1\mc{V}_k f$, we have by the extension correction formula \eqref{e:hForm}, 
\begin{equation}
    \tilde\gamma_e\mc{H}_M^{\Omega} \chi_1\mc{V}_k f=\chi_1\mc{V}_k\tilde\gamma_e\mc{H}_M^{\Omega}\gamma_e \tilde\chi f-\tilde\gamma_e\mc{R}_M^{\Omega}[-\Delta,\chi_1\mc{V}_k]\tilde\gamma_e\mc{H}_M^{\Omega}\gamma_e \tilde\chi f
\end{equation}
and thus by $[-\Delta, \chi_1\mc{V}_k]=\oxy\chi_2$ from Lemma~\ref{l:algCommute}, 
\begin{equation}
    \chi\mc{V}_1\cdots \mc{V}_k\tilde\gamma_e\mc{H}_M^{\Omega}\gamma_e\tilde\chi f=\chi\mc{V}_1\cdots \mc{V}_{k-1}\tilde\gamma_e\mc{H}_M^{\Omega}\chi_1\mc{V}_k f+\chi\mc{V}_1\cdots \mc{V}_{k-1}\tilde\gamma_e\mc{R}_M^{\Omega}\oxy\chi_2\mc{H}_M^{\Omega}\gamma_e \tilde\chi f.
\end{equation}
By the iterative assumption, the first term can be bounded by
\begin{equation}
    \|\chi\mc{V}_1\cdots \mc{V}_{k-1}\tilde\gamma_e\mc{H}_M^{\Omega}\chi_1\mc{V}_k \gamma_ef\|_{X_m^l}\le C\sum_{|\alpha|\leq k-1}\|\chi_1 (\langle x\rangle D)^\alpha \mc{V}_k f\|_{H^{\frac{1}{2},l}_\partial} \le C\sum_{|\alpha|\leq k}\|\tilde\chi (\langle x\rangle D)^\alpha f\|_{H^{\frac{1}{2},l}_\partial}.
\end{equation}
With the iterative assumption with $\chi_2$ in place of $\chi$, the second term can be bounded by Lemma~\ref{l:longDistanceResolve}, 
\begin{equation}
    \|\chi\mc{V}_1\cdots \mc{V}_{k-1}\tilde\gamma_e\mc{R}_M^{\Omega}\oxy\chi_2\mc{H}_M^{\Omega}\gamma_e\tilde\chi  f\|_{X_m^l}\le C\|\chi_2\mc{H}_M^{\Omega}\gamma_e \tilde\chi f\|_{X_m^{l,k-1}}\le C\sum_{|\alpha|\leq k-1}\|\tilde\chi (\langle x\rangle D)^\alpha f\|_{H^{\frac{1}{2},l}_\partial}.
\end{equation}
This yields \eqref{e:4l14} in the case $k$ as desired. 

2. We claim in this step the second term of \eqref{e:4l32} can be bounded by 
\begin{equation}
    \|\chi\mc{V}_1\cdots \mc{V}_k\chi_1\mc{H}^{\Omega}_M\gamma_e(1-\tilde\chi)f\|_{X_m^l(\Omega)}\le C\|\gamma_e f\|_{H^{\frac12-N, m-\frac12+0}_\partial}. 
\end{equation}
Let $\varphi\in S^0(\Omega)$ such that $\supp\gamma_e(1-\tilde\chi)\cap \supp(1-\varphi)=\supp\varphi\cap\supp\chi_1=\emptyset$. Let $\varphi_1, \varphi_2,\dots$ satisfy the same requirements for $\varphi$ while $\supp\varphi_n\cap\supp(1-\varphi_{n+1})=\emptyset$ for all $n$, $\supp\varphi\cap\supp(1-\varphi_1)=\emptyset$. We first note
\begin{equation}
    \chi_1\mc{H}_M^{\Omega}\varphi=\chi_1\mc{R}_M^{\Omega}[\Delta, \varphi_{N+1}]\mc{R}_M^{\Omega}[\Delta, \varphi_{N}]\cdots\mc{R}_M^{\Omega}[\Delta, \varphi_1]\varphi_2\mc{H}_M^{\Omega}\varphi. 
\end{equation}
Lemma~\ref{l:smoothing} implies $\chi_1\mc{R}_M^{\Omega}[\Delta, \varphi_{N+1}]: X_m^{l, -N}\rightarrow X_m^{l, N}$ for $m,l$ within $|m|<\rardyr_{\Omega}$, $|m-l|<\frac12$. Lemma~\ref{l:oxybddxmlk} implies for each $n$, $\mc{R}_M^{\Omega}[\Delta, \varphi_{n}]:X_m^{l-1,-N}\rightarrow X_m^{l,-N}$ for any $m, l$ within the same range. Then by choosing a large $N$, we have 
\begin{equation}\label{e:4l39}
    \|\chi_1\mc{H}_M^{\Omega}\varphi f\|_{X_m^{l,N}}\le C\|\varphi_2 \mc{H}_M^{\Omega}\varphi f\|_{X_{m}^{m-\frac12+0, -N}}\le C\|\varphi_2 \mc{H}_M^{\Omega}\varphi f\|_{X_{m,-N}^{m-\frac12+0}}
\end{equation}
for any $f$. Thus for $N$ large we have
\begin{equation}
    \|\chi\mc{V}_1\cdots \mc{V}_k\chi_1\mc{H}^{\Omega}_M\gamma_e(1-\tilde\chi)f\|_{X_m^l}\le C\|\varphi_2\mc{H}^{\Omega}_M \gamma_e (1-\tilde\chi)f\|_{X_{m,-N}^{m-\frac12+0}}\le C\|\gamma_e f\|_{H^{\frac12-N, m-\frac12+0}_\partial}
\end{equation}
from Lemma~\ref{l:hnegbreg}.

3. For the last term in \eqref{e:4l32}, we have
\begin{equation}
    \chi\mc{V}_1\cdots \mc{V}_k\mc{H}^\Omega_M\gamma_i f=\chi\mc{V}_1\cdots \mc{V}_k\mc{R}^\Omega_M[\Delta, \tilde\gamma_i]\mc{H}_M^{\Omega}\gamma_i f.
\end{equation}
Lemma~\ref{l:smoothing} and Lemma~\ref{l:ellregunbdd} then imply
\begin{equation}
    \|\chi\mc{V}_1\cdots \mc{V}_k\mc{R}^\Omega_M[\Delta, \tilde\gamma_i]\mc{H}_M^{\Omega}\gamma_i f\|_{X_m^l(\Omega)}\le C\|\kappa \mc{H}_M^{\Omega}\gamma_i f\|_{L^2(\Omega)}\le C\|\gamma_i f\|_{H^{0,-N}(\partial\Omega)}. 
\end{equation}
since $\supp\chi\cap\supp\kappa=\emptyset$ and $\kappa$ is compactly supported. 
\end{proof}

\begin{lemma}[Nonlocal effect of $\mc{H}_M^\Omega$ on smooth part is $b$-smoothing]\label{l:4t2}
    Let $\Omega$ be an \acl domain. Let $m,l\in\mathbb{R}$ that $|m|<\rardyr_\Omega$, $|m-l|<\frac12$. Let $\gamma\in S^0(\partial\Omega)$ with $\supp\gamma\subset\partial\oreg$. Let $\chi\in S^0(\Omega)$ with $\supp\gamma\cap(\supp\chi\cap\partial\Omega)=\emptyset$. Then for all $k\in\mathbb{N}_0$, for all $\mc{V}_1,\dots, \mc{V}_k\in \Diff_b^{1,1}(\Omega)$, for any $N>0$ large, we have
    \begin{equation}\label{e:4l3}
        \chi\mathcal{H}_M^\Omega \gamma \mathcal{V}_1\cdots\mathcal{V}_k: H^{-N,m-\frac12+0}(\partial\Omega)\rightarrow \dot X_m^{l,N}(\Omega)
    \end{equation}
    uniformly in $M\in [0,M_0]$. See Figure~\ref{f:4f1}(A). 
\end{lemma}
\begin{proof}
    1. We confuse $\gamma\in S^0(\partial\Omega)$ with one of its extensions $\gamma\in S^0(\Omega)$ such that $\supp\gamma\cap \supp\chi=\emptyset$. Let $\tilde\gamma\in S^0(\Omega)$ that $(\supp\tilde\gamma\cap \partial\Omega)\subset\partial\oreg$, $\supp\gamma\cap\supp(1-\tilde\gamma)=\emptyset$ while $\supp\tilde\gamma\cap \supp\chi=\emptyset$. Let $\gamma_1, \gamma_2, \dots$ satisfy the same conditions for $\tilde\gamma$ with $\supp\gamma_j\cap\supp(1-\gamma_{j+1})=\supp\gamma_j\cap\supp(1-\tilde\gamma)=\emptyset$ for each $j$, and $\supp\gamma\cap\supp(1-\gamma_1)=\emptyset$. Without loss of generality picking $N$ large that $\frac12+k-2N<-N$, we have from Lemma~\ref{l:hnegbreg} that
    \begin{multline}\label{e:4l4}
        \|\mathcal{H}_M^\Omega \gamma \mc{V}_1\cdots \mc{V}_k f\|_{X_m^{m-\frac12+0,-2N}}\le C\|\mathcal{H}_M^\Omega \gamma \mc{V}_1\cdots \mc{V}_k f\|_{X_{m,k-2N}^{m-\frac12+0,-k}}\le C\|\tilde\gamma f\|_{H^{\frac12+k-2N,m-\frac12+0}_\partial}\\
        \le C\|\tilde\gamma f\|_{H^{-N,m-\frac12+0}_\partial},
    \end{multline}
    since $X_{m,-s}^{m-\frac12+0,-k}\subset X_{m}^{m-\frac12+0,-k-s}$ for all $s,k\in \mathbb{N}_0$.

   2. The extension correction formula \eqref{e:hForm} implies
   \begin{equation}
       \chi\mathcal{H}_M^\Omega\gamma\mc{V}_1\cdots \mc{V}_kf=\chi\mathcal{R}_M^\Omega \tilde\gamma[\Delta, \gamma_{N+1}]\mathcal{R}_M^\Omega [\Delta, \gamma_{N}]\cdots\mathcal{R}_M^\Omega [\Delta, \gamma_1]\gamma_2\mathcal{H}_{M}^\Omega\gamma\mc{V}_1\cdots \mc{V}_kf.
   \end{equation}
    Since $\chi\mathcal{R}_M^\Omega \tilde\gamma$ is bounded per Lemma~\ref{l:smoothing}, Lemma~\ref{l:oxybddxmlk} now implies
    \begin{multline}
        \|\chi\mathcal{H}_M^\Omega\gamma\mc{V}_1\cdots \mc{V}_kf\|_{X_{m}^{l,N}}\le C\|\mathcal{R}_M^\Omega [\Delta, \gamma_{N}]\cdots\mathcal{R}_M^\Omega [\Delta, \gamma_1]\gamma_2\mathcal{H}_{M}^\Omega\gamma\mc{V}_1\cdots \mc{V}_kf\|_{X^{l,-N}_{m}}\\
        \le C\|\mathcal{H}_{M}^\Omega\gamma\mc{V}_1\cdots \mc{V}_kf\|_{X^{m-\frac12+0,-N}_{m}},
    \end{multline}
    similar to Step 2 of Lemma~\ref{l:longDistanceHarmonic}. Revisit \eqref{e:4l4} to conclude the proof. 
\end{proof}

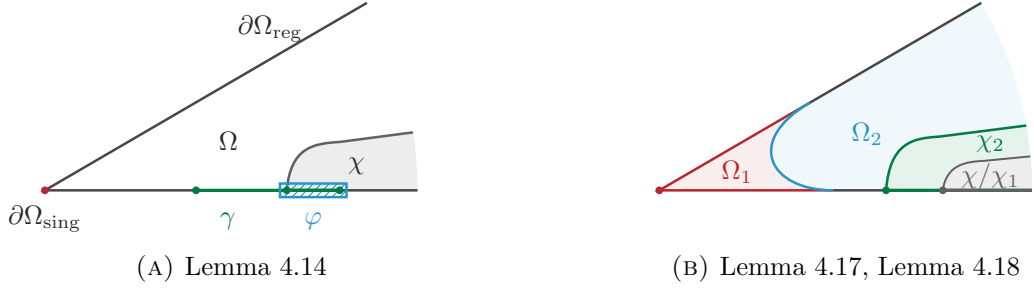
\begin{figure}[t]
  \centering
  \begin{subfigure}[t]{0.5\textwidth}
        \centering
        \begin{tikzpicture}[x=1cm,y=1cm]
  \usetikzlibrary{patterns}
  \definecolor{colorAKANE}{RGB}{183,40,46}
  \definecolor{colorTOKIWA}{RGB}{0,123,67}
  \definecolor{colorUSUAI}{RGB}{40,148,200}
  \colorlet{plotPrimary}{colorAKANE}
  \colorlet{plotSecondary}{colorTOKIWA}
  \colorlet{plotTertiary}{colorUSUAI}
  \colorlet{plotAxis}{black!72}

  \tikzset{
    >=latex, font=\small, mark size=1.5pt, mark options={thin},
    line cap=butt, line join=round,
    wedge edge/.style={draw=plotAxis,line width=0.95pt},
    inner support/.style={draw=plotSecondary,line width=1.2pt,solid},
    overlap support/.style={draw=plotSecondary,line width=1.2pt,solid},
    cutoff box/.style={draw=plotTertiary,line width=0.95pt},
    chi edge/.style={draw=black!62,line width=0.95pt,line cap=round},
    figure label/.style={text=black!84,inner sep=0.8pt}
  }

  \pgfmathsetmacro{\figureWidth}{7.5}
  \pgfmathsetmacro{\figureHeight}{3.2}
  \pgfmathsetmacro{\frameBottom}{-0.28}
  \pgfmathsetmacro{\frameTop}{\frameBottom+\figureHeight}
  \pgfmathsetmacro{\edgeDx}{4.26}
  \pgfmathsetmacro{\edgeDy}{2.48}
  \pgfmathsetmacro{\edgeLength}{veclen(\edgeDx,\edgeDy)}
  \pgfmathsetmacro{\tipX}{(\figureWidth-\edgeLength)/2}
  \pgfmathsetmacro{\tipY}{0.34}
  \pgfmathsetmacro{\wedgeAngle}{atan2(\edgeDy,\edgeDx)}

  \path[use as bounding box]
    (0,\frameBottom) rectangle (\figureWidth,\frameTop);
  \fill[fill=white]
    (0,\frameBottom) rectangle (\figureWidth,\frameTop);
  \coordinate (vertex) at (\tipX,\tipY);
  \coordinate (upperEnd) at ({\tipX+\edgeDx},{\tipY+\edgeDy});
  \coordinate (lowerEnd) at ({\tipX+\edgeLength},\tipY);
  \coordinate (cropBeyond) at ({\tipX+\edgeLength+0.5},\tipY);

  \coordinate (gammaStart) at ({\tipX+2.00},\tipY);
  \coordinate (gammaEnd) at ({\tipX+3.90},\tipY);
  \coordinate (chiStart) at ({\tipX+3.20},\tipY);
  \coordinate (varphiLowerLeft) at
    ([xshift=-0.09cm,yshift=-0.09cm]chiStart);
  \coordinate (varphiUpperRight) at
    ([xshift=0.09cm,yshift=0.09cm]gammaEnd);
  \path (chiStart) -- (gammaEnd) coordinate[midway] (overlapMiddle);
  \coordinate (chiControlOne) at ({\tipX+3.25},{\tipY+0.57});
  \coordinate (chiControlTwo) at ({\tipX+3.55},{\tipY+0.62});
  \coordinate (chiShoulder) at ({\tipX+3.90},{\tipY+0.64});
  \coordinate (chiEnd) at
    ({\tipX+\edgeLength+0.5},{\tipY+0.84});

  \begin{scope}
    \clip (vertex) -- (lowerEnd)
      arc[start angle=0,end angle=\wedgeAngle,radius=\edgeLength cm]
      -- cycle;
    \fill[fill=black!7]
      (chiStart) .. controls (chiControlOne) and (chiControlTwo)
      .. (chiShoulder) -- (chiEnd) -- (cropBeyond) -- cycle;
    \draw[chi edge]
      (chiStart) .. controls (chiControlOne) and (chiControlTwo)
      .. (chiShoulder) -- (chiEnd);
  \end{scope}

  \path[pattern=north east lines,pattern color=plotTertiary]
    (varphiLowerLeft) rectangle (varphiUpperRight);

  \draw[wedge edge] (upperEnd) -- (vertex) -- (lowerEnd);
  \draw[inner support] (gammaStart) -- (gammaEnd);
  \draw[overlap support] (chiStart) -- (gammaEnd);
  \fill[fill=plotSecondary] (gammaStart) circle[radius=1.45pt];
  \foreach \point in {chiStart,gammaEnd}
    \fill[fill=plotSecondary] (\point) circle[radius=1.45pt];
  \fill[fill=plotPrimary] (vertex) circle[radius=1.45pt];

  \draw[cutoff box] (varphiLowerLeft) rectangle (varphiUpperRight);

  \path (vertex) ++(\wedgeAngle:3.55)
    ++({\wedgeAngle+90}:0.17) coordinate (regularLabel);
  \node[figure label,anchor=south]
    at (regularLabel) {$\partial\Omega_{\rm reg}$};
  \node[figure label,anchor=base]
    at (\tipX,-0.10) {$\partial\Omega_{\rm sing}$};
  \node[figure label] at ({\tipX+2.43},{\tipY+0.65}) {$\Omega$};
  \node[figure label] at ({\tipX+4.12},{\tipY+0.32}) {$\chi$};
  \node[figure label,text=plotSecondary,anchor=base]
    at ({\tipX+2.42},-0.10) {$\gamma$};
  \node[figure label,text=plotTertiary,anchor=base]
    at ([yshift=-0.44cm]overlapMiddle) {$\varphi$};
\end{tikzpicture}%
        \caption{Lemma~\ref{l:longDistanceHarmonicRight}}
  \end{subfigure}\hfill
    \begin{subfigure}[t]{0.5\textwidth}
        \centering
        \begin{tikzpicture}[x=1cm,y=1cm]
  \definecolor{colorAKANE}{RGB}{183,40,46}
  \definecolor{colorTOKIWA}{RGB}{0,123,67}
  \definecolor{colorUSUAI}{RGB}{40,148,200}
  \colorlet{plotPrimary}{colorAKANE}
  \colorlet{plotSecondary}{colorTOKIWA}
  \colorlet{plotTertiary}{colorUSUAI}
  \colorlet{plotAxis}{black!72}

  \tikzset{
    >=latex, font=\small, mark size=1.5pt, mark options={thin},
    line cap=butt, line join=round,
    common boundary/.style={draw=plotAxis,line width=0.95pt},
    first boundary/.style={draw=plotPrimary,line width=0.95pt},
    second boundary/.style={draw=plotTertiary,line width=0.95pt},
    outer cutoff edge/.style={draw=plotSecondary,line width=0.95pt,
      line cap=round},
    inner cutoff edge/.style={draw=black!62,line width=0.75pt,
      line cap=round},
    outer boundary support/.style={draw=plotSecondary,line width=1.2pt},
    inner boundary support/.style={draw=black!62,line width=1.2pt},
    figure label/.style={text=black!84,inner sep=0.8pt}
  }

  \pgfmathsetmacro{\figureWidth}{7.5}
  \pgfmathsetmacro{\figureHeight}{3.2}
  \pgfmathsetmacro{\frameBottom}{-0.28}
  \pgfmathsetmacro{\frameTop}{\frameBottom+\figureHeight}
  \pgfmathsetmacro{\edgeDx}{4.26}
  \pgfmathsetmacro{\edgeDy}{2.48}
  \pgfmathsetmacro{\edgeLength}{veclen(\edgeDx,\edgeDy)}
  \pgfmathsetmacro{\tipX}{(\figureWidth-\edgeLength)/2}
  \pgfmathsetmacro{\tipY}{0.34}
  \pgfmathsetmacro{\wedgeAngle}{atan2(\edgeDy,\edgeDx)}
  \def\commonStartRadius{2.30}

  \path[use as bounding box]
    (0,\frameBottom) rectangle (\figureWidth,\frameTop);
  \fill[fill=white]
    (0,\frameBottom) rectangle (\figureWidth,\frameTop);
  \coordinate (vertex) at (\tipX,\tipY);
  \coordinate (upperEnd) at ({\tipX+\edgeDx},{\tipY+\edgeDy});
  \coordinate (lowerEnd) at ({\tipX+\edgeLength},\tipY);
  \path (vertex) ++(\wedgeAngle:\commonStartRadius)
    coordinate (commonUpper);
  \path (vertex) ++(0:\commonStartRadius)
    coordinate (commonLower);

  \path (commonUpper) ++(\wedgeAngle:-1.00)
    coordinate (roundingUpperControl);
  \coordinate (roundingLowerControl) at ({\tipX+1.45},\tipY);
  \def\roundingPath{(commonUpper)
    .. controls (roundingUpperControl) and (roundingLowerControl)
    .. (commonLower)}
  \def\secondDomainPath{\roundingPath -- (lowerEnd)
    arc[start angle=0,end angle=\wedgeAngle,radius=\edgeLength cm]
    -- cycle}

  \coordinate (chiTwoStart) at ({\tipX+3.00},\tipY);
  \coordinate (chiOneStart) at ({\tipX+3.75},\tipY);
  \coordinate (chiTwoControlOne) at ({\tipX+3.05},{\tipY+0.63});
  \coordinate (chiTwoControlTwo) at ({\tipX+3.35},{\tipY+0.69});
  \coordinate (chiTwoShoulder) at ({\tipX+3.70},{\tipY+0.71});
  \coordinate (chiTwoEnd) at
    ({\tipX+\edgeLength+0.50},{\tipY+0.93});
  \coordinate (chiOneControlOne) at ({\tipX+3.80},{\tipY+0.34});
  \coordinate (chiOneControlTwo) at ({\tipX+4.05},{\tipY+0.38});
  \coordinate (chiOneShoulder) at ({\tipX+4.30},{\tipY+0.39});
  \coordinate (chiOneEnd) at
    ({\tipX+\edgeLength+0.50},{\tipY+0.50});
  \coordinate (cropBeyond) at ({\tipX+\edgeLength+0.50},\tipY);
  \def\cutoffPath{(chiTwoStart)
    .. controls (chiTwoControlOne) and (chiTwoControlTwo)
    .. (chiTwoShoulder) -- (chiTwoEnd)}
  \def\innerCutoffPath{(chiOneStart)
    .. controls (chiOneControlOne) and (chiOneControlTwo)
    .. (chiOneShoulder) -- (chiOneEnd)}

  \begin{scope}
    \clip (vertex) -- (lowerEnd)
      arc[start angle=0,end angle=\wedgeAngle,radius=\edgeLength cm]
      -- cycle;
    \fill[fill=plotPrimary!8]
      (vertex) -- (lowerEnd)
      arc[start angle=0,end angle=\wedgeAngle,radius=\edgeLength cm]
      -- cycle;
    \fill[fill=plotTertiary!7] \secondDomainPath;
    \fill[fill=plotSecondary!10]
      \cutoffPath -- (cropBeyond) -- cycle;
    \fill[fill=black!7]
      \innerCutoffPath -- (cropBeyond) -- cycle;
    \draw[outer cutoff edge] \cutoffPath;
    \draw[inner cutoff edge] \innerCutoffPath;
  \end{scope}

  \draw[first boundary] (commonUpper) -- (vertex) -- (commonLower);
  \draw[second boundary] \roundingPath;
  \draw[common boundary] (commonUpper) -- (upperEnd);
  \draw[common boundary] (commonLower) -- (lowerEnd);
  \fill[fill=plotPrimary] (vertex) circle[radius=1.45pt];
  \draw[outer boundary support] (chiTwoStart) -- (lowerEnd);
  \draw[inner boundary support] (chiOneStart) -- (lowerEnd);
  \fill[fill=plotSecondary] (chiTwoStart) circle[radius=1.45pt];
  \fill[fill=black!62] (chiOneStart) circle[radius=1.45pt];

  \node[figure label,text=plotPrimary]
    at ({\tipX+1.03},{\tipY+0.24}) {$\Omega_1$};
  \node[figure label,text=plotTertiary]
    at ({\tipX+2.75},{\tipY+0.75}) {$\Omega_2$};
  \node[figure label,text=plotSecondary]
    at ({\tipX+4.38},{\tipY+0.62}) {$\chi_2$};
  \node[figure label,text=black!62]
    at ({\tipX+4.38},{\tipY+0.19}) {$\chi/\chi_1$};
\end{tikzpicture}%
        \caption{Lemma~\ref{l:DomChange}, Lemma~\ref{l:sameNearInfinity}}
  \end{subfigure}
  \caption{Illustration for Lemma~\ref{l:longDistanceHarmonic}, Lemma~\ref{l:DomChange} and Lemma~\ref{l:sameNearInfinity}.}
  \label{f:4f2}
\end{figure}

\begin{lemma}
\label{l:longDistanceHarmonicRight}
Let $\Omega$ be an \acl domain. Let $m,l\in\mathbb{R}$ that $|m|<\rardyr_\Omega$, $|m-l|<\frac12$. Let $\gamma\in S^0(\partial\Omega)$ with $\supp\gamma\subset\partial\oreg$, $\chi\in S^0(\Omega)$ with $(\supp\chi\cap\partial\Omega)\subset\partial\oreg$. Let $\tilde\gamma\in S^0(\partial\Omega)$ that $\supp\tilde\gamma\subset\partial\oreg$ and $\supp\gamma\cap\supp(1-\tilde\gamma)=\emptyset$. Let $\varphi\in S^0(\partial\Omega)$ that $\supp\gamma\chi\cap\supp(1-\varphi)=\emptyset$, $\supp\varphi\subset\partial\oreg$. Then for all $k,k'\in \mathbb{N}_0$, vector fields $\mc{V}_1,\dots, \mc{V}_k, \mc{V}_1',\dots, \mc{V}_{k'}'\in\Diff^{1,1}_b(\Omega)$, for each $N>0$, we have
\begin{equation}
\| \chi\mc{V}_1'\dots\mc{V}_{k'}'\mathcal{H}_{M}^\Omega\gamma\mc{V}_1\dots \mc{V}_k f\|_{X_{m}^l(\Omega)}\leq C\|\tilde{\gamma} f\|_{H^{-N,m-\frac12+0}(\partial\Omega)}
+C\sum_{|\alpha|\leq k+k'}\| \varphi(\langle x\rangle D)^\alpha f\|_{H^{\frac{1}{2},l}(\partial\Omega)},
\end{equation}
uniformly for all $f\in\mc{S}'(\partial\Omega)$ and $M\in[0,M_0]$. See Figure~\ref{f:4f2}(A). 
\end{lemma}
\begin{proof}
We do not distinguish $\gamma\in S^0(\partial\Omega)$ from its extension $\gamma\in S^0(\Omega)$. Let $\chi_1,\chi_2,\chi_3,\chi_4\in S^0(\Omega)$ with $\supp\chi\cap\supp(1-\chi_1)=\emptyset$, $\supp\chi_j\cap\supp(1-\chi_{j+1})=\emptyset$ for $j=1,2,3$. Without loss of generality assume $\supp\gamma\chi_2\cap\supp(1-\varphi)=\emptyset$. 
By decomposing $\mc{H}_M^\Omega=\mc{H}_M^\Omega\chi_2+\mc{H}_M^\Omega(1-\chi_2)$, we have
\begin{equation}\label{e:4l6}
    \| \chi\mc{V}_1'\cdots\mc{V}_{k'}'\mathcal{H}_{M}^\Omega \gamma\mc{V}_1\cdots \mc{V}_k f\|_{X_m^l} \le C\| \chi_1\mathcal{H}_{M}^\Omega \gamma\chi_2\mc{V}_1\cdots \mc{V}_k f\|_{X_{m}^{l,k'}}+C\| \chi_1\mathcal{H}_{M}^\Omega \gamma(1-\chi_2)\mc{V}_1\cdots \mc{V}_k f\|_{X_{m}^{l,k'}}
\end{equation}
With that the support of $\chi_1$ and $1-\chi_2$ are disjoint, Lemma~\ref{l:4t2} implies
\begin{equation}
    \| \chi_1\mathcal{H}_{M}^\Omega\gamma(1-\chi_2)\mc{V}_1\cdots \mc{V}_k f\|_{X_{m}^{l,k'}}\le C\|\tilde\gamma f\|_{H^{-N,m-\frac12+0}_\partial},
\end{equation}
that the second term of \eqref{e:4l6} is an error term.

2. We iterate in $k\in\mathbb{N}_0$ to show for all $k'\in\mathbb{N}_0$, 
\begin{equation}\label{e:4l5}
\|\chi_1\mathcal{H}_M^\Omega\gamma\chi_2\mc{V}_1\cdots \mc{V}_k f\|_{X_{m}^{l,k'}}\le C\|\gamma f\|_{H^{-N,m-\frac12+0}_\partial}+C\sum_{|\alpha|\le k'+k}\|\varphi(\langle x\rangle D)^\alpha f\|_{H^{\frac12, l}_\partial}.
\end{equation}
Lemma~\ref{l:longDistanceHarmonic} and Lemma~\ref{l:4t2} give case $k=0$ and all $k'\in\mathbb{N}_0$. Let \eqref{e:4l5} holds for $k-1$ and all $k'\in\mathbb{N}_0$. Note 
\begin{equation}
(-\Delta+M^2)\gamma\chi_2\mc{V}_1\chi_3 \mc{H}_M^\Omega \tilde{\gamma}\chi_4\mc{V}_2\cdots \mc{V}_kf=\chi_3[-\Delta,\gamma\chi_2\mc{V}_1]\chi_3\mc{H}_M^\Omega\tilde{\gamma}\chi_4\mc{V}_2\cdots \mc{V}_kf,
\end{equation}
where $\gamma\chi_2\mc{V}_1 \chi_3\mc{H}_M^\Omega \tilde{\gamma}\chi_4\mc{V}_2\cdots \mc{V}_kf$ has trace $\gamma\chi_2 \mc{V}_1\cdots \mc{V}_k f$ on $\partial\Omega$. This implies
\begin{equation}
    \chi_1\mathcal{H}_M^\Omega\gamma\chi_2\mc{V}_1\dots \mc{V}_k f=\gamma\chi_1\mc{V}_1 \chi_3\mc{H}_M^\Omega \tilde{\gamma}\chi_4\mc{V}_2\cdots \mc{V}_kf-\chi_1\mathcal{R}_M^\Omega\chi_3[-\Delta,\gamma\chi_2\mc{V}_1]\chi_3\mc{H}_M^\Omega\tilde{\gamma}\chi_4\mc{V}_2\cdots \mc{V}_kf.
\end{equation}
Now note $\chi_1\mathcal{R}_M^\Omega\chi_3:Y_{m}^{l,k'}\rightarrow \dot X_{m}^{l,k'}$ and $[-\Delta, \gamma\chi_2\mc{V}_1]=\oxy: \dot X_{m}^{l,k'}\rightarrow Y_{m}^{l,k'}$ are bounded per Lemma~\ref{l:lowRegularityAction}, Lemma~\ref{l:oxybddxmlk} and Lemma~\ref{l:algCommute}. Apply the iterative assumption to $\chi_3,\chi_4$ instead of $\chi_1,\chi_2$ to see
\begin{multline}
\|\chi_1\mathcal{H}_M^\Omega\gamma\chi_2\mc{V}_1\cdots \mc{V}_k f\|_{X_{m}^{l,k'}}\le C\|\chi_3\mc{H}_M^\Omega \tilde{\gamma}\chi_4\mc{V}_2\cdots \mc{V}_kf\|_{X_{m}^{l,k'+1}}+C\|\chi_3\mc{H}_M^\Omega\tilde{\gamma}\chi_4\mc{V}_2\cdots \mc{V}_kf\|_{X_{m}^{l,k'}}\\
    \le C\|\tilde{\gamma} f\|_{H^{-N,m-\frac12+0}_\partial}+C\sum_{|\alpha|\le (k'+1)+(k-1)}\|\tilde\varphi(\langle x\rangle D)^\alpha f\|_{H^{\frac12, l}_\partial},
\end{multline}
for some $\tilde\varphi$ with slightly larger support than $\varphi$. Iterate and rename $\tilde\varphi$ to $\varphi$ to get \eqref{e:4l5} for all $k,k'\in\mathbb{N}_0$. Revisit \eqref{e:4l6} to conclude the proof.
\end{proof}

\begin{lemma}
Let $\Omega$ be an \acl domain. Let $m,l\in\mathbb{R}$ that $|m|<\rardyr_\Omega$, $|m-l|<\frac12$. Let $k,s\in\mathbb{N}_0$. Let $N\in\mathbb{N}_0$ be large. Let $\alpha$ be a multi-index. Let $\chi\in S^0(\Omega)$. Let $\gamma,\tilde\gamma\in S^0(\partial\Omega)$ with support inside $\partial\oreg$, and $\supp\gamma\cap\supp(1-\tilde\gamma)=\emptyset$. Let $\psi(\xi)\in S^0(\mathbb{R}^{d-1})$ with $0\notin \supp\psi$. Let $\varphi\in S^0(\Omega)$, with properties to be specified. Let $\mc{U}_1,\dots, \mc{U}_{k}\in \Diff^{1,0}_b(\Omega)$. Then the following are true uniformly for all $f\in \mc{S}'(\Omega)$, $\fbdy\in \mc{S}'(\partial\Omega)$:
\begin{enumerate}
\item If $\supp\varphi\cap\supp\chi=\emptyset$, then
\begin{equation}
\label{e:thing2}
\|\chi\rr^{l+k}\dd^{m-l}\mc{U}_1\cdots \mc{U}_{k}\nabla \mc{R}^\Omega_M\varphi f\|_{L^2(\Omega)}\leq C\|f\|_{Y^{l,-N}_{m}(\Omega)} .
\end{equation}
\item  We have
\begin{equation}
\label{e:thing3}
\|\chi\rr^{l+k}\dd^{m-l}\mc{U}_1\cdots \mc{U}_{k}\nabla \mc{R}^\Omega_M f\|_{L^2(\Omega)}\leq C\sum_{|\alpha|\leq k} \|(\rr \mc{W})^{\alpha} f\|_{Y^{l}_{m}(\Omega)}.
\end{equation}
\item If $(\supp\varphi\cap\partial\Omega)\cap \supp\gamma=\emptyset$, then
\begin{equation}
\label{e:thing6}
\|\varphi\rr^{l+k}\dd^{m-l}\mc{U}_1\cdots \mc{U}_{k}\nabla \mc{H}^\Omega_M\gamma \fbdy\|_{L^2(\Omega)}\le C\|\gamma \fbdy\|_{H^{-N,m-\frac12+0}(\partial\Omega)}.
\end{equation}
\item We have
\begin{equation}
\label{e:thing7}
\|\chi\rr^{l+k}\dd^{m-l}\mc{U}_1\dots \mc{U}_{k}\nabla \mc{H}^\Omega_M\gamma \fbdy\|_{L^2(\Omega)}\leq C\sum_{|\alpha|\leq k} \|(\langle x\rangle D)^{\alpha}\gamma \fbdy\|_{H^{\frac{1}{2},l}(\partial\Omega)}+C\|\gamma \fbdy\|_{H^{-N, m-\frac12+0}}.
\end{equation}
\item We have
\begin{equation}\label{e:thing4}
    \|\mc{H}^\Omega_M\gamma D^\alpha \fbdy\|_{X_{m,s}^{l,-|\alpha|}(\Omega)}\leq C \|\tilde\gamma\fbdy\|_{H^{s+\frac{1}{2},l-|\alpha|}(\partial\Omega)}.
\end{equation}
\item If $(\supp\varphi\cap\partial\Omega)\cap \supp\gamma=\emptyset$, then
\begin{equation}
\label{e:thing5}
\|\varphi \rr^{k}\mc{U}_1\cdots \mc{U}_{k}\mc{H}^\Omega_M\gamma D^\alpha \fbdy\|_{X^{l}_m(\Omega)}\leq C\|\tilde\gamma \fbdy\|_{H^{-N,m-\frac12-|\alpha|+0}(\partial\Omega)}.
\end{equation}

\item We have
\begin{equation}
\label{e:thing1}
    \|\mc{H}_M^\Omega\gamma \psi(D)\tilde\gamma \fbdy\|_{X_{m-2k, s}^{l-2k, 0}}\le C\|\tilde\gamma \fbdy\|_{H^{s+\frac12, l-2k}(\partial\Omega)}.  
\end{equation}

\item If $(\supp \varphi\cap \partial\Omega)\cap \supp \gamma =\emptyset$, then

\begin{equation}
\label{e:thing8}
\|\varphi\mc{H}_M^\Omega\gamma \psi(D)\tilde\gamma \fbdy\|_{X^{l,N}_{m,N}(\Omega)}\le C\|\tilde\gamma \fbdy\|_{H^{-N,-N}(\partial\Omega)}.
\end{equation}

\item If $\supp\varphi$ is compact with $(\supp \varphi\cap \partial\Omega)\cap \supp \gamma=\emptyset$, then
\begin{equation}
\label{e:thing8comp}
\|\varphi\mc{H}_M^\Omega\gamma \psi(D)\tilde\gamma \fbdy\|_{H^{\frac32-0,N}(\Omega)}\leq C\|\tilde\gamma\fbdy\|_{H^{-N,-N}(\partial\Omega)}.
\end{equation}
\item If $(\supp \varphi\cap \partial\Omega)\cap \supp \gamma=\emptyset$, then
\begin{equation}\label{e:thing9}
    \|\gamma\psi(D)\gamma\partial_\nu \mathcal{R}_M^\Omega\varphi f\|_{H^{N,N}(\partial\Omega)}\le C\|\varphi f\|_{Y_{-m, -N}^{-l, -N}(\Omega)}.
\end{equation}
\end{enumerate}
\end{lemma}
\begin{proof}
We set up some cutoffs. Let $\tilde\chi$, $\tilde\varphi$, $\tilde{\tilde{\gamma}}$ satisfy the same conditions for $\chi$, $\varphi$, $\tilde{\gamma}$ with $\supp\chi\cap\supp(1-\tilde\chi)=\supp\varphi\cap\supp(1-\tilde\varphi)=\supp\tilde\gamma\cap\supp(1-\tilde{\tilde{\gamma}})=\emptyset$. We confuse $\gamma,\tilde\gamma,\tilde{\tilde\gamma}\in S^0(\partial\Omega)$ with their extensions in $S^0(\Omega)$ that $\supp\gamma\cap\supp(1-\tilde\gamma)=\supp\tilde\gamma\cap\supp(1-\tilde{\tilde{\gamma}})=\emptyset$. 

(1) Note $\tilde\chi$ and $\tilde\varphi$ have disjoint support. Lemma~\ref{l:smoothing} implies
\begin{equation}
    \|\chi\rr^{l+k}\dd^{m-l}\mc{U}_1\cdots \mc{U}_{k}\nabla \mc{R}^\Omega_M\varphi f\|\le C\|\tilde\chi \mathcal{R}_M^\Omega \tilde\varphi\varphi f\|_{X_{m}^{l,k}}\le C\|\varphi f\|_{Y_{m}^{l,-N}}.
\end{equation}

(2) Lemma~\ref{l:longDistanceResolve} implies
\begin{equation}
    \|\chi\rr^{l+k}\dd^{m-l}\mc{U}_1\cdots \mc{U}_{k}\nabla \mc{R}^\Omega_M f\| \le C\|f\|_{Y_{m}^{l,k}}\le C\sum_{|\alpha|\leq k} \|(\rr \mc{W})^{\alpha} f\|_{Y_m^l}.
\end{equation}

(3) Lemma~\ref{l:4t2} implies
\begin{equation}
    \|\varphi \rr^{l+k}\dd^{m-l}\mathcal{U}_1\cdots \mathcal{U}_k \nabla \mathcal{H}_M^\Omega \gamma \fbdy\|\le \|\varphi \rr^k\mathcal{U}_1\cdots \mathcal{U}_k \mathcal{H}_M^\Omega \tilde\gamma\gamma \fbdy\|_{X^{l}_m}\le C\|\gamma \fbdy\|_{H^{-N, m-\frac12+0}_\partial}. 
\end{equation}

(4) Lemma~\ref{l:longDistanceHarmonicRight} implies
\begin{multline}
    \|\chi \rr^{l+k}\dd^{m-l}\mathcal{U}_1\cdots \mathcal{U}_k \nabla \mathcal{H}_M^\Omega \gamma \fbdy\|\le \|\chi \rr^k\mathcal{U}_1\cdots \mathcal{U}_k \mathcal{H}_M^\Omega \tilde\gamma\gamma \fbdy\|_{X^{l}_m}\le C\sum_{|\alpha|\leq k} \|(\langle x\rangle D)^{\alpha}\gamma \fbdy\|_{H^{\frac{1}{2},l}_\partial}\\
    +C\|\gamma \fbdy\|_{H^{-N, m-\frac12+0}}. 
\end{multline}

(5) Let $s\in\mathbb{N}_0$. We claim that by induction in $k\in\mathbb{N}_0$,
\begin{equation}
    \mathcal{H}_M^\Omega\gamma\mathcal{U}_1\cdots \mathcal{U}_k: H^{s+\frac12, l-k}_{\partial}\rightarrow X^{l,-k}_{m,s}
\end{equation}
is bounded. Lemma~\ref{l:longDistanceHarmonic} and Proposition~\ref{p:harmonic} implies
\begin{equation}
    \|\mc{H}_M^\Omega\gamma\fbdy\|_{X_{m,s}^l}\le C\sum_{|\beta|\le s}\| D^\beta \gamma\fbdy\|_{H^{\frac{1}{2},l}_\partial}\le C\|\gamma\fbdy\|_{H^{s+\frac12,l}_\partial}, 
\end{equation}
and that the case $k=0$ holds. Assume that the cases up to $k-1$ hold, then consider that
\begin{equation}\label{e:4l24}
    \mathcal{H}_M^\Omega\gamma\mathcal{U}_1\cdots \mathcal{U}_k=\mathcal{H}_M^\Omega\gamma \rr\mathcal{U}_1 \mathcal{U}_2 \cdots \mathcal{U}_k \rr^{-1} 
    -
    \mathcal{H}_M^\Omega\gamma [\rr, \mathcal{U}_1\cdots \mathcal{U}_k]\rr^{-1},
\end{equation}
the second term of which is desirably bounded from $H^{s+\frac12, l-k}_\partial$ to $X_{m,s}^{l,-k+1}$ because of the inductive assumption and that $[\rr, \mathcal{U}_1\dots \mathcal{U}_k]\in \Diff_b^{k-1,0}$. For the first term, since $\mathcal{U}_1$ is tangential to $\partial\Omega$, 
\begin{equation}
    \gamma \rr \mathcal{U}_1 \mathcal{H}_{M}^\Omega \tilde\gamma\mathcal{U}_2 \cdots \mathcal{U}_k \rr^{-1}\fbdy= \mathcal{H}_M^\Omega\gamma \rr\mathcal{U}_1 \mathcal{U}_2 \cdots \mathcal{U}_k \rr^{-1} \fbdy \ \text{ on }\partial\Omega
\end{equation}
and thus the extension correction formula \eqref{e:hForm} gives
\begin{equation}\label{e:4l23}
    \mathcal{H}_{M}^\Omega \gamma \rr \mathcal{U}_1 \tilde\gamma \mathcal{U}_2 \cdots \mathcal{U}_k \rr^{-1} =\gamma \rr \mathcal{U}_1 \mathcal{H}_{M}^\Omega \tilde\gamma\mathcal{U}_2 \cdots \mathcal{U}_k \rr^{-1} -\mathcal{R}_M^\Omega [-\Delta, \gamma \rr\mathcal{U}_1] \mathcal{H}_{M}^\Omega \tilde\gamma \mathcal{U}_2 \cdots \mathcal{U}_k \rr^{-1}.
\end{equation}
Now note that by the inductive assumption, the first term here $\gamma \rr\mathcal{U}_1\mathcal{H}_{M}^\Omega \tilde\gamma \mathcal{U}_2 \cdots \mathcal{U}_k \rr^{-1}$ is bounded from $H^{s+\frac12, l-k}_\partial$ to $X_{m,s}^{l,-k}$. Now note $[-\Delta, \gamma \rr\mathcal{U}_1]=\oxy$ maps $X_{m,s}^{l,-k}$ to $Y_{m,s}^{l,-k}$ and $\mathcal{R}_M^\Omega: Y_{m,s}^{l,-k}\rightarrow \dot X_{m,s}^{l,-k}$ per Lemma \ref{l:lowRegularityAction}, we have the second term of \eqref{e:4l23} being bounded from $H^{s+\frac12, l-k}_\partial$ to $\dot X_{m,s}^{l,-k}$. Revisit \eqref{e:4l24} to observe the case $k$ holds and this concludes the induction. 

(6) Note $\varphi\gamma=0$ and 
\begin{equation}
    \tilde\varphi \mc{H}^\Omega_M\gamma D^\alpha \fbdy=\tilde\varphi \mc{H}^\Omega_M\gamma (\rr D)^\alpha \rr^{-|\alpha|}\fbdy+\tilde\varphi \mc{H}^\Omega_M\gamma O(\Diff^{|\alpha|-1, |\alpha|-1}_b)\rr^{-|\alpha|}\fbdy.
\end{equation}
Applying Lemma~\ref{l:4t2} on each term then gives
\begin{equation}
    \|\varphi \rr^{k}\mc{U}_1\cdots \mc{U}_{k}\mc{H}^\Omega_M\gamma D^\alpha \fbdy\|_{X_{m}^l}\le C\|\tilde\gamma\rr^{-|\alpha|}\fbdy\|_{H^{-N,m-\frac12+0}_\partial}.
\end{equation}

(7) Consider that
\begin{equation}
\mathcal{H}_M^\Omega\gamma\psi(D)\tilde{\tilde{\gamma}}=\mathcal{H}_M^\Omega\gamma (-\Delta)^k\tilde\gamma (-\Delta)^{-k}\psi(D)\tilde{\tilde{\gamma}}.
\end{equation}
Here $\tilde\gamma (-\Delta)^{-k}\psi(D)\tilde{\tilde{\gamma}}$ is bounded from $H^{s+\frac12, l-2k}(\partial\Omega)$ to $H^{s+\frac12+2k, l-2k}(\partial\Omega)$. Now invoke part (5) to see $\mathcal{H}_M^\Omega\gamma\psi(D)\tilde\gamma \fbdy\in X_{m,s+2k}^{l,-2k}(\Omega)\subset X_{m-2k,s}^{l-2k,0}(\Omega)$ provided $\tilde\gamma \fbdy\in H^{s+\frac12, l-2k}(\partial\Omega)$.

(8) Consider that
\begin{equation}
\varphi\mathcal{H}_M^\Omega\gamma\psi(D)\tilde{\tilde{\gamma}}=\varphi\mathcal{H}_M^\Omega\gamma (-\Delta)^k\tilde\gamma (-\Delta)^{-k}\psi(D)\tilde{\tilde{\gamma}}. 
\end{equation}
Given $\tilde\gamma \fbdy\in H^{-N, -N}(\partial\Omega)$, we have $\tilde\gamma (-\Delta)^{-k}\psi(D)\tilde\gamma \fbdy\in H^{2k-N, -N}_\partial$. Take $k$ sufficiently large that $l-2k\le-N$. Part (6) implies that $\varphi\mathcal{H}_M^\Omega\gamma\psi(D)\tilde\gamma \fbdy\in X^{l, 2N}_{m, 0}\subset X^{l, N}_{m, N}$. 

(9) Note $\varphi\mathcal{H}_M^\Omega \gamma\psi(D)\tilde\gamma \fbdy|_{\partial\Omega}=0$. This implies
\begin{equation}
    \varphi\mathcal{H}_M^\Omega \gamma\psi(D)\tilde\gamma \fbdy=
    \mathcal{R}_M^\Omega [-\Delta, \varphi]\tilde\varphi\mathcal{H}_M^\Omega \gamma\psi(D)\tilde\gamma \fbdy
\end{equation}
and thus part (8) implies
\begin{equation}
    \|\tilde\varphi\mathcal{H}_M^\Omega \gamma\psi(D)\tilde\gamma \fbdy\|_{H^{\frac12}}\le \|\tilde\varphi\mathcal{H}_M^\Omega \gamma\psi(D)\tilde\gamma \fbdy\|_{X^{0}_{0}}\le \|\tilde\gamma \fbdy\|_{H^{-N,-N}_\partial},
\end{equation}
with $[-\Delta, \varphi]\tilde\varphi\mathcal{H}_M^\Omega \gamma\psi(D)\tilde\gamma \fbdy\in Y^0_0$. Note $\mathcal{R}_M^\Omega: Y^0_0\rightarrow \dot X_0^0$ and thus Lemma~\ref{l:nearBoundary} implies
\begin{multline}
    \|\varphi\mathcal{H}_M^\Omega \gamma\psi(D)\tilde\gamma \fbdy\|_{H^{\frac32-0}}\le C\|\rr^{-N}\psi_{\dd}(\dd)\mathcal{R}_M^\Omega [-\Delta, \varphi]\tilde\varphi\mathcal{H}_M^\Omega \gamma\psi(D)\tilde\gamma \fbdy\|_{H^{\frac32-0}}\\\le C\|[-\Delta, \varphi]\tilde\varphi\mathcal{H}_M^\Omega \gamma\psi(D)\tilde\gamma \fbdy\|_{H^{-\frac12-0}}+C\|\mathcal{R}_M^\Omega [-\Delta, \varphi]\tilde\varphi\mathcal{H}_M^\Omega \gamma\psi(D)\tilde\gamma \fbdy\|_{X^0_0}\le C\|\tilde\gamma \fbdy\|_{H^{-N,-N}_\partial}
\end{multline}
as desired, for some $\psi_{\dd}\in C^\infty_c(\mathbb{R})$ that $\psi_{\dd}(\dd)\varphi=\varphi$.

(10) Part (8) implies for any $\phi_\partial\in \mc{S}'(\partial\Omega)$, we have $\tilde\varphi\mc{H}_M^\Omega\gamma \psi(D)\gamma \phi_\partial\in\dot X_{m,N}^{l,N}$ and 
\begin{equation}
    \|\tilde\varphi\mc{H}_M^\Omega\gamma \psi(D)\gamma \phi_\partial\|_{X_{m,N}^{l,N}}\le C\|\gamma\phi_\partial\|_{H^{-N,-N}_\partial}.
\end{equation}
By noting the duality $(\mathcal{H}_M^\Omega)^*=-\partial_\nu \mathcal{R}_M^\Omega$, we have 
\begin{equation}
    |\langle\gamma\psi(D)\gamma\partial_\nu \mathcal{R}_M^\Omega\varphi f, \phi_\partial\rangle_{\partial\Omega}|=|\langle \varphi f, \tilde\varphi\mc{H}_M^\Omega\gamma \psi(D)^*\gamma \phi_\partial\rangle|\le C\|\varphi f\|_{Y_{-m, -N}^{-l,-N}}\|\phi_\partial\|_{H^{-N,-N}_\partial},
\end{equation}
yielding the desired bounds. 
\end{proof}

\subsection{The \dtn map on \acl domains}
\label{s:parametrix}
Next, we record local elliptic estimates for $\Lambda_M^\Omega$.
\begin{lemma}[Elliptic estimate for $\Lambda^\Omega_M$ near the Lipschitz boundary]\label{l:EllEst_ACLip}
Let $\Omega$ be an \acl domain. Let $l\in\mathbb{R}$ that $|l|<\rardyr_\Omega+\frac12$. Let $\omega>0$. Let $\chi,\tilde\chi,\tilde{\tilde{\chi}} \in C_c^\infty(\partial\Omega)$ that $\supp\chi\cap\supp(1-\tilde\chi)=\supp\tilde\chi\cap\supp(1-\tilde{\tilde{\chi}})=\emptyset$, $\supp(1-\chi)\subset\partial\oreg$. Let $\psi_0(\xi)\in C^\infty_c(\mathbb{R}^{d-1})$ that $0\notin\supp(1-\psi_0)$. For each $0\le s\le 1$, $|l|<\rardyr_{\Omega}+\frac12$, and $N\in\mathbb{R}$, we have
\begin{multline}
    \|\chi u\|_{H^{s}(\partial \Omega)}\leq C\|\tilde\chi (\Lambda_M^\Omega-\sqrt{\omega^2+M^2}) u\|_{H^{s-1}(\partial \Omega)}+C\|\tilde{\tilde{\chi}} u\|_{H^{0}(\partial \Omega)}+C\|(1-\tilde\chi) u\|_{H^{-N,-N}(\partial \Omega)}\\
    +C\|(1-\tilde{\tilde{\chi}})\psi_0(D)(1-\tilde\chi)u\|_{H^{-N,l}(\partial\Omega)}.
\end{multline}
uniformly for all $u\in \mc{S}'(\partial\Omega)$ and $M\in[0,M_0]$. 
\end{lemma}
\begin{proof}
We confuse the cutoffs with their extension onto $\Omega$ with 0 Neumann data. By \eqref{e:4l28} from Lemma~\ref{l:ellregunbdd}, we have
\begin{multline}
    \|\chi \mc{H}^\Omega_M u\|_{H^{s+\frac12}}\le C\|\tilde\chi(\Lambda_M^\Omega-\sqrt{\omega^2+M^2}) u\|_{H^{s-1}_\partial}+C\|\tilde\chi\mc{H}^\Omega_M\big(\tilde{\tilde{\chi}}+(1-\tilde{\tilde{\chi}})\psi_0(D)(1-\tilde\chi)\big) u\|_{H^{s-\frac12}}
    \\+C\|\tilde\chi u\|_{H^{s-1}_\partial}
    +C\|\tilde\chi \mc{H}^\Omega_M(1-\tilde{\tilde{\chi}})(1-\psi_0(D))(1-\tilde\chi) u\|_{H^{s-\frac12}}. 
\end{multline}
It remains to estimate the last two terms. Lemma~\ref{l:ellregunbdd} implies
\begin{equation}
    \|\tilde\chi \mc{H}^\Omega_M\big(\tilde{\tilde{\chi}}+(1-\tilde{\tilde{\chi}})\psi_0(D)(1-\tilde\chi)\big) u\|_{H^{s-\frac12}}
    \le C\|\tilde{\tilde{\chi}} u\|_{H^{0}_\partial}+C\|(1-\tilde{\tilde{\chi}})\psi_0(D)(1-\tilde\chi)u\|_{H^{-N,l}_\partial}+C\|(1-\tilde\chi)u\|_{H^{-N,-N}_\partial}.
\end{equation}
On another hand, \eqref{e:thing8comp} implies
\begin{equation}
    \|\tilde\chi \mc{H}^\Omega_M(1-\tilde{\tilde{\chi}})(1-\psi_0(D))(1-\tilde\chi) u\|_{H^{s-\frac12}}
    \le C\|(1-\tilde\chi)u\|_{H^{-N,-N}_\partial}. 
\end{equation}
This concludes the proof. 
\end{proof}

\begin{lemma}[Change of domain for $\Lambda_M^\Omega$]\label{l:DomChange}
Let $\Omega_1, \Omega_2\subset \mathbb{R}^d$ be \acl domains. Let $\chi\in C^\infty(\Omega_1)$ that $\supp\chi\subset\Omega_1\cap \Omega_2$, $\Omega_1, \Omega_2$ share the common boundary near $\supp\chi$ and $\partial_\nu
\chi\equiv 0$. Then we have
\begin{equation}
\label{e:genCommutator}
\chi\Lambda^{\Omega_1}_M= \Lambda_M^{\Omega_2} \chi+ \tilde{\mc{E}}_M, \ \tilde{\mc{E}}_M=\partial_\nu \mc{R}_M^{\Omega_2}[-\Delta, \chi]\mc{H}^{\Omega_1}_M,
\end{equation}
uniformly for all $M\in [0,M_0]$. See Figure~\ref{f:4f2}(B). 
\end{lemma}
\begin{proof}
Noting $\chi \mc{H}_M^{\Omega_1}u$ extends $\chi u$ into $\Omega_2$, we have
\begin{equation}
    \mc{H}^{\Omega_2}_M\chi u=\chi \mc{H}^{\Omega_1}_Mu-\mc{R}^{\Omega_2}_M(-\Delta+M^2)\chi \mc{H}^{\Omega_1}_Mu=\chi \mc{H}^{\Omega_1}_Mu-\mc{R}^{\Omega_2}_M[-\Delta,\chi] \mc{H}^{\Omega_1}_Mu,
\end{equation}
by the extension correction formula \eqref{e:hForm}. Thus
\begin{equation}
    \Lambda^{\Omega_2}_M\chi u=\partial_\nu\mc{H}^{\Omega_2}_M\chi u=\chi\Lambda^{\Omega_1}_M u-\partial_\nu \mc{R}^{\Omega_2}_M[-\Delta,\chi] \mc{H}^{\Omega_1}_Mu,
\end{equation}
as desired. 
\end{proof}

We can now prove our main result on the structure of the \dtn map near infinity on an \acl domain. Given $\psi(\xi)\in S^0(\mathbb{R}^{d-1})$, we interpret $\psi(D)$ as a Fourier multiplier in $\Psi^{0,0}(\mathbb{R}^{d-1})$. Given $\chi,\tilde\chi\in S^0(\partial\oreg)$, we have $\chi \psi(D) \tilde\chi\in \Psi^{0,0}(\partial\oreg)$. 

\begin{lemma}[Comparing $\Lambda_M^\Omega$ near infinity]
\label{l:sameNearInfinity}
Let $\Omega_1,\Omega_2\subset\mathbb{R}^d$ be \acl domains. Let $\chi_1,\chi_2\in S^0(\partial\Omega_1\cap\partial\Omega_2)$ with $\supp\chi_1\cap\supp(1-\chi_2)=\emptyset$. Assume $\Omega_1, \Omega_2$ share the common smooth boundary near $\supp\chi_2$. Then we have
\begin{equation}
    \chi_1\Lambda^{\Omega_1}_M= \chi_1\Lambda_M^{\Omega_2}\chi_2+ \chi_1\tilde{\mathcal{E}}_M,
\end{equation}
where we confuse $\chi_1,\chi_2$ with extensions specified in the proof to write
\begin{equation}\label{e:4l42}
    \chi_1 \tilde{\mc{E}}_M=\chi_1\partial_\nu \mc{R}_M^{\Omega_2}[-\Delta, \chi_2]\mc{H}^{\Omega_1}_M.
\end{equation}
Moreover, let $\psi(\xi)\in S^0(\mathbb{R}^{d-1})$ with $0\notin \supp \psi$, $N\geq 0$, $l_j\in (-\rardyr_{\Omega_j}-\frac12,\rardyr_{\Omega_j}+\frac12)$ for $j=1,2$ such that
\begin{equation}\label{e:4l34}
    l_2< l_1+1 , \ l_2<\rardyr_{\Omega_1}+\frac12, \ l_1>-\rardyr_{\Omega_2}-\frac12.
\end{equation}
We can estimate the error $\tilde{\mc{E}}_M$ by
\begin{gather}\label{e:up1}
 \|\chi_1\psi(D) \chi_1\tilde{\mathcal{E}}_M\chi_1\psi(D)\chi_1 f\|_{H^{N,N}(\partial\Omega_1)}\leq C\|f\|_{H^{-N,-N}(\partial\Omega_1)} \\\label{e:up2}
  \|\chi_1\psi(D)\chi_1\tilde{\mathcal{E}}_Mf\|_{H^{N,N}(\partial\Omega_1)}\leq C\|f\|_{H^{0,l_1}(\partial\Omega_1)}, \\\label{e:up3}
    \|\chi_1\tilde{\mathcal{E}}_M\chi_1\psi(D)\chi_1 f\|_{H^{N,l_2}(\partial\Omega_1)}\leq C\|f\|_{H^{-N,-N}(\partial\Omega_1)},\\\label{e:up4}
     \|\chi_1\tilde{\mathcal{E}}_M f\|_{H^{N,l_2}(\partial\Omega_1)}\leq C\|f\|_{H^{0,l_1}(\partial\Omega_1)},
\end{gather}
uniformly for all $f\in \mc{S}'(\partial\Omega_1)$ and $M\in [0,M_0]$. See Figure~\ref{f:4f2}(B). 
\end{lemma}
\begin{proof}
0. Throughout the proof, we let $m_j$ that $|m_j|<\rardyr_{\Omega_j}$, $|m_j-l_j|<\frac12$ for $j=1, 2$. We further assume without loss of generality that $m_2\le m_1$, whose existence is guaranteed by \eqref{e:4l34}. Note $\Omega_1\cap\Omega_2$ admits a conic subset $U$ that $\supp\chi_2\subset U$. We confuse $\chi_2$ with its extension in $S^0(U)$ that $\partial_{\nu}\chi_2=0$. 

1. We have by~\eqref{e:genCommutator}
\begin{equation}
\chi_2\Lambda^{\Omega_1}_M=\Lambda_M^{\Omega_2}\chi_2+\partial_{\nu}\mc{R}_M^{\Omega_2}[-\Delta,\chi_2]\mc{H}_M^{\Omega_1}.
\end{equation}
Hence
\begin{equation}
\chi_1\Lambda^{\Omega_1}_M=\chi_1\Lambda_M^{\Omega_2}\chi_2+\chi_1\partial_\nu \mc{R}_M^{\Omega_2}[-\Delta,\chi_2]\mc{H}_M^{\Omega_1}=\chi_1\Lambda_M^{\Omega_2}\chi_2+\chi_1\tilde{\mathcal{E}}_M.
\end{equation}
It remains to establish the four desired estimates for 
\begin{equation}
    \chi_1\tilde{\mathcal{E}}_{M}=\chi_1\partial_\nu \mc{R}_M^{\Omega_2}\varphi[-\Delta,\chi_2]\varphi\mc{H}_M^{\Omega_1},
\end{equation}
where $\varphi\in S^0(U)$ with $\varphi\equiv 1$ near $\supp \nabla\chi_2$, $\supp \varphi\cap \supp \chi_1=\emptyset$.

2. We show \eqref{e:up1}. By \eqref{e:thing8} and \eqref{e:thing9}, we have
\begin{equation}
    \varphi\mc{H}^{\Omega_1}_M\chi_1\psi(D)\chi_1: H^{-N, -N}_\partial\rightarrow X^{l_1}_{m_1}(\Omega_1), \quad \chi_1\psi(D)\chi_1\partial_\nu \mc{R}_M^{\Omega_2}\varphi: \dot Y_{m_2}^{l_2}(\Omega_2)\rightarrow H^{N,N}_\partial,
\end{equation}
bounded uniformly in $M$. Lemma~\ref{l:oxybddxmlk} implies $[-\Delta, \chi_2]: X^{l_1}_{m_1}(\Omega_1)\rightarrow \varphi \dot Y^{l_1+1}_{m_1}(\Omega_1)\subset \dot Y^{l_2}_{m_2}(\Omega_2)$ and this concludes \eqref{e:up1}. Note here we used $l_2\le l_1+1$ and $m_2\le m_1$.   

3. For \eqref{e:up2}, it suffices to show $[-\Delta,\chi_2]\varphi\mc{H}_M^{\Omega_1}: H^{0,l_1}_\partial\rightarrow \dot Y^{l_1+1}_{m_1}(\Omega_1)$ is bounded. Take $\kappa\in C^\infty_c(\overline{\Omega_1})$ that $\partial\Omega_1$ is smooth over $\supp(1-\kappa)$. Let $\kappa_1,\tilde \kappa$ satisfy the same requirement for $\kappa$ while $\supp\kappa\cap\supp(1-\kappa_1)=\supp\kappa_1\cap\supp(1-\tilde\kappa)=\emptyset$. For $f\in H^{0,l_1}_\partial$, Lemma~\ref{l:ellregunbdd} implies
\begin{equation}\label{e:4l43}
    \|[-\Delta,\chi_2]\varphi\tilde\kappa\mc{H}_M^{\Omega_1}f\|_{H^{-1}(\Omega_1)}\le C\|\varphi\tilde\kappa \mathcal{H}_M^{\Omega_1}f\|_{H^{\frac12}(\Omega_1)}\le C\|f\|_{H^{0,l_1}_\partial}.
\end{equation}
To estimate $[-\Delta,\chi_2]\varphi(1-\tilde\kappa)\mc{H}_M^{\Omega_1}f$, we decompose it into
\begin{equation}
    [-\Delta,\chi_2]\varphi(1-\tilde\kappa)\mc{H}_M^{\Omega_1}f=[-\Delta,\chi_2]\varphi(1-\tilde\kappa)\mc{H}_M^{\Omega_1}(1-\kappa) f+[-\Delta,\chi_2]\varphi(1-\tilde\kappa)\mc{H}_M^{\Omega_1}\kappa f.
\end{equation}
Lemma~\ref{l:oxybddxmlk} and Lemma~\ref{l:hnegbreg} allow us to bound the first term by 
\begin{equation}
    \|[-\Delta,\chi_2]\varphi(1-\tilde\kappa)\mc{H}_M^{\Omega_1}(1-\kappa) f\|_{Y^{l_1+1}_{m_1}}\le C\|\mc{H}_M^{\Omega_1}(1-\kappa) f\|_{X^{l_1,0}_{m_1,-1}}\le C\|f\|_{H^{-\frac12, l_1}_\partial}.
\end{equation}
Noting $(1-\tilde\kappa)\mc{H}_M^{\Omega_1}\kappa=(1-\tilde\kappa)\mc{R}_M^{\Omega_1}[\Delta, \kappa_1]\mc{H}_M^{\Omega_1}\kappa$ we have
\begin{multline}
    \|[-\Delta,\chi_2]\varphi(1-\tilde\kappa)\mc{H}_M^{\Omega_1}\kappa f\|_{Y^{l_1+1}_{m_1}}\le C\|(1-\tilde\kappa)\mc{R}_M^{\Omega_1}[\Delta, \kappa_1]\mc{H}_M^{\Omega_1}\kappa f\|_{X^{l_1,0}_{m_1,-1}}\\\le C\|[\Delta, \kappa_1]\mc{H}_M^{\Omega_1}\kappa f\|_{H^{-1}}
    \le C\|f\|_{H^{0, l_1}_\partial}
\end{multline}
again by Lemma~\ref{l:ellregunbdd}. The above give $[-\Delta,\chi_2]\varphi\mc{H}_M^{\Omega_1}f\in \dot{Y}^{l_1+1}_{m_1}(\Omega_1)$ as desired. One concludes \eqref{e:up2} like in Step 2. 

4. For \eqref{e:up3}, Lemma~\ref{l:4t2} implies $\varphi\mc{H}_M^{\Omega_2}\chi_1: H^{-N, -l_2}_\partial\rightarrow X^{-l_2}_{-m_2}(\Omega_2)$ is bounded,
and thus its adjoint $-\chi_1\partial_\nu \mathcal{R}_M^{\Omega_2}\varphi: \dot Y^{l_2}_{m_2}(\Omega_2)\rightarrow H^{N, l_2}_\partial$ is bounded. Equation \eqref{e:up4} follows similarly. 
\end{proof}

\begin{lemma}[$M$-dependency of error]\label{l:errordependsonM}
Denote $\rardy=\min(\rardyr_{\Omega_1}, \rardyr_{\Omega_2})$. Let $|l|< \rardy+\frac12$, $|l'|< \rardy+\frac12$ with $l'\ge l$. In 
the setting of Lemma~\ref{l:sameNearInfinity}, we can estimate the error $\tilde{\mc{E}}_M$ by
\begin{gather}
    \|\chi_1\psi(D)\chi_1(\tilde{\mc{E}}_M-\tilde{\mc{E}}_0)\chi_1\psi(D)\chi_1f\|_{H^{N, N}(\partial\Omega_1)}\le CM^{\min(2\rardy-0,2)}\|f\|_{H^{-N, -N}(\partial\Omega_1)}\\
    \|\chi_1\psi(D)\chi_1(\tilde{\mc{E}}_M-\tilde{\mc{E}}_0)f\|_{H^{N, N}(\partial\Omega_1)}\le CM^{\min(2\rardy-0,\rardy+\frac12+l'-0,2)}\|f\|_{H^{0, l'}(\partial\Omega_1)}\\
    \|\chi_1(\tilde{\mc{E}}_M-\tilde{\mc{E}}_0)\chi_1\psi(D)\chi_1f\|_{H^{N, l}(\partial\Omega_1)}\le CM^{\min(2\rardy-0,\rardy+\frac12-l-0,2)}\|f\|_{H^{-N, -N}(\partial\Omega_1)}\\
    \|\chi_1(\tilde{\mc{E}}_M-\tilde{\mc{E}}_0)f\|_{H^{N, l}(\partial\Omega_1)}\le CM^{\min(2\rardy-0,\rardy+\frac12+\min(l',-l)-0, 1+l'-l-0,2)}\|f\|_{H^{0, l'}(\partial\Omega_1)}\label{e:4l53}
\end{gather}
uniformly for all $f\in \mc{S}'(\partial\Omega_1)$ and $M\in [0,M_0]$.
\end{lemma}
\begin{remark}
    If one sets $l=-\rardy-\frac12+0$, $l'=\rardy+\frac12-0$, all four $M$-terms are simplified to $M^{\min(2\rardy-0,2)}$. 
\end{remark}
\begin{proof}
    0. In this proof, notations are inherited from Lemma~\ref{l:sameNearInfinity} by making $l_1=l_2=l$, and all compatibility conditions there are satisfied. We have some $\epsilon>0$ and some $m_1'>m_1+\epsilon$, $m_2'>m_2+\epsilon$ such that $m_2'\le m_1'$ and $|m_1'|<\rardyr_{\Omega_1}$, $|m_1'-l'|<\frac12, |m_2'|<\rardyr_{\Omega_2}$, $|m_2'-l'|<\frac12$ as well. Note that given fixed $l,l'$, one can always choose $m_1, m_2, m_1', m_2'$ that $\e=\min(2\rardy-0,\rardy+\frac12+\min(l',-l)-0, 1+l'-l-0,2)$: this is similar to Proposition~\ref{p:harmonic}. We confuse $\chi_1$ with its extension $\chi_1\in S^0(\Omega_1\cap\Omega_2)$ be that $\partial_\nu \chi_1\equiv 0$. Let $\tilde\chi_1\in S^0(\Omega_1\cap\Omega_2)$ that $\supp\chi_1\cap\supp(1-\tilde\chi_1)=\emptyset$, $\supp\tilde\chi_1\cap\supp\varphi=\emptyset$. Decompose
\begin{equation}\label{e:4l49}
    \chi_1(\tilde{\mc{E}}_M-\tilde{\mc{E}}_0)=\chi_1\partial_\nu \mc{R}_0^{\Omega_2}[-\Delta, \chi_2](\mc{H}^{\Omega_1}_M-\mc{H}^{\Omega_1}_0)+\chi_1\partial_\nu (\mc{R}_M^{\Omega_2}-\mc{R}_0^{\Omega_2})[-\Delta, \chi_2]\mc{H}^{\Omega_1}_M. 
\end{equation}
We estimate those two terms separately to prove \eqref{e:4l53}.

1. For the first term in \eqref{e:4l49}, we claim
\begin{equation}\label{e:4l46}
    \|\chi_1\partial_\nu \mc{R}_0^{\Omega_2}[-\Delta, \chi_2](\mc{H}^{\Omega_1}_M-\mc{H}^{\Omega_1}_0) f\|_{H^{N,l}_\partial}\le C
    M^\e\|f\|_{H^{0,l'}_\partial}.
\end{equation}
Considering $\mc{H}^{\Omega_1}_M-\mc{H}^{\Omega_1}_0=-M^2\mc{R}_{0}^{\Omega_1}\mc{H}_M^{\Omega_1}$, we can estimate
\begin{multline}
     \|\chi_1\partial_\nu \mc{R}_0^{\Omega_2}[-\Delta, \chi_2]\mc{R}_{0}^{\Omega_1}\mc{H}_M^{\Omega_1} f\|_{H^{N,l}_\partial} \le C\|[-\Delta, \chi_2]\mc{R}_{0}^{\Omega_1}\mc{H}_M^{\Omega_1}f\|_{Y_{m_2}^{l}}\le C\|\varphi\mc{R}_{0}^{\Omega_1}\mc{H}_M^{\Omega_1}f\|_{X_{m_1,-1}^{l,0}}\\\le C\|\mc{H}_M^{\Omega_1}f\|_{Y_{m_1,-1}^{l,0}},
\end{multline}
using $\chi_1\partial_\nu \mathcal{R}_M^{\Omega_2}\varphi: \dot Y^{l}_{m_2}(\Omega_2)\rightarrow H^{N-\frac12, l}_\partial$ as in Step 4 in Lemma~\ref{l:sameNearInfinity}. Now Corollary~\ref{c:HimproveM} implies
\begin{equation}
    M^{2-\e}\|\mc{H}_M^{\Omega_1}f\|_{Y_{m_1,-1}^{l,0}}\le C\|\langle M\dd\rangle^{2-\e}\mc{H}_M^{\Omega_1}f\|_{Y_{m_1'-2,-1}^{l',0}}\le C\|f\|_{H^{0, l'}_\partial},
\end{equation}
giving \eqref{e:4l46} as desired. 

2. We address the second term in \eqref{e:4l49}. We claim 
\begin{equation}\label{e:4l45}
    \|\chi_1\partial_\nu(\mc{R}_M^{\Omega_2}-\mc{R}_0^{\Omega_2})[-\Delta, \chi_2]\mc{H}^{\Omega_1}_M f\|_{H^{N,l}_\partial}\le CM^\e\|f\|_{H^{0,l'}_\partial}. 
\end{equation}
Denote by $v=[-\Delta, \chi_2]\mc{H}^{\Omega_1}_M f$. Consider
\begin{equation}\label{e:4l52}
    \chi_1\partial_\nu(\mc{R}_M^{\Omega_2}-\mc{R}_0^{\Omega_2})\varphi=-M^2\chi_1\partial_\nu\mc{R}_0^{\Omega_2}\mc{R}_M^{\Omega_2}\varphi=-M^2\chi_1\partial_\nu\mc{R}_0^{\Omega_2}(1-\tilde\chi_1)\mc{R}_M^{\Omega_2}\varphi-M^2\chi_1\partial_\nu\mc{R}_0^{\Omega_2}\tilde\chi_1\mc{R}_M^{\Omega_2}\varphi,
\end{equation}
Corollary~\ref{c:RimproveM} implies $\chi_1\partial_\nu\mc{R}^{\Omega_2}_0(1-\tilde\chi_1): Y^{l,0}_{m_2, -N}(\Omega_2)\rightarrow H^{N-\frac12, l}(\partial\Omega_1)$ since $\chi_1$ and $1-\tilde\chi_1$ have disjoint support. Then the first term of \eqref{e:4l52} can be estimated by
\begin{multline}
    M^{2-\e}\|\chi_1\partial_\nu\mc{R}^{\Omega_2}_0(1-\tilde\chi_1)\mc{R}_M^{\Omega_2}\varphi v\|_{H^{N-\frac12, l}_\partial}\le CM^{2-\e}\|\mc{R}_M^{\Omega_2}\varphi v\|_{Y^{l}_{m_2}}\le CM^{2-\e}\|\mc{R}_M^{\Omega_2}\varphi v\|_{X^{l}_{m_2+2}}\\
    \le C\|\langle M\dd\rangle^{2-\e}\mc{R}_M^{\Omega_2}\varphi v\|_{X^{l'}_{m_2'}}\le C\|\langle M\dd\rangle^{2-\e}\varphi v\|_{Y^{l'}_{m_2'}}\le C\|\langle M\dd\rangle^{2-\e}\varphi v\|_{Y^{l'}_{m_1'}}.
\end{multline}
Noting $\chi_1\partial_\nu\mc{R}^{\Omega_2}_0: Y^{l,0}_{m_2,N}(\Omega_2)\rightarrow H^{N-\frac12, l}(\partial\Omega_1)$, we can estimate the other term by
\begin{multline}
    M^{2-\e}\|\chi_1\partial_\nu\mc{R}_0^{\Omega_2}\tilde\chi_1\mc{R}_M^{\Omega_2}\varphi v\|_{H^{N-\frac12, l}_\partial}\le CM^{2-\e}\|\tilde\chi_1\mc{R}_M^{\Omega_2}\varphi v\|_{Y^{l, 0}_{m_2, N}}\le CM^{2-\e}\|\tilde\chi_1\mc{R}_M^{\Omega_2}\varphi v\|_{X^{l, 0}_{m_2+2, N}}\\
    \le C\|\langle M\dd\rangle^{2-\e}\tilde\chi_1\mc{R}_M^{\Omega_2}\varphi v\|_{X^{l', 0}_{m_2', N}}\le C\|\langle M\dd\rangle^{2-\e}\varphi v\|_{Y^{l'}_{m_2'}}\le C\|\langle M\dd\rangle^{2-\e}\varphi v\|_{Y^{l'}_{m_1'}},
\end{multline}
since $\tilde\chi_1$ and $\varphi$ are disjointly supported, again using Corollary~\ref{c:RimproveM}. Let $\kappa\in C_c^\infty(\overline{\Omega_1})$ be that $\supp\kappa\cap\supp\varphi=\emptyset$, $\supp(1-\kappa)\cap\partial\osing=\emptyset$. Now Step 1 of Corollary~\ref{c:HimproveM} and Lemma~\ref{l:ellregunbdd} imply
\begin{multline}\label{e:4l54}
    \|\langle M\dd\rangle^{2-\e}\varphi v\|_{Y^{l'}_{m_1'}}\le C\|\langle M\dd\rangle^{2-\e}\varphi\mc{H}^{\Omega_1}_M f\|_{X^{l', 0}_{m_1', -1}}\le C\|\langle M\dd\rangle^{2-\e}\varphi\mc{H}^{\Omega_1}_M(1-\kappa) f\|_{X^{l', 0}_{m_1', -1}}\\
    +C\|\langle M\dd\rangle^{2-\e}\varphi\mc{H}^{\Omega_1}_M \kappa f\|_{X^{l', 0}_{m_1', -1}} \le C\|(1-\kappa)f\|_{H^{0, l'}_\partial}+C\|\kappa f\|_{H^{0, l'}_\partial}\le C\|f\|_{H^{0, l'}_\partial},
\end{multline}
giving \eqref{e:4l45} as desired. This concludes the proof of \eqref{e:4l53}. 

3. For the other three claims, it suffices to pay attention to \eqref{e:4l46} and \eqref{e:4l54}. Use \eqref{e:thing9} and \eqref{e:thing8} to observe the claims analogously.
\end{proof}

\begin{lemma}\label{l:errordependsonMcpt}
Let $\chi_1,\chi_2\in C^\infty_c(\partial\Omega)$ that $\supp\chi_1\cap\supp(1-\chi_2)=\emptyset$. Then we have
\begin{equation}
    \|\chi_2(\Lambda_M^\Omega-\Lambda_0^\Omega)\chi_1 f\|_{L^2(\partial\Omega)}\le CM^{\min(2\rardyr_{\Omega}-0,2)}\|\chi_1f\|_{L^2(\partial\Omega)},
\end{equation}
uniformly for all $f\in \mc{S}'(\partial\Omega)$ and $M\in [0,M_0]$. 
\end{lemma}

\begin{proof}
    We confuse $\chi_1,\chi_2$ with their extension in $C^\infty_c(\overline{\Omega})$ that $\supp\chi_1\cap\supp(1-\chi_2)=\emptyset$. Let $\tilde\chi_2\in C^\infty_c(\overline{\Omega})$ that $\supp\chi_2\cap\supp(1-\tilde\chi_2)=\emptyset$ and $\partial_\nu\tilde\chi_2=0$. Let $\tilde\Omega$ be a bounded uniformly Lipschitz domain such that $\tilde\Omega,\Omega$ share the common boundary near $\supp \tilde\chi_2$. Lemma \ref{l:DomChange} implies
    \begin{gather}\label{e:4l47}
        \chi_2(\Lambda_M^\Omega-\Lambda_0^\Omega)\chi_1=\chi_2(\Lambda_M^{\tilde\Omega}-\Lambda_0^{\tilde\Omega})\chi_1+\chi_2(\tilde{\mc E}_M-\tilde{\mc E}_0)\chi_1,\\
        \tilde{\mc E}_M-\tilde{\mc E}_0= \partial_\nu \mc R^{\tilde\Omega}_M[-\Delta,\tilde\chi_2]\mc H^\Omega_M- \partial_\nu\mc R^{\tilde\Omega}_0[-\Delta,\tilde\chi_2]\mc H^\Omega_0.
    \end{gather}
    Note $\Lambda_M^{\tilde\Omega}-\Lambda_0^{\tilde\Omega}=-M^2 \partial_\nu \mc R_0^{\tilde\Omega}\mc H^{\tilde\Omega}_M$ and we can estimate the first term of \eqref{e:4l47} by
     \begin{equation}\label{e:DiffLambdaBdd}
        \|\chi_2(\Lambda_M^{\tilde\Omega}-\Lambda_0^{\tilde\Omega}) \chi_1 f\|_{L^2_\partial}\le CM^2 \|\mc H^{\tilde\Omega}_M \chi_1 f\|_{H^\frac12}\le  CM^2 \|\chi_1 f\|_{L^2_\partial}
    \end{equation}
    since $\partial_\nu \mc R_0^{\tilde\Omega}:L^2\rightarrow L^2_\partial$ and $\mc H^{\tilde\Omega}_M: L^2_\partial\rightarrow H^\frac12$ are bounded per Lemma~\ref{l:reg2}. On the other hand, 
    \begin{equation}\label{e:4l48}
        \chi_2(\tilde{\mc E}_M-\tilde{\mc E}_0)\chi_1=\chi_2\partial_\nu  (\mc{R}^{\tilde\Omega}_M-\mc R^{\tilde\Omega}_0)[-\Delta,\tilde\chi_2]\mc H^\Omega_M \chi_1
    +\chi_2 \partial_\nu\mc R^{\tilde \Omega}_0[-\Delta,\tilde\chi_2](\mc H^\Omega_M-\mc H^\Omega_0)\chi_1,
    \end{equation}
    whose first term can be estimated by $\mc{R}^{\tilde\Omega}_M-\mc R^{\tilde\Omega}_0=-M^2\mc{R}^{\tilde\Omega}_0\mc R^{\tilde\Omega}_M$ and
    \begin{multline}
        \|\chi_2\partial_\nu  (\mc{R}^{\tilde\Omega}_M-\mc R^{\tilde\Omega}_0)[-\Delta,\tilde\chi_2]\mc H^\Omega_M \chi_1 f\|_{L^2_\partial}\le CM^2\|\mc R^{\tilde\Omega}_M[-\Delta,\tilde\chi_2]\mc H^\Omega_M \chi_1 f\|_{L^2}\\
        \le CM^2\|[-\Delta,\tilde\chi_2]\mc H^\Omega_M \chi_1 f\|_{H^{-\frac32}}
        \le CM^2\|\chi_1 f\|_{L^2_\partial}
    \end{multline}
    from Lemma~\ref{l:ellregunbdd}. Let $m=\min(\rardyr_{\Omega}-0, 1)$. The second term of \eqref{e:4l48} can be estimated by
    \begin{multline}
        \|\chi_2 \partial_\nu\mc R^{\tilde \Omega}_0[-\Delta,\tilde\chi_2](\mc H^\Omega_M-\mc H^\Omega_0)\chi_1 f\|_{L^2_\partial}\le CM^2\|[-\Delta,\tilde\chi_2]\mc R_0^\Omega\mc H^\Omega_M\chi_1 f\|_{L^2}
        \le CM^{2}\|\mc H^\Omega_M\chi_1 f\|_{Y^{0,0}_{-m,-1}}\\
        \le CM^{2m}\|\langle M\dd\rangle^{2-2m}\mc H^\Omega_M\chi_1 f\|_{Y^{0,0}_{m-2,-1}}\le CM^{2m}\|\chi_1 f\|_{L^2_\partial}
    \end{multline}
    via Corollary~\ref{c:HimproveM}. This concludes the proof. 
\end{proof}

We now record our main description of the \dtn map on \acl domains in $\mathbb{R}^d$: 
\begin{proposition}
\label{p:dtN}
Let $\chi, \tilde\chi \in S^0(\partial\Omega)$ with $\supp\chi\cap\supp(1-\tilde\chi)=\emptyset$, $\supp\tilde\chi\subset\partial\oreg$. Then we have
\begin{equation}
    \chi \Lambda^{\Omega}_M= \chi P\tilde\chi+\chi\mathcal{E}_M, \ \Lambda^{\Omega}_M\chi= \tilde\chi P\chi+\mathcal{E}_M'\chi, 
\end{equation}
where $\chi P\tilde\chi, \tilde\chi P\chi\in \Psi^{0,-\delta}_{1,\delta}(\partial\Omega)$ are the left/right microlocal description of $\Lambda^{\Omega}_M$ near $\supp\chi$ with
\begin{equation}
    \sigma(P)= h^{-1}\langle x\rangle^{-\delta}\big\langle h\langle  x\rangle^\delta \sqrt{|\xi|^2_g+M^2}\big\rangle
\end{equation}
given in Lemma~\ref{l:factor}. We confuse $\chi,\tilde\chi$ with extensions specified in the proof, and the errors $\mathcal{E}_M, \mathcal{E}_M'$ are given by
\begin{gather}\label{e:4l36}
    \chi\mathcal{E}_M=\chi\mc{E}_\flat \tilde\chi+\chi\tilde{\mc{E}}_M=\chi(\Lambda^\Omega_M-P)\tilde\chi+\chi\partial_\nu \mc{R}_M^{\Omega}[-\Delta, \tilde\chi]\mc{H}^{\Omega}_M,\\
    \label{e:4l37}
    \mathcal{E}_M'\chi=\tilde\chi\mc{E}_\flat \chi+\tilde{\mc{E}}_M^*\chi=\tilde\chi(\Lambda^\Omega_M-P)\chi-\partial_\nu \mc{R}_M^{\Omega}[-\Delta, \tilde\chi]\mc{H}^{\Omega}_M\chi.
\end{gather}
Moreover, let $\psi(\xi)\in S^0(\mathbb{R}^{d-1})$ with $0\notin \supp \psi$, $N\geq 0$, $|l|<\rardyr_\Omega+\frac12$. We can estimate the error $\mc{E}_M$ by
\begin{gather}\label{e:errorEst}
 \|\chi\psi(D) \chi \mathcal{E}_M\chi\psi(D)\chi f\|_{H^{N,N}(\partial\Omega)}\le C\|f\|_{H^{-N,-N}(\partial\Omega)},\\
  \|\chi\psi(D) \chi\mathcal{E}_M  f\|_{H^{N,N}(\partial\Omega)}\le C\|f\|_{H^{0,l}(\partial\Omega)}, \label{e:errorEst2}\\
    \|\chi\mathcal{E}_M\chi\psi(D)\chi f\|_{H^{N,l}(\partial\Omega)}\le C\|f\|_{H^{-N,-N}(\partial\Omega)}, \label{e:errorEst3}\\
     \|\chi\mathcal{E}_M f\|_{H^{N,l}(\partial\Omega)}\le C\|f\|_{H^{0,\max(l-1, -\rardyr_\Omega-\frac12)+0}(\partial\Omega)}, \label{e:errorEst4}
\end{gather}
uniformly for all $f\in \mc{S}'(\partial\Omega)$ and $M\in [0,M_0]$. Same estimates hold for $\mc{E}_M'$ instead of $\mc{E}_M$. All of $P,\mc{E}_M, \mc{E}_M'$ depend on $M$ continuously. 
\end{proposition}
\begin{proof}
1. Let $U$ be a regular collar of $\Omega$ as in Remark~\ref{r:aclcollar} that $\supp\tilde\chi\subset U$. We confuse $\tilde\chi, \chi$ with their extensions in $S^0(U)$ that $\partial_{\nu}\tilde\chi=\partial_{\nu}\chi=0$ and $\supp\chi\cap\supp(1-\tilde\chi)=\emptyset$. We address the case $\chi\mc{E}_M$ by proving estimates \eqref{e:errorEst} to \eqref{e:errorEst4} with $\chi\mc{E}_M$ replaced by $\chi\tilde{\mc{E}}_M$ and $\chi\mc{E}_\flat\tilde\chi$ respectively. The case $\mc{E}_M'\chi$ then follows similarly. 

2. We will use Lemma~\ref{l:sameNearInfinity} for $\chi\tilde{\mc{E}}_M$. The compatibility condition reduces to $l_2< l_1+1$. Estimates \eqref{e:errorEst}, \eqref{e:errorEst2} and \eqref{e:errorEst3} follow immediately. For \eqref{e:errorEst4}, we take
    \begin{equation}
        l_2=l, \ l_1=\max(l-1,-\rardyr_\Omega-\frac12)+0,
    \end{equation}
    which satisfies the compatibility condition and we have the claim. 
    
    3. Proposition~\ref{p:dtNParametrix} implies for any $u\in\mc{S}'(\partial\Omega)$, 
    \begin{equation}
    \|\tilde\chi\Gamma^{N,l}_{N,\delta}\chi\mc{E}_\flat\tilde\chi u\|_{L^2(\partial\Omega)}\le C\|\Gamma^{-N,\max(l+\frac12,1-\rardyr_\Omega)-\frac32\delta+0}_{-N,\delta} u\|_{L^2(\partial\Omega)},
    \end{equation}
    where choosing $\delta<1$ close to $1$ makes
    \begin{equation}
        \max(l+\frac12,1-\rardyr_\Omega)-\frac32\delta+0\le\max(l-1, -\rardyr_\Omega-\frac12)+0.
    \end{equation}
    The estimate above then implies \eqref{e:errorEst} to \eqref{e:errorEst4} with $\chi\mc{E}_M$ replaced with $\chi\mc{E}_\flat\tilde\chi$. This concludes the proof.
\end{proof}
\begin{remark}
    Here the microlocal error $\mc{E}_\flat$ measures the nonlocality at the zero frequency of $\Lambda^{\Omega}_M$ due to its symbolic singularity in $\sqrt{|\xi|^2+M^2}$, and the nonlocal error $\tilde{\mc{E}}_M$ measures the nonlocality at the zero frequency of $\Lambda^\Omega_M$ due to the Lipschitz singularity in $x\in \partial\Omega$.  
\end{remark}

\begin{remark}
    All constants in this section can depend on $\Omega$ in the following way: the upper bound on uniform Lipschitz characters of $\partial\Omega$, the upper bound on the choice of $r_0\ge 0$, and the lower bound on the radial Hardy constant $\rardyr_{\Omega}$, the upper bound on the diameter of $\partial\Omega\setminus\partial\tilde\Omega$ in Remark~\ref{r:aclcollar}, for each multi-index $\alpha$, the upper bounds on $C_\alpha$ of derivatives of the metric $\gplus$ in \eqref{e:3l26}, and lower bounds on $|\det \gplus|$.
\end{remark}

\section{Microlocal analysis of the \dtn map on \acl domains}
\label{s:microlocalDtN}
The purpose of this section is to prove a limiting absorption principle for the DtN map on \acl domains. In fact, following the now standard ideas of Vasy and Melrose~\cite{mel94,vasy2000semiclassical,Vasy2021Resolvent} we show that the limiting absorption resolvent as an actual inverse on certain variable order Sobolev spaces. 

\noindent{\bf{Outline of the strategy}}

This section relies heavily on the microlocal description of the DtN map from Section~\ref{s:3}. With that model available the proof of the limiting absorption principle follows steps similar to those used when proving the limiting absorption principle for the Laplacian on an asymptotically conic manifold. 
\begin{enumerate}
    \item Prove elliptic estimates near $|x|=\infty$ and away from $|\xi|_g^2+M^2=\omega^2.$

    In this step we use the standard elliptic parametrix construction away from $\xi=0$. However, $\xi=0$ requires special treatment since $\Lambda_M^\Omega$ is \emph{not} pseudodifferential there. Instead, we show in Lemma~\ref{l:EllEst_Near0} that the difference between $\Lambda_M^\Omega$ and a pseudodifferential operator is sufficiently well behaved to be absorbed into the error. 

    \item Elliptic estimates on $|x|<\infty$, $|\xi|=\infty$. 
    
    Once again, away from the Lipschitz part of $\partial\Omega$, we are able to use standard elliptic parametrices. Near the Lipschitz part, we rely on the elliptic theory for the Laplacian as in Lemma~\ref{l:ellregunbdd}. 

    \item Propagation estimates on $\{|\xi|_g^2+M^2=\omega^2\}\cap \{|x|=\infty\}$. 

    With the microlocal description of the DtN map from Section~\ref{s:3}, we can apply the standard real principal type propagation estimates and radial point estimates (see e.g.~\cite[Appendix E]{DyZw:19}). It is in this step where, as usual, we are required to use variable order Sobolev spaces.

    \item Absence of non-$L^2$ outgoing states.

    Finally, we need to show that if $u$ is purely outgoing an an appropriate microlocal sense and satisfies $(\Lambda_M^\Omega-\omega)u=0$ for some $\omega>M$, then $u$ is in $L^2$. Together with the propagation estimates above, this will allow us to show that if $\omega$ is not an eigenvalue for the DtN map, then i) $(\Lambda_M^\Omega-\omega)$ is invertible on some appropriately defined spaces and ii) this inverse is also the limiting absorption resolvent.  
\end{enumerate}

\subsection{Geometric and microlocal preliminaries}

In this section, we work on an \acl domain $\Omega$, with its boundary $\partial\Omega$ being a scattering manifold. 

We compactify $\partial\Omega$ and $T^*\partial\Omega$ in the picture of radial compactification in both $x$ and $\xi$. We abuse the notation and may freely use the projective coordinates $(\rho, y, \tau,\eta)$ on $T^*\partial\Omega$, with trivialization $\tau\rho^{-2}d\rho+\eta\rho^{-1}dy$, or the local normal coordinates $(x,\xi)$ on $T^*\partial\Omega$ induced from Fermi normal coordinates $(x_0,x)$ locally near $\partial\Omega$ equipped with induced metric $g$, by restricting to $T^*\partial\oreg$ and $\partial\oreg$ respectively if necessary. One has the canonical transformation 
\begin{equation}
    \xi\cdot dx=\tau\rho^{-2}d\rho+\rho^{-1}\eta\cdot dy. 
\end{equation}

Note $\partial T^*\partial\oreg=\partial T^*\partial\Omega$. We say $a\in S^{s,l}_{k,\delta}(T^*\partial\Omega)$, if $a\in C^\infty(T^*\partial\oreg)$ and for each $K\subset\partial\oreg$ closed, and for each multi-indices $\alpha,\beta$, there is $C_{\alpha\beta K}>0$ that
\begin{equation}
\abs{(\langle x\rangle\partial_{x})^\alpha (\langle x\rangle^{-\delta} \langle \langle x\rangle^{\delta}\sqrt{|\xi|_g^2+M^2}\rangle\partial_\xi)^\beta a(x,\xi)}\le C_{\alpha\beta K}\langle \xi\rangle^{s}\langle x\rangle^{l}\langle \langle x\rangle^{\delta}\sqrt{|\xi|_g^2+M^2}\rangle^{k},
\end{equation}
uniformly for all $(x,\xi)\in T^*K$. 

Consider a finite family of charts covering $(\Phi_{j}, U_j)$ over $\partial\Omega$, associated with a partition of unity $\phi_j\in S^0(\partial\Omega)$ and $\supp\phi_j\subset U_j$ that $\sum_{j}\phi_j=1$. For those $U_j\cap\partial\osing=\emptyset$, let $\varphi_{j0}=\phi_j$ and $\varphi_{jk}=0$ for all $k\ge 1$.
For those $U_j\cap\partial\osing\neq\emptyset$ 
consider a dyadic partition of unity $\varphi_{jk}\in S^0(\partial\Omega)$ that is locally finite on $\partial\oreg$, $\supp\varphi_{j0}\subset\{x:d(x,\partial\osing)>\frac12\}$,
\begin{equation}
    \supp\varphi_{jk}\subset \{x:2^{-k-1}<d(x,\partial\osing)<2^{-k+1}\} \text{ for }k\ge 1, \ \sum_{k=0}^\infty \varphi_{jk}=\phi_{j}\mathbbm{1}_{U_j\cap\partial\oreg}. 
\end{equation}
Let $\tilde\varphi_{jk}\in S^0(\partial\Omega)$ with  $\supp(1-\tilde\varphi_{jk})\cap \supp \varphi_{jk}=\emptyset$ and $d(\supp\tilde\varphi_{j0}, \partial\osing)>\frac12$ and $\supp\tilde\varphi_{jk}\subset \{x:2^{-k-1}<d(x,\partial\osing)<2^{-k+1}\}$ for $k\ge 1$. Fix a localizer in distance $\psi\in C^\infty_c(\mathbb{R})$ with $\supp\psi\subset [-\frac14,\frac14]$ and $\psi\equiv 1$ on $[-\frac18,\frac18]$. We define the individual quantization map on each dyadic piece by
\begin{equation}
    \Op_{jk}(a)=\frac1{(2\pi)^{d-1}}\int_{\mathbb{R}^{d-1}_\xi} e^{i\langle x-y, \xi\rangle}a(x,\xi)\psi(2^k|x-y|)\,d\xi
\end{equation}
and thus the global quantization map
\begin{equation}
    \Op(a)=\sum_{j}\sum_{k=0}^\infty\tilde\varphi_{jk}\Phi_j^{*}\Op_{jk}(\Phi_j^{-*}(\varphi_{jk}a))\Phi_j^{-*}\tilde\varphi_{jk}.
\end{equation}

We say $A\in\Psi^{s,l}_{k,\delta}(\partial\Omega)$, if 
\begin{enumerate}
    \item $A: C^\infty_c(\partial\Omega)\rightarrow L^2(\partial\Omega)$.
    \item There is $\chi\in S^0(\partial\Omega)$ with $\supp\chi\subset\partial\oreg$ and $a\in S^{s,l}_{k,\delta}(T^*\partial\Omega)$ that 
    \begin{equation}
        A=\chi A\chi=\Op(a)+\chi O(\Psi^{-\infty,-\infty}(\partial\Omega))\chi. 
    \end{equation}
\end{enumerate}
We say $A\in\Psi^{s,l}(\partial\Omega)$ if the above condition is met with $a\in S^{s,l}(T^*\partial\Omega)$. 
We define the wavefront set and elliptic set by
\begin{gather}
    \WF A=\bigcup \{(x,\xi)\in\WF(\chi A\chi): \chi\in S^0(\partial\Omega), \supp\chi\subset\partial\oreg\},\\
    \Ell A=\bigcup \{(x,\xi)\in\Ell(\chi A\chi): \chi\in S^0(\partial\Omega), \supp\chi\subset\partial\oreg\}.
\end{gather}

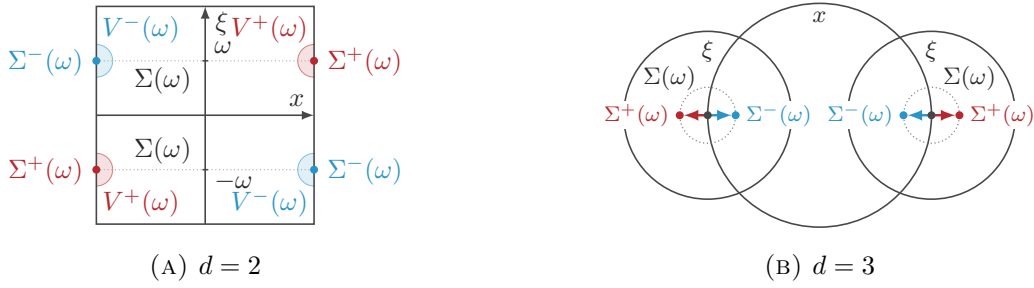
\begin{figure}[t]
  \centering
  \begin{subfigure}[t]{0.5\textwidth}
        \centering
        \begin{tikzpicture}[x=1cm,y=1cm]
  \definecolor{colorAKANE}{RGB}{183,40,46}
  \definecolor{colorUSUAI}{RGB}{40,148,200}
  \colorlet{plotPrimary}{colorAKANE}
  \colorlet{plotTertiary}{colorUSUAI}
  \colorlet{plotAxis}{black!72}
  \colorlet{plotGuide}{black!34}

  \tikzset{
    >=latex, font=\small, mark size=1.5pt, mark options={thin},
    line cap=round, line join=round,
    boundary/.style={draw=plotAxis,line width=0.6pt},
    axis/.style={draw=plotAxis,line width=0.6pt,->},
    characteristic guide/.style={draw=plotGuide,line width=0.4pt,
      densely dotted},
    plus neighborhood/.style={draw=plotPrimary!60,fill=plotPrimary!14,
      line width=0.35pt},
    minus neighborhood/.style={draw=plotTertiary!60,fill=plotTertiary!14,
      line width=0.35pt},
    plus point/.style={circle,draw=none,fill=plotPrimary,
      inner sep=0pt,minimum size=3pt},
    minus point/.style={circle,draw=none,fill=plotTertiary,
      inner sep=0pt,minimum size=3pt},
    figure label/.style={text=black!84,inner sep=0.8pt},
    colored label/.style={figure label,anchor=base,yshift=-0.65ex}
  }

  \path[use as bounding box] (-3.75,-1.6) rectangle (3.75,1.6);
  \begin{scope}[scale=0.72]
    \coordinate (A) at (-2,0);
    \coordinate (B) at (2,0);
    \coordinate (C) at (0,-2);
    \coordinate (D) at (0,2);
    \coordinate (upperLeft) at (-2,1);
    \coordinate (upperRight) at (2,1);
    \coordinate (lowerLeft) at (-2,-1);
    \coordinate (lowerRight) at (2,-1);

    \begin{scope}[radius=0.3cm]
      \clip (A |- C) rectangle (B |- D);
      \path[minus neighborhood] (upperLeft) circle[];
      \path[plus neighborhood] (upperRight) circle[];
      \path[plus neighborhood] (lowerLeft) circle[];
      \path[minus neighborhood] (lowerRight) circle[];
    \end{scope}

    \begin{scope}
      \draw[boundary,line join=miter] (A |- C) rectangle (B |- D);
      \draw[characteristic guide] (upperLeft) -- (upperRight);
      \draw[characteristic guide] (lowerLeft) -- (lowerRight);
      \draw[axis] (A) -- (B);
      \draw[axis] (C) -- (D);
      \draw[boundary] (-0.05,1) -- (0.05,1);
      \draw[boundary] (-0.05,-1) -- (0.05,-1);
    \end{scope}

    \begin{scope}
      \node[minus point] at (upperLeft) {};
      \node[plus point] at (upperRight) {};
      \node[plus point] at (lowerLeft) {};
      \node[minus point] at (lowerRight) {};
      \node[colored label,text=plotTertiary,anchor=base east,xshift=-4pt]
        at (upperLeft) {$\Sigma^-(\omega)$};
      \node[colored label,text=plotPrimary,anchor=base west,xshift=4pt]
        at (upperRight) {$\Sigma^+(\omega)$};
      \node[colored label,text=plotPrimary,anchor=base east,xshift=-4pt]
        at (lowerLeft) {$\Sigma^+(\omega)$};
      \node[colored label,text=plotTertiary,anchor=base west,xshift=4pt]
        at (lowerRight) {$\Sigma^-(\omega)$};
    \end{scope}

    \begin{scope}
      \node[colored label,text=plotTertiary] at (-1.18,1.60)
        {$V^-(\omega)$};
      \node[colored label,text=plotPrimary] at (1.18,1.60)
        {$V^+(\omega)$};
      \node[colored label,text=plotPrimary] at (-1.18,-1.60)
        {$V^+(\omega)$};
      \node[colored label,text=plotTertiary] at (1.18,-1.60)
        {$V^-(\omega)$};
      \node[figure label] at (-0.78,0.65) {$\Sigma(\omega)$};
      \node[figure label] at (-0.78,-0.65) {$\Sigma(\omega)$};
      \node[figure label,anchor=south west] at (0.12,1.10) {$\omega$};
      \node[figure label,anchor=north west] at (0.12,-1.01) {$-\omega$};
      \node[figure label,anchor=south east] at (1.87,0.12) {$x$};
      \node[figure label,anchor=north west] at (0.14,1.96) {$\xi$};
    \end{scope}
  \end{scope}
\end{tikzpicture}%
        \caption{$d=2$}
  \end{subfigure}\hfill
  \begin{subfigure}[t]{0.5\textwidth}
        \centering
        \begin{tikzpicture}[x=1cm,y=1cm]
  \definecolor{colorAKANE}{RGB}{183,40,46}
  \definecolor{colorUSUAI}{RGB}{40,148,200}
  \colorlet{plotPrimary}{colorAKANE}
  \colorlet{plotTertiary}{colorUSUAI}
  \colorlet{plotAxis}{black!72}

  \tikzset{
    >=latex, font=\small\fontsize{9}{11}\selectfont,
    mark size=1.5pt, mark options={thin},
    line cap=round, line join=round,
    boundary/.style={draw=plotAxis,line width=0.6pt},
    characteristic/.style={draw=black!52,line width=0.45pt,
      densely dotted},
    outward vector/.style={draw=plotPrimary,line width=0.9pt,->,
      shorten >=1.5pt},
    inward vector/.style={draw=plotTertiary,line width=0.9pt,->,
      shorten >=1.5pt},
    center point/.style={circle,draw=none,fill=plotAxis,
      inner sep=0pt,minimum size=3pt},
    plus point/.style={center point,fill=plotPrimary},
    minus point/.style={center point,fill=plotTertiary},
    figure label/.style={text=black!84,inner sep=0.8pt},
    radial label/.style={figure label,
      font=\footnotesize\fontsize{8}{10}\selectfont}
  }

  \path[use as bounding box] (-3.75,-1.6) rectangle (3.75,1.6);
  \begin{scope}[scale=1.48]
    \coordinate (A) at (-1,0);
    \coordinate (B) at (1,0);
    \coordinate (leftPlus) at (-1.25,0);
    \coordinate (leftMinus) at (-0.75,0);
    \coordinate (rightMinus) at (0.75,0);
    \coordinate (rightPlus) at (1.25,0);

    \begin{scope}
      \draw[boundary] (0,0) circle[radius=1cm];
      \draw[boundary] (A) circle[radius=0.75cm];
      \draw[boundary] (B) circle[radius=0.75cm];
      \draw[characteristic] (A) circle[radius=0.25cm];
      \draw[characteristic] (B) circle[radius=0.25cm];
    \end{scope}

    \begin{scope}
      \node[radial label,text=plotPrimary,fill=black!0,
        anchor=east,xshift=-3pt]
        at (leftPlus) {$\Sigma^+(\omega)$};
      \node[radial label,text=plotPrimary,fill=black!0,
        anchor=west,xshift=3pt]
        at (rightPlus) {$\Sigma^+(\omega)$};
      \node[radial label,text=plotTertiary,fill=black!0,
        anchor=west,xshift=3pt]
        at (leftMinus) {$\Sigma^-(\omega)$};
      \node[radial label,text=plotTertiary,fill=black!0,
        anchor=east,xshift=-3pt]
        at (rightMinus) {$\Sigma^-(\omega)$};
    \end{scope}

    \begin{scope}
      \draw[outward vector] (A) -- (leftPlus);
      \draw[outward vector] (B) -- (rightPlus);
      \draw[inward vector] (A) -- (leftMinus);
      \draw[inward vector] (B) -- (rightMinus);
      \node[center point] at (A) {};
      \node[center point] at (B) {};
      \node[plus point] at (leftPlus) {};
      \node[plus point] at (rightPlus) {};
      \node[minus point] at (leftMinus) {};
      \node[minus point] at (rightMinus) {};
    \end{scope}

    \begin{scope}
      \node[figure label] at (0,0.88) {$x$};
      \node[figure label] at (-1,0.57) {$\xi$};
      \node[figure label] at (1,0.57) {$\xi$};
      \node[figure label] at (-1.34,0.33) {$\Sigma(\omega)$};
      \node[figure label] at (1.34,0.33) {$\Sigma(\omega)$};
    \end{scope}
  \end{scope}
\end{tikzpicture}%
        \caption{$d= 3$}
  \end{subfigure}
  \caption{Sink/source structures.}
  \label{f:5f1}
\end{figure}

Let $g$ be the induced metric of $\partial\oreg\subset\mathbb{R}^d$ and denote the principal symbol of the approximation $P$ to $\Lambda^\Omega_M$ by
\begin{equation}
p(x,\xi;M)=\langle x\rangle^{-\delta}\langle \langle x\rangle^{\delta}\sqrt{|\xi|_g^2+M^2}\rangle\in S^{0,-\delta}_{1,\delta}(\partial\Omega),
\end{equation}
where the small parameter $h>0$ was fixed small in Remark~\ref{r:fixe}, and $\delta\in [0,1)$. In this section $\omega>0$. Denote the characteristics of the approximation $P$ to $\Lambda^\Omega_M-\om$ by the set of dimension $2d-4$, 
\begin{equation}
\Sigma(\omega)=\{(x,\xi)\in \partial T^*\partial\Omega: p(x,\xi;M)=\om\}.
\end{equation}
Near $\Sigma(\omega)$, we have $g(\rho,y)=g_\infty(y)+O(\rho S^0(\partial W_Y))$ and 
\begin{gather}
    p(\rho, y,\tau,\eta)=\sqrt{\tau^2+|\eta|^2_{g_\infty}+M^2}+ O(\rho S^{1,0})+ O(\rho^{2\delta}S^{-1,0})
\end{gather}
Without loss of generality assuming $\delta\ge\frac12$ and noting the frequency is bounded near $\Sigma(\omega)$, we can ignore the two error terms above as subprincipal. Let $\hat H_{p}=\rho^{-1}H_p$ be the rescaled Hamiltonian vector field associated to $p$ on $T^*\partial\Omega$, given by
\begin{equation}
    \hat H_p=(\tau^2+|\eta|^2_{g_\infty}+M^2)^{-\frac12}\Big( \tau\rho\partial_\rho+ \eta^t\cdot g_\infty \cdot \partial_y-|\eta|^2_{g_\infty} \partial_\tau
    +\sum_{j=1}^{d-2}(\tau\eta_j -\frac12(\eta^t\cdot(\partial_{y_j}g_\infty)\cdot\eta)\partial_{\eta_j}\Big),
\end{equation}
near $\Sigma(\omega)$. Let $e^{t \hat H_p}$ be the flow associated to $\hat H_{p}$ on $T^*\partial\Omega$. We let the set of dimension $d-2$, 
\begin{equation}
    \Sigma^\pm(\omega)=\{(\rho,y,\tau,\eta)\in \Sigma(\omega):\rho=0, \tau=\mp\omega, \eta=0\}
\end{equation}
as radial sink/source set of the approximation $P$ to $\Lambda^\Omega_M-\sqrt{\omega^2+M^2}$. See Figure~\ref{f:5f1}. 

The sink/source nature follows from the lemma below. 
\begin{lemma}
The radial set $\Sigma^-(\omega)\subset T^*\partial\Omega$ consists of sources, and $\Sigma^+(\omega)$ consists of sinks. 
\end{lemma}
\begin{proof}
    It suffices to show the claim on an \acs $\Omega_e$ that $\Omega_e\cap\{\rho\le\epsilon\}=\Omega\cap\{\rho\le\epsilon\}$ for some $\e>0$ small. Note $p|_{\rho=0}=\sqrt{|\xi|^2_g+M^2}$ and $2p H_{p}|_{\rho=0}=H_{-\Delta}|_{\rho=0}$. This implies the sink and source of $P$ are the same of $-\Delta$ on $\partial\Omega_e$ and the claim follows from \cite[Lemma 2]{MeZw:96}.
\end{proof}

\subsection{Microlocal estimates through zeros, sinks and sources}
We start by showing an elliptic estimate for the operator $\Lambda^\Omega_M-\om-i\mu$, on $\{|\xi|_g\neq\omega\}$. Let $\intw\subseteq (0,\infty)$ be a compact interval. We define the decay threshold at the zero frequency by
\begin{equation}
    L_0=\rardyr_{\Omega}+\frac12.
\end{equation}
Fix $\kappa\in C^\infty_c(\overline{\Omega})$ that $\supp(1-\kappa)\cap \partial\osing=\emptyset$. Fix $\kappa_0, \tilde\kappa$ that satisfy the same requirement for $\kappa$ with $\kappa_0\prec\kappa\prec\tilde\kappa$.

\begin{lemma}[Elliptic estimate at zero frequency]\label{l:EllEst_Near0}
Let $\Omega\subset \mathbb{R}^d$ be an \acl domain. Let $|l_0|<L_0$. Let $A_0, B_0\in \Psi^{0,0}(\partial\Omega)$, $\WF A_0\cap \Sigma(\omega)=\emptyset$, $\WF A_0\subset \Ell B_0$, and
\begin{equation}
    \{(x,\xi): x\in \supp(1-\kappa), \xi=0\}\subset\Ell B_0, \ \WF(B_0)\subset\{(x,\xi):|\xi|\le \e_0\}
\end{equation}
for some $\e_0>0$. Then for any $N>0$ we have
\begin{equation}
    \|A_0 u\|_{H^{N,l_0}}\le C\|B_0(\Lambda^\Omega_M-\om-i\mu)u\|_{H^{-N,l_0}}+C\|\tilde\kappa u\|_{H^{0}}+C\|u\|_{H^{-N, -N}}
\end{equation}
uniformly for all $u\in \mathcal{S}'$ with $\WF^{-\infty,-L_0+0}(u)\cap \{\xi=0\}=\emptyset$, $\omega\in\intw$ and $\mu\ge0$ small, $M\in[0, M_0]$. 
\end{lemma}
\begin{proof}
    1. We claim the lemma holds for some $A_0'$ in place of $A_0$, that is, 
    \begin{equation}\label{e:3l24}
    \|A_0' u\|_{H^{N,l_0}}\le C\|B_0(\Lambda^\Omega_M-\om-i\mu)u\|_{H^{-N,l_0}}+C\|\tilde\kappa u\|_{H^{0}}+C\|u\|_{H^{-N, -N}},
    \end{equation}
    where $A_0'\in\Psi^{0,0}(\partial\Omega)$ with $\WF A_0\subset\Ell A_0'$, $\Pi_x(\WF (\Id -A_0')\cap \{\xi=0\})\subset \supp\tilde\kappa$, $\WF A_0'\subset \Ell B_0$. Then the standard elliptic estimate implies the lemma. To show the claim, note since $B_0=\chi B_0\chi$ for some $\chi\in S^0$ whose support touches finitely many connected components of $\partial\Omega$, we can approximate $\Lambda_M^\Omega$ by $P\in \Psi^{0,-\delta}_{1,\delta}$ in Proposition~\ref{p:dtN}. Note $A_0',B_0\in\Psi^{-\infty,0}_{0,\delta}$ per the compact support in frequencies. We also have $B_0\chi(P-\om-i\mu)\tilde\chi \in \Psi^{-\infty,0}_{0,\delta}$ is elliptic on $\WF A_0'$ for some $\tilde\chi$ satisfying all conditions for $\chi$ with $\supp\chi\cap\supp(1-\tilde\chi)=\emptyset$. By Lemma \ref{l:MicroPar} there is $E\in \Psi^{-\infty,0}_{0,\delta}$ such that
    \begin{equation}
        \Gamma^{s,l_0}_{k,\delta}A_0'=E\Gamma^{s,l_0}_{k,\delta}B_0\chi (P-\om-i\mu)\tilde\chi+O(\Psi^{-\infty,-\infty}).
    \end{equation}
    Hence we have for any $s\in\mathbb{R}$, 
    \begin{multline}\label{e:EllEst_Tot}
        \|\Gamma^{s,l_0}_{k,\delta}A_0' u\|\leq C\|\Gamma^{s,l_0}_{k,\delta}B_0\chi(P-\om-i\mu)\tilde\chi u\|+C\|\Gamma^{-N,-N}_{-N,\delta}u\|\\
        \le C\|\Gamma^{s,l_0}_{k,\delta}B_0\chi(\Lambda^\Omega_M-\om-i\mu)u\|+C\|\Gamma^{s,l_0}_{k,\delta}B_0\chi\mc{E}_M u\|+C\|\Gamma^{-N,-N}_{-N,\delta}u\|,
    \end{multline}
    where $\mc{E}_M$ is from Proposition~\ref{p:dtN}. Taking $k=0$ gives on $\partial \Omega$, 
    \begin{equation}
        \|A_0' u\|_{H^{s,l_0}}\le C\|B_0(\Lambda^\Omega_M-\om-i\mu)u\|_{H^{s,l_0}}+C\|\chi\mc{E}_M u\|_{H^{s, l_0}}+C\|u\|_{H^{-N, -N}}. 
    \end{equation}
Let $\psi_0(\xi)\in C^\infty_c(\mathbb{R}^{d-1})$ be that $\supp(1-\psi_0)\cap\{\xi=0\}=\emptyset$ and $\supp ((1-\tilde\kappa(x))\psi_0(\xi))\subset \Ell A_0'$. With the decomposition $\Id=(1-\tilde\kappa)\psi_0(D)(1-\kappa)+(1-\tilde\kappa)(\Id-\psi_0(D))(1-\kappa)+\tilde\kappa$, we have from Proposition~\ref{p:dtN} and the standard elliptic estimate that
\begin{multline}\label{e:5l7}
    \|\chi\mc{E}_M u\|_{H^{s, l_0}}\le C\|(1-\tilde\kappa)\psi_0(D)(1-\kappa) u\|_{H^{-N, \max(l_0-1+0, -L_0+0)}}+C\|u\|_{H^{-N, -N}}+C\|\tilde\kappa u\|_{H^{0}}\\
    \le C\|A_0' u\|_{H^{-N, \max(l_0-1+0, -L_0+0)}}+C\|u\|_{H^{-N, -N}}+C\|\tilde\kappa u\|_{H^{0}}.
\end{multline}

2. It remains to reduce the $\|A_0' u\|_{H^{-N, \max(l_0-1+0, -L_0+0)}}$ term in \eqref{e:5l7}. When $l_0<-L_0+1$, one can interpolate
\begin{equation}
    \|A_0' u\|_{H^{-N, \max(l_0-1+0, -L_0+0)}}\le C\e^{-1}\|u\|_{H^{-N, -N}}+\e\|A_0' u\|_{H^{s, -L_0+0}}
\end{equation}
for small $\e>0$ to obtain \eqref{e:3l24} via absorption. For $-L_0+1\le l_0< -L_0+2$ with $l_0<L_0$, interpolation implies
\begin{multline}
    \|A_0' u\|_{H^{-N, \max(l_0-1+0, -L_0+0)}}\le C\e^{-1}\|u\|_{H^{-N, -N}}+\e\|A_0' u\|_{H^{s, l_0-1+0}}\\
    \le C\|B_0(\Lambda^\Omega_M-\om-i\mu)u\|_{H^{-N,l_0-1+0}}+C\|\tilde\kappa u\|_{H^{0}}+C\|u\|_{H^{-N, -N}}
\end{multline}
where we used \eqref{e:3l24} for $l_0-1<-L_0+1$. Iterate this procedure to obtain \eqref{e:3l24} for all $|l_0|<L_0$.  
\end{proof}

\begin{lemma}[Conormal regularity at zero frequency]\label{l:EllEst_Near0_breg}
Let $|l_0|< L_0$. Let $A_0, B_0\in \Psi^{0,0}(\partial\Omega)$, $\WF A_0\cap \Sigma(\omega)=\emptyset$, $\WF A_0\subset \Ell B_0$, and
\begin{equation}
    \{(x,\xi): x\in \supp(1-\kappa), \xi=0\}\subset\Ell B_0, \ \WF(B_0)\subset\{(x,\xi):|\xi|\le \e_0\}
\end{equation}
for some $\e_0>0$. Let $k\in \mathbb{N}_0$. Let $\mc{V}_1,\dots, \mc{V}_k\in \Diffhom_b^{1,1}(\partial\Omega)$ with support disjoint from $\supp\kappa$. Then for any $N>0$ we have
\begin{equation}
    \|A_0 \mc{V}_1\cdots\mc{V}_k u\|_{H^{N,l_0}}\le C\|B_0(\Lambda^\Omega_M-\om-i\mu)u\|_{H^{-N,l_0+k+0}}+C\|\tilde\kappa u\|_{H^{0}}+C\|u\|_{H^{-N, -N}}
\end{equation}
uniformly for all $u\in \mathcal{S}'$ with $\WF^{-\infty,-L_0+0}(u)\cap \{\xi=0\}=\emptyset$, $\omega\in\intw$ and $\mu\ge0$ small, $M\in[0, M_0]$. 
\end{lemma}
\begin{proof}
    0. Let $\tilde B_0,\tilde{\tilde{B}}_0$ satisfy the same requirement for $B_0$ with $\WF B_0\cap\WF(\Id-\tilde B_0)=\WF \tilde B_0\cap \WF(\Id-\tilde{\tilde{B}}_0)=\emptyset$. Without loss of generality assume $\chi,\tilde\chi\in S^0(\Omega)$ that $\nabla\chi, \nabla\tilde\chi$ is compactly supported, $\mc{V}_1=\mc{V}_1\chi$, $\tilde{\tilde{B}}_0=\chi \tilde{\tilde{B}}_0\chi$, $\supp\chi\cap\supp(1-\tilde\chi)=\emptyset$, $\supp\tilde\chi\cap\supp\kappa=\emptyset$. Let $\varphi\in C^\infty_c(\Omega)$ that $\supp\nabla\tilde\chi\cap \supp(1-\varphi)=\emptyset$, $\supp\chi\cap\supp\varphi=\emptyset$. 
    
    1. We iteratively claim
    \begin{equation}\label{e:5l13}
    \|A_0 \mc{V}_1\cdots\mc{V}_k u\|_{H^{N,l_0}}\le C\|B_0(\Lambda^\Omega_M-\om-i\mu)u\|_{H^{-N,l_0+k+0}}+C\|\tilde\kappa u\|_{H^{0}}+C\|u\|_{H^{-N, -N}}.
    \end{equation}
    The case $k=0$ follows from Lemma~\ref{l:EllEst_Near0}. Assume the case holds for up to $k-1$. Applying \eqref{e:EllEst_Tot} to $\mc{V}_1\cdots\mc{V}_k u$ implies
    \begin{multline}
        \|A_0 \mc{V}_1\cdots\mc{V}_k u\|_{H^{N,l_0}}\le C\|B_0\chi(P-\om-i\mu)\tilde\chi\mc{V}_1\mc{V}_2\cdots\mc{V}_ku\|_{H^{-N, l_0}}+C\|u\|_{H^{-N, -N}}\\
        \le C\|B_0\mc{V}_1\cdots\mc{V}_k\chi(\Lambda_M^\Omega-\om-i\mu)u\|_{H^{-N, l_0}}+C\|B_0[\chi P\tilde\chi, \mc{V}_1\cdots\mc{V}_k]u\|_{H^{-N, l_0}}\\
        +C\|\tilde B_0\mc{V}_1\cdots\mc V_k \mc{E}_M u\|_{H^{-N, l_0}}+C\|u\|_{H^{-N, -N}}\le C\|\tilde B_0(\Lambda^\Omega_M-\om-i\mu)u\|_{H^{-N,l_0+k+0}}\\
        +C\|B_0[\chi P\tilde\chi, \mc{V}_1\cdots\mc{V}_k]u\|_{H^{-N, l_0}}
        +C\|\tilde B_0\mc{V}_1\cdots\mc V_k \mc{E}_Mu\|_{H^{-N, l_0}}+C\|\tilde\kappa u\|_{H^{0}}+C\|u\|_{H^{-N, -N}}
    \end{multline}
    by the iterative assumption. 
    
    2. It remains to bound the commutator term and the $\mc{E}$-term. Decompose
    \begin{equation}
    [\chi P\tilde\chi,\mc{V}_1\cdots \mc{V}_k]=-\sum_{|\alpha|<k} \ad_\mc{V}^{\alpha^c}(\chi P\tilde\chi)\mc{V}_\alpha,
    \end{equation}
    as a summation over all proper subwords $\alpha$ of $(1,\dots, k)$, defined in \eqref{e:adword}. Noting $\ad_\mc{V}^{\alpha^c}(\chi P\tilde\chi)\in\Psi_{0,\delta}^{1,0}$ we have
    \begin{multline}
        \|B_0[\chi P\tilde\chi,\mc{V}_1\cdots \mc{V}_k] u\|_{H^{-N, l_0}}\le \sum_{|\alpha|<k}\|\tilde B_0 \mc{V}_\alpha u\|_{H^{-N, l_0}}+C\|u\|_{H^{-N, -N}}\\
        \le C\|\tilde{\tilde{B}}_0(\Lambda^\Omega_M-\om-i\mu)u\|_{H^{-N,l_0+k-1+0}}+C\|\tilde\kappa u\|_{H^{0}}+C\|u\|_{H^{-N, -N}}
    \end{multline}
    by the iterative assumption. For the other term, one can decompose 
    \begin{equation}
        B_0\mc{V}_1\cdots\mc{V}_k\chi \mc{E}_Mu=B_0\mc{V}_1\cdots\mc{V}_k\chi (\Lambda_M^{\Omega}-P)\tilde\chi u
        +B_0\mc{V}_1\cdots\mc{V}_k\chi \partial_\nu\mc{R}^{\Omega}_M[-\Delta, \tilde\chi]\mathcal{H}_M^\Omega u,
    \end{equation}
    as in Proposition~\ref{p:dtN}. Here Proposition~\ref{p:dtNParametrix} and the iterative assumption imply
    \begin{multline}
        \|B_0 \mc{V}_1\cdots\mc{V}_k\chi(\Lambda_M^{\Omega}-P)\tilde\chi u\|_{H^{-N,l_0}}\le C\|\Gamma_{k,\delta}^{-N, l_0+(1-\delta)k}(\Lambda_M^{\Omega}-P)\tilde\chi u\|\le C\|B_0\tilde\chi u\|_{H^{-N, l_0+0}}\\
        +C\|u\|_{H^{-N, -N}}\le C\|\tilde{\tilde{B}}_0(\Lambda^\Omega_M-\om-i\mu)u\|_{H^{-N,l_0+0}}+C\|\tilde\kappa u\|_{H^{0}}+C\|u\|_{H^{-N, -N}}
    \end{multline}
    by choosing $\delta\in [0,1)$ very close to $1$.  Lemma~\ref{l:smoothing} gives
    \begin{multline}\label{e:RH_Est}
        \|B_0\mc{V}_1\cdots\mc{V}_k\chi \partial_\nu \mc{R}^{\Omega}_M[-\Delta, \tilde\chi]\mc{H}_M^\Omega u\|_{H^{-N,l_0}}\le C\|\varphi\mc{H}_M^\Omega u\|_{X_{m,-N}^{l_0, -N}}\\
        \le C\|\varphi\mc{H}_M^\Omega \chi(1-\tilde\kappa) u\|_{X_{m,-N}^{l_0, -N}}+C\|\varphi\mc{H}_M^\Omega (1-\chi)(1-\tilde\kappa)u\|_{X_{m,-N}^{l_0, -N}}+C\|\varphi\mc{H}_M^\Omega \tilde\kappa u\|_{X_{m,-N}^{l_0, -N}}\\
        \le C\|B_0\chi u\|_{H^{-N, l_0}}+C\|u\|_{H^{-N, -N}}+C\|\tilde\kappa u\|_{H^{0, -N}},
    \end{multline}
    where we used Lemma~\ref{l:hnegbreg} and Lemma~\ref{l:longDistanceHarmonic} for the second and third terms on the last line. Invoke the iterative assumption and rename $\tilde{\tilde{B}}_0$ to $B_0$ to conclude \eqref{e:5l13}. 
\end{proof}

Note that $\Psi^{0,-\delta}_{1,\delta}\subset \Psi^{1,0}_{0,\delta}$
and $\Psi^{s,l}\subset \Psi^{s,l}_{0,\delta}$. When away from $\xi=0$, note $\Psi^{s,l}_{k,\delta}=\Psi^{s+k,l+k\delta}$. Here we record the following elliptic, propagation and source/sink estimates within the class $\Psi^{s,l}$ (see also \cite[(25),(26) and (27)]{DZ:19}). We denote the following error term, at the $0$-frequency or in the Lipschitz part, by 
\begin{equation}
    O_0=C(\|(1-\kappa)\psi_0(D)(1-\kappa_0)u\|_{H^{-N, -L_0+0}}+\|\kappa u\|_{H^{0}}+\|u\|_{H^{-N, -N}}),
\end{equation}
if for each $N>0$ and $\psi_0(\xi)\in C^\infty_c(\mathbb{R}^{d-1})$ with $\psi_0\equiv 1$ near $\xi=0$ there is such $C>0$ that does not depend on $u, \omega$ or $\mu$ in the following microlocal statements. 

\begin{lemma}[Elliptic estimate away from zero frequency]\label{l:EllEst_Away0}
Let $\Omega\subset \mathbb{R}^d$ be an \acl domain. Let $|l_0|<L_0$. Let $s,l\in\mathbb{R}$. For each $A, B\in \Psi^{0,0}(\partial\Omega)$, $\WF A\subset \Ell B$,
\begin{equation}
    \WF B\cap\Sigma(\omega)=\emptyset,\ \WF B\subset \Ell(1-\kappa_0), \ \WF B\cap\{(x,\xi):|\xi|_g\le \e_0\}=\emptyset,
\end{equation}
for some $\e_0>0$, for any $N>0$ we have
\begin{equation}\label{e:ResEst_Away0}
    \|A u\|_{H^{s,l}}\le C\|B(\Lambda^\Omega_M-\om-i\mu)u\|_{H^{s-1,l}}+O_0
\end{equation}
uniformly for all $u\in \mathcal{S}'$, $\omega\in\intw$ and $\mu\ge0$ small, $M\in[0, M_0]$. 
\end{lemma}
\begin{proof}
    Let $B=\chi B\chi$ for some $\chi,\tilde\chi\in S^0(\Omega)$ away from $\partial\osing$ with $\supp\chi\cap\supp(1-\tilde\chi)=\emptyset$. The standard elliptic estimate and Proposition~\ref{p:dtN} imply
    \begin{multline}
        \|A u\|_{H^{s,l}}\le C\|B\chi(P-\om-i\mu)\tilde\chi u\|_{H^{s-1,l}}+C\|u\|_{H^{-N, -N}}\\
        \le C\|B\chi(\Lambda^\Omega_M-\om-i\mu)u\|_{H^{s-1,l}}+C\|B\chi\mc{E}_M u\|_{H^{s-1,l}}+C\|u\|_{H^{-N, -N}},
    \end{multline}
    and $B\chi\mc{E}_M(1-\kappa)(\Id-\psi_0(D))(1-\kappa_0):H^{-N,-N}\rightarrow H^{N,N}$ is bounded since $\WF B\cap \{\xi=0\}=\emptyset$. Then
    \begin{equation}
        \|B\chi\mc{E}_M u\|_{H^{s-1,l}}\le C(\|(1-\kappa)\psi_0(D)(1-\kappa_0) u\|_{H^{-N, -L_0+0}}+\|u\|_{H^{-N, -N}}+\|\kappa u\|_{H^{0}})= O_0
    \end{equation}
    where we decomposed $\Id=(1-\kappa)\psi_0(D)(1-\kappa_0)+(1-\kappa)(\Id-\psi_0(D))(1-\kappa_0)+\kappa$.
\end{proof}
With exactly the same treatment of $\Lambda_M^\Omega-P$, we have
\begin{lemma}[Propagation Estimate]\label{l:PropEst}
    Let $\Omega\subset \mathbb{R}^d$ be an \acl domain. Let $s,l\in\mathbb{R}$. Let $A, B, Q\in \Psi^{0,0}$ with $\WF B\subset \Ell(1-\kappa_0)$, $\WF B\cap\{\xi=0\}=\emptyset$, $\WF A\cup\WF Q\subset \Ell B$. Moreover, assume there is $T>0$ such that
    \begin{equation}
        \WF A\subset \bigcup_{t\in[0,T]}e^{t\hat H_{p}}(\Ell Q), \ \bigcup_{t\in[0,T]}e^{t\hat H_{p}}(\WF Q)\subset \Ell(B).
    \end{equation}
    Then for each $N>0$, we have
    \begin{equation}
        \|Au\|_{H^{s,l}}\le C\|B (\Lambda^\Omega_M-\om-i\mu)u\|_{H^{s-1, l+1}}+C\|Q u\|_{H^{s,l}}+O_0,
    \end{equation}
    uniformly for all $u\in\mathcal{S}'$, $\omega\in\intw$, $M\in[0,M_0]$ and $\mu>0$ small. 
\end{lemma}

\begin{lemma}[Source estimate] \label{l:SourceEst}
    Let $\Omega\subset \mathbb{R}^d$ be an \acl domain. Let $s\in\mathbb{R}$, $l^->-\frac12$. Let $B_-\in \Psi^{0,0}$ such that $\WF B_-\subset \Ell(1-\kappa_0)$, $\WF B_-\cap\{\xi=0\}=\emptyset$, $\Sigma^-(\omega)\subset\Ell B_-$.
    There exists $A_-\in \Psi^{0,0}$ with $\Sigma^-(\omega)\subset \Ell A_-$ such that for any $N>0$ we have
    \begin{equation}
    \|A_-u\|_{H^{s, l^-}}\leq C\|B_-(\Lambda^\Omega_M-\om-i\mu) u\|_{H^{s-1, l^-+1}}+O_0, 
    \end{equation}
    uniformly for all $u\in \mathcal{S}'$ with $\WF^{s,-\frac12+0}(u)\cap \Sigma^-(\omega)=\emptyset$, $\omega\in\intw$, $M\in[0,M_0]$ and $\mu>0$ small. 
\end{lemma}

\begin{lemma}[Sink estimate]\label{l:SinkEst}
    Let $\Omega\subset \mathbb{R}^d$ be an \acl domain. Let $s\in\mathbb{R}$, $l^+<-\frac12$. Let $B_+\in \Psi^{0,0}$ such that $\WF B_+\subset \Ell(1-\kappa_0)$ $\WF B_+\cap\{\xi=0\}=\emptyset$, $\Sigma^+(\omega)\subset\Ell B_+$.
    There exist $A_+, Q\in \Psi^{0,0}$ with $\Sigma^+(\omega)\subset \Ell A_+$, $\WF Q\subset\Ell B_+\setminus \Sigma^+(\omega)$ such that for any $N>0$ we have 
    \begin{equation}
    \|A_+u\|_{H^{s, l^+}}\leq C\|B_+(\Lambda^\Omega_M-\om-i\mu) u\|_{H^{s-1, l^++1}}+C\|Qu\|_{H^{s,l^+}}+O_0,
    \end{equation}
    uniformly for all $u\in \mathcal{S}'$, $\omega\in\intw$, $M\in[0,M_0]$ and $\mu>0$ small. 
\end{lemma}

\subsection{Variable order Sobolev spaces and Fredholm estimates}\label{s:5-2}
In this subsection, we define the spaces on which the DtN map is Fredholm. Let $\omega>0$. Let $\kappa_1, \kappa_2, \kappa_3, \kappa_4 \in C_c^\infty(\overline{\Omega})$ with $\tilde\kappa\prec\kappa_1\prec\kappa_2\prec\kappa_3\prec\kappa_4$. Let $\psi_1, \tilde \psi_1$ in $C_c^\infty(\mathbb{R})$ such that $\psi_1(\xi)\equiv 1$ near $\xi=1$ and $\tilde\psi_1(\xi)\equiv 0$ near $\xi=0$, $\psi_1\prec\tilde\psi_1$. We rescale by $\psi_\omega, \tilde\psi_\omega\in C_c^\infty(\mathbb{R})$ that $\psi_\omega(\xi)=\psi_1(\omega^{-1}\xi)$, $\tilde\psi_\omega(\xi)=\tilde\psi_1(\omega^{-1}\xi)$. Let $\psi_{00}, \psi_0\in C_c^\infty(\mathbb{R})$ with $\supp \psi_0\cap \supp \psi_\omega=\emptyset$, $\psi_{00}\prec\psi_0$, $0\notin \supp(1-\psi_0)$. The choice of $\psi_{00}, \psi_0$ can be made independent of $\omega$ if we ask $\omega\ge \omega_0$ for some $\omega_0\ge 0$. Let $V^{\pm}(\omega)$ be neighborhoods of $\Sigma^{\pm}(\omega)$ in $T^*\partial\Omega$ and $\varphi\in C^\infty (\overline{T}^*\partial\Omega)$ such that $\supp (\pm 1-\varphi)\cap V^{\mp}(\omega)=\emptyset$, and $\hat H_{p}\varphi \leq 0$ on $\Sigma(\omega)$.

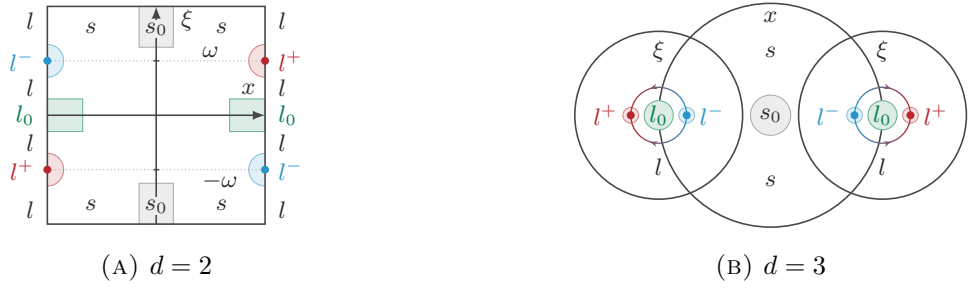
\begin{figure}[t]
  \centering
  \begin{subfigure}[t]{0.5\textwidth}
        \centering
        \begin{tikzpicture}[x=1cm,y=1cm]
  \definecolor{colorAKANE}{RGB}{183,40,46}
  \definecolor{colorTOKIWA}{RGB}{0,123,67}
  \definecolor{colorUSUAI}{RGB}{40,148,200}
  \colorlet{plotPrimary}{colorAKANE}
  \colorlet{plotSecondary}{colorTOKIWA}
  \colorlet{plotTertiary}{colorUSUAI}
  \colorlet{plotAxis}{black!72}
  \colorlet{plotGuide}{black!34}

  \tikzset{
    >=latex, font=\small\fontsize{9}{11}\selectfont,
    mark size=1.5pt, mark options={thin},
    line cap=round, line join=round,
    boundary/.style={draw=plotAxis,line width=0.6pt},
    box border/.style={line width=0.35pt,line join=miter},
    axis/.style={boundary,->},
    characteristic guide/.style={draw=plotGuide,line width=0.4pt,
      densely dotted},
    plus neighborhood/.style={draw=plotPrimary!60,fill=plotPrimary!14,
      line width=0.35pt},
    minus neighborhood/.style={draw=plotTertiary!60,fill=plotTertiary!14,
      line width=0.35pt},
    plus point/.style={circle,draw=none,fill=plotPrimary,
      inner sep=0pt,minimum size=3pt},
    minus point/.style={plus point,fill=plotTertiary},
    figure label/.style={text=black!84,inner sep=0.8pt},
    left weight/.style={figure label,anchor=east,xshift=-4pt},
    right weight/.style={figure label,anchor=west,xshift=4pt}
  }

  \path[use as bounding box] (-3.75,-1.6) rectangle (3.75,1.6);
  \begin{scope}[scale=0.72]
    \coordinate (A) at (-2,0);
    \coordinate (B) at (2,0);
    \coordinate (C) at (0,-2);
    \coordinate (D) at (0,2);
    \coordinate (upperLeft) at (-2,1);
    \coordinate (upperRight) at (2,1);
    \coordinate (lowerLeft) at (-2,-1);
    \coordinate (lowerRight) at (2,-1);
    \node[inner xsep=2pt,inner ysep=0pt] (compactExtent) at (0,0)
      {\phantom{$s_0$}};
    \coordinate (compactLeft) at (compactExtent.west |- A);
    \coordinate (compactRight) at (compactExtent.east |- A);
    \coordinate (compactTopInner) at (0,1.25);
    \coordinate (compactBottomInner) at (0,-1.25);

    \begin{scope}
      \path[box border,fill=black!7,draw=plotGuide]
        (compactLeft |- compactTopInner) rectangle (compactRight |- D)
        (compactLeft |- C) rectangle (compactRight |- compactBottomInner);
    \end{scope}

    \begin{scope}
      \path[box border,fill=plotSecondary!14,draw=plotSecondary!60]
        (-2,-0.30) rectangle (-1.35,0.30)
        (1.35,-0.30) rectangle (2,0.30);
    \end{scope}

    \begin{scope}[radius=0.30cm]
      \clip (A |- C) rectangle (B |- D);
      \path[minus neighborhood] (upperLeft) circle[];
      \path[plus neighborhood] (upperRight) circle[];
      \path[plus neighborhood] (lowerLeft) circle[];
      \path[minus neighborhood] (lowerRight) circle[];
    \end{scope}

    \begin{scope}
      \draw[boundary,line join=miter] (A |- C) rectangle (B |- D);
      \draw[characteristic guide] (upperLeft) -- (upperRight);
      \draw[characteristic guide] (lowerLeft) -- (lowerRight);
      \draw[axis] (A) -- (B);
      \draw[axis] (C) -- (D);
      \draw[boundary] (-0.05,1) -- (0.05,1);
      \draw[boundary] (-0.05,-1) -- (0.05,-1);
      \node[minus point] at (upperLeft) {};
      \node[plus point] at (upperRight) {};
      \node[plus point] at (lowerLeft) {};
      \node[minus point] at (lowerRight) {};
    \end{scope}

    \begin{scope}
      \foreach \yy in {-1.63,1.63} {
        \foreach \xx in {-1.20,1.20} {
          \node[figure label,anchor=mid] at (\xx,\yy) {$s$};
        }
        \node[figure label,anchor=mid,fill=black!7]
          at (0,\yy) {$s_0$};
      }
    \end{scope}

    \begin{scope}
      \foreach \yy in {-1.73,-0.50,0.50,1.73} {
        \node[left weight] at (-2,\yy) {$l$};
        \node[right weight] at (2,\yy) {$l$};
      }
      \node[left weight,text=plotSecondary] at (A) {$l_0$};
      \node[right weight,text=plotSecondary] at (B) {$l_0$};
      \node[left weight,text=plotTertiary] at (upperLeft) {$l^-$};
      \node[right weight,text=plotPrimary] at (upperRight) {$l^+$};
      \node[left weight,text=plotPrimary] at (lowerLeft) {$l^+$};
      \node[right weight,text=plotTertiary] at (lowerRight) {$l^-$};
    \end{scope}

    \begin{scope}
      \node[figure label,anchor=south west] at (0.80,1.02) {$\omega$};
      \node[figure label,anchor=north west] at (0.80,-1.01) {$-\omega$};
      \node[figure label,anchor=south east] at (1.87,0.35) {$x$};
      \node[figure label,anchor=north west,xshift=2pt]
        at (compactRight |- D) {$\xi$};
    \end{scope}
  \end{scope}
\end{tikzpicture}%
        \caption{$d=2$}
  \end{subfigure}\hfill
  \begin{subfigure}[t]{0.5\textwidth}
        \centering
        \begin{tikzpicture}[x=1cm,y=1cm]
  \definecolor{colorAKANE}{RGB}{183,40,46}
  \definecolor{colorTOKIWA}{RGB}{0,123,67}
  \definecolor{colorUSUAI}{RGB}{40,148,200}
  \colorlet{plotPrimary}{colorAKANE}
  \colorlet{plotSecondary}{colorTOKIWA}
  \colorlet{plotTertiary}{colorUSUAI}
  \colorlet{plotAxis}{black!72}
  \colorlet{plotGuide}{black!34}

  \pgfdeclarehorizontalshading{s5f4FlowLeft}{100bp}{
    color(0bp)=(plotPrimary);
    color(29.16667bp)=(plotPrimary);
    color(70.83333bp)=(plotTertiary);
    color(100bp)=(plotTertiary)}
  \pgfdeclarehorizontalshading{s5f4FlowRight}{100bp}{
    color(0bp)=(plotTertiary);
    color(29.16667bp)=(plotTertiary);
    color(70.83333bp)=(plotPrimary);
    color(100bp)=(plotPrimary)}

  \tikzset{
    >=latex, font=\small\fontsize{9}{11}\selectfont,
    mark size=1.5pt, mark options={thin},
    line cap=round, line join=round,
    boundary/.style={draw=plotAxis,line width=0.6pt},
    compact region/.style={draw=plotGuide,fill=black!7,line width=0.35pt},
    zero region/.style={draw=plotSecondary!60,fill=plotSecondary!14,
      line width=0.35pt},
    plus neighborhood/.style={draw=plotPrimary!60,fill=plotPrimary!14,
      line width=0.35pt},
    minus neighborhood/.style={draw=plotTertiary!60,fill=plotTertiary!14,
      line width=0.35pt},
    plus point/.style={circle,draw=none,fill=plotPrimary,
      inner sep=0pt,minimum size=3pt},
    minus point/.style={plus point,fill=plotTertiary},
    figure label/.style={text=black!84,inner sep=0.8pt},
    radial label/.style={figure label}
  }

  \path[use as bounding box] (-3.75,-1.6) rectangle (3.75,1.6);
  \begin{scope}[scale=1.48]
    \coordinate (spatialCenter) at (0,0);
    \coordinate (A) at (-1,0);
    \coordinate (B) at (1,0);
    \coordinate (leftPlus) at (-1.25,0);
    \coordinate (leftMinus) at (-0.75,0);
    \coordinate (rightMinus) at (0.75,0);
    \coordinate (rightPlus) at (1.25,0);

    \begin{scope}
      \draw[boundary] (spatialCenter) circle[radius=1cm];
      \draw[boundary] (A) circle[radius=0.75cm];
      \draw[boundary] (B) circle[radius=0.75cm];
    \end{scope}

    \begin{scope}
      \path[compact region] (spatialCenter) circle[radius=0.18cm];
      \path[zero region] (A) circle[radius=0.13cm];
      \path[zero region] (B) circle[radius=0.13cm];
    \end{scope}

    \begin{scope}
      \path[plus neighborhood] (leftPlus) circle[radius=0.07cm];
      \path[minus neighborhood] (leftMinus) circle[radius=0.07cm];
      \path[minus neighborhood] (rightMinus) circle[radius=0.07cm];
      \path[plus neighborhood] (rightPlus) circle[radius=0.07cm];
      \pgfmathsetlengthmacro{\flowOuterRadius}{0.25cm+0.25pt/1.48}
      \pgfmathsetlengthmacro{\flowInnerRadius}{0.25cm-0.25pt/1.48}
      \pgfmathsetlengthmacro{\flowHeadLength}{2.8pt/1.48}
      \pgfmathsetlengthmacro{\flowHeadHalfWidth}{1.1pt/1.48}
      \pgfmathsetlengthmacro{\flowHeadNotch}{0.675*\flowHeadLength}
      \foreach \flowCenter/\flowShading/\topAngle/\topHeading/\bottomAngle/\bottomHeading in {
        A/s5f4FlowLeft/108/198/-108/-198,
        B/s5f4FlowRight/72/-18/288/378} {
        \begin{scope}[shift={(\flowCenter)}]
          \begin{scope}[even odd rule]
            \clip
              (0,0) circle[radius=\flowOuterRadius]
              (0,0) circle[radius=\flowInnerRadius];
            \shade[shading=\flowShading]
              (-0.30,-0.30) rectangle (0.30,0.30);
          \end{scope}
          \foreach \flowAngle/\flowHeading in {
            \topAngle/\topHeading,\bottomAngle/\bottomHeading} {
            \begin{scope}
              \clip[shift={(\flowAngle:0.25)},rotate=\flowHeading]
                (0,0) -- (-\flowHeadLength,\flowHeadHalfWidth)
                -- (-\flowHeadNotch,0)
                -- (-\flowHeadLength,-\flowHeadHalfWidth) -- cycle;
              \shade[shading=\flowShading]
                (-0.30,-0.30) rectangle (0.30,0.30);
            \end{scope}
          }
        \end{scope}
      }
      \node[plus point] at (leftPlus) {};
      \node[minus point] at (leftMinus) {};
      \node[minus point] at (rightMinus) {};
      \node[plus point] at (rightPlus) {};
      \node[radial label,text=plotPrimary,anchor=east,xshift=-4pt]
        at (leftPlus) {$l^+$};
      \node[radial label,text=plotTertiary,anchor=west,xshift=4pt]
        at (leftMinus) {$l^-$};
      \node[radial label,text=plotTertiary,anchor=east,xshift=-4pt]
        at (rightMinus) {$l^-$};
      \node[radial label,text=plotPrimary,anchor=west,xshift=4pt]
        at (rightPlus) {$l^+$};
    \end{scope}

    \begin{scope}
      \node[figure label] at (spatialCenter) {$s_0$};
      \node[figure label,text=plotSecondary]
        at (A) {$l_0$};
      \node[figure label,text=plotSecondary]
        at (B) {$l_0$};
      \foreach \xx in {-1,1} {
        \node[figure label] at (\xx,-0.46) {$l$};
      }
      \node[figure label,anchor=mid] (upperOrder) at (0,0.58) {$s$};
      \node[figure label,anchor=mid] (lowerOrder) at (0,-0.58) {$s$};
    \end{scope}

    \begin{scope}
      \node[figure label] at (0,0.88) {$x$};
      \node[figure label] at (-1,0.57) {$\xi$};
      \node[figure label] at (1,0.57) {$\xi$};
    \end{scope}
  \end{scope}
\end{tikzpicture}%
        \caption{$d=3$}
  \end{subfigure}
  \caption{Variable order Sobolev spaces and regularity/decay indices.}
  \label{f:5f2}
\end{figure}

\begin{definition}[Variable-order Sobolev spaces on \acl domains]\label{d:Zspace}
Let $s,l,l_0,l^-,l^+\in\mathbb{R}$, and $|s_0|\le 1$. We define the logarithmic weight $W=\frac12(\Op(w)+\Op(w)^*)$ where
\begin{multline}\label{e:5l5}
w=(1-\kappa_4(x))\Big(\psi_\omega(|\xi|_g)\big(\tfrac{l^++l^-}{2}+\tfrac{l^--l^+}{2}\varphi\big) + l_0\psi_0(|\xi|_g) +l(1-\psi_0(|\xi|_g)-\psi_\omega(|\xi|_g))\Big)\log\langle x\rangle \\
+\big(s_0(1-\kappa_{0}(x))\kappa_4(x) +s(1-\kappa_4(x))\big)\log \langle \omega^{-1}\xi\rangle.
\end{multline}
Note $W$ is essentially self-adjoint on $L^2$, and $(e^W)^{-1}=e^{-W}$. Fix $N\ge\max(|s|,|s_0|, |l|,|l_0|,|l^-|,|l^+|)$ large. We define the variable-order Sobolev norm by 
\begin{equation}
\label{e:zSpace}
\|u\|_{Z^{s,l,s_0}_{l_0,l^-,l^+}}=\|\kappa_3 u\|_{H^{s_0}}+\| (1-\kappa_2)e^W(1-\kappa_2)u\|_{L^2}+\|(1-\kappa_1)u\|_{H^{-N,-N}}.
\end{equation}
For an \acl domain $\Omega$ with Lipschitz $\partial\Omega$, we define the variable-order Sobolev spaces by
\begin{equation}
    Z^{s,l,s_0}_{l_0,l^-,l^+}(\partial\Omega;\omega)=\Big\{u\in 
    \mathcal{S}'(\partial\Omega): \|u\|_{Z^{s,l,s_0}_{l_0,l^-,l^+}}<\infty \Big\}.
\end{equation}
We may drop $\partial\Omega$ and $\omega$ whenever there is no confusion. The definition of $Z^{s,l,s_0}_{l_0,l^-,l^+}$ only depends on the choice of $\kappa_4$, not on $\kappa_0$ nor $\kappa_1$. When $\Omega$ is \acs and $\partial\Omega$ is smooth, one can take $\kappa_0=\kappa_4\equiv 0$ and $w$ no longer depends on $s_0$. In this case, we abbreviate by
\begin{equation}
    Z^{s, l}_{l_0, l^-, l^+}(\partial\Omega)=Z^{s, l, 0}_{l_0, l^-, l^+}(\partial\Omega).
\end{equation}
See Figure~\ref{f:5f2}. 
\end{definition}

\begin{lemma}
    $(Z^{s,l,s_0}_{l_0,l^-,l^+}(\partial\Omega))^*=Z^{-s,-l,-s_0}_{-l_0,-l^-,-l^+}(\partial\Omega)$. 
\end{lemma} 
\begin{proof}
    Throughout this proof, $\langle\cdot, \cdot\rangle$ indicates the inner product in $L^2$ or $(L^2)^3$. For $t\in\mathbb{R}$, the cutoffs $\kappa_t\in C^\infty_c(\partial\Omega)$ are chosen compatible with $\kappa_0,\kappa_1, \kappa_2,\kappa_3, \kappa_4$ in the manner that $\supp\kappa_{t}\cap\supp(1-\kappa_{t'})=\emptyset$ if $t<t'$. 
    
    1. We show $Z^{-s,-l,-s_0}_{-l_0,-l^-,-l^+}\subset(Z^{s,l,s_0}_{l_0,l^-,l^+})^*$. For each $u\in Z^{-s,-l,-s_0}_{-l_0,-l^-,-l^+}$ and $\phi\in Z^{s,l,s_0}_{l_0,l^-,l^+}$, we claim
    \begin{equation}
    |\langle u,\phi\rangle|=|\langle \kappa_{2.2} u, \kappa_3 \phi\rangle + \langle (1-\kappa_{2.2}) u,\phi\rangle|\le C\|u\|_{Z^{-s,-l,-s_0}_{-l_0,-l^-,-l^+}}\|\phi\|_{Z^{s,l,s_0}_{l_0,l^-,l^+}}.
    \end{equation}
    Indeed, the first term can be estimated by
    \begin{equation}
        |\langle \kappa_{2.2} u, \kappa_3 \phi\rangle|\le C\|\kappa_3 u\|_{H^{-s_0}}\|\kappa_3 \phi\|_{H^{s_0}}\le C\|u\|_{Z^{-s,-l,-s_0}_{-l_0,-l^-,-l^+}}\|\phi\|_{Z^{s,l,s_0}_{l_0,l^-,l^+}}.
    \end{equation}
    The second term can be estimated by
    \begin{equation}
        |\langle (1-\kappa_{2.2}) u,\phi\rangle|=|\langle e^{-W}(1-\kappa_{2.2}) u,e^{W}(1-\kappa_{2.1}) \phi\rangle|\le C\|u\|_{Z^{-s,-l,-s_0}_{-l_0,-l^-,-l^+}}\|\phi\|_{Z^{s,l,s_0}_{l_0,l^-,l^+}},
    \end{equation}
    since $\WF(e^{-W}(1-\kappa_{2.2}))\subset\Ell((1-\kappa_{2})e^{-W}(1-\kappa_{2}))$ implies
    \begin{equation}
        \|e^{-W}(1-\kappa_{2.2}) u\|\le C\|(1-\kappa_{2})e^{-W}(1-\kappa_{2})u\|+C\|(1-\kappa_1)u\|_{H^{-N,-N}}\le C\|u\|_{Z^{-s,-l,-s_0}_{-l_0,-l^-,-l^+}},
    \end{equation}
    and similarly we have
    \begin{equation}
        \|e^{W}(1-\kappa_{2.1}) \phi\|\le C\|(1-\kappa_{2})e^{W}(1-\kappa_{2})\phi\|+C\|(1-\kappa_1)\phi\|_{H^{-N,-N}}\le C\|\phi\|_{Z^{s,l,s_0}_{l_0,l^-,l^+}}.
    \end{equation}
    This concludes the claim. 

    2. We show $(Z^{s,l,s_0}_{l_0,l^-,l^+})^*\subset Z^{-s,-l,-s_0}_{-l_0,-l^-,-l^+}$. Let $u\in(Z^{s,l,s_0}_{l_0,l^-,l^+})^*$ and $\phi\in Z^{s,l,s_0}_{l_0,l^-,l^+}$. Consider the embedding $\Phi: Z^{s,l,s_0}_{l_0,l^-,l^+}\leftrightarrow \mc{V}\subset \mc{H}=H^{s_0}\times L^2\times H^{-N,-N}$, where
    \begin{equation}
        \Phi(\phi)=(\kappa_3\phi, (1-\kappa_2)e^W(1-\kappa_2)\phi, (1-\kappa_1)\phi), \ \mc{V}=\operatorname{Ran}(\Phi).
    \end{equation}
    One can check $\Phi$ is injective and $\mc{V}=\operatorname{Ran}(\Phi)$ is closed: the invertibility of $\Phi$ then follows from the open mapping theorem. The functional $u\in(Z^{s,l,s_0}_{l_0,l^-,l^+})^*$ then admits a unique lift $\tilde u\in \mc{V}^*$ given by
    \begin{equation}
        \tilde u(\tilde \phi)=u(\Phi^{-1}(\tilde \phi)), \ \|\tilde u\|_{\mc{V}^*}\le C\|u\|_{(Z^{s,l,s_0}_{l_0,l^-,l^+})^*},
    \end{equation}
    for each $\tilde\phi\in \mc{V}$. Hahn-Banach theorem implies $\tilde u\in \mc{V}^*$ extends to an element in $\mc{H}^*$ represented by $\tilde{\tilde{u}}=(\tilde{\tilde{u}}_1, \tilde{\tilde{u}}_2, \tilde{\tilde{u}}_3)\in H^{-s_0}\times L^2 \times H^{N, N}$, such that for each $\tilde\phi\in \mc{V}$, 
    \begin{equation}\label{e:5l8}
        \langle \tilde{\tilde{u}}, \tilde\phi\rangle=\tilde u(\tilde \phi)=u(\Phi^{-1}(\tilde \phi)), \ \|\tilde{\tilde{u}}\|_{H^{-s_0}\times L^2 \times H^{N, N}}=\|\tilde u\|_{\mc{V}^*}\le C\|u\|_{(Z^{s,l,s_0}_{l_0,l^-,l^+})^*}. 
    \end{equation}
    Now we set the representation
    \begin{equation}
        w=\kappa_3\tilde{\tilde{u}}_1+(1-\kappa_2)e^W(1-\kappa_2)\tilde{\tilde{u}}_2+ (1-\kappa_1)\tilde{\tilde{u}}_3. 
    \end{equation}
    Since $\langle D\rangle^{-s_0}\kappa_3(1-\kappa_2)e^W(1-\kappa_2)\in \Psi^{0, -\infty}$, we have
    \begin{equation}
        \|\kappa_3 w\|_{H^{-s_0}}\le C\|\tilde{\tilde{u}}_1\|_{H^{-s_0}}+C\|\tilde{\tilde{u}}_2\|_{L^2}+\|\tilde{\tilde{u}}_3\|_{H^{N, N}}.
    \end{equation}
    We can similarly estimate $\|(1-\kappa_2)e^{-W}(1-\kappa_2) w\|$ and $\|(1-\kappa_1)w\|_{H^{-N,-N}}$ to observe 
    \begin{equation}
        \|w\|_{Z^{-s,-l,-s_0}_{-l_0,-l^-,-l^+}}\le C\|\tilde{\tilde{u}}_1\|_{H^{-s_0}}+C\|\tilde{\tilde{u}}_2\|_{L^2}+\|\tilde{\tilde{u}}_3\|_{H^{N, N}}.
    \end{equation}
    Setting $\tilde\phi=\Phi(\phi)$ in \eqref{e:5l8}, we end up with
    \begin{equation}
        u(\phi)=\langle \tilde{\tilde{u}}, \Phi(\phi)\rangle=\langle \tilde{\tilde{u}}_1, \kappa_3\phi\rangle+\langle \tilde{\tilde{u}}_2, (1-\kappa_2)e^W(1-\kappa_2)\phi\rangle+\langle \tilde{\tilde{u}}_3, (1-\kappa_1)\phi\rangle=\langle w, \phi\rangle,
    \end{equation}
    that is we have represented $u$ by $w\in Z^{-s,-l,-s_0}_{-l_0,-l^-,-l^+}$. This concludes the proof.
\end{proof}

\begin{lemma}\label{l:dtnfredbdd}
    Let $\Omega$ be an \acl domain. Let $M_0>0$. Let $s,l, l^-, l^+\in\mathbb{R}$. Let $s_0\in [0,1]$, $|l_0|<L_0$. Then 
    \begin{equation}
        \Lambda^\Omega_M: Z^{s,l,s_0}_{l_0,l^-,l^+}(\partial\Omega;\omega)\rightarrow Z^{s-1,l,s_0-1}_{l_0,l^-,l^+}(\partial\Omega;\omega)
    \end{equation}
    is continuous and uniformly bounded in $M\in[0,M_0]$ as bounded maps.
\end{lemma}
\begin{proof}
0. In this proof, for $t\in\mathbb{R}$, the cutoffs $\kappa_t\in C^\infty_c(\partial\Omega)$ are chosen compatible with $\kappa_0,\kappa_1, \kappa_2,\kappa_3, \kappa_4$ in the manner that $\supp\kappa_{t}\cap\supp(1-\kappa_{t'})=\emptyset$ if $t<t'$. 

1. Through Step 1--4, we show
\begin{equation}
    \|\Lambda_M^\Omega u\|_{Z^{s-1,l,s_0-1}_{l_0,l^-,l^+}}\le C\|u\|_{Z^{s,l,s_0}_{l_0,l^-,l^+}}
\end{equation}
uniformly for all $M\in [0, M_0]$. All estimates in this proof are uniform in this way. 

2a. We estimate $\kappa_3\Lambda^\Omega_M u$ in this step, by claiming
\begin{equation}
    \kappa_3 \Lambda^\Omega_M \kappa_{3.5}=\kappa_3\Lambda^{\tilde\Omega}_M\kappa_{3.5}-\kappa_{3}\partial_\nu\mc{R}^\Omega_M\varphi[-\Delta, \kappa_{3.6}]\mc{H}_M^{\tilde\Omega}\kappa_{3.5}: H^{s_0}\rightarrow H^{s_0-1},
\end{equation}
where we take a bounded uniformly Lipschitz domain $\tilde\Omega$ that $\tilde\Omega, \Omega$ share the common boundary near $\supp\kappa_4\cap \partial\Omega$, and $\varphi\in C^\infty_c(\overline{\tilde\Omega})$ that $\supp\nabla\kappa_{3.6}\cap\supp(1-\varphi)=\supp\varphi\cap\supp \kappa_3=\emptyset$. Only the second term needs to be justified here. Lemma~\ref{l:reg2} implies
\begin{equation}
    [-\Delta, \kappa_{3.6}]\mc{H}_M^{\tilde\Omega}: H^{s_0}(\partial\tilde\Omega)\rightarrow Y_{N}^N(\Omega). 
\end{equation}
Noting $\kappa_{3}\varphi=0$, we have
\begin{equation}
    \partial_\nu\kappa_3\mc{R}_M^\Omega\varphi=\partial_\nu\mc{R}_M^{\tilde\Omega}[-\Delta, \kappa_3]\mc{R}_M^\Omega\varphi: Y_N^N(\Omega)\rightarrow H^{0}(\partial\Omega),
\end{equation}
where we used Proposition~\ref{p:extendNoMass} and $\mc{R}_M^{\tilde\Omega}: L^2(\tilde\Omega)\rightarrow \dot H^{\frac32}(\tilde\Omega)$ from Lemma~\ref{l:reg2}. This concludes the step. 

2b. Let $\varphi'\in C^\infty_c(\overline{\Omega})$ that $\supp\nabla\kappa_{3.1}\cap\supp(1-\varphi')=\supp\varphi'\cap\supp\kappa_3=\emptyset$. Proposition~\ref{p:dtN} implies
\begin{equation}
     \kappa_3\Lambda^\Omega_M(1-\kappa_{3.5})=\kappa_3\tilde{\mc{E}}_M^*(1-\kappa_{3.5})=\kappa_3\partial_\nu\mc{R}^{\Omega}_M\varphi'[\Delta, 1-\kappa_{3.1}]\mc{H}^{\Omega}_M(1-\kappa_{3.5}).
\end{equation}
Lemma~\ref{l:ellregunbdd} implies $\varphi'\mc{H}_M^\Omega \kappa_3: H^{0, -N}(\partial\Omega)\rightarrow \dot H^{\frac32, N}(\Omega)$ bounded and so is its adjoint $-\kappa_3\partial_\nu\mc{R}_M^\Omega\varphi': H^{-\frac32, -N}(\Omega)\rightarrow H^{0, N}(\partial\Omega)$. Away from the zero-frequency, note
\begin{equation}
\kappa_3\partial_\nu\mc{R}^{\Omega}_M\varphi'[\Delta, 1-\kappa_{3.1}]\mc{H}^{\Omega}_M(1-\kappa_{3.2})(1-\psi_0(D))(1-\kappa_{3.5}): H^{-N, -N}\rightarrow H^{0,N}
\end{equation}
is bounded, since $[\Delta, 1-\kappa_{3.1}]\mc{H}^{\Omega}_M(1-\kappa_{3.2})(1-\psi_0(D))(1-\kappa_{3.5}): H^{-N, -N}(\partial\Omega)\rightarrow Y_{\rardyr_\Omega-0, N}^{L_0-0, N}(\Omega)$ bounded per \eqref{e:thing8}. At the zero-frequency, we have 
\begin{equation}
    \kappa_3\partial_\nu\mc{R}^{\Omega}_M\varphi'[\Delta, 1-\kappa_{3.1}]\mc{H}^{\Omega}_M(1-\kappa_{3.2})\psi_0(D)(1-\kappa_{3.5}): H^{-N, -L_0+0}\rightarrow H^{0,N},
\end{equation}
bounded since $\mc{H}_M^{\Omega}: H^{\frac12, -L_0+0}\rightarrow X_{-\rardyr_{\Omega}+0}^{-L_0+0}(\Omega)$ per Proposition~\ref{p:harmonic}, and $\varphi'\mc{H}_M^{\Omega}\kappa_3$. Note this implies
\begin{multline}
    \| \kappa_3\partial_\nu\mc{R}^{\Omega}_M\varphi'[\Delta, 1-\kappa_{3.1}]\mc{H}^{\Omega}_M(1-\kappa_{3.2})\psi_0(D)(1-\kappa_{3.5})u\|_{H^{0}}\le C\|(1-\kappa_{4.1})u\|_{H^{-N, -L_0+0}}\\
    +C\|(1-\kappa_1)u\|_{H^{-N, -N}}.
\end{multline}
Bringing the estimates together we have
\begin{equation}
    \|\kappa_3\Lambda^\Omega_M u\|_{H^{s_0-1}}\le C\|\kappa_{3.9} u\|_{H^{s_0}}+C\|(1-\kappa_{4.1})u\|_{H^{-N, -L_0+0}}
    +C\|(1-\kappa_1)u\|_{H^{-N, -N}}\le C\|u\|_{Z^{-N, -N, s_0}_{-L_0+0, -N, -N}}
\end{equation}
as desired. 

3a. Let the weights $W$ be with parameter $(s-1,l,s_0-1,l_0,l^-,l^+)$ and $W'$ be with parameter $(s,l,s_0,l_0,l^-,l^+)$ in \eqref{e:5l5}. Denote $A=(1-\kappa_2)e^W(1-\kappa_2)$ and $A'=(1-\kappa_2)e^{W'}(1-\kappa_2)$. We begin with the first term of 
\begin{equation}
    A \Lambda_M^\Omega=A P e^{-W'}e^{W'}(1-\kappa_{1.9})+A \mc{E},
\end{equation}
given by Proposition~\ref{p:dtN}. Note $P\in \Psi^{0, -\delta}_{1,\delta}\subset\Psi^{1,0}_{0,\delta}$, $APe^{-W'}\in \Psi^{0,0}_{0,\delta}$ to see 
\begin{equation}
    \|AP(1-\kappa_{1.9}) u\|\le C\|e^{W'}(1-\kappa_{1.9}) u\|\le 
    C\|A' u\|+C\|\kappa_3 u\|_{H^{s_0}}+C\|(1-\kappa_1)u\|_{H^{-N,-N}}.   
\end{equation}

3b. Denote $B_{00}=(1-\kappa_{2.4})\psi_{00}(D)(1-\kappa_{2.5})$. Then
\begin{equation}
    A\mc{E}=A\mc{E}\kappa_{2.5}+A\mc{E}(1-\kappa_{2.5}-B_{00})+A\mc{E}B_{00}.
\end{equation}
Proposition~\ref{p:dtN} allows us to bound the error terms by
\begin{gather}
    \|A \mc{E}\kappa_{2.5} u\|\le C\|\kappa_{2.5} u\|_{H^{0,-N}},\\
    \|A \mc{E}(1-\kappa_{2.5}-B_{00}) u\|\le C\|A'u\|_{H^{-N,-N}}+C\|(1-\kappa_{1}) u\|_{H^{-N,-N}},\\
    \|A\mc{E}B_{00}u\|\le C\|A' u\|_{H^{-N,0}}+C\|(1-\kappa_{1}) u\|_{H^{-N,-N}}.
\end{gather}
Together with Step 3a we have
\begin{equation}
    \|A\Lambda_M^\Omega u\|\le C\|A' u\|+C\|\kappa_3 u\|_{H^{s_0}}+C\|(1-\kappa_1)u\|_{H^{-N,-N}}\le C\|u\|_{Z^{s,l,s_0}_{l_0, l^-,l^+}}.
\end{equation}

4. Step 2 and 3 imply
\begin{equation}
    \|(1-\kappa_1)\Lambda_M^\Omega u\|_{H^{-N,-N}}\le \|(1-\kappa_1)\kappa_3\Lambda_M^\Omega u\|_{H^{-N,-N}}+\|(1-\kappa_3)\Lambda_M^\Omega u\|_{H^{-N,-N}}\le C\|u\|_{Z^{s,l,s_0}_{l_0, l^-,l^+}}.
\end{equation}
This concludes the claim in Step 1. 

5. That $\Lambda_M^\Omega$ is continuous in $M$ follows from the continuity of all objects we used in the proof, established in earlier sections.
\end{proof}

We summarize the microlocal estimates in terms of the variable-order Sobolev spaces.
\begin{lemma}\label{l:microfred}
    Let $s,l,l_0, l^\pm\in\mathbb{R}$, $s_0\in[0,1]$. Let $A, A_0, B_0, B, B^\pm, Q\in \Psi^{0,0}(\partial\Omega)$. There exist some $C, \e_0, T>0$ that the following statements are true for any $N>0$ and uniformly for all $u\in \mathcal{S}'$, $\omega\in\intw$ and $\mu\ge0$ small, $M\in[0, M_0]$:
\begin{enumerate}[wide]
    \item Let $|l_0|<L_0$. If $\{(x,\xi): x\in \supp(1-\kappa), \xi=0\}\subset\Ell B_0$, $\WF(B_0)\subset\{(x,\xi):|\xi|_g\le \e_0\}$, $\WF A_0\cap \Sigma(\omega)=\emptyset$, $\WF A_0\subset \Ell B_0$, then
\begin{equation}\label{e:microfred0}
    \|A_0 u\|_{Z^{N,l,s_0}_{l_0,l^-,l^+}}\le C\|B_0(\Lambda^\Omega_M-\om-i\mu)u\|_{Z^{-N,l,s_0-1}_{l_0,l^-,l^+}}+C\|\tilde\kappa u\|_{H^{0}}+C\|u\|_{H^{-N, -N}}.
\end{equation}
uniformly for all $u$ described above with $\WF^{-\infty,-L_0+0}(u)\cap \{\xi=0\}=\emptyset$.

    \item If $\WF B\subset \Ell(1-\kappa_0)$, $\WF B\cap(\Sigma(\omega)\cup \{(x,\xi):|\xi|_g\le \e_0\})=\emptyset$, $\WF A\subset \Ell B$, then
    \begin{equation}\label{e:microfredaway0}
    \|A u\|_{Z^{s,l,s_0}_{l_0,l^-,l^+}}\le C\|B(\Lambda^\Omega_M-\om-i\mu)u\|_{Z^{s-1,l,s_0-1}_{l_0,l^-,l^+}}+O_0. 
    \end{equation}

    \item If $\WF B\subset \Ell(1-\kappa_0)$, $\WF B\cap\{\xi=0\}=\emptyset$, $\WF A\cup\WF Q\subset \Ell B$, and
    \begin{equation}
        \WF A\subset \bigcup_{t\in[0,T]}e^{t\hat H_{p}}(\Ell Q), \ \bigcup_{t\in[0,T]}e^{t\hat H_{p}}(\WF Q)\subset \Ell(B)
    \end{equation}
    for some $T>0$, then
    \begin{equation}\label{e:microfredprop}
        \|Au\|_{Z^{s,l,s_0}_{l_0,l^-,l^+}}\le C\|B (\Lambda^\Omega_M-\om-i\mu)u\|_{Z^{s-1,l,s_0-1}_{l_0,l^-+1,l^++1}}+C\|Q u\|_{Z^{s,l,s_0}_{l_0,l^-,l^+}}+O_0.
    \end{equation}

    \item Let $l^->-\frac12$. If $\WF B^-\subset \Ell(1-\kappa_0)$, $\WF B^-\cap\{\xi=0\}=\emptyset$, $\Sigma^-(\omega)\subset\Ell B^-$, then there exists $A^-\in \Psi^{0,0}$ with $\Sigma^-(\omega)\subset \Ell A^-$ that
    \begin{equation}\label{e:microfredsource}
    \|A^-u\|_{Z^{s,l,s_0}_{l_0,l^-,l^+}}\leq C\|B^- (\Lambda^\Omega_M-\om-i\mu)u\|_{Z^{s-1,l,s_0-1}_{l_0,l^-+1,l^++1}}+O_0, 
    \end{equation}
    uniformly for all $u$ described above with $\WF^{s,-\frac12+0}(u)\cap \Sigma^-(\omega)=\emptyset$.

    \item Let $l^+<-\frac12$. If $\WF B^+\subset \Ell(1-\kappa_0)$, $\WF B^+\cap\{\xi=0\}=\emptyset$, $\Sigma^+(\omega)\subset\Ell B^+$ then there exist $A^+, Q^+\in \Psi^{0,0}$ with $\Sigma^+(\omega)\subset \Ell A^+$, $\WF Q^+\subset\Ell B^+\setminus \Sigma^+(\omega)$ that
    \begin{equation}\label{e:microfredsink}
    \|A^+u\|_{Z^{s,l,s_0}_{l_0,l^-,l^+}}\leq C\|B^+(\Lambda^\Omega_M-\om-i\mu) u\|_{Z^{s-1,l,s_0-1}_{l_0,l^-+1,l^++1}}+C\|Q^+u\|_{Z^{s,l,s_0}_{l_0,l^-,l^+}}+O_0.
    \end{equation}
\end{enumerate}
In the above estimates, one can replace $-i\mu$ by $+i\mu$ by switching $Z^{s,l,s_0}_{l_0,l^-,l^+}$ and $Z^{s,l,s_0}_{l_0,l^+,l^-}$, $\Sigma^-(\omega)$ and $\Sigma^+(\omega)$, that is, switching the source and sinks. 
\end{lemma}

Combining the estimates above, we can write down the Fredholm estimates in the variable order Sobolev spaces. 
\begin{proposition}[Fredholm estimates]\label{p:FredEst}
Let $\Omega$ be an \acl domain. Let $M_0>0$. Let $s,l\in\mathbb{R}$, $s_0\in [0,1]$. Let $|l_0|<L_0$. Then for any $N>0$ we have the following estimates:
\begin{enumerate}[wide]
    \item For $l^->-\frac12$, $l^+<-\frac12$, we have
    \begin{equation}
    \|u\|_{Z^{s,l, s_0}_{l_0,l^-,l^+}}\leq C\|(\Lambda^\Omega_M-\om-i\mu) u\|_{Z^{s-1,l, s_0-1}_{l_0,l^-+1,l^++1}}+C\|\tilde\kappa u\|_{L^2}+C\|u\|_{H^{-N,-N}},
    \end{equation}
    uniformly for all $u\in Z^{-\infty, -\infty, 0}_{-L_0+0, -\frac12+0, -\infty}(\partial\Omega)$, $\omega\in\intw$, $\mu\ge 0$ small, $M\in [0, M_0]$. 
    \item For $l^-<-\frac12$, $l^+>-\frac12$, we have
    \begin{equation}
        \|u\|_{Z^{s,l, s_0}_{l_0,l^-,l^+}}\leq C\|(\Lambda^\Omega_M-\om+i\mu) u\|_{Z^{s-1,l, s_0-1}_{l_0,l^-+1,l^++1}}+C\|\tilde\kappa u\|_{L^2}+C\|u\|_{H^{-N,-N}},\label{e:5l9}
    \end{equation}
    uniformly for all $u\in Z^{-\infty, -\infty, 0}_{-L_0+0, -\infty, -\frac12+0}(\partial\Omega)$, $\omega\in\intw$, $\mu\ge 0$ small, $M\in [0, M_0]$. 
    \item When $s_0\in(0,1]$, the $\|\tilde\kappa u\|_{L^2}$ terms can be dropped from the above estimates.
\end{enumerate}
\end{proposition}
\begin{proof}
    1. Let $f=(\Lambda^\Omega_M-\om-i\mu) u$. Lemma~\ref{l:microfred} implies
    \begin{equation}
        \|(1-\kappa)u\|_{Z^{s,l,s_0}_{l_0,l^-,l^+}}\le C\|(1-\kappa_0) f\|_{Z^{s-1,l,s_0-1}_{l_0,l^-+1,l^++1}}+O_0
    \end{equation}
    by applying \eqref{e:microfredsource}, \eqref{e:microfredprop}, \eqref{e:microfredsink}, \eqref{e:microfredaway0} and \eqref{e:microfred0} sequentially. Invoke Lemma~\ref{l:EllEst_Near0}, the elliptic estimate at the zero frequency, to see for $O_0$, 
    \begin{equation}
        \|(1-\kappa)\psi_0(D)(1-\kappa_0) u\|_{H^{N,l_0}}\le C\|(1-\kappa_0)f\|_{H^{-N, l_0}}+C\|\tilde\kappa u\|_{H^{0}}+C\|u\|_{H^{-N,-N}}.
    \end{equation}
    Thus we end up with
    \begin{equation}\label{e:3l25}
        \|(1-\kappa)u\|_{Z^{s,l,s_0}_{l_0,l^-,l^+}}\le C\|(1-\kappa_0) f\|_{Z^{s-1,l,s_0-1}_{l_0,l^-+1,l^++1}}+C\|\tilde\kappa u\|_{H^{0}}+C\|u\|_{H^{-N, -N}}.
    \end{equation}

    2. We now estimate $\|\tilde\kappa u\|_{Z^{s,l,s_0}_{l_0,l^-,l^+}}$. Note $\|\tilde\kappa u\|_{Z^{s,l,s_0}_{l_0,l^-,l^+}}\sim\|\tilde\kappa u\|_{H^{s_0}}$. By Lemma~\ref{l:EllEst_ACLip} and \eqref{e:3l25}, we have
    \begin{multline}
        \|\tilde \kappa u\|_{H^{s_0}}\le C\|\tilde\kappa f\|_{H^{s_0-1}}+C\|(1-\kappa_1)u\|_{H^{-N,-N}}+C\|(1-\kappa_2)\psi_0(D)(1-\kappa_1) u\|_{H^{-N, -L_0+0}}\\+C\|\kappa_2 u\|_{H^{0}}
        \le C\|f\|_{Z^{s-1,l, s_0-1}_{l_0,l^-+1,l^++1}}+C\|(\kappa_2-\tilde\kappa )u\|_{H^{0}}+C\e^{-1}\|\tilde\kappa u\|_{H^{-N}}+\e\|\tilde \kappa u\|_{H^{s_0}}
    \end{multline}
    where we used interpolation
    \begin{equation}
        \|\tilde\kappa u\|_{H^{0}}\le C\e^{-1}\|\tilde\kappa u\|_{H^{-N}}+\e\|\tilde\kappa u\|_{H^{s_0}}.
    \end{equation}
    Noting $\kappa_2-\tilde\kappa$ is supported away from $\supp\kappa$, we have from \eqref{e:3l25} and interpolation that 
    \begin{equation}
        \|(\kappa_2-\tilde\kappa )u\|_{H^{0}}\le C\|(1-\kappa_0) f\|_{Z^{s-1,l, s_0-1}_{l_0,l^-+1,l^++1}}+C\e^{-1}\|\tilde\kappa u\|_{H^{-N}}+\e\|\tilde\kappa u\|_{H^{s_0}}.
    \end{equation}
    Choose small $\e>0$ to allow the absorption to see
    \begin{equation}
        \|\tilde \kappa u\|_{H^{s_0}}\le C\|f\|_{Z^{s-1,l, s_0-1}_{l_0,l^-+1,l^++1}}+C\|u\|_{H^{-N, -N}},
    \end{equation}
    and $O_0$ can be bounded by the same right hand side. Revisit \eqref{e:3l25} to conclude the proof. 

    3. Equation \eqref{e:5l9} follows similarly by switching the source and the sink. 
\end{proof}

\begin{lemma}[Absorption resolvent]\label{l:absorbres}
Let $\Omega$ be an \acl domain. Let $M_0>0$. Let $s,l,l^-,l^+\in\mathbb{R}$, $s_0\in (0,1)$, $|l_0|<L_0$. Fix $\mu>0$. Then for any $N>0$ we have
    \begin{equation}\label{e:5l4}
        \|u\|_{Z^{s,l, s_0}_{l_0,l^-,l^+}}\leq C_\mu\|(\Lambda^\Omega_M-\om\mp i\mu) u\|_{Z^{s-1,l, s_0-1}_{l_0,l^-,l^+}}+C_\mu\|u\|_{H^{-N,-N}},
    \end{equation}
    uniformly for all $u\in Z^{-\infty, -\infty, 0}_{-L_0+0, -\frac12+0, -\infty}$, $\omega\in\intw$, $M\in [0, M_0]$. Moreover, for each $\mu>0$, 
    \begin{equation}
        (\Lambda^\Omega_M-\om\mp i\mu)^{-1}: Z^{s-1,l, s_0-1}_{l_0,l^-,l^+}\rightarrow Z^{s,l, s_0}_{l_0,l^-,l^+}
    \end{equation}
    is bounded. 
\end{lemma}
\begin{proof}
    1. We first briefly discuss \eqref{e:5l4}. For fixed $\mu>0$, $P-\om\mp i\mu$ is elliptic. One can prove a similar version of Lemma~\ref{l:EllEst_Away0} away from the zero frequency. Together with Lemma~\ref{l:EllEst_Near0}, Lemma~\ref{l:EllEst_ACLip},
    one can show \eqref{e:5l4} by repeating the proof of Proposition~\ref{p:FredEst} without invoking any propagation or radial point estimates. 

    2. We show $\Lambda^\Omega_M-\om- i\mu$ is invertible. Equation~\eqref{e:5l4} implies it has closed range. It suffices to check it has trivial kernel and cokernel. Let $(\Lambda^\Omega_M-\om-i\mu)u=0$ for $u\in Z^{s,l,s_0}_{l_0,l^-, l^+}$. Then \eqref{e:5l4} implies $u\in L^2$. Noting $\Lambda_M^\Omega$ is self-adjoint on $L^2$, we have
    \begin{equation}
        \mu\|u\|^2=-\Im\langle(\Lambda^\Omega_M-\om-i\mu) u, u\rangle=0,
    \end{equation}
    yielding $u=0$ and $\Lambda^\Omega_M-\om-i\mu$ has trivial kernel. The trivial cokernel and the case of $\Lambda^\Omega_M-\om+ i\mu$ follow similarly.
\end{proof}

\subsection{Limiting absorption principle}
We show that any outgoing resonant state is an eigenfunction for $\Lambda^\Omega_M$. 
\begin{lemma}[No purely outgoing states]\label{l:EigenReg}
    Let $\omega\in \intw$, $s_0\in[0,1]$. Assume on $\partial\Omega$ we have
    \begin{equation}
        u\in Z^{-\infty, -\infty, 0}_{-L_0+0, -\frac12+0, -\infty}, \ (\Lambda^\Omega_M-\om)u\in Z^{\infty, \infty, s_0-1}_{L_0-0, \infty, \infty}, \ \Im\langle (\Lambda^\Omega_M-\om)u, u\rangle\geq 0.
    \end{equation}
    Then we have
    \begin{equation}
        u\in Z^{\infty, \infty, s_0}_{L_0-0, \infty,\infty}(\partial\Omega)\subset L^2(\partial\Omega).
    \end{equation}
\end{lemma}
\begin{proof}
    1. The Fredholm estimate, Proposition~\ref{p:FredEst}, together with the assumptions on $u$ and $(\Lambda_M^{\Omega}-\sqrt{\omega^2+M^2})u$, implies $u\in Z^{\infty, \infty, s_0}_{L_0-0, \infty, -\frac12-0}$. It remains to show $\WF u\cap \Sigma^+(\omega)=\emptyset$ under the extra sign condition. To do this, we follow the idea in \cite[Lemma~2.3]{Dyatlov2017RuelleSurfaces}. Throughout the proof, we will use the coordinate $(\rho, y, \xi)$ where $\rho$ is the boundary defining function of the scattering manifold and $y$ the cross-section variable, and $\xi$ is the dual variable to $x$. 
    
    Let $\varphi_0, \varphi, \tilde\varphi\in C^\infty_c((0, \infty)_\rho)$ that $[\frac19, 9]\cap\supp(1-\varphi)=\supp\varphi_0\cap\supp(1-\varphi)=\supp\varphi\cap\supp(1-\tilde\varphi)=\emptyset$. Let $\chi,\tilde\chi, \tilde{\tilde\chi}\in C^\infty_c([0, \infty)_\rho)$ with 
    $\supp\tilde\varphi\cap\supp(1-\chi)=\supp\chi\cap\supp(1-\tilde\chi)=\supp\tilde\chi\cap\supp(1-\tilde{\tilde\chi})=\emptyset$. Fix $\e_0>0$ and denote by $S^{s,\comp}$ those symbols in $S^{s,-\infty}$ with compact support in $x$, and $S^{\comp, l}$ analogously. We will prove that for all $a(x,\xi)\in S^{0,0}$, with $\supp a\subset  \{(x,\xi): x\notin\supp\kappa, |\xi|> \e_0\}$, there exists $b\in S^{0, 0}$ with $\supp b\subset  \{(x,\xi): x\notin\supp\kappa, |\xi|> \e_0\}$ such that for any $N>0$, uniformly for $h>0$ small, 
    \begin{equation} \label{e:SmoothEst}
        \|\varphi_0(h^{-1}\rho)\Op (a)u\|\leq Ch^{\frac12}\|\tilde\varphi(h^{-1}\rho)\Op (b)u\|+O(h^N\|(\Lambda^\Omega_M-\om)u\|_{Z^{\infty, \infty, s_0}_{L_0-0, \infty, \infty}}).
    \end{equation}
    Once we showed it, the iterative use of \eqref{e:SmoothEst} gives for any $N>0$,
    \begin{equation}
        \|\varphi_0(h^{-1}\rho)\Op (a)u\|\le C h^N (\|u\|_{H^{-N, -N}}+\|(\Lambda^\Omega_M-\om)u\|_{Z^{\infty, \infty, s_0}_{L_0-0, \infty, \infty}})
    \end{equation}
    uniformly in $h>0$ small. This implies $\|\Op (a)u\|_{H^{0, N}}$ is bounded by the same bounds on the right and $\WF u\cap \Sigma^+(\omega)=\emptyset$ as desired.

    2. We remind the reader that the neighborhood $V^+(\omega)$ of $\Sigma^+(\omega)$ was defined at the beginning of Section~\ref{s:5-2}. Consider a function $f_1\in S^{\operatorname{comp}, 1}$ such that
    \begin{equation}
        C_0^{-1}\rho^{-1}\le f_1\le C_0\rho^{-1}, \  H_p f_1\ge C_0^{-1}\rho\tilde\varphi \text{ on } V^+(\omega). 
    \end{equation}
    Let $f_2\in C^\infty_c(\mathbb{R})$ with
    \begin{equation}
        0\notin\supp(1-f_2), \ f_2'\le 0 \text{ on } [0,\infty), \ f_2'<0 \text{ on }[\frac12C_0^{-1},2C_0], \ \supp f_2'\cap [0,\infty)\subset[\frac14C_0^{-1},4C_0]. 
    \end{equation}
    Now construct the escape function
    \begin{equation}
        f=(1-f_2)(h f_1)\psi_\omega(\xi)\in S^{0,0}, \ \supp f\subset\{0\le \rho \le 4C_0^2h, |\xi|>\frac12\e_0\}.
    \end{equation}
    We have uniformly in $h>0$ small that 
    \begin{equation}
        H_p f=-h\psi_\omega(\xi)f_2'(h f_1)H_p f_1+(1-f_2)(hf_1)H_p\psi_\omega(\xi),
    \end{equation}
    the first term of which is supported inside $\{\frac14C_0^{-2}h\le \rho \le 4C_0^2h, |\xi|>\frac12\e_0\}$, disjoint from $\supp(1-\varphi_0(h^{-1}\rho))$, and the second term of which is supported away from $V^+(\omega)$. We have
    \begin{equation}\label{e:EscapeProp}
        H_p f>0\text{ on }V^+(\omega)\cap\{\rho\in [\frac12h, 2h]\}, \ H_p f\ge 0 \text{ on } V^+(\omega). 
    \end{equation}
    Fix $C_0=\sqrt{2}$. Let
    \begin{equation}
        F=\frac12 \big(\Op f+(\Op f)^*\big), \ F^*=F, \ \WF F\subset\{0\le \rho \le 8h, |\xi|>\frac12\e_0\}, \ F=\chi(h^{-1}\rho) F\chi(h^{-1}\rho),
    \end{equation}
    for $h>0$ small since the quantization procedure does not allow off-diagonal communication further than distance $\frac14$. Note $\rho\sim h$ on $\supp\varphi(h^{-1}\rho)$. There then exists $q\in S^{0,0}$ with $Q=\Op(q)$, 
    \begin{equation}
        \supp q\subseteq \{x\notin\supp\kappa, |\xi|>\e_0 \}, \ \supp q \cap V^+(\omega)=\emptyset, 
    \end{equation}
    such that for some $C>0$,
    \begin{equation}\label{e:FlowEst}
        \frac12 H_{p} f
        +\rho\varphi(h^{-1}\rho)^2|q(\rho, y,\xi)|^2
        \geq C^{-1}\rho \varphi_0(h^{-1}\rho)^2|a(\rho, y,\xi)|^2.
    \end{equation}
    Note that $\WF Q$ is away from $\WF u$ implies $Qu\in \mc{S}$ and $\varphi(h^{-1}\rho)Q u=O_{H^{N,N}}(h^\infty\|u\|_{H^{-N,-N}})$. Now we choose $b\in S^{0,0}$ with $B=\Op(b)$ such that
    \begin{equation}\label{e:WF-B_Disj}
        \supp H_{p}f\cup \supp q\cup\supp a\subset \{(x,\xi): b(x,\xi)=1\}, \ \supp b\cap \supp\kappa=\emptyset. 
    \end{equation}
    Note $[P, \Id-F]\in \Psi^{0,-1}$ even when $\Id-F$ is in general not in our calculus. Apply the sharp G\r{a}rding inequality with $\frac{1}{2i} [ P,\Id-F]+Q^*\varphi^2(h^{-1}\rho)\rho Q-C^{-1}A^*\varphi^2_0(h^{-1}\rho) \rho A$ on $\tilde\varphi(h^{-1}\rho)Bu$ to obtain
    \begin{multline}
        \|\rho^\frac12\varphi_0(h^{-1}\rho) A\tilde\varphi(h^{-1}\rho)Bu\|^2\leq
        C\|\rho^\frac12\varphi(h^{-1}\rho) Q\tilde\varphi(h^{-1}\rho)Bu\|^2+\frac{C}{2i}\langle [P ,\Id-F]\tilde\varphi(h^{-1}\rho)Bu, \tilde\varphi(h^{-1}\rho)B u\rangle \\ +C\| \tilde\varphi(h^{-1}\rho) Bu \|_{H^{-\frac12,-1}}^2
        \leq
        -\frac{C}{2i}\langle [P,F]u, u\rangle+Ch^2\|\tilde\varphi(h^{-1}\rho)Bu\|^2
         +O(h^\infty\|u\|_{H^{-N,-N}}^2)
    \end{multline}
    since $[P,F](\Id -\tilde\varphi(h^{-1}\rho)B)\in
    h^\infty\Psi^{-\infty,-\infty}$ uniformly . This implies
    \begin{equation}\label{e:3l16}
        \|\varphi_0(h^{-1}\rho) Au\|^2\le -\frac{C}{2i}\langle [P,F]u, u\rangle+Ch\|\tilde\varphi(h^{-1}\rho)Bu\|^2\\
         +O(h^\infty\|u\|_{H^{-N,-N}}^2)
    \end{equation}
    uniformly in $h>0$ small. 
    
    3. It remains to estimate $\langle [P,F] u, u\rangle$. Note $F=\chi(h^{-1}\rho)F\chi(h^{-1}\rho)$ and $(\Lambda^\Omega_M-\om)u\in Z^{\infty, \infty, s_0-1}_{L_0-0, \infty, \infty}$. This implies $\|F(\Lambda^\Omega_M-\om)u\|=O(h^\infty)$. On another hand, noting $\Im\langle (\Lambda^\Omega_M-\om)u, u\rangle\ge 0$, we have
\begin{multline}\label{e:5l2}
    O(h^\infty)=\Im \langle F(\Lambda^\Omega_M-\om)u ,u \rangle=\Im \langle (\Lambda^\Omega_M-\om)u ,u \rangle- \Im\langle (\Id-F)(\Lambda^\Omega_M-\om)u ,u \rangle\\
    \ge -\frac{1}{2i} \big(\langle [P,F]u ,u \rangle+\langle  Fu ,(\mc{E}')^*u \rangle-\langle \mc{E}u , F u \rangle\big)
    =-\frac{1}{2i} \langle [P,F]u ,u \rangle+O(h^\infty),
\end{multline}
where we used $(\Id-F)u\in H^{s_0, L_0-0}\subset L^2$ implies 
\begin{equation}
    -\Im\langle (\Id-F)(\Lambda^\Omega_M-\om)u ,u \rangle=-\frac1{2i}\langle [\Lambda^\Omega_M,F]u ,u \rangle
\end{equation}
and we used Proposition~\ref{p:dtN} for
\begin{multline}
    \Lambda^\Omega_M F= \Lambda^\Omega_M \chi(h^{-1}\rho) F\chi(h^{-1}\rho)=\tilde\chi(h^{-1}\rho) P \chi(h^{-1}\rho) F\chi(h^{-1}\rho)+ \mc{E}' \chi(h^{-1}\rho)F\chi(h^{-1}\rho)\\
    =P F+ \mc{E}' F+O(h^\infty\Psi^{-\infty, -\infty}),
\end{multline}
since $(1-\tilde\chi)(h^{-1}\rho) P \chi(h^{-1}\rho) F\chi(h^{-1}\rho)\in h^\infty\Psi^{-\infty,-\infty}$. The $F\Lambda^\Omega_M$ term is addressed similarly. Here the error $O(h^\infty)$ comes from 
   \begin{multline}\label{e:5l1}
   |\langle \mc{E} u ,Fu \rangle|
   \leq C\|Fu\|_{H^{-N,-N}} \|\chi(h^{-1}\rho)\tilde\psi_\omega(D)\chi(h^{-1}\rho)\mc{E} u\|_{H^{N,N}}+Ch^N\|u\|_{H^{-N,-N}}\|\chi(h^{-1}\rho)\mc{E} u\|_{H^{-N,-N}}\\
   \leq Ch^{N} \|u\|_{Z^{\infty, \infty, s_0}_{-L_0+0, \infty, -\infty}}^2,
   \end{multline}
   where $F=\chi(h^{-1}\rho)\tilde\psi_\omega(D)\chi(h^{-1}\rho) F+O(h^\infty\Psi^{-\infty,-\infty})$. For $\langle  Fu ,(\mc{E}')^*u \rangle$, note in  Proposition~\ref{p:dtN}, $(\mc{E}')^*=\mc{E}_0^*\tilde\chi(h^{-1}\rho)+\mc{E}_0'$ where $\langle Fu, \mc{E}_0' u\rangle$ can be estimated in a way similar to \eqref{e:5l1}. It remains to estimate $\langle  Fu ,\mc{E}_0^*\tilde\chi(h^{-1}\rho)u \rangle$. Proposition~\ref{p:dtNParametrix} implies $\Gamma^{N,0}_{N,\delta}\tilde\chi(h^{-1}\rho)\mc{E}_0\chi(h^{-1}\rho)\tilde\psi_\omega(D)\chi(h^{-1}\rho)\Gamma^{N,0}_{N,\delta}$ is bounded on $L^2(\partial\Omega)$ and so is $\Gamma^{N,0}_{N,\delta}\chi(h^{-1}\rho)\tilde\psi_\omega(D)\chi(h^{-1}\rho)\mc{E}_0^*\tilde\chi(h^{-1}\rho)\Gamma^{N,0}_{N,\delta}$. Thus we have
\begin{multline}
    \|\Gamma^{N,0}_{N,\delta}\chi(h^{-1}\rho)\tilde \psi_\omega(D)\chi(h^{-1}\rho)\mc{E}_0^*\tilde\chi(h^{-1}\rho) u\|\le C\|\psi_\omega(D)\tilde{\tilde{\chi}}(h^{-1}\rho)u\|_{H^{-N,-N}}\\
    +C\|(1-\psi_\omega(D))\tilde{\tilde{\chi}}(h^{-1}\rho)u\|_{L^2}\le C\|u\|_{Z^{0, 0, 0}_{0, -\infty, -\infty}}
\end{multline}
and thus
\begin{equation}
    |\langle F u, \mc{E}_0^*\tilde\chi(h^{-1}\rho) u\rangle|\le C\|\Gamma^{-N, 0}_{-N,\delta}Fu\|\|\Gamma^{N,0}_{N,\delta}\chi(h^{-1}\rho)\tilde \psi_\omega(D)\chi(h^{-1}\rho)\mc{E}_0^*\tilde\chi(h^{-1}\rho) u\|\le Ch^{N\delta}\|u\|_{Z^{0, 0, 0}_{0, -\infty, -\infty}}^2. 
\end{equation}
This justifies $|\langle  Fu ,(\mc{E}')^*u \rangle-\langle \mc{E}u , F u \rangle|=O(h^\infty)$ in \eqref{e:5l2}, and thus $\langle [P,F]u ,u \rangle=O(h^\infty)$. Bring it back in \eqref{e:3l16} to see \eqref{e:SmoothEst} as desired. 
\end{proof}

Now we show that there are finitely many eigenvalues for $\Lambda^\Omega_M$ in any given $\intw$.
\begin{lemma}[Eigenvalues are isolated] \label{l:FiniteEigenVal}
    Let $M_0>0$. Denote the point spectrum of $\Lambda^\Omega_M$ by $\operatorname{Spec}_{\mathrm{pp}}(\Lambda^\Omega_M)$. Then
    \begin{equation}
        |\{\sqrt{\omega^2+M^2}: \omega\in\intw\}\cap \operatorname{Spec}_{\mathrm{pp}}(\Lambda^\Omega_M)|<\infty
    \end{equation}
    for each $M\in [0, M_0]$. Without loss of generality, assume 
    \begin{equation}
        \{\sqrt{\omega^2+M^2}: \omega\in\partial\intw\}\cap \operatorname{Spec}_{\mathrm{pp}}(\Lambda^\Omega_M)=\emptyset
    \end{equation}
    for all $M\in[0, M_0]$. Then
    \begin{equation}
        \Big\{(\sqrt{\omega^2+M^2}, M):\sqrt{\omega^2+M^2}\in\operatorname{Spec}_{\mathrm{pp}}(\Lambda^\Omega_M), M\in [0, M_0], \omega\in \intw\Big\}\subset \mathbb{R}\times [0, M_0]
    \end{equation}
    is closed. We can define $\Pi_M$ as the projection operator
    \begin{equation}
    \Pi_M: Z^{-\infty, -\infty, 0}_{-L_0+0, -\infty, -\infty}\rightarrow \bigoplus_{\omega}\{\ker_{L^2}(\Lambda^\Omega_M-\sqrt{\omega^2+M^2}): \omega\in\intw\}\subset Z^{\infty,\infty, 1}_{L_0-0, \infty, \infty},
    \end{equation}
    and $\Pi_M^\perp=\Id-\Pi_M: Z^{s,l,s_0}_{l_0, l^-, l^+}\rightarrow Z^{s,l,s_0}_{l_0, l^-, l^+}$. 
\end{lemma}
\begin{remark}
    For $M_j\rightarrow M$, one can prove $\Pi_M^\perp \Pi_{M_j}\rightarrow 0$, but to show $\Pi_{M_j}^\perp \Pi_{M}\rightarrow0$ and thus $\Pi_{M_j}\rightarrow\Pi_{M}$ would require more. We will not pursue this direction in the paper.  
\end{remark}
\begin{proof}
    1. Let $\omega\in \intw$, $M\ge 0$. Assume against contradiction that there is a sequence of distinct $\omega_j\in\intw$ with $\omega_j\rightarrow \omega$, $\|u_j\|_{L^2}\equiv 1$, $(\Lambda^\Omega_M-\sqrt{\omega_j^2+M^2})u_j=0$. Then $\{u_j\}$ is an orthonormal sequence, and $u_j\rightharpoonup 0$ weakly in $L^2$, hence $u_j\rightarrow 0$ in the norm topology in $H^{0-0,0-0}$ via some subsequences. On another hand, Proposition~\ref{p:FredEst} implies
    \begin{equation}\label{e:5l11}
        \|u_j\|_{Z^{N,N, 1}_{L_0-0,N,-\frac12-0}}\leq C\|u_j\|_{H^{-N,-N}},\ \|u_j\|_{Z^{N,N, 1}_{L_0-0,-\frac12-0,N}}\leq C\|u_j\|_{H^{-N,-N}}.
    \end{equation}
    Note here we have both estimates since Lemma~\ref{l:EigenReg} implies $u_j\in Z^{\infty, \infty, 1}_{L_0-0, \infty, \infty}$. Then we have $\|u_j\|_{L^2}=o(1)$, yielding the desired contradiction for the finiteness of point spectrum. 

    2. Let $\sqrt{\omega^2+M^2}\notin\operatorname{Spec}_{\mathrm{pp}}(\Lambda^\Omega_M)$. Assume against contradiction that there is a sequence of distinct $M_j\rightarrow M$, $\omega_{j}\rightarrow \omega$ inside $\intw$, $\|u_j\|_{L^2}\equiv 1$, $(\Lambda^\Omega_{M_j}-\sqrt{\omega_j^2+M_j^2})u_j=0$. There then exists $u\in L^2$ that $u_j\rightharpoonup u$ weakly in $L^2$ and $u_j\rightarrow u$ in $H^{0-0,0-0}$ via some subsequences. Since $\Lambda^\Omega_{M}$ is continuous in $M$, we have
    \begin{equation}
        (\Lambda^\Omega_{M}-\sqrt{\omega^2+M^2})u_j=(\Lambda^\Omega_{M}-\Lambda^\Omega_{M_j})u_j-(\sqrt{\omega^2+M^2}-\sqrt{\omega_j^2+M_j^2})u_j=o_{Z^{\infty, \infty, 0}_{L_0-0, \infty, \infty}}(1). 
    \end{equation}
    This implies $(\Lambda^\Omega_{M}-\sqrt{\omega^2+M^2})u=0$ and thus $u=0$. Similar to \eqref{e:5l11} we have $\|u_j\|_{L^2}=o(1)$. The contradiction gives the stated closed graph property of $(\sqrt{\omega^2+M^2}, M)$.
\end{proof}

\begin{lemma}[Uniform absorption resolvent]\label{l:UAR}
Let $\Omega$ be an \acl domain. Let $M_0>0$. Let $s,l\in\mathbb{R}$, $s_0\in (0,1]$. Let $|l_0|< L_0$. For each fixed $M\in [0, M_0]$, we have the following estimates:
\begin{enumerate}[wide]
    \item For $l^->-\frac12$, $l^+<-\frac12$, we have
    \begin{equation}\label{e:5l12}
    \|u\|_{Z^{s,l, s_0}_{l_0,l^-,l^+}}\leq C\|(\Pi_M^\perp\Lambda^\Omega_M\Pi_M^\perp-\om-i\mu) u\|_{Z^{s-1,l, s_0-1}_{l_0,l^-+1,l^++1}}
    \end{equation}
    uniformly for all $u\in Z^{-\infty, -\infty, 0}_{-L_0+0, -\frac12+0, -\infty}(\partial\Omega)$, $\omega\in\intw$, $\mu\ge 0$ small.

    \item For $l^-<-\frac12$, $l^+>-\frac12$, we have
    \begin{equation}
    \|u\|_{Z^{s,l, s_0}_{l_0,l^-,l^+}}\leq C\|(\Pi_M^\perp\Lambda^\Omega_M\Pi_M^\perp-\om+i\mu) u\|_{Z^{s-1,l, s_0-1}_{l_0,l^-+1,l^++1}}
    \end{equation}
    uniformly for all $u\in Z^{-\infty, -\infty, 0}_{-L_0+0,-\infty,  -\frac12+0}(\partial\Omega)$, $\omega\in\intw$, $\mu\ge 0$ small.

    \item If additionally
    \begin{equation}
    \{\sqrt{\omega^2+M^2}: \omega\in\intw\}\cap \operatorname{Spec}_{\mathrm{pp}}(\Lambda^\Omega_M)=\emptyset, \text{ for all }M\in [0, M_0],
    \end{equation}
    then the above estimates are also uniform in $M\in [0, M_0]$. 
\end{enumerate}
\end{lemma}
\begin{proof}
Assume \eqref{e:5l12} does not hold and there are $\|u_j\|_{Z^{s,l, s_0}_{l_0,l^-,l^+}}\equiv 1$ with 
\begin{equation}
    (\Pi_M^\perp\Lambda^\Omega_M\Pi_M^\perp-\sqrt{\omega_j^2+M^2}-i\mu_j) u_j=o_{Z^{s-1,l, s_0-1}_{l_0,l^-+1,l^++1}}(1), \ \omega_j\rightarrow \omega\text{ within }\intw, \ \mu_j\rightarrow \mu
\end{equation}
This implies $\Pi_M u_j=o(1)$ in $Z^{\infty, \infty, 1}_{L_0-0, \infty, \infty}$ and thus
\begin{equation}
    (\Lambda^\Omega_M-\sqrt{\omega_j^2+M^2}-i\mu_j) \Pi_M^\perp u_j=o_{Z^{s-1,l, s_0-1}_{l_0,l^-+1,l^++1}}(1).
\end{equation}
Via a subsequence there is $u\in Z^{s,l, s_0}_{l_0,l^-,l^+}$ that $\Pi_M^\perp u_j\rightharpoonup u$ weakly in $Z^{s,l, s_0}_{l_0,l^-,l^+}$ and $\Pi_M^\perp u_j\rightarrow u$ in the norm of $Z^{s-0,l-0, s_0-0}_{l_0-0,l^--0,l^+-0}$. The limit $u$ solves $(\Lambda^\Omega_M-\sqrt{\omega^2+M^2}-i\mu) u=0$ and $\Pi_M u=0$, and thus $u=0$. Now the Fredholm estimates, Proposition~\ref{p:FredEst} implies $\Pi_M^\perp u_j=o(1)$ in $Z^{s,l, s_0}_{l_0,l^-,l^+}$. This contradicts the normalization of $u_j$, and yields the desired estimates at fixed $M$. For uniformity in $M$, note the extra assumption yields $\Pi_M=0$ and $\Pi_M^\perp=\Id$ and the Fredholm estimates are uniform in $M$ and the proof is similar. 
\end{proof}

We show the limiting absorption principle. 
\begin{proposition}[Limiting absorption principle] \label{p:LAP}
Let $\Omega$ be an \acl domain. Let $\omega\in \intw$, $M_0>0$. Let $M\in [0, M_0]$. Let $s,l\in\mathbb{R}$, $s_0\in (0,1)$. Let $|l_0|<L_0$, $l^->-\frac12$, 
$l^+<-\frac12$. Then the following are true for each $f\in Z^{s-1,l, s_0-1}_{l_0, l^-+1, l^++1}(\partial\Omega)$:
\begin{enumerate}[wide]
    \item The limit
\begin{equation}\label{e:LSALimit3}
    u=(\Pi_M^\perp\Lambda^\Omega_M\Pi_M^\perp-\om-i0)^{-1}f=\lim_{\mu\rightarrow 0^+}(\Pi_M^\perp\Lambda^\Omega_M\Pi_M^\perp-\om-i\mu)^{-1}f
\end{equation}
as a limit in $Z^{s-0,l-0, s_0-0}_{l_0-0, l^--0, l^+-0}$ exists. 
\item The limit $u\in Z^{s,l, s_0}_{l_0, l^-, l^+}$, and is the unique solution in $Z^{s,l, s_0}_{l_0, l^-, l^+}$ to
\begin{equation}\label{e:LSALimit4}
    (\Pi_M^\perp\Lambda^\Omega_M\Pi_M^\perp-\om)u=f.
\end{equation}
\item The limit $u$ depends continuously in $\omega\in\intw$ for each fixed $M\in[0, M_0]$.
\item If additionally,
\begin{equation}
    \{\sqrt{\omega^2+M^2}: \omega\in\intw\}\cap \operatorname{Spec}_{\mathrm{pp}}(\Lambda^\Omega_M)=\emptyset, \text{ for all }M\in [0, M_0],
\end{equation}
then the limit $u$ depends continuously in both $\omega\in\intw$ and $M\in [0,M_0]$. 
\item For any $k\in\mathbb{N}$, 
    \begin{equation}\label{e:LAP_Derivative}
        \partial^k_{\omega_M} (\Pi_M^\perp\Lambda^\Omega_M\Pi_M^\perp-\omega_M-i0)^{-1}f=k!
        (\Pi_M^\perp\Lambda^\Omega_M\Pi_M^\perp-\omega_M-i0)^{-k-1}f\in Z^{s,l,s_0}_{l_0,l^-,l^+-k},
    \end{equation}
for any $f\in Z^{s-1-k,l, s_0-1}_{l_0, l^-+1+k,l^++1}$ and $\omega_M=\sqrt{\omega^2+M^2}$. 
\end{enumerate}
\end{proposition} 
\begin{proof}
    1. We show $u\in Z^{s,l, s_0}_{l_0, l^-, l^+}$, the solution to $(\Pi_M^\perp\Lambda^\Omega_M\Pi_M^\perp-\om)u=f$, is unique in $Z^{s,l, s_0}_{l_0, l^-, l^+}$. Let $v\in Z^{-\infty, -\infty, 0}_{-L_0+0, -\frac12+0, -\infty}$, the difference between any two such solutions, solve $(\Pi_M^\perp\Lambda_M^\Omega\Pi_M^\perp-\om)v=0$. Lemma~\ref{l:UAR} implies $v=0$ and $u$ must be unique in $Z^{s,l, s_0}_{l_0, l^-, l^+}$. 

    2a. Let a sequence $\mu\rightarrow 0^+$ and $\omega_\mu\rightarrow \omega$. Lemma~\ref{l:absorbres} implies the $\mu$-sequence
    \begin{equation}
        u_\mu=(\Pi_M^\perp\Lambda_M^\Omega\Pi_M^\perp-\sqrt{\omega_\mu^2+M^2}-i\mu)^{-1}f=(\Lambda^\Omega_{M}-\sqrt{\omega_\mu^2+M^2}-i\mu)^{-1}\Pi_{M}^\perp f
        -(\sqrt{\omega_\mu^2+M^2}+i\mu)^{-1}\Pi_{M} f
    \end{equation}
    exists in $Z^{s, l, s_0}_{l_0, l^-+1, l^++1}$ and are uniformly bounded in $Z^{s, l, s_0}_{l_0, l^-, l^+}$ per Lemma~\ref{l:UAR}. There then exists $u\in Z^{s, l, s_0}_{l_0, l^-, l^+}$ that $u_\mu\rightharpoonup u$ weakly in $Z^{s, l, s_0}_{l_0, l^-, l^+}$ via a subsequence. 

    2b. We claim $(\Pi_M^\perp\Lambda^\Omega_M\Pi_M^\perp-\om)u=f$. Indeed, for any $\phi\in C^\infty_c$ we have $(\Pi_M^\perp\Lambda^\Omega_M\Pi_M^\perp-\om)\phi\in Z^{\infty,\infty, 0}_{L_0-0, \infty,\infty}$ and
    \begin{multline}
        \langle(\Pi_M^\perp\Lambda^\Omega_M\Pi_M^\perp-\omega_M)u, \phi\rangle=\langle u, (\Pi_M^\perp\Lambda^\Omega_M\Pi_M^\perp-\omega_M)\phi\rangle=
        \lim_{\mu\rightarrow 0^+}\langle u_{\mu}, (\Pi_M^\perp\Lambda^\Omega_M\Pi_M^\perp-\omega_M)\phi\rangle\\
        =\lim_{\mu\rightarrow 0^+}\langle (\Lambda^\Omega_M-\omega_M)u_{\mu}, \phi\rangle-\langle \Lambda^\Omega_M\Pi_M u_{\mu}, \phi\rangle=\langle f, \phi\rangle,
    \end{multline}
    as desired. Then $u$ is the unique limit of $u_\mu$ with respect to any sequence of $\mu\rightarrow 0$ via Step 1. This justifies \eqref{e:LSALimit3}, \eqref{e:LSALimit4} and the continuity in $\omega$. 

    2c. The continuity in $M\in[0,M_0]$ is immediate since the additional assumption gives uniform estimates in Lemma~\ref{l:UAR}. 

    3. We show \eqref{e:LAP_Derivative}. For $\mu>0$, the resolvent identity gives
    \begin{multline}
        (\Pi_M^\perp\Lambda^{\Omega}_M\Pi_M^\perp-\omega_M'-i\mu)^{-1}=(\omega_M'-\omega_M)(\Pi_M^\perp\Lambda^{\Omega}_M\Pi_M^\perp-\omega_M-i\mu)^{-1}(\Pi_M^\perp\Lambda^{\Omega}_M\Pi_M^\perp-\omega_M'-i\mu)^{-1}\\
        +(\Pi_M^\perp\Lambda^{\Omega}_M\Pi_M^\perp-\omega_M-i\mu)^{-1}.
    \end{multline}
    Taking $\omega_M'\rightarrow \omega_M$ gives $\partial_{\omega_M} (\Pi_M^\perp\Lambda^{\Omega}_M \Pi_M^\perp-\omega_M-i\mu)^{-1}=(\Pi_M^\perp \Lambda^{\Omega}_M \Pi_M^\perp-\omega_M-i\mu)^{-2}$. Take $\mu\rightarrow 0$ and iterate for \eqref{e:LAP_Derivative}. 
\end{proof}
\begin{remark}
    All constants in this section can depend on $\Omega$ in the following way: the upper bound on uniform Lipschitz characters of $\partial\Omega$, the upper bound on the choice of $r_0\ge 0$, and the lower bound on the radial Hardy constant $\rardyr_{\Omega}$, the upper bound on the diameter of $\partial\Omega\setminus\partial\tilde\Omega$ in Remark~\ref{r:aclcollar}, for each multi-index $\alpha$, the upper bounds on the derivatives of the metric $\gplus$ in \eqref{e:3l26}, and lower bounds on $|\det \gplus|$.
\end{remark}

\section{The scattering matrix}
\label{s:scatMat}
This section concerns the scattering matrix for the Dirichlet-to-Neumann map $\Lambda^{W_Y}_M$ on the boundary $\partial W_Y$, boundary of the \acl domain $W_Y$ with smooth section $Y$. We remind the readers that $L_0=\rardyr_{W_Y}+\frac12$. 

The construction of the scattering matrix relies heavily on the microlocal parametrix for the DtN map from Section~\ref{s:3}. Indeed, once we have this structure, the construction of the scattering matrix follows familiar lines (see e.g.~\cite{DZ:19, mel94}).  In lemma~\ref{l:PlaneWave}, we show that outgoing functions solving  $(\Lambda^{W_Y}_M-\sqrt{\omega^2+M^2})u=0$ near the radial sink have the form $a(y,\rho)e^{i\sqrt{\omega^2+M^2}/\rho}$ near $|x|=\infty$. To do this, we commute pseudodifferential operators defining $\Sigma_+$ through the DtN map, eventually showing that $u$ has Lagrangian regularity and $\Sigma^+$. We are then able to extract the asymptotic expansion by applying $\Lambda^{W_Y}_M-\sqrt{\omega^2+M^2}$ to $u$ and computing its symbol as a Lagrangian distribution. 

To construct the scattering matrix, given incoming data $a_-(y)$, we then construct an incoming function $u^-=\tilde{a}_-(y,\rho)e^{-i\sqrt{\omega^2+M^2}/\rho}$ supported near $\rho=0$ that solves $)(\Lambda_{M}^{W_Y}-\sqrt{\omega^2+M^2})u^-=0$ microlocally near the source. An application of the outgoing resolvent together with the asymptotic expansion from Lemma~\ref{l:PlaneWave} then yields the scattering matrix.

\subsection{Asymptotic expansions of outgoing solutions}

We show that we can decompose resonant states of $\Lambda^{W_Y}_M$ into plane waves.
\begin{lemma}[Plane wave expansion]\label{l:PlaneWave}
Let $W_Y\subset \R^d$ be an \acl domain with smooth cross-section $Y$. Let $M\geq 0$ and $\omega>0$. Let microlocalizer $B^+\in \Psi^{0,0}(\partial W_Y)$ be that $\WF(B^+)\subset \{|\xi|>c\omega\}$ for some $c>0$, $\Sigma^+(\omega)\subset\Ell(B^+)$. Let $\omega$-family $u\in Z^{-\infty, -\infty, 0}_{-L_0+0, -\frac12+0, -\infty}(\partial W_Y; \omega)$ with
\begin{equation}
    (\Lambda^{W_Y}_M-\om)u\in Z^{\infty, \infty, -1}_{-L_0+0, \infty, \infty}(\partial W_Y;\omega), \ B^+ (\Lambda^{W_Y}_M-\om)u=O(\omega^{-\infty}\mc{S}),
\end{equation}
the second condition is uniform in $\omega>0$ large. Then there exist $\rho_0>0$, constant $\alpha(\omega)\ge 1$, and $\tilde a^+(\rho,y;\omega) \in S^{0}(\partial W_Y)$ that, for each $A\in \Psi^{0,0}(\partial W_Y)$ with $\WF A\subset\{\rho<\rho_0, |\xi|> c\omega\}$ for some $c>0$, we have
\begin{equation}\label{e:PlaneWaveExpSink}
        (A u)(\rho,y)=\alpha(\omega)\big(A\rho^{\frac{d}{2}-1} \tilde a^+(\rho,y;\omega) e^{i \omega/\rho}+O(\omega^{-\infty}\mathcal{S})\big),
\end{equation}
in local coordinates $y$ on $\partial Y$, uniformly in $\omega>0$ large. Given $\alpha(\omega)$, such $\tilde a^+\in S^0(\partial W_Y)$ is unique in $\{\rho<\rho_0\}$ modulo $\omega^{-\infty}\mc{S}(\partial W_Y)$. Moreover, there exist $a^{+}(y;\omega)\in S^0(\partial Y)$, $a^+_{k}(\rho, y; \omega)\in S^{0}(\partial W_Y)$ for $k\ge 1$ such that
\begin{equation}
    \tilde a^{+}(\rho, y;\omega)=a^{+}(y;\omega)+\sum_{k=1}^\infty \omega^{-k}\rho^ka^+_{k}(\rho,y;\omega)+O(\omega^{-\infty}\mc{S}).
\end{equation}
\end{lemma}
\begin{proof}
We mostly follow the strategy of \cite[Lemma~4.1]{DZ:19} and \cite{mel94} with a small modification in Step 2, where we commute $\rho \Delta_y$ with the operator, rather than $\Delta_y$ in coordinates where $\Sigma^+$ is given by $\{(0,y,-\omega,0)\}$. 
    
    1. We set up some notations. Since $A\in\Psi^{0,0}(\partial W_Y)$, there is $\chi\in S^0(\partial W_Y)$ where $\partial W_Y$ is smooth over $\supp\chi$ and $A=A\chi$. Let $\tilde\chi, \tilde{\tilde{\chi}}\in S^0(\partial W_Y)$ such that $\chi\prec\tilde\chi\prec\tilde{\tilde{\chi}}$, and $\partial W_Y$ is smooth over $\supp\tilde{\tilde{\chi}}$. Let $\sigma(B^+)=b^+$ and we assume without loss of generality $\WF (\Id-B^+)\cap \Sigma^+(\omega)=\emptyset$. Throughout the proof, we ask any microsupport separation to be of distance at least $c\omega$ in $\xi$ for some $c>0$. Consider microlocalizer $B\in \Psi^{0,0}$ near $\Sigma^+(\omega)$,
    \begin{equation}
        \WF(\Id-B)\cap \Sigma^+(\omega)=\emptyset, \ \WF(\Id-B^+)\cap \WF B=\emptyset.
    \end{equation}
     For convenience, let the microlocal representation of $\Lambda_M^{W_Y}-\om$ near the sink $\Sigma^+(\omega)$ be 
    \begin{equation}
        P_\omega^+=B^+\chi(P-\om)\tilde\chi\in \Psi^{1,0}, \ p_\omega^+=\sigma(P_\omega^+)=b^+\chi(\sqrt{|\xi|_g^2+M^2}-\om).
    \end{equation}
    We know by assumption $\WF^{-\infty, -\frac12+0} u\cap \Sigma^-(\omega)=\emptyset$. The Fredholm estimate, Proposition~\ref{p:FredEst}, then implies $u\in Z^{\infty, \infty, 0}_{-L_0+0, \infty, -\frac12-0}$, and thus $B \tilde{\tilde{\chi}}u\in H^{\infty, -\frac12-0}$, for each $\omega$ not uniformly. We normalize them by $\alpha(\omega)=\max(\|B \tilde{\tilde{\chi}}u\|_{H^{0, -\frac12-0}}, 1)$, and denote $u^+=\alpha(\omega)^{-1}B \tilde{\tilde{\chi}}u\in H^{\infty, -\frac12-0}$. Note Proposition~\ref{p:dtN} implies
    \begin{equation}
        P_\omega^+ u^+=\alpha(\omega)^{-1}\big(B^+\chi(\Lambda^{W_Y}_M-\om)u-B^+\chi\mc{E}_Mu
        -B^+\chi (P-\om)\tilde\chi (\Id-B)\tilde{\tilde{\chi}} u\big) =O(\omega^{-\infty}\mc{S}), 
    \end{equation}
    where we used that $\WF(B^+\chi (P-\om)\tilde\chi (\Id-B)\tilde{\tilde{\chi}})$ is disjoint from $\Sigma^+(\omega)$ and $\{\xi=0\}$. From this point onwards, all estimates in this proof are uniform in $\omega>0$ large. 

    2. We show that $u$ near $\Sigma^+(\omega)$ is a Lagrangian distribution associated with the sink $\Sigma^+(\omega)=\{(0,y,-\omega,0)\}$, that is to claim for all $k\ge 0$ and all $L_1$, \dots, $L_k\in \Psi^{1,1}$ such that $\omega^{-1}L_j$ is uniformly bounded in $\omega$, $\sigma(L_j)|_{\Sigma^+(\omega)}=0$ for $1\le j\le k$, we have
    \begin{equation}\label{e:LagrReg}
        L_{1}\cdots L_{k} u^+ \in H^{\infty, -\frac12-0}.
    \end{equation}
    Consider in local coordinates near $(0,y)$, the testing operators in $\Psi^{1,1}$ for $\Sigma^+(\omega)$, 
    \begin{equation}
        V=(V_j)_{j=0}^{d-2}, \ V_0=\rho^{-1}P_\omega^+, \ V_j=D_{y_j},\ 1\le j\le d-2.
    \end{equation}
    They correspond to $S^{1,1}$-defining functions of $\Sigma^+(\omega)$,
    \begin{equation}
        v_0=\sigma(V_0)=\rho^{-1}(\sqrt{\tau^2+|\eta|_{g_{\infty}}^2+M^2}-\om), \ v_j=\sigma(V_j)=\rho^{-1}\eta_j.
    \end{equation}
    In Step 2, we iterate in $(k,n)\in \mathbb{N}_0\times\mathbb{N}_0$ to claim for any $k'\le k,n'\le n$, any word $\beta$ of length $k'$, 
    \begin{equation}
        \beta=(\beta_1,\dots, \beta_{k'}), \ 0\le \beta_l\le d-2 \text{ for each } 1\le l\le k',
    \end{equation}
    we have, with the notation $V_\beta=V_{\beta_1}\cdots V_{\beta_{k'}}$, 
    \begin{equation}\label{e:6l7}
        (\rho\Delta_y)^{n'} V_\beta u^+\in H^{\infty, -\frac12-0}.
    \end{equation} 
    Once we have shown \eqref{e:6l7}, \eqref{e:LagrReg} follows immediately. Indeed, noting $L_j$ has principal symbols vanishing on $\Sigma^+(\omega)$, we have
    \begin{equation}
        L_j=\sum_{l=0}^{d-2}O(\Psi^{0,0})V_l+O(\Psi^{0,0}),
    \end{equation}
    and thus
    \begin{equation}
        L_1\cdots L_k u^+=\sum_{|\beta|\le k}O(\Psi^{0,0})V_\beta u^+\in H^{\infty, -\frac12-0},
    \end{equation}
    by the case $(k,0)$ of \eqref{e:6l7}. The rest of Step 2 is devoted to showing \eqref{e:6l7}. 

    2a. We provide a priori estimates on commutators. We check near $\Sigma^+(\omega)$,
    \begin{gather}
        H_{p} v_0=if_0 v_0,\ f_0=\frac{i\rho\tau}{\sqrt{\tau^2+|\eta|_{g_{\infty}}^2+M^2}},\\
        H_{p} v_j=-\sum_{l,m=1}^{d-2} \frac{\eta_l(\partial_{y_j}(g_{\infty})^{lm})\eta_m}{2\sqrt{\tau^2+|\eta|_{g_{\infty}}^2+M^2}},\ 1\le j\le d-2,
    \end{gather}
    and thus we have for $1\le j\le d-2$, modulo $O(\Psi^{1,0})$ supported away from $\Sigma^+(\omega)$, 
    \begin{gather}
    [P_\omega^+,V_0]=F_0 V_0+O(\Psi^{0,-1}), \ \sigma(F_0)=f_0\in S^{0,-1},\\
    \label{e:6l10}
        [P_\omega^+, V_j]=\sum_{l}O(\Psi^{0,-1})V_l+O(\Psi^{0,-1})=\sum_{l,m=1}^{d-2}O(\Psi^{-1,0})\rho^2V_lV_m+O(\Psi^{0,-1}).
    \end{gather}
    Noting $V_0=\rho^{-1}P_\omega^+$, we have for each $1\le j,m\le d-2$, 
    \begin{equation}
        [V_0, V_j]=O(\Psi^{0,0})V_0+\sum_{l=1}^{d-2}O(\Psi^{0,0})V_l+O(\Psi^{0,0}), \ [V_j,V_m]=\sum_{l=1}^{d-2}O(\Psi^{0,0})V_l
    \end{equation}
    modulo $O(\Psi^{1,1})$ supported away from $\Sigma^+(\omega)$, and implies with the same modulo class,
    \begin{equation}\label{e:6l17}
        \ad_V^{\gamma}V_0=O(\Psi^{0,0})V_0+\sum_{l=1}^{d-2}O(\Psi^{0,0})V_l+O(\Psi^{0,0}).
    \end{equation}
    Moreover, if $\gamma_j\neq 0$ for all $j$, then modulo $O(\Psi^{1,0})$ supported away from $\Sigma^+(\omega)$, 
    \begin{equation}
        \ad_{V}^\gamma P_\omega^+=\sum_{l,m=1}^{d-2}O(\Psi^{-1,0})\rho^2V_lV_m+O(\Psi^{0,-1}).
    \end{equation}
    We have near $\Sigma^{+}(\omega)$, 
    \begin{equation}
        H_{p} (\rho^{-1}|\eta|_{g_\infty}^2)=\frac{\tau|\eta|_{g_\infty}^2+\tau\rho \eta^t\cdot \partial_\rho g_\infty\cdot \eta }{\sqrt{\tau^2+|\eta|_{g_\infty}^2+M^2}}=\frac{\tau|\eta|_{g_\infty}^2+O(\rho\tau|\eta|_{g_\infty}^2)}{\sqrt{\tau^2+|\eta|_{g_\infty}^2+M^2}}.
    \end{equation}
    This implies, modulo $O(\Psi^{2,0})$ supported away from $\Sigma^+(\omega)$, 
    \begin{equation}
        [\rho \Delta_y, P_\omega^+]=F_1\rho \Delta_y+O(\Psi^{1,-1}), \ F_1\in \Psi^{0,-1}, \ \sigma(F_1)=\frac{i\rho\tau}{\sqrt{\tau^2+|\eta|_{g_\infty}^2+M^2}}.
    \end{equation}
    Then for all $n\in\mathbb{N}$, we have modulo $O(\Psi^{n+1,0})$ supported away from $\Sigma^+(\omega)$,
    \begin{equation}
        \ad_{\rho\Delta_y}^n P_\omega^+=O(\Psi^{n-1,-1})\rho\Delta_y+O(\Psi^{n,-1}). 
    \end{equation}
    
    2b. We bootstrap by showing the case $(0,n)$ holds for every $n$. Step 1 implies the case $(0,0)$ holds. Assume the case $(0,n-1)$ holds and consider
    \begin{equation}
        P_\omega^+(\rho\Delta_y)^n u^+=(\rho\Delta_y)^n P_\omega^+ u^+-\sum_{n'<n}\binom{n}{n'}(\ad_{\rho\Delta_y}^{n-n'}P_\omega^+)(\rho\Delta_y)^{n'}u^+
    \end{equation}
    For $F_1\in \Psi^{0,-1}$, we have
    \begin{multline}
        \sum_{n'<n}\binom{n}{n'}(\ad_{\rho\Delta_y}^{n-n'}P_\omega^+)(\rho\Delta_y)^{n'} u^+=\sum_{n'=0}^{n-1}O(\Psi^{n-n',-1})(\rho\Delta_y)^{n'} u^++nF_1(\rho\Delta_y)^n u^+\\
        =nF_1(\rho\Delta_y)^n u^++O(H^{\infty, \frac12-0}),
    \end{multline}
    by the case $(0,n-1)$. This gives
    \begin{equation}
        (P_\omega^++n F_1)(\rho\Delta_y)^n u^+\in H^{\infty, \frac12-0}. 
    \end{equation}
    Applying the standard radial sink estimate to $(P_\omega^++n F_1)$ gives for any $s, N>0$, 
\begin{equation}
    \|(\rho\Delta_y)^n u^+\|_{H^{s,-\frac12-0}}\le C\|(P_\omega^++n F_1)(\rho\Delta_y)^nu^+\|_{H^{s+1, \frac12-0}}+C\|Q_+(\rho\Delta_y)^n u^+\|_{H^{s, -\frac12-0}}+C\|u^+\|_{H^{-N, -N}}
\end{equation}
for some $Q_+\in \Psi^{0,0}$ supported near $\Sigma^+(\omega)$ with $\WF Q_+\cap \Sigma^+(\omega)=\emptyset$. One can choose $Q_+$ that $\WF Q_+\cap\WF B=\emptyset$ and this yields $(\rho\Delta_y)^n u^+\in H^{\infty, -\frac12-0}$ as desired, and validates the case $(0,n)$ iteratively.

    2c. We claim if $V_{\beta_j}=V_0$ for some $1\le j\le k$, then \eqref{e:6l7} holds for case $(k,n)$ for any $n$. Indeed if so, then 
    \begin{equation}
        (\rho\Delta_y)^n V_\beta u^+=(\rho\Delta_y)^nV_{\beta_1}\cdots V_{\beta_{j-1}}\big(V_{\beta'} \rho^{-1}P_\omega^+u^+-\sum_{\gamma\subsetneq\beta'}\binom{\beta'}{\gamma}(\ad_V^{\beta'\setminus\gamma} V_0)V_\gamma u^+\big)
    \end{equation}
    is in $H^{\infty, -\frac12-0}$ as desired, for $\beta'=(\beta_{j+1}, \dots, \beta_{k})$, since
    \begin{equation}
        (\rho\Delta_y)^nV_{\beta_1}\cdots V_{\beta_{j-1}}(\ad_V^{\beta'\setminus\gamma} V_0)V_\gamma u^+=\sum_{n'\le n}\sum_{|\beta''|<k} O(\Psi^{n-n',0})(\rho\Delta_y)^{n'} V_{\beta''} u^+\in H^{\infty, -\frac12-0},
    \end{equation}
    per iterative assumption \eqref{e:6l7} of the case $(k-1,n)$. Here we used $\ad_{V}^\gamma O(\Psi^{0,0})=O(\Psi^{0,0})$ and $\ad_{\rho\Delta_y}^{n'} O(\Psi^{0,0})=O(\Psi^{n',0})$. 

    2d. It then suffices to show \eqref{e:6l7} for case $(k,n)$ under the extra condition that $\beta_j\neq 0$ for all $1\le j\le k$. We proceed by induction on $n$. We first show the case $(k,0)$ holds. Consider
    \begin{equation}
        P_\omega^+ V_\beta u^+=V_\beta P_\omega^+ u^+-\sum_{\gamma\subsetneq \beta}\big(\ad_V^{\beta\setminus\gamma}P_\omega^+\big)V_\gamma u^+=O(H^{\infty,\frac12-0}),
    \end{equation}
    where we used 
    \begin{equation}
        \sum_{\gamma\subsetneq \beta}\big(\ad_V^{\beta\setminus\gamma}P_\omega^+\big)V_\gamma u^+=\sum_{\gamma\subsetneq \beta}\sum_{j,l=1}^{d-2}O(\Psi^{-1,0})\rho^2 V_jV_l V_\gamma u^++O(\Psi^{0,-1}V_\gamma u^+)=O(H^{\infty,\frac12-0}),
    \end{equation}
    since $\Delta_y$ is elliptic over $V_jV_l$ and $\rho\Delta_yV_\gamma u^+\in H^{\infty,-\frac12-0}$ for $|\gamma|\le k-1$ per the case $(k-1,1)$. This implies the case $(k,0)$. 
    
    2e. Now assume the case $(k,n-1)$ holds. Consider
    \begin{equation}
        P_\omega^+(\rho\Delta_y)^n V_\beta u^+=(\rho\Delta_y)^n V_\beta P_\omega^+ u^+-\sum_{n'<n}\binom{n}{n'}(\ad_{\rho\Delta_y}^{n-n'}P_\omega^+)(\rho\Delta_y)^{n'}V_\beta u^+-\sum_{\gamma\subsetneq \beta}(\rho\Delta_y)^{n} \big(\ad_V^{\beta\setminus\gamma}P_\omega^+\big)V_\gamma u^+.
    \end{equation}
    For $F_1\in \Psi^{0,-1}$, we have
    \begin{multline}
        \sum_{n'<n}\binom{n}{n'}(\ad_{\rho\Delta_y}^{n-n'}P_\omega^+)(\rho\Delta_y)^{n'}V_\beta u^+= \sum_{n'=0}^{n-1}O(\Psi^{n-n',-1})(\rho\Delta_y)^{n'}V_\beta u^++n F_1(\rho\Delta_y)^nV_\beta u^+\\=n F_1(\rho\Delta_y)^nV_\beta u^++O(H^{\infty, \frac12-0}),
    \end{multline}
    by the case $(k,n-1)$. On another hand, 
    \begin{multline}
        \sum_{\gamma\subsetneq \beta}(\rho\Delta_y)^n \big(\ad_V^{\beta\setminus\gamma}P_\omega^+\big)V_\gamma u^+=\sum_{\gamma\subsetneq \beta}\sum_{n'\le n}\sum_{j,l=1}^{d-2}O(\Psi^{n-n'-1,0})(\rho\Delta_y)^{n'} \rho^2 V_jV_l V_\gamma u^+\\+O(\Psi^{n-n',-1})(\rho\Delta_y)^{n'}V_\gamma u^+
        =O(H^{\infty,\frac12-0}),
    \end{multline}
    where we used $(\Delta_y)^{n'+1}$ being elliptic over $(\rho\Delta_y)^{n'}V_jV_l$ and $(\rho\Delta_y)^{n'+1}V_\gamma u^+\in H^{\infty,-\frac12-0}$ for $|\gamma|\le k-1$ and $n'\le n$ per the case $(k-1,n+1)$. We ended up with 
    \begin{equation}
        (P_\omega^++n F_1)(\rho\Delta_y)^n V_\beta u^+\in H^{\infty, \frac12-0}. 
    \end{equation}
    Similar to Step 2b, applying the standard radial sink estimate to $(P_\omega^++n F_1)$ yields $(\rho\Delta_y)^n V_\beta u^+\in H^{\infty, -\frac12-0}$ as desired. By induction we have shown \eqref{e:6l7} for all $(k,n)$. 

3. Denote $h=\omega^{-1}$. Now we invoke \cite[Theorem 25.2.4]{HorV4} to see the asymptotic structure at the sink: \eqref{e:LagrReg} implies there is $\tilde a_+(\rho,y) \in S^{0}$ such that
    \begin{equation}\label{e:6l18}
        u^+=\alpha(\omega)^{-1}B\tilde{\tilde{\chi}} u=\rho^{\frac{d}{2}-1} \tilde a^+(\rho,y) e^{i /h\rho}+O(h^\infty\mathcal{S}). 
    \end{equation}
    It remains to show that there exist $a^{+}\in S^0(\partial Y)$, $a^+_{k}\in S^{0,0}(\partial W_Y)$ that 
    \begin{equation}\label{e:6l5}
        \tilde a^+(\rho, y)=a^{+}(y)+\sum_{k =1}^\infty (h\rho)^k a^+_{k}(\rho, y)+O(h^\infty\mathcal{S}).
    \end{equation}
    Note \eqref{e:6l18} implies
    \begin{equation}\label{e:6l12}
        O(h^\infty\mc{S})=P_\omega^+ u^+=-\frac{i\omega }{\sqrt{\omega^2+M^2}}\rho^{\frac d2}\Big(-\rho\partial_\rho \tilde a^+(\rho,y)+\sum_{k=2}^\infty (h\rho)^{k-1} L^+_{k} \tilde a^+(\rho,y)\Big)e^{i /h\rho}+O(h^{\infty}\mc{S}),
    \end{equation}
    for $L^+_{k}$ being b-differential operators locally involving $\rho\partial_\rho$ and $\partial_y$ up to order $2k-1$. Note $L^+_{k}\tilde a_+\in S^{0}$. We then have the governing transport equation:
    \begin{equation}\label{e:6l3}
        \rho\partial_\rho \tilde a^+(\rho,y)-\sum_{k=0}^\infty (h\rho)^{k+1} L^+_{k+2} \tilde a^+(\rho,y)=O(h^{\infty}\mc{S}). 
    \end{equation}

4a. We proceed by induction to show \eqref{e:6l5}. It is equivalent to showing for any $K\ge 1$, there exist $a^+=a_0^+, a^+_{1},\dots, a^+_{K-1}\in S^{0}$ and $\tilde a^+_{K}$ that
\begin{equation}\label{e:6l2}
    \tilde a^+(\rho,y)=a^{+}(y)+\sum_{k=1}^{K-1} (h\rho)^k a^+_{k}(\rho,y)+(h\rho)^K\tilde a^+_{ K}, \ \rho\partial_\rho ((h\rho)^K \tilde a^+_{K})=O((h\rho)^K).
\end{equation}
We address the case $K=1$ first. From \eqref{e:6l3} we have
\begin{equation}\label{e:6l6}
    \rho\partial_\rho \tilde a^+(\rho, y)=h\rho L^+_{2} \tilde a^+(\rho,y)+O(h^2\rho^2). 
\end{equation}
We can solve the transport equation $\rho\partial_\rho a^{+}(y)=0$ for $a^+_{0}(y)\in S^{0}$ by $a^{+}(y)=\tilde a^+(0, y)$. Plugging
\begin{equation}
    \tilde a^{+}_{1}(\rho, y)=(h\rho)^{-1}(\tilde a^+(\rho, y)-a^{+}(y))
\end{equation}
into \eqref{e:6l6} gives $\rho\partial_\rho (h\rho\tilde a_{+,1})=O(h\rho)$. This validates the case $K=1$. 

4b. We show if \eqref{e:6l2} holds for some $K\ge 1$, then there exist $a^+_{K}\in S^{0,0}$ and $\tilde a^+_{K+1}$ that 
\begin{equation}
    \tilde a^+(\rho,y)=a^{+}(y)+\sum_{k=1}^{K} (h\rho)^k a^+_{k}(\rho,y)+(h\rho)^{K+1}\tilde a^+_{K+1}, \ \rho\partial_\rho ((h\rho)^{K+1} \tilde a^+_{K+1})=O((h\rho)^{K+1}).
\end{equation}
Indeed, substituting \eqref{e:6l2} into the governing equation \eqref{e:6l3} gives
\begin{equation}\label{e:6l4}
    \rho\partial_\rho ((h\rho)^K \tilde a^+_{K})=\sum_{k=0}^{K-1}(h\rho)^{K-k} L^+_{K+1-k} ((h\rho)^k a^+_{k})+O((h\rho)^{K+1}),
\end{equation}
where we used $\rho\partial_\rho ((h\rho)^K \tilde a^+_{K})$ is of size $(h\rho)^K$. Solving the transport equation
\begin{equation}\label{e:6l19}
    \rho\partial_\rho ((h\rho)^K a^+_{K})=\sum_{k=0}^{K-1}(h\rho)^{K-k} L^+_{K+1-k} ((h\rho)^k a^+_{k})
\end{equation}
gives the desired $a^+_{K}\in S^{0}$. Now letting
\begin{equation}
    \tilde a^+_{K+1}=(h\rho)^{-1}(\tilde a^+_{K}-a^+_{K})
\end{equation}
gives $\rho\partial_\rho ((h\rho)^{K+1} \tilde a^+_{K+1})=O((h\rho)^{K+1})$ in \eqref{e:6l4}. One can then iterate and this shows \eqref{e:6l5} at the sink. To observe \eqref{e:PlaneWaveExpSink}, it suffices to note for any $A=A\tilde{\tilde{\chi}}$ described in the statement of the lemma, one has that $\alpha(\omega)^{-1}Au=\alpha(\omega)^{-1}AB\tilde{\tilde{\chi}}u+O(h^\infty\mc{S})$ by noting $\alpha(\omega)^{-1}A(\Id-B)\tilde{\tilde{\chi}}u=O(h^\infty\mc{S})$. 

5. To show uniqueness, it suffices to show for $u^+\in H^{-\frac12-0}$ with $P_\omega^+u^+=O(h^\infty\mathcal{S})$,
\begin{equation}
    u^+=\Big(\rho^{\frac{d}{2}-1} \tilde a^+(\rho,y) e^{i /h\rho}\Big)+O(h^\infty\mathcal{S}), \ a_+(y)=0
\end{equation}
implies $u^+=O(h^\infty \mc{S})$. Indeed, given $a^+_0, \dots, a^+_{K-1}=0$, \eqref{e:6l19} implies $(h\rho)^Ka^+_K$ does not depend on $\rho$ while $a^+_K\in S^0$, hence $a^+_K=0$. Induction concludes the uniqueness claim.  
\end{proof}

\subsection{Unitarity and existence of the scattering matrix}
We prove a pairing formula. 
\begin{lemma}[Green's second identity for $\Lambda^{W_Y}_M$]\label{l:green}
Let $M\ge 0$, $\omega>0$ fixed. Let $\chi\in S^0(\partial W_Y)$ with $\partial W_Y$ smooth over $\supp\chi$, $\supp(1-\chi)$ compact. Let $u,v\in Z^{\frac12, 0, \frac12}_{0, -\frac12-0, -\frac12-0}(\partial W_Y)$ where 
    \begin{gather}
        u=\chi\rho^{\frac{d-2}{2}}\Big(u^+ e^{i\omega/\rho}+u^- e^{-i\omega/\rho}\Big)+O(Z^{\frac12, 0, \frac12}_{0,-\frac12+0, -\frac12+0}),\\
        v=\chi\rho^{\frac{d-2}{2}}\Big(v^+ e^{i\omega/\rho}+v^- e^{-i\omega/\rho}\Big)+O(Z^{\frac12, 0, \frac12}_{0,-\frac12+0, -\frac12+0}),
    \end{gather}
    for $u^\pm, v^\pm\in L^2(\partial Y)$. Assume 
    \begin{equation}
        (\Lambda^{W_Y}_M-\om)u, (\Lambda^{W_Y}_M-\om)v\in Z^{-\frac12,0,-\frac12}_{0, \frac12+0, \frac12+0}(\partial W_Y)
    \end{equation}
    Then we have
    \begin{multline}\label{e:6l13}
        \langle (\Lambda^{W_Y}_M-\om)u,v\rangle_{\partial W_Y}-\langle u,(\Lambda^{W_Y}_M-\om)v\rangle_{\partial W_Y}\\
        =\frac{\omega}{i\om}\Big(\langle u^+,v^+\rangle_{\partial Y}-\langle u^-,v^-\rangle_{\partial Y}\Big).
    \end{multline}
\end{lemma}
\begin{proof}
1. We first turn the left of \eqref{e:6l13} into its microlocal representation. Let $\tilde\chi\in S^0(\partial W_Y)$ that $\supp\chi\cap\supp(1-\tilde\chi)=\emptyset$ with $\partial W_Y$ smooth over $\supp\tilde\chi$. Denote the microlocal representation of $\Lambda_M^{W_Y}-\om$ by 
    \begin{equation}
        P_\omega=\chi(P-\om)\tilde\chi\in \Psi^{1,0}, \ p_\omega=\sigma(P_\omega)=\chi(\sqrt{|\xi|_g^2+M^2}-\om).
    \end{equation}

1a. Let $\varphi(\rho) \in S^0(\partial W_Y)$ that only depends on $\rho$, with $\supp\varphi\cap\supp(1-\chi)=\emptyset$, $\varphi(\rho)\equiv 1$ near $\rho=0$ and $\supp\varphi\subset\{\rho<1\}$. We remind the reader that $\psi_\omega \in C^\infty_c(\mathbb{R}_\xi)$ defined in Definition~\ref{d:Zspace} satisfies $\omega\notin \supp(1-\psi_\omega)$ and $0\notin \supp\psi_\omega$. Let $\varphi_0(\rho), \varphi_{00},\tilde\varphi \in S^0(\partial W_Y)$ satisfy the same requirement for $\varphi$, with $\varphi_{00}\prec \varphi_0\prec \varphi\prec\tilde\varphi$. Let $0<h<\frac12$. Let $Q\in \Psi^{0,0}$ be that
\begin{equation}
Q=\Op(\varphi_{00}(h^{-1}\rho)\psi_\omega(|\xi|_g)),\ \sigma(Q)=q(\rho, y;\tau, \eta)=\varphi_{00}(h^{-1}\rho)\psi_\omega(|\xi|_g).
\end{equation}
Note $Q=\varphi_0(h^{-1}\rho)Q\varphi_0(h^{-1}\rho)$ since the quantization procedure does not move the support more than $\frac14$. We have 
\begin{equation}
    Q=Q^*+O_{\Psi^{-1,0}}(h), \ \WF(\Id-Q)\cap \Sigma(\omega)=\emptyset, \ \WF Q\cap \{(\tau,\eta)=0\}=\emptyset.
\end{equation}
Noting $\langle (\Lambda_M^{W_Y}-\om) u, Q v\rangle=o_h(1)$ and $(\Id-Q^*)v\in L^2$, we have
\begin{equation}
    \langle (\Lambda_M^{W_Y}-\om) u, v\rangle=\langle(\Id-Q)(\Lambda_M^{W_Y}-\om) u, v\rangle+o_h(1),
\end{equation}
since $\Lambda_M^{W_Y}$ is self-adjoint on $L^2$. Similar arguments around $\langle Q u, (\Lambda_M^{W_Y}-\om) v\rangle$ give
\begin{equation}
    \langle(\Lambda_M^{W_Y}-\om) u, v\rangle-\langle u, (\Lambda_M^{W_Y}-\om)v\rangle=\langle[\Lambda_M^{W_Y}, Q]u, v\rangle+o_h(1). 
\end{equation}

1b. We further reduce the $[\Lambda^{W_Y}_M,Q]$ into $[P_{\omega}^{\varphi},Q]$ for $P_{\omega}^{\varphi}=\varphi(h^{-1}\rho) P_\omega\tilde\varphi(h^{-1}\rho)$. By Proposition~\ref{p:dtN} we have
\begin{equation}
    \varphi(h^{-1}\rho)\Lambda^{W_Y}_M=\varphi(h^{-1}\rho)P\tilde\varphi(h^{-1}\rho)+\varphi(h^{-1}\rho)\mc{E}_M,\ \Lambda^{W_Y}_M\varphi(h^{-1}\rho)=\tilde\varphi(h^{-1}\rho) P\varphi(h^{-1}\rho)+\mc{E}_M'\varphi(h^{-1}\rho). 
\end{equation}
With $Q\mc{E}_M:H^{0, -L_0+0}\rightarrow H^{N,N}$ by \eqref{e:errorEst2}, this implies
\begin{equation}
    Q\Lambda_M^{W_Y}u=Q\varphi(h^{-1}\rho)P \tilde\varphi(h^{-1}\rho)u+O_{H^{N,N}}(h^\infty). 
\end{equation}
On another hand, note $Q=\varphi_0(h^{-1}\rho)Q=\tilde\varphi(h^{-1}\rho)Q$ implies
\begin{equation}
    \tilde\varphi(h^{-1}\rho) P\varphi(h^{-1}\rho)Q=\varphi(h^{-1}\rho) P\tilde\varphi(h^{-1}\rho)Q+(\tilde\varphi-\varphi)(h^{-1}\rho) P\varphi_0(h^{-1}\rho)Q,
\end{equation}
the second term of which is $O_{\Psi^{-\infty, -\infty}}(h^\infty)$. With $\mc{E}_M'Q:H^{-N, -N}\rightarrow H^{0,\rardyr_{W_Y}+\frac12-0}$, we have
\begin{equation}
    \Lambda_M^{W_Y}Qu=\varphi(h^{-1}\rho)P\tilde\varphi(h^{-1}\rho)Qu+ O_{H^{0,\rardyr_{W_Y}+\frac12-0}}(h^\infty).
\end{equation}
This implies
\begin{equation}
    \langle(\Lambda_M^{W_Y}-\om) u, v\rangle-\langle u, (\Lambda_M^{W_Y}-\om)v\rangle=\langle[P_\omega^{\varphi}, Q]u, v\rangle+o_h(1).
\end{equation}

2. It remains to look at $[P_\omega^{\varphi}, Q]\in \Psi^{0,-1}$. We remind the reader that the neighborhood $V^+(\omega)$ of $\Sigma^+(\omega)$ inside $T^*\partial  W_Y$ was defined near \eqref{e:5l5}. Noting $\rho\sim h$ on $\WF[P_\omega^{\varphi},Q]$, we have in $V^+(\omega)\cap \{(\tau,\eta)=(-\omega, 0)\}$, 
\begin{equation}
    H_p q=-\frac{\omega \rho^2\partial_\rho q}{\om}+O(S^{-2,-3+0})=-\frac{\omega h^{-1}\rho^2(\partial_\rho\varphi_{00})(h^{-1}\rho)\psi_\omega(|\xi|_g)}{\om}+O(S^{-2,-3+0})\in S^{0, -1}
\end{equation}
with remainder terms up to size $h^{1-0}S^{-1,-1-0}$. Acting on a Lagrangian distribution, we have
\begin{multline}
    [P_\omega^{\varphi}, Q]\varphi\rho^{\frac{d-2}2} u^+(y)e^{i\omega/\rho}=\frac{i\omega}{\om}\rho^{\frac{d-2}2}h^{-1}\rho^2(\partial_\rho\varphi_{00})(h^{-1}\rho)u^+(y)e^{i\omega/\rho}\\
    +O_{\Psi^{-1,-1-0}}(h^{1-0})\big(\varphi\rho^{\frac{d-2}2} u^+(y)e^{i\omega/\rho}\big),
\end{multline}
the last remainder term of which is of size $O_{H^{1,\frac12+0}}(h^{1-0})$. Now
\begin{multline}
    \langle[P_\omega^{\varphi}, Q]\varphi\rho^{\frac{d-2}2} u^+(y)e^{i\omega/\rho}, \varphi\rho^{\frac{d-2}2} v^+(y)e^{i\omega/\rho}\rangle\\
    =\frac{i\omega}{\om}\langle u^+, v^+\rangle_{\partial Y} \int_{0}^1 \rho^{d-2}h^{-1}\rho^2(\partial_\rho\varphi_{00})(h^{-1}\rho)\,\rho^{-d}d\rho+O(h^{1-0})\\=\frac{\omega}{i\om}\langle u^+, v^+\rangle_{\partial Y}+O(h^{1-0}),
\end{multline}
where that integral in $\rho$ is just $\varphi_{00}(h^{-1})-\varphi_{00}(0)=-1$ for small $h$. Similar arguments with $H_p q$ in $V^-(\omega)\cap \{(\tau,\eta)=(\omega, 0)\}$ give
\begin{multline}
    \langle (\Lambda^{W_Y}_M-\om)u,v\rangle_{\partial W_Y}-\langle u,(\Lambda^{W_Y}_M-\om)v\rangle_{\partial W_Y}\\
    =\frac{\omega}{i\om}\Big(\langle u^+,v^+\rangle_{\partial Y}-\langle u^-,v^-\rangle_{\partial Y}\Big)+o_h(1).
\end{multline}
Note all terms above except $o_h(1)$ are independent of $h$ to conclude the proof.
\end{proof}

Now we show the existence of the scattering matrix. 
\begin{proposition}[The scattering matrix]\label{p:existscmatrix}
    Let $W_Y$ be an \acl domain with smooth cross-section $Y$. Let $M\ge 0$ and $\omega>0$. Then for all $a^{-}(y)\in C^\infty(\partial Y)$, there exists
    \begin{equation}
        u\in Z^{\infty,\infty,1-0}_{L_0-0,-\frac12-0,-\frac12-0}(\partial W_Y;\omega), \ (\Lambda^{W_Y}_M-\om)u=0
    \end{equation}
    unique modulo $Z^{\infty, \infty, 1}_{L_0-0, \infty, \infty}(\partial W_Y)$, such that there exist $\rho_0>0$, $\tilde a^\pm(\rho,y) \in S^{0}(\partial W_Y)$ that for each $A\in \Psi^{0,0}(\partial W_Y)$ with $\WF A\subset\{\rho<\rho_0, |\xi|> c\omega\}$ for some $c>0$, we have
\begin{equation}\label{e:PlaneWaveExp}
        A u(\rho,y)=A\rho^{\frac{d}{2}-1} \sum_{\pm}\tilde a^{\pm}(\rho,y;\omega) e^{\pm i \omega/\rho}+O(\omega^{-\infty}\mathcal{S}),
\end{equation}
in local coordinates $y$ on $\partial Y$, uniformly in $\omega>0$. Such $\tilde a^\pm(\rho, y;\omega)$ in $\{\rho<\rho_0\}$ modulo $\omega^{-\infty}\mc{S}(\partial W_Y)$ are uniquely determined by $a^-(y)$ modulo $\omega^{-\infty}\mc{S}(\partial W_Y)$. Moreover, there exist $a^{+}(y;\omega)\in S^0(\partial Y)$, $a^\pm_{k}(\rho, y; \omega)\in S^{0}(\partial W_Y)$ for $k\ge 1$ such that
\begin{equation}\label{e:6l20}
    \tilde a^{\pm}(\rho, y;\omega)=a^{\pm}(y)+\sum_{k=1}^\infty \omega^{-k}\rho^ka^{\pm}_{k}(\rho,y;\omega)+O(\omega^{-\infty}\mc{S}).
\end{equation}
Moreover, the map $\Sigma^{W_Y}_{M,\omega}:a^{-}(y)\mapsto a^{+}(y;\omega)$ is $L^2(\partial Y)$-unitary. 
\end{proposition}
\begin{proof}
     1. Let $\chi,\tilde\chi,\tilde{\tilde{\chi}}\in S^0(\partial W_Y)$ that $\partial W_Y$ is smooth over $\supp\tilde{\tilde{\chi}}$, $\supp(1-\chi)$ is compact and $\chi\prec\tilde\chi\prec\tilde{\tilde{\chi}}$. Pick $B^-\in \Psi^{0,0}$ with $\WF (\Id-B^-)\cap \Sigma^-(\omega)=\emptyset$, $\WF B^-\subset\{\rho<\rho_0, |\xi|> c\omega\}$ for some $c>0$. Let
     \begin{equation}
         P_\omega^-=B^-\chi(P-\om)\tilde\chi\in \Psi^{1,0}. 
     \end{equation}
     In this step, we build $\tilde a^-(\rho, y;\omega), a^-_k(\rho, y;\omega)\in S^0$ such that $P_\omega^- v^i=O (\omega^{-\infty} \mc{S})$ with
     \begin{gather}
         v^i=\tilde{\tilde{\chi}}\rho^\frac{d-2}{2}\tilde a^{-}(\rho, y;\omega)e^{-i\omega/\rho} \in 
        Z^{\infty, \infty, 1}_{L_0-0, -\frac12-0, \infty}, \\ \tilde a^-(\rho, y;\omega)=a^{-}(y)+\sum_{k=1}^\infty (\omega^{-1}\rho)^ka^-_{k}(\rho, y;\omega)+O(\omega^{-\infty}\mc{S}).
     \end{gather}
     Denote $h=\omega^{-1}$. Indeed, one notes the Lagrangian structure of $v^i$ with respect to $\Sigma^-(\omega)$, and via computations similar to \eqref{e:6l12} we have
     \begin{equation}
         P_\omega^-\rho^\frac{d-2}{2}\tilde a^-(\rho, y;h)e^{-i/h\rho}=\frac{i\omega}{\om}\rho^{\frac d2}\Big(-\rho\partial_\rho \tilde a^-(\rho, y;h)+\sum_{k=2}^\infty (h\rho)^{k-1} L^-_{k} \tilde a^{-}(\rho,y;h)\Big)e^{-i/h\rho}+O(h^\infty\mc{S}),
     \end{equation}
     where $L^-_{k}$ are b-differential operators locally involving $\rho\partial_\rho$ and $\partial_y$ up to order $k$. We solve the governing equation
     \begin{equation}
         \rho\partial_\rho \tilde a^-(\rho,y;h)-\sum_{k=0}^\infty (h\rho)^{k+1} L^-_{k+2} \tilde a^-(\rho,y;h)=O(h^\infty\mc{S}), \ \tilde a^-(\rho, y;h)\sim a^{-}(y)+\sum_{k=1}^\infty (h\rho)^ka^-_{k}(\rho, y;h).
     \end{equation}
     At the heights of $h\rho$ and $(h\rho)^2$, the transport equations are
     \begin{equation}
         \rho\partial_\rho(h\rho a^{-}_{1}(\rho,y;h))=h\rho L^-_{2}a^{-}(y), \ \rho\partial_\rho(h^2\rho^2  a^-_{2}(\rho,y;h))=(h\rho)^2 L^-_{3}a^{-}(y)+h\rho L^-_{2}(h\rho a^-_{1}(\rho, y;h)).
     \end{equation}
     One obtains $a^-_{1}$, $a^-_{2}$. The transport equations always give $a^-_{k}$ based on $a^-, a^{-}_{1},\dots, a^{-}_{k-1}$ for $k\ge 2$. One can iteratively solve for all $a^-_{k}\in S^{0}$ and this concludes the step. Note $a^-(y)$ uniquely determines each coefficient and hence $\tilde a^-(\rho ,y;\omega)$, modulo $\omega^{-\infty}\mc{S}$, by arguments similar to Step 5 of Lemma~\ref{l:PlaneWave}. 

    2. For some $B^+$ stated in Lemma~\ref{l:PlaneWave}, we claim
     \begin{equation}
         f=(\Lambda^{W_Y}_M-\om)v^i\in Z^{\infty, \infty, 0}_{L_0-0, \infty, \infty}, \ B^+ f=O(\omega^{-\infty}\mc{S}).  
     \end{equation}
     Here the first bound on $f$ is not uniform in $\omega$. Note $v^i\in Z^{\infty, \infty, 1}_{L_0-0, -\frac12-0, \infty}$ immediately implies that $f\in Z^{\infty, \infty, 0}_{L_0-0, -\frac12-0, \infty}$ in a non-uniform way. It remains to show $B^-\chi f\in \mc{S}$. Indeed, Proposition~\ref{p:dtN} implies
     \begin{equation}
         B^-\chi f=P_\omega^- v^i+B^-\chi \mc{E}_M v^i=O(\omega^{-\infty}\mc{S}), 
     \end{equation}
     since $B^-\chi \mc{E}_M: H^{0, -L_0+0}\rightarrow \omega^{-N}H^{N, N}$ is bounded. On the other hand, we have
     \begin{equation}
         B^+\chi f=B^+\chi(P-\om)\tilde\chi v^i+B^+\chi \mc{E}_M v^i=O(\omega^{-\infty}\mc{S}),
     \end{equation}
     as claimed. 
    
    3. Without loss of generality assume $\ran(\Pi_M)=\ker_{L^2}(\Lambda^{W_Y}_M-\om)\subset Z^{\infty, \infty, 1}_{L_0-0, \infty, \infty}$ by shrinking $\intw$ if necessary. This implies $f=\Pi_M^\perp f$. The limiting absorption principle, Proposition~\ref{p:LAP}, implies there is
   \begin{equation}
       v^o=-(\Pi_M^\perp \Lambda^{W_Y}_M\Pi_M^\perp -\om-i0)^{-1}f\in Z^{\infty, \infty, 1-0}_{L_0-0, \infty, -\frac12-0}
   \end{equation}
   that solves $(\Lambda^{W_Y}_M-\om) v^o=-f$. A priori, we do not have uniform control of $v^o$. Apply Lemma~\ref{l:PlaneWave} to see with $\WF A$ away from the zero frequency and $\Sigma^-(\omega)$, 
   \begin{equation}
       A v^o(\rho,y)=\alpha(\omega)\big(A\rho^\frac{d-2}{2}\tilde a^+(\rho, y) e^{i \omega/\rho}+O(\omega^{-\infty}\mathcal{S})\big),
   \end{equation}
   where
   \begin{equation}
       \tilde a^+(\rho, y)=a^{+}(y)+\sum_{k=1}^\infty \omega^{-k}\rho^ka^+_{k}(\rho,y) +O(\omega^{-\infty}\mathcal{S}),
   \end{equation}
   for some $a^+(y;\omega), a^+_{k}(\rho, y;\omega)\in S^0$, uniquely determined modulo $\omega^{-\infty}\mc{S}$ by $a^-(y)$ once $\alpha(\omega)$ is fixed. Now set $u=v^o+v^i$. We have $(\Lambda^{W_Y}_M-\om)u=0$ with
   \begin{equation}
        A u(\rho,y)=A\rho^{\frac{d}{2}-1} \Big(\tilde a^-(\rho,y;\omega)e^{-i \omega/\rho}+\alpha(\omega)\tilde a^+(\rho,y;\omega) e^{i \omega/\rho}\Big)+O((1+\alpha(\omega))\omega^{-\infty}\mathcal{S}).
    \end{equation}
   We remind the reader that we defined in Step 1 of Lemma~\ref{l:PlaneWave} that uniformly for $\omega>0$ large,
   \begin{equation}
       \alpha(\omega)=\max(\|B\tilde{\tilde{\chi}}v^o\|_{H^{0, -\frac12-0}},1)\le \max(K\alpha(\omega)\|a^+(y;\omega)\|_{L^2(\partial Y)},1)
   \end{equation}
   for some constant $K>0$. Now a posteriori, the Green's identity, Lemma~\ref{l:green} implies
   \begin{equation}
       \alpha(\omega)^2\|a^+(y;\omega)\|_{L^2(\partial Y)}^2=\|a^-(y)\|_{L^2(\partial Y)}^2.
   \end{equation}
   The above implies
   \begin{equation}
       \alpha(\omega)\le \max(K\|a^-(y)\|_{L^2(\partial Y)},1)
   \end{equation}
   is uniformly bounded for $\omega>0$ large. This allows us to renormalize $\tilde a^+(\rho, y;\omega)$ such that $\alpha(\omega)\equiv 1$. This gives \eqref{e:PlaneWaveExp} and $\Sigma^{W_Y}_{M,\omega}:a^{-}\mapsto a^{+}$ is $L^2$-unitary. 
\end{proof}

As a corollary, we show that there is no meaningful reflection on an \acs domain of $d=2$. 
\begin{corollary}[Scattering matrix on a smooth curve]\label{c:noreflection}
    Let $\Omega$ be a $2$-dimensional \acs domain whose boundary $\partial\Omega$ is a connected smooth curve. Let $M\ge 0$ and $\omega>0$. Then for all $\mathbf{a}^i=(a_+^i, a_-^i)\in \mathbb{C}^2$, there exists a unique $\mathbf{a}^o=(a_+^o, a_-^o)\in \mathbb{C}^2$ depending on $M, \omega$ that \eqref{e:PlaneWaveExp}, \eqref{e:6l20} hold with
    \begin{equation}
        a^+(\pm\infty)=a^o_\pm, \ a^-(\pm\infty)=a^i_\mp. 
    \end{equation}
    Moreover, there exist $\theta_\pm(M,\omega)\in[0, 2\pi)$ that
    \begin{equation}\label{e:6l15}
        \begin{pmatrix}
            a^{o}_-\\
            a^{o}_+
        \end{pmatrix}=\Sigma^\Omega_{M,\omega}\begin{pmatrix}
            a^{i}_+\\
            a^{i}_-
        \end{pmatrix}, \ \Sigma^\Omega_{M,\omega}=\begin{pmatrix}
            0 & e^{i\theta_{-}}\\
            e^{i\theta_{+}} & 0
        \end{pmatrix}+O(\omega^{-\infty}),
    \end{equation}
    uniformly for $\omega>0$ large. 
\end{corollary}
\begin{proof}
    Proposition~\ref{p:existscmatrix} gives the desired $u\in Z^{\infty,\infty}_{L_0-0, -\frac12-0, -\frac12-0}$, $(\Lambda^{\Omega}_M-\om)u=0$. It suffices to show \eqref{e:6l15}. We claim
    \begin{equation}\label{e:6l16}
        \psi_\omega(D)(\Lambda^{\Omega}_M-\om)\psi_\omega(D) u=O(\omega^{-\infty}\mc{S}),
    \end{equation}
    where $\psi_\omega\in C^\infty_c((0,\infty))$ defined in Definition~{d:Zspace} satisfies $\omega\in \supp\psi_\omega$ and $0\notin \supp\psi_\omega$. Indeed, 
    \begin{equation}
        [\Lambda_M^\Omega, \psi_\omega(D)]=[P, \psi_\omega(D)]+\mc{E}_\flat\psi_\omega(D)-\psi_\omega(D)\mc{E}_\flat.
    \end{equation}
    Note $\WF u,\WF[P, \psi_\omega(D)]$ are separated at scale $\omega$ in $\xi$, and $\mc{E}_\flat$ is in $\Psi^{-\infty,-\infty}$ away from the zero frequency in Proposition~\ref{p:dtNParametrix}. They imply \eqref{e:6l16}. Now Green's identity, Lemma~\ref{l:green} implies
    \begin{equation}
        O(\omega^{-\infty})=2i\Im \langle (\Lambda^{\Omega}_M-\om)\psi_\omega(D)u,\psi_\omega(D)u\rangle_{\partial \Omega}=\frac{\omega}{i\om}\Big(|a^o_+|^2-|a^i_+|^2\Big).
    \end{equation}
    The proof giving $|a_-^o|\sim|a_-^i|$ is similar. 
\end{proof}
\begin{remark}
    All constants in this section can depend on $\Omega$ in the following way: the upper bound on uniform Lipschitz characters of $\partial\Omega$, the upper bound on the choice of $r_0\ge 0$, and the lower bound on the radial Hardy constant $\rardyr_{\Omega}$, the upper bound on the diameter of $\partial\Omega\setminus\partial\tilde\Omega$ in Remark~\ref{r:aclcollar}, for each multi-index $\alpha$, the upper bounds on the derivatives of the metric $\gplus$ in \eqref{e:3l26}, and lower bounds on $|\det \gplus|$.
\end{remark}

\part{The Steklov spectrum on curvilinear polygons}

\section{Convergence of the scattering matrix in two dimension}
\label{s:scatteringatrixConvergence}

This section considers scattering matrices on simply connected \acl domains in $\mathbb{R}^2$ with a single point singularity in the high energy limit. We call these domains \emph{curved wedges with angle $\alpha$} if the two smooth components of the boundary meet at angle $\alpha\in(0,2\pi)$. The main goal if this section is to determine the extent to which the scattering matrix is determined by $\alpha$ together with the derivatives of the smooth parts of the boundary at the singular point. 

Specifically, letting $\rardya=\rardy^0_{W_\alpha}=\frac\pi\alpha$ and $\mc{W}_{\alpha},\tilde{\mc{W}}_{\alpha}\subset\mathbb{R}^2$ curved wedges of corner angle $\alpha$ with strong radial Hardy constants $\rardy_{\mc{W}_\alpha}^0, \rardy_{\tilde{\mc{W}}_\alpha}^0$. Without loss of generality, we assume that the singular point for $\mc{W}_{\alpha}$ and $\tilde{\mc{W}}_{\alpha}$ is at $0$. Then, setting $\min(\rardy_{\mc{W}_\alpha}^0,\rardy_{\tilde{\mc{W}}_\alpha}^0)=\rardya-0$, we obtain 
\begin{proposition}\label{p:scconverge}
Let $M_0>0$, $\omega\in[\frac12, 2]$, $k_0\geq 0$ such that the first $(k_0+1)$ derivatives of the boundary of $\mc{W}_{\alpha}$ agreeing with those for $\tilde{\mc{W}}_{\alpha}$ at $0$. Then there exists $h_0>0$ such that uniformly for all $M\in[0, M_0]$, all $h\in(0, h_0]$, we have
    \begin{equation}\label{e:7l26}
        \Sigma_{M, h^{-1}\omega}^{\tilde{\mc{W}}_\alpha}=\Sigma_{0, \omega}^{\mc{W}_{\alpha}}
        +O(h^{\min(2\rardya,1+k_0)-0})+O((hM)^{2-0}).
    \end{equation}
\end{proposition}
\begin{remark}
    In Section~\ref{s:grushin}, for any $\e>0$, we will apply Proposition~\ref{p:scconverge} to domains that satisfy $\rardy_{\mc{W}_\alpha}^0>\rardya-\e$, $\rardy_{\tilde{\mc{W}}_\alpha}^0>\rardya-\e$. 
\end{remark}

\noindent{\bf Outline of the strategy}

We prove Proposition~\ref{p:scconverge} in two steps. We first compare $\Sigma^{\tilde{\mc{W}}_\alpha}_{M,h^{-1}\omega}$ with $\Sigma^{\mc{W}_{\alpha}}_{M,h^{-1}\omega}$ and then compare $\Sigma^{\mc{W}_{\alpha}}_{M,h^{-1}\omega}$ with $\Sigma^{\mc{W}_{\alpha}}_{0,h^{-1}\omega}$.  The first comparison modifies the underlying geometry, while the second only modifies the mass. These correspond to Lemmas~\ref{l:scstraighten}, and ~\ref{l:schfmzero}, respectively, which together imply the Proposition. 

\noindent Geometric comparison:

When comparing scattering matrices on two different wedges, we consider a solution, $u$, to
$$
(\Lambda_{M}^{\tilde{\mc{W}}_{\alpha}}-h^{-1}\ohm )u=0.
$$
Our goal is then to transplant $u$ onto $\mc{W}_\alpha$, as the function $v'$, in such a way that it satisfies 
$$
(\Lambda_{M}^{\mc{W}_{\alpha}}-h^{-1}\ohm )v'=\text{small}.
$$

More precisely, we find $F:B(0,\e)\cap \mc{W}_{\alpha}\to \tilde{\mc{W}}_{\alpha}$ a diffeomorphism that is as close to the identity as possible. We then zoom in to scale $h$, considering $\mc{W}^h_{\alpha}=:h^{-1}\mc{W}_\alpha$ and $\tilde{\mc{W}}^h_{\alpha}:=h^{-1}\tilde{\mc{W}}_{\alpha}$, the map $F_h:= h^{-1}F(hx)$ and consider $u_h:= u(hx)$ satisfying
\begin{equation}
\label{e:rescaledWave}
(\Lambda_{hM}^{\tilde{\mc{W}}^h_{\alpha}}-\ohm )u_h=0.
\end{equation}
Recall from Proposition~\ref{p:existscmatrix} that such a $u_h$ exists and is can be decomposed as 
\begin{equation}
    u=u_{h,\text{low}}+u_{h,\text{wave}},
\end{equation}
where $u_{h,\text{low}}$ is low frequency near infinity  and $u_{h,\text{wave}}$ is an approximate plane wave. Moreover, $u_{h,\text{wave}}$ is a sum of approximate plane waves that is determined by only the scattering data. 

We then build $v'_h$ out of three pieces:
\begin{equation}
\label{e:extensionType}
v'(h\cdot):=v'_h= \chi^\delta F_h^*u + (\chi-\chi^\delta)F_h^*u_{h,\text{low}}+(1-\chi^\delta)(w_{\text{right}}+w_{\text{left}}).
\end{equation}
Here, $\chi^\delta$ is a localizer to $|x|\lesssim h^{-\delta}$ for $0<\delta\ll 1$, and $\chi$ to $|x|\lesssim h^{-1}$, $u_{h,\text{low}}$ denotes the low frequencies of $u$, and $w_{\text{right}/\text{left}}$ are approximate plane waves on the right and left of the singularity of $\mc{W}_{\alpha}$. The reason for the choice of scales is as follows:
\begin{enumerate}
\item  The asymptotic properties of $u_{h,\text{wave}}$ set in (modulo $h^\infty$) at scale $h^{-\delta}$ so, at this scale, we can construct $w_{\text{right}/\text{left}}$ on $\mc{W}^h_{\alpha}$ that are very close to $F^*u_{h,\text{wave}}$ at scale $h^{-\delta}$ and are very accurate quasimodes on $\mc{W}_{\alpha}^h$.  This allows us to match $F^*u_{h,\text{wave}}$ at scale $h^{-\delta}$ and appropriately extend to infinity.
\item We have very little control over $u_{low}$ and, in addition, we expect its properties to be global. Therefore, although we must truncate these low frequencies somewhere, we want to wait as long as possible to do so. Since we can only accurately compare the geometries to scale $h^{-1}$, this is the scale we true.
\end{enumerate}

Estimating
\begin{equation}
    (\Lambda_{hM}^{\mc{W}^h_{\alpha}}-\ohm )v_h'
\end{equation}
is technically involved, so we do not discuss details here. We instead point out one important feature that is crucial in proving Proposition~\ref{p:scconverge}. From the Fredholm estimates in Section~\ref{s:microlocalDtN}, one apriori knows that $u_{h,\text{low}}$ decays like $|x|^{-\frac{\pi}{\alpha}-\frac{1}{2}+0}$. At first, it seems that this is not sufficient to obtain control better than $h^{\frac{\pi}{\alpha}+\frac{1}{2}-0}$. However, by measuring the error in a space that allows growth at low frequency like $|x|^{\frac{\pi}{\alpha}+\frac{1}{2}-0}$, we are able to recover the difference. 

Once we show that the right $(\Lambda_{hM}^{\mc{W}^h_{\alpha}}-\ohm)v_h'$ is small enough, an application of the outgoing resolvent will finish the proof.

\noindent Mass comparison:

The mass comparison is significantly simpler than the geometric comparison, but the process is similar. We begin we a solution of 
\begin{equation}
    (\Lambda^{\mc{W}_\alpha}_M-h^{-1}\sqrt{\omega^2+h^2M^2}) u=0 
\end{equation}
and show that $u$ satisfies
\begin{equation}
\label{e:merror}
(\Lambda^{\mc{W}_\alpha}_0-h^{-1}\omega) u=\text{small}. 
\end{equation}
Applying the outgoing resolvent then finishes the proof.

The keys to estimating the error in~\eqref{e:merror} are, again, the apriori structure of $u$ and an accurate microlocal description of $\Lambda_M^{\mc{W}_\alpha}$ near infinity. 

\subsection{Rescaling \acl domains with exactly one corner}\label{s:8-1}

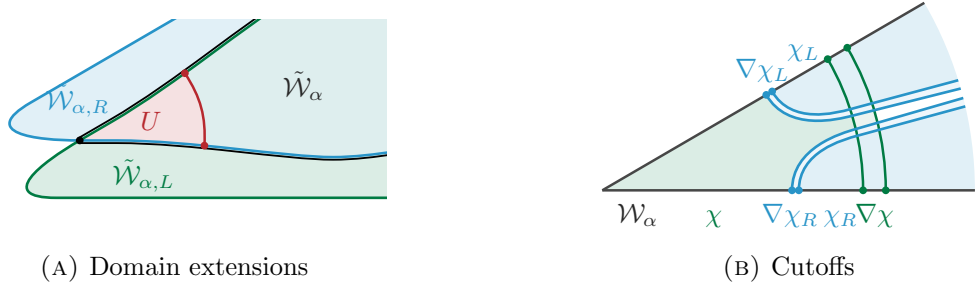
\begin{figure}[t]
  \centering
  \begin{subfigure}[t]{0.5\textwidth}
        \centering
        \begin{tikzpicture}[x=1cm,y=1cm]
  \usetikzlibrary{decorations.markings}
  \definecolor{colorAKANE}{RGB}{183,40,46}
  \definecolor{colorTOKIWA}{RGB}{0,123,67}
  \definecolor{colorUSUAI}{RGB}{40,148,200}
  \colorlet{plotPrimary}{colorAKANE}
  \colorlet{plotSecondary}{colorTOKIWA}
  \colorlet{plotTertiary}{colorUSUAI}
  \colorlet{plotInterpolated}{plotSecondary!50!plotTertiary}
  \colorlet{plotAxis}{black!72}

  \tikzset{
    >=latex,
    font=\small,
    line cap=butt,
    line join=round,
    wedge edge/.style={draw=plotAxis,line width=0.95pt},
    loop point/.style={fill=plotSecondary,draw=none},
    edge point/.style={fill=plotTertiary,draw=none},
    pair arc/.style={draw=plotSecondary,line width=0.95pt,line cap=round},
    escape curve/.style={
      draw=plotTertiary,
      line width=1pt,
      line cap=round,
      line join=round
    },
    nearby edge/.style={draw=black,line width=0.7pt,line cap=round},
    inner region/.style={fill=plotSecondary!14,draw=none},
    pair region/.style={fill=plotSecondary!8,draw=none},
    blue pair region/.style={fill=plotTertiary!8,draw=none},
    upper blue region/.style={fill=plotTertiary!14,draw=none},
    lower blue region/.style={fill=plotTertiary!14,draw=none},
    green lower region/.style={fill=plotSecondary!14,draw=none},
    interpolated region/.style={fill=plotInterpolated!14,draw=none},
    red sector region/.style={fill=plotPrimary!14,draw=none},
    red region label/.style={text=plotPrimary,inner sep=0.8pt},
    blue curve label/.style={text=plotTertiary,inner sep=0.8pt},
    green curve label/.style={text=plotSecondary,inner sep=0.8pt},
    figure label/.style={text=black!84,inner sep=0.8pt},
    intersection point/.style={fill=black,draw=none},
    one-third coordinate/.style args={#1}{
      decoration={
        markings,
        mark=at position {1/3} with {
          \coordinate (#1) at (0,0);
        }
      },
      postaction={decorate}
    }
  }

  \pgfmathsetmacro{\figureWidth}{7.5}
  \pgfmathsetmacro{\figureHeight}{3.2}
  \pgfmathsetmacro{\frameBottom}{-0.28}
  \pgfmathsetmacro{\frameTop}{\figureHeight+\frameBottom}
  \pgfmathsetmacro{\edgeDx}{4.26}
  \pgfmathsetmacro{\edgeDy}{2.48}
  \pgfmathsetmacro{\tipY}{\frameBottom+2*\figureHeight/5}
  \pgfmathsetmacro{\edgeLength}{veclen(\edgeDx,\edgeDy)}
  \pgfmathsetmacro{\tipX}{\figureWidth/3}
  \pgfmathsetmacro{\wedgeAngle}{atan2(\edgeDy,\edgeDx)}
  \pgfmathsetmacro{\innerClipLeft}{0.125*\figureWidth}
  \pgfmathsetmacro{\innerClipRight}{0.875*\figureWidth}
  \pgfmathsetmacro{\innerClipBottom}{\frameBottom+0.1*\figureHeight}
  \pgfmathsetmacro{\innerClipTop}{\frameBottom+0.9*\figureHeight}
  \pgfmathsetmacro{\extensionShoulderX}{-0.65}
  \pgfmathsetmacro{\extensionShoulderY}{0.50}
  \pgfmathsetmacro{\extensionFarX}{
    \extensionShoulderX+3*\figureWidth*cos(\wedgeAngle)
  }
  \pgfmathsetmacro{\extensionFarY}{
    \extensionShoulderY+3*\figureWidth*sin(\wedgeAngle)
  }

  \path[use as bounding box]
    (0,{\frameBottom}) rectangle ({\figureWidth},{\frameTop});
  \begin{pgfinterruptboundingbox}

  \coordinate (tip) at ({\tipX},{\tipY});
  \coordinate (domainLabel) at
    ({0.75*(\innerClipRight-0.12)+0.25*0.5*\figureWidth},
     {0.75*(\innerClipTop-0.10)+0.25*0.5*(\frameBottom+\frameTop)-0.26});
  \coordinate (blueLeftLabel) at
    ({\tipX-0.03},{\tipY+0.28});
  \coordinate (greenLeftLabel) at
    ({\tipX+0.40},{\tipY-0.22});

  \def\curvedWedgeEdgeFromTip{%
    -- ({0.1*\edgeLength},0)
    .. controls
      ({0.30*\edgeLength},0)
      and ({0.55*\edgeLength},-0.22)
      .. ({0.68*\edgeLength},-0.22)
    .. controls
      ({0.79*\edgeLength},-0.22)
      and ({0.92*\edgeLength},-0.08)
      .. ({\edgeLength},0)%
  }
  \def\curvedWedgeEdge{%
    (0,0) \curvedWedgeEdgeFromTip
  }
  \def\nearbyWedgeEdge{%
    (0,-0.035)
    -- ({0.1*\edgeLength},-0.035)
    .. controls
      ({0.30*\edgeLength},-0.035)
      and ({0.55*\edgeLength},-0.255)
      .. ({0.68*\edgeLength},-0.255)
    .. controls
      ({0.79*\edgeLength},-0.255)
      and ({0.92*\edgeLength},-0.115)
      .. ({\edgeLength},-0.035)%
  }
  \def\curvedWedgeExtension{%
    ({\extensionFarX},{\extensionFarY})
    -- ({\extensionShoulderX},{\extensionShoulderY})
    .. controls
      ({\extensionShoulderX-0.80*cos(\wedgeAngle)},
       {\extensionShoulderY-0.80*sin(\wedgeAngle)})
      and (-0.60,0)
      .. (0,0)%
  }
  \def\extendedWedgeRegionPath{%
    \curvedWedgeExtension
    \curvedWedgeEdgeFromTip
    -- ({3*\figureWidth},0)
    -- cycle%
  }

  \begin{scope}
  \clip
    (0,{\frameBottom}) rectangle ({\figureWidth},{\frameTop});
  \clip
    ({\innerClipLeft},{\innerClipBottom})
    rectangle
    ({\innerClipRight},{\innerClipTop});

  \begin{scope}[shift={(tip)}]
    \path[one-third coordinate={bottomRedDot}]
      \curvedWedgeEdge;
    \begin{scope}[rotate=\wedgeAngle,yscale=-1]
      \path[one-third coordinate={topRedDot}]
        \curvedWedgeEdge;
    \end{scope}
  \end{scope}

  \begin{scope}[shift={(tip)}]
    \fill[upper blue region] \extendedWedgeRegionPath;
  \end{scope}

  \begin{scope}[shift={(tip)}]
    \begin{scope}[rotate=\wedgeAngle,yscale=-1]
      \fill[green lower region] \extendedWedgeRegionPath;
    \end{scope}
  \end{scope}

  \begin{scope}[shift={(tip)}]
    \clip \extendedWedgeRegionPath;
    \begin{scope}[rotate=\wedgeAngle,yscale=-1]
      \fill[interpolated region] \extendedWedgeRegionPath;
    \end{scope}
  \end{scope}

  \node[figure label,anchor=north east]
    at (domainLabel) {$\tilde{\mathcal{W}}_\alpha$};

  \fill[red sector region]
    (tip) -- (bottomRedDot)
    to[bend right=18] (topRedDot) -- cycle;
  \draw[draw=plotPrimary,line width=0.95pt,line cap=round]
    (bottomRedDot) to[bend right=18] (topRedDot);
  \path
    (bottomRedDot) --
    coordinate[midway] (redChordMid)
    (topRedDot);
  \path
    (tip) --
    coordinate[pos=0.62] (redRegionLabel)
    (redChordMid);
  \node[red region label] at (redRegionLabel) {$U$};

  \begin{scope}[shift={(tip)}]
    \draw[escape curve,draw=plotTertiary] \curvedWedgeExtension;
    \begin{scope}[rotate=\wedgeAngle,yscale=-1]
      \draw[escape curve,draw=plotSecondary] \curvedWedgeExtension;
    \end{scope}
    \draw[nearby edge] \nearbyWedgeEdge;
    \draw[wedge edge,draw=plotTertiary]
      \curvedWedgeEdge;
    \begin{scope}[rotate=\wedgeAngle,yscale=-1]
      \draw[nearby edge] \nearbyWedgeEdge;
      \draw[wedge edge,draw=plotSecondary]
        \curvedWedgeEdge;
    \end{scope}
  \end{scope}
  \path[intersection point] (tip) circle[radius=1.45pt];
  \path[fill=plotPrimary,draw=none]
    (bottomRedDot) circle[radius=1.45pt]
    (topRedDot) circle[radius=1.45pt];
  \node[blue curve label,anchor=south]
    at (blueLeftLabel) {$\tilde{\mathcal{W}}_{\alpha,R}$};
  \node[green curve label,anchor=north west]
    at (greenLeftLabel) {$\tilde{\mathcal{W}}_{\alpha,L}$};
  \end{scope}
  \end{pgfinterruptboundingbox}
\end{tikzpicture}%
        \caption{Domain extensions}
  \end{subfigure}\hfill
  \begin{subfigure}[t]{0.5\textwidth}
        \centering
        \includestandalone{../images/s7f2}
        \caption{Cutoffs}
  \end{subfigure}
  \caption{Illustration of domain extensions and cutoffs}
  \label{f:7l1}
\end{figure}

Consider an \acl domain $\tilde{\mc{W}}_\alpha\subset\mathbb{R}^2$ with exactly one corner of central angle $\alpha\in (0, 2\pi)$ positioned at $z=0$, and suppose that we have an $h$-family of $u_h$ that $(h\Lambda_M^{\tilde{\mc{W}}_{\alpha}}-\omega)u_h=0$ for $\omega\in [\frac12, 2]$ fixed. Its boundary $\partial\tilde{\mc{W}}_\alpha$ can be parametrized by the arclength variable $x$, where $x=0$ at the corner. Let $\partial\tilde{\mc{W}}_\alpha=\Gamma_L\cup\{x=0\}\cup\Gamma_R$, where each of $\Gamma_{L}=\{x\in \partial\tilde{\mc{W}}_\alpha: x<0\}$ and $\Gamma_{R}=\{x\in \partial\tilde{\mc{W}}_\alpha: x>0\}$ is a smooth ray that originates from the corner and is asymptotically linear as $|x|\rightarrow \infty$. There then exist \acs extension $\tilde{\mc{W}}_{\alpha, L/R}\subset \mathbb{R}^2$ such that $\Gamma_{L/R}\subset\partial \tilde{\mc{W}}_{\alpha, L/R}$ with $\rardy^\infty_{\tilde{\mc{W}}_{\alpha, L/R}}\ge \rardya$. See Figure~\ref{f:7l1}(A). 

Let $k_0\ge 0$. Our goal in this subsection is to relate $\tilde{\mc{W}}_\alpha$ to $\mc{W}_\alpha$, another curved wedge in $\mathbb{R}^2$ cornered at $z=0$ with corner angle $\alpha$, with first $(k_0+1)$ derivatives of the boundary at $z=0$ agreeing with those for $\tilde{\mc{W}}_\alpha$. In this section, we use $\cl(\cdot)$ to denote the closure of an open set to avoid confusion with bar and tilde notations. Let $U= \cl\big(\tilde{\mc{W}}_\alpha\cap B(0, R)\big)$ for some $R>0$ fixed large, and $F:U\to F(U)\subset\cl(\mc{W}_{\alpha})$ be a smooth map such that $F|_{\partial\tilde{\mc{W}}_\alpha \cap \partial U}$ is an isometry, $dF(0)=\Id$, with $dF(\partial_{\nu_{\tilde{\mc{W}}_\alpha}})=\partial_{\nu_{\mc{W}_{\alpha}}}$ for $\partial_{\nu_{\tilde{\mc{W}}_\alpha}}$ and $\partial_{\nu_{\mc{W}_{\alpha}}}$ being respectively the outward pointing normal to $\partial\tilde{\mc{W}}_\alpha$ and $\partial\mc{W}_\alpha$. We ask 
\begin{equation}
    F(z)=z+\sum_{k=2}^\infty F_k(z)=z+O(|z|^{k_0+2}), \ F_k=\frac{1}{k!}\nabla_z^kF(0)(z, \dots, z), F_k\equiv 0\text{ for }2\le k\le k_0+2, 
\end{equation}
where $F_k(\lambda z)=\lambda^kF_k(z)$ is homogeneous of degree $k$ for $\lambda\in \mathbb{R}$. The existence of such $F$ follows from Lemma~\ref{l:Fprop}. Its pushforward is
\begin{equation}
    DF(z)=\Id+\sum_{k=2}^\infty DF_k(z)=\Id+O(|z|^{k_0+1}),\ DF_k(z)z'=\frac{1}{(k-1)!}\nabla^k_z F(0)(z, \dots, z, z'), 
\end{equation}
with $DF_k$ homogeneous of degree $k-1$ in $z$. We denote $k_1=k_0+1$.

We now zoom in near $z=0$ and rescale the domains by 
\begin{equation}
    \tilde{\mc{W}}^{h}_\alpha=\{z\in\mathbb{R}^2: hz\in \tilde{\mc{W}}_{\alpha}\}, \ \mc{W}^{h}_\alpha=\{z\in\mathbb{R}^2: hz\in \mc{W}_{\alpha}\}.
\end{equation}
Then, setting $u(z)=u_h(hz)$, we have
\begin{equation}
    (\Lambda_{hM}^{\tilde{\mc{W}}^{h}_\alpha}-\omega)u=0.
\end{equation}
Let $U^h=\{z: hz\in U\}\subset \tilde{\mc{W}}_{\alpha}^h$, and $F_h(z)=h^{-1}F(hz)$ as a map from $U^h\subset \tilde{\mc{W}}_\alpha^h$ to  $\mc{W}^h_{\alpha}$. This implies that $dF_h(z)=\Id+O(h^{k_1}|z|^{k_1})$. For the convenience of our analysis, we will introduce another scale $0<\delta\ll 1$, by 
\begin{equation}
    \tilde{\mc{W}}^{\delta}_\alpha=\{z\in\mathbb{R}^2: h^{1-\delta}z\in \tilde{\mc{W}}_{\alpha}\}, \ \mc{W}^{\delta}_\alpha=\{z\in\mathbb{R}^2: h^{1-\delta}z\in \mc{W}_{\alpha}\},
\end{equation}
and $U^\delta=\{z:h^{1-\delta}z\in U\}$. We also have $F_\delta=h^{-1+\delta}F(h^{1-\delta}z)$ and
\begin{equation}
    (h^\delta\Lambda_{h^{1-\delta}M}^{\tilde{\mc{W}}^{\delta}_\alpha}-\omega) u_\delta=0, \ u_\delta(z)=u_h(h^{1-\delta} z).
\end{equation}
We see $\tilde{\mc{W}}_{\alpha}, \tilde{\mc{W}}_{\alpha}^\delta, \tilde{\mc{W}}_{\alpha}^h$ are the domains that we see at the physical scale, matching scale, and semiclassical scale. We rescale $\tilde{\mc{W}}_{\alpha, L/R}$ to obtain $\tilde{\mc{W}}_{\alpha, L/R}^{h}, \tilde{\mc{W}}_{\alpha, L/R}^{\delta}$ accordingly. 

We denote the unitary rescaling operator
\begin{gather}
    T_h: L^2(\partial\mc{W}_\alpha)\rightarrow L^2(\partial\mc{W}^h_\alpha), \ (T_h u)(z)=h^{\frac12}u(hz),\\
    T_\delta: L^2(\partial\mc{W}_\alpha)\rightarrow L^2(\partial\mc{W}^\delta_\alpha), \ (T_\delta u)(z)=h^{\frac{1-\delta}2}u(h^{1-\delta}z).
\end{gather}
These rescaling operators naturally extends to other rescaled domains. 

We now set up some cutoff functions we will need for Lemma~\ref{l:scstraighten}. Throughout this subsection and Lemma~\ref{l:scstraighten}, we use superscripts $h$ and $\delta$ to denote functions that lives on $\mc{W}_\alpha$ at a scale $|x|\sim h$ and $|x|\sim h^{1-\delta}$ respectively, and to denote domains in which a unit ball looks likes a ball on $\mc{W}_\alpha$ at a scale $|x|\le h$ and $|x|\le h^{1-\delta}$ respectively. On $\mc{W}_\alpha$, consider $\chi\in C^\infty_c(F(U))$ that $0\notin\supp(1-\chi)\cap\partial\mc{W}_\alpha$ and $\partial_\nu \chi|_{\partial\mc{W}_\alpha}\equiv 0$. Let $\tilde\chi,\tilde{\tilde{\chi}}$ satisfy the same conditions as $\chi$, and $\chi\prec \tilde\chi\prec \tilde{\tilde\chi}$. Consider the following rescaled functions
\begin{equation}
    \chi^h(z)=\chi(h^{-1}z)\in C_c^\infty(F^h(U^h)), \ \chi^\delta(z)=\chi(h^{\delta-1}z)\in C_c^\infty(F^\delta(U^\delta)),
\end{equation}
We understand $\chi\in C^\infty_c(U)$ via the pullback induced by $F$, and $\chi^h\in C^\infty_c(U^h)$, $\chi^\delta\in C^\infty_c(U^\delta)$ similarly. All functions with superscripts $h$ or $\delta$ are scaled accordingly, such as $\tilde\chi^h,\tilde\chi^\delta,\tilde{\tilde\chi}^h,\tilde{\tilde\chi}^\delta$. On $\mc{W}_\alpha$, we let $\chi_{L/R}\in S^0(\mc{W}_\alpha)$ such that $0\notin\supp\chi_{L/R}$, $\chi_{L/R}\equiv 0$ on $\Gamma_{R/L}$, $\chi_L+\chi_R\equiv 1$ on $\partial \mc{W}_\alpha\cap \supp(1-\chi)$, and $\partial_\nu \chi_{L/R}|_{\partial\mc{W}_\alpha}\equiv 0$. We ask $\nabla\chi_{L/R}$ is supported inside $\{C^{-1}<|x|,x_0>C^{-1}\langle x\rangle-C\}$ in $\tilde{\mc{W}}_{\alpha, L/R}$. Rescale $\chi_{L/R}$ to obtain $\chi_{L/R}^h$, $\chi_{L/R}^\delta$. See Figure~\ref{f:7l1} (B).

\begin{figure}[t]
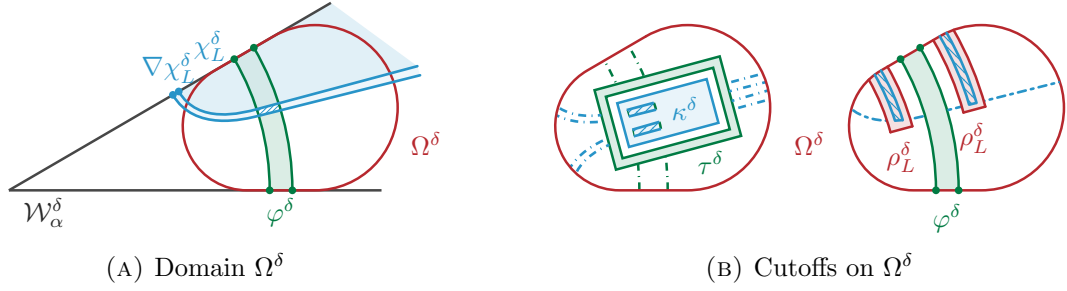

  \centering
  \begin{subfigure}[t]{0.5\textwidth}
        \centering
        \includestandalone{../images/s7f1}
        \caption{Domain $\Omega^{\delta}$}
  \end{subfigure}\hfill
  \begin{subfigure}[t]{0.5\textwidth}
        \centering
        \includestandalone{../images/s7f4}
        \caption{Cutoffs on $\Omega^{\delta}$}
  \end{subfigure}
  \caption{Illustration of $\Omega^{\delta}$ and cutoffs on $\Omega^{\delta}$.}
  \label{f:7l2}
\end{figure}

On $\mc{W}^\delta_\alpha$ we fix $h$-independent $\varphi^\delta$, $\tilde\varphi^\delta\in C^\infty(\cl(\mc{W}^\delta_\alpha))$ such that $(\nabla \chi^\delta)\varphi^\delta=\nabla \chi^\delta$, $\varphi^\delta\tilde\varphi^\delta=\varphi^\delta$. Note $\supp\varphi^\delta$ is compact on $\mc{W}_\alpha^\delta$ and away from 0, there exists an $h$-uniform compact smooth domain $\Omega^\delta\subset \mc{W}_\alpha^\delta$ such that $\supp\tilde\varphi^\delta\subset \Omega^\delta$, $\partial\Omega^\delta\cap \supp\tilde\varphi^\delta=\partial\mc{W}_\alpha\cap\supp\tilde\varphi^\delta$. Then $\varphi^\delta$, $\tilde\varphi^\delta$ can be viewed as functions in $C^\infty(\Omega^\delta)$. We can pull back $\Omega^\delta$ to obtain $\tilde\Omega^\delta=F_\delta^{-1}(\Omega^\delta)\subset \tilde{\mc{W}}_\alpha^\delta$, whose shapes then depends on $h$. Making $\Omega^\delta$ larger if necessary, we let $\kappa^\delta, \tilde\kappa^\delta\in C_c^\infty(\Omega^\delta)$ be that $\tilde\varphi^\delta (\nabla \chi^\delta_\pm)\prec\kappa^\delta$ and $\kappa^\delta\prec \tilde\kappa^\delta$, and let $\tau^\delta\in C_c^\infty(\Omega^\delta)$ be that $\nabla \tilde\kappa^\delta\prec\tau^\delta$ and $\supp\tau^\delta\cap\supp\kappa^\delta=\emptyset$. Note $\tilde\kappa^\delta$ and $\tau^\delta$ are supported away from $\partial\Omega^\delta$. Let $\rho^\delta_{L/R}\in C^\infty(\bar\Omega^\delta)$ be that $\nabla \tilde\varphi^\delta\chi_{L/R}^\delta\prec\rho^\delta_{L/R}$ and $\supp\rho^\delta_{L/R}\cap\supp\varphi^\delta=\emptyset$. See Figure~\ref{f:7l2}. 

\subsection{Uniform resolvent estimates along family of curves}
Denote $L_0=\rardya+\frac12-0$. We also record the 2D resolvent estimates in terms of the variable-order Sobolev spaces:
\begin{lemma}[Small-mass 2D uniform absorption resolvent]\label{l:uniformres}
    Let $M_0>0$. Let $s,l\in\mathbb{R}$, $s_0\in(0,1]$, $l^->-\frac12$, $l^+<-\frac12$, $|l_0|<L_0$. Then there are $C>0$ and $h_0>0$ that we have
    \begin{equation}
        \|u\|_{Z^{s,l, s_0}_{l_0, l^-, l^+}(\partial\tilde{\mc{W}}_\alpha^h)}\le C\|(\Lambda^{\tilde{\mc{W}}_\alpha^h}_{hM}-\omega-i\mu)u\|_{Z^{s-1, l, s_0-1}_{l_0, l^-+1, l^++1}(\partial\tilde{\mc{W}}_\alpha^h)}
    \end{equation}
    uniformly for $u\in Z^{-\infty,-\infty, 0}_{-L_0+0, -\frac12+0, -\infty}(\partial\tilde{\mc{W}}_\alpha^h; \omega)$, $h\in(0,h_0]$, $\mu\ge 0$, $\omega\in [\frac14, 4]$, $M\in [0, M_0]$. 
\end{lemma}
\begin{proof}
    1. Observe that $\partial\tilde{\mc{W}}_\alpha^h$ can be isometrically identified by $\partial\mc{W}_\alpha^h$. Denote 
    \begin{equation}
        \mc{X}=Z^{s, l, s_0}_{l_0, l^-, l^+}(\partial\mc{W}_\alpha^h), \ \mc{Y}=Z^{s-1, l, s_0-1}_{l_0, l^-+1, l^++1}(\partial\mc{W}_\alpha^h), \ \mc{K}=H^{-N,-N}(\partial\mc{W}_\alpha^h),
    \end{equation}
    which then identify their counterparts on $\partial\tilde{\mc{W}}_\alpha^h$. Note $\mc{X}\hookrightarrow \mc{K}$ compactly. Assume against contradiction that there are $u_k$, $\mu_k\to 0^+$ with $\|u_k\|_{\mc{X}}=1$, $h_k\to 0$ , such that $\|(\Lambda^{\tilde{\mc{W}}_\alpha^{h_k}}_{h_kM}-\omega-i\mu_k)u_{k}\|_{\mc{Y}}\rightarrow 0$. From Proposition~\ref{p:FredEst}, we have 
    \begin{equation}
        \|u_k\|_{\mc{X}}\le C\|(\Lambda^{\tilde{\mc{W}}_\alpha^{h_k}}_{h_kM}-\omega-i\mu_k)u_k\|_{\mc{Y}}+C_N\|u_k\|_{\mc{K}},
    \end{equation}
    since $\{\tilde{\mc{W}}_\alpha^{h_k}\}$ have uniformly controlled constants needed for Proposition~\ref{p:FredEst}. Observe that all estimates we used are uniform in $\omega,M$ in the prescribed range, and only the limit $\mu_k\rightarrow 0^+$ is interesting here. Since $u_k$ are bounded in $\mc{X}$, $u_k\rightharpoonup v$ weakly in $\mc{X}$ and $u_k\to v$ for some $v\in \mc{X}$ in norm in $\mc{K}$ via a subsequence.
    
    2. We claim $(\Lambda_0^{\mc{W}_\alpha^h}-\omega)v=0$. Indeed, for a test function $\phi\in C_c^\infty(\partial \mc{W}_\alpha^h)$, we have
\begin{multline}
\langle (\Lambda_0^{\mc{W}_\alpha^h} -\omega)v,\phi\rangle 
=\lim_{k\to \infty}\langle u_k,(\Lambda_0^{\mc{W}_\alpha^h} -\omega)\phi\rangle
= \lim_{k\to \infty}\langle u_k,(\Lambda_0^{\mc{W}_\alpha^h}-\Lambda_{h_kM}^{\tilde{\mc{W}}^{h_k}_\alpha}) \phi\rangle +\langle (\Lambda_{h_kM}^{\tilde{\mc{W}}^{h_k}_\alpha}-\omega)u_k,\phi\rangle \\
=\lim_{k\to \infty}\langle u_k,(\Lambda_0^{\mc{W}_\alpha^{h_k}}-\Lambda_{h_kM}^{\tilde{\mc{W}}^{h_k}_\alpha})\phi\rangle.
\end{multline}
Therefore, since $u_k$ is bounded in $\mc{X}$, it is enough to show that $\Lambda_{hM}^{\tilde{\mc{W}}^{h}_\alpha}\phi\to \Lambda_0^{\mc{W}_\alpha^h} \phi$ in $\mc{X}^*$ as $h\rightarrow 0$. Now noting $\phi=\chi^\delta \phi$ for small $h$ and $\tilde\chi^\delta F^*_{h} \mathcal{H}_{0}^{\mc{W}_\alpha^{h}} \phi$ extends $\phi$ from $\partial\tilde{\mc{W}}^{h}_\alpha$ to $\tilde{\mc{W}}^{h}_\alpha$, we have
\begin{equation}
    \mc{H}_{h M}^{\tilde{\mc{W}}^{h}_\alpha}\phi=\tilde\chi^\delta F^*_{h} \mc{H}_{0}^{\mc{W}_\alpha^{h}} \phi-\mc{R}_{hM}^{\tilde{\mc{W}}^{h}_\alpha}(-\Delta+h^2M^2)\tilde\chi^\delta F^*_{h} \mathcal{H}_{0}^{\mc{W}_\alpha^{h}} \phi,
\end{equation}
and thus by $dF_h(\partial_\nu)=\partial_\nu$ on $\supp\tilde\chi^\delta$, we have
\begin{equation}
    \Lambda_{h M}^{\tilde{\mc{W}}^{h}_\alpha}\phi=\tilde\chi^\delta \Lambda_{0}^{\mc{W}_\alpha^{h}} \phi+\partial_\nu\mc{R}_{hM}^{\tilde{\mc{W}}^{h}_\alpha}(-\Delta+h^2M^2)\tilde\chi^\delta F^*_{h} \mathcal{H}_{0}^{\mc{W}_\alpha^{h}} \phi.
\end{equation}
One can check $[-\Delta, \tilde\chi^\delta]=O(h^\delta)\nabla+O(h^{2\delta})$ and
\begin{equation}\label{e:8l1}
    (F_h^{-1})^*(-\Delta)(F_h)^*=-\Delta+L, \ L=O(h^{k_1}|z|^{k_1})\nabla^2+O(h^{k_1}|z|^{{k_1}-1})\nabla,
\end{equation}
which implies
\begin{equation}
    \partial_\nu\mc{R}_{hM}^{\tilde{\mc{W}}^{h}_\alpha}(-\Delta+h^2M^2)\tilde\chi^\delta F^*_{h} \mathcal{H}_{0}^{\mc{W}_\alpha^{h}} \phi=\partial_\nu\mc{R}_{hM}^{\tilde{\mc{W}}^{h}_\alpha}F^*_{h}([-\Delta,\tilde\chi^\delta]+(h^2M^2+L)\tilde\chi^\delta)\mathcal{H}_{0}^{\mc{W}_\alpha^{h}} \phi=O_{\mc{X}^*}(h^\delta) 
\end{equation}
decays, since $|z|\le Ch^{-\delta}$ on $\supp\tilde\chi^\delta$ in $\mc{W}^h_\alpha$. Now we have
\begin{equation}
    (\Lambda_{hM}^{\tilde{\mc{W}}_\alpha^h}-\Lambda_{0}^{\mc{W}_\alpha^h})\phi=-(1-\tilde\chi^\delta)\Lambda_{0}^{\mc{W}_\alpha^h}\chi^\delta \phi+O_{\mc{X}^*}(h^\delta)=O_{\mc{X}^*}(h^\delta),
\end{equation}
goes to zero as desired, where we used \eqref{e:genCommutator} to obtain
\begin{equation}
    (1-\tilde\chi^\delta)\Lambda_{0}^{\mc{W}_\alpha^h}\chi^\delta \phi=\partial_\nu \mc{R}_{0}^{\mc{W}_\alpha^h}[-\Delta, 1-\tilde\chi^\delta]\mc{H}_{0}^{\mc{W}_\alpha^h}\chi^\delta \phi=O_{\mc{X}^*}(h^\delta),
\end{equation}
since $[-\Delta, 1-\tilde\chi^\delta]=O(h^\delta)\nabla+O(h^{2\delta})$. This concludes the claim. 

3. Now set $\mc{W}_\alpha=W_\alpha$, which is scaling-invariant. Then $(\Lambda_0^{W_\alpha^h}-\omega)v=0$ and $v\in \mc{X}$. The Fredholm estimate, Proposition~\ref{p:FredEst} then implies $v\in Z^{\infty,\infty,1}_{L_0-0, \infty, -\frac12-0}$, and thus $\WF v\subset\Sigma^+(\omega)\cup \{\xi=0\}$. Lemma~\ref{l:EigenReg} gives $\WF v\subset\{\xi=0\}$ and thus $v\in L^2$, and $\omega$ is an eigenvalue to $\Lambda^{W_\alpha^h}_0$. Note $W^h_\alpha$ is a straight wedge and is invariant under rescaling. For each $\lambda\in\mathbb{R}_{>0}$, we have
\begin{equation}
    (\Lambda_{0}^{W_\alpha^h}-\lambda\omega)v(\lambda x)=\lambda((\Lambda_{0}^{W_\alpha^h}-\omega)v)(\lambda x)=0.
\end{equation}
This implies $(0,\infty)\subset\operatorname{Spec}_{\mathrm{pp}}(\Lambda_{0}^{W_\alpha^h})$, contradicting with its discreteness per Lemma~\ref{l:FiniteEigenVal}.
\end{proof}

\subsection{Comparison between curved wedges}
In this subsection, we will show the scattering matrices of $\tilde{\mc{W}}_j$ and $\mc{W}_{\alpha_j}$ are close to each other.
\begin{lemma}\label{l:scstraighten}
Let $M_0>0$. Then there exists $h_0>0$ that uniformly for all $M\in[0, M_0]$, $\omega\in[\frac12, 2]$, $h\in(0, h_0]$, we have 
\begin{equation}\label{e:7l32}
    \Sigma^{\tilde{\mc{W}}_{\alpha}}_{M, h^{-1}\omega}=\Sigma^{\mc{W}_{\alpha}}_{M,h^{-1} \omega}+O(h^{\min(2\rardya,1+k_0)-0})+O(M^2h^{2-0}),
\end{equation}
given $k_0\ge 0$ that the first $(k_0+1)$ derivatives of the boundaries of $\mc{W}_\alpha,\tilde{\mc{W}}_\alpha$ agree at the corner.
\end{lemma}
\begin{proof}
1. We begin with some generalized eigenfunctions $u$ on $\partial\mc{W}^h_\alpha$. Let $a^i=(a^i_+, a^i_-)\in \mathbb{C}^2$ with $|a^i|=1$. Proposition~\ref{p:existscmatrix} implies the existence of $u\in Z^{-\infty,-\infty, 1-0}_{L_0-0, -\frac12-0, -\frac12-0}(\partial\mc{W}_\alpha^h; \omega)$ such that $(\Lambda_{hM}^{\mc{W}_{\alpha}^h}-\ohm)u=0$, $u\sim \tilde a^i_\mp(x)e^{-i\omega|x|}+\tilde a^o_\pm(x)e^{i\omega|x|}$ away from the zero frequency as $x\rightarrow \pm\infty$, where
    \begin{equation}
        \tilde a_\pm^i(x;h)=a_\pm^i+O_{S^0}(\langle x\rangle^{-1}), \ \tilde a_{\pm}^o(x;h)=a_\pm^o(h)+O(\langle x\rangle^{-1}), \ a^o=(a^o_+, a^o_-)=\Sigma^{\mc{W}_\alpha}_{M, h^{-1}\omega} a^i=\Sigma^{\mc{W}_\alpha^h}_{hM, \omega} a^i.
    \end{equation}
Denote its harmonic extension
\begin{equation}
    \bar{u}=\mc{H}^{\mc{W}_{\alpha}^h}_{hM} u\in H^{1}_{\loc}(\mc{W}_\alpha^h).
\end{equation}
In Step 1, we do some a priori analysis of $\bar u$. 

1a. Lemma~\ref{l:2dquasimodeappendix} gives the existence of quasimodes
$\tilde e_{\pm,L/R,h}(x_0, x;h, \omega, M)$ on $\tilde{\mc{W}}_{L/R}$ that 
\begin{equation}
    (-\Delta+M^2)\tilde e_{\pm,L/R,h}=O_{\mc{S}(\tilde{\mc{W}}_{L/R})}(h^\infty), \ (h\partial_\nu-\om)\tilde e_{\pm,L/R,h}=O_{\mc{S}(\partial\tilde{\mc{W}}_{L/R})}(h^\infty),
\end{equation}
whose boundary restriction $e_{\pm,L/R,h}(x;h, \omega, M)=\tilde e_{\pm,L/R,h}(0, x;h, \omega, M)$ satisfies
\begin{equation}
    e_{\pm,L/R,h}(x;h,\omega, M)=e^{\pm i\omega x/h}(1+O_{S^0(\partial\tilde{\mc{W}}_{L/R})}(h^2M^2)).
\end{equation}
Rescale by $\tilde e_{\pm,L/R}=h^{-\frac12}T_h \tilde e_{\pm,L/R,h}$ to obtain
\begin{gather}\label{e:7l72}
    (-\Delta+h^2M^2)\tilde e_{\pm,L/R}=O_{\mc{S}(\tilde{\mc{W}}_{L/R}^h)}(h^\infty), \ (\partial_\nu-\ohm)\tilde e_{\pm,L/R}=O_{\mc{S}(\partial\tilde{\mc{W}}_{L/R}^h)}(h^\infty),\\
    e_{\pm,L/R}(x;h,\omega, M)=e^{\pm i\omega x}(1+O_{S^0(\partial\tilde{\mc{W}}_{L/R}^h)}(h^2M^2)).
\end{gather}

1b. Parametrize $\partial\mc{W}^h_\alpha$ by $x$ where $\partial\mc{W}^h_\alpha$ is singular at 0, and let $(x_0,x)$ be the normal coordinate in the support of $\chi^h_{L/R}\in S^0(\mc{W}_\alpha^h)$ and $x_0=0$ on $\partial \mc{W}_{\alpha}^h$. We craft the principal oscillator by
\begin{equation}\label{e:7l71}
    \tilde u_{\omega}(x_0, x)=a_{+}^{o}\chi_R^h\tilde e_{+, R}+a_{+}^{i}\chi_L^h\tilde e_{+, L}+a_{-}^{i}\chi_R^h\tilde e_{-, R}+a_{-}^{o}\chi_L^h\tilde e_{-, L}.
\end{equation}

Let the zero-frequency residue be $\tilde u_0=\bar u-\tilde{u}_\omega\in X_m^l(\mc{W}_{\alpha}^h)$ shown in Step 1c below. In this proof, we use bar to denote harmonic extension like $\bar u$, and tilde to denote non-harmonic extensions like $\tilde u_\omega$. Let $u_0=\tilde u_0|_{\partial\mc{W}_\alpha^h}$, $u_\omega=\tilde u_\omega|_{\partial\mc{W}_\alpha^h}$. Lemma~\ref{l:EllEst_Near0_breg} implies in $\partial\mc{W}_\alpha^h$,
\begin{equation}\label{e:8l4}
    \|\langle x\rangle^{\rardya+\frac12-0}(\langle x\rangle\partial_x)^k (1-\chi^h) u_0(x)\|_{H^N}\le C, \ |\partial^k (1-\chi^h) u_{0}(x)|\le C\langle x\rangle^{-\rardya-1+0-k}. 
\end{equation}
Heuristically, $u_\omega$ is a mode oscillating at frequencies $\pm \omega$ and at distance at least $1$ away from the corner in the semiclassical scale, and $u_0$ contains modes oscillating at frequency $0$ and at distance at least $1$ away from the corner and at all frequencies near the corner. 

1c. We claim $\tilde u_0\in X_{\rardya-0}^{\rardya+\frac12-0,N}(\mc{W}_\alpha^h)$. We take note of the key relation
\begin{multline}\label{e:7l75}
    (-\Delta+h^2M^2)\tilde u_0=-(-\Delta+h^2M^2)\tilde u_\omega=a_{+}^{o}[\Delta,\chi_R^h]\tilde e_{+, R}+a_{+}^{i}[\Delta,\chi_L^h]\tilde e_{+, L}+a_{-}^{i}[\Delta,\chi_R^h]\tilde e_{-, R}\\
    +a_{-}^{o}[\Delta,\chi_L^h]\tilde e_{-, L}+O_{\mc{S}}(h^\infty),
\end{multline}
since it is supported inside $\{C^{-1}<|x|,x_0>C^{-1}\langle x\rangle-C\}$ in $\mc{W}_{\alpha}^h$. That $\tilde e_{\pm,L/R}$ are exponentially decaying in $x_0$ then implies $(-\Delta+h^2M^2)\tilde u_0\in X^{N,N}_N(\mc{W}_{\alpha}^h)\subset Y^{N,N}_{N-2}(\mc{W}_{\alpha}^h)$ for any $N>0$ large. Note
\begin{equation}
    \tilde u_0=\mc{H}_{hM}^{\mc{W}_\alpha^h} u_0+\mc{R}_{hM}^{\mc{W}_\alpha^h}(-\Delta+h^2M^2)\tilde u_0=\mc{H}_{hM}^{\mc{W}_\alpha^h} u_0+O(X^{\rardya+\frac12-0,N}_{\rardya-0})\in X_{\rardya-0}^{\rardya+\frac12-0,N}(\mc{W}_\alpha^h),
\end{equation}
by applying Lemma~\ref{l:longDistanceHarmonic} per \eqref{e:8l4}. 

2. Our goal in Step 2 and 3, is to use the scattering data of generalized eigenfunctions $u$ on $\partial\mc{W}_{\alpha}^{h}$, to construct an $h$-family quasimodes $v'$ on $\partial\tilde{\mc{W}}_{\alpha}^{h}$. Let $\delta\in(0,1)$ small be a scale to be chosen later. Consider $F_h^*\chi^\delta w$ as a function localized in an $h^{-\delta}$-neighbourhood of 0 in $h^{-1}U\subset \tilde{\mc{W}}_{\alpha}^{h}$. We confuse the notation by $\chi^\delta F_h^* u=F_h^*\chi^\delta u$ via the pullback $F_h^*$. We want to extend $\chi^\delta F_h^*u$ by quasimodes $e_{\pm,L/R}'$ on $\tilde{\mc{W}}_{L/R}^h$:
\begin{gather}
    w_{R}=a_{+}^o(h)e_{+, R}'+a_{-}^ie_{-, R}'=a_{+}^o(h)\big(1+M^2h^2O(S^0)\big)e^{i \omega x}+a_{-}^i\big(1+M^2h^2O(S^0)\big)e^{-i \omega x},\\
    w_{L}=a_{+}^ie_{+, L}'+a_{-}^o(h)e_{-, L}'=a_{+}^i\big(1+M^2h^2O(S^0)\big)e^{i \omega x}+a_{-}^o(h)\big(1+M^2h^2O(S^0)\big)e^{-i \omega x}.
\end{gather}
Denote their harmonic extensions by $\bar w_{L/R}=\mc{H}^{\tilde{\mc{W}}_{L/R}^h}_{hM} w_{L/R}\in C^\infty(\tilde{\mc{W}}_{L/R}^h)$.

We then construct the quasimodes on $\partial\tilde{\mc{W}}_\alpha^h$, by
\begin{equation}\label{e:7l28}
    v'=\chi^\delta F^*_h u+(\chi-\chi^\delta) F_h^* u_0+(1-\chi^\delta)(\chi^\delta_Rw_R+\chi^\delta_Lw_L).
\end{equation}
The rest of Step 2 and Step 3 is devoted to showing that for each $\delta\ll 1$, we have
\begin{equation}\label{e:7l45}
(\Lambda^{\tilde{\mc{W}}^h_{\alpha}}_{hM}-\ohm)v'=O_{Z^{-\frac12, N, -\frac12}_{-L_0+0, N, N}}(h^{\min(2\rardya,k_1-K\delta)-0}+M^2h^{2-K\delta})
\end{equation}
for some $K>0$ only depending on $k_1, N,\rardya$, not depending on $\delta$.

2a. We begin with some a priori analysis of \eqref{e:7l28}. First note 
\begin{gather}
\label{e:7l30}
    (-\Delta +h^2M^2)F_h^*\bar{u}=F_h^*(-\Delta +h^2M^2+L)\bar{u}=F_h^*L \bar{u},\\
\label{e:7l31}
    \partial_\nu F_h^*\bar u=F_h^*DF_h(\partial_\nu) \bar u=F_h^*\partial_\nu \bar u=\ohm F_h^* \bar u,
\end{gather}
where $L$ was defined back in \eqref{e:8l1}. We will use them later. The first two terms of \eqref{e:7l28} give
\begin{multline}
    (\Lambda^{\tilde{\mc{W}}^h_{\alpha}}_{hM}-\ohm)(\chi^\delta F^*_h u+(\chi-\chi^\delta) F_h^* u_0)=(\Lambda^{\tilde{\mc{W}}^h_{\alpha}}_{hM}-\ohm)(\chi^\delta F^*_h u_\omega+\chi F_h^* u_0)\\
    =(\partial_\nu-\ohm)(\chi^\delta F^*_h \tilde u_\omega+\chi F_h^* \tilde u_0)-\partial_\nu\mc{R}^{\tilde{\mc{W}}^h_{\alpha}}_{hM}(-\Delta+h^2M^2)(\chi^\delta F^*_h \tilde u_\omega+\chi F_h^* \tilde u_0)\\
    =(\partial_\nu-\ohm)(\chi^\delta F^*_h \tilde u_\omega+\chi F_h^* \tilde u_0)-\partial_\nu\mc{R}^{\tilde{\mc{W}}^h_{\alpha}}_{hM}\chi^\delta(-\Delta+h^2M^2)F^*_h \tilde u_\omega\\-\partial_\nu\mc{R}^{\tilde{\mc{W}}^h_{\alpha}}_{hM}\chi(-\Delta+h^2M^2)F^*_h \tilde u_0-\partial_\nu\mc{R}^{\tilde{\mc{W}}^h_{\alpha}}_{hM}([-\Delta,\chi^\delta]F^*_h \tilde u_\omega+[-\Delta,\chi]F^*_h \tilde u_0),
\end{multline}
where the first term, via \eqref{e:7l31}, can be further written as
\begin{multline}
    (\partial_\nu-\ohm)(\chi^\delta F^*_h \tilde u_\omega+\chi F_h^* \tilde u_0)\\
    =\chi F^*_h (\partial_\nu-\ohm)\bar u-(\chi-\chi^\delta) F_h^* (\partial_\nu-\ohm) \tilde u_\omega=O(h^\infty),
\end{multline}
from the quasimode nature of $\tilde u_\omega$ per \eqref{e:7l71}, $\eqref{e:7l72}$. On the other hand, $(1-\chi^\delta)\chi_{L/R}^\delta\bar w_{L/R}$ is an extension of $(1-\chi^\delta)\chi_{L/R}^\delta w_{L/R}$, and the last term of \eqref{e:7l28} give
\begin{multline}
    (\Lambda^{\tilde{\mc{W}}^h_{\alpha}}_{hM}-\ohm)\sum_{L/R}(1-\chi^\delta)\chi^\delta_{L/R} w_{L/R}=\sum_{L/R}\Big((\partial_\nu-\ohm)(1-\chi^\delta)\chi_{L/R}^\delta\bar w_{L/R}\\
    -\partial_\nu \mathcal{R}^{\tilde{\mc{W}}^h_{\alpha}}_{hM}[-\Delta, (1-\chi^\delta)\chi_{L/R}^\delta]\bar w_{L/R}\Big)
    =\sum_{L/R} \partial_\nu \mathcal{R}^{{\tilde{\mc{W}}^h_{\alpha}}}_{hM}[-\Delta, \chi^\delta]\chi_{L/R}^{\delta}\bar w_{L/R}\\
    -\sum_{L/R} \partial_\nu \mathcal{R}^{{\tilde{\mc{W}}^h_{\alpha}}}_{hM}(1-\chi^\delta)[-\Delta, \chi_{L/R}^{\delta}]\bar w_{L/R}
    +O(h^\infty \mc{S}).
\end{multline}
Thus together we have
\begin{multline}\label{e:7l29}
    (\Lambda^{\tilde{\mc{W}}^h_{\alpha}}_{hM}-\ohm)v'=-\partial_\nu\mc{R}^{\tilde{\mc{W}}^h_{\alpha}}_{hM}(\chi-\chi^\delta)(-\Delta+h^2M^2)F^*_h \tilde u_0-\partial_\nu\mc{R}^{\tilde{\mc{W}}^h_{\alpha}}_{hM}\chi^\delta(-\Delta+h^2M^2)F^*_h \bar u\\-\partial_\nu\mc{R}^{\tilde{\mc{W}}^h_{\alpha}}_{hM}[-\Delta,\chi]F^*_h \tilde u_0-\sum_{L/R} \partial_\nu \mathcal{R}^{{\tilde{\mc{W}}^h_{\alpha}}}_{hM}(1-\chi^\delta)[-\Delta, \chi_{L/R}^{\delta}]\bar w_{L/R}
    \\
    +\partial_\nu\mc{R}^{\tilde{\mc{W}}^h_{\alpha}}_{hM}[-\Delta,\chi^{\delta}]\Big(-F^*_h \tilde u_\omega+\chi_{L}^\delta\bar w_{L}+\chi_{R}^\delta\Big)+O(h^\infty \mc{S}).
\end{multline}
We will estimate the above term by term, except $O(h^\infty \mc{S})$, in the rest of Step 2 and Step 3.

2b. We estimate the first term in \eqref{e:7l29} by claiming
\begin{equation}
\|\partial_\nu\mc{R}^{\tilde{\mc{W}}^h_{\alpha}}_{hM}(\chi-\chi^\delta)(-\Delta+h^2M^2)F^*_h \tilde u_0\|_{Z^{-\frac12, N, -\frac12}_{-L_0+0, N, N}(\partial \tilde{\mc{W}}^h_{\alpha})}\le C h^{\min(2\rardya-0, k_1)}. 
\end{equation}
Note we follow \eqref{e:7l75} for
\begin{multline}
    (\chi-\chi^\delta)(-\Delta+h^2M^2)\tilde{u}_\omega=(\chi-\chi^\delta)\Big(a_{+}^{o}[\Delta,\chi_R^h]\tilde e_{+, R}+a_{+}^{i}[\Delta,\chi_L^h]\tilde e_{+, L}+a_{-}^{i}[\Delta,\chi_R^h]\tilde e_{-, R}
    \\
    +a_{-}^{o}[\Delta,\chi_L^h]\tilde e_{-, L}+O_{\mc{S}}(h^\infty)\Big)=O_{\mc{S}}(h^\infty),
\end{multline}
since the support of the right hand side is inside $\{C^{-1}h^{-\delta}<|x|<Ch^{-1},x_0>C^{-1}\langle x\rangle\}$ in $\mc{W}_{\alpha}^h$. This with $\tilde{u}_0=\bar u-\tilde u_\omega$ implies $(\chi-\chi^\delta)(-\Delta+h^2M^2)\tilde{u}_0=O_{\mc{S}}(h^\infty)$ as well. Then
\begin{equation}
(\chi-\chi^\delta) (-\Delta+h^2M^2)F_h^*\tilde u_{0}=(\chi-\chi^\delta) F_h^*L\tilde u_0+O_{\mc{S}}(h^\infty)
\in \chi h^{k_1}\rr^{k_1}Y_{\rardya-0}^{\rardya+\frac12-0,N}\subset h^{\min(2\rardya-0, k_1)}Y_{-\rardya+0}^{-\rardy_{\alpha}-\frac12+0,N}
\end{equation}
uniformly as desired, since $\rr^{-1}Y_{\rardya-1-0}^{\rardya+\frac12-0,N}\subset Y_{\rardya-0}^{\rardya+\frac12-0,N}$ and 
\begin{equation}
    h^{k_1}\rr^{k_1}Y_{\rardya-0}^{\rardya+\frac12-0,N}\subset h^{k_1}\rr^{{k_1}-(2\rardya-0)}Y_{-\rardya+0}^{-\rardy_{\alpha}-\frac12+0,N}\subset h^{\min(2\rardya-0, {k_1})}Y_{-\rardya+0}^{-\rardy_{\alpha}-\frac12+0,N},
\end{equation}
Apply $\partial_\nu \mc{R}^{\tilde{\mc{W}}^h_{\alpha}}_{hM}$ on the left to see the estimates we need for the zero frequency, 
\begin{equation}
\|\partial_\nu\mc{R}^{\tilde{\mc{W}}^h_{\alpha}}_{hM}(\chi-\chi^\delta)(-\Delta+h^2M^2)F^*_h \tilde u_0\|_{H^{-\frac12, -\rardy_{\alpha}-\frac12+0}(\partial \tilde{\mc{W}}^h_{\alpha})}\le C h^{\min(2\rardya-0, {k_1})}. 
\end{equation}
Let $\gamma^h, \tilde\gamma^h\in S^0(\partial\tilde{\mc{W}}^h_{\alpha})$ supported in the smooth part and $\supp\gamma^h\cap\supp(1-\tilde\gamma^h)=\emptyset$. Note \eqref{e:thing4} gives $\tilde\gamma^h D^N\gamma^h\partial_\nu\mc{R}^{\tilde{\mc{W}}^h_{\alpha}}_{hM}: \dot Y^{-\rardya-\frac12+0, N}_{-\rardya+0}\rightarrow H^{-\frac12, -\rardya-\frac12+N+0}$ bounded, and
\begin{equation}
    \tilde\gamma^h (1-\psi_0)(D)\gamma^h=\tilde\gamma^h (1-\psi_0)(D)D^{-N}\tilde\gamma^h D^N\gamma^h.
\end{equation}
We then have
\begin{multline}\label{e:8l6}
    \|\tilde\gamma^h (1-\psi_0)(D)\gamma^h\partial_\nu\mc{R}^{\tilde{\mc{W}}^h_{\alpha}}_{hM}(\chi-\chi^\delta)(-\Delta+h^2M^2)F^*_h \tilde u_0\|_{H^{-\frac12, -\rardya-\frac12+N-0}}\\\le C\|(\chi-\chi^\delta)(-\Delta+h^2M^2)F^*_h \tilde u_0\|_{Y^{-\rardya-\frac12+0, N}_{-\rardya+0}}
    \le Ch^{\min(2\rardya-0, {k_1})}.
\end{multline}
Picking large $N>0$ gives desired estimates away from zero frequencies and thus the claim. 

2c. We estimate the second term in \eqref{e:7l29} by claiming
\begin{equation}
    \|\partial_\nu\mc{R}^{\tilde{\mc{W}}^h_{\alpha}}_{hM}\chi^\delta(-\Delta+h^2M^2)F^*_h \bar u\|_{Z^{-\frac12, N, -\frac12}_{-L_0+0, N, N}(\partial \tilde{\mc{W}}^h_{\alpha})}\le C h^{{k_1}-K\delta-0}.
\end{equation}
Note $\rr\le Ch^{-\delta}$ over the support of $\chi^\delta $. With $(-\Delta+h^2M^2)\bar u=0$, we have
\begin{equation}
    \chi^\delta(-\Delta+h^2M^2)F^*_h \bar u=F^*_h\chi^\delta L \bar u=F^*_h\chi^\delta L \tilde u_0+F^*_h\chi^\delta L \tilde u_\omega,
\end{equation}
where arguments similar to Step 2b give
\begin{equation}
    F^*_h\chi^\delta L \tilde u_0\in\chi^\delta h^{k_1}\rr^{k_1}Y_{\rardya-0}^{\rardya+\frac12-0}\subset \chi^\delta h^{k_1}\rr^{{k_1}-(2\rardya-0)}Y_{-\rardya+0}^{-\rardy_{\alpha}-\frac12+0}\subset h^{{k_1}-\delta\max({k_1}-(2\rardya-0),0)}Y_{-\rardya+0}^{-\rardy_{\alpha}-\frac12+0}(\tilde{\mc{W}}_\alpha^h).
\end{equation}
For the second term, we have
\begin{multline}
    F^*_h\chi^\delta L \tilde u_\omega\in \chi^\delta h^{k_1}\rr^{k_1} H^{-1,-\frac12-0}\subset \chi^\delta h^{k_1}\rr^{{k_1}-\rardya+0} \dd^{\frac32+0}(\rr^{\rardya+\frac12-0}\dd^{-\frac32-0}H^{-1})\\
    \subset h^{{k_1}-\delta\max({k_1}-\rardya+\frac32+0,\frac32+0)}Y_{-\rardya+0}^{-\rardy_{\alpha}-\frac12+0}(\tilde{\mc{W}}_\alpha^h)
\end{multline}
uniformly. Together we get 
\begin{equation}
    \|\partial_\nu\mc{R}^{\tilde{\mc{W}}^h_{\alpha}}_{hM}\chi^\delta(-\Delta+h^2M^2)F^*_h \bar u\|_{H^{-\frac12, -\rardy_{\alpha}-\frac12+0}(\partial \tilde{\mc{W}}^h_{\alpha})}\le C h^{k_1-\delta\max({k_1}-\rardya+\frac32,\frac32)-0} 
\end{equation}
for the zero frequency. For the nonzero frequencies, note \eqref{e:thing1} gives $\gamma^h (1-\psi_0)(D)\gamma^h\partial_\nu\mc{R}^{\tilde{\mc{W}}^h_{\alpha}}_{hM}: \dot Y^{-\rardya-\frac12+N+0}_{-\rardya+0}\rightarrow H^{-\frac12, -\rardya-\frac12+N+0}$ bounded and thus
\begin{multline}\label{e:8l7}
    \|\gamma^h (1-\psi_0)(D)\gamma^h\partial_\nu\mc{R}^{\tilde{\mc{W}}^h_{\alpha}}_{hM}\chi^\delta(-\Delta+h^2M^2)F^*_h \bar u\|_{H^{-\frac12, -\rardya-\frac12+N-0}}\\\le C\|\chi^\delta(-\Delta+h^2M^2)F^*_h \bar u\|_{Y^{-\rardya-\frac12+N+0}_{-\rardya+0}}
    \le C h^{k_1-\delta\max(N+{k_1}-\rardya+\frac32,N+\frac32)-0},
\end{multline}
similar to \eqref{e:8l6}. This completes the claim.

2d. On the other hand, we have
\begin{equation}
    [-\Delta,\chi]F_h^*\tilde u_{0}=-(\Delta\chi)F_h^*\tilde u_{0}-2(\nabla\chi)F_h^*(\nabla+O(h^{k_1}\rr^{k_1})\nabla)\tilde u_{0}.
\end{equation}
Noting $\rr\sim h^{-1}$ on $\supp\nabla \chi$, we have $|\Delta\chi|\le Ch^{2}$, $|\nabla\chi|\le Ch$ and
\begin{multline}
[-\Delta,\chi]F_h^*\tilde u_{0}\in h^2 Y_{\rardya-2-0}^{\rardya+\frac12-0,N}+hY_{\rardya-1-0}^{\rardya+\frac12-0,N} \subset (h^2\rr^{-(2\rardya-2-0)}+h \rr^{-(2\rardya-1-0)}) Y_{-\rardya+0}^{-\rardy_{\alpha}-\frac12+0,N}\\
\subset h^{2\rardya-0}Y_{-\rardya+0}^{-\rardy_{\alpha}-\frac12+0,N}(\tilde{\mc{W}}_\alpha^h),
\end{multline}
allowing us to estimate the third term in \eqref{e:7l29} by
\begin{equation}
    \|\partial_\nu\mc{R}^{\tilde{\mc{W}}^h_{\alpha}}_{hM}[-\Delta,\chi]F^*_h \tilde u_0\|_{H^{-\frac12, -\rardy_{\alpha}-\frac12+0}(\partial \tilde{\mc{W}}^h_{\alpha})}\le Ch^{2\rardya-0},
\end{equation}
and thus by estimating the nonzero frequencies similarly to \eqref{e:8l6} we have
\begin{equation}
    \|\partial_\nu\mc{R}^{\tilde{\mc{W}}^h_{\alpha}}_{hM}[-\Delta,\chi]F^*_h \tilde u_0\|_{Z^{-\frac12, N, -\frac12}_{-L_0+0, N, N}(\partial \tilde{\mc{W}}^h_{\alpha})}\le Ch^{2\rardya-0}.
\end{equation}

2e. For the fourth term in \eqref{e:7l29}, note $\supp(1-\chi^\delta)[-\Delta, \chi^\delta_{L/R}]$ is disjoint from $\partial\tilde{\mc{W}}^\delta_{L/R}$ and is in $\{|x|\ge C^{-1} h^{-\delta}\}$ for some $C>0$, $w_{L/R} =\psi_\omega(D) w_{L/R}+O(\mc{S})$. Then \eqref{e:thing8} implies
\begin{equation}
    \|(1-\chi^\delta)[-\Delta, \chi^\delta_{L/R}]\mc{H}^{\tilde{\mc{W}}_{L/R}^h}_{hM} w_{L/R}\|_{H^{N,N}}\le Ch^N.
\end{equation}
Repeat Step 2c to conclude
\begin{equation}
    \|\partial_\nu \mathcal{R}^{{\tilde{\mc{W}}^h_{\alpha}}}_{hM}(1-\chi^\delta)[-\Delta, \chi_{L/R}^{\delta}]\bar w_{L/R}\|_{Z^{-\frac12, N, -\frac12}_{-L_0+0, N, N}(\partial \tilde{\mc{W}}^h_{\alpha})}\le C h^{N}.
\end{equation}

3. It remains to estimate the last term in \eqref{e:7l29}. We claim
\begin{equation}\label{e:8l8}
    \partial_\nu\mc{R}^{\tilde{\mc{W}}^h_{\alpha}}_{hM}[-\Delta,\chi^\delta]\Big(F^*_h \tilde u_\omega-\chi_{L}^\delta\bar w_{L}-\chi_{R}^\delta\bar w_{R}\Big)=O_{Z^{-\frac12, N, -\frac12}_{0,N,N}}(h^{{k_1}-K\delta}+M^2h^{2-K\delta}).
\end{equation}
We want to estimate instead on $\tilde\Omega^\delta$, the pullback of $F^*_h \tilde u_\omega-\chi_{L}^\delta\bar w_{L}-\chi_{R}^\delta\bar w_{R}$, that is,
\begin{equation}
    f=\tilde\varphi^\delta(F_{\delta h}^* \tilde u_\omega-\chi_{L}^\delta\bar w_{L,\delta}-\chi_{R}^\delta\bar w_{R,\delta}), \ F_\delta(z)=h^{-1+\delta}F(h^{1-\delta}z), \ F_{\delta h}(z)=h^{-\delta} F_\delta(z)=F_h(h^{-\delta} z).
\end{equation}
Here we rescaled $w_{L/R}$ by $w_{L/R,\delta}=F_{\delta h}^*w_{L/R}$ and $\tilde e_{+, R,\delta}'=F_{\delta h}^*\tilde e_{+, R}'$, $\tilde e_{+, R,\delta}=F_{\delta h}^*\tilde e_{+, R}$:
\begin{gather}
    w_{R,\delta}=a_{+}^o(h)e_{+, R,\delta}'+a_{-}^ie_{-, R,\delta}'=a_{+}^o(h)\big(1+M^2h^2O(S^0)\big)e^{i h^{-\delta}\omega x}+a_{-}^i\big(1+M^2h^2O(S^0)\big)e^{-i h^{-\delta}\omega x},\\
    w_{L,\delta}=a_{+}^ie_{+, L,\delta}'+a_{-}^o(h)e_{-, L,\delta}'=a_{+}^i\big(1+M^2h^2O(S^0)\big)e^{i h^{-\delta}\omega x}+a_{-}^o(h)\big(1+M^2h^2O(S^0)\big)e^{-i h^{-\delta}\omega x}.
\end{gather}
We inform the reader this process is similar to addressing Lemma~\ref{l:blowup} in the next section. 
Then $f$ extends $\tilde\varphi^\delta (F^*_{\delta h} u_\omega-\chi_{L}^\delta w_{L,\delta}-\chi_{R}^\delta w_{R,\delta})$ from $\partial\tilde\Omega^\delta$ to $\tilde\Omega^\delta$, and thus
\begin{equation}\label{e:7l50}
    \varphi^\delta f=\varphi^\delta \mathcal{H}_{h^{1-\delta}M}^{\tilde\Omega^\delta}\tilde\varphi^\delta (F^*_{\delta h} u_\omega-\chi_{L}^\delta w_{L,\delta}-\chi_{R}^\delta w_{R,\delta})+\varphi^\delta \mathcal{R}_{h^{1-\delta}M}^{\tilde\Omega^\delta}(-\Delta+h^{2-2\delta}M^2)\tilde\varphi^\delta(F^*_{\delta h} \tilde u_\omega-\chi_{L}^\delta\bar w_{L,\delta}-\chi_{R}^\delta\bar w_{R,\delta}).
\end{equation}
We will estimate it term by term. First note on $\supp\tilde\varphi^\delta\chi^\delta_{L/R}$,
\begin{equation}
    F_{\delta h}(x_0, x)=(h^{-\delta}x_0,h^{-\delta}x)+h^{{k_1}-({k_1}+1)\delta}O(|(x_0,x)|^{k_1+1}),
\end{equation}
where we understand coordinate $x_0$ as the normal variable from $\Gamma_{L/R}^\delta$. Note $\tilde\varphi^\delta$ has uniform compact support in the $\delta$-scale. Hence it follows from \eqref{e:7l71} that
\begin{equation}\label{e:7l36}
    \tilde\varphi^\delta F_{\delta h}^*\tilde u_\omega(x_0,x)=\tilde\varphi^\delta \Big(a_{+}^{o}\chi_R^h\tilde e_{+, R,\delta}+a_{+}^{i}\chi_L^h\tilde e_{+, L,\delta}+a_{-}^{i}\chi_R^h\tilde e_{-, R,\delta}+a_{-}^{o}\chi_L^h\tilde e_{-, L,\delta}\Big)+O_{H^{1}}(h^{{k_1}-({k_1}+2)\delta})
\end{equation}
and over the boundary we have
\begin{multline}
    \tilde\varphi^\delta (F^*_{\delta h} u_\omega-\chi_{L}^\delta w_{L,\delta}-\chi_{R}^\delta w_{R,\delta})=\tilde\varphi^\delta O_{H^{\frac12}_\partial}(h^{k_1-(k_1+\frac32)\delta})+O_{S^0_\partial}(M^2h^{k_1+2-k_1\delta})+O(M^2h^2)\\
    =O_{H^{\frac12}_\partial}(h^{k_1-(k_1+\frac32)\delta}+M^2h^{2-\frac\delta2}), 
\end{multline}
where we used $|x|\sim 1$ on $\supp\tilde\varphi^\delta$. Thus 
\begin{equation}
    \varphi^\delta \mathcal{H}_{h^{1-\delta}M}^{\tilde \Omega^\delta}\tilde \varphi^\delta (F^*_{\delta h} u_\omega-\chi_{L}^\delta w_{L,\delta}-\chi_{R}^\delta w_{R,\delta})=O_{H^{1}}(h^{k_1-(k_1+\frac32)\delta}+M^2h^{2-\frac\delta2}).
\end{equation}
It remains to estimate 
\begin{multline}\label{e:7l35}
    \varphi^\delta \mathcal{R}_{h^{1-\delta}M}^{\tilde\Omega^\delta}(-\Delta+h^{2-2\delta}M^2)\tilde\varphi^\delta(F^*_{\delta h} \tilde u_\omega-\chi_{L}^\delta\bar w_{L,\delta}-\chi_{R}^\delta\bar w_{R,\delta})
    =\varphi^\delta \mathcal{R}_{h^{1-\delta}M}^{\tilde\Omega^\delta}\tilde\varphi^\delta(-\Delta+h^{2-2\delta}M^2)F^*_{\delta h} \tilde u_\omega \\-\sum_{L/R}\varphi^\delta \mathcal{R}_{h^{1-\delta}M}^{\tilde{\Omega}^\delta}\tilde\varphi^\delta [-\Delta,\chi_{L/R}^\delta]\bar w_{L/R,\delta}+\varphi^\delta \mathcal{R}_{h^{1-\delta}M}^{\tilde{\Omega}^\delta}[-\Delta,\tilde\varphi^\delta](F^*_{\delta h} \tilde u_\omega-\chi_{L}^\delta\bar w_{L,\delta}-\chi_{R}^\delta\bar w_{R,\delta}).
\end{multline}

3a. We show the following a priori elliptic estimates on $\tilde \Omega^\delta$ uniformly: for each $N>0$, there is $C>0$ that
\begin{gather}
    \|\kappa^\delta \mathcal{R}^{\tilde\Omega^\delta}_{h^{1-\delta}M}\tau^\delta\|_{H^{-N}(\tilde\Omega^\delta)\rightarrow \dot H^N(\tilde\Omega^\delta)}\le C,\label{e:7l54}\\
    \|\varphi^\delta \mathcal{R}^{\tilde\Omega^\delta}_{h^{1-\delta}M}\rho_{L/R}^\delta\|_{Y^{0,-N}_{0, -N}(\tilde\Omega^\delta)\rightarrow \dot X^{0,N}_{0,N}(\tilde\Omega^\delta)}\le C,\label{e:7l55}
\end{gather}
uniformly in $h>0$ small. Note although $\tilde\Omega^\delta$ depends on $h$, $\Omega^\delta=F_\delta(\tilde\Omega^\delta)$ is a fixed compact domain with smooth boundary inside $\mc{W}_\alpha$, that does not depend on $h$. 

First note $dF_\delta(z)=\Id+O(h^{k_1(1-\delta)}z^{k_1})$ and so is $dF_\delta^{-1}$. We can pullback the massive Laplacian between $\tilde\Omega^\delta$ and $\Omega^\delta$ by
\begin{gather}
    (F_\delta^{-1})^*(-\Delta+h^{2-2\delta}M^2) F_\delta^*=-\Delta+h^{2-2\delta}M^2+L^\delta, \\
    F_\delta^*(-\Delta+h^{2-2\delta}M^2) (F_\delta^{-1})^*=-\Delta+h^{2-2\delta}M^2+\tilde L^\delta, 
\end{gather}
where
\begin{equation}
    L^\delta=O(h^{{k_1}(1-\delta)}|z|^{k_1})\nabla^2+O(h^{{k_1}(1-\delta)}|z|^{k_1-1})\nabla, \ \tilde L^\delta=O(h^{{k_1}(1-\delta)}|z|^{k_1})\nabla^2+O(h^{{k_1}(1-\delta)}|z|^{k_1-1})\nabla.  
\end{equation}
Now note
\begin{equation}
    (-\Delta+h^{2-2\delta}M^2+L^\delta)\mc{R}_{h^{1-\delta}M}^{\Omega^\delta}(\Id-(\Id+L^\delta\mathcal{R}_{h^{1-\delta}M}^{\Omega^\delta})^{-1} L^\delta\mc{R}_{h^{1-\delta}M}^{\Omega^\delta})=\Id.
\end{equation}
This implies
\begin{equation}
    \mc{R}_{h^{1-\delta}M}^{\tilde\Omega^\delta}=F_\delta^*\mc{R}_{h^{1-\delta}M}^{\Omega^\delta}(\Id-(\Id+L^\delta\mathcal{R}_{h^{1-\delta}M}^{\Omega^\delta})^{-1}L^\delta\mc{R}_{h^{1-\delta}M}^{\Omega^\delta})(F^{-1}_\delta)^*
\end{equation}
maps from $H^{-1}(\tilde\Omega^\delta)$ to $\dot H^{1}(\tilde\Omega^\delta)$ with bounds uniform in $h$ small, and all $M\le M_0$.

By Lemma~\ref{l:smoothing}, we have $\varphi^\delta \mathcal{R}^{\Omega^\delta}_{h^{1-\delta}M}\rho_{L/R}^\delta: Y^{0,-N}_{0, -N}(\Omega^\delta)\rightarrow \dot X^{0,N}_{0,N}(\Omega^\delta)$ uniformly bounded with arbitrary gain in the b-regularity because of the disjoint support of $\varphi^\delta$ and $\rho_{L/R}^\delta$. Noting $L^\delta: \dot X^{0,N}_{0,N}(\Omega^\delta) \rightarrow Y^{0,N}_{0, N}(\Omega^\delta)$ is a local differential operator of size $O(h^{{k_1}(1-\delta)})$, we have for small $h$ that
\begin{equation}
    (F_\delta^{-1})^* \varphi^\delta \mathcal{R}^{\tilde\Omega^\delta}_{h^{1-\delta}M}\rho_{L/R}^\delta F_\delta^*=\varphi^\delta \Big(\mathcal{R}^{\Omega^\delta}_{h^{1-\delta}M}-\sum_{k=0}^\infty \mathcal{R}^{\Omega^\delta}_{h^{1-\delta}M}(-L^\delta\mathcal{R}^{\Omega^\delta}_{h^{1-\delta}M})^k L^\delta \mathcal{R}^{\Omega^\delta}_{h^{1-\delta}M}\Big)\rho_{L/R}^\delta
\end{equation}
is uniformly bounded from $Y^{0,-N}_{0, -N}(\Omega^\delta)\rightarrow \dot X^{0,N}_{0,N}(\Omega^\delta)$, and thus \eqref{e:7l55}. Similarly noting $\kappa^\delta$ and $\tau^\delta$ have disjoint support both away from $\partial\Omega^\delta$, we have $\kappa^\delta\mc{R}^{\Omega^\delta}_{h^{1-\delta}M}\tau^\delta:H^{-N}(\Omega^\delta)\rightarrow \dot H^{N}(\Omega^\delta)$ uniformly bounded and thus \eqref{e:7l54}.

3b. We estimate \eqref{e:7l35} in the rest of Step 3. By \eqref{e:7l36} we can estimate
\begin{equation}
    \tilde\varphi^\delta(-\Delta+h^{2-2\delta}M^2)F^*_{\delta h} \tilde u_\omega=O_{H^{-1}}(h^{{k_1}-({k_1}+2)\delta}),
\end{equation}
and the first term of \eqref{e:7l35} uniformly by
\begin{equation}
    \varphi^\delta \mathcal{R}_{h^{1-\delta}M}^{\tilde\Omega^\delta}\tilde\varphi^\delta(-\Delta+h^{2-2\delta}M^2)F^*_{\delta h} \tilde u_\omega=O_{\dot H^{1}}(h^{{k_1}-({k_1}+2)\delta}).
\end{equation}
Note  $\tilde\varphi^\delta [-\Delta,\chi_\pm^\delta]=\tilde\varphi^\delta [-\Delta,\chi_{L/R}^\delta]\kappa^\delta$ and $[-\Delta,\tilde\kappa^\delta]=\tau^\delta[-\Delta,\tilde\kappa^\delta]$ for the second term of \eqref{e:7l35}. We want to understand $\kappa^\delta \bar w_{L/R}$. Indeed, we have
\begin{equation}
    \tilde\kappa^\delta \bar w_{{L/R},\delta}=\mc{R}^{\tilde\Omega^\delta}_{h^{1-\delta}M}(-\Delta+h^{2-2\delta}M^2)\tilde\kappa^\delta \bar w_{{L/R},\delta}=\mc{R}^{\tilde\Omega^\delta}_{h^{1-\delta}M}\tau^\delta[-\Delta,\tilde\kappa^\delta]\tau^\delta \bar w_{{L/R},\delta}
\end{equation}
and thus $\kappa^\delta \bar w_{{L/R},\delta}=\kappa^\delta\mc{R}^{\tilde\Omega^\delta}_{h^{1-\delta}M}\tau^\delta[-\Delta,\tilde\kappa^\delta] \bar w_{{L/R},\delta}$. Apply \eqref{e:7l54} and \eqref{e:thing8} to see
\begin{equation}\label{e:7l53}
    \|\kappa^\delta \bar w_{L/R,\delta}\|_{H^1(\tilde\Omega^\delta)}\le C\|\tau^\delta \bar w_{L/R,\delta}\|_{H^{-2N}(\bar\Gamma^\delta_\pm)}\le C\|w_{L/R,\delta}\|_{H^{-2N+\frac12,-N}_\partial} \le C(h^{N\delta}+h^{1+{k_1}-{k_1}\delta})=O(h^{1+{k_1}-{k_1}\delta})
\end{equation}
for large $N>0$, since $\tau^\delta $ is compactly supported and we have 
\begin{equation}
    e^{\pm ih^{-\delta}\omega x }=\Big(\frac{\pm h^\delta D}{\omega}\Big)^N e^{\pm ih^{-\delta}\omega x}
\end{equation}
is away from the zero frequency. This implies we can estimate the second term of \eqref{e:7l35} by
\begin{equation}
    \|\varphi^\delta \mathcal{R}_{h^{1-\delta}M}^{\tilde{\Omega}^\delta}\tilde\varphi^\delta [-\Delta,\chi_{L/R}^\delta]\bar w_{L/R,\delta}\|_{\dot H^1(\tilde\Omega^\delta)}=O(h^{1+{k_1}-{k_1}\delta}).
\end{equation}

3c. For the second piece of the third term of \eqref{e:7l35}, note $\rho_{L/R}^\delta[-\Delta, \tilde\varphi^\delta]\chi_{L/R}^\delta=[-\Delta, \tilde\varphi^\delta]\chi_{L/R}^\delta$. Similarly to \eqref{e:7l53}, we can use Lemma~\ref{l:hnegbreg}, \eqref{e:7l55}, \eqref{e:thing8} to estimate for any $N>0$ large
\begin{equation}\label{e:7l49}
    \|\varphi^\delta \mc{R}_{h^{1-\delta}M}^{\tilde\Omega^\delta}[-\Delta, \tilde\varphi^\delta]\chi_{L/R}^\delta \bar w_{{L/R},\delta}\|_{H^1(\tilde\Omega^\delta)}\le C\|\rho_{L/R}^\delta \bar w_{{L/R},\delta}\|_{X^{-N,0}_{0,-N}(\tilde{\mc{W}}^\delta_{L/R})}\le C\|w_{{L/R},\delta}\|_{H^{-2N+\frac12,-N}_\partial}=O(h^{1+{k_1}-{k_1}\delta})
\end{equation}
for large $N>0$, since $\rho^\delta_{L/R}$ is again compactly supported.

3d. For the first piece of the last term of \eqref{e:7l35}, we want to understand $\rho_{L/R}^\delta\chi_{L/R}^\delta F_{\delta h}^*\tilde u_\omega$. On $\tilde\Omega^\delta$, by replacing $\tilde\varphi^\delta$ by $\rho_{L/R}^\delta\chi_{L/R}^\delta$ in \eqref{e:7l36} we have
\begin{multline}
    \rho_R^\delta\chi_R^\delta F_{\delta h}^*\tilde u_\omega=\rho_R^\delta\chi_R^\delta \Big(a_{+}^oe_{+,R,\delta}+a_{-}^ie_{-,R,\delta}\Big)+O_{H^1}(h^{{k_1}-({k_1}+2)\delta})\\
    =O_{X^{0,0}_{0,-N}}(h^{N\delta})+O_{H^1}(h^{{k_1}-({k_1}+2)\delta}),
\end{multline}
and $\rho_L^\delta\chi_L^\delta F_{\delta h}^*\tilde u_\omega=O_{X^{0,0}_{0,-N}}(h^{N\delta})+O_{H^1}(h^{k_1-(k_1+2)\delta})$ accordingly. Similarly to \eqref{e:7l49}, apply \eqref{e:7l55} to see
\begin{equation}
    \|\varphi^\delta \mathcal{R}_{h^{1-\delta}M}^{\tilde{\Omega}^\delta}[-\Delta,\tilde\varphi^\delta]\chi_\pm^\delta F_{\delta h}^*\tilde u_\omega\|_{\dot H^1}\le Ch^{k_1-(k_1+2)\delta}.
\end{equation}
 
3e. We now finish proving the claim \eqref{e:7l45}. Steps 3a--3d imply \eqref{e:7l35} can be estimated by 
\begin{equation}
    \varphi^\delta \mathcal{R}_{h^{1-\delta}M}^{\tilde\Omega^\delta}(-\Delta+h^{2-2\delta}M^2)\tilde\varphi^\delta(F^*_{\delta h} \tilde u_\omega-\chi_L^\delta\bar w_{L,\delta}-\chi_R^\delta\bar w_{R,\delta})=O_{\dot H^1(\tilde\Omega^\delta)}(h^{k_1-(k_1+2)\delta}), 
\end{equation}
and thus \eqref{e:7l50} by 
\begin{equation}
    \varphi^\delta f=O_{H^1(\tilde\Omega^\delta)}(h^{k_1-(k_1+2)\delta}+M^2h^{2-\frac\delta2}).
\end{equation}
Pulling $\varphi^\delta f$ back to $\tilde{\mc{W}}_\alpha^h$, we have
\begin{gather}
    \varphi^\delta f=O_{L^2(\tilde{\mc{W}}_\alpha^h)} (h^{k_1-(k_1+3)\delta}+M^2h^{2-\frac\delta2}), \ \nabla \varphi^\delta f=O_{L^2(\tilde{\mc{W}}_\alpha^h)} (h^{k_1-(k_1+2)\delta}+M^2h^{2-\frac\delta2}), \\ \varphi^\delta f=O_{X_0^0(\tilde{\mc{W}}_\alpha^h)} (h^{k_1-(k_1+3)\delta}+M^2h^{2-\frac\delta2}).
\end{gather}
Revisit \eqref{e:7l29} to give the estimates for the zero frequencies
\begin{equation}
    \partial_\nu\mc{R}^{\tilde{\mc{W}}^h_{\alpha}}_{hM}[-\Delta,\chi^\delta]\Big(F^*_h \tilde u_\omega-\chi_L^\delta\bar w_{L,\delta}-\chi_R^\delta\bar w_{R,\delta}\Big)=\partial_\nu\mc{R}^{\tilde{\mc{W}}^h_{\alpha}}_{hM}[-\Delta,\chi^\delta]\varphi^\delta f=O_{H^{-\frac12, 0}_\partial}(h^{k_1-(k_1+\frac72)\delta}+M^2h^{2-\frac\delta2}).
\end{equation}
For the nonzero frequencies, similarly to \eqref{e:8l7} we have
\begin{multline}
    \|\gamma^h (1-\psi_0)(D)\gamma^h\partial_\nu\mc{R}^{\tilde{\mc{W}}^h_{\alpha}}_{hM}[-\Delta,\chi^\delta]\Big(F^*_h \tilde u_\omega-\chi_L^\delta\bar w_{L,\delta}-\chi_R^\delta\bar w_{R,\delta}\Big)\|_{H^{-\frac12, N}}\le C\|\varphi^\delta f\|_{X^{0, N}_{0}}
    \\
    \le C h^{k_1-(N+k_1+\frac72)\delta}+CM^2h^{2-N\delta},
\end{multline}
similar to \eqref{e:8l6}. This gives the claim \eqref{e:8l8} stated at the beginning of Step 3. Together with Step 2 we have \eqref{e:7l45} as desired. 

4. In the end, by setting $0<\delta\ll 1$ and renumbering $N$, we have
\begin{equation}
(\Lambda^{\tilde{\mc{W}}^h_{\alpha}}_{hM}-\ohm)v'=O_{Z^{-\frac12, N, -\frac12}_{-L_0+0, N+1, N+1}}(h^{\min(2\rardya,k_1)-0}+M^2h^{2-0}), \ \|v'\|_{Z^{\frac12,N,\frac12}_{-L_0+0,-\frac12-0,-\frac12-0}}\sim 1,
\end{equation}
where near $\pm\omega$ frequencies we have
\begin{equation}
    v'(x)\sim (a^i_\mp+O(M^2h^{2})) e^{-i\omega|x|}+(a^o_\pm+O(M^2h^{2})) e^{i\omega|x|},
\end{equation}
as $x\rightarrow \pm\infty$. Let $v=v'-(\Lambda^{\tilde{\mc{W}}^h_{\alpha}}_{hM}-\ohm-i0)^{-1}(\Lambda^{\tilde{\mc{W}}^h_{\alpha}}_{hM}-\ohm)v'$. Proposition~\ref{p:LAP} implies $(\Lambda^{\tilde{\mc{W}}^h_{\alpha}}_{hM}-\ohm-i0)^{-1}: Z^{-\frac12, N, -\frac12}_{-L_0+0, N+1, N+1}\rightarrow Z^{\frac12, N, \frac12}_{-L_0+0, N, -\frac12-0}$ is bounded, and we have
\begin{equation}
    \|(\Lambda^{\tilde{\mc{W}}^h_{\alpha}}_{hM}-\ohm-i0)^{-1}(\Lambda^{\tilde{\mc{W}}^h_{\alpha}}_{hM}-\ohm)v'\|_{Z^{\frac12,N,\frac12}_{-L_0+0, N, -\frac12-0}}\le Ch^{\min(2\rardya,k_1)-0}+CM^2h^{2-0}
\end{equation}
and $(\Lambda^{\tilde{\mc{W}}^h_{\alpha}}_{hM}-\ohm)v=0$. We have
\begin{equation}
    v(x)\sim (a^i_\mp+O(M^2h^{2-0})+O(h^{\min(2\rardya,k_1)-0})) e^{-i\omega|x|}+(a^o_\pm+O(M^2h^{2-0})+O(h^{\min(2\rardya,k_1)-0})) e^{i\omega|x|}
\end{equation}
near $\pm\omega$ frequencies. This implies
\begin{equation}
    \Sigma^{\tilde{\mc{W}}_\alpha}_{M, h^{-1}\omega}=\Sigma^{\tilde{\mc{W}}_\alpha^h}_{hM, \omega}=\Sigma^{\mc{W}_\alpha^h}_{hM, \omega}+O(M^2h^{2-0})+O(h^{\min(2\rardya,k_1)-0})=\Sigma^{\mc{W}_\alpha}_{M, h^{-1}\omega}+O( M^2h^{2-0})+O(h^{\min(2\rardya,k_1)-0})
\end{equation}
via rescaling. 
\end{proof}

\subsection{Comparison with $M=0$}
We can now work on the straight wedges $\mc{W}_{\alpha_j}$. In this subsection, we compare the scattering matrices at high frequencies and fixed intermediate frequencies, on straight wedges. 
\begin{lemma}\label{l:schfmzero}
Let $M_0>0$, $\omega\in[\frac12, 2]$, $\alpha\in (0,2\pi)$. Then there exists $h_0>0$ that uniformly for all $M\in[0, M_0]$, all $h\in(0, h_0]$, we have
\begin{equation}
    \Sigma^{\mc{W}_\alpha}_{M, h^{-1}\omega}=\Sigma^{\mc{W}_\alpha}_{0, \omega}+O((hM)^{\min(2\rardya,2)-0}+h^\infty).
\end{equation}
\end{lemma}
\begin{proof}
1. Let $a^i=(a^i_+, a^i_-)\in \mathbb{C}^2$ with $|a^i|=1$. In this proof, we adopt the geometry set up in Section~\ref{s:8-1} with $\tilde{\mc{W}}_\alpha=\mc{W}_\alpha$, and drop the tilde as long as there is no confusion. Proposition~\ref{p:existscmatrix} implies there exists $u\in Z^{\infty, \infty, 1-0}_{L_0-0,-\frac12-0, -\frac12-0}(\partial\mc{W}_\alpha^h; \omega)$ that $(\Lambda_{hM}^{\mc{W}_\alpha^h}-\ohm)u(x)=0$ with $u\sim \tilde a^i_\mp(x)e^{-i\omega|x|}+\tilde a^o_\pm(x)e^{i\omega|x|}$ away from the zero frequency as $x\rightarrow \pm\infty$, where
    \begin{equation}\label{e:7l65}
        \tilde a_\pm^i(x)=a_\pm^i+O(|x|^{-1}), \ \tilde a_{\pm}^o(x)=a_\pm^o+O(|x|^{-1}), \ a^o=(a^o_+, a^o_-)=\Sigma^{\mc{W}_\alpha}_{M, h^{-1}\omega} a^i=\Sigma^{\mc{W}_\alpha^h}_{hM, \omega} a^i.
    \end{equation}

2. We show $\tilde a_{\pm}^{i/o}(x)$ are $h^2M^2$-close to $a_{\pm}^{i/o}$ as $|x|\rightarrow\infty$. In order to examine the gains in decay at different frequencies, we consider cutoff $\chi_0, \chi,\chi_1, \tilde\chi\in C_c^\infty(\mathbb{R})$ with $\chi_0\equiv 1$ near $0$, $\chi_0\prec\chi\prec\chi_1\prec\tilde\chi$. For any of these cutoffs, $\psi$, set $\psi^h(x)=\psi(h^{\delta}x)$ for some $\delta>0$ to be chosen arbitrarily small. Let also $\psi,\tilde\psi,\tilde{\tilde{\psi}}\in C_c^\infty( -2\omega,2\omega)$ with $\supp(1-\psi)\cap\{\pm \omega\}=\emptyset$, $\supp{\tilde{\tilde{\psi}}}\cap\{0\}=\emptyset$ and $\psi\prec\tilde\psi\prec\tilde{\tilde{\psi}}$.

By Lemma~\ref{l:sameNearInfinity}, we can compare $\Lambda^{\mc{W}_{\alpha}^h}_{hM}$ to the model operator
\begin{equation}
    \Lambda^{L/R}_{hM}= 1_{x>0}\Lambda_{hM}^{\mc{W}^{h}_{\alpha,R}}1_{x>0}+1_{x<0}\Lambda_{hM}^{\mc{W}^{h}_{\alpha,L}}1_{x<0}. 
\end{equation}
We have
\begin{gather}
\label{e:rep1}
    (1-\chi^h_{1})\Lambda^{L/R}_{hM}=(1-\chi^h_{1})\Lambda^{L/R}_{hM}(1-\chi_{0}^h)+(1-\chi^h_{1})\tilde{\mc{E}}_{hM},\\
    \Lambda^{L/R}_{hM}(1-\chi^h_{1})=(1-\chi^h_{0})\Lambda^{L/R}_{hM}(1-\chi^h_{1})+\tilde{\mc{E}}_{hM}^*(1-\chi^h_{1}),\label{e:rep2}
\end{gather}
where $\tilde{\mc{E}}_{hM}$ is the error in Lemma~\ref{l:sameNearInfinity} given by
\begin{equation}
    \tilde{\mathcal{E}}_{hM}=\sum_{L/R}\partial_\nu \mc{R}_{hM}^{\mc{W}^{h}_{\alpha,L/R}}[\Delta, \chi_{0}^h]\mc{H}^{\mc{W}_\alpha^h}_{hM}.    
\end{equation}
Notice that $u$ satisfies
\begin{equation}
    (1-\chi^h_1)(\Lambda^{L/R}_{hM}(1-\chi_0^h) +\tilde{\mc{E}}_{hM}-\ohm )u=0.
\end{equation}
Therefore, using Lemma~\ref{l:1dnormalform} on $\mc{W}^{h}_{\alpha,L/R}$ and then rescaling, we obtain $b_0\in S^{0,-1}$ such that 
\begin{multline}
(1-\chi^h_1)\Big(\big(\sqrt{D^2+h^2M^2}\tilde\psi(D)+h^3M^2b_0(hx,D)\big)(1-\chi^h_0)-\ohm \Big)\psi(D)u\\
=-(1-\chi^h_1)\big(\Lambda_{hM}^{L/R}(1-\chi^h_0) +\tilde{\mc{E}}_{hM}-\ohm \big)(1-\psi(D))u \\
\qquad -(1-\chi^h_1)E_{hM}^{L/R}(1-\chi^h_0)\psi(D)u-(1-\chi_1^h)\tilde{\mc{E}}_{hM}\psi(D)u.
\end{multline}
Now, by Proposition~\ref{p:existscmatrix} $(1-\psi(D))u\in Z^{\infty,\infty,1-0}_{L_0-0,\infty,\infty}$. Therefore, using Lemma~\ref{l:sameNearInfinity} and Proposition~\ref{p:dtN} and the estimate on $E_{hM}^{L/R}$ from Lemma~\ref{l:1dnormalform}, 
\begin{multline} \label{e:highFreqStuff}
\|\tilde\psi(D)(1-\chi_1^h)\Big(\big(\sqrt{D^2+h^2M^2}\tilde\psi(D)+h^3M^2b_0(hx,D)\big)(1-\chi_0^h) -\ohm \Big)\psi(D)u\|_{H^{N,N}}\\
\leq C_Nh^N\|u\|_{Z^{s,l,s_0}_{L_0-0,-\frac{1}{2}-0,-\frac{1}{2}-0}}.
\end{multline}
Since we already know that 
\begin{equation}
    (1-\chi_0^h)\psi(D)u= (1-\chi_0^h)\psi(D)\sum_{\pm}\big(\tilde a_{\pm}^{i}(x)e^{-i|x|\omega}+\tilde a_{\pm}^{o}(x)e^{i|x|\omega}\big)+O(\mc{S}),
\end{equation}
the calculations in Lemma~\ref{l:2dquasimodeappendix} imply that
\begin{equation}
\label{e:uForm}
(1-\chi^h)(\tilde a_{\pm}^{i/o}(x)-( a_{\pm}^{i/o}+h^2M^2b_{\pm}^{i/o}(x)))\in h^\infty\mc{S}
\end{equation}
with $b_{\pm}^{i/o}(x)$ satisfying
\begin{equation}
    |\partial^\beta_x b_{\pm}^{i/o}(x)|\leq C\langle x\rangle^{-2}\|u\|_{Z^{s,l,s_0}_{L_0-0,-\frac{1}{2}-0,-\frac{1}{2}-0}}, \ \beta\geq 1.
\end{equation}

3. We consider
\begin{equation}
    F_{hM}=\Lambda_{hM}^{\mc{W}^h_\alpha}-\ohm-\Lambda_0^{\mc{W}^h_\alpha}+\omega.
\end{equation}
We claim for each $s,l\in\mathbb{R}$, $s_0\in(0,1)$, we have
\begin{equation}\label{e:7l64}
    \|F_{hM}u\|_{Z^{-\frac{1}{2},  -\rardya ,\frac{1}{2}}_{-\rardya-\frac{1}{2}+0, \frac32-0, \frac32-0}}\le C ((hM)^{\min(2\rardya,2)-0}+h^\infty)\|u\|_{Z^{\frac{1}{2},\rardya+\frac{1}{2}-0, \frac{1}{2}}_{\rardya+\frac{1}{2}-0, -\frac12-0, -\frac12-0}}.
\end{equation}
Once we have the claim, its application gives
\begin{equation}
    \|(\Lambda_0^{\mc{W}^h_\alpha}-\omega-i0)^{-1}(\Lambda_0^{\mc{W}^h_\alpha}-\omega) u\|_{Z^{\frac{1}{2},-\rardya,\frac{1}{2}}_{-\rardya-\frac{1}{2}+0, \frac12-0, -\frac12-0}}\le C (hM)^{\min(2\rardya,2)-0}+h^\infty,
\end{equation}
uniformly in $M\in[0, M_0]$ and $h$ small, since $(\Lambda_0^{\mc{W}^h_\alpha}-\omega) u=-F_{hM}u$. We then construct a sequence of $v$ by
\begin{equation}
    v=u-(\Lambda_0^{\mc{W}^h_\alpha}-\omega-i0)^{-1}(\Lambda_0^{\mc{W}_\alpha^h}-\omega) u\in Z^{\infty, \infty, 1-0}_{-L_0+0, -\frac12-0, -\frac12-0}, \ (\Lambda_0^{\mc{W}_\alpha^h}-\omega)v=0,
\end{equation}
with $u-v=O_{Z^{\frac{1}{2}, -\rardya, \frac{1}{2}}_{-\rardya-\frac{1}{2}+0, -\frac12+0, -\frac12-0}}((hM)^{\min(2\rardya,2)-0}+h^\infty)$. This implies $u$, $v$ differ by that size in the scattering data, that is
\begin{equation}
\begin{aligned}
    \Sigma^{\mc{W}_\alpha}_{0, h^{-1}\omega}a^i=\Sigma^{\mc{W}_\alpha^h}_{0, \omega}a^i
    &=\big(\Sigma^{\mc{W}^h_\alpha}_{hM, \omega}+O((hM)^{\min(2\rardya,2)-0}+h^\infty)\big) a^i\\
    &=\big(\Sigma^{\mc{W}_\alpha}_{M, h^{-1}\omega}+O((hM)^{\min(2\rardya,2)-0}+h^\infty)\big) a^i
    \end{aligned}
\end{equation}
as desired. The rest of the proof will be devoted to \eqref{e:7l64}. We can decompose
\begin{equation}\label{e:7l27}
    F_{hM}=\tilde\chi^h F_{hM} \chi^h+(1-\tilde\chi^h)F_{hM} (1-\chi^h)+\tilde\chi^h F_{hM} (1-\chi^h)+(1-\tilde\chi^h) F_{hM} \chi^h,
\end{equation}
which we estimate term by term. 

3a. For arbitrary eligible $\chi,\tilde\chi$, Lemma~\ref{l:errordependsonMcpt} implies
\begin{equation}
    \|\tilde\chi^h F_{hM} \chi^h \|_{H^{0, -N}\rightarrow H^{0,N}}\le C(hM)^{\min(2\rardya-0,2)}h^{-C\delta}.
\end{equation}
This specifically implies
\begin{equation}
    \|\tilde\chi^h F_{hM}\chi^h u\|_{Z^{N,  N,0}_{N, N, N}}\le C (hM)^{\min(2\rardya-0,2)}h^{-C\delta}\|u\|_{Z^{-N, -N, 0}_{-N, -N, -N}},
\end{equation}
which allows us to estimate the first term of \eqref{e:7l27}. 

3b. We now consider the second term of \eqref{e:7l27}. By~\eqref{e:rep1} it is enough to estimate
\begin{multline}
(1-\tilde\chi^h)F_{hM}(1-\chi^h)=(1-\tilde\chi^h)(\Lambda^{L/R}_{hM}-\ohm-\Lambda^{L/R}_0+\omega+\tilde{\mc{E}}_{hM}-\tilde{\mc{E}}_0)(1-\chi^h)\\
= (1-\tilde\chi^h)(\Lambda_{hM}^{L/R}-\ohm-\Lambda_0^{L/R}+\omega)(1-\chi^h) +(1-\tilde\chi^h)(\tilde{\mc{E}}_{hM}-\tilde{\mc{E}}_0)(1-\chi^h).
\end{multline}
We first show that 
\begin{equation}
    \label{e:7l70}
\|\tilde\psi(D)(1-\tilde{\chi}^h)(\Lambda_{hM}^{L/R}-\ohm-\Lambda_0^{L/R}+\omega)(1-\chi^h)u\|_{Z^{N,N,1}_{N,\frac{3}{2}-0,\frac{3}{2}-0}}\leq C\big((hM)^2+h^N\big)\|u\|_{Z^{s,l,s_0}_{L_0-0,-\frac{1}{2}-0,-\frac{1}{2}-0}}.
\end{equation}
Using Lemma~\ref{l:1dnormalform} again and rescaling, we have both
\begin{multline}
\label{e:leftCut}
\tilde\psi(D)(1-\tilde{\chi}^h)(\Lambda_{hM}^{L/R}-\Lambda_0^{L/R})(1-\chi^h)
= \tilde\psi(D)(1-\tilde{\chi}^h)\Big(\sqrt{D^2+h^2M^2}\tilde{\tilde{\psi}}(D)+h^3M^2b_0(hx,D)\\-|D|\tilde{\tilde{\psi}}(D)\Big)(1-\chi^h)+O_{H^{-N,-\rardya-\frac12+0}\rightarrow H^{N,N}}(h^\infty)
\end{multline}
and 
\begin{multline}\label{e:rightCut}
(1-\tilde{\chi}^h)(\Lambda_{hM}^{L/R}-\Lambda_0^{L/R})(1-\chi^h)\psi(D)
= (1-\tilde{\chi}^h)\Big(\sqrt{D^2+h^2M^2}\tilde{\tilde{\psi}}(D)+h^3M^2b_0(hx,D)\\
-|D|\tilde{\tilde{\psi}}(D)\Big)(1-\chi^h)\psi(D)+O_{H^{-N,-N}\rightarrow H^{N,\rardya+\frac12-0}}(h^\infty).
\end{multline}
Now, by~\eqref{e:highFreqStuff} we have
\begin{multline}
\label{e:7l69}
\tilde\psi(D)(1-\tilde{\chi}^h)(\sqrt{D^2+h^2M^2}\tilde{\tilde{\psi}}(D)+h^3M^2b_0(hx,D)-\ohm-|D|\tilde{\tilde{\psi}}(D)+\omega)(1-\chi^h)u\\
=-\tilde\psi(D)(1-\tilde{\chi}^h)(|D|\tilde{\tilde{\psi}}(D)-\omega)(1-\chi^h)u+O_{\mc{S}}(h^\infty). 
\end{multline}
By \eqref{e:uForm} we have 
\begin{multline}\label{e:7l66}
    \|\tilde\psi(D)(1-\tilde{\chi}^h)(|D|\tilde{\tilde{\psi}}(D)-\omega)(1-\chi^h)u\|_{Z^{s-1,l}_{N, \frac32-0, \frac32-0}}\le C(hM)^2\|(1-\chi)u\|_{Z^{s,l}_{L_0-0, -\frac12-0, -\frac12-0}}\\
    +Ch^N\|u\|_{Z^{s,l,s_0}_{-L_0+0, -N, -N}},
\end{multline}
since
\begin{equation}
    |D|\tilde{\tilde{\psi}}(D)e^{\pm i\omega x}\tilde a^{i/o}_\pm(x)\mp  D e^{\pm i \omega x}\tilde a^{i/o}_{\pm}(x) \in \mc{S}.
\end{equation}
The above implies
\begin{equation}
    \|(1-\tilde\psi(D))(1-\tilde{\chi}^h)(\Lambda_{hM}^{L/R}-\ohm-\Lambda_0^{L/R}+\omega)(1-\chi^h)\psi(D) u\|_{Z^{N,N,1}_{L_0-0,N,N}}\leq C_Nh^N\| u\|_{Z^{s,l,s_0}_{-L_0+0,-N,-N}}.
\end{equation}
On another hand, we can use Proposition~\ref{p:harmonic} and Remark~\ref{r:pickm} to see
\begin{multline}
\|(1-\tilde\psi(D))(1-\tilde{\chi}^h)(\Lambda_{hM}^{L/R}-\ohm-\Lambda_0^{L/R}+\omega)(1-\chi^h)(1-\psi(D))u\|_{Z^{-\frac{1}{2},-\rardya+\frac{1}{2}-0,-\frac{1}{2}}_{-\rardya+\frac{1}{2}+0,N,N}}\\
\leq C\|(1-\tilde\psi(D))(1-\tilde{\chi}^h)(\Lambda_{hM}^{L/R}-\ohm-\Lambda_0^{L/R}+\omega)(1-\chi^h)(1-\psi(D))u\|_{H^{-\frac12,-\rardya+\frac{1}{2}+0}}\\
\leq C(hM)^{\min(2\rardya,2)-0}\|(1-\chi^h)(1-\psi(D))u\|_{H^{\frac12,\rardya+\frac{1}{2}-0}}\leq C(hM)^{\min(2\rardya,2)-0}\|u\|_{Z^{\frac{1}{2},\rardya+\frac{1}{2}-0,\frac{1}{2}}_{\rardya+\frac{1}{2}-0,-N,-N}}. 
\end{multline}
This demonstrates \eqref{e:7l70}. 

 For the error terms, one notes also we have from Lemma~\ref{l:sameNearInfinity} and Lemma~\ref{l:errordependsonM} that
 \begin{equation}
 \label{e:errorM}
\begin{gathered}
    (1-\tilde\chi^h)(1-\psi_{0}(D))(1-\tilde\chi^h)(\tilde{\mc{E}}_{hM}-\tilde{\mc{E}}_0)(1-\chi^h)(1-\psi_{0}(D))(1-\chi^h)=O_{H^{-N, -N}\rightarrow H^{N, N}}((hM)^{\min(2\rardya,2)-0}),\\
    (1-\tilde\chi^h)(1-\psi_{0}(D))(1-\tilde\chi^h)(\tilde{\mc{E}}_{hM}-\tilde{\mc{E}}_0)(1-\chi^h)=O_{H^{0, L_0-0}\rightarrow H^{N, N}}((hM)^{\min(2\rardya,2)-0}),\\
    (1-\tilde\chi^h)(\tilde{\mc{E}}_{hM}-\tilde{\mc{E}}_0)(1-\chi^h)(1-\psi_{0}(D))(1-\chi^h)=O_{H^{-N, -N}\rightarrow H^{N, -L_0+0}}((hM)^{\min(2\rardya,2)-0}),\\
    (1-\tilde\chi^h)(\tilde{\mc{E}}_{hM}-\tilde{\mc{E}}_0)(1-\chi^h)=O_{H^{0, L_0-0}\rightarrow H^{N, -L_0+0}}((hM)^{\min(2\rardya,2)-0}).
\end{gathered}
\end{equation}
This implies 
\begin{equation}
    \|(1-\tilde\chi^h)(\tilde{\mc{E}}_{hM}-\tilde{\mc{E}}_0)(1-\chi^h)u\|_{Z^{N,  N,1}_{-L_0+0, N, N}}\le C (hM)^{\min(2\rardya,2)-0}\|u\|_{Z^{-N,-N, -1}_{L_0-0, -N, -N}}. 
\end{equation}

3c. For the third term of \eqref{e:7l27}, we can use~\eqref{e:rep2} to see that it is enough to estimate
\begin{equation}
    \tilde{\chi}^h(1-\chi^h_{0})(\Lambda_{hM}^{L/R}-\ohm-\Lambda_0^{L/R}+\omega)(1-\chi^h_{1})u+\tilde{\chi}^h(\tilde{\mc{E}}_{M}^*-\tilde{\mc{E}}_0^*)(1-\chi^h_{1})u.
\end{equation}
The first term can be estimated as in Step 3b and the last term using Lemma~\ref{l:errordependsonM} as in~\eqref{e:errorM}. 

3d. For the fourth term of \eqref{e:7l27}, noting $\supp\chi^h\cap \supp(1-\tilde{\chi}^h)=\emptyset$ and using the representation~\eqref{e:leftCut} and~\eqref{e:rep1}, we have
\begin{gather}
(1-\tilde\chi^h) F_M\chi^h u=O(h^\infty)_{\Psi^{-\infty,-\infty}}u+(1-\tilde\chi^h)(\tilde{\mc{E}}_M-\tilde{\mc{E}}_0)\chi^h
\end{gather}
and the estimate follows similar to~\eqref{e:errorM}. Bring all together to observe \eqref{e:7l64} as desired. 
\end{proof}

\subsection{Derivatives of scattering matrices}

It will be convenient in Section~\ref{s:spectralAsymptotics} to know that the scattering matrix on a any curved wedge is $h^{2\pi/\alpha-0}$ close to one with bounded derivatives.  (In fact, it is possible to show that the matrices themselves have bounded derivatives, but since it would complicate the proof and we do not need it, we do not pursue this generalization here.)

Throughout this subsection, we work on a curved wedge whose boundary consists of two half-lines outside of a compact set. We then identify the boundary of a curved wedge, $\mc{W}_\alpha$ with 
 $\gamma:\mathbb{R}\to \partial\mc{W}_{\alpha}$ such that $\gamma$ is smooth away from $0$, has unit speed, and satisfies $\gamma''(s)=0$ for $|s|\geq L$. 

\begin{lemma}
\label{l:approxResolvent}
    Suppose that $\mc{W}_\alpha$ is a curved wedge whose boundary consists of two half-lines outside of a compact set. Let $\chi\in C^\infty(\mathbb{R})$ with $\supp \chi\subset (L,\infty)$, $\supp (1-\chi)\subset (-\infty, L+1)$. Then for all $j\geq 0$
    \begin{equation}
        (h\Lambda_0^{\mc{W}_\alpha}-\omega-i0)^{-1} x^j\chi_\pm e^{\pm i x\omega/h}\mp \tfrac{i}{h(j+1)}x^{j+1}\chi_\pm e^{\pm ix\omega/ h}=:e_{j,\pm}\in  h^{-1}Z_{\min(\frac{3}{2},\frac{\pi}{\alpha}+\frac{1}{2})-0,\infty,-\frac{1}{2}-0}^{\infty,\infty,1-0},
    \end{equation}
    where $\chi_{\pm}(x)=\chi(\pm x)$. Moreover, there are $a_1,a_2\in\mathbb{C}$ such that for $\psi\in C_c^\infty((-\omega,\omega))$ with $0\notin \supp 1- \psi$ and $\phi\in C_c^\infty$ with $\supp (1-\phi)\cap [-L,L]=\emptyset$,
    \begin{equation}
        (1-\psi(hD))\psi(x)e_{\pm,j}= h^{-1} (a_1\chi_+e^{ix\omega/h}+a_2\chi_- e^{-ix/h})+h^{-1}\mc{S}.
    \end{equation}
\end{lemma}
\begin{proof}

    We consider the case of $+$, the other being similar
    First, notice that, with $u_j:=x^{j+1}\chi_{\pm}e^{ix\omega/h}$, 
    \begin{equation}
        \|u_j\|_{Z^{\infty, \infty, 1}_{\infty,\infty,-\frac{1}{2}-j-0}}\leq C.
    \end{equation}
    Therefore, we may apply $\Lambda_0^{\mc{W}_\alpha}$ to $u_j$. To prove the lemma, let $\tilde{\chi}_+\in C^\infty(\mathbb{R})$ with $\supp \tilde{\chi}_+\subset (L,\infty)$,$\supp \chi_+\cap \supp(1-\tilde{\chi}_+)=\emptyset$. Then using Lemma~\ref{l:sameNearInfinity} we compute
\begin{multline}\label{e:ladeeda}
        (h\Lambda_0^{\mc{W}_\alpha}-\omega)u_j=\tilde{\chi}_+(h|D|-\omega)\chi_+x^{j+1}e^{ix\omega/h}+\mc{E}_0'\chi_+x^{j+1} e^{ix\omega/h}
     = \tilde{\chi}_+(h|D|\psi(hD)-\omega)\chi_+x^{j+1}e^{ix\omega/h}\\+\tilde{\chi}_+h|D|(1-\psi(hD))\chi_+x^{j+1}e^{ix\omega/h}+\mc{E}_0'\chi_+x^{j+1} e^{ix\omega/h}
     = \chi_+(h|D|\psi(hD)-\omega)x^{j+1}e^{ix\omega/h}\\
     +\tilde{\chi}_+[h|D|\psi(hD),\chi_+]x^{j+1}e^{ix\omega/h}
     +\tilde{\chi}_+h|D|(1-\psi(hD))\chi_+x^{j+1}e^{ix\omega/h}+\mc{E}_0'\chi_+ x^{j+1}e^{ix\omega/h}\\
          =-ih\chi_+(j+1)x^je^{ix\omega/h}+\tilde{\chi}_+[h|D|\psi(hD),\chi_+]x^{j+1}e^{ix\omega/h}+\tilde{\chi}_+h|D|(1-\psi(hD))\chi_+x^{j+1}e^{ix\omega/h}\\
          +\mc{E}_0'\chi_+x^{j+1} e^{ix\omega/h}
          =:-ih\chi_+(j+1)x^je^{ix\omega/h}+R_1+R_2+R_3.
    \end{multline}
    Now, we estimate the last three terms in~\eqref{e:ladeeda}. First, notice that standard scattering calculus estimates imply
    \begin{equation}
        R_1=\tilde{\chi}_+[h|D|\psi(hD),\chi_+]x^{j+1}e^{ix\omega/h}\in h Z^{\infty,\infty,\infty}_{\infty,\infty,\infty}.
    \end{equation}
    Next, since $(1-\psi(hD))\chi_+x^{j+1}e^{ix\omega/h}\in h^\infty \mc{S}$, 
    \begin{equation}
        R_2=\tilde{\chi}_+h|D|(1-\psi(hD))\chi_+x^{j+1}e^{ix\omega/h}\in h^\infty Z^{\infty,\infty, \infty}_{\frac{3}{2}-0,\infty,\infty}
    \end{equation}
    and, using Lemma~\ref{l:sameNearInfinity}
    \begin{equation}
        R_3=\mc{E}_0'\chi_+x^{j+1} e^{ix\omega/h}\in h^\infty Z^{\infty,\infty, \infty}_{\frac{\pi}{\alpha}+\frac{1}{2}-0,\infty,\infty}.
    \end{equation}
    Hence, since $u_j$ is outgoing
    \begin{equation}
        (h\Lambda_0^{\mc{W}_\alpha}-\omega-i0)^{-1}x^j\chi_{+}e^{ix\omega/h} = \frac{i}{h(j+1)}\Big(u_j -(h\Lambda_0^{\mc{W}_\alpha}-\omega-i0)^{-1}(R_1+R_2+R_3)\Big).
    \end{equation}
    The lemma now follows from Lemmas~\ref{l:uniformres} and Lemma~\ref{l:PlaneWave}.
\end{proof}
\begin{lemma}
\label{l:scatteringDerivatives}
    Suppose that $\mc{W}_\alpha$ is a curved wedge whose boundary consists of two half-lines outside of a compact set. Then,
    \begin{equation}
        |h^j\partial_\omega^j \Sigma^{\mc{W}_\alpha}_{0,h^{-1}\omega}|\leq C_j,\qquad j=1,2. 
    \end{equation}
\end{lemma}
\begin{proof}
    Let $\chi_-,\tilde{\chi}_-\in C^\infty(\mathbb{R})$ with $\supp \tilde{\chi}_-\subset (-\infty,-L)$,$\supp \chi_-\cap \supp(1-\tilde{\chi}_-)=\emptyset$, and $\supp (1-\tilde{\chi}_-)\cap (-\infty,-L-1)= \emptyset$. Define $u^-:= \chi_-(x)e^{ix\omega/h}$. 

    Then, observe that by Lemma~\ref{l:sameNearInfinity},
    \begin{align*}
    (h\Lambda_{0}^{\mc{W}_\alpha}-\omega)u^-&= (h\Lambda_{0}^{\mc{W}_\alpha}-\omega)\chi_-e^{ix\omega/h}\\
     &= \tilde{\chi}_-(h|D|-\omega)\chi_-e^{ix\omega/h}+\mc{E}_0'\chi_- e^{ix\omega/h}\\
     &= \tilde{\chi}_-(h|D|\psi(hD)-\omega)\chi_-e^{ix\omega/h}+\tilde{\chi}_-h|D|(1-\psi(hD))\chi_-e^{ix\omega/h}+\mc{E}_0'\chi_- e^{ix\omega/h}\\
     &= \chi_-(h|D|\psi(hD)-\omega)e^{ix\omega/h}+\tilde{\chi_-}[h|D|\psi(hD),\chi_-]e^{ix\omega/h}\\
     &\qquad \tilde{\chi}_-h|D|(1-\psi(hD))\chi_-e^{ix\omega/h}+\mc{E}_0'\chi_- e^{ix\omega/h}\\
          &=\tilde{\chi}_-[h|D|\psi(hD),\chi_-]e^{ix\omega/h}+\tilde{\chi}_-h|D|(1-\psi(hD))\chi_-e^{ix\omega/h}+\mc{E}_0'\chi_- e^{ix\omega/h}.
    \end{align*}
    
    Therefore, 
    $$    (hD_\omega)^\ell(h\Lambda_{0}^{\mc{W}_\alpha}-\omega)u^-=\big(\tilde{\chi}_-[h|D|\psi(hD),\chi_-]+\tilde{\chi}_-h|D|(1-\psi(hD))\chi_-+\mc{E}_0'\chi_- \big)x^\ell e^{i\omega x/h}.
    $$
    In particular, arguing as in Lemma~\ref{l:approxResolvent}, we have 
    $$
(hD_\omega)^\ell(h\Lambda_{0}^{\mc{W}_\alpha}-\omega)u^-\in hZ^{\infty,\infty,\infty}_{\min(\frac{3}{2},\frac{\pi}{\alpha}+\frac{1}{2})-0,\infty,\infty}
    $$
    Thus, 
    $$
    (h\Lambda_0^{\mc{W}_\alpha} -\omega-i0)^{-1}(hD_\omega)^\ell(h\Lambda_{0}^{\mc{W}_\alpha}-\omega)u^-=a_{+,\ell}\chi_+e^{ix\omega/h}+a_{-,\ell}\chi_-e^{-ix\omega/h} +Z^{\infty,\infty,1-0}_{\frac{\pi}{\alpha}+\frac{1}{2}-0,\infty,\infty}
    $$
    Applying Lemma~\ref{l:approxResolvent}, we thus obtain
    \begin{multline}\label{e:resolveFormula}
(h\Lambda_0^{\mc{W}_\alpha} -\omega-i0)^{-j-1}(hD_\omega)^\ell(h\Lambda_{0}^{\mc{W}_\alpha}-\omega)u^-=(ih^{-1})^ja_{+,\ell}\frac{x^j}{j!}\chi_+e^{ix\omega/h}+(ih^{-1})^ja_{-,\ell}\chi_-(-x)^je^{-ix\omega/h}\\
+h^{-j}Z^{\infty,\infty,1-0}_{\frac{\pi}{\alpha}+\frac{1}{2}-0,\infty,\infty}.
    \end{multline}
    Now, consider 
    $$
    u= u^--(h\Lambda_0^{\mc{W}_\alpha}-\omega-i0)^{-1}(h\Lambda_0^{\mc{W}_\alpha}-\omega)u^-.
    $$
    We study first $e^{-i\omega x/h}\chi_+u$.

    Define $f(\omega):= (h\Lambda_0^{\mc{W}_\alpha}-\omega )u^-$. Then, using~\eqref{e:resolveFormula} in the last line, we have
    \begin{align*}
    &hD_\omega e^{-i\omega x}\chi_+u\\
    &= -hD_{\omega}e^{-i\omega x}\chi_+(h\Lambda_0^{\mc{W}_\alpha}-\omega-i0)^{-1}f(\omega)\\
    &= -( e^{-i\omega x}\chi_+ (-x-ih(h\Lambda_0^{\mc{W}_\alpha}-\omega-i0)^{-1})(h\Lambda_0^{\mc{W}_\alpha}-\omega-i0)^{-1}f(\omega)\\
    &\qquad +e^{-i\omega x} \chi_+ (h\Lambda_0^{\mc{W}_\alpha}-\omega-i0)^{-1}hD_\omega f(\omega)\\
    &=b\chi_+^2 +e^{-i\omega x/h}Z_{\frac{\pi}{\alpha}+\frac{1}{2}-0,\infty,1-0}^{\infty,\infty,\infty},
    \end{align*}
for some $|b|\leq C$. Together with the same argument but considering 
    $
    e^{i\omega x/h}\psi(hD)\chi_-u
    $
    for some $\psi\in C_c^\infty (\mathbb{R}\setminus 0)$ with $-\omega\notin \supp (1-\psi)$, 
    instead, we obtain 
    $$
\|hD_\omega\Sigma_0^{\mc{W}_\alpha}\|\leq C.
    $$
    
    Next,
$$
\begin{aligned}
&(hD_\omega)^2\Big(
e^{-ix\omega/h}
(h\Lambda_0^{\mathcal W_\alpha}-\omega-i0)^{-1}f(\omega)
\Big)
\\
&\quad =
e^{-ix\omega/h}\Big[
(h\Lambda_0^{\mathcal W_\alpha}-\omega-i0)^{-1}
(hD_\omega)^2f(\omega)
\\
&\qquad
-2x(h\Lambda_0^{\mathcal W_\alpha}-\omega-i0)^{-1}
hD_\omega f(\omega)
-2ih(h\Lambda_0^{\mathcal W_\alpha}-\omega-i0)^{-2}
hD_\omega f(\omega)
\\
&\qquad
+x^2(h\Lambda_0^{\mathcal W_\alpha}-\omega-i0)^{-1}f(\omega)
+2ihx(h\Lambda_0^{\mathcal W_\alpha}-\omega-i0)^{-2}f(\omega)
\\
&\qquad
-2h^2(h\Lambda_0^{\mathcal W_\alpha}-\omega-i0)^{-3}f(\omega)
\Big].
\end{aligned}
$$
Again, using~\eqref{e:resolveFormula} (together with a similar argument for $
    e^{i\omega x/h}\psi(hD)\chi_-u
    $, this implies the lemma.
\end{proof}

\begin{remark}
    All constants in this section can depend on $\Omega$ in the following way: the upper bound on uniform Lipschitz characters of $\partial\Omega$, the upper bound on the diameter of $\partial\Omega\setminus\partial\tilde\Omega$ in Remark~\ref{r:aclcollar}, for each multi-index $\alpha$, the upper bounds on the derivatives of the metric $\gplus$ in \eqref{e:3l26}, and lower bounds on $|\det \gplus|$.
\end{remark}

\section{Steklov eigenfunctions on curvilinear polygons}
\label{s:aprioriSteklov}
In this section, we study the structure of Steklov eigenfunctions, $u$, on curvilinear polygons, $\mc{P}$. 

\noindent{\bf{ Outline of the strategy}}
We study the structure of $u$ in two steps. 

First, we provide a relatively rough description of $u$ on the interior of any smooth edge. More precisely we show in Lemma~\ref{l:fullDescription} that on the interior of a smooth edge $u$ consists of a sum of (approximate) plane waves and a low frequency function (i.e. a function with uniformly bounded derivatives as the energy tends to infinity). This is done by exploiting the fact that, near smooth points, the DtN map is an elliptic first order pseudodifferential operator and hence that, away from the zero frequency $h\Lambda_M^{\mc{P}}-\ohm$ is a \emph{semiclassical} psuedodifferential operator of real principal type. 

Once the rough structure from Lemma~\ref{l:fullDescription} is in place. We both improve it and determine the structure of $u$ near corners. The overall strategy is to transplant $u$ near a corner onto a curvilinear wedge, $\mc{W}_\alpha$ that agrees with the curvilinear polygon near a corner. Using a strategy very close to that in the geometric comparison of Section~\ref{s:scatteringatrixConvergence}, we are able to show that $u$, or rather an extension and rescaling, $v_h'$ of $u$ of the form~\eqref{e:extensionType} is nearly a solution to 
\begin{equation}
    (\Lambda_{hM}^{\mc{W}^h_\alpha}-\ohm)v_h'=0,
\end{equation}
and hence that $v_h'$ and therefore $u_h(x)=u(hx)$ inherits many of the properties of distorted plane waves on the infinite wedge. 

\subsection{\dtn Eigenfunctions on bounded Lipschitz domains}
We first show we can decompose a quasimode $u$ that $(h\Lambda_{M}^{\Omega}-\ohm)u=f$ into plane waves $\sum_{\pm}\tilde a_{\pm}(x; h)e^{\pm i\omega x/h}$ plus $u_e, u_0$, pieces at non-characteristic non-zero and zero frequencies respectively.  
\begin{lemma}[Eigenfunction structure on bounded Lipschitz domains]
\label{l:fullDescription} Let ${\Omega}\subset\mathbb{R}^{2}$
be a bounded uniformly Lipschitz domain. Fix $\omega\in[\frac12, 2]$. Let an $h$-dependent family of $u\in L^2(\partial{\Omega})$ with $\|u\|_{L^2}\le C$ satisfy $(h\Lambda_{M}^{\Omega}-\ohm)u=f$ with the global decomposition
\begin{equation}
    f=f_0+f_e+f_\omega, \ f_0\in \mc{D}'(\partial\Omega), \ f_e\in \mc{D}'_{c}(\partial\oreg), \ f_\omega \in C^\infty_{c}(\partial\oreg),
\end{equation}
where we ask
\begin{equation}
    \WFh f_e\cap (\{\xi=0\}\cup\{|\xi|=\omega\})=\WFh f_\omega\cap \{\xi=0\}=\emptyset.
\end{equation}
Let $\Gamma\subset\partial{\Omega}$ be open, smooth and simply connected. Let cutoffs $\gamma,\tilde\gamma\in C_c^{\infty}(\Gamma)$ that $\supp\gamma\cap\supp(1-\tilde\gamma)=\emptyset$ and both have connected support. Near $\supp\gamma$, we further ask $f_0\in S^0$ and and there exist $\tilde b_\pm(x;h)\in S^0(\Gamma)$ that
\begin{equation}
     f_\omega=h\sum_{\pm}\tilde b_{\pm}(x; h)e^{\pm i\omega x/h}+O_{S^0}(h^\infty)\text{ near }\supp\gamma. 
\end{equation}
Then there exist $\tilde a_{\pm}(x;h)\in S^0(\Gamma)$, $u_e(x;h)\in\mc{D}'(\Gamma)$, $u_0(x;h)\in S^0(\Gamma)$ such that 
\begin{equation}\label{e:AsymExp_Dim2}
    \gamma(x)u(x)= \gamma(x)\Big(\sum_{\pm}\tilde a_{\pm}(x; h)e^{\pm i\omega x/h}+u_e(x;h)+u_0(x;h)\Big)+O_{S^0}(h^\infty),
\end{equation}
in local arclength coordinates $x$ on $\Gamma$. Moreover, there exist $a_{\pm}(h)\in\mathbb{C}$ bounded in $h$, $a_{\pm,k}(x;h)\in S^0(\Gamma)$ that
\begin{equation}\label{e:7l51}
    \tilde a_\pm(x;h)\sim a_\pm(h)\pm i \int_{x_0}^{x}\tilde b_\pm(t;h)\,dt+h^2M^2\sum_{k=2}^\infty h^{k-2}a_{\pm,k}(x;h),
\end{equation}
where the values $\tilde a_\pm(x_0;h)$ at any point $x_0\in \supp \gamma$ uniquely determine $\tilde a_\pm(x;h)$ for all $x\in \supp\gamma$ modulo $O_{S^0}(h^\infty)$. For all $k,N\ge 0$, we can estimate
\begin{gather}\label{e:8l13}
    \|u_e\|_{H^k_h}\le C\|\tilde\gamma f_e\|_{H^{k-1}_h}+Ch^N\|u\|_{L^2}, \ \WFh(u_e)\subset \WFh(f_e), \\
    \|u_0\|_{H^k}\le C\|\tilde\gamma f_0\|_{S^0}+Ch\|u\|_{L^2}, \ \WFh(u_0)\subset \{\xi=0\}.\label{e:7l60}
\end{gather}
\end{lemma}

\begin{proof}
We devote Steps 1 to 4 to show the plane wave structure at the characteristic frequencies, Step 5 to the zero frequencies, and Step 6 to the non-characteristic non-zero frequencies. 

1. We make an a priori decomposition. Take $\tilde{\gamma}\in C_c^\infty(\Gamma)$, $\supp\gamma\cap \supp(1-\tilde\gamma)=\emptyset$, $\WF_h(f_0+f_e)\cap (\{|\xi|=\omega\}\cap \supp \tilde\gamma)=\emptyset$ and $\Omega$, $\tilde \Omega$ share the common smooth boundary near $\supp\tilde\gamma$ for some compact smooth domain $\tilde \Omega$. Assume $\partial_\nu\tilde\gamma=\partial_\nu \gamma=0$. Then by Lemma~\ref{l:DomChange} we have
\begin{equation}
    \gamma\Lambda^{\Omega}_M=\gamma\Lambda^{\tilde\Omega}_M\tilde\gamma +\gamma\partial_\nu\mathcal{R}_M^{\tilde\Omega}[-\Delta, \tilde\gamma]\mathcal{H}_M^{\Omega}.
\end{equation}
Note $\gamma$ and $[-\Delta,\tilde\gamma]$ have disjoint support 
to see for all $N>0$ that 
\begin{align}
\gamma\Lambda_{M}^{\Omega}&=\gamma B\tilde\gamma +\gamma\mc{E}, \ \gamma \mc{E}:H_\partial^{0}\rightarrow H^{N}_\partial, \label{eq:TransformationProperty}
\end{align}
where $B=\Lambda^{\tilde\Omega}_M  \in \Psi_{\textup{phg}}^1(\partial\tilde\Omega)$ is self-adjoint with the classical principal and subprincipal symbols $\sigma_1(B)=|\xi|$ and $\sigma_{0}(B)=0$. See \cite[Proposition C.1.]{Ta:23} for further details.

2a. Let $Q\in\Psi^0_h(\partial\tilde\Omega)$ to be fixed later with $\WF_h(Q)\cap (\{|\xi|=0\}\cup \supp (1-\tilde \gamma))=\emptyset$ and $\WF_h(Q)\cap \WF_h(f_0+f_e)=\emptyset$. We will prove that for any $\gamma,\tilde\gamma$ described in Step 1, for all $k\ge0$, all $A_{1},\ldots,A_{k}\in\Psi_h^{1}(\partial\tilde\Omega)$ with $\WFh(A_j)\subset\Ell_h(Q)$,  $\sigma_{h}(A_{j})|_{\xi=\pm\omega}=0$, 
\begin{equation}
\left\Vert A_{k}\cdots A_{1}\gamma u\right\Vert _{L^{2}} = O(h^{k}\|u\|_{L^2})+O(h^{k}\|\tilde b_\pm\|_{S^0}) =O(h^k).\label{eq:LagrangianCond}
\end{equation}
The nature of $f_\omega$ being a Lagrangian state implies
\begin{equation}
    \|A_{k}\cdots A_{1}\gamma f_\omega\|_{L^2}=O(h^{k+1}\|\tilde b_\pm\|_{S^0}). 
\end{equation}
To prove \eqref{eq:LagrangianCond}, we first recall that if $B\in \Psi^m$ and $A\in \Psi_h^l$ with $\WF_h(A)\cap \{\xi=0\}=\emptyset$, then $Ah^mB\in \Psi^{m+l}_h$ and $\sigma_h(Ah^mB)=\sigma_h(A)(x,\xi)\sigma(B)(x,\xi)$. 

2b. We proceed by induction on $k\geq0$. For $k=0$, \eqref{eq:LagrangianCond} is true by the assumption. Suppose that it holds up to the case $k-1$. Now, with $a_1=\sigma(A_1)$, we may write 
\begin{equation}
    a_1\gamma =q_0\gamma (\sqrt{\xi^2+h^2M^2}-\ohm)\tilde\gamma,\ Q_0=\Op_h(q_0)
\end{equation}
for some $q_0\in S^0$ with $\supp q_0\subset \Ell_h(Q)$, $\supp q_0\cap \supp (1-\tilde{\gamma})=\emptyset$. In particular, 
\begin{equation}
    A_1\gamma=Q_0\gamma(hB-\ohm)\tilde\gamma+hR\tilde\gamma,\ R\in \Psi^0_h, \ \WF_h(R)\subset \Ell_h(Q).
\end{equation}
Then we have
\begin{multline}
A_{k}\cdots A_1 \gamma u 
= A_k\cdots A_2 Q_0 \gamma (hB-\ohm) \tilde\gamma u+ hA_k\dots A_2 R\tilde\gamma u\\
= A_k\cdots A_2 Q_0\gamma(h\Lambda_M^{\Omega} -\ohm) u
-hA_k\cdots A_2 Q_0 \gamma\mc{E} u
+ hA_k\cdots A_2 R\tilde\gamma u+O_{\Psi^{-\infty}}(h^\infty)u.
\end{multline}
Noting $\WFh(Q_0)$ is supported away from $\{\xi=0\}$, we have for any $N>0$ large that
\begin{equation}
    \|A_k\dots A_2 Q_0 \gamma \mc{E} u\|_{L^2}\leq C\|Q_0 \gamma \mc{E} u\|_{H_h^{k+N-1}}\leq C \|(hD)^{k+N-1}\gamma \mc{E} u \|_{L^2}\leq Ch^{k+N-1}\|u\|_{L^2},\label{e:7l56}
\end{equation}
and thus by noting $A_k\cdots (A_2R)\tilde\gamma u=O_{L^2}(h^{k-1})$ from the case $k-1$ applied to $\tilde\gamma u$, we have
\begin{equation}
    \|A_k\cdots A_1\gamma u\|_{L^2}= O(h^{k}\|u\|_{L^2})+O(h^k\|\tilde b_\pm\|_{S^0})=O(h^k)
\end{equation}
as desired. The induction on $k$ then gives \eqref{eq:LagrangianCond}. 

3. We fix $\psi\in C^{\infty}(\mathbb{R},[0,1])$ with $0\notin\supp\psi$, $\psi\equiv1$ near $\pm\omega$, and $\WFh (f_e)\cap \supp\psi\cap \supp\gamma=\emptyset$. We ask $\tilde\psi$ to satisfy the same conditions for $\psi$ while $\supp\psi\cap\supp(1-\tilde\psi)=\emptyset$. Let $\varphi\in C_c^\infty(\Gamma)$ that $\supp\gamma\cap \supp(1-\varphi)=\supp\varphi\cap\supp(1-\tilde\gamma)=\emptyset$. Fixing $Q$ that $\WF_h(\psi(hD)\gamma)\subset\Ell_h(Q)$, the Lagrangian structure result \eqref{eq:LagrangianCond} applied to $\tilde\gamma u$ implies there are $\tilde a_{\pm}(x;h)\in S^{0}(\Gamma)$ that
\begin{equation}
    \psi(hD)\tilde\gamma(x)u(x)= \tilde\gamma(x)(\tilde a_{+}(x;h)e^{i\omega x/h}+\tilde a_{-}(x;h)e^{-i\omega x/h})+O_{C^\infty}(h^\infty). 
\end{equation}
See further details in Propositions 25.1.5 and 25.1.3 in \cite{HorV4}. It remains to show that on $\supp\gamma$, \eqref{e:7l51} holds for some $a_\pm(h)\in \C$ and some $a_{\pm,k}(x;h)\in S^0(\Gamma)$. Indeed, note 
\begin{multline}
    h\sum_{\pm}\tilde b_\pm \gamma e^{\pm i\omega x/h}\sim \psi(hD)\gamma(h\Lambda_M^\Omega-\ohm)u=\psi(hD)\gamma(hB-\ohm)\tilde\gamma u\\+h\psi(hD)\gamma\mc{E}u+O_{C^\infty}(h^\infty),
\end{multline}
where $h\psi(hD)\gamma\mc{E}u=O_{S^0}(h^\infty\|u\|_{L^2})$. On another hand, we claim 
\begin{multline}\label{e:7l67}
    \psi(hD)\gamma(hB-\ohm)\tilde\gamma u
    =\psi(hD)\gamma(hB-\ohm)\tilde\gamma (\tilde a_{+}(x;h)e^{i\omega x/h}+\tilde a_{-}(x;h)e^{-i\omega x/h})\\
    +O_{C^\infty}(h^\infty)\sim \gamma \sum_\pm\Big(\pm \frac{h}{i}\tilde a_\pm'(x;h)+h^2M^2\sum_{k=2}^\infty h^{k-1} L_{\pm, k}\tilde a_\pm(x;h)\Big)e^{\pm i\omega x/h},
\end{multline}
where $L_{\pm,k}$ are $h$-uniform differential operators in $\partial_x$ up to order $k$ with smooth, bounded coefficients. Indeed, from Lemma~\ref{l:1dnormalform} we have
\begin{equation}
    hB\psi(hD)=h\sqrt{D_x^2+M^2}\psi(hD)+O_{S^0}(h^3M^2)+O_{S^0}(h^\infty).
\end{equation}
Note $\WFh \tilde a_+\subset\{\xi=0\}$ and we have
\begin{equation}
    h\sqrt{D_x^2+M^2}\tilde a_+(x) e^{i\omega x/h}=(\mc{F}^{-1}(\sqrt{(h\xi+\omega)^2+h^2M^2}(\mc{F}\tilde a_+)(\xi)))e^{i\omega x/h},
\end{equation}
and near $h\xi=0$, we have
\begin{equation}
   \sqrt{(h\xi+\omega)^2+h^2M^2}-\sqrt{\omega^2+h^2M^2}=h\xi-\frac{h^3M^2\xi}{2\omega(h\xi+\omega)}+O(\frac{h^4M^4}{(h\xi+\omega)^3}).
\end{equation}
This yields \eqref{e:7l67} as claimed. Then we have
\begin{equation}\label{e:7l58}
    h\tilde b_\pm(x;h)\sim \pm \frac{h}{i}\tilde a_\pm'(x;h)+h^2M^2\sum_{k=2}^\infty h^{k-1} L_{\pm, k}\tilde a_\pm(x;h)
\end{equation}
The transport equation at size $h$ reads
\begin{equation}
    h\tilde b_\pm(x;h)\sim \pm \frac{h}{i}\tilde a_\pm'(x;h)+O(h^3M^2)
\end{equation}
and implies $\tilde a_{\pm}(x;h)=a_{\pm,0}(h)\pm i \int_{x_0}^{x}\tilde b_\pm(t;h)\,dt+h^2M^2 \tilde{a}_{\pm,2}(x;h)$ with $\tilde{a}_{\pm,2}(x;h), a_{\pm,0}(h)\in S^0$. Then at $h^3M^2$,
\begin{equation}
    0 \sim h^3M^2\big(\pm \frac1i\tilde{a}'_{\pm,2}(x;h)+L_{\pm,2}(a_{\pm,0}(h)\pm i \int_{x_0}^{x}\tilde b_\pm(t;h)\,dt)\big)+O(h^4M^2)
\end{equation}
implies $\tilde{a}_{\pm,2}'=\mp i L_{\pm,2}(a_{\pm,0}(h) \pm i\int_{x_0}^{x}\tilde b_\pm(t;h)\,dt)+O(h)$ and thus $\tilde a_{\pm}(x;h)=a_{\pm,0}(h)\pm i\int_{x_0}^{x}\tilde b_\pm(t;h)\,dt+h^2M^2 a_{\pm,2}(x;h)+h^3\tilde{a}_{\pm,3}(x;h)$ for $a_{\pm,2}, \tilde a_{\pm,3}\in S^0$. Solving recursively for $a_{\pm,k}(x;h)$ and renaming $a_{\pm, 0}(h)$ to $a_{\pm}(h)$ then proves \eqref{e:7l51}.

4. We show that the values of $\tilde a_{\pm}(x_0;h)$ at a point $x_0\in \supp \gamma$ determine $\tilde a_{\pm}(x;h)$ uniquely modulo $O(h^\infty)$. This is to show $\tilde a_{\pm}(x;h)=O(h^\infty)$, if we assume $\tilde a_{\pm}(x_0;h)=0$ for some $x_0\in\supp \gamma$ and $f=0$. Indeed, \eqref{e:7l58} and $\tilde a_\pm(x;h)=O_{S^0}(1)$ implies $\tilde a'_\pm(x;h)=\sum_{k=2}^\infty O(h) L_{\pm,k}\tilde a_\pm(x;h)=O_{S^0}(h)$, and $\tilde a_\pm(x_0;h)=0$ improves the size to $\tilde a_\pm(x;h)\in O_{S^0}(h)$. Now \eqref{e:7l58} implies $\tilde a'_\pm(x;h)=O_{S^0}(h^2)$ and one iterates to see $\tilde a_\pm(x;h)=O_{S^0}(h^\infty)$ as claimed. This concludes the description of the characteristic frequencies. 

5a. Now we turn to the zero frequencies. We show that there is $\e>0$ small such that for all $k\in\mathbb{N}_0$ and all $\psi_{0}\in C_{c}^{\infty}((-\e,\e);[0,1])$ with $\psi_0\equiv 1$ near $0$,
\begin{equation}
\left\Vert \psi_{0}(hD_{x})(hD_{x})^{k}\gamma u\right\Vert _{L^{2}}\le C h^{k}\big( \|\gamma f_0\|_{H^k}+h\|u\|_{L^2} \big).\label{eq:SemiclassIneq}
\end{equation}
The estimate \eqref{eq:SemiclassIneq} applied to $\tilde\gamma u$ instead of $\gamma u$ then gives the desired $u_0(x;h)=\psi_{0}(hD_{x})\tilde\gamma u+O_{S^0}(h^\infty)$. Without loss of generality we assume $\{|\xi|\le\e\}\cap \WFh(f_e)=\emptyset$, and let $\tilde\gamma, \tilde{\tilde\gamma}\in C^\infty(\Gamma)$ such that $\supp\gamma\cap \supp(1-\tilde\gamma)=\supp\tilde\gamma\cap \supp(1-\tilde{\tilde{\gamma}})=\emptyset$ and $\tilde\gamma f=\tilde\gamma(f_0+f_e+f_\omega)$. We start from the base case $k=0$, let $\tilde{\psi}_{0}\in C_{c}^{\infty}((-\e,\e),[0,1])$ with $\supp \psi_0\cap \supp (1-\tilde{\psi}_0)=\emptyset$. We again utilize \eqref{eq:TransformationProperty} that $\gamma \Lambda_M^{\Omega}=\gamma B \tilde\gamma +\gamma \mc{E}$ to see
\begin{multline}
    (\tilde{\psi}_{0}(hD)hB-\ohm)\psi_{0}(hD)\gamma u = \psi_{0}(hD)(hB-\ohm)\gamma u-\tilde{\psi}_0(hD)[hB,(1-\psi_0(hD))]\gamma u\\
= \psi_{0}(hD)\gamma(hB-\ohm)\tilde\gamma u+\psi_0(hD)[hB,\gamma]\tilde\gamma u+O_{\Psi^{0}_h }(h)\gamma u\\
= \psi_{0}(hD)\gamma(h\Lambda_M^{\Omega}-h\mc{E}-\ohm) u+O_{\Psi^{0}}(h) u
=  O_{L^2}(\|\gamma f_0\|_{L^2}+h\|u\|_{L^2}),
\end{multline}
Now for some $C>0$ independent of $h$, we have
\begin{equation}
    \|\tilde{\psi}_{0}(hD)hB\|_{L^2\to L^2}\leq h\|\tilde{\psi}_0(hD)\|_{H^{-1}\to L^2}\|B\|_{L^2\to H^{-1}}\leq C\e \|B\|_{L^2\to H^{-1}}.
\end{equation}
Now picking a small $\e$ such that $\|(\tilde{\psi}_{0}(hD)hB-\ohm)^{-1}\|$ is uniformly bounded in $h$ small gives
\begin{equation}
    \|\psi_0(hD)\gamma u\|_{L^2}\leq C(\|\gamma f_0\|_{L^2}+h\|u\|_{L^2})
\end{equation}
as desired. 

5b. Now we prove the inductive step. Assume \eqref{eq:SemiclassIneq} holds for $k\le K-1$ and consider the $k=K$ case. Note
\begin{multline}
    (\tilde{\psi}_0(hD)hB-\ohm)\psi_{0}(hD) (hD)^{K}\gamma u\\=(hD)^{K}(\tilde{\psi}_0(hD)hB-\ohm)\psi_{0}(hD)\gamma u-\sum_{j=0}^{K-1} (hD)^{j}\tilde{\psi}_0(hD)h[hD, B]\psi_{0}(hD)(hD)^{K-1-j}\gamma u\\
    =(hD)^{K}(\tilde{\psi}_0(hD)hB-\ohm)\psi_{0}(hD)\gamma u-\sum_{j=0}^{K-1}\sum_{l=0}^j\binom{j}{l}\tilde{\psi}_0(hD)h\ad^{j-l+1}_{hD}(B)\psi_{0}(hD)(hD)^{K-1-j+l}\gamma u\\=
    (hD)^{K}(\tilde{\psi}_0(hD)hB-\ohm)\psi_{0}(hD)\gamma u+O_{L^2}(h^{K}),
\end{multline}
the last line of which comes from $\ad^{k}_{hD}(B)\in h^k\Psi^1$ and the inductive hypothesis, that
\begin{equation}
    \|h\ad^{j-l+1}_{hD}(B)\psi_{0}(hD_{x})(hD)^{K-1-j+l}\gamma u\|\le Ch^{j-l+1}\|\psi_{0}(hD)(hD)^{K-1-j+l}\gamma u\|\le Ch^{K}.
\end{equation}
Now similar to the case $k=0$, we have 
\begin{multline}\label{e:7l1}
    (\tilde{\psi}_0(hD)hB-\ohm)\psi_{0}(hD) (hD)^{K}\gamma u=(hD)^{K}\psi_0(hD)\gamma(h\Lambda^{\Omega}_M-\ohm)u\\+(hD)^{K}\psi_0(hD)h[B,\gamma]u
    -(hD)^{K}\tilde\psi_0(hD)[hB,1-\psi_0(hD)]\gamma u+(hD)^{K}\psi_0(hD)\gamma h\mc{E} u +O_{L^2}(h^{K})\\
    =(hD)^{K}\psi_0(hD)\gamma(h\Lambda^{\Omega}_M-\ohm)u+(hD)^{K}\psi_0(hD)h[B,\gamma]u\\
    -(hD)^{K}\tilde\psi_0(hD)[hB,1-\psi_0(hD)]\gamma u+O_{L^2}(h^{K}).
\end{multline}
That $\gamma f_0\in S^0$ implies the first term on the right is of size $h^{K}$: 
\begin{equation}
    \|(hD)^{K}\psi_0(hD)\gamma(h\Lambda^{\Omega}_M-\ohm)u\|\le Ch^{K}\|\gamma f_0\|_{H^{K}}+Ch^N\|f_e\|_{H^{-N}_h}+Ch^N\|f_\omega\|_{H^{-N}_h},
\end{equation}
and we will estimate the other terms. Note
\begin{multline}
    (hD)^{K}\psi_0(hD)[hB,\gamma]u=\sum_{l=0}^{K-1}\binom{K-1}{l}hD \psi_0(hD)\ad_{hD}^{K-1-l}([hB,\gamma])(hD)^{l}u\\
    =\sum_{l=0}^{K-1}\binom{K-1}{l}hD \psi_0(hD)\ad_{hD}^{K-1-l}([hB,\gamma])\Big(\tilde\psi_0(hD)(hD)^{l}\tilde\gamma u+(hD)^{l}(1-\tilde\gamma)u\\+(1-\tilde\psi_0(hD))(hD)^{l}\tilde\gamma u\Big).
\end{multline}
Here for the first term, by the inductive hypothesis and that $\ad_D^{K-1-l}[B, \gamma]\in \Psi^0$, we have
\begin{equation}
    \|\ad_{hD}^{K-1-l}([hB,\gamma])\tilde\psi_0(hD)(hD)^{l}\tilde\gamma u\|\le C h^{K-l}\|\tilde\psi_0(hD)(hD)^{l}\tilde\gamma u \|\le C h^{K}(\|\tilde{\tilde{\gamma}} f_0\|_{H^{l}}+h\|u\|). 
\end{equation}
For the second term, since $\gamma$ and $1-\tilde\gamma$ have disjoint support, $D\ad_{D}^{K-1-l}[B, \gamma](1-\tilde\gamma)\in \Psi^{-\infty}$ and 
\begin{equation}
\|h D\ad_{hD}^{K-1-l}([hB,\gamma])(hD)^{l}(1-\tilde\gamma)u\|\le C h^{K+1}\|u\|_{L^2}. 
\end{equation}
For the last term, note that $\WF_h(1-\tilde\psi_0(hD))(hD)^{l}\tilde\gamma$ is disjoint from $\{\xi=0\}$ and $\supp \psi_0$ to see
\begin{equation}
    \|\psi_0(hD)\ad_{hD}^{K-1-l}([hB,\gamma])(1-\tilde\psi_0(hD))(hD)^{l}\tilde\gamma u\|=O( h^\infty)\|u\|_{L^2}. 
\end{equation}
We then obtain in \eqref{e:7l1}
\begin{equation}
    (hD_{x})^{K}\psi_0(hD)[hB,\gamma]u=O_{L^2}(h^{K}). 
\end{equation}
We move on to estimate $(hD)^{K}\tilde\psi_0(hD)[hB,1-\psi_0(hD)]\gamma u$. Let $\tilde{\tilde{\psi}}_0$ be supported on $(-\epsilon, \epsilon)$ that $\supp\tilde\psi_0\cap\supp(1-\tilde{\tilde{\psi}}_0)=\emptyset$. On $\supp\gamma$, $[hB, 1-\psi_0(hD)]\in h\Psi^0_h$ and $\WFh(\tilde\psi_0(hD)[hB, 1-\psi_0(hD)])$ is supported inside $\{\tilde{\tilde{\psi}}_0\equiv 1\}$. Thus by the elliptic estimate and the inductive hypothesis we have
\begin{align}
    \|(hD)^{K}\tilde\psi_0(hD)[hB,1-\psi_0(hD)]\gamma u\|
    \le & Ch\|hD\tilde{\tilde{\psi}}_0(hD) (hD)^{K-1}\gamma u\|+Ch^N\|\gamma u\|_{H_h^{-N}} \nonumber \\
    \le & C h^{K}(\|\gamma f_0\|_{H^{K-1}}+h\|u\|_{L^2}). 
\end{align}
Revisit \eqref{e:7l1} to see
\begin{equation}
    \|\psi_0(hD)(hD)^K\gamma u\|\le C\|(\tilde{\psi}_0(hD)hB-\ohm)\psi_{0}(hD) (hD)^{K}\gamma u\|
    \le C h^{K}\big( \|\gamma f_0\|_{H^K}+h\|u\|_{L^2} \big),
\end{equation}
which verifies the inductive case $k=K$ of \eqref{eq:SemiclassIneq} and thus gives \eqref{e:7l60}. 

6. It remains to invert $f_e$ at non-characteristic non-zero frequencies. Let $\psi_e=1-\psi_0-\psi$. Let $\tilde\psi_e\in C^\infty(\mathbb{R}_\xi,[0,1])$ that $\supp\psi_e\cap\supp(1-\tilde\psi_e)=\emptyset$ and $0\notin \supp\tilde\psi_e$. Let $u_e=\tilde\gamma\psi_e(hD)\gamma_1 u$. Let $\gamma_1,\gamma_2\in C^\infty_c(\Gamma)$ with connected support and $\supp\tilde\gamma\cap\supp(1-\gamma_1)=\supp\gamma_1\cap\supp(1-\gamma_2)=\emptyset$. Noting $\WFh u_e\subset\Ell_h(\tilde \psi_e(hD)\gamma_1(hB-\ohm)\gamma_2)$ we have for large $N>0$, 
\begin{equation}
    \|u_e\|_{H^k_h}\le C\|\tilde \psi_e(hD)\gamma_1(hB-\ohm)\gamma_2 u\|_{H^{k-1}_h}+Ch^N\|\gamma_2 u\|_{H^{-N}_h}.
\end{equation}
Note $\psi_e\equiv 1$ over $\WFh(f_e)$. We have $\gamma_1 \Lambda^{\Omega}_M=\gamma_1 B\gamma_2+\gamma_1\mc{E}$ and for large $N>0$, 
\begin{equation}
    \tilde\psi_e(hD)\gamma_1 (hB-\ohm)\gamma_2 u=\tilde\psi_e(hD)\gamma_1 f-h\tilde\psi_e(hD)\gamma_1\mc{E} u=\gamma_1 f_e+O_{S^0}(h^N\|u\|_{L^2}+h^N\|f_e\|_{H^{-N}_h}),
\end{equation}
since $\tilde\psi_e$ is supported away from 0, similar to \eqref{e:7l56}. This gives
\begin{equation}
    \|u_e\|_{H^k_h}\le C\|\gamma_1 f_e\|_{H^{k-1}_h}+Ch^N\|u\|_{L^2}.
\end{equation}
Rename $\gamma_1$ to $\tilde\gamma$ to see \eqref{e:8l13}.
\end{proof}

\subsection{Blow-up near the corners of curvilinear polygons}\label{s:8-2}
Let $\omega\in\mathbb{R}$, $\omega>0$, and $z\in\mathbb{C}$. We first focus on one connected component $\Gamma^l=\cup_{j=1}^{J}\Gamma_j$ of $\partial\mc{P}$, where $\Gamma_j$ connects $x_{j}$ to $x_{j+1}$. We abuse the notation and treat $x_{0}=x_{J}$, $x_{J+1}=x_1$ cyclically, and similarly for $\Gamma_0$ and $\Gamma_{J+1}$. 

Near each corner $x_j$, we build a \acl domain $\mc{W}_j\subset\mathbb{R}^2$ that $\mc{W}_j$ and $\mc{P}$ share the common boundary near $x_j$. Indeed, for small $\e>0$, we truncate $\partial\mc{P}\cap\{\dist(x, x_j)<\e\}$ and extend the two truncated endpoint tangentially by two straight lines to obtain $\partial\mc{W}_j$, and thus $\mc{W}_j$. Such an $\mc{W}_j$ has corner angle $\alpha_j$ at $x_j$ and is bi-Lipschitz with small $\e$-dependent error to $W_{\alpha_j}$, a straight wedge of corner angle $\alpha_j$. Lemma~\ref{l:DomComp} then implies for each small $\delta>0$, there is $\e>0$ such that $\rardy^0_{\mc{W}_j}>\frac{\pi}{\alpha_j}-\delta$. We fix an $\mc{W}_j$ with sufficiently small $\e>0$ and abuse the notation to say $\rardy^0_{\mc{W}_j}>\frac{\pi}{\alpha_j}-0$, through Lemma~\ref{l:conerardy}. See Figure~\ref{f:8l1}. 

\begin{figure}[t]
  \centering
        \begin{tikzpicture}[x=1cm,y=1cm]
  \usetikzlibrary{calc,plotmarks}
  \definecolor{colorAKANE}{RGB}{183,40,46}
  \definecolor{colorUSUAI}{RGB}{40,148,200}
  \colorlet{plotPrimary}{colorAKANE}
  \colorlet{plotTertiary}{colorUSUAI}
  \colorlet{plotAxis}{black!72}
  \colorlet{plotGuide}{black!34}
  \tikzset{
    >=latex,
    font=\small,
    mark size=1.5pt,
    mark options={thin},
    line cap=round,
    line join=round,
    boundary edge/.style={draw=black!88,line width=0.7pt},
    arclength axis/.style={draw=plotAxis!75,line width=0.4pt},
    origin guide/.style={draw=plotGuide,line width=0.3pt,
      dash pattern=on 1pt off 2pt},
    wave arrow/.style={line width=0.65pt,->},
    incoming wave/.style={wave arrow,draw=plotTertiary,text=plotTertiary},
    outgoing wave/.style={wave arrow,draw=plotPrimary,text=plotPrimary},
    figure label/.style={font=\fontsize{9}{10.8}\selectfont,
      text=black,inner sep=1pt},
    wave label/.style={font=\fontsize{8}{9.6}\selectfont,inner sep=1pt},
    axis label/.style={font=\fontsize{8}{9.6}\selectfont,
      text=plotAxis,inner sep=1pt}
  }

  \path[use as bounding box] (0,0) rectangle (7.5,3.2);
  \begin{pgfinterruptboundingbox}
  \begin{scope}[shift={(0.10,-0.175)}]
    \coordinate (vertex j) at (3.90,1.60);
    \coordinate (vertex next) at (6.30,1.60);
    \coordinate (left extent) at (3.13,2.93);
    \coordinate (right extent) at (7.07,2.93);
    \coordinate (axis origin) at ($(vertex j)+(0,-0.85)$);
    \coordinate (axis left limit) at (0.64,0);
    \coordinate (axis right limit) at (6.86,0);
    \coordinate (axis left) at (axis left limit |- axis origin);
    \coordinate (axis right) at (axis right limit |- axis origin);
    \coordinate (axis next) at (vertex next |- axis left);
    \coordinate (edge center) at ($(vertex j)!0.5!(vertex next)$);
    \coordinate (edge outgoing center) at ($(edge center)+(0,0.32)$);
    \coordinate (edge incoming center) at ($(edge center)+(0,-0.32)$);

    \coordinate (left flow upper) at ($(left extent)!0.24!(vertex j)$);
    \coordinate (left flow lower) at ($(left extent)!0.80!(vertex j)$);
    \coordinate (left normal) at
      ($(left flow upper)!0.27cm!90:(left flow lower)-(left flow upper)$);
    \coordinate (left incoming upper) at ($(left flow upper)+(left normal)$);
    \coordinate (left incoming lower) at ($(left flow lower)+(left normal)$);
    \coordinate (left outgoing upper) at ($(left flow upper)-(left normal)$);
    \coordinate (left outgoing lower) at ($(left flow lower)-(left normal)$);

    \draw[origin guide] (vertex j) -- (axis origin);
    \draw[origin guide] (vertex next) -- (axis next);

    \begin{scope}
      \draw[boundary edge]
        (left extent) -- (vertex j) -- (vertex next) -- (right extent);
      \path[fill=black,draw=none]
        (vertex j) circle[radius=1.25pt]
        (vertex next) circle[radius=1.25pt];
      \node[figure label,anchor=east] at ($(left extent)+(-0.10,0.04)$)
        {$\Gamma_{j-1}$};
      \node[figure label,anchor=south]
        at ($(vertex j)!0.90!(vertex next)+(0,0.08)$)
        {$\Gamma_j$};
      \node[figure label,anchor=east]
        at ($(vertex next)!0.50!(right extent)+(-0.10,0)$)
        {$\Gamma_{j+1}$};
      \node[figure label,anchor=north east] at ($(vertex j)+(-0.10,-0.10)$)
        {$x_j$};
      \node[figure label,anchor=north west] at ($(vertex next)+(0.10,-0.10)$)
        {$x_{j+1}$};
    \end{scope}

    \begin{scope}
      \draw[incoming wave] (left incoming upper) -- (left incoming lower)
        node[wave label,pos=0,right=4pt]
          {$a_{+,j-1}^{i}\,e^{iz(x+L_{j-1})/h}$};
      \draw[outgoing wave] (left outgoing lower) -- (left outgoing upper)
        node[wave label,midway,left=4pt]
          {$a_{-,j-1}^{o}\,e^{-iz(x+L_{j-1})/h}$};
    \end{scope}

    \begin{scope}
      \draw[outgoing wave]
        ($(edge outgoing center)+(-0.70,0)$) -- ($(edge outgoing center)+(0.70,0)$)
        node[wave label,midway,above=1.5pt]
          {$a_{+,j}^{o}\,e^{izx/h}$};
      \draw[incoming wave]
        ($(edge incoming center)+(0.70,0)$) -- ($(edge incoming center)+(-0.70,0)$)
        node[wave label,midway,below=1pt]
          {$a_{-,j}^{i}\,e^{-izx/h}$};
    \end{scope}

    \begin{scope}
      \draw[arclength axis,<->] (axis left) -- (axis right);
      \draw[arclength axis]
        ($(axis origin)+(0,-0.055)$) -- ($(axis origin)+(0,0.055)$);
      \draw[arclength axis]
        ($(axis next)+(0,-0.055)$) -- ($(axis next)+(0,0.055)$);
      \node[axis label,anchor=north] at ($(axis origin)+(0,-0.07)$)
        {$x=0$};
      \node[axis label,anchor=north] at ($(axis next)+(0,-0.07)$)
        {$x=L_j$};
      \node[axis label,anchor=north]
        at ($(axis left)!0.35!(axis origin)+(0,-0.07)$) {$x<0$};
      \node[axis label,anchor=west] at ($(axis right)+(0.08,0)$) {$x$};
    \end{scope}
  \end{scope}
  \end{pgfinterruptboundingbox}
\end{tikzpicture}%
  \caption{Extension of corner $x_j$ of $\mc{P}$ into curved wedge $\mc{W}_j$}
  \label{f:8l1}
\end{figure}

Now let $\chi_j,\tilde{\chi}_j\in C_c^\infty(\cl\mc{P}\cap\cl\mc{W}_j)$ with the vertex $x_j\notin \supp (1-\chi_j)$ and $\chi_j\prec \tilde\chi_j$. Parametrize $\partial\mc{W}_j=(-\infty, 0)\cup\{0\}\cup(0, \infty)$ by arclength, where $x>0$ on $\Gamma_j\cap\partial\mc{W}_j$ and $x< 0$ on $\Gamma_{j-1}\cap\partial\mc{W}_j$. We extend $(0, \infty)$ tangentially by a straight line on the left to obtain the eventually conic smooth domain $\mc{W}_{j,R}\subset\mathbb{R}^2$ that $\mc{W}_{j,R}$, $\mc{W}_j$ share the common smooth boundary near $[\frac\e2,\infty)$. We similarly obtain $\mc{W}_{j,L}\subset\mathbb{R}^2$ that $\mc{W}_{j,L}$, $\mc{W}_j$ share the common smooth boundary near $(-\infty,-\frac\e2]$. Similar to Step 1a of Lemma~\ref{l:scstraighten}, Lemma~\ref{l:2dquasimodeappendix} then yields interior quasimodes $\tilde e_{\pm,j,L/R,h}(x;h,\omega)$ on $\mc{W}_{j, L/R}$ that
\begin{equation}
    (-\Delta+M^2)\tilde e_{\pm,j, L/R, h}=O_{\mc{S}(\mc{W}_{j, L/R})}(h^\infty), \ (h\partial_\nu-\om)\tilde e_{\pm, j, L/R,h}=O_{\mc{S}(\partial\mc{W}_{j, L/R})}(h^\infty),
\end{equation}
whose boundary restriction $e_{\pm,j,L/R,h}(x;h, \omega, M)=\tilde e_{\pm,j,L/R,h}(0, x;h, \omega, M)$ satisfies
\begin{equation}
    e_{\pm,j,L/R,h}(x;h,\omega, M)=e^{\pm i\omega x/h}(1+O_{S^0(\partial\mc{W}_{j, L/R})}(h^2M^2)).
\end{equation}

\begin{lemma}[Blow-up near the corner]\label{l:blowup}
Let $\Omega=\mc{P}\subset\mathbb{R}^{2}$ be a curvilinear polygon with angles $\{\alpha_j\}_{j=1}^J$ and edges $\{\Gamma_j\}_{j=1}^J$. Let $f,f_0, f_e, f_\omega$ be as defined in Lemma~\ref{l:aprioriSteklovDescription}, and $u$ be that $(h\Lambda_{M}^{\Omega}-\ohm)u=f$. Let $\gamma_j(x)\in C_c^\infty(\Gamma_j)$ with the vertices disjoint from $\supp\gamma_j$. Let $\chi_j\in C^\infty(\mc{P})$ that $x_j\notin \supp(1-\chi_j)$ with $\supp(1-\gamma_j)\cap \Gamma_j\cap\supp\chi_j\cup\supp\chi_{j+1}$. For each $j$, Lemma~\ref{l:aprioriSteklovDescription} gives
\begin{equation}\label{e:7l5b}
\gamma_j u= \gamma_j\Big(\sum_{\pm}\tilde a_{\pm,j}(x;h) e^{\pm i\omega x/h}+u_e(x;h)+u_0(x;h)\Big) +O_{S^0}(h^\infty).
\end{equation}
It allows us to read off $a_{\pm, j}^{o/i}(h)\in \mathbb{C}$ from $\gamma_ju$ from the right of $x_j$, by
\begin{gather}\label{e:9l9b}
    \tilde a_{+,j}(x;h)\sim a_{+,j}^o(h)+O_{S^0}(h^2M^2\|u\|_{L^2}), \text{ near }\supp\gamma_j,\\ \label{e:9l10b}
    \tilde a_{-,j}(x;h)\sim a_{-,j}^i(h)+O_{S^0}(h^2M^2\|u\|_{L^2}), \text{ near }\supp\gamma_j, 
\end{gather}
and also to read off $a_{\mp, j}^{o/i}(h)\in \mathbb{C}$ from the left of $x_j$, by
\begin{gather}
    \tilde a_{+,j-1}(x;h)\sim a_{+,j-1}^i(h)+O_{S^0}(h^2M^2\|u\|_{L^2}), \text{ near }\supp\gamma_{j-1},\\ 
    \tilde a_{-,j-1}(x;h)\sim a_{-,j-1}^o(h)+O_{S^0}(h^2M^2\|u\|_{L^2}), \text{ near }\supp\gamma_{j-1}. 
\end{gather}
Then we have
\begin{multline}\label{e:8l14}
    \begin{pmatrix}a_{-,j-1}^o\\a_{+,j}^o\end{pmatrix}=\Big(\Sigma_{M, h^{-1}\omega}^{\mc{W}_j}+O(M^2h^{2-0})\Big)\begin{pmatrix}a_{+,j-1}^i\\a_{-,j}^i\end{pmatrix}+O(h^{\frac{\pi}{\alpha_j}-1-0}\|u_0\|_{H^{\frac12}_h})+O(h^\infty)\\
    +O(h^{-1-0}\|\chi_j f\|_{H^{-\frac12}_h}),
\end{multline}
and
\begin{equation}\label{e:8l11}
    \|\chi_j u\|_{H^{\frac12}_h}\le C|a^i|+h^{\frac{\pi}{\alpha_j}-1-0}\|u_0\|_{H^{\frac12}_h}+h^{-1-0}\|\chi_j f\|_{H^{-\frac12}_h}+Ch^N\|u\|_{H^{\frac12}_h}.
\end{equation}
\end{lemma}
\begin{remark}
    Lemma~\ref{l:aprioriSteklovDescription} below improves the remainder $h^{\frac{\pi}{\alpha_j}-1-0}u_0$ in \eqref{e:8l14} to $h^{\frac{2\pi}{\alpha_+}-0}u$:
    \begin{multline}
    \begin{pmatrix}a_{-,j-1}^o\\a_{+,j}^o\end{pmatrix}=\Big(\Sigma_{M, h^{-1}\omega}^{\mc{W}_j}+O(M^2h^{2-0})\Big)\begin{pmatrix}a_{+,j-1}^i\\a_{-,j}^i\end{pmatrix}+O(h^{\frac{2\pi}{\alpha_+}-0}\|u\|_{L^2})+O(h^\infty)\\
    +O\Big(\sum_{j=1}^{J}h^{\frac{2\pi}{\alpha_+}-1-0}\|\chi_j f\|_{H^{-\frac12}_h}+h^{\frac{\pi}{\alpha_+}-1-0}\|\tilde\gamma f_0\|_{S^0}\Big),
    \end{multline}
    where $\alpha_+=\max_j(\alpha_j)$. 
\end{remark}
\begin{proof}
    To prove the lemma, we write  $\mc{P}^h=h^{-1}\mc{P}$ and define the unitary rescaling operator
    \begin{equation}
        T_h:L^2(\partial\mc{P})\to L^2(\partial\mc{P}^h), \ T_hu(x)=h^{\frac12}u(hx).
    \end{equation}
    Then under the rescaling we have
    \begin{equation}\label{e:7l73}
        (\Lambda_{hM}^{\mc{P}^h}-\sqrt{\omega^2+h^2M^2})T_hu=T_hf.
    \end{equation}
    We adopt the superscript rescaling notations in Section~\ref{s:scatteringatrixConvergence}. We now work near each corner $x_j$. Let $\tilde\chi_j, \tilde{\tilde{\chi_j}}$ satisfy the requirement for $\chi_j$ with $\chi_j\prec\tilde\chi_j\prec \tilde{\tilde{\chi_j}}$. As in Step 1 of Lemma~\ref{l:scstraighten}, we construct the principal oscillator by
\begin{equation}
    \tilde u_{\omega}(x_0, x)=\tilde{\tilde{\chi_j}}\Big(a_{+}^{o}\chi_R^h\tilde e_{+, R,h}+a_{+}^{i}\chi_L^h\tilde e_{+, L,h}+a_{-}^{i}\chi_R^h\tilde e_{-, R}+a_{-}^{o}\chi_L^h\tilde e_{-, L}\Big).
\end{equation}
    as in \eqref{e:7l71}, and the zero-frequency residue by $\tilde u_0=\tilde{\tilde{\chi_j}}\bar u-\tilde{u}_\omega\in X_m^l(\mc{W}_{j})$. We then construct 
    \begin{equation}\label{e:9l11}
        v'_j=\chi_j^{\delta}T_hu+(\chi_j-\chi_j^{\delta})T_hu_0+(1-\chi_j^\delta)(\chi_R^\delta w_R+\chi_L^\delta w_L)
    \end{equation}
    on $\mc{W}_j$ with
    \begin{equation}
        w_R=h^{\frac12}\mathbbm{1}_{\Gamma_j^h}(a_{+,j}^oe_{+, j}^R+a_{-,j}^ie_{-, j}^R), \ w_L=h^{\frac12}\mathbbm{1}_{\Gamma_{j-1}^h}(a_{+,j-1}^oe_{+, j-1}^L+a_{-,j-1}^ie_{-, j-1}^L).
    \end{equation}
We can then identify \eqref{e:9l11} with the setting \eqref{e:7l28} via
\begin{equation}
    \mc{W}_\alpha=\tilde{\mc{W}}_\alpha=\mc{W}_j, \ \tilde{\mc{W}}_\alpha^R=\mc{W}_j^R, \ \tilde{\mc{W}}_\alpha^L=\mc{W}_j^L, \ \chi=\chi_j, \ F=\Id, \ u=T_h u, \ u_0=T_h u_0,\ v'=v_j'.
\end{equation}
Equation~\eqref{e:7l29} now reads
\begin{multline}
    (\Lambda^{\mc{W}^h_{j}}_{hM}-\ohm)v_j'=\chi_j T_h f-\partial_\nu\mc{R}^{\mc{W}^h_{j}}_{hM}(\chi_j-\chi_j^\delta)(-\Delta+h^2M^2)T_h\tilde u_0-\partial_\nu\mc{R}^{\mc{W}^h_{j}}_{hM}\chi_j^\delta(-\Delta+h^2M^2) T_h\bar u\\
    -\partial_\nu\mc{R}^{\mc{W}^h_{j}}_{hM}[-\Delta,\chi_j]T_h\tilde u_0-\sum_{L/R} \partial_\nu \mathcal{R}^{{\mc{W}^h_{j}}}_{hM}(1-\chi_j^\delta)[-\Delta, \chi_{L/R}^{\delta}]\bar w_{L/R}\\    +\partial_\nu\mc{R}^{\mc{W}^h_{j}}_{hM}[-\Delta,\chi_j^{\delta}]\Big(-T_h\tilde u_\omega+\sum_{L/R} \chi_{{L/R}}^\delta\bar w_{L/R}\Big) + O(h^\infty \mc{S}).
\end{multline}
Steps 2b, 2c, 2d, 2e and 3 of Lemma~\ref{l:scstraighten} remain the same:
\begin{gather}
    \|\partial_\nu\mc{R}^{\mc{W}^h_{j}}_{hM}(\chi_j-\chi_j^\delta)(-\Delta+h^2M^2)T_h\tilde u_0\|_{Z^{-\frac12, N,-\frac12}_{-L_0+0, N, N}}\le Ch^{\frac{2\pi}{\alpha_j}-0}\|\tilde\chi T_h u_0\|_{H^{\frac12, \frac{\pi}{\alpha_j}+\frac12-0}},\\
    \|\partial_\nu\mc{R}^{\mc{W}^h_{j}}_{hM}\chi_j^\delta(-\Delta+h^2M^2) T_h\bar u\|_{Z^{-\frac12, N,-\frac12}_{-L_0+0, N, N}}\le Ch^N\|\tilde\chi^\delta T_h u\|_{H^{\frac12}},\\
    \|\partial_\nu\mc{R}^{\mc{W}^h_{j}}_{hM}[-\Delta,\chi_j]T_h\tilde u_0\|_{Z^{-\frac12, N,-\frac12}_{-L_0+0, N, N}}\le Ch^{\frac{2\pi}{\alpha_j}-0}\|\tilde\chi T_h u_0\|_{H^{\frac12, \frac{\pi}{\alpha_j}+\frac12-0}},\\
    \|\partial_\nu \mathcal{R}^{{\mc{W}^h_{j}}}_{hM}(1-\chi_j^\delta)[-\Delta, \chi_{L/R}^{\delta}]\bar w_{L,R}\|_{Z^{-\frac12, N,-\frac12}_{-L_0+0, N, N}}\le C h^N(|a^i|+|a^o|),\\
    \|\partial_\nu\mc{R}^{\mc{W}^h_{j}}_{hM}[-\Delta,\chi_j^{\delta}]\Big(-T_h\tilde u_\omega+\sum_{L/R} \chi_{L/R}^\delta\bar w_{L/R}\Big)\|_{Z^{-\frac12, N,-\frac12}_{-L_0+0, N, N}}\le Ch^{\frac12}(M^2h^{2-0}+h^N)(|a^i|+|a^o|). 
\end{gather}
Note unitary $T_h: L^2(\partial\mc{W}_j)\rightarrow L^2(\partial\mc{W}_j^h)$ is also an isometry from $H^{\frac12}_h(\partial\mc{W}_j)$ to $H^{\frac12}(\partial\mc{W}_j^h)$ for
\begin{equation}
    h^{\frac{2\pi}{\alpha_j}-0}\|\tilde\chi_j T_h u_0\|_{H^{\frac12, \frac{\pi}{\alpha_j}+\frac12-0}}\le Ch^{\frac{\pi}{\alpha_j}-\frac12-0}\|\tilde\chi_j T_h u_0\|_{H^{\frac12}}\le Ch^{\frac{\pi}{\alpha_j}-\frac12-0}\|u_0\|_{H^{\frac12}_h}.
\end{equation}
On the other hand, we have
\begin{equation}
    \|\chi_j T_h f\|_{Z^{-\frac12, \frac12+0,-\frac12}_{\frac12+0, \frac12+0, \frac12+0}}\le Ch^{-\frac12-0}\|\chi_j f\|_{H^{-\frac12}_h}.
\end{equation}
The above implies
\begin{multline}
    \|(\Lambda^{\mc{W}^h_{j}}_{hM}-\ohm)v_j'\|_{Z^{-\frac12, \frac12+0,-\frac12}_{-L_0+0, \frac12+0, \frac12+0}}\le Ch^{\frac{\pi}{\alpha_j}-\frac12-0}\|u_0\|_{H^{\frac12}_h}+Ch^{\frac12}(M^2h^{2-0}+h^N)(|a^i|+|a^o|)\\
    +Ch^N\|u\|_{H^{\frac12}_h}+Ch^{-\frac12-0}\|\chi_j f\|_{H^{-\frac12}_h}.
\end{multline}
We correct the small eigenerror by building generalized eigenfunctions
\begin{multline}
    v_j=v_j'-(\Lambda^{\mc{W}^h_{j}}_{hM}-\ohm-i0)^{-1}(\Lambda^{\mc{W}^h_{j}}_{hM}-\ohm)v'\\
    =v_j'+O_{Z^{\frac12, \frac12+0,\frac12}_{-L_0+0, -\frac12+0, -\frac12-0}}\Big(h^{\frac{\pi}{\alpha_j}-\frac12-0}\|u_0\|_{H^{\frac12}_h}+h^{\frac12}(M^2h^{2-0}+h^N)(|a^i|+|a^o|)+h^N\|u\|_{H^{\frac12}_h}\\
    +O(h^{-\frac12-0}\|\chi_j f\|_{H^{-\frac12}_h})\Big).
\end{multline}
We here have $(\Lambda^{\mc{W}^h_{j}}_{hM}-\ohm) v_j=0$ and the incoming scattering amplitudes of $v$ must agree with $h^{\frac12}a^i$. The scattering matrix result Proposition~\ref{p:existscmatrix} implies
\begin{equation}
     \Sigma_{M, h^{-1}\omega}^{\mc{W}_j}a^i = a^o+O(M^2h^{2-0})(|a^i|+|a^o|))+O(h^{\frac{\pi}{\alpha_j}-1-0}\|u_0\|_{H^{\frac12}_h})+O(h^\infty)+O(h^{-1-0}\|\chi_j f\|_{H^{-\frac12}_h}),
\end{equation} 
and
\begin{equation}
    a^o=\Big(\Sigma_{M, h^{-1}\omega}^{\mc{W}_j}+O(M^2h^{2-0})\Big)a^i+O(h^{\frac{\pi}{\alpha_j}-1-0}\|u_0\|_{H^{\frac12}_h})+O(h^{-1-0}\|\chi_j f\|_{H^{-\frac12}_h}),
\end{equation}
as desired. Moreover we have
\begin{equation}
    v_j\sim h^{\frac12}a^ie^{-i\omega|x|}+h^{\frac12}a^oe^{-i\omega|x|}+O_{Z^{N, N,1-0}_{L_0-0,N, N}}(h^{\frac12}|a^i|),
\end{equation}
and thus
\begin{gather}
    \|v_j\|_{Z^{\frac12, 0,\frac12}_{L_0-0, -\frac12-0, -\frac12-0}}\le Ch^{\frac12}|a^i|,\\
    \chi_j T_h^* v_j\sim a^ie^{-ih^{-1}\omega|x|}+a^oe^{-ih^{-1}\omega|x|}+O_{H^{1-0}_h}(h^{\frac12}|a^i|).
\end{gather}
This implies
\begin{gather}\label{e:7l74}
    \|v_j'\|_{Z^{\frac12, 0,\frac12}_{-L_0+0, -\frac12-0, -\frac12-0}}\le Ch^{\frac{\pi}{\alpha_j}-\frac12-0}\|u_0\|_{H^{\frac12}_h}+Ch^{\frac12}|a^i|+Ch^N\|u\|_{H^{\frac12}_h}+O(h^{-\frac12-0}\|\chi f\|_{H^{-\frac12}_h}),\\
    \chi_j T_h^* v_j'\sim a^ie^{-ih^{-1}\omega|x|}+a^oe^{-ih^{-1}\omega|x|}+O_{H^{\frac12}_h}(h^{\frac12}|a^i|+h^{0-0}\|u_0\|_{H^{\frac12}_h}+h^{-\frac12-0}\|\chi_j f\|_{H^{-\frac12}_h}+h^\infty).
\end{gather}
We record \eqref{e:7l74} for later use. Now note $\chi_j u=\chi_j T_h^* v_j'+O(M^2h^{2-0})$ to conclude \eqref{e:8l11}. 
\end{proof}

\subsection{\dtn eigenfunctions on curvilinear polygons}
Our final goal in this section is to show, on a curvilinear polygon $\mc{P}$, the zero-frequency remainder $u_0(x;h)$ from Lemma~\ref{l:fullDescription} is actually smaller than one a priori expects, where the improvement depends on the upper bound of angles of $\mc{P}$. In this section, we use the results in Section~\ref{s:harm} and Section~\ref{s:3} based on $X_m^l(\mc{P}^h; \mc{V}^h)$ with the set of vertices $\mc{V}^h=\{h^{-1}x_j: 1\le j\le J\}\in \mc{P}^h=h^{-1}\mc{P}$. Denote unscaled set of vertices by $\mc{V}=\{h^{-1}x_j: 1\le j\le J\}\in \mc{P}$ and the distance to vertices by $r(z)=\dist(z;\mc{V})$ for $z\in \mc{P}$. 
\begin{lemma}[Improved remainders on polygons I]
\label{l:aprioriSteklovDescriptionb1}
Let $\Omega=\mc{P}\subset\mathbb{R}^{2}$
be a curvilinear polygon with angles $\{\alpha_j\}_{j=1}^J$. Let $\alpha_{+}=\max_j(\alpha_j)$. Then we have
\begin{equation}
        \|u_0\|_{H^k}\le Ch^{0-0}\|\tilde\gamma f_0\|_{H^{k}}+Ch^{\frac{\pi}{\alpha_+}+\frac{1}{2}-0}\|\langle h^{-1}r\rangle^{-N}u\|_{L^2},
\end{equation}
as a further improvement to \eqref{e:7l60}.
\end{lemma}
\begin{proof}
    Let $\chi,\tilde{\chi},\tilde{\tilde{\chi}}\in S^0(\mc{P}^h)$ with $\chi\prec\tilde\chi\prec\tilde{\tilde{\chi}}$ and $A_0, \tilde A_0\in \Psi^{0,0}(\partial\mc{P}^h)$ with $\WF(A_0)\cap \{|\xi|=\omega\}=\emptyset$ and $\WF A_0\cap\WF (\Id-\tilde A_0)=\emptyset$. We have from \eqref{e:7l73},
    \begin{equation}
        \tilde A_0\tilde\chi (P-\sqrt{\omega^2+h^2M^2})\tilde{\tilde{\chi}}T_h u=\tilde A_0\tilde\chi T_h f+\tilde A_0\tilde\chi(\Lambda_{hM}^{\mc{P}^h})(1-\tilde{\tilde{\chi}})T_h u -\tilde A_0\tilde\chi(\Lambda_{hM}^{\mc{P}^h}-P)\tilde{\tilde{\chi}}
    \end{equation}
    We then use the elliptic parametrix construction in Lemma~\ref{l:EllEst_Near0} together with Proposition~\ref{p:dtNParametrix} and Lemma~\ref{l:sameNearInfinity} adopted to $X_m^l(\mc{P}^h; \mc{V}^h)$ setting, to obtain that for $l<\frac{\pi}{\alpha_+}+\frac{1}{2}$ and any $N$,
    \begin{equation}
        \|A_0\chi T_h u\|_{H_k^{s,l}}\leq C\|\tilde A_0\tilde\chi T_hf\|_{H_k^{-N,l}}+\|T_h u\|_{H^{0,-N}}
    \end{equation}
    In particular, 
    \begin{equation}
        \|A_0\chi T_h u\|_{H_k^{s,0}}\leq C\|\tilde\chi T_h f_0\|_{H_k^{-N,0}}+Ch^l \|T_h u\|_{H^{0,-N}}. 
    \end{equation}
    Hence, 
    \begin{equation}
    \label{e:upgradeEstonia}
    \|u_0\|_{H^{k,0}}\leq Ch^{0-0}\|\tilde\chi f_0\|_{H^{k,0}}+Ch^{l}\|\langle h^{-1}r\rangle^{-N}u\|_{L^2}.
    \end{equation}
    Making $l=\frac{\pi}{\alpha_+}+\frac12=0$ concludes the proof. 
\end{proof}

\begin{lemma}[Improved remainders on polygons II]
\label{l:aprioriSteklovDescription}
Let $\Omega=\mc{P}\subset\mathbb{R}^{2}$
be a curvilinear polygon with angles $\{\alpha_j\}_{j=1}^J$. Let $\alpha_{+}=\max_j(\alpha_j)$. Then we have
\begin{equation}
        \|u_0\|_{H^k}\le Ch^{\frac{\pi}{\alpha_+}+1-0}\|u\|_{L^2}+Ch^{0-0}\|\tilde\gamma f_0\|_{S^0}+Ch^{\frac{\pi}{\alpha_+}-0}\sum_{j}\|\chi_j f\|_{H^{-\frac12}_h},
\end{equation}
as a further improvement to \eqref{e:7l60}. 
\end{lemma}
\begin{proof}
Lemma~\ref{l:aprioriSteklovDescriptionb1} implies \eqref{e:7l74} can be improved to
\begin{multline}
    \|v_j'\|_{Z^{\frac12, 0,\frac12}_{-L_0+0, -\frac12-0, -\frac12-0}}\le Ch^{\frac{\pi}{\alpha_+}-0}\|u_0\|_{H^{\frac12}_h}+Ch^{\frac12}(|a^i|+|a^o|)+Ch^N\|u\|_{H^{\frac12}_h}+Ch^{-\frac12-0}\|\chi_j f\|_{H^{-\frac12}_h}\\
    \le Ch^{\frac{2\pi}{\alpha_+}-\frac12-0}\|\langle h^{-1}r\rangle^{-N}u\|_{L^2}+Ch^{\frac{\pi}{\alpha_+}-0}\|\tilde\gamma f_0\|_{H^{\frac12}}+Ch^{\frac12}\|u\|_{L^2}+Ch^N\|u\|_{H^{\frac12}_h}+Ch^{-\frac12-0}\|\chi_j f\|_{H^{-\frac12}_h}
\end{multline}
Note $u=T_h^* v'$ in $r\le C^{-1}h$ and thus
\begin{multline}
    \|\langle h^{-1}r\rangle^{-N}u\|_{L^2}\le C\sum_{j=1}^J\|\langle r\rangle^{-N} v_j'\|_{L^2}+Ch^N\|u\|_{L^2}\le Ch^{\frac{2\pi}{\alpha_+}-\frac12-0}\|\langle h^{-1}r\rangle^{-N}u\|_{L^2}+Ch^{\frac12}\|u\|_{L^2}\\
    +Ch^N\|u\|_{H^{\frac12}_h}+Ch^{\frac{\pi}{\alpha_+}-0}\|\tilde\gamma f_0\|_{H^{\frac12}}+Ch^{-\frac12-0}\sum_{j}\|\chi_j f\|_{H^{-\frac12}_h}. 
\end{multline}
Absorption gives
\begin{equation}
    \|\langle h^{-1}r\rangle^{-N}u\|_{L^2}\le Ch^{\frac12}\|u\|_{L^2}+Ch^N\|u\|_{H^{\frac12}_h}+Ch^{\frac{\pi}{\alpha_+}-0}\|\tilde\gamma f_0\|_{H^{\frac12}}+Ch^{-\frac12-0}\sum_{j}\|\chi_j f\|_{H^{-\frac12}_h},
\end{equation}
as desired. Revisit Lemma~\ref{l:aprioriSteklovDescriptionb1} to conclude. 
\end{proof}
\begin{remark}
    All constants in this section can depend on $\mc{P}$ in the following way: the upper bound on the number of corners $J$, the upper and lower bounds over the angles $\alpha_j$, the upper and lower bounds on all $|\Gamma_j|$, the derivatives of the defining functions of $\Gamma_j$. 
\end{remark}

\section{The Grushin problem for the \dtn map}
\label{s:grushin}
In this section, we will use the geometric setup in Section~\ref{s:8-2}. We now set up a Grushin problem that will allow us to identify the spectrum of $h\Lambda_M^{\mc{P}}-z$: see \cite[Appendix C]{DZ:19} for further introduction to Grushin problems. Our goal is to relate the invertibility of $h\Lambda_M^{\mc{P}}-z$ to the effective Hamiltonian $E_{-+}(z;h)$ defined below. To do this, we augment $h\Lambda_M^{\mc{P}}-z$ by observation/control operators $R_\pm$ which determine the amplitudes of waves on various smooth segments of the boundary.

We define the observation operator $R_+$. Note $\chi_j,\tilde{\chi}_j\in C_c^\infty(\partial\mc{P})$ localize functions near the corner $x_j$. 
Within $\supp\tilde\chi_j$ we adopt the arclength coordinates where $x=0$ at $x_j$, $x<0$ on $\Gamma_{j-1}$, $x>0$ on $\Gamma_j$. For $1\le j\le J$, we define $R_+(z):H^{\frac12}_h(\partial\mc{P})\to \mathbb{C}^{2J}$ by 
\begin{gather}
    (R_+(z) u)_{2j-1}(z)=\int_{\Gamma_{j-1}} u(x)\chi_j'(x) e^{-iz (x+L_{j-1})/h}\,dx,\\
    (R_+(z) u)_{2j}(z)=-\int_{\Gamma_{j}} u(x)\chi_j'(x) e^{iz x/h}\,dx,
\end{gather}
within the arclength coordinate $x\in [-L_{j-1}, L_j]$ around $x_j$. Note for $\omega\ge 0$, 
\begin{equation}
    (R_+(\omega)e_{+, j-1})_{2j-1}=1+O(h^2M^2), \ (R_+(\omega)e_{-, j})_{2j}=1+O(h^2M^2).
\end{equation}
Here $R_+$ decomposes the incoming profile $u$ on $\partial\mc{P}$ into incoming data in $\mathbb{C}^{2J}$, where $(R_+)_{2j-1}$ and $(R_+)_{2j}$ read off the amplitudes of $u$ incoming towards $x_j$ from $\Gamma_{j-1}$ and $\Gamma_{j}$ respectively. 

We define the control operator $R_-:\mathbb{C}^{2J}\to H^{-\frac12}_h(\partial\mc{P})$. For $u_-=(u_{-, k})_{k=1}^{2J}\in \mathbb{C}^{2J}$ we let
\begin{equation}
    R_-(z)u_-=\sum_{j=1}^{J}-iu_{-,2j-1}\mathbbm{1}_{\Gamma_{j-1}}\tilde{\chi}_j'(x)e^{-iz(x+L_{j-1})/h}+iu_{-,2j}\mathbbm{1}_{\Gamma_{j}}\tilde{\chi}_j'(x)e^{iz x/h}.
\end{equation}
within the arclength coordinate $x\in [-L_{j-1}, L_j]$ around $x_j$.
Here $R_-(z)$ reconstructs outgoing profiles from outgoing data $u_-\in\mathbb{C}^{2J}$, where $u_{-, 2j-1}$, $u_{-,2j}$ are the amplitudes of reconstructed waves outgoing from $x_j$ towards $\Gamma_{j-1}$ and $\Gamma_{j}$ respectively. We give a schematic display of what $\mathbb{C}^{2J}$-data means:
\begin{equation}
    \begin{matrix}
        R_+(z)\\
        \text{ maps to }\\
        \text{ amplitudes }
    \end{matrix}\begin{pmatrix}
        \text{incoming to }x_1\text{ on }\Gamma_{J}\\
        \text{incoming to }x_1\text{ on }\Gamma_{1}\\
        \text{incoming to }x_2\text{ on }\Gamma_{1}\\
        \text{incoming to }x_2\text{ on }\Gamma_{2}\\
        \cdots\\
        \text{incoming to }x_J\text{ on }\Gamma_{J-1}\\
        \text{incoming to }x_J\text{ on }\Gamma_{J}
    \end{pmatrix},\ 
    \begin{matrix}
        R_-(z)\\
        \text{ maps from }\\
        \text{ amplitudes }
    \end{matrix}
    \begin{pmatrix}
        \text{outgoing from }x_1\text{ on }\Gamma_{J}\\
        \text{outgoing from }x_1\text{ on }\Gamma_{1}\\
        \text{outgoing from }x_2\text{ on }\Gamma_{1}\\
        \text{outgoing from }x_2\text{ on }\Gamma_{2}\\
        \cdots\\
        \text{outgoing from }x_J\text{ on }\Gamma_{J-1}\\
        \text{outgoing from }x_J\text{ on }\Gamma_{J}
    \end{pmatrix}.
\end{equation}

We now define the scattering matrix and the transfer matrix. Let
\begin{equation}
\Sigma^{\mc{P}}_{M,h^{-1}\omega}=\begin{pmatrix}\Sigma^{\mc{W}_1}_{M,h^{-1}\omega}&0&\dots &0\\
0&\Sigma^{\mc{W}_2}_{M,h^{-1}\omega}&\dots &0\\
\vdots&\vdots&\ddots&\vdots\\
0&0&\dots&\Sigma^{\mc{W}_J}_{M,h^{-1}\omega}
\end{pmatrix}, \ 
T_{z/h}=T_-\begin{pmatrix}E_{L_1}^{z/h}&0&\dots &0\\
0&E_{L_2}^{z/h}&\dots &0\\
\vdots&\vdots&\ddots&\vdots\\
0&0&\dots&E_{L_J}^{z/h}
\end{pmatrix} T_+
\end{equation}
where $T_\pm$ cyclically shifts up/down the data by one row, and $E_{L_j}^{z/h}$ is the twisted phase shift on $\Gamma_j$: 
\begin{equation}
    (T_+ \mathbf{a})_{k}=(\mathbf{a})_{k+1},\  (T_- \mathbf{a})_{k}=(\mathbf{a})_{k-1},\ E_{L_j}^{z/h}=\begin{pmatrix}
    0 & e^{-iL_j z/h}\\
    e^{-iL_j z/h} & 0
\end{pmatrix}. 
\end{equation}
for vector $\mathbf{a}\in\mathbb{C}^{2J}$. 

The above construction are on each individual component $\Gamma^l$ labeled by $l$. We can then glue up them in $l$ to obtain the global construction on all components: 
\begin{gather}
    R_+(z)=\bigoplus_{l=1}^L R_+^l(z): H^{\frac12}_h(\partial\mc{P})\rightarrow \mathbb{C}^{2J}, \ R_-(z)=\bigoplus_{l=1}^L R_-^l(z): \mathbb{C}^{2J}\rightarrow H^{-\frac12}_h(\partial\mc{P}),\\
    \Sigma^{\mc{P}}_{M,h^{-1}\omega}=\operatorname{diag}(\Sigma^{\mc{P}, l}_{M,h^{-1}\omega}; 1\le l\le L), \ T_{z/h}=\operatorname{diag}(T_{z/h}^l; 1\le l\le L). 
\end{gather}
Note here $J$ is now the total number of vertices on all components of $\partial\mc{P}$. 

In this paper, $z\mapsto\sqrt{z}$ has the branch cut along the negative real axis. Here is our main result of this section: 
\begin{proposition}[The Grushin problem for $\Lambda_M^{\mc{P}}$]\label{p:grushin}
There exists $h_0<\frac12$, $r_1>0$ such that for $z\in [\frac12, 2]+i[-r_1 h, r_1h]$, and $h\in(0,h_0]$, the Grushin problem $\mc{G}(z;h):H^{\frac{1}{2}}(\partial\mc{P})\oplus \mathbb{C}^{2J}\to H^{-\frac12}(\partial\mc{P})\oplus \mathbb{C}^{2J}$ defined by
\begin{equation}
    \mc{G}(z;h):=\begin{pmatrix}
\Lambda_M^{\mc{P}}-\sqrt{h^{-2}z^2+M^2}&R_-(z)\\R_+(z)&0
\end{pmatrix}
\end{equation}
is well-posed with inverse
\begin{equation}
    \mc{G}^{-1}(z;h)=\begin{pmatrix} E(z;h)&E_+(z;h)\\ E_-(z;h)& E_{-+}(z;h)\end{pmatrix}, 
\end{equation}
where
\begin{equation}
    \|\mc{G}^{-1}(z;h)\|_{H_h^{-\frac{1}{2}}(\partial\mc{P})\oplus \mathbb{C}^{2J}\to H_h^{\frac12}(\partial\mc{P})\oplus \mathbb{C}^{2J}}\le C,
\end{equation}
and
\begin{equation}\label{e:7l14}
    E_{-+}(z;h)=\Sigma^{\mc{P}}_{M,h^{-1}\omega}-T_{z/h}+O(h^{\frac{2\pi}{\alpha_+}-0}+M^2h^{2-0})+O(|\Im z/h|)
\end{equation}
is analytic in $z\in [\frac12, 2]+i[-r_1 h, r_1h]$. 
\end{proposition}
\begin{remark}
    This implies for $z\in\mathbb{C}$ with $\Re z>0$ sufficiently large and $|\Im z|\le r_1$, we can define
    \begin{equation}
        \mc{Q}^{\Lambda}_{z}=E_{-+}(z)=\Sigma^{\mc{P}}_{M,\Re z}-T_{z}+O(|\Re z|^{\frac{2\pi}{\alpha_+}-0}+M^2|\Re z|^{2-0})+O(|\Im z|),
    \end{equation}
    analytic in $z$, by picking a small $h$ and set $E_{-+}(z)=E_{-+}(hz;h)$. 
\end{remark}
\begin{proof}
1. For readability, we only address the case $\partial\mc{P}$ has only one connected component. The same proof works for cases when $\partial\mc{P}$ has more than one connected component. In Steps 1 to 7, we prove the invertibility of $\mc{G}(\omega;h)$ in the case $z=\omega\in[\frac12, 2]$. Suppose that 
\begin{equation}
    \mc{G}(z;h)\begin{pmatrix}u\\u_-\end{pmatrix}=\begin{pmatrix}v\\v_+\end{pmatrix}.
\end{equation}
Since $\Lambda_{M}^\mc{P}:H^{\frac12}\to H^{-\frac12}$ is Fredholm with index 0, the same holds for $\mc{G}(\omega;h)$. Therefore to show invertibility of $\mc{G}(\omega;h)$, it is enough to check its injectivity. 

2a. We set up some analytic tools we need later. Let $\gamma_{j, 0}, \gamma_j, \tilde\gamma_j,\gamma_{j,0.5},\tilde{\tilde{\gamma}}_{j},\gamma_{j,1}, \gamma_{j,2} \in C_c^\infty(\Gamma_j)$ be cutoffs that $\gamma_{j,0}\prec \gamma_j\prec \tilde\gamma_j\prec\gamma_{j,0.5}\prec\tilde{\tilde{\gamma}}_j\prec\gamma_{j,1}\prec\gamma_{j,2}$ such that $\chi_j'|_{\Gamma_j}\prec\gamma_{j,0}$, $\chi_{j+1}'|_{\Gamma_{j}}\prec \gamma_{j,0}$, $\supp (1-\gamma_{j, 0})\cap\Gamma_j\subset\supp{\chi_j}\cup\supp{\chi_{j+1}}$ on each edge $\Gamma_j$. Let $\tilde \gamma_j=\tilde \gamma_j(X_\omega+X_0+X_\infty)$, where $X_\omega\in \Psi^{\comp}_h(\Gamma_j)$ with $\WF_h(X_\omega)\cap \{\xi=0\}=\emptyset$ and $\WF_h(I-X_\omega)\cap \{|\xi|=\omega\}\cap (\supp \chi_j\cup \supp \chi_{j+1})=\emptyset$, $X_0\in \Psi^{\operatorname{comp}}_h$ with $\WF_h(X_0)\subset \{|\xi|<\frac \e 4\}$ and $\WF_h(I-X_0)\cap \{|\xi|=0\}\cap (\supp \chi_j\cup \supp \chi_{j+1})=\emptyset$, for some $\e>0$ to be chosen small. Let the standard semiclassical parametrix $E_\infty\in\Psi_h^{-1}$ be that $(hB-\ohm)\tilde{\tilde{\gamma}}_j E_\infty=\tilde\gamma_j X_\infty+h^\infty\Psi^{-\infty}$, where $B\in\Psi^1$ is the local structure of $\Lambda_M^{\mc{P}}$ in the sense of 
\begin{equation}\label{e:transform2}
    \Lambda^{\mc{P}}_M\gamma_{j,1}=\gamma_{j,2} B\gamma_{j,1}+\mc{E}\gamma_{j,1}, \ \mc{E}: H^{-N}\rightarrow H^0
\end{equation}
by taking the adjoint of \eqref{eq:TransformationProperty}.  

2b. The rest of Step 2 is devoted to show the existence of $E_0(h)\in\Psi^{0}$ uniformly bounded for small $h>0$, with $(hB-\ohm)\tilde{\tilde{\gamma}}_j E_0(h)=\tilde\gamma_j X_0+O(h \Psi^{-\infty})$. For the rest of this proof, we may abbreviate $\omega_{hM}=\ohm$ if it makes the proof easier to read. Indeed, let $\tilde X_0\in \Psi^{\operatorname{comp}}_h$ be that $\tilde X_0 X_0=X_0+O(h^\infty\Psi^{-\infty})$ with $\WF_h(\tilde X_0)\subset \{|\xi|_g<\frac {\e}{ 2}\}$. Then, since $B:H^1\to L^2$ and 
\begin{equation}
    \|h\tilde X_0\|_{L^2\to H^1}\leq C(h\|\tilde{X}_0\|_{L^2\to L^2}+C\|h\nabla \tilde{X}_0\|_{L^2\to L^2})\leq C\e,
\end{equation}
we have $\|hB\tilde{X}_0\|_{L^2\to L^2}\leq C\e$. Thus, choosing $\e>0$ small enough, uniformly for small $h$, $(h B\tilde X_0-\omega)$ is invertible on $L^2$ and
\begin{equation}
    (h B\tilde X_0-\omega_{hM})^{-1}=-\omega_{hM}^{-1}\Big(\Id+\sum_{k=1}^\infty (\omega_{hM}^{-1}h B\tilde X_0)^k\Big). 
\end{equation}

2c. We next show that for $k\in \mathbb{Z}$, there is $C$ such that 
\begin{equation}
\label{e:regularityClaim}
\|
(hB\tilde{X}_0-\omega_{hM})^{-1}\|_{H^k\to H^k}\leq C,
\end{equation}
uniformly in $h$ small. For this, observe that for $A\in \Psi^k$, 
\begin{equation}\label{e:9l12}
[A,(hB\tilde{X}_0-\omega_{hM})^{-1}]= -(hB\tilde{X}_0-\omega_{hM})^{-1}[A, hB\tilde{X}_0](hB\tilde{X}_0-\omega_{hM})^{-1},
\end{equation}
where
\begin{equation}
\label{e:commuteForm}
[A,hB\tilde{X}_0]=h[A,B]\tilde{X_0}+hB[A,\tilde{X}_0]\in \Psi^{k-1},
\end{equation}
and $h\tilde{X}_0\in\Psi^{-1}$, $[A,B]\in \Psi^{k}$, and $[A,h\tilde{X}_0]\in \Psi^{k-2}$. 
Hence, for any $s\in\mathbb{R}$ we have uniformly in $h>0$ small that
\begin{equation}
\label{e:commuteEstimate}
\begin{aligned}
\|[A,hB\tilde{X}_0]\|_{H^{s+k-1}\to H^s}&\leq C.
\end{aligned}
\end{equation}
In particular, using this with $A=\langle D\rangle$ implies 
\begin{equation}
    \langle D\rangle(hB\tilde{X}_0-\omega_{hM})^{-1}\langle D\rangle^{-1}=(hB\tilde{X}_0-\omega_{hM})^{-1}-(hB\tilde{X}_0-\omega_{hM})^{-1}[\langle D\rangle,hB\tilde{X}_0](hB\tilde{X}_0-\omega_{hM})^{-1}\langle D\rangle^{-1}
\end{equation}
is bounded on $L^2$ and hence $(hB\tilde{X}_0-\omega_{hM})^{-1}:H^1\to H^1$ is uniformly bounded. Now take $A=\langle D\rangle^2$ gives the boundedness of $(hB\tilde{X}_0-\omega_{hM})^{-1}:H^2\to H^2$. Iterative arguments and duality show \eqref{e:regularityClaim}.

2d. Next, we claim that for $A_1,\cdots, A_N\in \Psi^0$, for any $k\in\mathbb{Z}$, 
\begin{equation}
\label{e:semiclassicalClaim}
\|\ad_{A_N}\ad_{A_{N-1}}\cdots \ad_{A_1}(hB\tilde{X}_0-\omega_{hM})^{-1}\|_{H^{k}\to H^{k+N-1}}\leq Ch.
\end{equation}
To do this, we argue as above. This time using \eqref{e:commuteForm} for $A\in \Psi^0$, $[A,B]\in \Psi^0$, and $[A,\tilde{X}_0]\in \Psi^{-1}$, we have
\begin{equation}
\label{e:commutator1}
[A,hB\tilde{X}_0]\in h\Psi^{0}
\end{equation}
Together with~\eqref{e:regularityClaim}, this implies~\eqref{e:semiclassicalClaim} with $N=1$. We now claim by induction that for any $N\geq 1$ and $A_1,\dots, A_N\in \Psi^0$, 
\begin{equation}
\label{e:inductionFormula}
\ad_{A_N}\ad_{A_{N-1}}\dots \ad_{A_1}(hB\tilde{X}_0-\omega_{hM})^{-1}=h\sum_{l=1}^{L_N} (hB\tilde{X}_0-\omega_{hM})^{-1}\prod_{m=1}^{M_l} Q_{lm}(hB\tilde{X}_0-\omega_{hM})^{-1},
\end{equation}
for some $L_N\in\mathbb{N}$, $Q_{lm}\in \Psi^{n_{lm}}$, $n_{lm}\leq 0$ and $\sum_{m=1}^{M_l} n_{lm}=-N+1$. Notice that~\eqref{e:semiclassicalClaim} then follows from~\eqref{e:inductionFormula} and~\eqref{e:regularityClaim} so it remains to prove~\eqref{e:inductionFormula}. The equality~\eqref{e:inductionFormula} with $N=1$ follows from 
\eqref{e:9l12} and~\eqref{e:commutator1}. Suppose by induction~\eqref{e:inductionFormula} holds for some $N\geq 1$. Let $A_1,\cdots, A_N, A_{N+1}\in \Psi^0$. Then,
\begin{multline}
\ad_{A_{N+1}}\ad_{A_{N}}\dots \ad_{A_1}(hB\tilde{X}_0-\omega_{hM})^{-1}=h\sum_{l=1}^{L_N} [A_{N+1},(hB\tilde{X}_0-\omega_{hM})^{-1}\prod_{m=1}^{M_l} Q_{lm}(hB\tilde{X}_0-\omega_{hM})^{-1}]\\
=-h\sum_{l=1}^{L_N} (hB\tilde{X}_0-\omega_{hM})^{-1}[A_{N+1},hB\tilde{X}_0](hB\tilde{X}_0-\omega_{hM})^{-1}\prod_{m=1}^{M_l} Q_{lm}(hB\tilde{X}_0-\omega_{hM})^{-1}\\
+h\sum_{l=1}^{L_N} (hB\tilde{X}_0-\omega_{hM})^{-1}\sum_{k=1}^{M_l}\Big(\prod_{m=1}^{k-1} Q_{lm}(hB\tilde{X}_0-\omega_{hM})^{-1}\Big)\Big([A_{N+1},Q_{lk}](hB\tilde{X}_0-\omega_{hM})^{-1}\\
 -Q_{lk}(hB\tilde{X}_0-\omega_{hM})^{-1}[A_{N+1},hB\tilde{X}_0](hB\tilde{X}_0-\omega_{hM})^{-1}\Big)
 \prod_{m=k+1}^{M_l}Q_{lm}(hB\tilde{X}_0-\omega_{hM})^{-1}. 
\end{multline}
The equality~\eqref{e:inductionFormula} with $N$ replaced by $N+1$ now follows from the fact that $[A_{N+1},hB\tilde{X}_0]\in \Psi^{-1}$ and $[A_{N+1},Q_{lm}]\in \Psi^{n_{lm}-1}$.

2e. We can now complete Step 2b. Let $E_0=(h B\tilde X_0-\omega_{hM})^{-1}\tilde\gamma_j X_0\in\Psi^{0}$. We then have by~\eqref{e:semiclassicalClaim} that for any $N$, 
\begin{equation}
    (1-\tilde{\tilde{\gamma}}_j)(h B\tilde X_0-\omega_{hM})^{-1}\tilde\gamma_j X_0=(1-\tilde{\tilde{\gamma}}_j)\ad_{(1-\gamma_{j, 0.5})}^{N}(h B\tilde X_0-\omega_{hM})^{-1}\tilde\gamma_j X_0\in h\Psi^{-N+1},
\end{equation}
Similarly, using that $\WF_h(1-\tilde{X}_0)\cap \WF_h(X_0)=\emptyset$, and letting $X_{0}'\in\Psi_h^0$ with $\WFh(I-\tilde{X}_0)\cap \WFh(X_{0}')=\WFh(1-X_{0}')\cap \WFh(X_0)=\emptyset$, we obtain for any $N>0$ that 
\begin{multline}
    (hB-\omega_{hM})\tilde{\tilde{\gamma}}_j E_0=\tilde\gamma_j X_0+hB(1-\tilde X_0)(h B\tilde X_0-\omega_{hM})^{-1}\tilde\gamma_jX_0 +O(h\Psi^{-N}) \\
    =\tilde\gamma_j X_0+hB(1-\tilde{X}_0)(\ad_{(1-X_0')}^{N}(h B\tilde X_0-\omega_{hM})^{-1})\tilde\gamma_jX_0+O(h\Psi^{-N}) 
    =\tilde\gamma_j X_0+O(h\Psi^{-N}),
\end{multline}
as desired. We also have
\begin{multline}
    (1-\tilde X_0)E_0=(1-\tilde X_0)(h B\tilde X_0-\omega_{hM})^{-1}\tilde\gamma_j X_0=(1-\tilde X_0)(\ad_{(1-X_0')}^N(h B\tilde X_0-\omega_{hM})^{-1})\tilde\gamma_j X_0\\+O(h^\infty\Psi^{-\infty})
    =O(h\Psi^{-N}),
\end{multline}
which gives $\|(1-\tilde X_0)E_0 u\|_{H^N}\le Ch\|u\|_{H^{-\frac12}_h}$. On the other hand, we have
\begin{equation}
    \|\tilde X_0 E_0 u\|_{H^N_h}\le C\|E_0 u\|_{L^2}\le C\|\tilde\gamma_j X_0 u\|_{L^2}\le C \|u\|_{H^{-\frac12}_h}.
\end{equation}
Hence, for any $N>0$, there is $C>0$ that for all $u\in H^{-\frac12}_h$, 
\begin{equation}\label{e:parametrixzeroestimate}
    \|E_0 u\|_{H^{N}_h}\le C\|u\|_{H^{-\frac12}_h}.
\end{equation}

3. We now approximate $u$ on edges away from the corners. Let $T>\max_j L_j$. Our strategy is to construct quasimodes $\tilde u$ that characterize the influence of $v$, and $u-\tilde u$ will then behave like plane waves. Indeed, we can construct $\tilde u_j$ localised on each edge $\Gamma_j$, that approximate $u$, via 
\begin{equation}
    \tilde{u}_j:= -i\tilde \gamma_j \int_0^T e^{it(B-\omega_{hM}/h)}X_\omega \tilde\gamma_j v dt +h\tilde{\tilde{\gamma}}_{j} (E_\infty +E_0)\tilde\gamma_j v, \ \tilde u=\sum_{j=1}^J \tilde u_{j},
\end{equation}
with $E_\infty$ and $E_0$ given in Step 2. We can estimate that
\begin{equation}
    \|\tilde{u}_j\|_{H^{-\frac12}_h}\le CT\|\langle hD\rangle^{\frac12}X_\omega \tilde \gamma_j v\|_{H^{-\frac12}_h}+Ch\|\tilde X_\infty \tilde \gamma_j v\|_{H^{-\frac12}_h}+Ch\|\tilde X_0 \tilde \gamma_j v\|_{H^{-\frac12}_h}\le CT\|\tilde \gamma_j v\|_{H^{-\frac12}_h},
\end{equation}
for some $\tilde X_\infty\in\Psi^{0}_h$ with $\WFh(1-\tilde X_\infty)\cap \WFh X_\infty=\emptyset$, since $e^{it(B-\omega_{hM}/h)}$ is unitary on $L^2$ and we have \eqref{e:parametrixzeroestimate}. Note that from \eqref{e:transform2}, we have
\begin{multline} \label{e:QuasiModeApprox}
    (h\Lambda_M^{\mc{P}}-\ohm)\tilde u_j=-i\gamma_{j,2}[hB, \tilde \gamma_j]\int_0^T e^{it(B-\omega_{hM}/h)} X_\omega \tilde \gamma_j v\,dt+ h \tilde \gamma_j X_\omega\tilde \gamma_j v\\
    -h\tilde \gamma_j e^{iT(B-\omega_{hM}/h)} X_\omega\tilde \gamma_j v
    +\gamma_{j,2}(hB-\omega_{hM})h\tilde{\tilde{\gamma}}_j E_\infty \tilde \gamma_j v
+\gamma_{j,2}(hB- \omega_{hM})h\tilde{\tilde{\gamma_j}} E_0 \tilde \gamma_j v+O_{H^{0}}(h\|\tilde \gamma_j v\|_{H^{-N}})\\
=-i\gamma_{j,2}[hB, \tilde \gamma_j]\int_0^T e^{it(B-\omega_{hM}/h)} X_\omega \tilde \gamma_j v\,dt+ h \tilde \gamma_j^2 v+O_{H^{0}}(h\|\tilde \gamma_j v\|_{H^{-N}}),
\end{multline}
since $\WFh(\tilde \gamma_j e^{iT(B-\omega_{hM}/h)}\tilde \gamma_j X_\omega)=\emptyset$ via propagation. It is interpreted that $\tilde u_j$ is a good quasimode on $\Gamma_j$ that $(h\Lambda_M^{\mc{P}}-\ohm)\tilde u_j$ is close to $h\tilde \gamma_j^2 v$. Here,
\begin{equation}
    \|e^{it(B-\omega_{hM}/h)} X_\omega \tilde \gamma_j v\|_{H^{-\frac12}_h} \leq \|e^{it(B-\omega_{hM}/h)} X_\omega \tilde \gamma_j v\|_{L^2} \leq C\|X_\omega \tilde \gamma_j v\|_{L^2} \leq C\|\tilde \gamma_j v\|_{H^{-\frac12}_h},
\end{equation}
and we can estimate 
\begin{equation}
    \|(h\Lambda_M^{\mc{P}}-\ohm)\tilde u_j\|_{H^{-\frac12}_h}\le CTh
    \|\tilde \gamma_j v\|_{H^{-\frac12}_h}
\end{equation}
and then by that $\Lambda^{\mc{P}}_M$ is classically elliptic, we have
\begin{equation}\label{e:7l8}
    \|\tilde{u}_j\|_{H^{\frac12}_h}\le C \|\tilde{u}_j\|_{H^{-\frac12}_h}+C\|h\Lambda^{\mc{P}}_M\tilde{u}_j\|_{H^{-\frac12}_h}\le C(1+\omega)T(1+h)\|\tilde \gamma_j v\|_{H^{-\frac12}_h} \le C\|\tilde \gamma_j v\|_{H^{-\frac12}_h}. 
\end{equation}
Now
\begin{multline}\label{e:7l22}
    (h\Lambda_M^{\mc{P}}-\ohm)(u-\tilde u)=\sum_j i\gamma_{j,2} [hB, \tilde \gamma_j]\int_0^T e^{it(B-\omega_{hM}/h)} X_\omega \tilde \gamma_j v~dt + h (1-\sum_j\tilde \gamma_j^2) v\\-hR_-(\omega) u_-
    +O_{H^{0}}(h\|v\|_{H^{-N}})=f.
\end{multline}
Note here we have the global expression of $f$ on $\partial\mc{P}$:
\begin{equation}
    f=-hR_-(\omega) u_-+O_{H^{-\frac12}_h}(h \|v\|_{H^{-\frac12}_h})=\sum_j \Big(-ih(u_{-})_{2j}\tilde\chi_{j}' e_{+, j}+ih(u_{-})_{2j+1}\tilde\chi_{j+1}' e_{-, j}\Big)+O_{H^{-\frac12}_h}(h \|v\|_{H^{-\frac12}_h}).
\end{equation}
We have
\begin{gather}\label{e:9l1}
    \|u-\tilde u\|_{H^{\frac12}_h}\le C\|u-\tilde u\|_{H^{-\frac12}_h}+C\|f\|_{H^{-\frac12}_h}
    \le C\|u\|_{H^{-\frac12}_h}+C\|v\|_{H^{-\frac12}_h}+Ch|u_-|,\\ \label{e:9l7}
    \|u\|_{H^{\frac12}_h}\le C\|u\|_{H^{-\frac12}_h}+C\|v\|_{H^{-\frac12}_h}+Ch|u_-|.
\end{gather}

4. We want to understand how the scattering data near two ends of the same edge communicate via propagation. In this step, we are in the global coordinates on $\Gamma_j$ where $x=0$ at $x_j$, $x=L_j$ at $x_{j+1}$. See Figure~\ref{f:9l1}. Lemma~\ref{l:aprioriSteklovDescription} applied to \eqref{e:7l22} give
\begin{equation}\label{e:7l5}
\gamma_j(u-\tilde{u})= \gamma_j\Big(\tilde a_{+,j}(x;h)e_{+,j,h}+\tilde a_{-,j}(x;h)e_{-,j,h}\Big) +O_{S^0}(h^{\frac{\pi}{\alpha_+}+1-0}\|u\|_{L^2}+h^{1-0}\|v\|_{H_h^{-\frac12}}),
\end{equation}
and allows us to read off $a_{\pm, j}^{o/i}(h)\in \mathbb{C}$ from $\gamma_j(u-\tilde{u})$ from the left of $\Gamma_j$, by
\begin{gather}\label{e:9l9}
    \tilde a_{+,j}(x;h)\sim a_{+,j}^o(h)+(u_{-})_{2j}(\tilde\chi_j(x)-1)+O_{S^0}(h^2M^2\|\tilde a_{+,j}\|_{S^0}), \text{ near }\supp\gamma_j,\\ \label{e:9l10}
    \tilde a_{-,j}(x;h)\sim a_{-,j}^i(h)+(u_{-})_{2j+1}\tilde\chi_{j+1}(x)+O_{S^0}(h^2M^2\|\tilde a_{-,j}\|_{S^0}), \text{ near }\supp\gamma_j. 
\end{gather}
We can also use Lemma~\ref{l:aprioriSteklovDescription} to read off $a_{\mp, j}^{o/i}(h)\in \mathbb{C}$ from the right of $\Gamma_j$, by
\begin{gather}
    \tilde a_{+,j}(x;h)\sim e^{-i\omega L_j/h}a_{+,j}^i(h)+e^{-i\omega L_j/h}(u_{-})_{2j}\tilde\chi_j(x)+O_{S^0}(h^2M^2\|\tilde a_{+,j}\|_{S^0}), \text{ near }\supp\gamma_j,\\ 
    \tilde a_{-,j}(x;h)\sim e^{i\omega L_j/h}a_{-,j}^o(h)+e^{i\omega L_j/h}(u_{-})_{2j+1}(\tilde\chi_{j+1}(x)-1)+O_{S^0}(h^2M^2\|\tilde a_{-,j}\|_{S^0}), \text{ near }\supp\gamma_j. 
\end{gather}
Here the phase shifts come from the fact that $a_{\mp, j}^{o/i}(h)$ are read near $x_{j+1}$ but normalized near $x_j$, in bases $e^{\mp i\omega(x+L_j)}$ in coordinates where $x=0$ at $x_j$. We remind the reader that $v_+=R_+(\omega)u=R_+(\omega)(u-\tilde u)+R_+(\omega)\tilde u$, and
\begin{equation}
    |R_+(\omega)\tilde u|_{\mathbb{C}^{2J}}\le C\|\tilde u\|_{H^{\frac12}_h}\le C\|v\|_{H^{-\frac12}_h}
\end{equation}
from \eqref{e:7l8} implies
\begin{gather}
    (v_+)_{2j-1}=(R_+(\omega)(u))_{2j-1}\sim a_{+, j-1}^i+O(\|v\|_{H^{-\frac12}_h}+h^2M^2\|\tilde a_{-,j-1}\|_{S^0}),\\ 
    (v_+)_{2j}=(R_+(\omega)u)_{2j}\sim a_{-, j}^i+O(\|v\|_{H^{-\frac12}_h}+h^2M^2\|\tilde a_{-,j}\|_{S^0}).
\end{gather}
Note the unique determination of $\tilde a_{\pm,j}$ in Lemma~\ref{l:aprioriSteklovDescription} then gives
\begin{equation}\label{e:9l5}
    \|\tilde a_{\pm,j}\|_{S^0}\le C\|\tilde a_{\pm,j}\|_{L^\infty}\le C|v_+|+C|u_-|+C\|v\|_{H^{-\frac12}_h}.
\end{equation}
We then use $O(h^2M^2)$ to abbreviate $O(h^2M^2(|v_+|+|u_-|+\|v\|_{H^{-\frac12}_h}))$. Evaluating the above reading at some $x\in\Gamma_j$ where $\tilde\chi_j(x)=\tilde\chi_{j+1}(x)=0$ gives
\begin{gather}
    a_{+, j}^o=e^{-i\omega L_j/h}a_{+, j}^i+(u_-)_{2j}+O(h^2M^2)+O(h^\infty),\\
    a_{-, j}^o=e^{-i\omega L_j/h}a_{-, j}^i+(u_{-})_{2j+1}+O(h^2M^2)+O(h^\infty)
\end{gather}

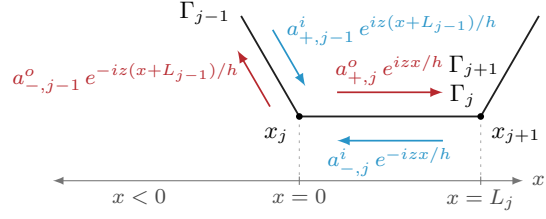
\begin{figure}[t]
  \centering
        \begin{tikzpicture}[x=1cm,y=1cm]
  \usetikzlibrary{calc,plotmarks}
  \definecolor{colorAKANE}{RGB}{183,40,46}
  \definecolor{colorUSUAI}{RGB}{40,148,200}
  \colorlet{plotPrimary}{colorAKANE}
  \colorlet{plotTertiary}{colorUSUAI}
  \colorlet{plotAxis}{black!72}
  \colorlet{plotGuide}{black!34}
  \tikzset{
    >=latex,
    font=\small,
    mark size=1.5pt,
    mark options={thin},
    line cap=round,
    line join=round,
    boundary edge/.style={draw=black!88,line width=0.7pt},
    arclength axis/.style={draw=plotAxis!75,line width=0.4pt},
    origin guide/.style={draw=plotGuide,line width=0.3pt,
      dash pattern=on 1pt off 2pt},
    wave arrow/.style={line width=0.65pt,->},
    incoming wave/.style={wave arrow,draw=plotTertiary,text=plotTertiary},
    outgoing wave/.style={wave arrow,draw=plotPrimary,text=plotPrimary},
    figure label/.style={font=\fontsize{9}{10.8}\selectfont,
      text=black,inner sep=1pt},
    wave label/.style={font=\fontsize{8}{9.6}\selectfont,inner sep=1pt},
    axis label/.style={font=\fontsize{8}{9.6}\selectfont,
      text=plotAxis,inner sep=1pt}
  }

  \path[use as bounding box] (0,0) rectangle (7.5,3.2);
  \begin{pgfinterruptboundingbox}
  \begin{scope}[shift={(0.10,-0.175)}]
    \coordinate (vertex j) at (3.90,1.60);
    \coordinate (vertex next) at (6.30,1.60);
    \coordinate (left extent) at (3.13,2.93);
    \coordinate (right extent) at (7.07,2.93);
    \coordinate (axis origin) at ($(vertex j)+(0,-0.85)$);
    \coordinate (axis left limit) at (0.64,0);
    \coordinate (axis right limit) at (6.86,0);
    \coordinate (axis left) at (axis left limit |- axis origin);
    \coordinate (axis right) at (axis right limit |- axis origin);
    \coordinate (axis next) at (vertex next |- axis left);
    \coordinate (edge center) at ($(vertex j)!0.5!(vertex next)$);
    \coordinate (edge outgoing center) at ($(edge center)+(0,0.32)$);
    \coordinate (edge incoming center) at ($(edge center)+(0,-0.32)$);

    \coordinate (left flow upper) at ($(left extent)!0.24!(vertex j)$);
    \coordinate (left flow lower) at ($(left extent)!0.80!(vertex j)$);
    \coordinate (left normal) at
      ($(left flow upper)!0.27cm!90:(left flow lower)-(left flow upper)$);
    \coordinate (left incoming upper) at ($(left flow upper)+(left normal)$);
    \coordinate (left incoming lower) at ($(left flow lower)+(left normal)$);
    \coordinate (left outgoing upper) at ($(left flow upper)-(left normal)$);
    \coordinate (left outgoing lower) at ($(left flow lower)-(left normal)$);

    \draw[origin guide] (vertex j) -- (axis origin);
    \draw[origin guide] (vertex next) -- (axis next);

    \begin{scope}
      \draw[boundary edge]
        (left extent) -- (vertex j) -- (vertex next) -- (right extent);
      \path[fill=black,draw=none]
        (vertex j) circle[radius=1.25pt]
        (vertex next) circle[radius=1.25pt];
      \node[figure label,anchor=east] at ($(left extent)+(-0.10,0.04)$)
        {$\Gamma_{j-1}$};
      \node[figure label,anchor=south]
        at ($(vertex j)!0.90!(vertex next)+(0,0.08)$)
        {$\Gamma_j$};
      \node[figure label,anchor=east]
        at ($(vertex next)!0.50!(right extent)+(-0.10,0)$)
        {$\Gamma_{j+1}$};
      \node[figure label,anchor=north east] at ($(vertex j)+(-0.10,-0.10)$)
        {$x_j$};
      \node[figure label,anchor=north west] at ($(vertex next)+(0.10,-0.10)$)
        {$x_{j+1}$};
    \end{scope}

    \begin{scope}
      \draw[incoming wave] (left incoming upper) -- (left incoming lower)
        node[wave label,pos=0,right=4pt]
          {$a_{+,j-1}^{i}\,e^{iz(x+L_{j-1})/h}$};
      \draw[outgoing wave] (left outgoing lower) -- (left outgoing upper)
        node[wave label,midway,left=4pt]
          {$a_{-,j-1}^{o}\,e^{-iz(x+L_{j-1})/h}$};
    \end{scope}

    \begin{scope}
      \draw[outgoing wave]
        ($(edge outgoing center)+(-0.70,0)$) -- ($(edge outgoing center)+(0.70,0)$)
        node[wave label,midway,above=1.5pt]
          {$a_{+,j}^{o}\,e^{izx/h}$};
      \draw[incoming wave]
        ($(edge incoming center)+(0.70,0)$) -- ($(edge incoming center)+(-0.70,0)$)
        node[wave label,midway,below=1pt]
          {$a_{-,j}^{i}\,e^{-izx/h}$};
    \end{scope}

    \begin{scope}
      \draw[arclength axis,<->] (axis left) -- (axis right);
      \draw[arclength axis]
        ($(axis origin)+(0,-0.055)$) -- ($(axis origin)+(0,0.055)$);
      \draw[arclength axis]
        ($(axis next)+(0,-0.055)$) -- ($(axis next)+(0,0.055)$);
      \node[axis label,anchor=north] at ($(axis origin)+(0,-0.07)$)
        {$x=0$};
      \node[axis label,anchor=north] at ($(axis next)+(0,-0.07)$)
        {$x=L_j$};
      \node[axis label,anchor=north]
        at ($(axis left)!0.35!(axis origin)+(0,-0.07)$) {$x<0$};
      \node[axis label,anchor=west] at ($(axis right)+(0.08,0)$) {$x$};
    \end{scope}
  \end{scope}
  \end{pgfinterruptboundingbox}
\end{tikzpicture}%
  \caption{Scattering and propagation near vertex $x_j$ on $\partial\mc{P}$}
  \label{f:9l1}
\end{figure}

5. In order to understand how the scattering data communicate through the corner $x_j$, we use Lemma~\ref{l:blowup} to blow up near the corner $x_j$. See Figure~\ref{f:9l1}. Note
\begin{equation}
    \chi_j f=O_{H^{-\frac12}_h}(h\|v\|_{H^{-\frac12}_h}), \ \tilde\gamma f_0=O_{S^0}(h\|v\|_{H^{-\frac12}_h}).
\end{equation}
Lemma~\ref{l:blowup} then implies
    \begin{equation}
    \begin{pmatrix}a_{-,j-1}^o\\a_{+,j}^o\end{pmatrix}=\Big(\Sigma_{M, h^{-1}\omega}^{\mc{W}_j}+O(M^2h^{2-0})\Big)\begin{pmatrix}a_{+,j-1}^i\\a_{-,j}^i\end{pmatrix}+O(h^{\frac{2\pi}{\alpha_+}-0}\|u\|_{L^2})+O(h^\infty)
    +O(h^{\frac{\pi}{\alpha_+}-0}\|v\|_{H^{-\frac12}_h}\Big),
    \end{equation}
where \eqref{e:8l11} and \eqref{e:9l5} imply
    \begin{equation}\label{e:9l13}
    \|\chi_j u\|_{H^{\frac12}_h} \le Ch^{\frac{2\pi}{\alpha_+}-0}\|u\|_{L^2}+C|v_+|+C|u_-|+C\|v\|_{H^{-\frac12}_h}.
    \end{equation}

6. We now glue up the scattering relations. Now note Steps 4 and 5 imply the scattering around corners and propagation on edges must agree:
\begin{multline}
    \begin{pmatrix}
    a_{-,J}^o\\ a_{+,1}^o\\a_{-,1}^o\\\vdots\\a_{+, J}^o
\end{pmatrix}=(\Sigma^{\mc{P}}_{M, h^{-1} \omega}+O(M^2h^{2-0}))\begin{pmatrix}
    a_{+,J}^i\\ a_{-,1}^i\\a_{+,1}^i\\\vdots\\a_{-, J}^i
\end{pmatrix}+O(h^{\frac{2\pi}{\alpha_+}-0}\|u\|_{L^2})+O(h^\infty)
    +O(h^{\frac{\pi}{\alpha_+}-0}\|v\|_{H^{-\frac12}_h}\Big)\\
    =T_{\omega/h}\begin{pmatrix}
    a_{+,J}^i\\ a_{-,1}^i\\a_{+,1}^i\\\vdots\\a_{-, J}^i
\end{pmatrix}+u_-+O(h^2M^2)+O(h^\infty),
\end{multline}
respectively. Thus we have
\begin{equation}\label{e:7l4}
    (\Sigma^{\mc{P}}_{M, h^{-1}\omega}-T_{\omega/h})v_+ + O(h^{\frac{2\pi}{\alpha_+}-0}\|u\|_{H^{-\frac12}_h}+h^{\frac{\pi}{\alpha_+}-0}\|v\|_{H^{-\frac12}_h})+O(M^2h^{2-0})+O(h^\infty)=(1+O(h^{\frac{2\pi}{\alpha_+}+1-0}))u_-.
\end{equation}

7. We now show $\mc{G}^{-1}(\omega;h)$ is uniformly bounded and obtain $E_{-+}(\omega;h)$. Note $\Sigma^{\mc{P}}_{M, h^{-1}\omega}-T_{\omega/h}$ is bounded and \eqref{e:7l4} implies
\begin{equation}\label{e:9l6}
    |u_-|\le C|v_+|+Ch^{\frac{2\pi}{\alpha_+}-0}\|u\|_{H^{-\frac12}_h}+Ch^{\frac{\pi}{\alpha_+}-0}\|v\|_{H^{-\frac12}_h}.
\end{equation}
Now note by \eqref{e:7l5} we have
\begin{multline}
    \|\gamma_j(u-\tilde u)\|_{H^{\frac12}_h}\le C\|\tilde a_{+,j}(x;h)\|_{S^0}+C\|\tilde a_{-,j}(x;h)\|_{S^0}+Ch^{\frac{\pi}{\alpha_+}+1-0}\|u\|_{H^{-\frac12}_h}+Ch^{1-0}\|v\|_{H^{-\frac12}_h}\\
    \le C|v_+|+Ch^{\max(\frac{2\pi}{\alpha_+},\frac{\pi}{\alpha_+}+1)-0}\|u\|_{H^{-\frac12}_h}+C\|v\|_{H^{-\frac12}_h}+O(h^\infty),
\end{multline}
where we used \eqref{e:9l5} and \eqref{e:9l6} to estimate $\tilde a_{\pm,j}(x;h)$. Now note that by \eqref{e:9l13}, \eqref{e:7l8} we have
\begin{equation}
    \|u\|_{H^{\frac12}_h}\le \sum_{j=1}^J\|\chi_j u\|_{H^{\frac12}_h}+\|\gamma_j \tilde u\|_{H^{\frac12}_h}+\|\gamma_j (u-\tilde u)\|_{H^{\frac12}_h}\le Ch^{\frac{\pi}{\alpha_+}}\|u\|_{L^2}+C\|v\|_{H^{-\frac12}_h}+C|v_+|,
\end{equation}
the first term on the right of which can be absorbed. This gives
\begin{equation}
    \|u\|_{H^{\frac12}_h}\le C\|v\|_{H^{-\frac12}_h}+C|v_+|.
\end{equation}
Then \eqref{e:7l4} reduces to 
\begin{multline}
    u_-=(\Sigma^{\mc{P}}_{M, h^{-1}\omega}-T_{\omega/h})v_+ + O(h^{\frac{2\pi}{\alpha_+}-0}|v_+|)+O(h^{\frac{\pi}{\alpha_+}-0}\|v\|_{H^{-\frac12}_h})+O(M^2h^{2-0})+O(h^\infty)\\
    =(\Sigma^{\mc{P}}_{M, h^{-1}\omega}-T_{\omega/h}+O(h^{\frac{2\pi}{\alpha_+}-0}+M^2h^{2-0}))v_++O(h^{\frac{\pi}{\alpha_+}-0}+M^2h^{2-0})\|v\|_{H^{-\frac12}_h})
\end{multline}
Here we used \eqref{e:9l6}. We then can see $\mc{G}^{-1}(\omega;h)$ given by 
\begin{equation}\label{e:7l25}
    \begin{pmatrix}
        u\\u_-
    \end{pmatrix}=\begin{pmatrix}
        O_{H^{\frac12}_h}\Big(\|v\|_{H^{-\frac12}_h}+|v_+|\Big)\\
        O\Big(h^{\frac{\pi}{\alpha_+}-0}+M^2h^{2-0}\Big)\|v\|_{H^{-\frac12}_h}+\Big(\Sigma^{\mc{P}w}_{M, h^{-1}\omega}-T_{\omega/h}+O(h^{\frac{2\pi}{\alpha_+}-0}+M^2h^{2-0})\Big)|v_+|
    \end{pmatrix}
\end{equation}
is bounded and
\begin{equation}
    E_{-+}(\omega; h)v_+= \Big(\Sigma^{\mc{P}}_{M, h^{-1}\omega}-T_{\omega/h}+O(h^{\frac{2\pi}{\alpha_+}-0}+M^2h^{2-0})\Big)v_+=u_-
\end{equation}
as desired. This gives the invertibility of $\mc{G}(\omega;h)$ for all $\omega\in[\frac12,2]$. 

8. We now show the invertibility of $\mc{G}(z;h)$ holds in a complex neighbourhood of $\omega\in [\frac12, 2]$. Note one can check for real $\mu$ close to $0$ we have
\begin{equation}
    \mc{G}(z;h)=\mc{G}(\omega;h)-\mc{A}(z;h), \ \mc{A}(z;h)=O_{H_h^{\frac12}\oplus \mathbb{C}^{2J}\rightarrow H_h^{-\frac12}\oplus \mathbb{C}^{2J}}(|\mu|),
\end{equation}
from the definition of $R_\pm(z)$ and
\begin{equation}
    \mc{A}(z;h)=\begin{pmatrix}
        \sqrt{h^{-2}(\omega+ih\mu)^2+M^2}-\sqrt{h^{-2}\omega^2+M^2} & R_-(\omega)-R_-(z)\\
        R_+(\omega)-R_+(z) & 0
    \end{pmatrix}.
\end{equation}
Now \eqref{e:7l25} implies $\|\mc{G}^{-1}(\omega;h)\|$ is bounded uniformly in $\omega\in [\frac12,2]$, $M\in[0, M_0]$, $h\in (0, h_0]$. This implies there is $r_1>0$ small that uniformly for any $z=\omega+ih\mu\in [\frac12, 2]+i[-r_1 h, r_1 h]$ we have $\|\mc{A}(z;h)\mc{G}^{-1}(\omega;h)\|<1$. Now we invert $\mc{G}(z;h)$ with the Neumann series:
\begin{equation}
    \mc{G}^{-1}(z;h)=\mc{G}^{-1}(\omega;h)\Big(\Id+\sum_{k=1}^{\infty}(\mc{A}(z;h)\mc{G}^{-1}(\omega;h))^{k}\Big)=\mc{G}^{-1}(\omega;h)+O_{H_h^{-\frac12}\oplus \mathbb{C}^{2J}\rightarrow H_h^{\frac12}\oplus \mathbb{C}^{2J}}(|\mu|).
\end{equation}
This implies
\begin{equation}
    E_{-+}(z;h)=E_{-+}(\omega;h)+O(|\mu|)=\Sigma^{\mc{P}}_{M,h^{-1}\omega}-T_{\omega/h}+O(h^{\frac{2\pi}{\alpha_+}-0}+M^2h^{2-0})+O(|\mu|).
\end{equation}
That $E_{-+}(z;h)$ is analytic Fredholm in $z$ follows from $\mc{G}(z;h)$ is so.
\end{proof}

\begin{remark}[Scattering matrices of curvature order $k_0\ge 0$]\label{r:curvedscmatrix}
    Let $\mc{W}_{\alpha_j}^{k_0}$ be some 2D curved wedge of corner angle $\alpha_j$ centered at $x_j$, with first $(k_0+1)$ derivatives of the boundary at $x_j$ agreeing with those for $\mc{P}$ at $x_j$. Proposition~\ref{p:existscmatrix} implies the existence of scattering matrices $\Sigma_{M,h^{-1}\omega}^{\mc{W}_{\alpha_j}^{k_0}}$. If $\mc{P}$ has all derivatives of order $2\le k\le (k_0+1)$ of the boundary at $x_j$ vanishes, then 
\begin{equation}
    \mc{W}_{\alpha_j}^{k}=W_{\alpha_j}, \ \Sigma_{M,h^{-1}\omega}^{\mc{W}_{\alpha_j}^{k_0}}=\Sigma_{M,h^{-1}\omega}^{W_{\alpha_j}}=\Sigma_{hM,\omega}^{W_{\alpha_j}}, \text{ for each } 2\le k\le (k_0+1)
\end{equation}
by the scaling invariance of the straight wedge $W_{\alpha_j}$. 
\end{remark}

Proposition~\ref{p:scconverge} allows us to only use fixed scattering matrices on $\mc{W}_{\alpha_j}^{k_0}$ to define the effective Hamiltonian. Indeed, let for $\omega>0$ large.
\begin{equation}
\Sigma^{k_0}(\omega)=\begin{pmatrix}\Sigma^{\mc{W}_{\alpha_1}^{k_0}}_{0,\omega}&0&\dots &0\\
0&\Sigma^{\mc{W}_{\alpha_2}^{k_0}}_{0,\omega}&\dots &0\\
\vdots&\vdots&\ddots&\vdots\\
0&0&\dots&\Sigma^{\mc{W}_{\alpha_J}^{k_0}}_{0,\omega}
\end{pmatrix},\ \Sigma^0=\begin{pmatrix}\Sigma^{\mc{W}_{\alpha_1}}_{0,1}&0&\dots &0\\
0&\Sigma^{\mc{W}_{\alpha_2}}_{0,1}&\dots &0\\
\vdots&\vdots&\ddots&\vdots\\
0&0&\dots&\Sigma^{\mc{W}_{\alpha_J}}_{0,1}
\end{pmatrix}.
\end{equation}

\begin{corollary}\label{c:effhamdtn}
There exist $C, h_0>0$ such that uniformly for $\omega\in [\frac12, 2]$, and $h\in(0,h_0]$, we have $\dim \ker (\Lambda_M^\mc{P}-h^{-1}\omega)=\dim \ker\mc{Q}_{h^{-1}\omega}^{\Lambda}$, where
\begin{equation}
    \mc{Q}_{h^{-1}\omega}^\Lambda
        =E_{-+}(\omega;h)=\Sigma^{k_0}(h^{-1}\omega)-T_{\omega/h}+O\Big(h^{\min(\frac{2\pi}{\alpha_+},1+k_0)-0}+(Mh)^{\min(\frac{2\pi}{\alpha_+},2)-0}\Big).
\end{equation}
\end{corollary}
\begin{proof}
    Proposition~\ref{p:grushin} implies that $\dim \ker(\Lambda^{\mc{P}}_M-h^{-1}z)=\dim \ker E_{-+}(z; h)$. It suffices to note Proposition~\ref{p:scconverge} implies
    \begin{equation}
        \Sigma^{k_0}(h^{-1}\omega)-\Sigma^{\mc{P}}_{M, h^{-1}\omega}=O\Big(h^{\min(\frac{2\pi}{\alpha_+},1+k_0)-0}+(Mh)^{\min(\frac{2\pi}{\alpha_+},2)-0}\Big)
    \end{equation}
    and thus
    \begin{multline}
        \mc{Q}_{h^{-1}\omega}^\Lambda=E_{-+}(\omega; h)=\Sigma^{\mc{P}}_{M, h^{-1}\omega}-T_{\omega/h}+O(h^{\frac{2\pi}{\alpha_+}-0})+O(h^2M^2)\\
        =\Sigma^{k_0}(h^{-1}\omega)-T_{\omega/h}+O\Big(h^{\min(\frac{2\pi}{\alpha_+},1+k_0)-0}+(Mh)^{\min(\frac{2\pi}{\alpha_+},2)-0}\Big),
    \end{multline}
    as desired. 
\end{proof}

\begin{corollary}
\label{c:derivatives}
    There exist $ h_0,r_1>0$ such that for all $\e>0$, $\ell \geq 0$, there is $C>0$ such that the uniformly for $z\in [\frac34, \frac 32]$, and $h\in(0,h_0]$,
    $$
    |\partial^\ell_z (Q_{h^{-1}z}^\Lambda+T_{z/h})|\leq Ch^{2-\ell-\e} 
    $$
\end{corollary}
\begin{proof}
    Observe that $\Sigma^{0}(z)=\Sigma^0$ is independent of $z$ and thus, by Corollary~\ref{c:effhamdtn},
    $$
    |f(h,z)|\leq Ch^{1-0}+C|\Im z|/h,\qquad f(h,z):=Q_{h^{-1}z}^\Lambda+T_{z/h}-\Sigma^0
    $$
    In particular since $f$ is analytic, for $z\in[\frac{1}{2},2]+i[-r_1h,r_1h]$, letting 
    $$
    \Gamma_\e= 2+[-\e h,\e h]\cup [2,\tfrac{1}{2}]+i\e h\cup \tfrac{1}{2}+i[\e h,-\e h]\cup [\tfrac{1}{2},2]-i\e h
    $$
    \begin{align*}
    \big|\partial_z^\ell f(z)\big|&=\frac{\ell!}{2\pi }\Big|\int_{\Gamma_\e}\frac{f(s)}{(s-z)^{\ell+1}}ds\Big|\\
    &\leq C \int_{-\e h}^{\e h} (h^{1-0}+s/h)ds+C\int_{-2}^2\frac{h^{1-0}+\e}{(\sqrt{s^2+\e^2h^2})^{\ell+1}}ds.\\
    &\leq C(\e h^{2-0}+\e^2h^{1-0})+\begin{cases}C(h^{1-0} +\e)\log h&\ell=1\\
    C(h^{1-0}+\e)h^{1-\ell}&\ell \geq 2.\end{cases}
    \end{align*}
    Hence, choosing $\e=h$, we obtain
    $$
    |\partial^{\ell}_z f(z)|\leq Ch^{2-\ell-0},
    $$
    which implies the corollary.
\end{proof}
\begin{remark}
    All constants in this section can depend on $\mc{P}$ in the following way: the upper bound on the number of corners $J$, the upper and lower bounds over the angles $\alpha_j$, the upper and lower bounds on all $|\Gamma_j|$, the derivatives of the defining functions of $\Gamma_j$. 
\end{remark}

\section{The Steklov Dirac operator and approximation on curvilinear polygons}
\label{s:spectralAsymptotics}

Using the Grushin problem from Section~\ref{s:grushin}, we prove Theorem~\ref{t:boundedSteklov} in this section. The key remaining issues are i) to show that the values of $\omega$ that $\mc{Q}^\Lambda_\omega=E_{-+}(\omega)$ from~\eqref{e:7l14} is not invertible are close, with multiplicity, to those of $\mc{Q}^{\mathscr{D}}(\omega)=\Sigma^{k_0}(\omega)-T_{\omega/h}$ from Corollary~\ref{c:effhamdtn} (which in turn is the Steklov Driac operator with matching from Theorem~\ref{t:boundedSteklov} and ii) to prove that the labeling of the eigenvalues for the Dirac operator is identical to that for the Steklov problem.

\noindent{\bf {Outline of the strategy}}

\noindent Closeness of the spectra:

To prove that the spectra of $\mc{Q}^{\Lambda}_\omega$ and $\mc{Q}^{\mathscr{D}}_\omega$ are close, with multiplicity, we extend $\mc{Q}^{\mathscr{D}}_\omega$ off of the real axis so that $\bar{\partial}\mc{Q}^{\mathscr{D}}_z|_{|\Im z|=0}=0$ and use a version of Rouche's theorem valid for smooth functions. The key features that make this possible are that $\Sigma^{k_0}$ is unitary and has small derivatives. This allows us to show that the eigenvalues of $(\Sigma^{k_0}(\omega))^*T_\omega$ wind monotonically and with strictly positive speed around the unit. Thus has two important consequences. First, there are at most $2J$ roots of $\det \mc{Q}^{\mathscr{D}}_\omega$ in fixed size intervals and hence the spectrum has many gaps, and second, that the extension of $\mc{Q}_\omega^{\mathscr{D}}$ off the real axis is well behaved. The closeness of $\mc{Q}^{\mathscr{D}}_z$ and the extension of $\mc{Q}^{\Lambda}_\omega$ then persists off the real axis, making the Rouché argument possible.

\noindent Eigenvalue labeling:

After proving that the roots of $\det \mc{Q}^{\Lambda}_z$ are close to the eigenvalues of the Steklov problem, we have established that, up to a possible shift of the index, the eigenvalues, $\lambda_n$ of the Dirac operator are close to those of the Steklov problem, $\sigma_n$ i.e. there is $N\in\mathbb{Z}$ such that
\begin{equation}
    |\sigma_n-\lambda_{n+N}|\leq C\sigma_n^{-\min(\frac{2\pi}{\alpha_+},1+k_0)+\e}+CM^2\sigma_n^{-2+\e},\qquad \sigma_n>\sigma_0.
\end{equation}
It thus remains to show that $N=0$.  

To prove that $N=0$ we make a deformation argument. Specifically, we observe that $\sigma_0$ is uniform over $M\in[0,M_0]$, for curvilinear polygons whose angles are uniformly bounded away from $0$ and $2\pi$, and for which the $C^k$ norms of the the edges of the polygon are uniformly bounded. Thus, we may deform any curvilinear polygon with $M\geq 0$  to a smooth domain with $M=0$ whose boundary is linear near each (fictitious) corner and without changing $\sigma_0$. This implies that $L$ is constant among all domains with the same topology and for all $M\geq 0$. It thus remains to observe that the results of~\cite{GPPS2014} (see also~\cite{Rozenblyum1986} for the simply connected case) imply $L=0$ for smooth domains. 

\subsection{The graph Dirac operator on curvilinear polygons}
Given $\mc{P}\Subset \mathbb{R}^2$, a connected Lipschitz domain with connected piecewise smooth boundary, we parametrize $\partial\mc{P}$ by $\mathbb{S}^1$ and the vertices by $x_j\in \mathbb{S}^1$ for $1\le j\le J$: here we understand it as a graph with edges $\Gamma_{j}$ connecting $x_j$ to $x_{j+1}$ of length $L_j$. Let $\mc{W}_{\alpha_j}^{k_0}$ be defined in Remark~\ref{r:curvedscmatrix}. Let $\psi_0\in C^\infty_c(\mathbb{R})$ be even and $\psi_0(\omega)\equiv 1$ on $[-\frac12,\frac12]$, and $\psi_0(\omega)\equiv 1$ on $|\omega|\ge 1$. Let $\delta>0$ to be chosen small and define the regularized scattering matrix
\begin{equation}
\label{e:scatDefine}
    \Sigma_{j}^{k_0}(\omega)=\exp\Big((1-\psi_0)(\delta\omega)\log\big(\Sigma^{\mc{W}_{\alpha_j}^{k_0}}_{0,\omega}(\Sigma^{W_{\alpha_j}}_{0,1})^*\big)\Big)\Sigma^{W_{\alpha_j}}_{0,1}, \ \omega\ge0, 
\end{equation}
and $\Sigma_{j}^{k_0}(\omega)=-\sigma_3\overline{\Sigma_{j}^{k_0}(-\omega)}\sigma_3$ for $\omega<0$, where
\begin{equation}
    \sigma_3=\begin{pmatrix}1&0\\0&-1\end{pmatrix}
\end{equation}
is the third Pauli matrix. Note $\Sigma_{j}^{k_0}(\omega)$ is unitary on all of $\omega\in\mathbb{R}$, and $\Sigma^0_j(\omega)=\Sigma^{W_{\alpha_j}}_{0,1}$ on $\omega\in\mathbb{R}$, where
\begin{equation}\label{e:scmatrixalpha}
    \Sigma^{W_{\alpha_{j}}}_{0,1} =\begin{pmatrix}
    -i\cos(\frac{\pi^2}{2\alpha})&\sin(\frac{\pi^2}{2\alpha})\\ \sin(\frac{\pi^2}{2\alpha})&-i\cos(\frac{\pi^2}{2\alpha})
\end{pmatrix},
\end{equation}
given by Lemma~\ref{l:scatMatrix}. One can check $\Sigma_{0,1}^{W_{\alpha_j}}=-\sigma_3\overline{\Sigma_{0,1}^{W_{\alpha_j}}}\sigma_3$ and thus $\Sigma_{j}^{k_0}(\omega)$ is identically $\Sigma_{0,1}^{W_{\alpha_j}}$ near $\omega=0$. 

Consider the Steklov Dirac operator $\mathscr{D}_\omega^{k_0}=\begin{pmatrix} D_x&0\\0&-D_x\end{pmatrix}$ with domain
\begin{equation}\label{e:7l11}
    \dom(\mathscr{D}_\omega^{k_0})=\Big\{ u=\begin{pmatrix}
    u_+\\ 
    u_-
\end{pmatrix}\in H^1(\partial\mc{P};\mathbb{C}^2)\,:\, \Sigma_{j}^{k_0}(\omega)\begin{pmatrix}u_+(x_j-0)\\u_-(x_j+0)\end{pmatrix}=\begin{pmatrix}u_-(x_j-0)\\u_+(x_j+0)\end{pmatrix}\Big\}\subset L^2(\partial\mc{P}; \mathbb{C}^2),
\end{equation}
where the Dirichlet-type boundary condition records the scattering happening at vertices $x_j$'s. We may drop the index $\omega$ and $k_0$ whenever there is no confusion. 

For $\omega\in\mathbb{C}$, we define the effective Hamiltonian $\mc{Q}_\omega^{\mathscr{D}}$ by, 
\begin{equation}
\mc{Q}_\omega^{\mathscr{D}}=\Sigma^{k_0}(\omega)-T_\omega, \ \Sigma^{k_0}(\omega)=\begin{pmatrix}\Sigma_{1}^{k_0}(\omega)&0&\dots &0\\
0&\Sigma_{2}^{k_0}(\omega)&\dots &0\\
\vdots&\vdots&\ddots&\vdots\\
0&0&\dots&\Sigma_{J}^{k_0}(\omega)
\end{pmatrix},
\end{equation}
where $\Sigma^{k_0}(\omega)$ tracks the scattering at vertices and the transfer operator $T_z$ for $z\in\mathbb{C}$ tracks the propagation of states at complex frequency $z$ along each edge, 
\begin{equation}
    T_z=T_-\begin{pmatrix}E_{L_1}^z&0&\dots &0\\
0&E_{L_2}^z&\dots &0\\
\vdots&\vdots&\ddots&\vdots\\
0&0&\dots&E_{L_J}^z
\end{pmatrix} T_+, \ 
E_{L_j}^z=\begin{pmatrix}
    0 & e^{-izL_j}\\
    e^{-izL_j} & 0
\end{pmatrix}.
\end{equation}
We remind the readers that $\Sigma^{k_0}(\omega)$ and $\mc{Q}_\omega^{\mathscr{D}}$ are not necessarily analytic in $\omega$ unless $k_0=0$. 

Set $S_3=\operatorname{diag}(\sigma_3,\dots, \sigma_3)\in \mathbb{C}^{2J\times 2J}$. Note
\begin{equation}
    S_3T_+S_3=-T_+, \ S_3T_-S_3=-T_-, \ S_3 T_z S_3=-T_z.
\end{equation}
The above implies
\begin{equation}
    S_3 \overline{\mc{Q}^{\mathscr{D}}_\omega} S_3=S_3 (\overline{\Sigma^{k_0}(\omega)-T_\omega})S_3=-(\Sigma^{k_0}(-\omega)-T_{-\omega})=-\mc{Q}^{\mathscr{D}}_{-\omega}.
\end{equation}
Taking the determinant gives
\begin{equation}
    \overline{\det \mc{Q}^{\mathscr{D}}_\omega}=\det S_3 \overline{\mc{Q}^{\mathscr{D}}_\omega} S_3=\det (-\mc{Q}^{\mathscr{D}}_{-\omega})=\det \mc{Q}^{\mathscr{D}}_{-\omega}.
\end{equation}
This implies the zeros of $\det \mc{Q}^{\mathscr{D}}_\omega$ are symmetric about $\omega=0$, and if $\det \mc{Q}_0^{\mathscr{D}}=0$, $\omega=0$ is a zero with even multiplicity. We define the zero set of $\det \mc{Q}_{\omega}^{\mathscr{D}}$ and $\det \mc{Q}_{\omega}^{\Lambda}$ by
\begin{equation}
    \mc{Z}^{\mathscr{D}}:=\{\omega\in\mathbb{R}\,:\, \det \mc{Q}_{\omega}^{\mathscr{D}}=0\},\ \mc{Z}^{\Lambda}:=\{\omega\in\mathbb{R}\,:\, \det \mc{Q}_{\omega}^{\Lambda}=0\}
\end{equation}
and their associated counting functions by 
\begin{equation}
    \mc{N}^{\mathscr{D}}_{I}:=\sum_{\omega\in I}\dim \ker \mc{Q}_{\omega}^{\mathscr{D}},\  \mc{N}^{\Lambda}_{I}:=\sum_{\omega\in I}\dim \ker \mc{Q}_{\omega}^{\Lambda}, \ I\subset\mathbb{R}.
\end{equation}

In the following lemma, we show that $\ker \mc{Q}_\omega^{\mathscr{D}}$ characterizes $\operatorname{char} (\mathscr{D}_\omega^{k_0})$.
\begin{lemma}[The Steklov Dirac operator]\label{l:graphdirac}
    The Steklov Dirac operator $\mathscr{D}^{k_0}_\omega$ is self-adjoint. Moreover, let $\omega\in\mathbb{R}$ and 
    \begin{equation}
        \mathbf{a}=(a^i_{+,J}, a^{i}_{-, 1}, a^{i}_{+,1}, \cdots, a^i_{-, J})^T\in \mathbb{C}^{2J}.
    \end{equation}
    Then we have $\mc{Q}_{\omega}^{\mathscr{D}}\mathbf{a}=0$ if and only if $u_{\mathbf{a}}\in H^1(\partial\mc{P};\mathbb{C}^2)$ defined by
    \begin{equation}
    u_{\mathbf{a}}|_{(x_j, x_{j+1})}=\begin{pmatrix}
        a^i_{+, j} e^{iz(x-x_{j+1})}\\
        a^i_{-, j} e^{-iz(x-x_j)}
    \end{pmatrix}
    \end{equation}
    is a solution to $(\mathscr{D}^{k_0}_\omega-\omega)u_{\mathbf{a}}=0$. 
\end{lemma}
\begin{proof}
    1. Let $z=\omega+i\mu$ and we formally solve $(\mathscr{D}-z)u=f$. Given $f\in L^2(\partial\mc{P};\mathbb{C}^2)$, we want to find $u\in \dom(\mathscr{D})$ that
\begin{equation}
    (\mathscr{D}-z)u=(\mathscr{D}-z)\begin{pmatrix}
        u_+\\u_-
    \end{pmatrix}=\begin{pmatrix}
        (D_x-z)u_+\\(-D_x-z)u_-
    \end{pmatrix} =\begin{pmatrix}
        f_+\\f_-
    \end{pmatrix}=f.
\end{equation}
We will denote the end data of $u$ near vertex $x_j$ by 
\begin{equation}
    a^i_{+, j-1}=u_+(x_j-0), \ a^o_{-, j-1}=u_-(x_j-0), \ a^o_{+, j}=u_+(x_j+0), \ a^i_{-, j}=u_-(x_j+0)
\end{equation}
and will solve for those data that matches the scattering relation governed by the boundary condition. Note that on edge $\Gamma_j$, we have two ordinary differential equations for $u_\pm$ and they are solved by
\begin{gather}\label{e:7l16}
    u_+(x)=a^i_{+, j}e^{iz(x-x_{j+1})}-i\int_x^{x_{j+1}} e^{iz(x-s)}f_+(s)\,ds, \\
    u_-(x)=a^i_{-,j}e^{-iz(x-x_j)}-i\int_{x_j}^x e^{-iz(x-s)}f_-(s)\,ds,\label{e:7l17}
\end{gather}
for $x\in \Gamma_j=(x_j, x_{j+1})$, where $a^i_{\pm, j}$ unique determines $u_\pm$. Notably we have 
\begin{equation}\label{e:7l10}
    a_{+, j}^i=e^{iz L_j}(a_{+,j}^o-F_{+,j}), \ a_{-, j}^o=e^{-iz L_j}a_{-, j}^i+F_{-, j},
\end{equation}
where 
\begin{equation}
    F_{+, j}=-i\int_{x_j}^{x_{j+1}} e^{iz (x_{j}-s)}f_+(s)\,ds,\ F_{-, j}=-ie^{-iz L_j}\int_{x_j}^{x_{j+1}} e^{-iz(x_j-s)}f_-(s)\,ds
\end{equation}
are purely determined by the data $f$. Now note $u\in\dom(\mathscr{D})$ if and only if the scattering relations per \eqref{e:7l11} holds, that is
\begin{equation}\label{e:7l18}
    \begin{pmatrix}
    a_{-,J}^o\\ a_{+,1}^o\\a_{-,1}^o\\\vdots\\a_{+, J}^o
\end{pmatrix}=\Sigma^{k_0}(\omega)\begin{pmatrix}
    a_{+,J}^i\\ a_{-,1}^i\\a_{+,1}^i\\\vdots\\a_{-, J}^i
\end{pmatrix}=T_z\begin{pmatrix}
    a_{+,J}^i\\ a_{-,1}^i\\a_{+,1}^i\\\vdots\\a_{-, J}^i
\end{pmatrix}+\begin{pmatrix}
    F_{-,J}\\
    F_{+,1}\\
    F_{-,1}\\
    \vdots\\
    F_{+, J}
\end{pmatrix}=T_z\mathbf{a}+\mathbf{F}.
\end{equation}
That is, $(\mathscr{D}-z)u=f$ for $u\in\dom(\mathscr{D})$, if and only if, $u$ is given by \eqref{e:7l16} and \eqref{e:7l17} and the scattering relation 
$(\Sigma^{k_0}(\omega)-T_z)\mathbf{a}=\mathbf{F}$ holds. For $z=\omega\in\mathbb{R}$, this means $\mc{Q}_\omega^{\mathscr{D}} \mathbf{a}=(\Sigma^{k_0}(\omega)-T_\omega)\mathbf{a}=\mathbf{F}$.

2. We show $(\Sigma^{k_0}(\omega)-T_z)^{-1}$ exists for $z=\omega+i\mu\in\mathbb{C}\setminus \mathbb{R}$. Note $\Sigma^{k_0}(\omega)$ is unitary. When $\mu>0$, we have $\|T_z^{-1}\|=e^{-\mu\min_{j}\{L_j\}}<1$, and thus for each $\mathbf{a}\in\mathbb{C}^{2J}$, we have
\begin{equation}
    |(\Sigma^{k_0}(\omega)-T_z)\mathbf{a}|\ge (e^{\mu\min_{j}\{ L_j\}}-1)|\mathbf{a}|\ge C^{-1}\min(\mu\min_{j}(L_j),1)e^{\mu\min_j(L_j)}|\mathbf{a}|
\end{equation}
as desired. When $\mu<0$, we have $\|T_z\|=e^{\mu\min_{j}\{L_j\}}<1$ and then
\begin{equation}
    |(\Sigma^{k_0}(\omega)-T_z)\mathbf{a}|\ge (1-e^{\mu\min_{j}\{ L_j\}})|\mathbf{a}|\ge C^{-1}\min(|\mu|\min_{j}(L_j),1)|\mathbf{a}|.
\end{equation}
The above shows for $z\in\mathbb{C}\setminus \mathbb{R}$, $(\Sigma^{k_0}(\omega)-T_z)^{-1}$ exists and by Step 1, $\mathscr{D}-z$ is surjective. 

3. We show $\mathscr{D}$ is self-adjoint. First note $\mathscr{D}$ is symmetric, since for $u,v\in \dom(\mathscr{D})$, 
\begin{multline}
    \langle \mathscr{D}u,v\rangle-\langle u,\mathscr{D}v\rangle=-i\sum_{j=1}^J \Big(u_+(x_{j+1}-0)\overline{v_+(x_{j+1}-0)}-u_+(x_{j}+0)\overline{v_+(x_{j}+0)}-u_-(x_{j+1}-0)\overline{v_-(x_{j+1}-0)}\\+u_-(x_{j}+0)\overline{v_-(x_{j}+0)}\Big)
    =-i\sum_{j=1}^J \Big(u_+(x_{j}-0)\overline{v_+(x_{j}-0)}-u_+(x_{j}+0)\overline{v_+(x_{j}+0)}-u_-(x_{j}-0)\overline{v_-(x_{j}-0)}\\
    +u_-(x_{j}+0)\overline{v_-(x_{j}+0)}\Big)
    =-i\sum_{j=1}^J \Big(\Big\langle \begin{pmatrix}
        u_+(x_j-0)\\u_-(x_j+0)
    \end{pmatrix},\begin{pmatrix}
        v_+(x_j-0)\\v_-(x_j+0)
    \end{pmatrix}\Big\rangle-\Big\langle \begin{pmatrix}
        u_-(x_j-0)\\u_+(x_j+0)
    \end{pmatrix},\begin{pmatrix}
        v_-(x_j-0)\\v_+(x_j+0)
    \end{pmatrix}\Big\rangle\Big)\\
    =-i\sum_{j=1}^J \Big(\Big\langle \begin{pmatrix}
        u_+(x_j-0)\\u_-(x_j+0)
    \end{pmatrix},\begin{pmatrix}
        v_+(x_j-0)\\v_-(x_j+0)
    \end{pmatrix}\Big\rangle-\Big\langle \Sigma_{j}^{k_0}(\omega)\begin{pmatrix}
        u_+(x_j-0)\\u_-(x_j+0)
    \end{pmatrix},\Sigma_{j}^{k_0}(\omega)\begin{pmatrix}
        v_+(x_j-0)\\v_-(x_j+0)
    \end{pmatrix}\Big\rangle\Big)\\
    =0,
\end{multline}
the last line of which follows from $\Sigma_{j}^{k_0}(\omega)$ being unitary. By Step 2, $\operatorname{Ran}(\mathscr{D}\pm i)=L^2(\partial\mc{P};\mathbb{C}^2)$ and thus $\mathscr{D}$ is self-adjoint. 

4. Let $\mc{Q}_z^{\mathscr{D}}\mathbf{a}=0$, then Step 1 immediately implies that $u_{\mathbf{a}}$ is a solution to $(\mathscr{D}-z)u_{\mathbf{a}}=0$. Let $u_{\mathbf{a}}\in\dom(\mathscr{D})$ be that $(\mathscr{D}-z)u_{\mathbf{a}}=0$, then by \eqref{e:7l16} and \eqref{e:7l17}, 
\begin{equation}
    u_{+}=a^i_{+, j}e^{iz(x-x_{j+1})},\ u_{-}=a^i_{-, j}e^{-iz(x-x_{j})},
\end{equation}
and $a^i_{\pm, j}$ must satisfy the scattering relation \eqref{e:7l18}, that is, $\mc{Q}_z^{\mathscr{D}}\mathbf{a}=0$. 
\end{proof}

\subsection{Properties of $\mc{Q}^{\mathscr{D}}_\omega$}

In this subsection, we collect some important facts about $\mc{Q}^{\mathscr{D}}_\omega$ that will be used in the proof of Theorem~\ref{t:boundedSteklov}. All estimates in this section are uniform in $0<\delta<1$ unless explicitly specified.

By Proposition~\ref{p:scconverge}, we have
\begin{equation}
    \Sigma^{k_0}(\omega)-\Sigma^0=O(|\omega|^{-1+0}),
\end{equation}
and, by Lemma~\ref{l:scatteringDerivatives}, we have
\begin{equation}
    \|\partial_{\omega}^2\Sigma^{k_0}\|+\|\partial_{\omega}\Sigma^{k_0}\|\leq C. 
\end{equation}
Therefore,
\begin{equation} \label{e:diracDerivative}
\|\partial_\omega\Sigma^{k_0}\|\le C(|\omega|^{-\frac{1}{2}+0}). 
\end{equation}
Denote for convenience
\begin{equation}
    U(\omega)=\Sigma^{k_0}(\omega),\ T(\omega)=T_\omega, \ W(\omega)=U^*(\omega)T(\omega), \ A_{W}(\omega)=-iW^*(\omega)W'(\omega).
\end{equation}
Here $U, T, W$ are unitary matrices in $\mathbb{C}^{2J\times 2J}$. We then have
\begin{equation}
    \mc{Q}_{\omega}^{\mathscr{D}}= U(\omega)(I-W(\omega)).
\end{equation}
Specifically we have
\begin{equation}
    A_T(\omega):=-iT^*(\omega)T'(\omega)=-iT_- \operatorname{diag}_j((E_{L_j}^\omega)^*(E_{L_j}^\omega)')T_+=-T_- \operatorname{diag}_j\begin{pmatrix}
        L_j & 0\\
        0 & L_j
    \end{pmatrix}T_+,
\end{equation}
since $(E_{L_j}^\omega)'=-iL_j E_{L_j}^\omega$. Then we have
\begin{equation}
    - L_+\Id\le A_T(\omega)\le - L_-\Id, \  L_-=\min_j(L_j), \ L_+=\max_j(L_j).
\end{equation}
Note that for $\omega$ large, 
\begin{equation}
\label{e:positiveRotation}
    -2L_+\Id\le A_W=A_T-iT^*U(U')^*T\le -\frac{L_-}2\Id,
\end{equation}
since $U'(\omega)=o(1)$ by~\eqref{e:diracDerivative}.

Since $W(\omega)$ is unitary, it has on orthonormal basis of eigenvectors $\{v_j(\omega)\}_{j=1}^{2J}$, with eigenvalues $\{e^{i\theta_j(\omega)}\}_{j=1}^{2J}$, for $\theta_j(\omega)\in\mathbb{R}$ and, moreover, \eqref{e:positiveRotation} implies that $-2L_+\leq \theta_j'(\omega)\leq -\frac{L_-}{2}$. In particular, there is $c>0$ that 
\begin{equation} \label{e:eigenvalueControl}\theta_j(\omega_0)\in 2\pi \mathbb{Z}\quad\Rightarrow\quad  
\dist(\theta_j(\omega),2\pi\mathbb{Z})\geq c|\omega-\omega_0|\text{ for }|\omega-\omega_0|<c.
\end{equation}
\begin{lemma}[Root density for $\det\mc{Q}_z^{\mathscr{D}}$]\label{l:rootnumber}
There is $r_0>0$ such that for any $\omega_0\in\mathbb{R}$, 
\begin{equation}
    \mc{N}^{\mathscr{D}}_{(\omega_0-r_0,\omega_0+r_0)}=\#\{\omega\in(\omega_0-r_0,\omega_0+r_0): \det\mc{Q}_\omega^{\mathscr{D}}=0\}\le 2J,
\end{equation}
counting multiplicities. 
\end{lemma}
\begin{proof}
    Observe that, since $\mc{Q}_{\omega}^{\mathscr{D}}= \Sigma^{k_0}(\omega)(I-W(\omega))$,
    \begin{equation}
        \mc{N}^{\mathscr{D}}_{(\omega_0-r_0,\omega_0+r_0)}=\#\{(\omega,j)\in (\omega_0-r_0,\omega_0+r_0)\times \{1,\dots,2J\}\,:\, \theta_j(x)\in 2\pi \mathbb{Z}\}.
    \end{equation}
    Thus, the lemma follows from~\eqref{e:eigenvalueControl}.
\end{proof}

\begin{lemma}[Control on $(\mc{Q}_{\omega}^{\mathscr{D}})^{-1}$]
    There is $C>0$ such that for all 
    \begin{equation}
    \label{e:distanceInverseEstimate}
    \|(\mc{Q}_{\omega}^{\mathscr{D}})^{-1}\|\leq C\dist\big(\omega,\mc{Z}^{\mathscr{D}}\big)^{-1}.
    \end{equation}
\end{lemma}
\begin{proof}
    Observe that 
    \begin{equation}
        (\mc{Q}_{\omega}^{\mathscr{D}})^{-1}=(\Id-W(\omega))^{-1}U(\omega)^*.
    \end{equation}
    Therefore, 
    \begin{equation}
        \|(\mc{Q}_{\omega}^{\mathscr{D}})^{-1}\|\leq C\inf_{j}\dist(\theta_j,2\pi \mathbb{Z})^{-1}.
    \end{equation}
    The lemma now follows from~\eqref{e:eigenvalueControl}.
\end{proof}

We build $\tilde{\mc{Q}}_z^{\mathscr{D}}$, the almost analytic extension of $\mc{Q}_\omega^{\mathscr{D}}$ of order 1 to $z\in\mathbb{C}$ by
\begin{equation}\label{e:Qextend}
    \tilde{\mc{Q}}_{\omega+i\mu}^{\mathscr{D}}:=\mc{Q}_\omega^{\mathscr{D}}+i\mu \partial_\omega\mc{Q}_\omega^{\mathscr{D}}.
\end{equation}
Notice, in particular, that $\bar{\partial}_{\omega+i\mu}\mc{Q}^{\mathscr{D}}_{\omega}=0$ for $\omega\in\mathbb{R}$. 
\begin{lemma}[Structure of $ \tilde{\mc{Q}}_{z}^{\mathscr{D}}$ near a root]
\label{l:identify}
    Let $\omega_0\in\mathbb{R}$ and define $m:=\dim \ker \mc{Q}_{\omega_0}^{\mathscr{D}}$.
    Then, there is $a\in\mathbb{C}$, $a\neq 0$ such that 
    \begin{equation}
        \det \tilde{\mc{Q}}_{z}^{\mathscr{D}}= a(z-\omega_0)^m+O(|z-\omega_0|^{m+1}).
    \end{equation}
\end{lemma}
\begin{proof}
Observe that for $|\mu|\leq C$,
\begin{equation}
    \|(\Sigma^{k_0}(\omega)+i\mu \partial_{\omega}\Sigma^{k_0}(\omega))^{-1}\|\leq C.
\end{equation}
Extend $W(\omega)$ by
\begin{equation}
    \tilde{W}(\omega+i\mu)=(\Sigma^{k_0}(\omega)+i\mu \partial_{\omega}\Sigma^{k_0}(\omega))^{-1}(T_{\omega}+i\mu\partial_\omega T_\omega), 
\end{equation}
so that $\tilde{W}(\omega)=W(\omega)$ and
\begin{equation}
    \det \tilde{\mc{Q}}_{\omega+i\mu}^{\mathscr{D}}=\det (\Sigma^{k_0}(\omega)+i\mu \partial_{\omega}\Sigma^{k_0}(\omega))\det (\Id-\tilde W(\omega+i\mu))=e(\omega,\mu) \det (\Id-\tilde{W}),
\end{equation}
where $e(\omega,\mu)$ is smooth and nonvanishing. 
    Let $V:=\ker (W(\omega_0)-\Id)$. Then, decomposing $\mathbb{C}^{2J}=V\oplus V^\perp$, we have
    \begin{equation}
        W(\omega_0)-\Id=\begin{pmatrix} 0&0\\0&W_{22}(\omega_0)-\Id\end{pmatrix},
    \end{equation}
    where $W_{22}(\omega_0)-\Id$ is invertible. Then, with $z=\omega+i\mu$, we write
    \begin{equation}
        \tilde{W}(z,\bar{z})-\Id=\begin{pmatrix} \tilde{W}_{11}(z,\bar{z})-\Id&\tilde{W}_{12}(z,\bar{z})\\W_{21}(z,\bar{z})&\tilde{W}_{22}(z,\bar{z})-\Id\end{pmatrix},
    \end{equation}
    with $\tilde{W}_{ij}(z,\bar{z})\in C^2$. Notice that 
    \begin{equation}
        -i\tilde{W}^*\partial_\omega \tilde{W}|_{\omega_0}=A_W(\omega_0)<0,\ -i\tilde{W}^*\partial_\mu \tilde{W}|_{\omega_0}=iA_W(\omega_0).
    \end{equation}
    Hence, since $W_{11}(\omega_0)=\Id$, 
    \begin{equation}
        \tilde{W}_{11}(z,\bar{z})= \Id+ i(A_W(\omega_0))_{11}(z-\omega_0)+O(|z-\omega_0|^2)
    \end{equation}
    with $(A_W(\omega_0))_{11}<0$. Thus, 
    \begin{equation}
        \tilde{W}(z,\bar{z})-\Id=\begin{pmatrix} i(z-\omega_0)(A_W(\omega_0))_{11}+O(|z-\omega_0|^2)&O(|z-\omega_0|)\\O(|z-\omega_0|)&\tilde{W}_{22}(z,\bar{z})-\Id\end{pmatrix}.
    \end{equation}
    The Schur complement formula then implies
    \begin{multline}
    \det (\tilde{W}(z,\bar{z})-\Id)=\det (\tilde{W}_{22}(z,\bar{z})-\Id)
\det \Big(i(z-\omega_0)(A_{W}(\omega_0))_{11}+O(|z-\omega_0|^2)\\
-O(|z-\omega_0|)(\tilde{W}_{22}(z,\bar{z})-\Id)^{-1}O(|z-\omega_0|)\Big)=e_1(z,\bar{z})(z-\omega_0)^m\det ((A_{W}(\omega)_0)_{11}+O(|z-\omega_0|))
    \end{multline}
    from which the lemma follows.
\end{proof}

For $\e>0$, we set up the target spectral scale
\begin{equation}\label{e:spectralscaler}
    \gamma_{\mc{P},M}(\omega):= \omega^{-\min(\frac{2\pi}{\alpha_+},1+k_0)+\e} +M^2\omega^{-2+\e},
\end{equation}
abbreviated $\gamma(\omega)$ when there is no confusion. 

\begin{lemma}[Counting zeros finely]
\label{l:nonUniformEps0}
    Let $r>0$, $\delta>0$, $\e>0$. There is $\omega_0>0$ (depending on $r,\e,\delta$) such that for any $\omega\geq \omega_0$ and any interval $[a,b]\subset [\frac12\omega,\infty)$ such that
    \begin{equation}
    \label{e:endpointDistance}
    \dist\big(\{a,b\},\mc{Z}^{\mathscr{D}}\big)>r\gamma(\omega)=r\big(\omega^{-\min(\frac{2\pi}{\alpha_+},1+k_0)+\e} +M^2\omega^{-2+\e}\big)
    \end{equation}
    we have
    \begin{equation}
        \mc{N}_{[a,b]}^{\mathscr{D}}=\sum_{\omega\in[a,b]}\dim \ker \mc{Q}_{\omega}^{\mathscr{D}}=\sum_{\omega\in[a,b]}\dim \ker \mc{Q}_\omega^\Lambda=\mc{N}_{[a,b]}^{\Lambda},
    \end{equation}
    where $\mc{Q}_z^\Lambda=E_{-+}(z)$ for $z\in\mathbb{C}$ from Proposition~\ref{p:grushin}.
\end{lemma}
\begin{proof}
Recall that for $\omega_0>0$ large enough $\mc{Q}_{z}^{\Lambda}$ is analytic in $z$ in $[\omega_0,\infty)+i[-r_1,r_1]$ for some $r_1>0$. Moreover, by Corollary~\ref{c:derivatives}, for $\beta\ge0$, 
\begin{equation}
\label{e:exactDerivative}
\|\partial_{\omega}^\beta (\mc{Q}^{\Lambda}_\omega+T_\omega)(\omega)\|\leq C|\omega|^{-2+0}.
\end{equation}
Recall also the smooth extension~\eqref{e:Qextend} of $\mc{Q}_{\omega}^{\mathscr{D}}$.

Observe that 
\begin{equation}
\tilde{\mc{Q}}_{\omega+i\mu}^{\mathscr{D}}=\Sigma^{k_0}(\omega)(\Id-W(\omega)(\Id-\mu A_W(\omega)))+i\mu \partial_\omega\Sigma^{k_0}(\omega)(\Id-W(\omega)).
\end{equation}
Note that the second term is of size $o(|\mu|)$ as $\omega\to \infty$. Since $A_W(\mu)\le-\frac{L_-}2\Id$, we have either
\begin{equation}
    \Id-\mu A_W(\omega)\ge (1+\frac{L_-\mu}2)\Id\text{ for }\mu>0, \text{ or }\Id-\mu A_W(\omega)\le (1-\frac{L_-|\mu|}2)\Id\text{ for }\mu<0.
\end{equation}
Since $W(\omega)$ is unitary, this implies for large enough $\omega_0$,
\begin{equation}
\label{e:extendedEstimate}
    \|(\tilde{\mc{Q}}_{\omega+i\mu}^{\mathscr{D}})^{-1}\|\le \frac{4}{L_-|\mu|}  ,\quad \omega+i\mu \in [\omega_0,\infty)+i[r_1,r_1]
\end{equation}

\begin{figure}[t]
  \centering
        \begin{tikzpicture}[x=1cm,y=1cm]
  \usetikzlibrary{calc,decorations.markings,plotmarks}
  \definecolor{colorAKANE}{RGB}{183,40,46}
  \definecolor{colorUSUAI}{RGB}{40,148,200}
  \colorlet{plotPrimary}{colorAKANE}
  \colorlet{plotTertiary}{colorUSUAI}
  \colorlet{plotAxis}{black!72}
  \colorlet{plotGuide}{black!34}
  \colorlet{plotField}{black!5}
  \tikzset{
    >=latex,
    font=\small,
    mark size=1.5pt,
    mark options={thin},
    line cap=round,
    line join=round,
    figure label/.style={font=\fontsize{8}{9.6}\selectfont,
      inner sep=1pt,text=black},
    small label/.style={figure label,font=\fontsize{7}{8.4}\selectfont},
    coordinate axis/.style={draw=plotAxis,line width=0.4pt,->},
    height guide/.style={draw=plotGuide,line width=0.3pt,
      dash pattern=on 1.2pt off 1.6pt},
    contour/.style={draw=plotTertiary,line width=0.85pt,
      postaction={decorate},decoration={markings,
        mark=at position 0.18 with {\arrow{latex}},
        mark=at position 0.435 with {\arrow{latex}},
        mark=at position 0.68 with {\arrow{latex}},
        mark=at position 0.935 with {\arrow{latex}}}},
    distance arrow/.style={draw=black!65,line width=0.35pt,<->},
    dirac zero/.style={fill=black,draw=black,line width=0.3pt},
    exact zero/.style={draw=plotPrimary,line width=0.7pt,
      mark=x,mark size=1.8pt,
      mark options={draw=plotPrimary,line width=0.7pt}}
  }

  \path[use as bounding box] (0,0) rectangle (7.5,3.2);
  \begin{pgfinterruptboundingbox}
  \begin{scope}[shift={(-0.06,-0.24)}]
    \coordinate (real left) at (2.30,1.60);
    \coordinate (real right) at (6.00,1.60);
    \coordinate (lower left) at ($(real left)+(0,-0.66)$);
    \coordinate (lower right) at ($(real right)+(0,-0.66)$);
    \coordinate (upper left) at ($(real left)+(0,0.66)$);
    \coordinate (upper right) at ($(real right)+(0,0.66)$);
    \coordinate (left gap left) at ($(real left)+(-0.66,0)$);
    \coordinate (left gap right) at ($(real left)+(0.66,0)$);
    \coordinate (right gap left) at ($(real right)+(-0.66,0)$);
    \coordinate (right gap right) at ($(real right)+(0.66,0)$);

    \path[fill=plotTertiary!4,draw=none]
      (lower left) rectangle (upper right);
    \path[fill=black!12,draw=none]
      ($(left gap left)+(0,-0.07)$)
        rectangle ($(left gap right)+(0,0.07)$)
      ($(right gap left)+(0,-0.07)$)
        rectangle ($(right gap right)+(0,0.07)$);

    \draw[coordinate axis] (1.10,1.60) -- (7.20,1.60);
    \draw[coordinate axis] (1.25,0.70) -- (1.25,2.67);
    \node[figure label,anchor=west] at (7.24,1.60) {$\omega$};
    \node[figure label,anchor=south] at (1.25,2.69) {$\mu$};
    \draw[height guide] (1.25,2.26) -- (upper left);
    \draw[height guide] (1.25,0.94) -- (lower left);
    \node[figure label,anchor=east] at (1.16,2.26) {$r\gamma(\lambda)$};
    \node[figure label,anchor=east] at (1.16,0.94) {$-r\gamma(\lambda)$};

    \draw[contour]
      (lower left) -- (lower right) -- (upper right) -- (upper left) -- cycle;
    \node[figure label,text=plotTertiary,anchor=south]
      at (4.45,2.32) {$\partial S$};
    \node[figure label] at (4.15,1.99) {$S$};
    \node[figure label,anchor=north]
      at ($(real left)+(0,-0.12)$) {$a$};
    \node[figure label,anchor=north]
      at ($(real right)+(0,-0.12)$) {$b$};

    \draw[height guide]
      (left gap right) -- ($(left gap right)+(0,0.24)$);
    \draw[height guide]
      (right gap left) -- ($(right gap left)+(0,0.24)$);
    \draw[distance arrow]
      ($(real left)+(0,0.24)$) -- ($(left gap right)+(0,0.24)$)
      node[small label,midway,above=2pt] {$r\gamma(\lambda)$};
    \draw[distance arrow]
      ($(right gap left)+(0,0.24)$) -- ($(real right)+(0,0.24)$)
      node[small label,midway,above=2pt] {$r\gamma(\lambda)$};

    \path[dirac zero]
      ($(left gap left)+(-0.20,0)$) circle[radius=1.35pt]
      (3.23,1.60) circle[radius=1.35pt]
      (4.12,1.60) circle[radius=1.35pt]
      (5.08,1.60) circle[radius=1.35pt]
      ($(right gap right)+(0.22,0)$) circle[radius=1.35pt];
    \draw[exact zero] plot[only marks] coordinates {
      (3.39,1.60) (4.29,1.60) (4.90,1.60)
    };
  \end{scope}
  \end{pgfinterruptboundingbox}
\end{tikzpicture}%
  \caption{Contour integral used in Lemma~\ref{l:nonUniformEps0}. Zeros of $\mc{Q}^{\mathscr{D}}_{\omega}$ and $\mc{Q}^{\Lambda}_{\omega}$ are marked as black dots and red crosses respectively.}
  \label{f:10l1}
\end{figure}
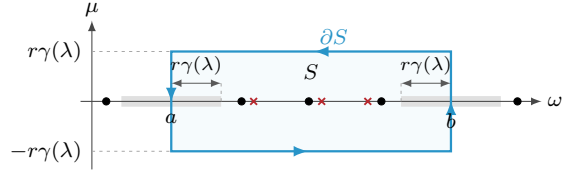

Let the box $S:= [a,b]+i[-r,r]\gamma(\omega)$, where $\gamma(\omega)$ was defined in \eqref{e:spectralscaler}. See Figure~\ref{f:10l1}. Let $z=\omega+i\mu\in \overline{S}$. By~\eqref{e:distanceInverseEstimate}, \eqref{e:endpointDistance}, \eqref{e:extendedEstimate}, and the fact that $\|\nabla \tilde{\mathcal{Q}}_{\omega+i\mu}^{\mathscr{D}}\|\leq C$ we have
\begin{equation}
    \|(\tilde{\mc{Q}}_{z}^{\mathscr{D}})^{-1}|_{\partial S}\|\le Cr^{-1}\gamma(\omega)^{-1}
\end{equation}
In addition, for $\omega \in [a,b]$ we have
\begin{equation}
    \|\mc{Q}_{\omega}^{\Lambda}-\mc{Q}_{\omega}^{\mathscr{D}}\|\leq C\gamma(\omega)\omega^{-\frac{\e}{2}},
\end{equation}
and by~\eqref{e:diracDerivative} and~\eqref{e:exactDerivative} we have
\begin{equation}
    \|\partial_\omega (\mc{Q}_{\omega}^{\Lambda}-\mc{Q}_{\omega}^{\mathscr{D}})\|\leq C\omega^{-\frac{1}{2}+0}.
\end{equation}
Hence, for $z=\omega+i\mu\in \overline{S}$, we have
\begin{equation}
    \|\mc{Q}_{\omega+i\mu}^{\Lambda}-\tilde{\mc{Q}}_{\omega+i\mu}^{\mathscr{D}}\|\leq C \gamma(\omega)\omega^{-\frac{\e}{2}}+C\omega^{-\frac{1}{2}+\e}|\mu|.
\end{equation}
Therefore, for $z\in\partial S$,
\begin{multline}
    \det \mc{Q}_z^{\Lambda}=\det \tilde{\mc{Q}}_z^{\mathscr{D}}\det (I+(\tilde{\mc{Q}}_z^{\mathscr{D}})^{-1}(\mc{Q}_{z}^{\Lambda}-\tilde{\mc{Q}}_{z}^{\mathscr{D}}))
    =\det \tilde{\mc{Q}}_z^{\mathscr{D}}\det (I+(\tilde{\mc{Q}}_z^{\mathscr{D}})^{-1}O(\gamma(\omega)\omega^{-\frac{\e}{2}}))\\
    =\det \tilde{\mc{Q}}_z^{\mathscr{D}}\det (I+O(\omega^{-\frac{\e}{2}}))=\det \tilde{\mc{Q}}_z^{\mathscr{D}} (1+O(\omega^{-\frac{\e}{2}})).
\end{multline}
Now, by homotopy invariance of the winding number~\cite[p. 134 property (c)]{Hatcher2002}, this implies 
\begin{equation}
\operatorname{wind}_{\partial S}\det\mc{Q}_z^{\Lambda}=\operatorname{wind}_{\partial S}\det\tilde{\mc{Q}}_z^{\mathscr{D}}=\sum_{\omega\in [a,b]}\dim \ker \mc{Q}^{\mathscr{D}}_\omega=\mc{N}^{\mathscr{D}}_{[a,b]},
\end{equation}
where $\operatorname{wind}_{\partial S}$ denotes the winding number along $\partial S$ and the second equality follows from Lemma~\ref{l:identify} and \eqref{e:extendedEstimate} implies that all zeros of $\det\tilde{\mc{Q}}_z^{\mathscr{D}}$ are real. 

Finally, the fact that $\|(\mc{Q}^{\Lambda}_z)^{-1}\|\leq C|\Im z|^{-1}$ implies that $\mc{Q}^{\Lambda}_z$ has no Jordan blocks with eigenvalue 0 and that all zeros of $\det \mc{Q}_z^\Lambda$ are real and hence that  
\begin{equation}
    \operatorname{wind}_{\partial S}\det\mc{Q}_z^{\Lambda}=\#\{\text{zeros of }\det\mc{Q}_z^{\Lambda}\text{ in }\Omega\}=\sum_{\omega\in [a,b]}\dim \ker \mc{Q}^{\Lambda}_\omega=\mc{N}^{\Lambda}_{[a,b]},
\end{equation}
since $\det \mc{Q}_z^{\Lambda}$ is analytic in $S$. This completes the proof.
\end{proof}

An almost identical proof yields the following lemma.
\begin{lemma}[Counting zeros loosely but uniformly in $\delta$]\label{l:uniformEps0} 
    Let $r>0$, $\e>0$. There is $\omega_0>0$ (depending on $r,\e$) such that for any $0<\delta<1$, any $\omega\geq \omega_0$ and any interval $[a,b]\subset [\frac12 \omega,\infty)$ such that 
    \begin{equation}
    \dist\big(\{a,b\},\mc{Z}^{\mathscr{D}}\big)>r\omega^{-1+\e}
    \end{equation}
    we have
    \begin{equation}
        \mc{N}_{[a,b]}^{\mathscr{D}}=\sum_{\omega\in[a,b]}\dim \ker \mc{Q}_{\omega}^{\mathscr{D}}=\sum_{\omega\in[a,b]}\dim \ker \mc{Q}_\omega^\Lambda=\mc{N}_{[a,b]}^{\Lambda},
    \end{equation}
\end{lemma}
We emphasize that the key to uniformity in $\delta$ is that $(\Sigma^{k_0}(\omega)-\Sigma^0)=O(\omega^{-1+\e})$ while omitting the proof. 

\subsection{Proof of Theorem \ref{t:boundedSteklov}}

\begin{proof}[Proof of Theorem~\ref{t:boundedSteklov}]
1. We remind the reader that our goal is to show first that, for each $\e>0$, there is $\omega_0>0$ such that eigenvalues in $\operatorname{Spec}(\Lambda^{\mathcal{P}})\cap [\omega_0,\infty)$ and $\operatorname{Spec}(\mathscr{D}^{k_0})\cap [\omega_0,\infty)$ are related pairwise, and each pair is separated by small distance. Specifically, in Steps 2--5, we show pairing happens at the scale
\begin{equation}
    \gamma_{\mc{P},M}(\omega)= \omega^{-\min(\frac{2\pi}{\alpha_+},1+k_0)+\e} +M^2\omega^{-2+\e},
\end{equation}
defined in \eqref{e:spectralscaler}, and then in Step 6 we show that the pairs of labeling of each pair agrees. To do the pairing, we will show $[\omega_0, \infty)$ can be written as a union of intervals of length $\sim \gamma(\omega)$, so that the numbers of eigenvalues must agree in each such interval. 

2. We show for each $r\le \frac14$, there is $\omega_0>0$ such that for any $\omega\ge \omega_0$ with $\det\mc{Q}_{\omega}^{\Lambda}=0$, we have
\begin{equation}\label{e:7l43}
    \mc{Z}^{\mathscr{D}}\cap[\omega-4r\gamma(\omega), \omega+4r\gamma(\omega)]=\{\omega'\in [\omega-4r\gamma(\omega), \omega+4r\gamma(\omega)]: \det\mc{Q}_{\omega'}^{\mathscr{D}}=0\}\neq \emptyset.
\end{equation}
Assume this is not true, then $\dist(\{\omega\pm 2r\gamma(\omega)\},\mc{Z}^{\mathscr{D}})\ge 2r\gamma(\omega)$, and the fine counting Lemma~\ref{l:nonUniformEps0} implies
\begin{equation}
    \mc{N}^{\mathscr{D}}_{[\omega-2r\gamma(\omega), \omega+2r\gamma(\omega)]}=\mc{N}^{\Lambda}_{[\omega-2r\gamma(\omega), \omega+2r\gamma(\omega)]}\ge 1,
\end{equation}
yielding a contradiction. This argument specifically implies that if \eqref{e:7l43} is empty, then 
\begin{equation}
    \{\omega\in [\omega-r\gamma(\omega), \omega+r\gamma(\omega)]: \det\mc{Q}_{\omega}^{\Lambda}=0\}=\emptyset.
\end{equation}

3. By Lemma~\ref{l:rootnumber} there are at most $2J$ zeros of $\det \mc{Q}^{\mathscr{D}}_\omega$ within intervals of length $2r_0$. Fix $r=\frac1{64(8J+1)}$, there is $\omega_0$ that allows the use of Step 2. Without loss of generality, let $\gamma(\omega_0)\le \frac{1}{4r_0}$. We claim for each $\omega>\omega_0$, there exists two intervals $I_\pm$ with center $a_\pm$ that
\begin{equation}
    I_\pm =\bar{B}(a_\pm,\frac{\gamma(\omega)}{32(8J+1)})\cap\mathbb{R},\ I_-\subset [\omega-\gamma(\omega), \omega-\frac{3}{4}\gamma(\omega)], \ I_+ \subset [\omega+\frac{3}{4}\gamma(\omega), \omega+\gamma(\omega)], 
\end{equation}
such that $I_\pm\cap (\mc{Z}^{\mathscr{D}}\cup\mc{Z}^{\Lambda})=\emptyset$. Indeed, the $2J$ upper bound on the zeros of $\det \mc{Q}^{\mathscr{D}}_\omega$ within
\begin{equation}
    [\omega-\gamma(\omega), \omega-\frac{3}{4}\gamma(\omega)]=\cup_{j=0}^{8J}I_j, \ I_j=[\omega-(1-\frac{j}{4(8J+1)})\gamma(\omega), \omega-(1-\frac{j+1}{4(8J+1)})\gamma(\omega)].
\end{equation}
Thus there must be one of those intervals centered at some $a_-$ that contains no zeros of $\det \mc{Q}^{\mathscr{D}}_\omega$, and Step 2 implies $I_-=\bar{B}(a_\pm,\frac{1}{32(8J+1)}\gamma(\omega))\cap\mathbb{R}$ contains neither zeros of $\det \mc{Q}^{\mathscr{D}}_\omega$ nor zeros of $\det \mc{Q}^{\Lambda}_\omega$. Similar argument gives $I_+$. 

4. We show $\det \mc{Q}^{\mathscr{D}}_\omega$ and $\det \mc{Q}^{\Lambda}_\omega$ have the same number of zeros in $\omega\in[a_-, a_+]$, where
\begin{equation}
    \omega-\gamma(\omega)\le a_-\le \omega-\frac34\gamma(\omega)\le \omega+\frac34\gamma(\omega)\le a_+\le \omega+\gamma(\omega).
\end{equation}
Note that $I_\pm\cap \mc{Z}^{\mathscr{D}}=\emptyset$ implies $\dist\big(\{a_-,a_+\},\mc{Z}^{\mathscr{D}}\big)>\frac{\gamma(\omega)}{64(8J+1)}$. Apply the fine counting Lemma~\ref{l:nonUniformEps0} to see
\begin{equation}
    \mc{N}^{\mathscr{D}}_{[a_-, a_+]}=\mc{N}^{\Lambda}_{[a_-, a_+]}. 
\end{equation}

5. Now we show $\mc{Z}^{\mathscr{D}}$ and $\mc{Z}^{\Lambda}$ have one-to-one correspondence after large $\omega\ge \omega_0$. Start with $\omega_0$, 
where $\mc{N}^{\mathscr{D}}_{[a_-, a_+]}=\mc{N}^{\Lambda}_{[a_-, a_+]}$. Let $\omega_1>\omega_0$ such that
\begin{equation}
    \omega_0+\frac34\gamma(\omega_0)\le \omega_1-(1-\frac{1}{32(8J+1)})\gamma(\omega_0)\le 
    a_+\le
    \omega_1-(\frac{3}{4}+\frac{1}{32(8J+1)})\gamma(\omega_0).
\end{equation}
Note $\dist\big(a_+,\mc{Z}^{\mathscr{D}}\big)>\frac{\gamma(\omega_1)}{32(8J+1)}$ satisfies the requirement in Step 3. Thus we can repeat Step 3 and get
some $a_+'\in [\omega_1+\frac{3}{4}\gamma(\omega_1), \omega_1+\gamma(\omega_1)]$ where we still have
\begin{equation}
    \mc{N}^{\mathscr{D}}_{[a_+, a_+']}=\mc{N}^{\Lambda}_{[a_+, a_+']}.
\end{equation}
We can iterate this process to cover $\mathbb{R}_{\ge \omega_0}$. 

6. We now show that the labels of the eigenvalues for the Dirac operator and the Steklov problem agree. For a curvilinear polygon $\mathcal{P}\subset \mathbb{R}^2$ and $M\ge 0$, we denote the Steklov eigenvalues by
\begin{equation}
    0=\sigma_0(\mathcal{P},M)\le \sigma_1(\mathcal{P},M)\leq \sigma_2(\mathcal{P},M)\le \cdots,
\end{equation}
and the eigenvalues of $\mathscr{D}$ by
\begin{equation}
    \cdots\lambda_{-2}(\mathcal{P})\le \lambda_{-1}(\mathcal{P})\le0\le\lambda_0(\mathcal{P})\leq \lambda_1(\mathcal{P})\le \lambda_2(\mathcal{P})\leq \cdots. 
\end{equation}
By the symmetry of $\mc{Q}^{\mathscr{D}}_\omega$, we may choose $\lambda_n$ so that 
$\lambda_{n}(\mathcal{P})=-\lambda_{-n-1}(\mathcal{P})$ for $n\in\mathbb{N}_0$. Let $M_0>0$, $J\ge 1$, $\{C_j\}_{j=0}^\infty$ be positive constants. 

Using the argument above in the previous steps, but appealing to Lemma~\ref{l:uniformEps0} instead of Lemma~\ref{l:nonUniformEps0}, there is $\omega_0>0$ such that for all $0<\delta<1$, all curvilinear polygons $\mathcal{P}$ with $J$ vertices with angles $\alpha_j\in(\e,2\pi -\e)$ for all $1\le j\le J$, and with edges
\begin{equation}
    \e \leq |\Gamma_j|\leq \e^{-1}, \ \|\Gamma_j\|_{C^k}\leq C_k,\ k\in \mathbb{N}_0,
\end{equation}
and all $M\in[0,M_0]$ there is $N=N(\mathcal{P},M)$ such that for all $n$ with $\sigma_n(\mathcal{P},M)>\omega_0$,
    \begin{equation}
    \label{e:matching}
    |\lambda_{n+N}(\mathcal{P})-\sigma_n(\mathcal{P},M)|\leq C\sigma_n(\mc{P},M)^{-1+\e}.
    \end{equation}
   
    Now, let $\mathscr{D}^0$ denote the basic Steklov Dirac operator on $\mc{P}$ and let $\delta<\frac{1}{2}\omega_0^{-1}$. Then, if $n\geq 0$,
    $
    \lambda_n(\mathscr{D}^0)\leq 2\omega_0
    $
    we have $\lambda_n(\mathscr{D}^0)=\lambda_n(\mathscr{D})$. In particular, the value of $N$ in~\eqref{e:matching} does not depend on $k_0$. Thus, from now on, we argue with $k_0=0$. 

6a. We show that $N(\mc{P},M)$ in~\eqref{e:matching} is unique. Indeed, suppose that~\eqref{e:matching} holds with both $N$ and $N'$ instead of $N$ with $N\neq N'$ for $n=n_0,n_0+1,\dots n_0+2J$. Suppose without loss of generality that $N'>N$. Then for all $0\leq k\leq 2J+1$ and $n_0$ large enough, 
    \begin{multline}  
    |\lambda_{n_0+k(N'-N)}(\mathcal{P})-\lambda_{n_0}(\mathcal{P})|
    \leq \sum_{m=0}^{k-1} |\lambda_{n_0+m(N'-N)}(\mathcal{P})-\sigma_{n_0+m(N'-N)-N}(\mathcal{P})|\\+\sum_{m=0}^{k-1} |\sigma_{n_0+m(N'-N)-N}(\mathcal{P})-\lambda_{n_0+(m+1)(N'-N)}
(\mathcal{P})|\leq Ck\sigma_{n_0-N}^{-1+\e}.
    \end{multline}
    Let $k=2J+1$. In particular since $\sigma_j\to \infty$, this implies that for any $\e>0$ there is large $C>0$ such that 
    \begin{equation}
        \#\{\lambda_j\in [C,C+\e]\}\geq k=2J+1,
    \end{equation}
    which contradicts Lemma~\ref{l:rootnumber}. In particular, for any $n_0$ large enough, there is a unique $N$ such that~\eqref{e:matching} holds for $n\in\{n_0,n_0+1,\dots, n_0+2J\}$. 

    6b. We next claim that if $\partial\mathcal{P}$ is smooth, $M=0$, and $\partial\mathcal{P}$ is flat near the (fictitious) corners, then $N(\mc{P},0)=0$. Let $\partial\mathcal{P}$ be smooth, flat near (fictitious) corners with connected components of lengths $\{L^1,\dots, L^L\}$. One can check that the Dirac operator orthogonally decomposes into $\mathscr{D}=\oplus_{l=1}^L\mathscr{D}^l$ where $\mathscr{D}^l$ is the Dirac operator on the $l$-th connected component of the boundary. Then since
    \begin{equation}
        \Sigma_{0,1}^{W_\pi}=\begin{pmatrix} 0&1\\1&0\end{pmatrix}, \ \lambda_{2n}(\mathscr{D}^l)=\lambda_{2n-1}(\mathscr{D}^l)=\frac{2\pi n}{L^l}, \ n\geq 0.
    \end{equation}
   Then, \cite[Theorem 1.4]{GPPS2014} (see also Rozenblum~\cite{Rozenblyum1986} for the simply connected version) proves the claim.

    6c. To prove that $N(\mc{P}, M)=0$ for all $M$ and curvilinear polygons $\mathcal{P}$, we use a deformation argument. Fix $M>0$ and $\mathcal{P}$ a curvilinear polygon with $J$ vertices with angles $\alpha_j\in(2\e,2\pi -2\e)$ for all $1\le j\le J$, and with edges
\begin{equation}
    2\e \leq |\Gamma_j|\leq \frac12\e^{-1}, \ \|\Gamma_j\|_{C^k}\leq \frac12C_k,\ k\in \mathbb{N}_0.
\end{equation}
Then there is a continuous family of curvilinear polygons $t\mapsto \mathcal{P}_t$ such that $\mathcal{P}_0=\mc{P}$ and $\mathcal{P}_1$ has smooth boundary and is flat near the corners $x_j$, such that $\mc{P}_t$ has the same set of edges and vertices, with angles 
\begin{equation}
    \e<\inf_{t\in[0,1]}(\alpha_j(\mathcal{P}_t))\le\sup_{t\in[0,1]}(\alpha_j(\mathcal{P}_t))<2\pi -\e,
\end{equation}
    and 
\begin{equation}
    \e \leq \inf_{t\in[0,1]}|\Gamma_j(\mathcal{P}_t)|\leq\sup_{t\in[0,1]}|\Gamma_j(\mathcal{P}_t)|\leq  \e^{-1}, \ \sup_{t\in[0,1]}\|\Gamma_j(\mathcal{P}_t)\|_{C^k}\leq C_k,
\end{equation}
for each $1\le j\le J$. In particular, there is $\lambda_0>0$ such that for all $t,t'\in[0,1]$, there is a unique $N_{t,t'}\in\mathbb{Z}$ such that for $\sigma_n(\mathcal{P}_t,(1-t')M)>\lambda_0$, 
    \begin{equation}
    \label{e:matching2}
    |\lambda_{n+N_{t,t'}}(\mathcal{P}_t)-\sigma_n(\mathcal{P}_t,(1-t')M)|<Cr_{\mc{P}_t,(1-t')M}(\sigma_n(\mathcal{P}_t,(1-t')M)).
    \end{equation}
Now, observe that since $\Sigma^{\mc{W}_{\alpha_t}}_{0,1}$ is continuous in $t$, $t\mapsto \lambda_n(\mathcal{P}_t)$ is continuous. Moreover, the variational characterization of $\sigma_n(\mathcal{P}_t,(1-t')M)$ implies that $t'\mapsto \sigma_n(\mathcal{P}_t,(1-t')M)$ is continuous and, by~\cite[Theorem 4.1]{BuGiTr:20}, $t\mapsto \sigma_n(\mathcal{P}_t,0)$ is also continuous. 

6d. We show that $N_{0,1}=0$. Since $\mathcal{P}_{1}$ is smooth, we have $N_{1,1}=0$. Let 
    \begin{equation}
        t_0=\inf\{t\in [0,1]\,:\, N_{s,1}=0,\,s\in [t,1]\}.
    \end{equation}
    If $t_0=0$, then by continuity of $t\mapsto \lambda_n(\mathcal{P}_t)$ and $t\mapsto\sigma_n(\mathcal{P}_t,0)$, $N_{0,1}=0$. Assume instead $t_0>0$. Then, by construction there are $s_k\rightarrow t_0-0$ and $t_k\rightarrow t_0+0$ such that 
    \begin{equation}
        N_{s_k,1}\neq 0,\ N_{t,1}=0,\ t\in[t_k,1].
    \end{equation}
    Hence, using continuity of $t\mapsto \lambda_n(\mathcal{P}_t)$ and $t\mapsto\sigma_n(\mathcal{P}_t,0)$ again, together with the fact that $N_{t_0,1}$ is determined by finitely many values of $\lambda_n$ and $\sigma_n$, $N_{t_0,1}=\inf_{t>t_0}N_{t,1}=0$ and $N_{t_0, 1}=\lim_{k\to \infty}N_{s_k,1}\neq 0$,
    yielding the desired contradiction. Note that the limit $\lim_{k\to \infty}N_{s_k,1}$ exists since $N_{s_k,1}$ is integer valued.

   6e. Now, we conclude with showing that $N_{0,0}=0$. To do this, we use that $N_{0,1}=0$. Let 
    \begin{equation}
        t'_0=\inf\{t'\in [0,1]\,:\, N_{0,s}=0,\,s\in [t',1]\}. 
    \end{equation}
    If $t'_0=0$, then by continuity of $s\mapsto\sigma_n(\mathcal{P}_0, (1-s)M)$, $N_{0,0}=0$. Assume instead $t'_0> 0$. Then, by construction there are $s_k\rightarrow t'_0-0$ and $t'_k\rightarrow t'_0+0$ such that 
    \begin{equation}
        N_{0,s_k}\neq 0,\ N_{0,t'}=0,\ t'\in[t'_k,1].
    \end{equation}
    Hence, using continuity of $s\mapsto\sigma_n(\mathcal{P}_0, (1-s)M)$, $N_{0,t'_0}=\inf_{t'>t'_0}N_{0,t'}=0$ and $N_{0,t'_0}=\lim_{k\to \infty}N_{0,s_k}\neq 0$,
    yielding the desired contradiction. 
\end{proof}

\part{Appendices}
\appendix

\section{A second-microlocal pseudodifferential calculus}
\label{a:1}
Given smooth scattering manifold $\partial\Omega$, consider its half-space collar $\Omega_+=\{(x_0, x): x_0\in [0,\infty), x\in \partial\Omega\}$. We denote $\bar x=(x_0,x)$. We define $(T^*\partial\Omega)_+=\{(x_0, x,\xi): x_0\in [0,\infty), (x,\xi)\in T^*\partial\Omega\}$ similarly. For $\delta\in[0,1)$ we consider the following symbol class:
\begin{definition}[Symbols in $S^{s,l}_{k,\delta}((T^*\partial\Omega)_+)$] \label{def:Symbols}
Let $M_0>0$. Let $h_0>0$ small. We say $a(x_0,x,\xi;M, h)\in C^\infty([0,\infty)_{x_0}\times \partial\Omega_x\times \mathbb{R}^{d-1}_\xi\times[0, M_0]_M\times (0, h_0]_h)$ is in $S^{s,l}_{k,\delta}((T^*\partial\Omega)_+)$, if for every multi-indices $\alpha,\beta$ we have
\begin{equation}
\abs{(\langle x\rangle\partial_{\bar x})^\alpha (\langle x\rangle^{-\delta} \langle \langle x\rangle^{\delta}\sqrt{|\xi|^2+M^2}\rangle\partial_\xi)^\beta a(x_0,x,\xi; M, h)}\le C_{\alpha\beta}\langle \xi\rangle^{s}\langle x\rangle^{l}\langle \langle x\rangle^{\delta}\sqrt{|\xi|^2+M^2}\rangle^{k},
\end{equation}
uniformly for all $x_0\in[0,\infty)$, $x\in \partial\Omega$, $\xi\in \mathbb{R}^{d-1}$, $h\in(0, h_0]$ and $M\in [0,M_0]$.
\end{definition}

For brevity, we may omit the dependency on $M, h$ whenever there is no confusion. Our goal in this section is to build a pseudodifferential calculus on $\Omega_+$, associated with the symbol class $S^{s,l}_{l,\delta}((T^*\partial\Omega)_+)$. To do it we need to look at the half-space case $\Omega_+=\mathbb{R}^{d}_+=\{(x_0,x): x_0\in[0,\infty), x\in \mathbb{R}^{d-1}\}$ with $\partial\Omega=\partial \mathbb{R}^{d}_+=\mathbb{R}^{d-1}$, $(T^*\partial\Omega)_+=\mathbb{R}_{x_0}\times T^*\mathbb{R}^{d-1}$, and then later transplant them to a general $\Omega_+$.  
\begin{definition}[A second-microlocal calculus on $\mathbb{R}^{d}_+$]
Let $a(x_0, x,\xi)\in S^{s,l}_{k,\delta}((T^*\mathbb{R}^{d-1})_+)$. We define its left quantization as an operator bounded from $\mc{S}(\mathbb{R}^d_+)$ to $\mc{S}'(\mathbb{R}^d_+)$ whose kernel is given by
\begin{equation}
    \oph(a)(x_0,x; y_0, y)=\frac{\delta(y_0-x_0)}{(2\pi h)^{d-1}}\int_{\mathbb{R}^{d-1}_{\xi}}e^{\frac{i}{h}\langle x-y,\xi\rangle}a(x_0,x,\xi)\,d\xi.
\end{equation}
We define the residual class by $h^\infty\Psi^{-\infty,-\infty}_{-\infty,\delta}(\mathbb{R}^d_+)=h^\infty\Psi^{-\infty, -\infty}(\mathbb{R}^d_+)$, consisting of all operators with kernel in $h^\infty \delta(x_0-y_0) C^\infty_b([0,\infty)_{x_0}; \mc{S}(\mathbb{R}^{d-1}_x\times\mathbb{R}^{d-1}_y)$. We define $\Psi^{s,l}_{k,\delta}(\mathbb{R}^d_+)=\oph(S^{s,l}_{k,\delta}((T^*\mathbb{R}^{d-1})_+))$. 
\end{definition}

\subsection{Change of coordinates}
Now we study the invariance of the calculus under change of coordinates.
We start from proving the invariance of $S^{s,l}_{k,\delta}((T^*\mathbb{R}^{d-1})_+)$.
\begin{lemma}[Change of variables for symbols]
    Let $\Phi:\R^{d-1}\to \R^{d-1}$ be a diffeomorphism that 
    \begin{equation}
        |\partial^\alpha \Phi(x)|\leq C_\alpha\langle x\rangle^{1-|\alpha|}, \ |\partial^\alpha (\Phi^{-1}(x'))|\leq C_\alpha\langle x'\rangle^{1-|\alpha|}
    \end{equation}
    for every multi-index $\alpha$. Then, for any $a(x_0, x,\xi)\in S^{s,l}_{k,\delta}((T^*\mathbb{R}^{d-1})_+)$, we have
    \begin{equation}
        (\Phi^{*}a)(x_0,x',\xi')=a(x_0,\Phi^{-1}(x'),D\Phi|_{\Phi^{-1}(x')}^T\xi')\in S^{s,l}_{k,\delta}((T^*\mathbb{R}^{d-1})_+).
    \end{equation}
\end{lemma}
\begin{proof}
    We follow the approach of \cite[Theorem 9.4]{Zw:12}. Let $x=\Phi^{-1} (x')$ with $\xi=D\Phi|_{\Phi^{-1}(x')}^T\xi'$ given by the symplectomorphism induced by $\Phi$. Compute
    \begin{equation}
        \nabla_{x'}=D\Phi|_x^{-T}(\nabla_x+((\nabla_x (D\Phi|_x^T))\xi')^T\nabla_\xi)=O(S^0)\nabla_x+O(S^{-1})|\xi|\nabla_\xi, \ \nabla_{\xi'}=D\Phi|_x\nabla_{\xi}=O(S^0)\nabla_\xi.
    \end{equation}
    Then iteratively for all multi-indices $\alpha, \beta$ and $j\ge 0$ we have
    \begin{multline}
        \partial_{x_0}^j\partial_{x'}^\alpha\partial_{\xi'}^\beta (\Phi^{*}a)(x_0,x',\xi')=\sum_{\alpha',\beta'}O(S^{|\alpha'|-|\alpha|})|\xi|^{|\beta'|-|\beta|}(\partial_{x_0}^j\partial_{x}^{\alpha'}\partial_{\xi}^{\beta'}a)(x_0, x,\xi)\\
        =O(\langle \xi\rangle^{s}\langle x\rangle^{l-j-|\alpha|+\delta|\beta|}\langle \langle x\rangle^\delta\sqrt{|\xi|^2+M^2}\rangle^{k-|\beta|}),
    \end{multline}
    summing over $|\alpha'|\le |\alpha|$, $|\beta'|\ge |\beta|$ and $|\beta'|+|\alpha'|\le|\beta|+|\alpha|$. Note here we used $\langle \langle x\rangle^\delta\sqrt{|\xi|^2+M^2}\rangle\ge \langle x\rangle^\delta|\xi|$. Note $\langle x\rangle\sim \langle x'\rangle$, $|\xi|\sim |\xi'|$ to conclude the proof.
    \end{proof}

\begin{lemma}\label{l:ChangeCoorOperator}
    Let $\Phi:\R^{d-1}\to \R^{d-1}$ be a diffeomorphism that 
    \begin{equation}
        |\partial^\alpha \Phi(x)|\leq C_\alpha\langle x\rangle^{1-|\alpha|}, \ |\partial^\alpha (\Phi^{-1}(x'))|\leq C_\alpha\langle x'\rangle^{1-|\alpha|}
    \end{equation}
    for every multi-index $\alpha$. Then for all $a\in S^{s,l}_{k,\delta}((T^*\mathbb{R}^{d-1})_+)$, there exists $a'\in S^{s,l}_{k,\delta}((T^*\mathbb{R}^{d-1})_+)$ such that $\oph(a')=\Phi^{-*}\oph(a)\Phi^*$ and
    \begin{equation}
        a'(x_0, x', \xi')=a(x_0,\Phi^{-1}(x'),D\Phi|_{\Phi^{-1}(x')}^T\xi')+O(h S^{s,l-1+\delta}_{k-1,\delta}((T^*\mathbb{R}^{d-1})_+)).
    \end{equation}
\end{lemma}

\begin{proof}
    1. Let $x'=\phi(x)$ and $y'=\phi(y)$. We follow the approach in the proof of \cite[Theorem 9.3]{Zw:12}. Assume for now $a\in \mc S$.
    $$
    \oph(a)(u)=\int K_a(\bar x, \bar y) u(\bar y)d\bar y,
    $$
    with $K_a(\bar x,\bar y):=\frac{\delta(x_0-y_0)}{(2h\pi)^{d-1}}\int a(\bar x, \xi ) e^{\frac{i}{h}\langle  x- y, \xi\rangle}\, d\xi$. Notice that the tangential kernel 
    $$
    k_a:=\frac{1}{(2h\pi)^{d-1}}\int a(\bar x, \xi ) e^{\frac{i}{h}\langle  x- y, \xi\rangle}\, d\xi \in \mc S
    $$
    since $a\in \mc S$. Moreover, the tangential kernel of $(\Phi^{-1})^*\oph(a)\Phi^*=\oph(\tilde a)$ is
    $$
    \tilde k(x_0,x',y')=k_a(x_0,\Phi^{-1}(x'),\Phi^{-1}(y'))|\det \partial(\Phi^{-1})(y')|.
    $$
    One can show by a direct estimate using the hypothesis on the growth of $\Phi$ and $\Phi^{-1}$ that $\tilde k\in \mc S$
    Hence, by Fourier inversion there is $\tilde a\in \mc S$ such that $(\Phi^{-1})^*\oph(a)\Phi^*=\oph(\tilde a)$.  Thus,
    $$
    \oph(\tilde a)(u')=\int K_{\tilde a}(\bar x', \bar y') u(\bar y')d\bar y',
    $$
    with $K_{\tilde a}(\bar x',\bar y'):=\frac{\delta(x_0-y_0)}{(2h\pi)^{d-1}}\int \tilde a(\bar x', \xi') e^{\frac{i}{h}\langle  x'- y', \xi' \rangle}\, d \xi'$.
    Namely,
    $$
    K_a(x_0,x,y)=K_{\tilde a}(x_0,\Phi(x),\Phi(y))|\det \partial\Phi(y)|.
    $$
    Now, 
    \begin{multline*}
        \int K_{\tilde a}(x',y') u(\bar y')d\bar y'= \int K_{\tilde a}(x_0,\Phi(x),y_0,\Phi(y)) u(y_0,\Phi(y)) dy_0 d\Phi( y)\\
        = \int \tilde a( x_0, \Phi(x), \xi') e^{\frac{i}{h}\langle  \Phi(x)- \Phi(y), \xi' \rangle} u(x_0,\Phi(y)) |\det \partial\Phi(y)| \, dy d\xi'.
    \end{multline*}
    
    2. We can assume that $| x- y|$ is small. In fact, a nonstationary phase argument in the integral of $K_{\tilde a}$ shows that for any fixed $\e>0$, the component for $|x-y|>\e$ is $O_{\mc S' \to \mc S}(h^\infty)$. Hence, by Taylor expanding $\Phi(y)$ and $\Phi(x)$ at the point $\frac{x+y}{2}$,
    $$
    \Phi(x)- \Phi(y)=(\partial \Phi)\Big(\frac{x+y}{2} \Big)(x-y)+O(
    \langle x\rangle^{-1}
    |x-y|^2)=G(x,y)(x-y),
    $$
    with $G(x,y)$ invertible for all $|x-y|\leq\e$ ($\e$ small) since $G(x,y)=\partial\Phi(x)+O(\langle x\rangle^{-1}|x-y|)=\partial\Phi(x)+O(\e)$, and $\partial\Phi(x)^{-1}$ is uniformly bounded.
    Denoting $\xi=G(x,y)^T  \xi'$,
    
    \begin{multline}
        \int K_{\tilde a}( x', y') u( y')d y' \\
        =\frac{\delta(x_0-y_0)}{(2h\pi)^{d-1}}  \int \tilde a( x_0, \Phi(x), (G(x,y)^T)^{-1} \xi) e^{\frac{i}{h}\langle   x-  y,  \xi \rangle} u(x_0,\Phi(y)) \, |\det G(x,y)|^{-1} d \xi |\det \partial \Phi(y)| d y \\
        =\frac{\delta(x_0-y_0)}{(2h\pi)^{d-1}}  \int \tilde a( x_0, \Phi(x), (G(x,y)^T)^{-1} \xi) e^{\frac{i}{h}\langle   x-  y,  \xi \rangle} u(x_0,\Phi(y)) \, |\det G(x,y)|^{-1} d \xi |\det (\partial \Phi)(y)| d y \\
        =\frac{\delta(x_0-y_0)}{(2h\pi)^{d-1}} \int \Big(\tilde a\Big( x_0, \Phi(x), \Big(\partial \Phi\Big(\frac{x+y}{2}\Big)^T\Big)^{-1} \xi\Big)
        +O_{S^{s,l}_{k,\delta}}
        (
        \langle x\rangle^{-1}
        | x -  y|)\Big) e^{\frac{i}{h}\langle   x-  y,  \xi 
        \rangle} u(x_0,\Phi(y)) \times \\
        \times \, \Big|\det G(x,y)\Big|^{-1} d \xi |\det (\partial \Phi)(y)| d y,
    \end{multline}
    The $O_{S^{s,l}_{k,\delta}}$ error follows from the Taylor expansion of 
    $$
    b(x,y,\xi):=\tilde a( x_0, \Phi(x), (G(x,y)^T)^{-1} \xi) |\det G(x,y)|^{-1} |\det (\partial \Phi)(y)|
    $$
    at $y=x$.
    In fact, 
    $$
    b(x,y,\xi)=b(x,x,\xi)+\sum_j(y_j-x_j)r_j(x,y,\xi),\quad r_j\in S^{s,l-1}_{k,\delta}.
    $$
    Hence, the integration by part using $(x_j-y_j)e^{\frac{i}{h}\langle x-  y, \xi \rangle}=hD_{\xi_j}e^{\frac{i}{h}\langle x-y,  \xi \rangle}$ gives a remainder in $hS^{s,l-1+\delta}_{k-1,\delta}$.
    Now we complete the proof of the claim by using that for $|x-y|$ small, again by a Taylor expansion
    $$
    \Big|\det G(x,y)\Big|^{-1} |\det \partial \Phi(y)|= 1+O(\langle x\rangle^{-1} |x-y|),
    $$
    and since
    $$
    (x-y)^\alpha e^{\frac ih \langle  x- y,  \xi \rangle}=(-hD_\xi)^\alpha  e^{\frac ih \langle  x- y,  \xi \rangle}.
    $$
    3. We conclude the proof, that is we treat the case $a\in S^{s,l}_{k,\delta}$, by a density argument, see for instance the proof of \cite[Theorem 5.2]{Hi:25} for details.
\end{proof}

\subsection{The calculus on manifolds}
Lemma~\ref{l:ChangeCoorOperator} allows us to quantize symbols inside charts and then we can define a global calculus by a finite sum of such quantization in individual charts. See further in \cite{DyZw:19,Hi:25}. 
\begin{definition}[A second-microlocal calculus on $\Omega_+$]
Consider a finite family of charts covering $(\Phi_{j}, U_j)$ over $\partial\Omega$, associated with a partition of unity $\phi_j\in S^0(\partial\Omega)$ and $\supp\phi_j\subset U_j$ that $\sum_{j}\phi_j=1$. Let $\tilde\phi_j\in S^0(\partial\Omega)$ with  $\supp \phi_{j}\cap\supp(1-\tilde\phi_{j})=\emptyset$, $\supp\tilde\phi_j\subset U_j$. Let $a(x_0, x,\xi)\in S^{s,l}_{k,\delta}((T^*\partial\Omega)_+)$. We
define its global left quantization by
\begin{equation}
    \oph(a)=\sum_{j}\tilde\phi_{j}\Phi_j^{*}\oph(\Phi_j^{-*}(\phi_j(x,\xi) a(x_0, x,\xi)))\Phi_j^{-*}\tilde\phi_{j}, 
\end{equation} 
where we confuse $\Phi$ with the symplectomorphism on $T^*\partial\Omega$ induced by $\Phi$ on $\partial\Omega$. We define the residual class by $h^\infty\Psi^{-\infty,-\infty}_{-\infty,\delta}(\Omega_+)=h^\infty\Psi^{-\infty, -\infty}(\Omega_+)$, consisting of all operators with kernel in $h^\infty \delta(x_0-y_0) C^\infty_b([0,\infty)_{x_0}; \mc{S}(\partial\Omega_x\times\partial\Omega_y)$. We define $\Psi^{s,l}_{k,\delta}(\Omega_+)=\oph(S^{s,l}_{k,\delta}((T^*\partial\Omega)_+))+O(h^\infty\Psi^{-\infty, -\infty}(\Omega_+))$. 
\end{definition}
As the corollary of Lemma~\ref{l:CompositionSymbols}, we have the following definition:
\begin{definition}
    We define the principal symbol map by
    \begin{equation}
        \sigma_h: \Psi^{s,l}_{k,\delta}(\Omega_+)\rightarrow S^{s,l}_{k,\delta}((T^*\partial\Omega)_+)/hS^{s,l-(1-\delta)}_{k-1,\delta}((T^*\partial\Omega)_+),
    \end{equation}
    such that $\sigma_h(\oph(a))=a$ for each $a\in S^{s,l}_{k,\delta}((T^*\partial\Omega)_+)$. Let $b\in S^{s',l'}_{k',\delta}((T^*\partial\Omega)_+)$. We have:
    \begin{gather}
        \sigma_h(\oph(a)\oph(b))=ab, \\
        \sigma_h([\oph(a), \oph(b)])=\{a, b\}=\langle\partial_{\xi}a,\partial_{x}b\rangle-\langle\partial_{x}a,\partial_{\xi}b\rangle.
    \end{gather}
\end{definition}

Moreover, we can as well define the calculus on $\partial\Omega$ by restricting to $x_0=0$. In fact the calculus introduced in this section is uniform in the $x_0$ which is, in fact, treated as a parameter.

\begin{lemma}[Off-chart information is residual]
    Given $\chi,\tilde\chi\in S^0(\Omega_+)$ with $\supp\chi\cap\supp(1-\tilde\chi)=\emptyset$ and $\supp\chi\subset\{(x_0,x): x_0\le C\langle x\rangle\}$ for some $C>0$. Then for any $A\in\Psi^{s,l}_{k,\delta}(\Omega_+)$, we have
    \begin{equation}
        (1-\tilde\chi)A\chi\in 
        h^\infty
        \langle x_0\rangle^{-\infty}\Psi^{-\infty,-\infty}(\Omega_+).
    \end{equation}
\end{lemma}
\begin{proof}
    Disjoint $1-\tilde\chi$, $\chi$ give $(1-\tilde\chi)A\chi\in h^{\infty}\Psi^{-\infty,-\infty}$. Conclude with $x_0\le C\langle x\rangle$ on $\supp\chi$. 
\end{proof}

Now we want to show that the calculus is well defined, that is $\Psi^{s_1,l_1}_{k_1,\delta}\circ \Psi^{s_2,l_2}_{k_2,\delta}\subseteq \Psi^{s_1+s_2,l_1+l_2}_{k_1+k_2,\delta}$. We start with a result about the boundedness of $e^{i\langle D_\xi,D_y\rangle}$.
\begin{lemma} \label{l:Continuity_Asym_Exp}
    Let $s_j,l_j,k_j\in\R$ for $j\in\{1,2\}$, and $M_0>0$. Then, 
    if $X:=(x,\xi)\in \partial\Omega 
    $, $Y:=(y,\eta) \in \R^{2(d-1)}$, and 
    $\chi\in C^\infty_c(\R^{d-1})$, $\chi=1$ near $0$, $\supp \chi\subseteq B(0,1)$,
    \begin{align*}
    \chi \Big(\frac{y-x}{\langle y\rangle+\langle x\rangle}\Big) e^{ih \langle D_\xi,D_y\rangle}: S^{s_1,l_1}_{k_1,\delta}((T^*\partial\Omega_{X})_+)\otimes S^{s_2,l_2}_{k_2,\delta}((T^*\partial\Omega_{Y})_+)
    \to S^{s_1,l_1}_{k_1,\delta}((T^*\partial\Omega_{X})_+)\otimes S^{s_2,l_2}_{k_2,\delta}((T^*\partial\Omega_{Y})_+)
    \end{align*}
    is continuous (uniformly in $h$ and for all $M\in [0,M_0]$) and for all $N\in \mathbb{N}$
    \begin{align}\label{e:AsympExp}
        \chi \Big(\frac{y-x}{\langle y\rangle+\langle x\rangle}\Big) \Big( e^{ih \langle D_\xi,D_y\rangle}a(\bar x,\xi,\bar y,\eta)
        -\sum^{N-1}_{k=0}\frac{(ih)^k}{k!}\langle D_\xi,D_y\rangle^k a(\bar x,\xi,\bar y,\eta)\Big) \\
        \in h^N S^{s_1,l_1-(1-\delta) N}_{k_1-N,\delta}((T^*\partial\Omega_{X})_+) \otimes S^{s_2,l_2}_{k_2,\delta}((T^*\partial\Omega_{Y})_+). \nonumber
    \end{align}
\end{lemma}
\begin{proof}
    1. By \cite[Theorem 4.8 (ii)]{Zw:12}  
    \begin{equation}\label{e:ExpDefS}
        \chi\Big(\frac{y-x}{\langle y\rangle+\langle x\rangle}\Big) e^{ih \langle D_\xi,D_y\rangle}a(\bar x,\xi,\bar y,\eta)=\frac{1}{(2h\pi)^{d-1}}\int_{\R^{2d-2}}\chi \Big(\frac{y-x}{\langle y\rangle+\langle x\rangle}\Big)e^{-\frac{i}{h}\langle \tilde \xi,\tilde y\rangle}a(\bar x,\xi+\tilde \xi,y_0,y+\tilde y,\eta)\,d\tilde y d\tilde \xi,
    \end{equation}
    for $a\in \mathscr{S}$. We can give meaning to this integral as a tempered distribution when $a\in S^{s_1,l_1}_{k_1,\delta}\otimes S^{s_2,l_2}_{k_2,\delta}$ by a stationary phase approach. (See, for instace, the proof of \cite[Theorem 3.18]{Zw:12}.)
    
    2. Now we prove \eqref{e:AsympExp}. First of all, for $\gamma\in (0,1)$,
    \begin{align}
        e^{ih \langle D_\xi,D_y\rangle}a(\bar x,\xi,\bar y,\eta) \nonumber
        =& \frac{1}{(2h\pi)^{d-1}}\int
        e^{-\frac{i}{h}\langle \tilde \xi,\tilde y\rangle}a(\bar x,\xi+\tilde \xi,y_0,y+\tilde y,\eta)\tilde\chi\Bigg(\frac{\tilde y}{\langle y\rangle}, \frac{\tilde \xi}{\langle \xi\rangle^\gamma}\Bigg)\,d\tilde y d\tilde \xi \\ \nonumber
         &+ \frac{1}{(2h\pi)^{d-1}}\int
         e^{-\frac{i}{h}\langle \tilde \xi,\tilde y\rangle}a(\bar x,\xi+\tilde \xi,y_0,y+\tilde y,\eta) (1-\tilde\chi)\Bigg(\frac{\tilde y}{\langle y\rangle}, \frac{\tilde \xi}{\langle \xi\rangle^\gamma}\Bigg)\,d\tilde y d\tilde \xi \\ \nonumber
        =:&  
        \,I_1 (\bar x,\xi,\bar y,\eta)+I_2(\bar x,\xi,\bar y,\eta)
        .
     \end{align}
     
    2a. We start by considering $I_1$. Set $Z:=\frac{\tilde y}{\langle y\rangle}$,
    $\Xi:=\frac{\tilde \xi}{\langle \xi\rangle^\gamma}$,
    and define $\lambda:=\frac{\langle y\rangle\langle \xi\rangle^\gamma}{h}$.
    Then,
    $$
    \begin{aligned}
    I_1(\bar x,\xi,\bar y,\eta)
    &= \frac{
    \langle y\rangle^{d-1}\langle \xi\rangle^{\gamma (d-1)} }{(2h\pi)^{d-1}}
    \int
    e^{-i\lambda\langle \Xi,Z\rangle}
    A(Z,\Xi)
    \,dZ\,d\Xi,
    \end{aligned}
    $$
    where $A(Z,\Xi):=a(\bar x,\xi+\langle \xi\rangle^\gamma\Xi,
        y_0,y+\langle y\rangle Z,\eta)
    \tilde\chi(Z,\Xi)$.
    Now, $A$ is compactly supported
    in $(Z,\Xi)$ uniformly in
    $(\bar x,\xi,\bar y,\eta)$. Thus, stationary phase may be applied
    uniformly with large parameter $\lambda=\frac{\langle y\rangle\langle \xi\rangle^\gamma}{h}$.
    Namely, (see \cite[(3.5.11)]{Zw:12})
    $$I_1(\bar x,\xi,\bar y,\eta)=  \sum_{k=0}^{N-1}
    \frac{(ih)^k}{k!}
    \langle D_\xi,D_y\rangle^k
    a(\bar x,\xi,\bar y,\eta)
    +
     h^N R_N(\bar x,\xi,\bar y,\eta).
    $$
    with
    \begin{equation} \label{e:ResidualExp}
        R_N(\bar x,\xi,\bar y,\eta)=(\langle y \rangle \langle \xi \rangle^\gamma)^{-N}\int^1_0 \frac{(1-t)^{N-1}}{(N-1)!}\Big[e^{i\frac{h}{\langle y\rangle\langle \xi\rangle^\gamma}t\langle D_\Xi,D_Z\rangle}(i\langle D_\Xi,D_Z\rangle)^N  A\Big]\Big|_{Z=\Xi=0}\,dt.
    \end{equation}
    Hence, (see the last line of the proof of  \cite[Theorem 3.13]{Zw:12})
    $$
    |R_N|\leq C_N (\langle y\rangle \langle \xi\rangle^\gamma)^{-N}\sum_{|\alpha|\leq 2d-1} \sup_{Z,\Xi\in \supp \tilde\chi} |\partial^\alpha_{Z,\Xi} \langle D_\Xi,D_Z \rangle^N A(Z,\Xi)|.
    $$
    This is estimated by the symbol seminorms of $a$. Hence, for all $N$, we showed that $R_N\in S^{s_1,l_1+\delta( N+4d-1)}_{k_1-N,\delta}\otimes S^{s_2,l_2-N}_{k_2,\delta}$ and, thus,
    \begin{align*}
        \chi\Big(\frac{y-x}{\langle y\rangle+\langle x\rangle}\Big) &\bigg(I_1(\bar x,\xi,\bar y,\eta)-\bigg( \sum^{N-1}_{k=0} \frac{(ih)^k}{k!} \langle D_{\tilde \xi},D_{\tilde y}\rangle^k a(\bar x,\xi+\tilde \xi,y_0,y+\tilde y,\eta)|_{\tilde\xi=\tilde y=0}\bigg) \bigg)\\
        &\in h^{N} S^{s_1,l_1+(\delta-1)N+\delta(4d-1)}_{k_1-N,\delta}\otimes S^{s_2,l_2}_{k_2,\delta} 
    \end{align*}
    because of the support of $\chi$.
    
    Now, we fix $K\geq \frac{4d-1}{1-\delta}$  so that $S^{s_1,l_1+(\delta-1)K+\delta(4d-1)}_{k_1-N,\delta}\subseteq  S^{s_1,l_1}_{k_1-N,\delta}$ and, hence, for all $N$
    \begin{align*}
        \chi\Big(\frac{y-x}{\langle y\rangle+\langle x\rangle}\Big)& \bigg(I_1(\bar x,\xi,\bar y,\eta)-\bigg( \sum^{N+K-1}_{k=0} \frac{(ih)^k}{k!}
        \langle D_\xi,D_y\rangle^k
        a(\bar x,\xi,y_0,y,\eta)\bigg) \bigg) \\
        &\in h^N S^{s_1,l_1+(\delta-1)N}_{k_1-N,\delta}\otimes S^{s_2,l_2}_{k_2,\delta}.
    \end{align*}
     Moreover, $\chi\Big(\frac{y-x}{\langle y\rangle+\langle x\rangle}\Big)\langle D_\xi,D_y\rangle^k
        a(\bar x,\xi,y_0,y,\eta)\in S^{s_1,l_1+(\delta-1)k}_{k_1-k,\delta}\otimes S^{s_2,l_2}_{k_2,\delta}\subseteq S^{s_1,l_1+(\delta-1)N}_{k_1-N,\delta}\otimes S^{s_2,l_2}_{k_2,\delta}$ for $k\geq N$ which implies
     \begin{align*}
         \chi\Big(\frac{y-x}{\langle y\rangle+\langle x\rangle}\Big) &\bigg(I_1(\bar x,\xi,\bar y,\eta)-\bigg( \sum^{N-1}_{k=0} \frac{(ih)^k}{k!} \langle D_{\tilde \xi},D_{\tilde y}\rangle^k a(\bar x,\xi,\bar y,\eta)
         \bigg) \bigg) \\
        &\in h^N S^{s_1,l_1+(\delta-1)N}_{k_1-N,\delta}\otimes S^{s_2,l_2}_{k_2,\delta}.
     \end{align*}
    
    2b. Now we focus on $I_2$. By changing coordinates $(\tilde y, \tilde \xi)\mapsto(\langle y\rangle \tilde y, \langle \xi\rangle^\gamma \tilde \xi)$
    \begin{align} \label{e:c_2-NonStatPhase}
    I_2= 
    &\langle y\rangle^{d-1}\langle \xi\rangle^{\gamma(d-1)} \int
    e^{-\frac{i}{h}\langle y\rangle \langle \xi\rangle^\gamma\langle \tilde \xi,\tilde y\rangle}a(\bar x,\xi+\langle \xi\rangle^\gamma\tilde \xi,y_0,y+\langle y\rangle\tilde y,\eta) (1-\tilde\chi)(\tilde y,\tilde \xi)\,d\tilde y d\tilde \xi \nonumber \\
    =& \langle y\rangle^{d-1}\langle \xi\rangle^{\gamma(d-1)} \int
    e^{-\frac{i}{h}\langle y\rangle \langle \xi\rangle^\gamma \langle \tilde \xi,\tilde y\rangle}a(\bar x,\xi+\langle \xi\rangle^\gamma \tilde \xi,y_0,y+\langle y\rangle \tilde y,\eta) (\psi_1+\psi_2)(\tilde y,\tilde \xi)\,d\tilde y d\tilde \xi
    \end{align}
    with $1-\tilde\chi=\psi_1+\psi_2$, $\supp \psi_1\subset\{|\tilde y|\geq c\}$, $\supp \psi_2\subset\{|\tilde y|\leq 2c,\,|\tilde\xi|\geq c\}$ for some $c>0$ small. We estimate
    $$
    I'_2= \langle y\rangle^{d-1}\langle \xi\rangle^{\gamma(d-1)} \int_{\R^{2(d-1)}}
    e^{-\frac{i}{h}\langle y\rangle \langle \xi\rangle^\gamma \langle \tilde \xi,\tilde y\rangle}a(\bar x,\xi+\langle \xi\rangle^\gamma \tilde \xi,y_0,y+\langle y\rangle \tilde y,\eta) \psi_1(\tilde y,\tilde \xi)\,d\tilde y d\tilde \xi,
    $$
    by a nonstationaryphase approach using the operator $L_{\tilde y}=-h\frac{\langle \tilde y,D_{\tilde\xi} \rangle}{|\tilde y|^2\langle y\rangle \langle \xi\rangle^\gamma}$ and $L_{\tilde \xi}=-h\frac{\langle \tilde \xi,D_{\tilde y} \rangle}{|\tilde \xi|^2\langle y\rangle \langle \xi\rangle^\gamma}$.
    In fact, 
    \begin{multline}
        D_{\tilde \xi}a(\bar x,\xi+\langle \xi\rangle^\gamma \tilde \xi,y_0,y+\langle y\rangle \tilde y,\eta)=\langle \xi\rangle^\gamma (D_\xi a)(\bar x,\xi+\langle \xi\rangle^\gamma \tilde \xi,y_0,y+\langle y\rangle \tilde y,\eta)\\
        \in \langle \xi\rangle^\gamma  \langle \langle x\rangle^\delta\sqrt{|\xi+\langle \xi \rangle^\gamma \tilde\xi|^2+M^2}\rangle^{k_1-1}  S^{K,K+\delta}_{0,\delta}\times S^{s_2,l_2}_{k_2,\delta}
    \end{multline}
    for some $K>0$. Hence,
    $L_{\tilde y}a\in h \langle y\rangle^{-1}\langle \tilde  y\rangle^{-1} \langle \langle x\rangle^\delta\sqrt{|\xi+\langle \xi \rangle^\gamma \tilde\xi|^2+M^2}\rangle^{k_1-1}  S^{K,K+\delta}_{0,\delta}\times S^{s_2,l_2}_{k_2,\delta}$ and, iterating, $L_{\tilde y}^Na\in h^N\langle y\rangle^{-N}\langle \tilde  y\rangle^{-N} \langle \langle x\rangle^\delta\sqrt{|\xi+\langle \xi \rangle^\gamma \tilde\xi|^2+M^2}\rangle^{k_1-N}  S^{K,K+\delta N}_{0,\delta}\times S^{s_2,l_2}_{k_2,\delta}$ with the same $K$. 
    Now, for large $N$, $L_{\tilde y}^Na\in h^N\langle y\rangle^{-N}\langle \tilde  y\rangle^{-N} \langle \xi\rangle^{-J_N}  \langle \tilde \xi\rangle^{-J_N}S^{K,K+\delta N}_{0,\delta}\times S^{s_2,l_2}_{k_2,\delta}$ when by a smooth partition we restrict to $\langle\xi+\langle \xi \rangle^\gamma \tilde\xi\rangle>\e (\langle \xi \rangle+\langle \xi \rangle^\gamma |\tilde\xi|)$, $\e$ small, with $J_N\to \infty$ as $N\to \infty$.
    For $\langle\xi+\langle \xi \rangle^\gamma \tilde\xi\rangle\leq \e (\langle \xi \rangle+\langle \xi \rangle^\gamma |\tilde\xi|)$ we iteratively integrate by part using the operator $L_{\tilde \xi}=-h\frac{\langle \tilde \xi,D_{\tilde y} \rangle}{|\tilde \xi|^2\langle y\rangle \langle \xi\rangle^\gamma}$. Namely, we get $L_{\tilde\xi}^NL_{\tilde y}^Na\in h^{2N}\langle y\rangle^{-N}\langle \tilde  y\rangle^{-N} \langle \xi\rangle^{-J_N}  \langle \tilde \xi\rangle^{-J_N}S^{K,K+\delta N}_{0,\delta}\times S^{s_2,l_2}_{k_2,\delta}$. 

    Now we estimate 
    $$
    I''_2= \langle y\rangle^{d-1}\langle \xi\rangle^{\gamma(d-1)} \int_{\R^{2(d-1)}}
    e^{-\frac{i}{h}\langle y\rangle \langle \xi\rangle^\gamma \langle \tilde \xi,\tilde y\rangle}a(\bar x,\xi+\langle \xi\rangle^\gamma \tilde \xi,y_0,y+\langle y\rangle \tilde y,\eta) \psi_2(\tilde y,\tilde \xi)\,d\tilde y d\tilde \xi,
    $$
    and this can be done as for $I'_2$ using the operator $L_{\tilde \xi}=-h\frac{\langle \tilde \xi,D_{\tilde y} \rangle}{|\tilde \xi|^2\langle y\rangle \langle \xi\rangle^\gamma}$.

    Hence, we have that for every $J>0$, choosing the number of integrations $N$ sufficiently large
    $$
    \chi \Big(\frac{y-x}{\langle y\rangle+\langle x\rangle}\Big)I_2 \in h^J S^{-2J,-2J}_{K_0,\delta}\otimes  S^{s_2,l_2}_{k_2,\delta}
    \subseteq  h^J S^{-J,-J}_{K_0-J,\delta}\otimes  S^{s_2,l_2}_{k_2,\delta},
    $$
    with $K_0$ uniform in $J$ as claimed.
\end{proof}

Now we show a preliminary result of boundedness for the quantization of a symbol in our class.
\begin{lemma} \label{l:BadBdd}
    Let $a\in S^{s,l}_{k,\delta}((T^*\R^{d-1})_+)$ for $s,l,k\in\R$. Then, there is $N\in\mathbb{N}$ such that $\Op_h(a):H^{2N,2N}((\partial\Omega)_+)\to H^{0,0}((\partial\Omega)_+)$ and $\Op_h(a):H^{0,0}((\partial\Omega)_+)\to H^{-2N,-2N}((\partial\Omega)_+)$ are bounded uniformly in $M$ and $h$.
\end{lemma}
\begin{proof}
    1. We want to prove that $\Op_h(a):H^{2N,2N}\to H^{0,0}$ is bounded for some $N$. Let $u\in H^{2N,2N}$ and $v:=\langle x\rangle^{2N} \langle h D_x\rangle^{2N}  u$. We follow a Schur test argument,
    $$
    \oph(a)u(x)=\int k(x,y)v(y)\,dy,\quad
    k(x,y)=(2\pi h)^{-(d-1)}\int e^{\frac{i}{h}\langle x-y,\xi \rangle}a(x,\xi)\langle \xi\rangle^{-2N} \langle y\rangle^{-2N}\,d\xi.
    $$
    By a stationary phase argument, for some $J>0$ and for $K>J+(1-\delta)^{-1}(|s|+|l|+|k|+d)$, $N>|s|+|l|+|k|+K+d$ 
    \begin{align}
        |k(x,y)|\leq& Ch^{-(d-1)}\langle h^{-1}(x-y) \rangle^{-2K}\langle y\rangle^{-2N}\int_{\R^{d-1}} |(1-\Delta_\xi)^K(a(x,\xi)\langle \xi \rangle^{-2N})|\,d\xi \nonumber\\
        \leq& Ch^{-(d-1)}\langle h^{-1}(x-y) \rangle^{-2K}\langle y\rangle^{-2N} \langle x \rangle^{J} \nonumber
        \leq Ch^{-(d-1)}\langle h^{-1}(x-y) \rangle^{J-2K}, \label{e:SchurKer}
    \end{align}
    where in the last estimate we used that by Peetre's estimate $\langle y\rangle^{-2N} \langle x \rangle^{J}\leq \langle x-y \rangle^{J}$.
    Now by the change of variables $x\mapsto hx+y$ and $y\mapsto hy+x$,
    $$
    \sup_x \int |K(x,y)|\,dy+\sup_y \int |K(x,y)|\,dx\leq C,
    $$
    uniformly in $h$. Hence, 
    $\|\oph(a)u\|\leq C\|v\|=C\|\langle x\rangle^{2N} \langle h D_x\rangle^{2N}  u\|$, as claimed.
    
    2. Now, we prove the boundedness of $\oph(a):H^{0,0}\to H^{-2N,-2N}$ by duality.
    Namely, the kernel of $\oph(a)^*$ is $(2\pi h)^{-(d-1)}\int e^{\frac{i}{h}\langle x-y,\xi \rangle}\overline{a(y,\xi)}\,d\xi$
    and we conclude similarly to Step 1.
\end{proof}

From these lemma a composition law for symbols in the class we are studying follows

\begin{lemma}[Composition of symbols in $S^{s_1,l_1}_{k_1,\delta}$] \label{l:CompositionSymbols}
Let $a\in S^{s_1,l_1}_{k_1,\delta}((T^*\partial\Omega)_+)$ and $b\in S^{s_2,l_2}_{k_2,\delta}((T^*\partial\Omega)_+)$. Then, there exists $c\in S^{s_1+s_2,l_1+l_2}_{k_1+k_2,\delta}((T^*\partial\Omega)_+)$ such that $\oph(a)\oph(b)=\oph(c)$ with
\begin{equation}\label{e:ComposAsym}
    c-\sum^{N-1}_{k=0}\frac{(ih)^k}{k!}\langle D_\xi,D_y\rangle^k a(x_0,x,\xi)b(x_0,y,\eta)|_{y=x,\,\eta=\xi}\in h^{N}S^{s_1+s_2,l_1+l_2-(1-\delta) N}_{k_1+k_2-N,\delta}((T^*\partial\Omega)_+)
\end{equation}
\end{lemma}
\begin{proof}
    1. We first consider $a$ and $b$ uniformly in $C^\infty_b([0,\infty)_{x_0};\mc S(\partial\Omega))$. By the composition formula \cite[(4.3.14)]{Zw:12}
    \begin{equation}\label{e:ComposSchw}
        c(x,\xi):=e^{ih\langle D_\xi, D_y\rangle}  a(x_0,x,\xi)  b(x_0,y,\eta)|_{y=x,\,\eta=\xi}.
    \end{equation}
    By Lemma \ref{l:Continuity_Asym_Exp},
    $$
    \chi\Big(\frac{y-x}{\langle x\rangle+\langle y \rangle}\Big)e^{ih \langle D_\xi,D_y \rangle}a(x_0,x,\xi)  b(y_0,y,\eta)\in S^{s_1,l_1}_{k_1,\delta}\otimes S^{s_2,l_2}_{k_2,\delta}.
    $$
    Restricting to $y_0=x_0$, $y=x$ and $\eta=\xi$, we have
    $c\in S^{s_1+s_2,l_1+l_2}_{k_1+k_2,\delta}$.

    2. We have an asymptotic expansion for $c$. In fact, by Lemma \ref{l:Continuity_Asym_Exp}
    $$
    \chi\Big(\frac{y-x}{\langle x\rangle+\langle y \rangle}\Big)\Big(e^{ih \langle D_\xi,D_y \rangle}a(x_0,x,\xi)  b(y_0,y,\eta)-\sum_{j=0}^{N-1}\frac{(ih)^j}{j!}\langle D_\xi,D_y \rangle^j a(x_0,x,\xi)  b(y_0,y,\eta) \Big)\in h^N S^{s_1,l_1-(1-\delta)N}_{k_1-N,\delta}\otimes S^{s_2,l_2}_{k_2,\delta}
    $$
    which implies
    $$
    c-\sum_{j=0}^{N-1}\frac{(ih)^j}{j!}\langle D_\xi,D_y \rangle^j a(x_0,x,\xi)  b(y_0,y,\eta) |_{y_0=x_0,\,y=x,\,\eta=\xi}\in h^N S^{s_1+s_2,l_1+l_2-(1-\delta)N}_{k_1+k_2-N,\delta}.
    $$
    
    3. Now, we treat the case $a,b\notin \mc S$ by a density argument.
    Let $\chi\in C^\infty_c(T^*\partial\Omega)$, $\chi=\chi(|x|,\xi)$,
    $\chi|_{|x|+|\xi|_g<1}=1$. Then, let $a_m(x_0,x,\xi):=a(x_0,x,\xi)\chi(|x|/m,\xi/m)$ and $b_m(y_0,y,\eta):=b(y_0,y,\eta)\chi(|y|/m,\eta/m)$. $a_m, b_m\in \mc S$ uniformly for all $x_0$, $y_0$. $a_m\to a$ in $S^{s_1+0,l_1+0}_{k_1+0,\delta}$ and $b_m\to b$ in $S^{s_2+0,l_2+0}_{k_2+0,\delta}$.
    Hence, by Step 1, $\oph(c_m)=\oph(a_m)\oph(b_m)$ with
    $$
    c_m(x,\xi):=e^{ih\langle D_\xi, D_y\rangle}  a_m(x_0,x,\xi)  b_m(x_0,y,\eta)|_{y=x,\,\eta=\xi}\in S^{s_1+s_2,l_1+l_2}_{k_1+k_2,\delta}.
    $$
    Then, by Lemma \ref{l:Continuity_Asym_Exp}, $c_m\to c$ in $S^{s_1+s_2+0,l_1+l_2+0}_{k_1+k_2+0,\delta}$ with 
    \begin{multline}
        c(x,\xi):=e^{ih\langle D_\xi, D_y\rangle}  a(x_0,x,\xi)  b(x_0,y,\eta)|_{y=x,\,\eta=\xi}=\\
    \sum^{N-1}_{k=0}\frac{(ih)^k}{k!}\langle D_\xi,D_y\rangle^k a(x_0,x,\xi)b(x_0,y,\eta)|_{y=x,\,\eta=\xi}+h^{N}S^{s_1+s_2,l_1+l_2-(1-\delta) N}_{k_1+k_2-N,\delta}((T^*\partial\Omega)_+).
    \end{multline}
    Hence, $c\in S^{s_1+s_2,l_1+l_2}_{k_1+k_2,\delta}((T^*\partial\Omega)_+)$.
    
    4. It remains to show that $\oph(c)=\oph(a)\oph(b)$. By Lemma \ref{l:BadBdd},
    $$
    \langle (\oph(a_m)\oph(b_m)-\oph(a)\oph(b))u,v \rangle\to  0,
    $$
    and
    $$
    \langle (\oph(c_m)-\oph(c))u,v \rangle\to  0.
    $$
    Passing to the limit, $\oph(c)=\oph(a)\oph(b)$ as $\mc S \to \mc S'$ bounded operators, as claimed. 
\end{proof}

We show that $\Psi^{0,0}_{0,\delta}(\Omega_+)$ are bounded on $L^2$. 
\begin{lemma} \label{l:GoodBdd}
    Let $a\in S^{0,0}_{0,\delta}((T^*\partial\Omega)_+)$. Then $\oph(a)$ is bounded on $C^\infty_b([0,\infty)_{x_0}; L^2(\partial\Omega))$, uniformly in $h\in(0, h_0]$ and $M\in [0, M_0]$. 
\end{lemma}
\begin{proof}
    1. $\Op(a):L^2\to L^2$ is bounded if $a\in S^{-N,-N}_{-N,\delta}$ for some large $N$
    by Lemma \ref{l:BadBdd} and Lemma \ref{l:CompositionSymbols}.
    
    2. Now, for all $\e>0$, $\oph(a):L^2\to L^2$ is bounded if $a\in S^{-\e,-\e}_{-\e,\delta}$. In fact, $\oph(a):L^2\to L^2$ is bounded if $\oph(a)^*\oph(a):L^2\to L^2$ is bounded. Hence, if $a \in S^{-N/2,-N/2}_{-N/2,\delta}$ then $\oph(a)$ is bounded onto $L^2$ since $\oph(a)^*\oph(a)\in \Psi^{-N,-N}_{-N,\delta}$ by Lemma \ref{l:CompositionSymbols}. By repeating this procedure $\oph(a):L^2\to L^2$ is bounded when 
    $a\in S^{-\e,-\e}_{-\e,\delta}$.\\
    3. Finally, let $a\in S^{0,0}_{0,\delta}$.  Let $c:=\sqrt{C-|a|^2}\in S^{0,0}_{0,\delta}$ with $C>\sup |a|^2+1$,
    \begin{multline}
        0\leq \langle \oph(c)^*\oph(c)u,u\rangle=\langle \oph(|c|^2+r)u,u\rangle=\langle \oph(C-|a|^2+r)u,u\rangle\\
        = C\|u\|^2-\|\oph(a)u\|^2+\langle \oph(\tilde r)u,u\rangle
        \leq C\|u\|^2-\|\oph(a)u\|^2+\tilde C \|u\|^2,
    \end{multline}
    where $r,\tilde r\in S^{0,\delta-1}_{-1,\delta}\subset S^{\frac{\delta-1}{2},\frac{\delta-1}{2}}_{\frac{\delta-1}{2},\delta}$ and the last inequality follows from Step 2. Hence, $\|\oph(a)u\|\leq C \|u\|$, as claimed.
\end{proof}

We use the properties of the calculus to show that we can define for all $s,l,k\in\R$ an invertible operator in $S^{s,l}_{k,\delta}$. 

\begin{lemma}\label{l:GammaInverse}
Let $N>0$. Let $s, l\in \mathbb{R}$ and $k\geq 0$. We define $\Gamma^{s,l}_{k,\delta}\in \Psi^{s,l}_{k,\delta}(\Omega_+)$ by
\begin{equation}
    \Gamma^{s,l}_{k,\delta}=(1-h^2\Delta_x)^{\frac s4}\langle x\rangle^{ {\frac l2}}\oph(\langle \langle x\rangle^\delta  \sqrt{|\xi|^2+M^2}\rangle^\frac{k}{2})^*\oph(\langle \langle x\rangle^\delta  \sqrt{|\xi|^2+M^2}\rangle^\frac{k}{2})\langle x\rangle^{ {\frac l2}}(1-h^2\Delta_x)^{\frac s4}
\end{equation}
Then there exists $h_0>0$ such that for $h\in(0,h_0]$, $k\le N$, $\Gamma^{s,l}_{k,\delta}:\mc S'((\partial\Omega)_+)\to \mc S'((\partial\Omega)_+)$ is invertible with $(\Gamma^{s,l}_{k,\delta})^{-1}\in \Psi^{-s,-l}_{-k,\delta}$. For $k<0$, we then define $\Gamma^{s,l}_{k,\delta}=(\Gamma^{-s,-l}_{-k,\delta})^{-1}\in \Psi^{s,l}_{k,\delta}((\partial\Omega)_+)$. 
\end{lemma}
\begin{proof}
    It suffices to show the invertibility of $A=\oph(\langle \langle x\rangle^\delta  \sqrt{|\xi|^2+M^2}\rangle^\frac{k}{2})$. Construct its principal parametrix $B_0=\oph(\langle \langle x\rangle^\delta  \sqrt{|\xi|^2+M^2}\rangle^{-\frac{k}{2}})$ and we have $hR=AB_0-\Id\in \Psi^{0,\delta-1}_{-1, \delta}$. Now for small $h$, $(\Id+hR)$ is invertible and we construct $
        B=B_0\sum_{j=0}^\infty (-hR)^{j}\in \Psi^{0,0}_{-\frac{k}2,\delta}$.
    This gives $AB=\Id$ and a right inverse as desired. The existence of a left inverse is similar. 
\end{proof}

Now, we construct the microlocal parametrices in this calculus.
\begin{lemma}\label{l:MicroPar}
    Let $a\in S^{s,l}_{k,\delta}$, $b\in S^{s',l'}_{k',\delta}$. Assume that $a$ is elliptic near $\supp b$, that is there exists $C>0$ such that
    $$
    |a(\bar x,\xi)|\geq C\langle \xi\rangle^s \langle x\rangle^l \langle \langle x\rangle^\delta \sqrt{|\xi|^2+M^2}\rangle^k \text{ near }\supp b,
    $$
    Then, there exists $E\in \Psi^{s'-s,l'-l}_{k'-k,\delta}$ such that
    $$
    \Op_h(b)=E\Op_h(a)+h^\infty\Psi^{-\infty,-\infty}_{-\infty,\delta}=\Op_h(a)E+h^\infty\Psi^{-\infty,-\infty}_{-\infty,\delta}.
    $$
    Moreover, the parametrices depend continuously on $a$ and $b$ in the corresponding symbol seminorms.
\end{lemma}
\begin{proof}
    1. We proceed by induction.
    Let $\chi\in S^{0,0}_{0,\delta}$, $a$ elliptic near $\supp \chi$, $\chi=1$ near $\supp b$.     Then, $e_0:=\frac{\chi}{a}\in S^{-s,-l}_{-k,\delta}$, $\supp e_0\subset \supp \chi$. Then, by Lemma \ref{l:CompositionSymbols}, $\Op_h(\chi)=\Op_h(e_0)\Op_h(a)+h\Op_h(r_1)+h^\infty\Psi^{-\infty,-\infty}$ with $r_1\in S^{0,\delta-1}_{-1,\delta}$, $\supp r_1\subset \supp \chi$. 
    
    2. Now, assume there exist $e_0,\dots,e_{N-1}$ such that 
    $$
    \Op_h(\chi)=\Big(\sum^{N-1}_{j=0} h^j\oph(e_j)\Big)\Op_h(a)+h^N\Op_h(r_N) +h^\infty\Psi^{-\infty,-\infty} 
    $$ 
    with $r_N\in S^{0,-N(1-\delta)}_{-N,\delta}$, $\supp r_N\subset \supp\chi$. Then, with $e_N=\frac{r_N}{a}\in S^{-s,-l-N(1-\delta)}_{-k-N,\delta}$
    $$
    \Op_h(\chi)=\Op_h\Big(\sum^{N}_{j=0} h^je_j\Big)\Op_h(a)+h^{N+1}\Op_h(r_{N+1})+h^\infty\Psi^{-\infty,-\infty},
    $$
    with $r_{N+1}\in S^{0,(N+1)(\delta-1)}_{-(N+1),\delta}$, $\supp r_{N+1}\subset \supp \chi$, concluding the induction.

    3. Now, let $e\sim \sum_{j\geq 0}h^j e_j$. Then, $E:=\oph(b)\oph(e)\in \Psi^{s'-s,l'-l}_{k'-k,\delta}$, as claimed.

    4. For the construction of the symbol of $E_R$ we can proceed similarly.
    The continuity statement follows from the recursive construction and the fixed ellipticity
    lower bound of $a$ on $\supp \chi$. 
\end{proof}

\newcommand{\hardr}{\mathfrak{r}}
\section{The scattering geometry and Hardy constants}\label{a:2}

\subsection{The scattering geometry}
We show that the boundary of any asymptotically conic domain is a scattering manifold. 
\begin{lemma}\label{l:AC=Scattering}
    The boundary of an \acs domain $\Omega\subset \R^d$ is a (Lipschitz) scattering manifold. 
\end{lemma}
\begin{proof}
    The case $d=2$ is trivial and we assume $d>2$. We just need to write what is the metric induced by $\R^d$ on $\partial\Omega$.
    In fact, 
    $\psi=\psi(\rho,\omega)\in C^\infty([0,1]\times  \partial Y; \mathbb{S}^{d-1})$ that $\psi(0,\omega)=\omega$, $\rho^{-1}\psi(\rho,\omega)$ parametrizes $\partial \Omega$. 
    Then in the local coordinates $\omega=(\omega_1, \dots, \omega_{d-2})$, we have
    \begin{equation}
            g=\rho^{-4}d\rho^2+\rho^{-2}\Big(|\partial_\rho \psi|^2 d\rho^2+2\sum_{j=1}^{d-2} \langle \partial_{\rho}\psi, \partial_{\omega_j}\psi \rangle d\rho d\omega^j+ \sum_{i,j=1}^{d-2} \langle \partial_{\omega_i}\psi, \partial_{\omega_j}\psi\rangle d\omega^i d\omega^j\Big).
    \end{equation}
    Note the smoothness of $\psi$ at $\rho=0$ implies
    \begin{equation}
        g(0,\omega)=\rho^{-4}d\rho^2+g_\infty(\omega), \ g_\infty(\omega)=\sum_{i,j=1}^{d-2} \langle \partial_{\omega_i}\psi(0, \omega), \partial_{\omega_j}\psi(0, \omega)\rangle d\omega^i d\omega^j
    \end{equation}
    where $g_\infty\in \operatorname{Sym}^2(T^*\partial Y))$ is positive-definite. Then \cite[Proposition 2.1]{JoSa:99} implies $\partial\Omega$ is a scattering manifold. 
\end{proof}

\subsection{Properties of the Hardy constant}
\begin{lemma}\label{l:fixhardy}
    Let $\Omega$ be an \acl domain. Then there is some $r_0\ge 0$ that for each $\e>0$, there exists $\gamma\in (0,\frac12)$ that $\hardy^{0}_{U^{r_0}_{\gamma}, \partial\Omega}\ge \frac12-\e$. If $\partial\Omega$ is only singular at up to one point, one can take $r_0=0$.
\end{lemma}
\begin{proof}
    There exists $r_0\ge 0$ such that $\partial\Omega\cap\{|z|> r_0\}$ is smooth. There then exists $\gamma\in(0, \frac12)$ that one can adopt the Fermi normal coordinates identifying $U^{r_0}_{\gamma}$ with a subset of $\{(x_0, x):x\in\partial\Omega\cap\{|z|> r_0\}, x_0\le 2\gamma |x|\}$ equipped with metric $g_+(x_0,x)=dx_0^2+g(x)+O(\gamma \operatorname{Sym}^2(T^*\partial\Omega))$. This implies
    \begin{equation}
        \sqrt{g_+(x_0,x)}=(1+O(\gamma))\sqrt{g(x)}, \ \dist(z,\partial\Omega)=x_0, \ |\partial_{x_0} u|^2\le |\nabla_{g_+} u|^2
    \end{equation}
    on $U^{r_0}_{\gamma}$. Let $u\in \dot H^1(U^{r_0}_{\gamma})$, extended trivially by $0$ onto the rest of $\{(x_0, x):x\in\partial\Omega\cap\{|z|> r_0\}, x_0\le 2\gamma |x|\}$. Now the one-dimensional Hardy inequality in $x_0$ implies
    \begin{equation}
        \int_{x_0=0}^\infty \Big|\frac{u(x_0, x)}{x_0}\Big|^2\,dx_0\le 4\int_{x_0=0}^\infty \Big|\partial_{x_0}u(x_0, x)\Big|^2\,dx_0
    \end{equation}
    for all $u$ with $u(0,x)=0$. We have
    \begin{multline}
        \Big\|\frac{u}{\dist(z, \partial\Omega)}\Big\|^2=(1+O(\gamma))\int_{\partial\Omega}\int_{x_0=0}^\infty \Big|\frac{u(x_0, x)}{x_0}\Big|^2\,\sqrt{g(x)}dx_0dx\\ 
        \le 4(1+O(\gamma)) \int_{\partial\Omega}\int_{x_0=0}^\infty \Big|\nabla_{g_+} u(x_0, x)\Big|^2\,\sqrt{g(x)}dx_0dx=4(1+O(\gamma))\|\nabla_{g_+} u\|^2.
    \end{multline}
    Make $\gamma$ small to conclude the proof for generic \acl domains. If $\partial\Omega$ is singular at up to one point, one can set that point $z=0$ and the above works with $r_0=0$. 
\end{proof}
Similarly one can show:
\begin{lemma}\label{l:fixhardypolygon}
    Let $\mc{P}$ be a curvilinear polygon, and let $\mc{V}=\{x_j\}_{j=1}^{J}$ contain all vertices of $\mc{P}$. Then for each $\e>0$, there exists $\gamma\in (0,\frac12)$ that $\hardy^{0}_{U^{0}_{\gamma}, \partial\Omega}\ge \frac12-\e$. 
\end{lemma}

\begin{lemma}\label{l:Fprop}
Suppose $\Omega\subset \mathbb{R}^2$ such that there are $\gamma_\pm\in C^\infty(\mathbb{R};\mathbb{R}^2)$ and $\alpha\in(0,2\pi)$ such that $\gamma_\pm(0)=0$, $|\gamma_\pm'(s)|=1$, and
\begin{equation}
    \partial \Omega\cap B(0,1)=\big(\gamma_+ ( [0,\infty))\cup \gamma_- ( [0,\infty)\big)\cap B(0,1).
\end{equation}
Let $\alpha\in (0, 2\pi)$ be the angle from $\gamma'_+(0)$ to $\gamma'_-(0)$. 
Then there is a neighborhood $U$ of $0$ and $F:\overline{U\cap \mathcal{W}_{\alpha}}\to \overline{\Omega}$, $F\in C^{1,1}(\overline{\mathcal{W}_{\alpha}})\cap C^\infty(\overline{\mathcal{W}_{\alpha}}\setminus\{0\})$ such that  
$$
F(s,0)=\gamma_+(s),\qquad F(s\cos \alpha,s\sin\alpha)=\gamma_-(s),
$$
and
$$
dF(s,0)\begin{pmatrix}0\\-1\end{pmatrix}=n_+(s),\qquad dF(s\cos \alpha,s\sin\alpha)\begin{pmatrix}-\sin \alpha\\ \cos \alpha\end{pmatrix}=n_-(s),
$$
where $n_\pm$ are the outward unit normals to $\Omega$ at $\gamma_\pm(s)$. 
\end{lemma}
\begin{proof}
1. Without loss of generality, we may assume that $\gamma_+'(0)=(1,0)$ and $\gamma_-'(0)=(\cos \alpha,\sin \alpha)$. We now use polar coordinates $\Theta_\pm:[0,\e)\to \R$ (for some $\e>0$) with $\Theta_-(0)=\alpha$, $\Theta_+(0)=0$ such that 
$$
\gamma_\pm=|\gamma_\pm|(\cos\Theta_\pm,\sin \Theta_\pm) \text{ on } [0,\e).
$$

2. Let $q_\pm (r)=\frac{|\gamma_\pm(r)|}{r}$ for $r>0$, $q_\pm(0)=1$ by smooth extension. We construct $\tilde{R}(r,\theta)$ and $\Theta(r,\theta)\in C^\infty$ such that 
\begin{gather}
\Theta(r,\alpha)=\Theta_-(r),\qquad\Theta(r,0)=\Theta_+(r),\qquad 
\tilde{R}(r,\alpha)=q_-(r),\qquad
\tilde{R}(r,0)=q_+(r),\nonumber\\
\Theta(0,\theta)=\theta,\qquad \tilde{R}(0,\theta)=1.\label{e:at0}
\end{gather}
 Let $\tilde R$, $\Theta\in C^\infty([0,\e)_r\times [0,\alpha]_\theta)$ that
$\tilde R(r,0)=q_+(r)$, $\tilde R(r,\alpha)=q_-(r)$, $\Theta(r,0)=\Theta_+(r)$, $\Theta(r,\alpha)=\Theta_-(r)$, $\Theta(0,\theta)=\theta$ and 
$$
\partial_\theta \tilde R(r,0)=-\langle n_+(r), (\cos (\Theta_+(r)),\sin(\Theta_+(r))) \rangle, \quad \partial_\theta \Theta(r,0)=-\frac{\langle n_+(r), (-\sin(\Theta_+(r),\cos(\Theta_+(r))) \rangle}{q_+(r)},
$$
$$
\partial_\theta \tilde R(r,\alpha)=\langle n_-(r), (\cos (\Theta_-(r)),\sin(\Theta_-(r))) \rangle, \quad \partial_\theta \Theta(r,\alpha)=\frac{\langle n_-(r), (-\sin(\Theta_-(r),\cos(\Theta_-(r))) \rangle}{q_-(r)},
$$

Then,
$$
F (r\cos \theta,r\sin \theta):= r\tilde{R}(r,\theta)\big(\cos\Theta(r,\theta),\sin\Theta(r,\theta)\big),
$$
satisfies the required properties.

3.  It remains to show that that $F\in C^{1,1}$. By Taylor expanding $\tilde R$ and $\Theta$ at $r=0$,
$$
F (r\cos \theta,r\sin \theta)= r(\cos(\theta),\sin(\theta))+r^2 \tilde F(r,\theta),
$$
with $\tilde F\in C^\infty$ bounded with all bounded derivatives. Since $DF=I+O(r)$ and $D^2 F=O(1)$, $F\in C^{1,1}$.
\end{proof}

\begin{remark}
    In fact, when $\alpha\neq \pi$, it is possible to find $F\in C^\infty(\overline{W_\alpha})$ mapping $\overline{W_\alpha}$ to $\overline{\Omega}$. However, when $\alpha=\pi$ there is an obvious obstruction to this if $\partial\Omega\notin C^\infty$ since $\partial\mathcal{W}_\pi \in C^\infty$. 
\end{remark}

\begin{lemma}\label{l:DomComp}
    Let $\Omega_1, \Omega_2\subset\mathbb{R}^d$ be two uniformly Lipschitz domains, with $\{z=0\}\subset\partial\Omega_1\cap\partial\Omega_2$. Let $F:\Omega_1\to \Omega_2$ be a bi-Lipschitz map such that $F(0)=0$. Then, for every $r_0\ge 0$, we have
    \begin{gather}
        \sqrt{\|\det DF\|_{L^\infty}\|\det DF^{-1}\|_{L^\infty}}\|DF\|_{L^\infty}\||F^{-1}(z)|/|z|\|_{L^\infty}\rardy^{r_0}_{\Omega_2}\geq \rardy^{r_0\||F^{-1}(z)|/|z|\|_{L^\infty}}_{\Omega_1},\\
        \sqrt{\|\det DF^{-1}\|_{L^\infty}\|\det DF\|_{L^\infty}}\|DF^{-1}\|_{L^\infty}\||F(z)|/|z|\|_{L^\infty}\rardy^{r_0}_{\Omega_1}\geq \rardy^{r_0\||F(z)|/|z|\|_{L^\infty}}_{\Omega_2}.
    \end{gather}
\end{lemma}
\begin{proof}
    Let $u\in \dot H^1(\Omega_2)$ and $\tilde r_0=\||F^{-1}(z)|/|z|\|_{L^\infty} r_0$. Then, 
\begin{multline}
    \Big\|\frac{u(z)}{|z|+r_0}\Big\|_{L^2(\Omega_2)}^2=\int_{\Omega_1} \frac{|u(F(\tilde z))|^2}{(|F(\tilde z)|+r_0)^2}|\det DF(\tilde z)|\,d\tilde z\\
    \leq \|\det DF\|_{L^\infty(\Omega_1)}\||F^{-1}(z)|/|z|\|_{L^\infty}^2\int_{\Omega_1} \frac{|u(F(\tilde z))|^2}{(|\tilde z|+\tilde r_0)^2}\,d\tilde z\\
    \leq (\rardy^{\tilde r_0}_{\Omega_1})^{-2}\|\det DF\|_{L^\infty(\Omega_1)}\||F^{-1}(z)|/|z|\|_{L^\infty}^2\int_{\Omega_1}|\nabla (u\circ F)(\tilde z)|^2\,d\tilde z\\
    \le (\rardy^{\tilde r_0}_{\Omega_1})^{-2}\|\det DF\|_{L^\infty(\Omega_1)}\||F^{-1}(z)|/|z|\|_{L^\infty}^2\int_{\Omega_1}|(\nabla u)(F(\tilde z))|^2|DF(\tilde z)|^2\,d\tilde z\\
    \leq (\rardy^{\tilde r_0}_{\Omega_1})^{-2}\|\det DF\|_{L^\infty(\Omega_1)}\|DF\|_{L^\infty(\Omega_1)}^2\||F^{-1}(z)|/|z|\|_{L^\infty}^2\int_{\Omega_2}|(\nabla u)(z)|^2|\det DF^{-1}(z)|\,dz\\
    \leq (\rardy^{\tilde r_0}_{\Omega_1})^{-2}\|\det DF\|_{L^\infty(\Omega_1)}\|DF\|_{L^\infty(\Omega_1)}^2\|\det DF^{-1}\|_{L^\infty(\Omega_2)}\||F^{-1}(z)|/|z|\|_{L^\infty}^2\|\nabla u\|_{L^2(\Omega_2)}^2. 
\end{multline}
Hence, the lemma follows after reversing the roles of $\Omega_1$ and $\Omega_2$.
\end{proof}

We determine the Hardy constant on exact cones in any dimensions. 
\begin{lemma}\label{l:conerardy}
     Let $Y\subset \mathbb{S}^{d-1}$ be open with smooth boundary. Let the exact cone with smooth cross section $Y$, $W_Y=[0,\infty)_r\times Y\subset\mathbb{R}^{d}$ be equipped with the metric $dr^2+r^2g_{\mathbb{S}^{d-1}}$. We then have
     \begin{equation}
         \rardy^0_{W_Y}= \sqrt{\lambda_1(-\Delta_Y)+\frac{(d-2)^2}{4}},
     \end{equation}
     where $-\Delta_Y$ is the Dirichlet Laplace--Beltrami operator on $Y$ with the round metric and $\lambda_1$ denotes the first (positive) eigenvalue. In particular, if $Y=(0,\alpha)$, then $\rardy^0_{W_Y}=\frac{\pi}{\alpha}$. 
\end{lemma}

\begin{proof}
    1. We show $\rardy^0_{W_Y}\ge \sqrt{\lambda_1(-\Delta_Y)+\frac{(d-2)^2}{4}}$. Let $(r, y)$ be the coordinates on $W_Y$. Let $u\in \dot H^1(W_Y)$. Consider $-\Delta_{W_Y}=\partial_r^*\partial_r-r^{-2}\Delta_Y$. This implies
    \begin{equation}\label{e:al1}
        \|\nabla u\|^2_{L^2(W_Y)}=\langle -\Delta_{W_Y} u, u\rangle_{L^2(W_Y)}=\|\partial_r u\|_{L^2(W_Y)}^2+\langle -r^{-2}\Delta_Y u, u\rangle_{L^2(W_Y)}.  
    \end{equation}
    We claim 
    \begin{equation}
        \|\nabla u\|^2_{L^2(W_Y)}\ge \Big(\lambda_1(-\Delta_Y)+\frac{(d-2)^2}{4}\Big)\|r^{-1}u\|_{L^2(W_Y)}^2.
    \end{equation}
    The second term of \eqref{e:al1} can be bounded from below by
    \begin{multline}
        \langle -r^{-2}\Delta_Y u, u\rangle_{L^2(W_Y)}=\int_0^\infty r^{-2}\langle-\Delta_Y u, u\rangle_{L^2(Y)}\,r^{d-1}dr\ge \lambda_1(-\Delta_Y)\int_0^\infty r^{-2}\|u\|_{L^2(Y)}^2\,r^{d-1}dr\\
        =\lambda_1(-\Delta_Y)\|r^{-1}u\|_{L^2(W_Y)}^2.
    \end{multline}
    The case $d=2$ immediately follows from the non-negativity of the first term of \eqref{e:al1}. For $d>2$, we first analyze the first term of \eqref{e:al1}. We note $\partial_r^*=-r^{-(d-1)}\partial_r r^{d-1}=-\partial_r-\frac{d-1}{r}$ gives a Rellich-type identity:
    \begin{equation}
        \langle \partial_r u, r^{-1}u\rangle_{L^2(W_Y)}=\langle u, \partial_r^*(r^{-1}u)\rangle_{L^2(W_Y)}=-\langle r^{-1}u, \partial_r u\rangle_{L^2(W_Y)}-(d-2)\|r^{-1}u\|_{L^2(W_Y)}^2. 
    \end{equation}
    This implies
    \begin{equation}
        \|r^{-1}u\|_{L^2(W_Y)}^2=-\frac{2}{d-2}\Re\langle \partial_r u, r^{-1}u\rangle_{L^2(W_Y)}\le \frac{2}{d-2}\|\partial_r u\|_{L^2(W_Y)}\|r^{-1}u\|_{L^2(W_Y)}
    \end{equation}
    and a Hardy inequality on $W_Y$:
    \begin{equation}
        \|r^{-1}u\|_{L^2(W_Y)}\le \frac{2}{d-2}\|\partial_r u\|_{L^2(W_Y)}. 
    \end{equation}
    The claim now follows from \eqref{e:al1}. 

    2. We now show that $\rardy^0_{W_Y}\leq \sqrt{\lambda_1(-\Delta_Y)+\frac{(d-2)^2}{4}}$. Let $\chi\in C_c^\infty((\frac12,4);[0,1])$ with $\chi\equiv 1 $ on $[1,2]$ and let $\phi\in \dot H^1(Y)\cap H^2(Y)$ satisfy
    \begin{equation}
        (-\Delta_Y-\lambda_1(-\Delta_Y))\phi(y)=0,\ \|\phi\|_{L^2(Y)}=1.
    \end{equation}
    For $\e>0$, define
    \begin{equation}
        u_\e(r,y)=\beta_\e^{-1}\chi(\e \log r)r^{-\frac{d-2}{2}}\phi(y), \ \beta_\e=\Big(\int_0^\infty r^{-1}\chi^2(\e \log r)\,dr\Big)^{\frac12}\geq \e^{-\frac12}, 
    \end{equation}
    normalized so that $\|r^{-1}u_\e\|_{L^2(W_Y)}\equiv 1$. We compute 
    \begin{equation}
        \partial_r u_\e= \beta_\e^{-1}\Big(\e \chi'(\e \log r)-\frac{d-2}{2}\chi(\e \log r)\Big)r^{-\frac{d}{2}}\phi(y)=-\frac{d-2}{2}r^{-1}u_\e+O_{L^2(W_Y)}(\e).
    \end{equation}
    Now we have
    \begin{equation}
        \|\nabla u_\e\|_{L^2(W_Y)}^2=\lambda_1(-\Delta_Y)\|r^{-1}u_\e\|_{L^2(W_Y)}^2+\|\partial_r u_\e\|_{L^2(W_Y)}^2=\lambda_1(-\Delta_Y)+\frac{(d-2)^2}{4}+O(\e). 
    \end{equation}
    Pushing $\e\rightarrow 0$ saturates the claimed bound. 
\end{proof}

\section{\dtn quasimodes for 2D \acs domains}\label{a:3}
Let $\Omega$ be a smooth \acl domain and $(x_0,x)$ denote Fermi normal coordinates in a neighborhood, 
$$
U=\{ (x_0,x)\,:\, 0\le x_0<c\langle x\rangle\}
$$
of $\partial\Omega$ with $x$ the arc length coordinates on $\partial\Omega$.  Recall that, in these coordinates,
$$
g_+(x_0,x)=dx_0^2+[R(x_0,x)]^2dx^2,\qquad R(x_0,x):=1+\kappa(x)x_0,
$$
where $\kappa(x)$ denotes the curvature of $\partial\Omega$ and 
$$
|\partial_x^\alpha \kappa(x)|\leq C_{\alpha}\langle x\rangle^{-1-|\alpha|}
$$

Then, in these coordinates,
$$
-h^2\Delta+h^2M^2=\big(h^2D_{x_0}^2-iR^{-1}\partial_{x_0}Rh^2D_{x_0}+R^{-2}h^2D_x^2+iR^{-1}\partial_xRh^2D_x+h^2M^2\big),
$$
and, conjugating away the volume form, we have
$$
R^{1/2}(-h^2\Delta+h^2M^2)R^{-1/2}=h^2D_{x_0}^2+(R^{-1/2}hD_xR^{-1/2})^2+h^2M^2-\frac{h^2\kappa^2}{4R^2}.
$$

\begin{lemma}
\label{l:carefulFactor}
    Let $\chi,\chi_+ \in C_c^\infty(\mathbb{R})$ with $0\notin \supp \chi_+$, $\supp \chi\cap \supp (1-\chi_+)=\emptyset$. Then there are $b_0$, $b_1$ such that, with
    $$
    p :=i\sqrt{\xi^2R^{-2}+h^2M^2}\chi_+(\xi)
    +i\frac{h}{2R}\kappa(x)+ h^2(hM^2b_0(x,\xi)+x_0b_1(x_0,x,\xi))\chi_+(\xi),
    $$
    we have
    $$
    R^{1/2}(-h^2\Delta+h^2M^2)R^{-1/2}\chi(hD_x)- (hD_{x_0}-p^w)(hD_{x_0} + p^w)\chi(hD_x)\in h^\infty \Psi^{-\infty,-\infty}(U).
    $$
    Moreover, 
\begin{equation}
\label{e:vanishing}
\begin{gathered}
 b_0\sim \sum_j h^jb_{0,j},\qquad 
 b_{0,j}\in S^{\comp,-1-j}\\
  b_1\sim \sum_j h^jb_{1,j},\qquad 
\partial_{x_0}^\ell b_{1,j}\in S^{\comp,-1-j-\ell}.
 \end{gathered}
\end{equation}
\end{lemma}
\begin{proof}
We proceed as in Section~\ref{s:3}. Observe that, we require
\begin{align*}
&(h^2D_{x_0}^2+(R^{-1/2}hD_xR^{-1/2})^2+h^2M^2-\frac{h^2\kappa^2}{4R^2})\chi(hD_x)\\
&=((hD_{x_0})^2 -(ih\partial_{x_0}p)^w-(p^w)^2)\chi(hD_x)+h^\infty\Psi^{-\infty,-\infty}.
\end{align*}

Thus, our goal is to solve 
\begin{multline*}
((p^w)^2
+(ih\partial_{x_0}p)^w)\chi(hD_x)
=(-(R^{-1/2}hD_xR^{-1/2})^2-h^2M^2+\frac{h^2\kappa^2}{4R^2})\chi(hD_x)+O(h^\infty)_{\Psi^{-\infty,-\infty}}.
\end{multline*}

Since $\supp (1-\chi_+)\cap \supp \chi=\emptyset$, we may ignore $\chi_+$ in all of the below calculations at the cost of an $h^\infty\Psi^{-\infty,-\infty}$ error. We also do not track decay in $\xi$ since every symbol will be compactly supported.

Recall that 
$$
(b_1^wb_2^w)=b_3^w,\qquad b_3=:b_1\sharp b_2\sim \sum_j \frac{h^ji^j}{j!}(A(D))^jb_1(y,\eta)b_2(x,\xi)\Big|_{\substack{y=x\\\eta=\xi}},\qquad A(D):=\frac{1}{2}\sigma(D_y,D_{\eta},D_x,D_\xi).
$$
In particular, if $b_1\sim \sum_jh^jb_{1,j}$ and $b_{2}\sim \sum_jh^jb_{2,j}$
$$
b_1\sharp b_2\sim \sum_jh^j\sum_{\ell=0}^j\sum_{m=0}^{\ell} \frac{i^m}{m!}(A(D))^{m}b_{1,j-\ell}(y,\eta)b_{2,\ell-m}(x,\xi)\Big|_{\substack{y=x\\\eta=\xi}},
$$
Therefore, if $b_{1}$ and $b_{2}$ satisfy~\eqref{e:vanishing}, then so do $b_1\sharp b_2$ and $h\partial_{x_0} b_1$. (Indeed, each power of $h$ in the asymptotic series comes exactly one $(x,y)$- derivative). 
Hence, with $p_0=i\sqrt{\xi^2R^{-2}+h^2M^2}+i\frac{h}{2R}\kappa$, 
\begin{align*}
(p^w)^2+(ih\partial_{x_0}p)^w&=(p_0\sharp p_0+h^2(p_0\sharp b+b\sharp p_0)+h^4b\sharp b)^w  +(ih\partial_{x_0} p_0+ih^3\partial_{x_0} b)^w.
\end{align*}
Now, 
\begin{multline}
p_0\sharp p_0
\sim -(\xi^2R^{-2}+h^2M^2)-\sqrt{\xi^2R^{-2}+h^2M^2}hR^{-1}\kappa -\tfrac{h^2\kappa^2}{4R^2}\\
-\sum_{j\geq 1}\frac{h^{2j}(-1)^j}{(2j)!2^{2j}}(D_{\eta}D_y-D_xD_\xi)^{2j}(\sqrt{\xi^2R(x)^{-2}+h^2M^2}+\tfrac{1}{2}hR(x)^{-1}\kappa(x)) \\\times(\sqrt{\eta^2R(y)^{-2}+h^2M^2}+\tfrac{1}{2}h
R(y)^{-1}
\kappa(y))|_{y=x,\eta=\xi}\\
\sim -(\xi^2R^{-2}+h^2M^2)-\sqrt{\xi^2R^{-2}+h^2M^2}hR^{-1}\kappa -\tfrac{h^2\kappa^2}{4R^2}+h^2\sum_{j}x_0c_{1,j}(x_0,x,\xi)h^j+ h^3M^2\sum_{j}c_{0,j}(x,\xi)h^j,
\end{multline}
with 
\begin{gather*}
|\partial_x^m c_{i,j}(x,\xi)|\leq C_m \langle x\rangle^{-2-j-m},
\end{gather*}
and 
$$
ih\partial_{x_0} p_0=\frac{h\xi^2(R^{-3}\kappa)}{\sqrt{\xi^2R^{-2}+h^2M^2}}+\frac{h^2\kappa^2}{2R^2}=h\sqrt{\xi^2R^{-2}+h^2M^2}R^{-1}\kappa-\frac{h^3M^2R^{-1}\kappa}{\sqrt{\xi^2R^{-2}+h^2M^2}}+\frac{h^2\kappa^2}{2R^2}.
$$

So our goal is to solve
\begin{multline*}
p_0\sharp (hM^2b_0+x_0b_1)+(hM^2b_0+x_0b_1)\sharp p_0+h^2(hM^2b_0+x_0b_1)\sharp (hM^2b_0+x_0b_1) +ih\partial_{x_0}(x_0b_1)\\
\sim -\sum_j x_0c_{1,j}h^j-hM^2\sum_j c_{0,j}+\frac{hM^2R^{-1}\kappa}{\sqrt{\xi^2R^{-2}+h^2M^2}}
\end{multline*}
Setting 
\begin{align*}
\tilde{b}_0&:=\frac{1}{2p_0(0,x,\xi)}\Big(- \sum_{j\geq 0}c_{0,j}h^j+\frac{R^{-1}\kappa}{\sqrt{\xi^2R^{-2}+h^2M^2}}\Big),\\
\tilde{b}_1&:=\frac{1}{2p_0}(- \sum_{j\geq 0}c_{1,j}h^j+\frac{hM^2}{x_0}\Big(1-
\frac{p_0(x_0,x,\xi)}{p_0(0,x,\xi)}\Big)\Big(\sum_{j\geq 0}c_{0,j}h^j+\frac{R^{-1}\kappa}{\sqrt{\xi^2R^{-2}+h^2M^2}}\Big),\\
\end{align*}
we obtain that $\tilde{b}_i$ satisfy~\eqref{e:vanishing}. Suppose by induction that $\tilde{b}_{K-1,i}$ satisfies~\eqref{e:vanishing} and, with $\tilde{b}_{K-1}=x_0\tilde{b}_{K-1,1}+hM^2\tilde{b}_{K-1,0}$, 
\begin{multline*}
p_0\sharp \tilde{b}_{K-1}+\tilde{b}_{K-1}\sharp p_0+h^2\tilde{b}_{K-1}\sharp \tilde{b}_{K-1} +ih\partial_{x_0}\tilde{b}_{K-1}+\sum_j x_0c_{1,j}h^j+hM^2\sum_j c_{0,j}- \frac{hM^2R^{-1}\kappa}{\sqrt{\xi^2R^{-2}+h^2M^2}}\\
=h^{K+1}M^2e_{K,0}+x_0h^{K}e_{K,1}.
\end{multline*}
with $e_{K,i}\in S^{\comp,-1-K}$.
Define 
\begin{align*}
\tilde{b}_{K,0}&:=\tilde{b}_{K-1,0}-\frac{h^K}{2p_0(0,x,\xi)}e_{K,0},\\
\tilde{b}_{K,1}&:=\tilde{b}_{K-1,1}-\frac{h^K}{2p_0}\Big(e_{K,1}+\frac{hM^2}{x_0}\Big(1-
\frac{p_0(x_0,x,\xi)}{p_0(0,x,\xi)}\Big)e_{K,0}\Big).
\end{align*}
Then, since $x_0^{-1}(p_0(0,x,\xi)-p_0(x_0,x,\xi))=\kappa(x)e(\kappa(x)x_0,\kappa(x),\xi)$ for some smooth function $e$, $\tilde{b}_{K,i}$ satisfies~\eqref{e:vanishing} and 
\begin{align*}
&p_0\sharp \tilde{b}_{K}+\tilde{b}_{K}\sharp p_0+h^2\tilde{b}_{K}\sharp \tilde{b}_{K} +ih\partial_{x_0}\tilde{b}_{K}+\sum_j x_0c_{1,j}h^j+hM^2\sum_j c_{0,j} -\frac{hM^2R^{-1}\kappa}{\sqrt{\xi^2R^{-2}+h^2M^2}}\\
&=h^{K+2}M^2e_{K+1,0}+x_0 h^{K+1}e_{K+1,1},
\end{align*}
with $e_{K+1,i}$ satisfying~\eqref{e:vanishing}
and hence the proof of~\eqref{e:vanishing} is complete after using the Borel summation lemma to set $b_\ell\sim\tilde{b}_{\ell,0}+ \sum_j \tilde{b}_{\ell,j}-\tilde{b}_{\ell,j-1}.$

\end{proof}

\begin{lemma}
\label{l:1dnormalform}
    Let $\chi \in C_c^\infty(\mathbb{R}\setminus \{0\})$ with $0\notin \supp \chi$. Then, with $p$ as in Lemma~\ref{l:carefulFactor}
    $$
    h\Lambda_{M}^\Omega \chi(hD_x)=-i(p^w-\tfrac{ih}{2}\kappa(x))\chi(hD_x)+E_M^\Omega,
    $$
    where
    $$
\|E_M^\Omega\|_{H_h^{-N,-N}\to H_h^{N,\rardy^{\Omega}+\frac{1}{2}-0}}\leq C_Nh^N,
$$
and
for any $\psi\in C_c^\infty(\mathbb{R})$ with $0\notin \supp (1-\psi)$,
$$
\|(1-\psi(D_x))E_M^\Omega\|_{H_h^{-N,-N}\to H_h^{N,N}}\leq C_Nh^N.
$$
\end{lemma}
\begin{proof}
In this proof, we denote by $\Psi_T^{-\infty,-\infty}$ the tangential residual class and $\Psi^{-\infty,-\infty}$, the full residual class.
    Let $\chi_0=\chi$ and $\chi_i\in C_c^\infty(\mathbb{R}\setminus\{0\})$, $i=1,2$  with $\supp (1-\chi_i)\cap \supp \chi_{i-1}=\emptyset$. Let  $U(x_0)$ satisfy, 
    \begin{equation} \label{e:usolve}
    (hD_{x_0}-p^w-i(1-\chi_+(hD))\langle hD\rangle)U(x_0)=0,\qquad U(0)=I.
    \end{equation}
    We first claim that 
    \begin{equation}
    \label{e:offDiagonal}
    (1-\chi_1(hD_x))U(x_0)\chi(hD_x)=O(e^{-cx_0/h}h^\infty)_{\Psi_T^{-\infty,-\infty}}.
    \end{equation}

    First, notice that 
    \begin{align*}
    h\partial_{x_0}\| U(x_0)f\|_{L^2}^2
    &= -2\Im \langle hD_{x_0}U(x_0)f,U(x_0)f\rangle\\
    &=-2\Im \langle (p^w+i(1-\chi_+(hD))\langle hD\rangle)U(x_0)f,U(x_0)f\rangle \\
    &\leq -c\|U(x_0)f\|_{L^2}^2.
    \end{align*}
    Hence, 
    $$
    \|U(x_0)f\|_{L^2}\leq Ce^{-cx_0/h}\|f\|_{L^2}.
    $$
    Furthermore, $a=a(x,\xi)\in S^{m,n}$, and $b\in S^{0,0}$ with $\supp(1-b)\cap \supp a=\emptyset$,  real valued,
    \begin{align*}
    h\partial_{x_0}\|a^w U(x_0)f\|_{L^2}^2
    &= -2\Im \langle a^whD_{x_0}U(x_0)f,a^wU(x_0)f\rangle\\
    &=-2\Im \langle a^w(p^w+i(1-\chi_+(hD))\langle hD\rangle)U(x_0)f,a^wU(x_0)f\rangle \\
    &=-2\Im \langle(p^w+i(1-\chi_+(hD))\langle hD\rangle) a^wU(x_0)f,a^wU(x_0)f\rangle\\
    &\qquad -2\Im \langle a^w[a^w,(p^w+i(1-\chi_+(hD))\langle hD\rangle) ]U(x_0)f,U(x_0)f\rangle\\
    &\leq -c\|a^wU(x_0)f\|_{L^2}^2+Ch^2\|b^wU(x_0)f\|^2_{H_h^{m-1,n-1}}+C_Nh^N\|U(x_0)f\|_{H_h^{-N,-N}}.
    \end{align*}
    We first use this estimate by letting $a\in\Psi^{m,n}$ elliptic and iterating to obtain 
    $$
    \|U(x_0)f\|_{H_{h}^{m,n}}\leq Ce^{-cx_0/h}\|f\|_{H_h^{m,n}}. 
    $$

    Next, for any $a\in S^{m,n}$ and $b\in S^{0,0}$ with $\supp(1-b)\cap \supp a=\emptyset$, we obtain
    \begin{align*}
    \|e^{\frac{cx_0}{2h}}a^wU(x_0)f\|_{L^2}^2&=\|a^wf\|_{L^2}^2+\int_0^{x_0}\partial_{s}e^{\frac{cs}{2h}}\|a^wU(s)f\|_{L^2}ds\\
    &\leq\|a^wf\|_{L^2}^2+ C\int_0^{x_0}h^2e^{\frac{cs}{2h}}\|b^wU(s)f\|_{H_h^{m-1,n-1}}^2ds+Ch^N\int_0^{x_0}e^{\frac{cs}{2h}}\|U(s)f\|_{H_h^{-N,-N}}ds\\
    &\leq \|a^wf\|_{L^2}^2+C\int_0^{x_0}h^2e^{\frac{cs}{2h}}\|b^wU(s)f\|_{H_h^{m-1,n-1}}^2ds+Ch^N\|f\|_{L^2}.
\end{align*}
Another iteration shows that 
$$
\|a^wU(x_0)f\|_{L^2}\leq Ce^{-cx_0/h}\big(\|a^wf\|_{L^2}^2+ch^N\|f\|_{L^2}\big),
$$
which implies~\eqref{e:offDiagonal}.

Now letting $\chi_0\in C_c^\infty(\mathbb{R})$ with $\supp (1-\chi_0)\cap \supp \chi=\emptyset$, $\supp(1-\chi_+)\cap \supp \chi_0=\emptyset$, and $\psi\in C_c^\infty(\mathbb{R})$ with $\psi\equiv 1$ near $0$, using~\eqref{e:usolve} to convert derivatives in $x_0$ into derivatives in $x$,  we have
\begin{align*}
(-h^2\Delta+h^2M^2)R^{-1/2}\psi(x_0/\langle x\rangle)U\chi(hD_x)&=(-h^2\Delta+h^2M^2)R^{-1/2}\chi_+(hD_x)\psi(x_0/\langle x\rangle)U\chi(hD_x)+O(h^\infty)_{\Psi^{-\infty,-\infty}}\\
&=O(h^\infty)_{\Psi^{-\infty,-\infty}}.
\end{align*}

Hence, solving away the error, we have 
$$
h\Lambda_{M}^{\Omega}\chi(hD_x)=-h\partial_{x_0} R^{-1/2}U(x_0)|_{x_0=0}\chi(hD_x)+E_M^{\Omega}, 
$$
where, 
$$
\|E_M^\Omega\|_{H_h^{-N,-N}\to H_h^{N,\rardy^{\Omega}+\frac{1}{2}-0}}\leq C_Nh^N
$$ and
for any $\psi\in C_c^\infty(\mathbb{R})$ with $0\notin \supp (1-\psi)$
$$
\|(1-\psi(D_x))E_M^\Omega\|_{H_h^{-N,-N}\to H_h^{N,N}}\leq C_Nh^N.
$$
This implies the lemma. 
\end{proof}

\begin{lemma}\label{l:2dquasimodeappendix}
Let $\psi_0\in C_c^\infty(\mathbb{R})$ with $0\notin \supp (1-\psi_0)$.
Then there is a solution, $u$ to 
    $$
    (h\Lambda^{\Omega}_M-\sqrt{\omega^2+M^2})e_\pm(s; h, \omega, M)=g_1+M^2g_2+\psi_0(D_x)g_3,
    $$
    such that $e_\pm(s; h, \omega, M)=e^{\pm i\omega x/h}(1+h^2M^2e_{\pm}(s))$ with
    $
    e_{\pm}\sim \sum_{j} h^je_{j,\pm},
    $ with $e_{j,\pm}(0)=0$ for all $j$ and
    \begin{equation} \label{e:properties}
\begin{gathered}
|\partial_x^{\ell }e_{j,\pm}(x)|\leq C_\ell \langle x\rangle^{-1-\ell-j},\qquad 0\leq j\leq 2,\, 1\leq \ell \\
|\partial_x^{\ell }e_{j,\pm}(x)|\leq C_\ell( \langle x\rangle^{-1-\ell-j}+M^2\langle x\rangle^{-1-\ell}),\qquad 3\leq j,\, 1\leq \ell.
\end{gathered}
\end{equation}
    and for all $N$ 
    $$
    \|g_1\|_{H_{h}^{N,N}}+\|g_2\|_{H_h^{N,\frac{3}{2}-0}}+\|g_2\|_{H_h^{N,\rardy^{\Omega}+\frac{1}{2}-0}}\leq C_Nh^N.
    $$
\end{lemma}

\begin{proof}
We consider the case of $u=e^{ix\sqrt{1-h^2M^2}/h}(1+h^2M^2e_+(x))$, the other being similar and drop the $\pm$ from the notation to simplify the writing.

We start by constructing $f\in S$ such that $f(0)=1$, and 
\begin{equation} \label{e:fSolve}
\begin{gathered}
 e^{-ix\sqrt{1-h^2M^2}/h}(p^w-i\tfrac{h}{2}\kappa-1)e^{ix\sqrt{1-h^2M^2}{h}}f=O(h^\infty)_{\mc{S}}+O(h^\infty M^2)_{S^{-2}},\\
 f\sim 1+h^2M^2\sum_{j\geq0} h^je_j(x),
 \end{gathered}
\end{equation}
with $e_j$ satisfying~\eqref{e:properties}.

Indeed, observe that, by Lemma~\ref{l:1dnormalform}
\begin{align*}
&e^{-is\sqrt{1-h^2M^2}/h}(-i( p^w-i\tfrac{h}{2}\kappa)-1)e^{ix\sqrt{1-h^2M^2}{h}}f\\
&=(\sqrt{(hD_x)^2+h^2M^2-1}+h^3M^2b_0^w)e^{ix\sqrt{1-h^2M^2}/h}f+O(h^\infty)_{\mathcal{S}}\\
&\sim \sum_j\frac{h^j}{j!i^j}(D_wD_\xi)^j \Big[\Big(\sqrt{(\xi+\sqrt{1-h^2M^2})^2+h^2M^2}+h^3M^2b_0(x-\tfrac{w}{2},\xi+\sqrt{1-h^2M^2})-1\Big)f(x-w))\Big]|_{w=\xi=0},
\end{align*}
and hence, differentiating in $s$ we have
\begin{multline}
\sum_j\frac{h^j}{j!i^j}D_w^\ell(D_wD_\xi)^j \Big[\Big(\sqrt{(\xi+\sqrt{1-h^2M^2})^2+h^2M^2}+h^3M^2b_0(x-\tfrac{w}{2},\xi+\sqrt{1-h^2M^2})-1\Big)(f-1)(x-w))\Big]|_{w=\xi=0}\\
\sim -\sum_j \frac{h^{j+3}M^2}{j!i^j}D_w^\ell(D_wD_\xi)^jb_0(x-\tfrac{w}{2},\xi+\sqrt{1-h^2M^2})|_{w=\xi=0}
\label{e:asymptotice}
\end{multline}
Now, reorganizing the sum, we obtain, with $e:=h^{-2}M^{-2}(f-1)\sim \sum_j h^je_j$ and $e_m(0)=0$
\begin{align*}
&-\sqrt{1-h^2M^2}D^{\ell+1}_xe_m(x)\\
&=-\sum_{j=0}^{m-1}\frac{1}{(j+2)!i^{j+2}}D_w^\ell(D_wD_\xi)^{j+2} \Big[\Big(\sqrt{(\xi+\sqrt{1-h^2M^2})^2+h^2M^2}-1\Big)e_{m-1-j}(x-w)|_{w=\xi=0}\\
&-\sum_{j=0}^{m-3}\frac{M^2}{j!i^j}D_w^\ell(D_wD_\xi)^j \Big[\Big(b_0(x-\tfrac{w}{2},\xi+\sqrt{1-h^2M^2})-1\Big)e_{m-j-3}(x-w))|_{w=\xi=0}\\
&-  \frac{1}{m!i^m}D_w^\ell(D_wD_\xi)^mb_0(x-\tfrac{w}{2},\xi+\sqrt{1-h^2M^2})|_{w=\xi=0}.
\end{align*}
We now proceed by induction. Observe that, since $b_0\in S^{-2}$ 
$$
|\partial_x^{\alpha}e_0|\leq C_{\alpha}\langle x\rangle^{-1-|\alpha|}, \, |\alpha|\geq 1.
$$

Now, it is easy to see that for $m\leq 2$, 
$$
 |\partial_x^\alpha e_m|\leq C_\alpha \langle x\rangle^{-1-m-|\alpha|},\quad |\alpha|\geq 1
$$
Suppose by induction that for $0\leq m\leq N$,  
$$
|\partial_x^\alpha e_m|\leq C_\alpha (\langle x\rangle^{-1-m-|\alpha|}+M^2\langle x\rangle^{-1-|\alpha|}),\quad  |\alpha|\geq 1,
$$
Then, using~\eqref{e:asymptotice}, we see that this holds for $0\leq m\leq N+1$ and hence for all $m$ by induction.

Using the Borel summation lemma, we have thus proved~\eqref{e:fSolve} and~\eqref{e:properties}.

Now, to conclude the proof, observe that for $\chi\in C_c^\infty(\mathbb{R})$ with $0\notin \supp \chi$ and $1\notin \supp (1-\chi)$, 
$
(1-\chi(hD_x))f=O(h^\infty)_{\mc{S}},
$ and hence by Lemma~\ref{l:1dnormalform}
\begin{align*}
(h\Lambda_M^\Omega-1)f&=(-i(p^w-i\tfrac{ih\kappa}{2})-1)f+O(h^\infty)_{\mc{S}}+E_M^\Omega f\\
&=g=O(h^\infty)_{\mc{S}}+e^{ix\sqrt{1-M^2h^2}/h}O(h^\infty M^2)_{S^{-2}}+\psi_0(D_x)O(h^\infty)_{H_h^{N,\rardy^{\Omega}+\frac{1}{2}-0}}.
\end{align*}
\end{proof}

\section{Explicit scattering matrices on 2D straight wedges}\label{a:4}
To determine the scattering matrix on an infinite wedge, we use parts of \cite[Appendix A]{LePaPoSh:22b} which in turn relies heavily on~\cite{Pe:50}. Given incoming data and $0<\alpha<2\pi$, we have shown that there is a unique solution $u\in Z_{-\frac{\pi}{\alpha}-\frac{1}{2}+,-\frac{1}{2}-,-\frac{1}{2}-}^{\infty,\infty,-1+}$ to 
\begin{equation}
    \label{e:planeWaveEquation}
(\Lambda_{\mathcal{W}_\alpha}-1)u=0,
\end{equation}
with that incoming data. Thus, since there is a two dimensional family of incoming data, it is enough to find two linearly independent
solutions to~\eqref{e:planeWaveEquation} and compute their asymptotic data.

\subsection{Solutions of the sloshing problem}
To find two linearly independent solutions of~\eqref{e:planeWaveEquation}, we first construct solutions $v\in L^\infty(\mathcal{W}_{\alpha/2})$ of the two problems
\begin{equation}
    \label{e:stekDir}
    \begin{cases}
        -\Delta v=0&\text{in }\mathcal{W}_{\alpha/2}=\{ z\in\mathbb{Z}\,:\, -\alpha/2<\arg z<0\}\\
        \partial_y v(s)=v(s)&\text{for }s>0\\
        v(s e^{-i\alpha/2})=0&\text{for }s>0.
    \end{cases}
\end{equation}
and
\begin{equation}
    \label{e:stekNeu}
    \begin{cases}
        -\Delta v=0&\text{in }\mathcal{W}_{\alpha/2}=\{ z\in\mathbb{Z}\,:\, -\alpha/2<\arg z<0\}\\
        \partial_y v(s)=v(s)&\text{for }s>0\\
        \partial_{\nu_{\alpha/2}} v(s e^{-i\alpha/2})=0&\text{for }s>0,
    \end{cases}
\end{equation}
where $\nu_{\alpha/2}$ is the normal vector to $\arg z=-\alpha/2$.

\begin{remark}The arguments in~\cite[Appendix A]{LePaPoSh:22b} apply directly for $0<\alpha<\pi$, but not for $\pi\leq \alpha<2\pi$. Therefore, we recall the parts of that appendix which we require here.
\end{remark}

Let $\mathcal{C}$ denote the logarithmic cover of $\mathbb{C}\setminus \{0\}$.  Define (at first $I_\alpha$ is defined for $\zeta>0$)
\begin{equation}
\label{e:IDef}
g_{\alpha}(\zeta):=\exp \Big(-I_\alpha(\zeta e^{-i\alpha/2})+\log \Big(\frac{\zeta+i}{\zeta}\Big)\Big),\qquad I_{\alpha}(\zeta):=\frac{1}{\pi}\int_0^{\infty}\log (1+w^{-2\pi /\alpha})\frac{\zeta}{w^2+\zeta^2}dw.
\end{equation}

 We show below (following~\cite{Pe:50}) that $I_\alpha$ has an analytic continuation to the logarithmic cover of 
$$\mathcal{C}\setminus \big( \bigcup_{j=0}^\infty e^{-i(\pi +(2j+1)\alpha)/2}[0,1]\cup\bigcup_{j=0}^{\infty} e^{i(\pi +(2j+1)\alpha)/2}[0,1]\big).$$
\begin{remark}
    In fact, it is not hard to see that it has a meromorphic continuation to a larger Riemann surface where each point $e^{\pm i(\pi +(2j+1)\alpha)/2}$ is a branch point.
\end{remark}
Following~\cite{LePaPoSh:22a,Pe:50} we will show that $v_d:=\Re f_d(z)$ and $v_n=\Re f_n(z)$ satisfy respectively~\eqref{e:stekDir} and~\eqref{e:stekNeu} with
$$
f_d(z):=i\int_P \frac{g_{\alpha}(\zeta)}{\zeta+i}\zeta^{-\mu}e^{z\zeta}d\zeta,\qquad f_n(z):=i\int_P \frac{g_{\alpha}(\zeta)}{\zeta+i}e^{z\zeta}d\zeta
$$
where $\mu:=\frac{\pi}{\alpha}$, $P$ is a contour consisting of the union of a circle of radius $2$ traversed counter clockwise and the two half lines  $$
\begin{aligned}
\Gamma_{out}&:=\{s\in(2,\infty)\,:\, se^{i(\pi +\alpha/4)}\},&
\Gamma_{in}&:=\{s\in(\infty,2)\,:\, se^{i(-\pi +\alpha/4)}\}.
\end{aligned}
$$

\subsection{Meromorphic continuation of $I_{\alpha}$ to the logarithmic cover \label{subsec:MeromorphicCont}}

\begin{lemma}
The function $I_{\alpha}(\zeta)$ defined by~\eqref{e:IDef} analytically continues to 
$$
\mathcal{C}\setminus \big( \bigcup_{j=0}^\infty e^{-i(\pi +(2j+1)\alpha)/2}[0,1]\cup\bigcup_{j=0}^{\infty} e^{i(\pi +(2j+1)\alpha)/2}[0,1]\big).
$$
Moreover, the function
$$
\tilde{g}(\zeta):=\exp(I_{\alpha}(\zeta))
$$
is meromorphic on $\mathcal{C}$ and its poles are given by $\{\lambda_k\}_{k=-\infty}^\infty$ with
\begin{equation}
\label{e:whereAreThePoles}
 \lambda_k =e^{-i(\pi/2-k\alpha)},\, k\leq 0,\qquad \lambda_k = e^{i(\pi/2 +k\alpha)},\,k\geq 1.
\end{equation}
\end{lemma}
\begin{proof}
Using that $\frac{\zeta}{v^2+\zeta^2}=\frac{1}{2i}(\frac{1}{v-i\zeta}-\frac{1}{v+i\zeta})$ and changing coordinates with $v=\pm iu$, it is easy to check that,
\begin{align}
I_{\alpha}(\zeta)= & \frac{1}{2\pi i}\int_{L}\log(1+  (-iu)^{-2\mu})\frac{1}{u-\zeta}\,du +\frac{1}{2\pi i}\int_{M}\log(1+(iu)^{-2\mu})\frac{1}{u-\zeta}\,du,\label{eq:I_alpha_Split}
\end{align}
where $L:=i(\infty,0)$ and $M:=-i(0,\infty)$ . 

We notice that (\ref{eq:I_alpha_Split}) represents $\zeta\mapsto I_{\alpha}(\zeta)$
in the half plane $-\pi/2<\arg\zeta<\pi/2$ since it defines an analytic
function on $-\pi/2<\arg\zeta<\pi/2$ which coincides with $I_{\alpha}$
when $\zeta>0$.

Fix $\theta_1<\theta_2\in \mathbb{R}$. Then, define
$$C_{R}:= \{u;\,\arg u\in(\theta_1,\theta_2),\,|u|=R\}$$
and
$$
\Gamma_{out,R}^{\theta}:=e^{i\theta}(R,\infty),\qquad \Gamma_{in,R}^{\theta}:= e^{i\theta}(\infty,R).
$$
Notice that for any $\theta_1,\theta_2$, 
\[
\lim_{R\rightarrow\infty}\int_{C_{R}}\log(1+(\pm iu)^{-2\mu})\frac{1}{u-\zeta}\,du=0
\]
Also, $\log(1+(-iu)^{-2\mu})$ is analytic on the logarithmic cover of 
$$
\mathbb{C}\setminus\Big( \{0\}\cup \bigcup_{k\in\mathbb{Z}}\{e^{i(\pi+(2k+1)\alpha)/2}\}\Big),
$$
since 
\[
(-iu)^{-2\mu}=(-iu)^{-2\mu}=-1,
\]
is equivalent to
\[
-2\mu(\log(u)-i\pi/2)=-i(\pi+2k\pi),\quad k\in\mathbb{Z}.
\]
Hence, since $u\mapsto\frac{\log(1+(-iu)^{-2\mu})}{u-\zeta}$ has
no singularity in the sector $\pi/2<\arg\zeta<\pi/2+\alpha/2$, 
\[
-\int_{\Gamma_{out,0}^{(\pi+\alpha)/2}}\log(1+(-iu)^{-2\mu})\frac{1}{u-\zeta}\,du=\int_{L}\log(1+(-iu)^{-2\mu})\frac{1}{u-\zeta}\,du.
\] Therefore,  for $\arg\zeta\in(-\pi/2,\pi/2+\alpha/2)$,
\begin{align*}
\zeta\mapsto & \frac{i}{2\pi }\int_{\Gamma_{out,0}^{(\pi+\alpha)/2}}\log(1+(-iu)^{-2\mu})\frac{1}{u-\zeta}\,du +\frac{1}{2\pi i}\int_{M}\log(1+(iu)^{-2\mu})\frac{1}{u-\zeta}\,du,
\end{align*}
is an analytic function which coincides with $\zeta\mapsto I_{\alpha}(\zeta)$
when $|\arg\zeta|<\pi/2$ and, hence, it is the holomorphic continuation
of $\zeta\mapsto I_{\alpha}(\zeta)$ to $\arg\zeta\in(-\pi/2,\pi/2+\alpha/2)$.
In a similar way we can shift the path $M$ into $M_{1}:=\Gamma_{in,0}^{-(\pi+\alpha)/2}$
in \eqref{eq:I_alpha_Split} since $u\mapsto\log(1+(iu)^{-2\mu})$
has singularities at $0$ and $e^{i(-\pi+(2k+1)\alpha)/2}$.

To continue $I_\alpha$ further, we repeat the
process of shifting, this time taking care of where the singularities of the
integrand are and computing the relevant residues. For
instance, if $-\pi/2-3\alpha/2<\theta<-\pi/2-\alpha/2$,
and if $-\pi/2-\alpha/2<\arg\zeta<\theta+2\pi$ then
\begin{align*}
  \frac{1}{2\pi i}\int_{\Gamma_{in,0}^{-(\pi+\alpha)/2}}\log(1+(iu)^{-2\mu})\frac{du}{u-\zeta}
= & \frac{1}{2\pi i}\int_{\Gamma_{in,0}^{\theta}}\log(1+(iu)^{-2\mu})\frac{du}{u-\zeta} +\frac{1}{2\pi i}\int_{\Gamma}\log(1+(iu)^{-2\mu})\frac{du}{u-\zeta},
\end{align*}
with $\Gamma$ being the loop path around the branch cut from $e^{-i(\pi/2+\alpha)}$
to $0$ (that is, $\Gamma:=e^{-i(\pi/2+\alpha/2-0)}(0,1)\cup(-e^{-i(\pi/2+\alpha/2+0)})(-1,0)$)
and one computes
\begin{align*}
 \frac{1}{2\pi i}\int_{\Gamma}\log(1+(iu)^{-2\mu})\frac{du}{u-\zeta}
=&\log\Big(\frac{\zeta-e^{-i(\pi/2+\alpha/2)}}{\zeta}\Big).
\end{align*}
 Hence, for $\arg\zeta\in(-\pi/2-3\alpha/2,\pi/2+\alpha/2)$, we have
 \begin{align*}
I_{\alpha}(\zeta)= & \frac{i}{2\pi }\int_{\Gamma_{out,0}^{\pi/2}}\log(1+(-iu)^{-2\mu})\frac{du}{u-\zeta}\\
 & +\frac{1}{2\pi i}\int_{\Gamma_{in,0}^{-(\pi+3\alpha)/2}}\log(1+(iu)^{-2\mu})\frac{du}{u-\zeta}+\log\Big(\frac{\zeta-e^{-i(\pi/2+\alpha/2)}}{\zeta}\Big),\\
 & 
\end{align*}
Repeating this argument we obtain for $\arg\zeta \in (-3\pi/2,\pi/2)$,
\begin{equation}
\begin{aligned}
I_{\alpha}(\zeta)= & \frac{i}{2\pi }\int_{\Gamma_{out,0}^{\pi/2}}\log(1+(-iu)^{-2\mu})\frac{du}{u-\zeta} \\
 & +\frac{i}{2\pi }\int_{\Gamma_{in,0}^{-3\pi/2}}\log(1+(iu)^{-2\mu})\frac{du}{u-\zeta}+\sum_{j=0}^{\lceil\pi/\alpha-\frac{3}{2}\rceil}\log\Big(\frac{\zeta-e^{-i(\pi/2+(2j+1)\alpha/2)}}{\zeta}\Big), \label{eq:I_alpha_AllPlaneExt1}
\end{aligned}
\end{equation}
To extend $I_\alpha$ further, we can now repeat the argument, rotating $\Gamma_{in}$ and $\Gamma_{out}$ in the logarithmic cover of $\mathbb{C}\setminus \{0\}$. 

We will also use the following formula for $\arg\zeta \in (-\pi/2,3\pi/2)$,
\begin{equation}
\begin{aligned}
I_{\alpha}(\zeta)= & \frac{i}{2\pi }\int_{\Gamma_{out,0}^{3\pi/2}}\log(1+(-iu)^{-2\mu})\frac{du}{u-\zeta} \\
 & +\frac{i}{2\pi }\int_{\Gamma_{in,0}^{-\pi/2}}\log(1+(iu)^{-2\mu})\frac{du}{u-\zeta}+\sum_{j=0}^{\lceil\pi/\alpha-\frac{3}{2}\rceil}\log\Big(\frac{\zeta-e^{i(\pi/2+(2j+1)\alpha/2)}}{\zeta}\Big), \label{eq:I_alpha_AllPlaneExt2}
\end{aligned}
\end{equation}
The process of rotating the contour can be continued to obtain~\eqref{e:whereAreThePoles}.
\end{proof}

It will be useful below to also record for $-\pi/2-\alpha/2<\arg \zeta <\pi/2-\alpha/2$
\begin{equation}
\begin{aligned}
I_{\alpha}(\zeta e^{-i\alpha/2})= & \frac{i}{2\pi }\int_{\Gamma_{out,0}^{\pi/2-\alpha}}\log(1+(-iu)^{-2\mu})\frac{du}{u-\zeta e^{-i\alpha/2}} \\
 & +\frac{i}{2\pi }\int_{\Gamma_{in,0}^{-\pi/2-\alpha}}\log(1+(iu)^{-2\mu})\frac{du}{u-\zeta e^{-i\alpha/2}}+\log \Big(\frac{\zeta +i}{\zeta}\Big)-\log \Big(\frac{\zeta-i}{\zeta}\Big), \label{eq:I_alpha_AllPlaneExt3}
\end{aligned}
\end{equation}

\begin{lemma}
\label{l:Ismall}
    For all $M>0$, there is $C>0$ such that
    $$
    |I_\alpha(\zeta)|\leq \frac{C}{|\zeta|},\qquad -M<\arg \zeta<M, |\zeta|>2.
    $$
\end{lemma}
\begin{proof}
    The logarithmic terms appearing, for example in~\eqref{eq:I_alpha_AllPlaneExt1} are bounded by $C/|\zeta|$ for $|\zeta|>1+\frac{1}{M}$. Therefore, we only need to consider the integrals. Notice that for $|\arg \zeta-\theta|>\e$, using that $2\mu>1$, 
    \begin{align*}
    &\Big|\int _{\Gamma^{\theta}_{out}}\log (1+(-iu)^{-2\mu})\frac{du}{u-\zeta}\Big|\\
    &\leq C\int_1^\infty (\log |1-t^{-1}|))\frac{dt}{|e^{i\theta}t-\zeta|}+C\int_0^1(1+|\log |1-t^{-1}|)\frac{dt}{|e^{i\theta} t-\zeta|}\leq \frac{C_\e}{|\zeta|}.
    \end{align*}
    A similar argument for the integral over $\Gamma^{\theta}_{in}$ proves the claim.
\end{proof}
\begin{lemma}\label{l:INear0}
    There is $C>0$ such that 
    $$
    |I_\alpha(\zeta)+\mu \log \zeta|<C,\quad -3\pi/2<\arg \zeta <3\pi/2,\, |\zeta|<\frac{1}{2}.
    $$
\end{lemma}
\begin{proof}
We consider $-\frac{3\pi}{2}<\arg \zeta<\frac{\pi}{2}$. The other case being similar. Using~\eqref{eq:I_alpha_AllPlaneExt1} 
    it is enough to study the function 
    $$
    \tilde{I}_{\theta}(\zeta):=\frac{i}{2\pi }\int_{\Gamma_{out,0}^{\pi/2}}\log(1+(-iu)^{-2\mu})\frac{du}{u-\zeta}  +\frac{i}{2\pi }\int_{\Gamma_{in,0}^{-3\pi/2}}\log(1+(iu)^{-2\mu})\frac{du}{u-\zeta}.
    $$
    For this, we argue as in Lemma~\ref{l:Ismall} to see that the integrals along $\Gamma_{in,3/4}^{-3\pi/2}$ and $\Gamma_{out,3/4}^{\pi/2}$ are bounded for $|\zeta|<\frac{1}{2}$. Therefore, it is enough to consider the remaining integral which is given by 
    \begin{align*}
    &\frac{i}{2\pi}\int_0^{\frac{3}{4}}\Big( \frac{i}{it-\zeta}\log (1+t^{-2\mu})-\frac{i}{it-\zeta}\log (1+( e^{-i\pi}t)^{-2\mu})\Big)dt\\
    &=\frac{i}{2\pi}\int_0^{\frac{3}{4}}\Big( \log (1+t^{-2\mu})-\log (1+e^{2i\mu \pi}t^{-2\mu})\Big)\frac{i}{it-\zeta}dt\\
    &=\frac{i}{2\pi}\int_0^{\frac{3}{4}}\Big( \log\Big(e^{-2i\mu \pi}\frac{1+t^{2\mu}}{1+e^{2i\mu \pi}t^{2\mu}}\Big)\frac{i}{it-\zeta}dt\\
     &=\frac{i}{2\pi}\log e^{-2i\mu \pi}\int_0^{\frac{3}{4}} \frac{i}{it-\zeta}dt+\frac{i}{2\pi}\int_0^{\frac{3}{4}}\log\Big(\frac{1+t^{2\mu}}{1+e^{2i\mu \pi}t^{2\mu}}\Big)\frac{i}{it-\zeta}dt\\
     &=\frac{i}{2\pi} \log e^{-2i\mu \pi}\log \Big(\frac{\frac{3}{4}+i\zeta}{i\zeta}\Big)+O(1)\\
     &=-\frac{i}{2\pi} \log e^{-2i\mu \pi}\log \zeta +O(1)\\
     &= \big(- \mu  +\lceil \tfrac{\pi}{\alpha}-\tfrac{3}{2}\rceil\big)\log \zeta+O(1).
    \end{align*}
    Inserting this into~\eqref{eq:I_alpha_AllPlaneExt1} completes the proof.
\end{proof}

\begin{lemma}
    For $0<\alpha<2\pi$, there is $r_\alpha>0$ such that 
\begin{equation}
\label{e:whatIsG}
g_{\alpha}(-i)=r_\alpha e^{i\frac{\pi}{4}(1-\frac{\pi}{\alpha})}.
\end{equation}
\end{lemma}
\begin{proof}
We compute using~\eqref{eq:I_alpha_AllPlaneExt3}
\begin{align*}
&\log g_\alpha(-i)\\
&=-\frac{i}{2\pi }\int_{\Gamma_{out,0}^{\pi/2-\alpha}}\log(1+(-iu)^{-2\mu})\frac{du}{u+ie^{-i\alpha/2}}-\frac{i}{2\pi}\int_{\Gamma_{in,0}^{-\pi/2-\alpha}}\log (1+(iu)^{-2\mu})\frac{du}{u+ie^{-i\alpha/2}}+\log 2\\
&=\frac{1}{2\pi i }\int_{0}^\infty\log(1+t^{-2\mu})\Big(\frac{1}{t+e^{i\alpha/2}}-\frac{1}{t-e^{i\alpha/2}}\Big)dt+\log 2\\
&=\frac{i}{\pi  }\int_{0}^\infty\log(1+t^{-2\pi/\alpha})\Big(\frac{e^{i\alpha/2}}{t^2-e^{i\alpha}}\Big)dt+\log 2\\
&=\frac{i}{\pi  }\int_{0}^\infty\Big(\frac{2\pi}{\alpha}\log t+\log(1+t^{-2\pi/\alpha})\Big)\Big(\frac{e^{-i\alpha/2}}{t^2-e^{-i\alpha}}\Big)dt+\log 2.
\end{align*}
Therefore, 
\begin{align*}
2\Im \log g_{\alpha}(-i)&=\frac{1}{\pi  }\int_{0}^\infty\log(1+t^{-2\pi/\alpha})\Big(\frac{e^{i\alpha/2}}{t^2-e^{i\alpha}}\Big)dt-\frac{1}{\pi  }\int_{0}^\infty\log(1+t^{-2\pi/\alpha})\Big(\frac{e^{-i\alpha/2}}{t^2-e^{-i\alpha}}\Big)dt\\
&=\frac{2}{\alpha}\int_0^\infty \log t\frac{e^{-i\alpha/2}}{t^2-e^{-i\alpha}}dt\\
&=\frac{2}{\alpha} \int_{e^{i(\pi/2-\alpha/2)}(0,\infty)} \log u\frac{e^{-i\alpha/2}}{u^2-e^{-i\alpha}}du\\
&=\frac{2}{\alpha}\int_0^\infty \big((\frac{\alpha-\pi}{2})-i\log t\big) \frac{1}{t^2+1}dt\\
&=-\frac{\pi(\alpha-\pi)}{2\alpha}.
\end{align*}
\end{proof}

\subsection{Verification that $v_d$ and $v_n$ satisfy~\eqref{e:stekDir} and~\eqref{e:stekNeu}}.

\subsubsection{Harmonicity}

By Lemma~\ref{l:Ismall} and the fact that for $\zeta$ on the infinite part of the contour $P$, $\arg \zeta = \pi +\alpha/4$ and hence, since $\alpha<2\pi$, for  $-\alpha/2\leq \arg z\leq 0$, 
$$
\Re (z\zeta) <-c_{\alpha}|z||\zeta|.
$$

Thus, the integrand defining $f_n$ and $f_d$ is rapidly decaying in $\zeta$ for all $-\alpha/2\leq \arg z\leq 0$, $z\neq 0$ and hence we may differentiate under the integral to see that $f_d$ and $f_n$ are holomorphic in $\mathcal{W}_{\alpha/2}$. In particular, 
$$
-\Delta v_d=-\Delta v_n=0,\text{ in }\mathcal{W}_{\alpha/2}. 
$$

\subsubsection{Steklov boundary conditions}
to obtain that for $z\in(0,\infty)$,
$$
\begin{gathered}
(\partial_y-1)v_d(z)=\Re (if_d'(z)-f_d(z))=-\Re\Big(\int_Pg_\alpha(\zeta)\zeta^{-\mu_\alpha}e^{z\zeta}d\zeta\Big),\\
(\partial_y-1)v_n(z)=\Re(if_n'(z)-f_n(z))=-\Re\Big(\int_Pg_\alpha(\zeta)e^{z\zeta}d\zeta\Big).
\end{gathered}
$$

By~\cite[A.1]{Pe:50}
$$
I(\zeta e^{-i\alpha/2})+\log\Big(\frac{\zeta+i}{\zeta}\Big)>0    
$$
for $\zeta>0$. Hence, $I_\alpha(\bar{\zeta})=\overline{I_\alpha(\zeta)}$ and, using that $I_\alpha$ is analytic on the 
$$
\{|z|>1\}\cap\mathcal{C},
$$
where $\mathcal{C}$ denotes the logarithmic cover of $\mathbb{C}\setminus 0$, $I_{\alpha}(\zeta e^{-i\zeta})$ is bounded on $|\zeta|>2$, and $|e^{z\zeta}|\leq e^{-cz}$ for $z>0$ and $-3\pi/2<\arg \zeta<3\pi/2$, may deform the contour $P$ to
$$
P_1:=\Gamma_{in,2}^{-\pi}\cup\{|z|=2\}\cup \Gamma_{out,2}^{\pi},
$$
and use that 
$$
g_{\alpha}(\bar{\zeta})=\overline{g_{\alpha}(\zeta)}
$$
to see that both integrals are pure imaginary and hence that 
$$
\partial_yv_d-v_d=\partial_yv_n-v_n=0,\quad\text{ on }z>0.
$$

\subsubsection{Neumann boundary conditions}
For the Neumann boundary condition, we change variables to $w=\zeta e^{-i\alpha/2}$ and deform the contour
\begin{align*}
\partial_{\nu_{\alpha/2}}v_n(se^{-i\alpha/2})&=\Re ie^{-i\alpha/2}f'(se^{-i\alpha/2})\\
&=\Re\Big(-e^{-i\alpha/2}\int_P\frac{g_{\alpha}(\zeta)\zeta }{\zeta+i}e^{se^{-i\alpha/2}\zeta}d\zeta \Big)\\
&=\Re\Big(-e^{-i\alpha/2}\int_P\exp(-I_{\alpha}(\zeta e^{-i\alpha/2})) e^{se^{-i\alpha/2}\zeta}d\zeta \Big)\\
&=\Re\Big(-\int_P\exp(-I_{\alpha}(w)) e^{sw}dw \Big)\\
&=\Re\Big(-\int_{P_1}\exp(-I_{\alpha}(w)) e^{sw}dw \Big).
\end{align*}
 Using that $I_{\alpha}(\bar\zeta)=\overline{I_{\alpha}(\zeta)}$ shows that the integral is pure imaginary and hence that $v_n$ satisfies the required boundary conditions.

 \subsubsection{Dirichlet boundary conditions}

 Finally, for the Dirichlet boundary condition, we again change variables to $w=\zeta e^{-i\alpha/2}$ and deform the contour to
\begin{align*}
v_d(se^{-i\alpha/2})&=\Re f(se^{-i\alpha/2})\\
&=\Re\Big(i\int_P\frac{g_{\alpha}(\zeta) }{\zeta+i}\zeta^{-\mu_{\alpha}}e^{se^{-i\alpha/2}\zeta}d\zeta \Big)\\
&=\Re\Big(\int_{P_1}\exp(-I_{\alpha}(w))(w)^{-\mu_{\alpha}-1} e^{sw}dw \Big)=0
\end{align*}
Using again that the integrand satisfies $I_{\alpha}(\bar{w})=\overline{I_{\alpha}(w)}$.

\subsubsection{Boundedness and asymptotics $f_d$ and $f_n$}

\begin{lemma}
\label{l:fasymptotic}
    Let $q\leq \mu$ and define
    $$
    f_q(z):=i\int_P\frac{g_{\alpha}(\zeta)}{\zeta+i}\zeta^{-q}e^{z\zeta}d\zeta.
    $$
    Then
    $$   \|f_q\|_{L^\infty(\mathcal{W}_{\alpha/2})}<\infty,
    $$
    and
    $$
    \lim_{x\to \infty}|f_q(x)+2\pi g_{\alpha}(-i)(-i)^{-q}e^{-ix}|=0.
    $$
\end{lemma}
\begin{proof}
For $-\alpha/2\leq \arg z\leq 0$, using Lemmas~\ref{l:Ismall} and~\ref{l:INear0}, we may deform the contour $P$ to  
$$
P_0:=\Gamma^{\pi+\theta}_{in,0}\cup \Gamma^{-\pi+\theta}_{out,0},
$$
where $|\theta-\alpha/4|$ is small enough so that $e^{z\zeta}$ is exponentially decreasing on $P_0$ for all $z\in\mathcal{W}_{\alpha/2}$ and $\{e^{i(\pi+\theta)},e^{i(-\pi+\theta)}\}\notin\{\lambda_k\}_{k=-\infty,\infty}$. 

Accounting for the poles of $g_{\alpha}$, we obtain
\begin{align*}
f_q(z)&= i\int_{P_0}\frac{g_{\alpha}(\zeta)}{\zeta+i}\zeta^{-q}e^{z\zeta}d\zeta-2\pi g_{\alpha}(-i)(-i)^{-q}e^{-iz}\\
&\qquad-2\pi\sum_{-\pi+\theta<\arg \lambda_k<\pi+\theta} e^{z\lambda_m} \operatorname{Res}_{\lambda_k}g_{\alpha}(\zeta)\zeta^{-q}/(\zeta+i)\\
&:= R_1(z)-2\pi g_{\alpha}(i)(i)^{-q}e^{iz}+R_2(z)\\
\end{align*}
It is easy to check that there are finitely many $k$ with 
$$
-\pi+\theta <\arg \lambda_k<\pi+\theta,
$$
and moreover, for such $k$,
$z\in\overline{\mathcal{W}_{\alpha/2}}$, $\Re z\lambda_k<-c|z|$.  This implies 
$$
|R_2(z)|\leq Ce^{-c|z|},\qquad z\in\mathcal{W}_{\alpha/2}.
$$

Using Lemma~\ref{l:INear0}, it is also easy to check that for $|q|\leq \mu$, 
$$
|R_1(z)|\leq C,\qquad z\in\overline{\mathcal{W}_{\alpha/2}}. 
$$
Moreover, for $x>0$ and $\zeta\in P_0$, $x\mapsto |e^{x\zeta}|$ is decreasing and $\lim_{x\to \infty}e^{x\zeta}=0$ for almost every $\zeta\in P_0$. In particular, the dominated convergence theorem implies that 
$$
\lim_{x\to \infty} |R_1(x)|=0. 
$$
this completes the proof of the lemma.
\end{proof}

\subsection{Determining the scattering matrix}

\begin{proof}[Proof of Lemma~\ref{l:scatMatrix}]
To determine the scattering matrix in $\mathcal{W}_{\alpha}=\{-\alpha<\arg z<0\}$, we set 
\begin{align*}
v_1(z)&:=\begin{cases}
    v_d(z)& -\alpha/2\leq \arg z\leq 0\\
    -v_d(ze^{2(-\arg z+\alpha/2)})&-\alpha\leq \arg z\leq \alpha/2,
\end{cases}\\
v_2(z)&:=\begin{cases}
    v_n(z)& -\alpha/2\leq \arg z\leq 0\\
    v_n(ze^{2(-\arg z+\alpha/2)})&-\alpha\leq \arg z\leq \alpha/2.
\end{cases}
\end{align*}
Then, since $v_d$ and $v_n$ satisfy respectively~\eqref{e:stekDir} and~\eqref{e:stekNeu}, we for $i=1,2$
$$
-\Delta v_i=0\, \text{in }\mathcal{W}_{\alpha},\qquad \partial_{\nu}v_i-v_i=0\, \text{on }\partial \mathcal{W}_{\alpha}.
$$
Moreover, since $v_i\in L^\infty$, 
$$
\partial_{\nu}v_i|_{\mathcal{W}_{\alpha}}=\Lambda_{\mathcal{W}_{\alpha}}v_i|_{\partial\mathcal{W}_{\alpha}}.
$$
Hence, setting $u_i:=v_i|_{\mathcal{W}_\alpha}$, we have
$$
(\Lambda_{\mathcal{W}_{\alpha}}-1)u_i=0.
$$
Therefore, to compute the scattering matrix, we only need to use Lemma~\ref{l:fasymptotic}. Indeed, 
$$
\lim_{x\to \infty}| v_n(x)+2\pi \Re (g_\alpha(-i))e^{-ix})|=0.
$$
Notice that by~\eqref{e:whatIsG}
$$
\Re g_{\alpha}(-i)e^{-ix}= \tfrac{r_\alpha}{2}(e^{i(\frac{\pi}{4}(1-\frac{\pi}{\alpha})-x)}+e^{-i(\frac{\pi}{4}(1-\frac{\pi}{\alpha})-x)})
$$
Therefore, 
$$
\Sigma^{\mathcal{W}_{\alpha}}_{0,1} \begin{pmatrix} 1\\1\end{pmatrix}  =-ie^{i\frac{\pi^2}{2\alpha}}\begin{pmatrix} 1\\1\end{pmatrix} .
$$
A similar computation using~\eqref{e:whatIsG} again shows that
$$
\Re g_{\alpha}(-i)(-i)^{-\mu}e^{-ix}= \tfrac{r_\alpha}{2}(e^{i(\frac{\pi}{4}(1+\frac{\pi}{\alpha})-x}+e^{-i(\frac{\pi}{4}(1+\frac{\pi}{\alpha})-x)}).
$$
Thus 
$$
\Sigma^{\mathcal{W}_{\alpha}}_{0,1} \begin{pmatrix}1\\-1\end{pmatrix}=-ie^{-i\frac{\pi^2}{2\alpha}}\begin{pmatrix}1\\-1\end{pmatrix}
$$

In particular,
$$
\Sigma^{\mathcal{W}_{\alpha}}_{0,1} =\begin{pmatrix}
    -i\cos(\frac{\pi^2}{2\alpha})&\sin(\frac{\pi^2}{2\alpha})\\ \sin(\frac{\pi^2}{2\alpha})&-i\cos(\frac{\pi^2}{2\alpha}).
\end{pmatrix}
$$
as claimed.
\end{proof}
\begin{lemma}
\label{l:scatMatrix}
The scattering matrix on the infinite wedge is given by
\begin{equation}
    \Sigma^{\mathcal{W}_{\alpha}}_{0,1} =\begin{pmatrix}
    -i\cos(\frac{\pi^2}{2\alpha})&\sin(\frac{\pi^2}{2\alpha})\\ \sin(\frac{\pi^2}{2\alpha})&-i\cos(\frac{\pi^2}{2\alpha}).
\end{pmatrix}
\end{equation}
\end{lemma}

\bibliography{biblio}
 \bibliographystyle{alpha}
\end{document}